\documentclass[12pt]{amsart}

\newif\ifdebug
\debugfalse

\newif \iffig
\figtrue

\newif \iftable
\tablefalse

\renewcommand{\bold}{\mathbf}

\usepackage{main_text_package_settings}
\usepackage{main_text_abbreviations}
\usepackage{main_text_thm_style_long}
\usepackage{pdfsync}

\newcommand{\tIP}[1]{{\left\langle\kern-0.40ex\left\langle\kern-0.40ex\left\langle #1 
    \right\rangle\kern-0.40ex\right\rangle\kern-0.40ex\right\rangle}}

\def\liessr{\Lie (\bbG^{ssr})}
\def\gB{\frakB}
\def\cross{\times}
\def\reals{\R}

\def\bfj{\bbj}

\renewcommand{\bH}{\bold H}

\title[Measure rigidity via Factorization]{Measure rigidity for generalized u-Gibbs states and stationary measures via the factorization method}

\hypersetup{
	pdftitle={Measure rigidity for generalized u-Gibbs states and stationary measures via the factorization method},	
	pdfauthor={Aaron Brown, Alex Eskin, Simion Filip, Federico Rodriguez Hertz},	
	pdfsubject={measure rigidity},		
	pdfkeywords={MSC},	
	pdfnewwindow=true,		
	colorlinks=true,		
	linkcolor=darkblue,		
	citecolor=darkred,		
	filecolor=darkblue,		
	urlcolor=darkblue,		
	pdfborder={0 0 0},
	breaklinks=true
}

\usepackage{imakeidx}
\makeindex[title = Index of Notation, columns = 3]

\thanks{Revised \textsc{\today} }

\date{February 2025}

\author[A.\ Brown]{
	Aaron Brown
}
\address{
	\parbox{0.5\textwidth}{
		Department of Mathematics\\
		Northwestern University\\}
}
\email{{awb@math.northwestern.edu}}

\author[A.\ Eskin]{
	Alex Eskin
}
\address{
	\parbox{0.5\textwidth}{
		Department of Mathematics\\
		University of Chicago\\
		5734 S University Ave\\
		Chicago, IL 60637\\}
}
\email{{eskin@math.uchicago.edu}}

\author[S.\ Filip]{
	Simion Filip
}
\address{
	\parbox{0.5\textwidth}{
		Department of Mathematics\\
		University of Chicago\\
		5734 S University Ave\\
		Chicago, IL 60637\\}
}
\email{{sfilip@math.uchicago.edu}}

\author[F.\ Rodriguez Hertz]{
	Federico Rodriguez Hertz
}
\address{
	\parbox{0.5\textwidth}{
		Department of Mathematics\\
		Pennsylvania State University\\
		328 McAllister Building\\
		State College, PA 16802\\}
}
\email{{fjr11@psu.edu}}
\newcommand{\td}{\tilde}
\newcommand{\wtd}{\widetilde}
\newcommand{\inv}{^{-1}}
\newcommand{\R}{\mathbb {R}}

\newcommand{\Z}{\mathbb {Z}}
\newcommand{\N}{\mathbb {N}}

\usepackage{enumitem}
\newtheorem{assumption}[theorem]{Assumption} 

\begin{document}

\begin{abstract}
We obtain measure rigidity results for stationary measures of random walks generated by diffeomorphisms, and for actions of $\SL(2,\R)$ on smooth manifolds.
Our main technical result, from which the rest of
the theorems are derived, applies also to the
case of a single diffeomorphism or $1$-parameter flow and establishes extra invariance of a class of measures that we call ``generalized u-Gibbs states''.
\end{abstract}

%

\maketitle
%

\ifdebug
   \listoffixmes
\fi

\tableofcontents



\section{Introduction}
	\label{sec:introduction}

In this paper, we study random dynamics and actions of $\SL(2,\R)$ on
smooth manifolds.
Our main results are stated in
\S\ref{sec:subsec:the:main:theorems}. We note that our main technical
result, \autoref{thm:inductive:step}, from which the rest of
the theorems are derived, applies also to the
situation where we are working with a single diffeomorphism or flow,
and may be of independent
interest. \autoref{thm:inductive:step} may be viewed as a
generalization of the main theorem of
\cite{Katz_Measure-rigidity-of-Anosov-flows-via-the-factorization-method}. 

As a special case, we prove a generalization of the results of
\cite{EskinMirzakhani_Invariant-and-stationary-measures-for-the-rm-SL2Bbb-R-action-on-moduli-space}
and \cite{EskinMirzakhaniMohammadi_Isolation} to actions on products
of moduli spaces of translation surfaces.

In this introductory section, we list some
easy to state consequences of the results of
\S\ref{sec:subsec:the:main:theorems}.

\subsection{Random Dynamics}
Let \index{$Q$}$Q$ be a connected smooth manifold.
Let \index{$D$@$\Diff_\infty(Q)$}$\Diff_\infty(Q)$ denote the group of $C^\infty$ diffeomorphisms of $Q$.
\begin{definition}\label{def:statandinvmsr}
Given a probability measure \index{$\mu$}$\mu$ on $\Diff_\infty(Q)$, a probability measure $\nu$ on $Q$ is said to be 
\begin{enumerate}
\item \emph{$\mu$-stationary} if $\int f_* \nu \, d \mu( f) = \nu$, and 
\item \emph{$\mu$-invariant } if $f_* \nu  = \nu$ for $\mu$-a.e.\ $f\in \Diff_\infty(Q)$.
\end{enumerate}
We observe the collection of $\mu$-stationary measures is a closed convex set.  
\begin{enumerate}[resume]
\item A $\mu$-stationary probability measure is said to be \emph{ergodic} if it is an extreme point of this collection.  
\end{enumerate}
Finally, 
\begin{enumerate}[resume]
\item we  say that \emph{$\mu$ is volume-preserving} if there exists a smooth volume form on $Q$ preserved  that is $\mu$-invariant.  
\end{enumerate}

\end{definition}
In this
version of the paper, we will be assuming that $\mu$ is finitely
supported. This assumption will be removed in a later version.

For $x \in Q$ let \index{$TQ(x)$}$TQ(x)$ denote the tangent space to $Q$.
Let \index{$G$@$\Gr_d(x)$}$\Gr_d(x)$ denote the Grassmannian of $d$-dimensional subspaces of $TQ(x)$.
Suppose that for each $x \in Q$ we are given a norm
\index{$\| \cdot \|_x$}$\| \cdot \|_x$ on the exterior power $\bigwedge^d(TQ(x))$. We then have a cocycle \index{$\sigma_d$}$\sigma_d \colon \Diff_\infty(Q) \cross TQ \to \reals$, given by
\begin{equation}
\label{eq:def:sigma:f:L}
  \sigma_d(f,L) = \log \frac{\| Df(l_1 \wedge \dots l_d) \|_{f(x)}}{\|l_1 \wedge \dots l_d \|_x}
\end{equation}
where $\{l_1, \dots, l_d\}$ is any basis for $L$. Then $\sigma_d(f,L)$ measures the change in the volume of the subspace $L$ under the diffeomorphism $f$. 


\begin{definition}[Uniform expansion]
\label{def:uniformly_expanding}
The measure $\mu$ is \emph{uniformly expanding in dimension $d$} if 
if there exist $C_d > 0$ and $N_d \in
\bN$ such that for all $x \in Q$ and all $L \in \Gr_d(x)$,
\begin{equation}
\label{eq:def:uniform:expansion}
	\int_{\Diff_\infty(Q)} \sigma_d(f,L) \, d\mu^{(N_d)}(f) > C_d > 0
\end{equation}
where \index{$\mu^{(N_d)}$}$\mu^{(N_d)}$ denotes the $N_d$-fold convolution of $\mu$ with itself.
\end{definition}

\begin{theorem}
\label{theorem:intro:stationary:is:invariant}
Suppose $\mu$ is volume preserving and satisfies uniform expansion in dimension $d$ for all $1 \le d < \dim(Q)$. Suppose also that either $Q$ is compact, or there exists an admissible Margulis function on $Q$ (see \autoref{def:random:margulis:function}). Then, any $\mu$-stationary measure is $\mu$-invariant.
\end{theorem}
A more technical result, but without assuming that $\mu$ is volume-preserving, is contained in \autoref{thm:random:inductive:step}.

\begin{definition}[Uniform expansion gaps]
	\label{def:uniform_expansion_gaps}
	The measure $\mu$ has a \emph{uniform expansion gap in dimension $d$} if 
	if there exist $C_d > 0, \delta_d\in\{\pm 1\}$, and  $N_d \in
	\bN$ such that for all $x \in Q$ and for all $L_0\subset L_1$ with $L_i \in \Gr_{d+i}(x)$, we have
	\begin{equation}
	\label{eq:def:uniform:expansion:gaps}
	\delta_d\cdot {\int_{\Diff_\infty(Q)} \big[\sigma_{d+1}(f,L_1)-\sigma_d(f,L_0)\big] \, d\mu^{(N_d)}(f) } > C_d > 0
	\end{equation}
	where $\mu^{(N_d)}$ denotes the $N_d$-fold convolution of $\mu$ with itself.
\end{definition}

\begin{remark}[On uniform expansion and gaps]
	\label{rmk:on_uniform_expansion_and_gaps}
	Suppose that $\nu$ is an ergodic $\mu$-stationary measure with Lyapunov spectrum $\lambda_1\geq \dots\geq \lambda_{\dim Q}$ on $TQ$.
	If $\mu$ is uniformly expanding in dimension $d$, then it follows from the same arguments as in the proof of \cite[Thm.~C]{Chung2020_Stationary-measures-and-orbit-closures-of-uniformly-expanding-random}, that the Lyapunov exponents satisfy $\lambda_1+\dots +\lambda_d>C_d/N_d$.
	Analogously, if $\mu$ has uniform expansion gap in dimension $d$ then $\delta_d \lambda_d>C_d/N_d$ where $\delta_d\in \{\pm 1\}$.
	
	Note also that both \autoref{def:uniformly_expanding} and \autoref{def:uniform_expansion_gaps} are properties that are open upon a $C^1$-perturbation.
\end{remark}

\begin{theorem}
\label{theorem:intro:ue:measure:classification}
Suppose $\mu$ is volume preserving and satisfies uniform expansion in dimension $d$ for all $1 \le d < \dim(Q)$.
Suppose \index{$\nu$}$\nu$ is an ergodic $\mu$-invariant measure on $Q$, and suppose that $\nu$ has no zero Lyapunov exponents.
Furthermore, suppose that either $Q$ is compact, or that there exist admissible Margulis functions on $Q$ and $Q\cross Q$ (see \autoref{def:random:margulis:function} and \autoref{def:random:QxQ:margulis:function}). 

Then $\nu$ is either absolutely continuous with respect to volume or is finitely supported.
If, in addition, $\mu$ has uniform expansion gap in dimension $d$ for all $1\leq d < \dim Q$, then the volume measure is $\mu$-ergodic and therefore $\nu$ is either volume or finitely supported.
\end{theorem}

\begin{definition}[$\mu$-invariant subbundle]
	\label{def:mu_invariant_subbundle}
	Suppose $\mu$ is a probability measure on $\Diff_{\infty}(Q)$ and $\nu$ is a $\mu$-stationary ergodic measure.
	A measurable subbundle $S\subset TQ$, defined $\nu$-a.e., is called \emph{$\mu$-invariant} if for $\mu$-a.e. $g$ and $\nu$-a.e. $q$ we have that $Dg(S(q))=S(g q)$.
\end{definition}
We emphasize that $S$ is defined on $Q$ and does not depend on the trajectory of the random walk.

\autoref{theorem:intro:stationary:is:invariant} is derived from
the following:
\begin{theorem}
\label{theorem:random:stationary:is:invariant}
Let $\mu$ be a finitely supported probability measure on $\Diff_\infty(Q)$. 
Let $\nu$ be an ergodic $\mu$-stationary measure on $Q$.
Suppose there is no
$\mu$-invariant $\nu$-measurable subbundle of $TQ$ on which the sum of
the Lyapunov exponents (see \autoref{ssec:random:skew:product}) is
negative. Then, $\nu$ is $\mu$-invariant.
\end{theorem}

Our main result for random actions is
\autoref{theorem:random:measure:classification}, from which
\autoref{theorem:intro:ue:measure:classification} is derived. 
The assumptions in \autoref{theorem:random:measure:classification}
are the same as in
\autoref{theorem:random:stationary:is:invariant}, and the
conclusion is a structure theorem for the stationary (and thus
invariant) measures, including a precise description of the
conditional measures along the stables and unstables.

We also state the following:
\begin{conjecture}
\label{conj:symplectic}
Suppose $Q$ is a compact symplectic manifold, and suppose $\mu$ is a probability measure on the group of symplectic diffeomorphisms of $Q$ which satisfies uniform expansion on all the isotropic subspaces of $TQ$.
Suppose $\nu$ is a $\mu$-stationary measure on $Q$ which has no $0$
Lyapunov exponents.  Then $\nu$ is a $\mu$-invariant measure on a smooth submanifold $M$ of $Q$. (The case $M=Q$ is allowed). 
\end{conjecture}

Our progress towards Conjecture~\ref{conj:symplectic} is given in
\autoref{theorem:symplectic}.

\subsection{\texorpdfstring{$\SL(2,\R)$}{SL(2,R)} actions}

Suppose that $\SL(2,\R)$ acts smoothly and locally freely on $Q$.
Let $\index{$g_t$}g_t = \begin{pmatrix} e^{t} & 0 \\ 0 & e^{-t} \end{pmatrix}$, and let
\begin{displaymath}
\index{$A$}A = \{ g_t  \in \R \}, \qquad \index{$N$}N = \begin{pmatrix} 1 & \ast \\ 0 &
  1 \end{pmatrix}, \qquad \index{$N$@$\bar{N}$}\bar{N} = \begin{pmatrix} 1 & 0 \\ \ast &
  1 \end{pmatrix},
\end{displaymath}
and let $\index{$P$}P=AN$.

Let \index{$SO(2)$}$SO(2)$ be a maximal compact subgroup of $\SL(2,\R)$.
Recall that the quotient $\SL(2,\R)/SO(2)$ is identified with the hyperbolic plane \index{$H$@$\mathbb{H}^2$}$\mathbb{H}^2$.
Let \index{$\Delta$}$\Delta$ denote the Laplace operator for the hyperbolic metric.

The analogue of uniform expansion is the following:
\begin{definition}[Nonnegative Laplacian]
\label{def:SL2R:nonnegative:laplacian}
Suppose $E\to Q$ is a cocycle over a smooth $\SL(2,\bR)$-action.
Given $d\in\lbrace 1,\dots,\dim E\rbrace$, we say that $E$ satisfies the \emph{non-negative Laplacian condition in rank $d$} if there exists a family of $SO(2)$-invariant norms
\index{$\| \cdot \|_q$}$\| \cdot \|_q$ on $\bigwedge^d E$ (for all $q \in Q$)
such that the following holds:

For any $q \in Q$, fix $\{v_1, \dots, v_d \}$ a linearly independent set in $E$.
Let $f\colon \SL(2,\R) \to \R$ be defined by
\begin{displaymath}
f(g) =  \log \| g v_1 \wedge \dots \wedge g v_d \|_{gq}. 
\end{displaymath}
Since $\| \cdot \|_q$ is $SO(2)$-invariant, $f$ descends to a function $h\colon \mathbb{H}^2 \to \R$. Then,
\begin{equation}
\label{eq:nonnegative:laplacian}
\Delta h \ge 0.   
\end{equation}
\end{definition}

The following is an analogue of \autoref{theorem:intro:stationary:is:invariant}:
\begin{theorem}
\label{theorem:intro:SL2R:P:invanriant:is:SL2R:invariant}
Suppose $SL(2,\reals)$ acts on $Q$ and satisfies the non-negative Laplacian condition on $TQ$, for all ranks $d\in\lbrace 1, \dots,\dim Q\rbrace$ (\autoref{def:SL2R:nonnegative:laplacian}).
Suppose $\nu$ is a $P$-invariant ergodic measure on $Q$.

Then, $\nu$ is $\SL(2,\R)$-invariant. 
\end{theorem}



\autoref{theorem:intro:SL2R:P:invanriant:is:SL2R:invariant} is
derived from the following:
\begin{theorem}
\label{theorem:SL2R:P:invariant:is:SL2R:invariant}
Suppose $SL(2,\reals)$ acts on $Q$ and 
let $\nu$ be an ergodic $P$-invariant measure on $Q$.
Suppose there is no $P$-invariant $\nu$-measurable subbundle of $TQ$
on which the sum of the Lyapunov exponents of $g_1$ is negative. Then, $\nu$ is $\SL(2,\R)$-invariant. 
\end{theorem}

Our main result for general $\SL(2,\R)$-actions is
\autoref{theorem:SL2R:measure:classification}.
The assumptions in \autoref{theorem:SL2R:measure:classification}
are the same as in
\autoref{theorem:SL2R:P:invariant:is:SL2R:invariant}, and the
conclusion is a structure theorem for the $P$-invariant (and thus
$\SL(2,\R)$-invariant) measures, including a precise description of the
conditional measures along the stable and unstable manifolds of $g_t$.


\subsection{Teichm\"uller dynamics}
	\label{ssec:teichmuller_dynamics_statements}

Suppose $g \ge 1$, \index{$g$} and 
let $\alpha = \alpha_1+ \dots + \alpha_n$ \index{$\alpha$} 
be a partition of $2g-2$, and 
let $\cH(\alpha)$ \index{$H($@$\cH(\alpha)$} be a stratum of Abelian differentials,
i.e.\ the space of pairs $(M,\omega)$ where $M$ \index{$M$} is a Riemann surface
and $\omega$ \index{$\omega$} is a holomorphic $1$-form on $M$ whose
zeroes have
multiplicities $\alpha_1 \dots \alpha_n$. The form $\omega$ defines a
canonical flat metric on $M$ with conical singularities at the zeros
of $\omega$. Thus we refer to points of $\cH(\alpha)$ as
{\em flat surfaces} or {\em translation surfaces}. For an introduction
to this subject, see the surveys \cite{Zorich2006_Flat-surfaces,Filip2024_Translation-surfaces:-Dynamics-and-Hodge-theory}.

The space $\cH(\alpha)$
admits an action of the group $\SL(2,\reals)$ which generalizes the
action of $SL(2,\reals)$ on the space $GL(2,\reals)/SL(2,\Z)$ of flat
tori. In
\cite{EskinMirzakhani_Invariant-and-stationary-measures-for-the-rm-SL2Bbb-R-action-on-moduli-space}
and \cite{EskinMirzakhaniMohammadi_Isolation}
some ergodic-theoretic rigidity properties of this action were proved. In this paper, we generalize the results to the diagonal action of $\SL(2,\reals)$ on products of the $\cH(\alpha)$. 

As in
\cite{EskinMirzakhani_Invariant-and-stationary-measures-for-the-rm-SL2Bbb-R-action-on-moduli-space}, 
we always replace $\cH(\alpha)$ by a finite cover which is a
manifold (for example the congruence cover with $\bZ/3$-homology trivialized).
However,
we suppress this from the notation.


\medskip

\paragraph{Products of strata}
Let $\index{$\alpha$@$\vec{\alpha}$}\vec{\alpha}=(\alpha_1, \dots, \alpha_n)$. Let $\index{$H$@$\cH(\vec{\alpha})$}\cH(\vec{\alpha}) =
\cH(\alpha_1) \cross \cH(\alpha_n)$. The moduli space
$\cH(\vec{\alpha})$ parametrizes a holomorphic $1$-form on a
disconnected surface $\index{$M$}M = M_1 \cup \dots \cup M_n$, which is an
$n$-tuple $\index{$\omega$}\omega=(\omega_1,\dots, \omega_n)$ where $\omega_i$ is a
holomorphic $1$-form on $M_i$.

Let $\index{$\Sigma_i$}\Sigma_i$ denote the set of zeroes of $\index{$\omega_i$}\omega_i$, and let $\index{$\Sigma$}\Sigma
= \bigcup_i \Sigma_i$. We have
\begin{displaymath}
H^1(M,\Sigma;\mathbb{C}) = \bigoplus_{i=1}^n H^1(M_i,\Sigma_i, \mathbb{C}).
\end{displaymath}
Let $\index{$\Phi$}\Phi: \cH^1(\vec{\alpha}) \to H^1(M,\Sigma;\mathbb{C})$ be the
period coordinate map. Then, $\Phi = (\Phi_1, \dots, \Phi_k)$ where
$\index{$\Phi_j$}\Phi_j: \cH(\alpha_j) \to H^1(M_j,\Sigma_j,\mathbb{C})$ is period map of
$\cH(\alpha_j)$ which takes $(X_j,\omega_j) \in \cH(\alpha_j)$ to the
cohomology class of $\omega_j$ in $H^1(M_j,\Sigma_j,\mathbb{C})$.

The area of a translation surface is given by 
\begin{displaymath}
\index{$a_j$}a_j(M_j,\omega_j) = \frac{i}{2} \int_{M_j} \omega_j \wedge \bar{\omega}_j.
\end{displaymath}
A ``unit hyperboloid'' $\cH_1(\alpha)$ \index{$H$@$\cH_1(\alpha)$}
is defined as a subset of translation surfaces in $\cH(\alpha)$ of
area one. We then define
\begin{displaymath}
\index{$H$@$\cH_1(\vec{\alpha}) $}\cH_1(\vec{\alpha}) = \cH_1(\alpha_1) \cross \dots \cross \cH_1(\alpha_n).
\end{displaymath}

Let $\index{$\lambda$@$\tilde{\lambda}_j$}\tilde{\lambda}_j$ denote the Masur-Veech measure on $\cH(\alpha_j)$,
and let \index{$\lambda_j$}$\lambda_j$ denote the Masur-Veech measure on
$\cH_1(\alpha_j)$ which is the restriction of $\tilde{\lambda}_j$,
so that $d\tilde{\lambda}_j = c d\lambda_j \, da_j$
where $c$ is a constant. We normalize $\lambda_j$ so that
$\lambda_j(\cH_1(\alpha_j)) = 1$. Let $\tilde{\lambda} =
\tilde{\lambda}_1 \cross \dots \cross \dots \tilde{\lambda}_n$ and let
$\lambda = \lambda_1 \cross \dots \lambda_n$, so that $\lambda$ is a
probability measure on $\cH_1(\vec{\alpha})$. 

For a subset $\cM_1 \subset \cH_1(\vec{\alpha})$ we write\index{$R$@$\reals^n \cM_1$}
\begin{multline*}
\reals^n \cM_1 = \{ ((M_1, t_1 \omega_1), \dots, (M_n,t \omega_n))
\;|\; ((M_1,\omega_1), \dots, (M_n,\omega_n) \in \cM_1, \\ t_j \in
\reals \setminus \{0 \} \} \subset \cH(\vec{\alpha}).
\end{multline*}
\begin{definition}
\label{def:affine:measure}
An ergodic $SL(2,\reals)$-invariant probability measure $\nu$ on
$\cH_1(\alpha)$ is called {\em affine} if the following conditions hold:
\begin{itemize}
\item[{\rm (i)}] The support $\cM_1$ of $\nu$ is an 
{\em immersed submanifold} of
  $\cH_1(\alpha)$, i.e.\
there exists a manifold $\cN$ and a proper continuous
  map $f: \cN \to \cH_1(\alpha)$ so that $\cM_1 =
  f(\cN)$. The self-intersection set of $\cM_1$, i.e.\ the set of
  points of $\cM_1$ which do not have a 
  unique preimage under $f$, is a closed subset of $\cM_1$ of
  $\nu$-measure $0$.   Furthermore, each point in $\cN$ has a
  neighborhood $U$ such 
  that locally $\reals^n f(U)$ is given by a complex linear subspace defined
  over $\reals$ in the period coordinates.
\item[{\rm (ii)}] Let $\tilde{\nu}$ be the measure supported on $\cM = \reals^n
  \cM_1$ so that $d\tilde{\nu} = 
    d\nu_1 da_1 \dots da_n$. Then each point in $\cN$ has a
      neighborhood $U$ such that the restriction of $\tilde{\nu}$ to $\reals^n
      f(U)$ is an affine  linear measure in the period
    coordinates on $\reals^n f(U)$, i.e.\ it is (up to normalization) the
    induced measure of the Masur-Veech measure $\tilde{\lambda}$ to the subspace
    $\reals^n f(U)$.
\end{itemize}
\end{definition}

\begin{definition}
\label{def:affine:invariant:submanifold}
We say that any suborbifold $\cM_1$ for which there exists a measure
$\nu$ such that the pair $(\cM_1, \nu)$ 
satisfies (i) and (ii) is an {\em affine invariant submanifold}. 
\end{definition}
We also consider all of $\cH(\vec{\alpha})$ to be an
(improper) affine invariant submanifold. 
It follows from the analogue of
\cite[Theorem~2.2]{EskinMirzakhaniMohammadi_Isolation} that
the self-intersection set of an affine invariant manifold is itself a
finite union of 
affine invariant manifolds of lower dimension. 
\medskip

Let
\begin{displaymath}
N = \left\{ \begin{pmatrix} 1 & t \\ 0 & 1 \end{pmatrix}, t \in \reals
  \right\}, \quad A =
\left\{ \begin{pmatrix} e^t & 0 \\ 0 & e^{-t} \end{pmatrix}, t \in
  \reals \right\}, \quad \bar{N} = \left\{ \begin{pmatrix} 1 & 0 \\ t & 1
  \end{pmatrix}, t \in \reals   \right\}
\end{displaymath}
Let $r_\theta = \begin{pmatrix} \cos \theta & \sin \theta
  \\ - \sin \theta & \cos \theta \end{pmatrix}$, \index{$r_\theta$}
and let $SO(2) = \{
r_\theta \mid \theta \in [0,2\pi) \}$. \index{$SO(2)$}
Then $N$, $\bar{N}$, $A$ and $SO(2)$ 
are subgroups of $\SL(2,\reals)$. Let $P = AN$
denote the set of upper triangular matrices of determinant
$1$, which is a subgroup of $SL(2,\reals)$.

The following generalizes \cite[Theorem~1.4]{EskinMirzakhani_Invariant-and-stationary-measures-for-the-rm-SL2Bbb-R-action-on-moduli-space}:
\begin{theorem}
\label{theorem:Teich:P:measures}
Let $\nu$ be any ergodic $P$-invariant probability measure on
$\cH_1(\vec{\alpha})$. Then $\nu$ is $\SL(2,\reals)$-invariant and affine.  
\end{theorem}

The following generalizes
\cite[Theorem~2.1]{EskinMirzakhaniMohammadi_Isolation}:  
\begin{theorem}
\label{theorem:Teich:closure:submanifold}
Suppose $x \in \cH_1(\vec{\alpha})$. Then the orbit closure
$\overline{Px}$ is equal to $\overline{\SL(2,\R)x}$ and
is an affine invariant
submanifold of $\cH_1(\vec{\alpha})$. 
\end{theorem}

We note that the other statements from Sections 2 and 3 of
\cite{EskinMirzakhaniMohammadi_Isolation}
also hold in this setting.



\subsection{Outline and General Remarks}
	\label{ssec:outline_and_general_remarks}

All the main results are deduced from \autoref{thm:inductive:step}, whose proof is outlined in \autoref{ssec:outline_proof_inductive_step}.

\subsubsection{Previous results}
	\label{sssec:previous_results}
The work on measure and topological rigidity questions in dynamics is vast and to survey it is beyond our capacity.

For homogeneous spaces, some landmark contributions in the rigidity of unipotent flows have been made by 
Furstenberg \cite{Furstenberg1973_The-unique-ergodicity-of-the-horocycle-flow},
Margulis \cite{Margulis1987_Formes-quadratriques-indefinies-et-flots-unipotents-sur-les-espaces-homogenes}
and Ratner \cite{Ratner1991_On-Raghunathans-measure-conjecture}.
Rigidity of measures invariant under abelian groups, under a positive entropy assumption, has been established in the work of Einsiedler, Katok, and Lindenstrauss \cite{EinsiedlerKatokLindenstrauss2006_Invariant-measures-and-the-set-of-exceptions-to-Littlewoods-conjecture}.
Still in the setting of Lie groups and homogeneous spaces, Benoist--Quint \cite{BenoistQuint_Mesures1} obtained the first rigidity results for stationary measures (see also the work of Bourgain--Furman--Lindenstrauss--Moses \cite{BFLM} for a precursor in the case of tori, and Eskin--Lindenstrauss \cite{Eskin:Lindenstrauss:short} for a view on this setting closer to the spirit of this work).

Beyond the homogeneous setting, in the setting of Teichm\"{u}ller dynamics rigidity of measures invariant under the upper-triangular subgroup of $\SL_2(\bR)$ was obtained by Eskin, Mirzakhani, and Mohammadi \cite{EskinMirzakhani_Invariant-and-stationary-measures-for-the-rm-SL2Bbb-R-action-on-moduli-space,EskinMirzakhaniMohammadi_Isolation}.
For random walks generated by diffeomorphisms of real $2$-dimensional surfaces, Brown and Rodriguez Hertz obtained a measure rigidity result
\cite{BrownRodriguez-Hertz2017_Measure-rigidity-for-random-dynamics-on-surfaces-and-related-skew} that was used and extended to the case of complex $2$-dimensional surfaces by Cantat and Dujardin 
\cite{CantatDujardin2023_Dynamics-of-automorphism-groups-of-projective-surfaces:-classification-examples} (see also the work of Roda \cite{Roda2024_Classifying-hyperbolic-ergodic-stationary-measures-on-compact}).

The QNI condition (see \autoref{sssec:qni_condition_proof_outline}) appears in the work of Katz \cite{Katz_Measure-rigidity-of-Anosov-flows-via-the-factorization-method}, as well as Eskin--Potrie--Zhang 
\cite{EskinPotrieZhang2023_Geometric-properties-of-partially-hyperbolic-measures-and-applications}.
Uniqueness of u-Gibbs states has been recently considered in the works of Alvarez, Leguil, Obata, Santiago, and Cantarino \cite{AlvarezLeguilObata2022_Rigidity-of-U-Gibbs-measures-near-conservative-Anosov, CantarinoSantiago2024_U-Gibbs-measure-rigidity-for-partially-hyperbolic-endomorphisms-on-surfaces}.

\subsubsection{The role of normal forms}
	\label{sssec:the_role_of_normal_forms}
The unstable manifolds of the flow are equipped with normal form coordinates, as we explain in \autoref{appendix:normal_forms_on_manifolds_and_cocycles}.
We view these as a class of charts on the unstable manifolds, with transition maps in the group of subresonant maps.
Among the consequences of this structure are several additional properties and constructions.
The first is that the unstable manifold acquires the action of the group of subresonant, and strictly subresonant maps.
The second is that we can linearize the action on the unstable manifold by polynomially embedding the manifold into a linear cocycle.

\subsubsection{Subgroups compatible with the measure}
	\label{sssec:subgroups_compatible_with_the_measure_proof_outline}
The subresonant structure on the unstable manifolds allows us to define generalized u-Gibbs states, see \autoref{def:generalized:u:Gibbs}.
These are $g_t$-invariant measures whose conditionals along unstables are invariant under a family of subgroups \index{$U^+(x)$}$U^+(x)$ of the strictly subresonant maps.

\subsubsection{QNI condition}
	\label{sssec:qni_condition_proof_outline}
The crucial extra invariance of \autoref{thm:inductive:step} is deduced from the assumption that the family $U^+$ and the center-stable manifolds satisfy a non-integrability condition, as explicated in \autoref{def:QNI}.

\subsubsection{Some notational conventions}
	\label{sssec:some_notational_conventions}
We will adopt the convention from
\cite{EskinMirzakhani_Invariant-and-stationary-measures-for-the-rm-SL2Bbb-R-action-on-moduli-space}
that \index{$W$@$\cW^s[q]$}
\begin{align*}
	\cW^s[q] & \text{ is the stable manifold}\\
	\index{$W^s(q)$}W^s(q) & \text{ is the stable tangent space}
\end{align*}
and more generally, we'll use \index{$E(q)$}$E(q)$ for the vector space which is a fiber of a vector bundle at $q$ and \index{$E[q]$}$E[q]$ for an actual manifold, or ``something inside the manifold''.

\subsubsection*{Acknowledgments}
The preparation of this work started almost 11 years ago.
In the intervening years, we have had many conversations and exchanges with a broad group of colleagues.
We are grateful to all of them for their patience and interest.

This material is based upon work supported by the National Science
Foundation under Grants No. DMS-2400191 (AWB), DMS-2020013 (AWB), DMS-1900778 (FRH), DMS-2005470 (SF), DMS-2305394 (SF), and a Simons Investigator grant from the Simons Foundation (AE); this research was partially conducted during the period SF served as a Clay Research Fellow.




\section{Technical preliminaries}
	\label{sec:technical_preliminaries}



\subsection{Assumptions on the total space}
	\label{ssec:assumptions_on_the_total_space}




\subsubsection{Notation}
	\label{sssec:notation_introduction}
We work on a manifold \index{$Q$}$Q$, denote its tangent bundle by \index{$TQ$}$TQ$, and assume that $Q$ is equipped with a smooth Riemannian metric.
The tangent space at $q\in Q$ is denoted \index{$TQ(q)$}$TQ(q)$ and the ball of radius $r>0$ at $0\in TQ(q)$ is denoted \index{$B(q;r)$}$B(q;r)$.
Later on, we will introduce other norms on $TQ(q)$ and such balls will be denoted $B_{\bullet}(q;r)$ where $\bullet$ indicates the flavor of norm.

\subsubsection{Norms on jet bundles}
	\label{sssec:norms_on_jet_bundles}
Denote by \index{$J_kQ$}$J_kQ$ the bundle of $k$-jets of functions on $Q$ and its dual \index{$J_k^\dual Q$}$J_k^\dual Q$.
For example $J_0Q$ is the trivial $1$-dimensional bundle, and $T^\dual Q = \ker\left(J_1 Q \onto J_0 Q\right)$.

The Riemannian metric on $Q$ yields:
\begin{enumerate}
	\item Norms on the $k$-jets $J_kQ$, and induced norms on the duals $J_k^\dual Q$.
	The norms are compatible, i.e. the embedding $J_k^\dual Q\into J_{k+1}^\dual Q$ is an isometry.
	\item The norm on $TQ$ gives a distance function on $Q$, by minimization along paths.
\end{enumerate}
We will denote these ``ambient'' norms by \index{$\norm$@$\tnorm{\bullet}$}$\tnorm{\bullet}$, and if necessary indicate with a subscript the order of the jets.

\subsubsection{Charts with bounds on jets}
	\label{sssec:charts_with_bounds_on_jets}
We assume also that there exists an ``ambient radius'' function $\index{$r_a$}r_a\colon Q\to \bR_{>0}$, such that for every $q\in Q$ there exist charts, which we will call 
``exponential maps'',  
\[
	\index{$e$@$\exp_q$}\exp_q\colon B(q;r_a(q))\to Q \text{ where } B(q;r_a(q))\text{ is the $\tnorm{\bullet}$-ball of radius }r_a(q).
\]
We assume the following compatibility between these charts, and the background splittings of higher order jets:
There exists a function \index{$C_k(q)$}$C_k(q)$, such that for any $v\in B(q;r_a(q))$ and $q':=\exp_q(v)$ the induced map between higher order jets
\[
	J_k(TQ(q))(v) \xleftarrow{J_k\exp_q} J_kQ(q')
\]
has norm bounded by $C_k(q)$.
Note that the jet bundles over a vector space equipped with a norm are naturally split and have natural metrics induced from the metric of the vector space.
 
We will put further restrictions on $r_a(q)$ and $C_k(q)$ in \autoref{sssec:systole_function} below.

\begin{remark}[On distances and norms]
	\label{rmk:on_distances_and_norms}
	If we have a Riemannian metric on $Q$, its Levi-Civita connection provides the above splittings.

	For every $n$, there exists a constant $c_n$ such that for every Finsler norm $\norm{-}_F$ on an $n$-dimensional vector space, there exists a positive-definite inner product $\norm{-}_J$ (provided by the ``John ellipsoid'' of the unit ball of the Finsler metric) such that $\norm{-}_F\leq \norm{-}_J\leq c_n\norm{-}_F$.
	So on tangent spaces, we can freely pass from a Finsler metric to a positive-definite quadratic form.
\end{remark}

\subsubsection{The flow $g_t$.}
Let \index{$g_t$}$g_t$ be a smooth flow on the manifold $Q$.
If we are considering an $\SL(2,\R)$-action, then $g_t$ denotes the action of
$\begin{pmatrix} e^{t} & 0 \\ 0 & e^{-t} \end{pmatrix} \in \SL(2,\R)$.

\subsubsection{Systole function for flows}
	\label{sssec:systole_function}
We assume there exists an (ambient) \emph{inverse systole function} \index{$\phi_a(q)$}$\phi_a(q)\geq  2$, that we think of as $y$ on the upper half-plane, or the inverse of the distance to a singular set of a subvariety, such that $\log \phi_a$ is uniformly continuous.
This holds if, for instance, $\log \phi_a$ is Lipschitz.
In particular, we assume $\phi_a$ is continuous and therefore bounded on compact sets.

We also require $\phi_a$ to control the dynamics in the following manner:
\begin{description}
	\item[Bounded growth] there exists $C>0$ such that 
	\[
		\phi_a(g_t q) \leq e^{C|t|}\phi_a(q)
	\]
	which follows from the infinitesimal bound: $|\cL_{g_t}\phi_a|\leq C \phi_a$, abusing notation for the Lie derivative along the vector field generating $g_t$.
	\item[Bounds on $k$-jet dynamics] There exists a constant $\index{$e$}e>0$, and a ``dynamical radius'' $\index{$r_d(q)$}r_d(q):=\phi_a(q)^{-e}r_a(q)$, such that for all $t\in[-1,1]$, the transformation $g_t$ takes $\exp_q(B(q;r_d(q)))$ into the image of $\exp_{g_t q}(B(q;r_a(g_t q)))$.
	Conjugating $g_t$ by these charts to measure the $C^k$-jets, we have the inequality:
	\[
	 	\sup_{t\in[-1,1]} \norm{g_t}_{C^k\left(B(q;r_d(q))\right)} \leq \phi_a(q)^{a_k}
	\]
	for some fixed constant $a_k>0$.
	\item[Bounds on $k$-jet transitions]: We assume that the size of charts $r_a(q)$ from the previous part is bounded below by $\phi_a(q)^{-e_0}$.
	Similarly, the constants $C_k(q)$ satisfy the bound $|C_k(q)|\leq \phi_a(q)^{e_k}$ for some $\index{$e_k$}e_k>0$.
\end{description}
We will denote by
\begin{align}
	\label{eqn:big_G_t_definition}
	\index{$G_{t,q}$}G_{t,q} := \exp_{g_t q}^{-1} \circ g_t \circ \exp_q 
\end{align}
the map between subsets of $TQ(q)$ and $TQ(g_t q)$ which are to be specified depending on the context.

\subsubsection{Some special cases}
	\label{sssec:some_special_cases}
For example, if $Q$ is compact we can take $\phi_a\equiv 10$, i.e. any constant larger than $1$.
Then the constants $C,e,a_k,e_k$ in \autoref{sssec:systole_function} can be chosen sufficiently large depending on the smooth vector field generating $g_t$.

\subsubsection{An elementary bound}
	\label{sssec:an_elementary_bound}
We record, for future use, the following elementary fact.
For any compact set $K\subset Q$ and any $t_0,\ve_0>0$ there exists $C:=C(K,t_0,\ve_0)$ such that if $q\in K$ and $d^Q(q,q')\leq \ve_0$ then
\begin{align}
	\label{eqn:elementary_compactness_2_point_estimate}
	d^Q(g_tq,g_tq')\leq C \cdot d^Q(q,q') \quad \forall t\in[0,t_0].
\end{align}

\subsubsection{Systole function for random dynamics}
	\label{sssec:random:systole_function}
Let \index{$S$@$\cS_\mu$}$\cS_\mu$ denote the support of $\mu$. Recall that we are assuming
that $\cS_\mu$ is a finite set. 
We assume that there exists a function $\phi_a: Q
\to \R$  such that all of the conditions of
\autoref{sssec:systole_function} are satisfied, where in place of
$g_t$ we consider any word in $\cS_\mu \cup \cS_\mu^{-1}$ of length at
most $t$.



\subsection{Assumptions on the measure and consequences}
	\label{ssec:assumptions_on_the_measure}

Suppose that $(Q,g_t)$ satisfies the conditions from \autoref{sssec:systole_function} with a systole function $q\mapsto \phi_a(q)$.  

\begin{assumption}
Throughout, \index{$\nu$}$\nu$ will be an ergodic, $g_t$-invariant 
Borel probability measure satisfying the integrability assumption:
\begin{equation}
\label{eq:systole:log:integrable}
\int \log^+|\phi_a(q)|d\nu(q)<+\infty.  
\end{equation}
\end{assumption}

\subsubsection{Lyapunov exponents}
As the systole function $\phi_a(\cdot)$ controls the norm of the derivative cocycle, we may apply the Oseledets theorem.  
We enumerate the Lyapunov exponents of the derivative cocycle as \index{$\lambda_i$}
\[
	\lambda_1>  \lambda_2 > \dots
\]
including possible zeros.


\begin{proposition}[Adapted Lyapunov norms]
	\label{prop:rescaling_metrics_to_size_1}
	For any $k\geq 1,\ve>0$, there exists a measurable ``dynamical systole function'' \index{$\phi_{k,\ve}(q)$}$\phi_{k,\ve}(q)\geq 1$, and a measurable family of $\ve$-Lyapunov metrics \index{$\norm{\bullet}_{k,\ve}$}$\norm{\bullet}_{k,\ve}$ on the tangent spaces, both defined $\nu$-a.e., with the following properties:
	\begin{enumerate}
		\item \label{Lnorm1} Slow variation:
	 	\[
	 		|\phi_{k,\ve}(g_t q)| \leq e^{|t|\ve} |\phi_{k,\ve}(q)|
	 	\]
	 	\item \label{Lnorm2}  Comparison to ambient metric: we have $\nu$-a.e. that
	 	\[
	 	 	\tnorm{\bullet} \leq \norm{\bullet}_{k, \ve} \leq \phi_{k,\ve}\tnorm{\bullet}
	 	\]
	 	\item  \label{Lnorm3}  Uniform hyperbolicity: the Lyapunov decomposition of the tangent cocycle $TQ$
	 	\[
	 		TQ = \oplus_{\lambda_i} TQ^{\lambda_i}
	 	\]
	 	is orthogonal for $\norm{\bullet}_{k,\ve}$ and satisfies on a set of full measure:
	 	\[
	 		e^{-t\lambda_i}\norm{Dg_t\vert_{TQ^{\lambda_i}}}_{k,\ve}
	 		\leq e^{|t|\ve} \quad 
	 		\text{ for all } t\in \bR.
	 	\]
	 	\item  \label{Lnorm4}  Dynamics bounds: for a constant $r_0>0$ that only depends on the Lyapunov spectrum, we have that (see \autoref{eqn:big_G_t_definition})
	 	\[
	 		G_{t,q}:=\exp_{g_tq}^{-1} \circ g_t \circ \exp_q\colon B_{k,\ve}(q;r_0)
	 		\to
	 		B_{k,\ve}(g_tq;1)
	 	\]
	 	is well-defined for $t\in[-1, 1]$, where $B_{k,\ve}(q';r)$ denotes the ball of radius $r$ in the metric $\norm{\bullet}_{k,\ve}$ centered at $0$ in $T_{q'}Q$.
	 	Furthermore we have the bound on derivatives:
	 	\[
	 		\sup_{t\in[-1,1]}\norm{G_{t,q}
	 		- Dg_t
	 		}_{C^k(B_{k,\ve}(q;r_0))}\leq \ve
	 	\]
	 	and the $C^k$-norm is for the norms $\norm{\bullet}_{k,\ve}$ on the respective vector spaces.
	 \end{enumerate} 
\end{proposition}
\noindent The proof is given in \autoref{sssec:proof_of_prop:rescaling_metrics_to_size_1}.

\begin{remark}[Basic $C^k$ bounds]
	\label{rmk:basic_c_k_bounds}
From (iii) and (iv) we obtain the following: for every $k$ and $\epsilon>0$ there is \index{$C_k$}$C_k>1$ and $\index{$\kappa_k$}\kappa_k>0$ 
such that for all $T$ and $\mu$-a.e. $x$,
the map $$\exp_{g_tq}^{-1} \circ g_t \circ \exp_q\colon B_{k,\ve}(q;C_k^{-1}e^{-\kappa_k T})\to 
B_{k,\ve}(g_tq;1)$$
is well-defined and
\[
	\sup_{t\in[-T,T]}\norm{\exp_{g_tq}^{-1} \circ g_t \circ \exp_q}
	_{C^k(B_{k,\ve}(q;C^{-1}e^{-\kappa_k T}))}\leq C_ke^{\kappa_k T}.
\]
Furthermore, there exists a tempered function $\index{$C_k$@$\wtilde{C}_k$}\wtilde{C}_k$ such that
\[
	\sup_{t\in[-T,T]}\tnorm{\exp_{g_tq}^{-1} \circ g_t \circ \exp_q}
	_{C^k(B_{k,\ve}(q;C^{-1}e^{-\kappa_k T}))}
	\leq \wtilde{C}_k(q) e^{\kappa_k T}.
\]
\end{remark}

\begin{definition}[Measurably good smooth dynamics]
	\label{def:measurably_good_smooth_dynamics}
	A system satisfying the conclusions of \autoref{prop:rescaling_metrics_to_size_1} will be said to have \emph{measurable good smooth dynamics}, or we'll say it is \emph{measurably good}.
\end{definition}

General properties of measurably good dynamical systems, such as the
existence of stable and unstable manifolds are developed in
\autoref{ssec:lyapunov_charts_and_stable_manifolds}. We denote the
stable manifold through $q \in Q$ by \index{$W$@$\cW^s[q]$}$\cW^s[q]$, and its tangent space
by $\index{$W^s(q)$}W^s(q)$. Similarly, \index{$W$@$\cW^u[q]$}$\cW^u[q]$ and $\index{$W^u(q)$}W^u(q)$ denote the
unstable manifold and its tangent space.

Normal forms coordinates on
stables and unstables are developed in
\autoref{appendix:normal_forms_on_manifolds_and_cocycles}. 



\subsubsection{Lyapunov charts}
	\label{sssec:lyapunov_charts_from_introduction}
We define
\[
	\index{$L$@$\cL_{k,\ve}[q]$}\cL_{k,\ve}[q]:=\exp_q(B_{k,\ve}(q;r_0))\subset Q
\]
where $B_{k,\ve}(q;r_0)$ is as in \autoref{prop:rescaling_metrics_to_size_1}.

For many purposes, it will be enough to take $k=2$ and $\ve(2)>0$ sufficiently small and work in such charts.
Later on, in \autoref{sssec:choice_of_k_epsilon}, we will fix once and for all $k,\ve$ (which will depend on the QNI condition).
Until then \index{$L$@$\cL[q]$}$\cL[q]$ will refer to $\cL_{k,\ve}[q]$ for $(k,\ve)$ to be specified.

\begin{remark}[On Lyapunov radius]
	\label{rmk:on_lyapunov_radius}
	From the construction of Lyapunov charts, it follows that there exists a measurable, $\ve$-slowly varying, a.e. nonzero function \index{$r_{k,\ve}(x)$}$r_{k,\ve}(x)$ such that if $d^Q(x',x)\leq r_{k,\ve}(x)$ then $x'\in \cL_{k,\ve}[x]$.
	Indeed, we take $r_{k,\ve}' \colon Q\to (0,1)$ to be $r_{k,\ve}'(q) = \frac 1{\phi_{k,\ve}(q)} {r_0}.$
	Then $$B(q, r_{k,\ve}'(q))\subset \cL_{k,\ve}[q].$$
	Since $r_{k,\ve}'$ is tempered, it can be bounded from below by an $\ve$-slowly varying $r_{k,\ve}$.
	Recall that we typically drop $(k,\ve)$ from the notation of Lyapunov charts.
\end{remark}

\subsection{Random Dynamics: The Skew Product Construction}
	\label{ssec:random:skew:product}

Let $\index{$S$@$\cS$}\cS \subset \Diff_\infty(Q)$ denote the support of $\mu$. 
We consider the two sided shift space \index{$S$@$\cS^\Z$}$\cS^\Z$. For $\omega \in \cS^\Z$, we have
$\omega = ( \dots, \omega_{-1}, \omega_0, \omega_1, \dots)$. 
We write $\omega = (\omega^-, \omega^+)$ where \index{$\omega^-$}$\omega^- = (\dots, \omega_{-1})$ is the
``past'', and \index{$\omega^+$}$\omega^+ = (\omega_0, \omega_1, \dots)$ is the ``future''. 
Let \index{$T$}$T: \cS^\Z \to \cS^\Z$ denote the left shift
$(T\omega)_i = \omega_{i+1}$ (which we are
thinking of as ``taking one step into the future'').
We use  \index{$T$@$\hat{T}$}$\hat{T}$ to denote the
``skew product'' map $\cS^\Z \cross Q \to \cS^\Z
\cross Q$ given by
\begin{equation}
\label{eq:def:T:on:SZ:cross:M}
\hat{T}(\omega, q) = (T\omega, \omega_0 \cdot q). 
\end{equation}

Let \index{$Q$@$\hat{Q}$}$\hat{Q} = (\cS^\Z \cross Q) \cross [0,1]$. 
Let \index{$T$@$\hat{T}^t$}$\hat{T}^t$
denote the suspension flow on
on $\hat{Q}$, i.e.\ $\hat{T}^t$ is obtained as a quotient of the flow $(x,s) \to
(x,t+s)$ on $(\cS^\Z \cross Q) \cross \reals$ by the equivalence relation
$(x,s+1) \sim (\hat{T}x,s)$, where $\hat{T}$ is as in (\ref{eq:def:T:on:SZ:cross:M}).

\subsubsection{Measures on skew-products.}
Suppose we are given an ergodic $\mu$-stationary measure $\nu$ on $Q$.
As in \cite{BenoistQuint_Mesures1}, 
for \index{$\omega$}
\begin{displaymath}
\omega = ( \dots, \omega_{-1}, \omega_0, \omega_1, \dots),
\end{displaymath}
let
\begin{displaymath}
\index{$\nu_{\omega^-}$}\nu_{\omega^-} = \lim_{n \to \infty}
\omega_{-1} \dots \omega_{-n} \, \nu.
\end{displaymath}
The fact that the limit exists follows from the Martingale
convergence theorem. 
Then $\nu_{\omega^-}$ is a measure on $Q$.
\medskip

\noindent
\textbf{Basic Fact:} Given a $\mu$-stationary measure $\nu$ on
$Q$, we get a $\hat{T}$-invariant measure $\hat{\nu}$
on $\cS^\Z \cross Q$ given by
\begin{equation}
\label{eq:def:hat:nu}
\index{$\nu$@$\hat{\nu}$}d\hat{\nu}(\omega^-,\omega^+,q) = d\mu^\Z(\omega^-,\omega^+) \, d\nu_{\omega^-}(q). 
\end{equation}

It is important that the measure $\hat{\nu}$ defined in \eqref{eq:def:hat:nu} is a product of a
measure depending on $(\omega^-,q)$ and a measure depending on $\omega^+$. 
(If instead of the two-sided shift space we use the 
one-sided shift $\index{$D$@$\Diff_\infty(Q)^\N$}\Diff_\infty(Q)^\N \cross Q$, then $\mu^\Z
\cross \nu$ would be an invariant measure for $\hat{T}$.)

We also use $\hat{\nu}$ to denote the $\hat{T}^t$-invariant measure on
$\hat{Q} = \cS^\Z \cross Q \cross [0,1]$ given by
\begin{displaymath}
\index{$\nu$@$\hat{\nu}$}d\hat{\nu}(\omega^-,\omega^+,q,s) =
d\mu^\Z(\omega^-,\omega^+) \, d\nu_{\omega^-}(q) \, dm(s), 
\end{displaymath}
where \index{$m(\cdot)$}$m(\cdot)$ denotes the Lebesgue measure on $\R$.

\begin{proposition}
\label{prop:hat:nu:ergodic}
If $\nu$ is an ergodic stationary measure on $Q$, then 
the $\hat{T}$-invariant measure $\hat{\nu}$ on $\Diff_\infty(Q)^\Z
\cross Q$ is $\hat{T}$-ergodic. Also the $\hat{T}^t$-invariant measure
$\hat{\nu}$ on $\hat{Q}$ is $\hat{T}^t$-ergodic. 
\end{proposition}

\begin{proof} This
follows from \cite[Lemma I.2.4, Theorem I.2.1]{Kifer_Ergodic_theory_of_random_transformations}
\end{proof}
\medskip

\noindent
\textbf{Stable and unstable manifolds.}
For $x =(\omega,s) \in \Diff_\infty(Q)^\Z \cross [0,1]$, let
\begin{multline*}
\index{$W$@$W^-[x]$}W^-[x] = \{ (\omega',s') \in
\Diff_\infty(Q)^\Z\cross [0,1] \st \\
\text{ $s' = s$ and for $n\in\N$
  sufficiently large, $\omega'_n = \omega_n$} \}. 
\end{multline*}
Then $W^-[x]$ consists of sequences $y$ which eventually agree with
$x$. We call $W^-[x]$ the ``stable leaf through $x$''.
We also have the subset
\begin{displaymath}
W^-_1[x] = \{ (\omega',m') \in \Diff_\infty(Q)^\Z \cross [0,1] \st m' = m \text{ and
  } (\omega')^+ = \omega^+ \} \subset W^-[x]. 
\end{displaymath}
Similarly, we define
\begin{multline*}
\index{$W^+[x]$}W^+[x] = \{ (\omega',s')  \in \Diff_\infty(Q)^\Z\cross [0,1] \st
\\ \text{ $s' = s$ and for
  $n \in \N$ sufficiently large, $\omega'_{-n} = \omega_{-n}$} \}, 
\end{multline*}
and we also have the subset
\begin{displaymath}
\index{$W^+_1[x]$}W^+_1[x] = \{ (\omega',s') \in \Diff_\infty(Q)^\Z\cross [0,1] \st (\omega')^- = \omega^- \} \subset W^+[x]. 
\end{displaymath}
For $x = (\omega,s) \in \Diff_\infty(Q)^\Z \cross [0,1]$, we define
$\index{$x^+$}x^+=(\omega^+,s)$ and $\index{$x^-$}x^- = (\omega^-,s)$.
Then $W_1^+[x]$ depends only on $x^-$ and $W_1^-[x]$ depends only on
$x^+$.

We now consider the skew-product map. 
The base dynamics is the flow $\hat{T}^t$  on $(\hat Q, \hat \nu)$.  
Associated to almost every $\index{$q$@$\hat{q}$}\hat{q} = (q,\omega,s)\in \hat{Q}$
is an unstable manifold $\widehat{\cW}^u[\hat{q}]$. It takes the form\index{$W$@$\widehat{\cW}^u[q,\omega,s]$}
\begin{displaymath}
\widehat{\cW}^u[q,\omega,s] = W^+[\omega,s] \cross \cW^u[q,\omega,s],
\end{displaymath}
where \index{$W$@$\cW^u[q,\omega,s]$}$\cW^u[q,\omega,s] \subset Q$. 

We write \index{$B$@$\gB_0$}$\gB_0$ for the measurable partition whose atom $\gB_0[\hat q]$ is contained in the (fiber) unstable manifold $\cW^u[\hat q]$ and satisfies $\hat T^t \gB_0\prec \gB_0$ for all $t\ge 0$ as constructed in \autoref{def:unstable_markov_partition}. 
Given $\hat q = (q, \omega, s)$, we write $\index{$B$@$\widehat{\gB}_0$}\widehat \gB_0= W^+_1[\omega,s] \cross \gB_0[\hat q]$.  Then $\widehat \gB_0$ is a measurable partition of $\hat Q$ and also satisfies 
$ \hat T^t \widehat  \gB_0\prec \widehat \gB_0$ for all $t\ge 0$.
We similarly define \index{$B$@$\gB_0^-$}$\gB_0^-$ and \index{$B$@$\widehat{\gB}_0^-$}$\widehat \gB_0^-$.

We use the notations
\index{$W$@$\cW^u[q,\omega,s]$}$\cW^u[q,\omega,s]$, \index{$G$@$\bbG^{sr}(\cW^u[q,\omega,s])$}$\bbG^{sr}(\cW^u[q,\omega,s])$,
\index{$G$@$\bbG^{ssr}(\cW^u[q,\omega,s])$}$\bbG^{ssr}(\cW^u[q,\omega,s])$, \index{$L\cW^u(q,\omega,s)$}$L\cW^u(q,\omega,s)$,
\index{$G$@$\bbG^{sr}(q,\omega,s)$}$\bbG^{sr}(q,\omega,s) \subset \GL(L \cW^u(q,\omega,s))$ 
as in
\autoref{sssec:setup_subgroups_compatible_with_the_measure}.

The flow $\hat{T}^t$ plays a role analogous to $g_t$ in the case of
$\SL(2,\R)$-actions. In particular, by ``the Lyapunov exponents of the
random dynamical system'' we mean the Lyapunov exponents of
$\hat{T}^t$.




\subsection{Margulis functions}
	\label{sec:subsec:Margulis:functions}

\subsubsection*{$\SL(2,\R)$ actions}
Let
$$\index{$g_t$}g_t = \begin{pmatrix} e^t & 0 \\ 0 & e^{-t} \end{pmatrix}
\in \SL(2,\R), \qquad
\index{$r_\theta$}r_\theta = \begin{pmatrix} \cos \theta & \sin \theta \\ -\sin \theta
                                        & \cos \theta \end{pmatrix}
                                      \in \SL(2,\R).$$ 
Let\index{$A$@$(A_\tau f)(x)$}
\begin{displaymath}
(A_\tau f)(x) = \frac{1}{2\pi} \int_0^{2\pi} f(g_\tau r_\theta x) \, d\theta. 
\end{displaymath}

\begin{definition}
\label{def:SL(2,R):admissible:margulis:function}
A continuous function $h: Q \to \reals^+$ is called an admissible
Margulis function on $Q$ if there exist $\tau > 0$, $c < 1$, $b > 0$
and $C > 0$ and $\delta > 0$ such that the following hold:
\begin{itemize}
\item[{\rm (a)}]
For all $q \in Q$,
\begin{displaymath}
(A_\tau h)(q) \le c \; h(q) + b.
\end{displaymath}
\item[{\rm (b)}] For all $q \in Q$,
\begin{displaymath}
h(q) \ge C \phi_a(q)^\delta,
\end{displaymath}
where $\phi_a$ is the systole function.
\end{itemize}
\end{definition}

\subsubsection*{Random dynamics}

\begin{definition}
\label{def:random:margulis:function}
A continuous function $h: Q \to \reals^+$ is called an admissible
Margulis function on $Q$ if there exist $c < 1$, $b > 0$
and $C > 0$ and $\delta > 0$ such that the following hold:
\begin{itemize}
\item[{\rm (a)}]
For all $q \in Q$,
\begin{equation}
\label{eq:def:random:margulis}
  \int_{\Diff_\infty(Q)} h(f(q)) \, d\mu(f) \le c \, h(q) + b.
\end{equation}
\item[{\rm (b)}] For all $q \in Q$,
\begin{displaymath}
h(q) > C \phi_a(q)^\delta,
\end{displaymath}
where $\phi_a$ is the systole function.
\end{itemize}
\end{definition}

\begin{proposition}
\label{prop:admissible:margulis:implies:good:measure}
Suppose that either $Q$ is compact or there exists an admissible
Margulis function on $Q$. Let $\nu$ be an ergodic $P$-invariant (or in the
random dynamics setting $\mu$-stationary) measure. Also assume in the
random dynamics setting that $\mu$ is finitely supported. 
Then, (\ref{eq:systole:log:integrable}) holds. 
\end{proposition}

\begin{proof}
Suppose $q \in Q$. Then, iterating (\ref{eq:def:random:margulis}) we get for any
$n \in \N$, 
\begin{displaymath}
\int_{\Diff_\infty(Q)} h(f(q)) \, d\mu^{(n)}(f) \le c^n \, h(q) + \frac{b}{1-c}.
\end{displaymath}
where $\mu^{(n)}$ is the $n$-fold convolution of $\mu$ with itself. 
Then, by the random ergodic theorem
\cite[Corollary~I.2.2]{Kifer_Ergodic_theory_of_random_transformations},
for $\nu$-a.e.\ $q \in Q$,
\begin{displaymath}
\frac{1}{n} \sum_{k =0}^{n-1} \int_{\Diff_\infty(Q)} h(f(q)) \,
d\mu^{(k)}(f) \to \int_Q h \, d\nu.
\end{displaymath}
Thus,
\begin{displaymath}
\int_Q h \, d\nu \le \frac{b}{1-c}.
\end{displaymath}
Since $f$ is admissible, this implies that for some $\delta > 0$,
\begin{displaymath}
\int_Q \phi_a(q)^\delta \, d\nu(q) \le \frac{Cb}{1-c}.
\end{displaymath}
This implies the statement of the lemma. 
\end{proof}

\subsubsection*{Avoiding the diagonal}

\begin{definition}
\label{def:random:QxQ:margulis:function}
A continuous function $h: Q \cross Q \to \reals^+ \cup \{\infty\}$ is
called an admissible 
Margulis function on $Q\cross Q$ if there exist $c < 1$, $b > 0$
and $C > 0$ and $\delta > 0$ such that the following hold:
\begin{itemize}
\item[{\rm (a)}] $h(x,y) < \infty$ unless $x = y$. 
\item[{\rm (b)}]
For all $x,y \in Q \cross Q$,
\begin{equation}
\label{eq:QxQ:margulis:inequality}
 \int_{\Diff_\infty(Q)} h(f(x),f(y))) \, d\mu(f) \le c \, h(x,y) + b.
\end{equation}
\item[{\rm (c)}] For all $(x,y) \in Q \cross Q$,
\begin{displaymath}
h(x,y) \ge C d(x,y)^{-\delta} + C \phi_a(x)^\delta + C \phi_a(y)^\delta.
\end{displaymath}
where $\phi_a$ is the systole function.
\end{itemize}
\end{definition}

\begin{lemma}
\label{lemma:Q:compact:plus:ue:implies:margulis:function}
Suppose $Q$ is compact and suppose $\mu$ satisfies uniform expansion
in dimension $1$ (see \autoref{def:uniformly_expanding}). Then,
there exists a Margulis function on $Q \cross Q$.
\end{lemma}

\begin{proof}
We define the length of a path $\gamma: [0,1] \to Q$ to be
\begin{displaymath}
\int_0^1 \|\gamma'(t) \| \, dt, 
\end{displaymath}
where the norm $\| \cdot \|$ is as in (\ref{eq:def:sigma:f:L}).
Given $x,y \in Q$, let $d(x,y)$ be the infimum of the lengths of paths
conecting $x$ and $y$. For $\delta > 0$, let $h(x,y) =
d(x,y)^{-\delta}$. Then, arguing is in the proof of \cite[Lemma~4.2]{EskinMargulis_Recurrence} we see that $h$ is a Margulis function on $Q \cross Q$. 
\end{proof}

\subsubsection{Dimension bounds}
We define, for $x \in Q$, the ``lower local dimension''
\begin{displaymath}
\index{$d$@$\underline{\dim}(\nu,x)$}\underline{\dim}(\nu,x) = \liminf_{r \to 0}
  \frac{\log \nu(B(x,r)}{\log r}, 
\end{displaymath}
where
\begin{displaymath}
\index{$B(x,r)$}B(x,r) = \{ y \in Q \st d(x,y) < r \}.
\end{displaymath}

\begin{proposition}
\label{prop:Margulis:function:implies:positive:dimension}
Suppose there exists an admissible
Margulis function on $Q \cross Q$. Let $\nu$ be  ergodic
$\mu$-stationary measure. Then, either $\nu$ is finitely supported, 
or for almost all $x \in Q$, $\underline{\dim}(\nu,x) > 0$. 
\end{proposition}

\begin{proof} This argument is taken from
\cite{Eskin:Lindenstrauss:short}. Suppose $\nu$ is not finitely
supported. We note that in view of the
ergodicity of $\nu$, $\nu$ has no atoms. Therefore the $\nu \cross \nu$
measure of the diagonal is $0$. 
  
By iterating (\ref{eq:QxQ:margulis:inequality}), for any
$x, y \in Q$ with $x \ne y$, 
\begin{equation}
\label{eq:result:margulis:iteration}
\limsup_{k \to \infty} \int_G h(f(x), f(y))) \,
d\mu^{(k)}(f) \le \frac{b}{1-c}.  
\end{equation}
By the random ergodic theorem \cite[Corollary~I.2.2]{Kifer_Ergodic_theory_of_random_transformations}, there
exists a function $\phi: Q \cross Q \to \reals$ such
that
\begin{displaymath}
\int_{Q\cross Q} \phi \, d(\nu \cross \nu) =
\int_{Q \cross Q} f \, d(\nu \cross \nu), 
\end{displaymath}
and for $\mu^\N$ a.e.\ $(f_1, \dots, f_n, \dots) \in \Diff_\infty(Q)^\N$ and $\nu \cross \nu$ a.e.\
$x$, $y \in Q$, 
\begin{equation}
\label{eq:random:ergodic:short}
\phi(x, y) = \lim_{k \to \infty} \frac{1}{k} \sum_{j=1}^k h( (f_{j} \dots f_1) x, (f_{j} \dots f_1) y).
\end{equation}
Then, integrating both sides of (\ref{eq:random:ergodic:short}) over
$\Diff_\infty(Q)^\N \cross Q \cross Q$, using Fatou's lemma to take
the limsup outside the integral, and then using
(\ref{eq:result:margulis:iteration}), we get
\begin{displaymath}
\int_{Q \cross Q} h \, d(\nu \cross \nu)\le \frac{b}{1-c}. 
\end{displaymath}
Therefore, for any $\eta > 0$ there exists $K'' \subset Q$ with
$\nu(K'') > 1-c(\eta)$ where $c(\eta) \to 0$ as $\eta \to
0$ and a constant $C = C(\eta)$ 
such that for any $x \in K''$, 
\begin{equation}
\label{eq:f:in:L1}
\int_{Q} h(x,y) \, d\nu(y) < C. 
\end{equation}
It follows from (\ref{eq:f:in:L1})
that for all $\epsilon > 0$ and all $x \in K''$,
\begin{displaymath}
\nu(B(x,\epsilon)) \le C(\eta) \epsilon^\delta,
\end{displaymath}
hence 
\begin{displaymath}
\frac{\log \nu(B(x,\epsilon))}{\log \epsilon} \ge \delta - \frac{|\log
  C(\eta)|}{|\log \epsilon|}. 
\end{displaymath}
This implies $\underline{\dim}(\nu,x) \ge \delta$ for $x \in
K''$.
Then, in view of the ergodicity of $\nu$,
the same holds for a.e.\ $x \in Q$. 
\end{proof}



\subsection{The finite measurable cover}
	\label{ssec:the_finite_measurable_cover}


\subsubsection{Basic setup}
The base dynamics is the flow $g_t$  on $(Q,\nu)$.  We assume the $g_t$-action on $(Q,\nu)$ is ergodic.  
Associated to almost every $q\in Q$ is an unstable manifold \index{$W$@$\cW^u[q]$}$\cW^u[q]$.
We have \index{$G$@$\bbG^{ssr}(q)$}$\bbG^{ssr}(q)$ as a 
subgroup of $\GL(L \cW^u(q))$. \index{$L\cW^u(q)$}
Recall for almost every $q\in Q$ and every $q'\in \cW^u[q]$, we have a well defined $L \cW^u(q')$ and $\bbG^{ssr}(q')$ and \index{$P^{ssr}_{GM}(q,q')$} identifications:
\begin{align}
	\label{eqn:defining_P_on_Gssr}
	\begin{split}
	P^{ssr}_{GM}(q,q')\colon \bbG^{ssr}(q) & \toisom  \bbG^{ssr}(q')\\
	g & \mapsto P^{L\cW^u}_{GM}(q,q') \circ g\circ  P^{L\cW^u}_{GM}(q',q)
	\end{split}	
\end{align}
where $P^{L\cW^u}_{GM}(q,q')\colon L\cW^u(q) \to L\cW^u(q')$ \index{$P^{L\cW^u}_{GM}(q,q')$} is as in
\autoref{thm:linearization_of_stable_dynamics_single_diffeo} (see also \autoref{sssec:holonomies_and_relating_different_groups}).

If needed, we also write \index{$G$@$\bbG^{sr}(\cW^u[q])$}$\bbG^{sr}(\cW^u[q])$ and
\index{$G$@$\bbG^{ssr}(\cW^u[q])$}$\bbG^{ssr}(\cW^u[q])$ to denote the groups of self-maps of
$\cW^u[q]$.  Thus with this notation, for $q'\in \cW^u[q]$ we
have $$\bbG^{*}(\cW^u[q])= \bbG^{*}(\cW^u[q'])$$ for $*\in
\{sr,ssr\}$.

Let \index{$L$@$\liessr(q)$}$\liessr(q)$ denote the Lie algebra of $ \bbG^{ssr}(q)$;  
we have
$$\mathrm {Ad}\left(P^{L\cW^u}_{GM}(q,q')\right) \liessr(q)=  \liessr(q').$$
Note that this map is the differential of the one defined in \autoref{eqn:defining_P_on_Gssr}.
The $g_t$-action intertwining $\cW^u$-manifolds induces an action intertwining the groups $\bbG^{ssr}$ which induces a linear cocycle  $$g_t\colon  \liessr(q) \to \liessr(g_tq).$$

\def\gF{\mathfrak F}
\subsubsection{The measurable cover}
	\label{sssec:the_measurable_cover}
Let \index{$L$}$L= \liessr$. We apply
\autoref{thm:jordan_normal_form} to the cocycle $L$ and the backwards
dynamics $g_{-t}$ to obtain a finite measurable cover $\index{$X$}X= Q\times F$.  Let $\index{$\sigma$}\sigma
\colon X\to Q$ be the finite-to-one map.  We let \index{$g_t$}$g_t$ denote the
induced cocycle on $X$.  We abuse notation and denote the measure on $X$ by \index{$\nu$}$\nu$.  

We view a point $x\in X$ as a pair $(q, \gF)$ where $q= \sigma (x)$
and $\index{$F$@$\gF$}\gF\in \gF_\bullet L(q)$ is a possible flag as in \autoref{thm:jordan_normal_form}.
We emphasize, that since we used the backwards dynamics to construct
$X$, $\gF$ is a backwards flag. 

The cocycles that will be key for our arguments, such as $\bbH$ from \autoref{ssec:the_bold_h_cocycle}, can be put in Jordan normal form on $X$ (see for instance \autoref{prop:LC_Jordan_normal_form}).

\def\gF{\mathfrak F}
\subsubsection{Stable and unstable manifolds on the measurable cover}
\label{sec:subsec:stable:manifold:measurable:cover}
Consider a biregular point $q\in Q$ and $\gF\in \gF_\bullet L(q)$.  Let $x= (q, \gF)\in \sigma^{-1}(q)$.  

We use \autoref{thm:jordan_normal_form}\ref{Jordanextention:iii} to define the measurable stable set $\index{$W$@$\cW^s[x]$}\cW^s[x]=\cW^s[q, \gF]$ of $x$ and  \autoref{thm:jordan_normal_form}\ref{Jordanextention:iiii} to define the (smooth) unstable manifold  $\index{$W$@$\cW^u[x]$}\cW^u[x]=\cW^u[q, \gF]$ of $x$ in $X$.



Recall from \autoref{ssec:conditional_measures}
the leaf-wise measure $\nu_q^s$  and $\nu_q^u$ defined on almost every $\cW^s[q]$ and $\cW^u[q]$, respectively.  
Given $q$ and $x\in \sigma^{-1}(q),$ we associate the leaf-wise measures $\index{$\nu_x^s$}\nu_x^s$ and \index{$\nu^u_x$}$\nu_x^u$ on $\cW^s[x]$ and $\cW^u[x]$, respectively via pull-back under the 1-1 map $\sigma \colon \cW^s[x]\to \cW^s[q]$ and $\sigma \colon \cW^u[x]\to \cW^u[q]$.

We note that $g_t$ intertwines stable sets and unstable manifolds on $X$.  Moreover, for a.e.\ $x$ and every $t\in \R$, the map $g_t\colon \cW^u[x]\to \cW^u[g_tx]$ is a diffeomorphism.

\subsubsection{Measurable connections on finite covers}
	\label{sssec:measurable_connections_on_finite_covers}
Recall from \autoref{sssec:standard_measurable_connection} that for any admissible cocycles $E$ on $Q$ we have measurable connections \index{$P^{\pm}(q,q')$}$P^{\pm}(q,q')$ for $q,q'$ on the same stable/unstable manifolds.
For such a cocycle $E$ pulled back to $X$, we will use the same notation and define by pullback the measurable connections \index{$P^{\pm}(x,x')$}$P^{\pm}(x,x')$ for $x'\in \cW^s[x]$ or $x'\in \cW^u[x]$ (with these manifolds defined as above).
We observe that on the finite cover $X$, for the subcocycle $\fraku^+\subset \frakg^{ssr}$, the measurable connections preserve it, by the Ledrappier invariance principle \autoref{thm:ledrappier_invariance_principle}.

\subsubsection{Pulling back Markov Partitions}
\label{sec:subsec:pulling:back:Markov}
Given a Markov partition $\gB_0$ of $Q$ subordinate to the unstables (see
\autoref{ssec:measurable_partitions})
we can pull it back to $X$. More precisely, we set
\begin{displaymath}
\index{$B$@$\gB_0[x]$}\gB_0[x] = \sigma^{-1}(\gB_0[\sigma(x)]) \cap \cW^u[x]. 
\end{displaymath}
Similarly, for the partition $\gB^-$ of $Q$ subordinate to the
stables, we set
\begin{displaymath}
\index{$B$@$\gB_0^-[x]$}\gB_0^-[x] = \sigma^{-1}(\gB_0^-[\sigma(x)]) \cap \cW^s[x]. 
\end{displaymath}
We then use all the notation of \autoref{ssec:measurable_partitions}.

\subsubsection{Distances and charts in $Q$}
	\label{sssec:summary_lyapunov_charts_and_distances_intro}
In addition to standard constructions of stable and unstable manifolds, we will make use of the following distances adapted to the dynamics:
\begin{description}
	\item[Ambient distance] \index{$d^Q$}$d^Q$ coming from the Riemannian metric on $Q$.
	\item[Lyapunov distance] \index{$d^{\cL}_q(q',q'')$}$d^{\cL}_q(q',q'')$ will only be defined for $q',q''\in \cL_{k,\ve}[q]$ (see \autoref{sssec:lyapunov_charts_and_distance}); the parameters $k,\ve$ will be suppressed from the notation.
	\item[Unstable distance] \index{$d^u(q,q')$}$d^u(q,q')$ will only be defined for $q'\in \cW^u[q]$ and will have a uniform contraction property under $g_{-t}$ for $t\geq 0$ (see \autoref{sssec:unstable_lyapunov_distance_and_charts}); the parameter $\ve>0$ is suppressed from the notation.
	\item[Induced ambient distance] $d^{Q,u}$ on $\cW^u[q]$ induced by the path metric from the ambient distance.
\end{description}
We will also construct the following types of charts adapted to the dynamics:
\begin{description}
	\item[Lyapunov charts] $\index{$L$@$\cL_{k,\ve}[q]$}\cL_{k,\ve}[q]\subset Q$ will denote the $(k,\ve)$-Lyapunov chart at $q$ (see \autoref{sssec:lyapunov_charts_and_distance}).
	\item[Local unstable manifold] (see \autoref{sssec:unstable_lyapunov_distance_and_charts})
	\[
		\index{$W$@$\cW^{u}_{loc}[q]$}\cW^{u}_{loc}[q]:=\{q'\colon d^u(q,q')<1\} \subset \cW^u[q].
	\]
\end{description}

\subsubsection{Distances in $X$}
	\label{sec:subsec:distances}
For $x_1, x_2 \in X$, we will
occasionally abuse notation and write \index{$d^Q(x_1,x_2)$}$d^Q(x_1,x_2)$ in place of
$d^Q(\sigma(x_1), \sigma(x_2))$.
We equip $X$ with the (non-Hausdorff) topology induced by this metric.  All notions of convergence in $X$ will be with respect to this topology.

For two sets $A,B\subset X$ we define
\[
	\index{$d^{Q}_{min}(A,B)$}d^{Q}_{min}(A,B):=\inf_{a\in A, b\in B}d^Q(a,b).
\]
On the finite cover $X$, the distance \index{$d^u(x_1,x_2)$}$d^u$ for points in the same unstable is defined such that the map $\sigma\colon  \cW^u[x]\to  \cW^u[\sigma(x)]$ is 
is a diffeomorphism and an isometry for the $d^u$-distance.
Note that $d^{\cL}$ is only defined for points on $Q$.



\begin{lemma}
\label{lemma:du:vs:dQ}
There exists a measurable function $\index{$r$}r\colon Q \to \reals^+$ such that the following holds:  
Suppose $x_1,x_2, y \in X$ and $x_1,x_2 \in \cW^u_{loc}[y]$.
Then,
\begin{displaymath}
r(y)d^Q(x_1, x_2) \le d^u(x_1,x_2) \le r(y)^{-1} d^Q(x_1, x_2).
\end{displaymath}
\end{lemma}

\begin{proof}
This follows from the construction of $d^u$ in \autoref{ssec:lyapunov_charts_and_stable_manifolds}.
\end{proof}



\begin{lemma}[Staying in the local unstable]
	\label{lem:lemma_stay_in_chart_stay_in_loc_new}
	Fix $\beta>0$.
	For any $\delta>0$ there exists a compact set $K\subset Q$ of measure at least $1-\delta$ with the following property.
	For $x_1,x_2,y\in X$ and $\tau>0$, suppose that $y,g_{-\tau}y\in \sigma\inv(K)$, $x_1\in \cW^u_{loc}[y]$, and that $g_{-\tau}x_1,g_{-\tau}x_2,\in \cW^u_{loc}[g_{-\tau}y]$.
	Suppose furthermore that $d^{Q}\left(g_{-t} x_1,g_{-t}x_2\right)\leq e^{-\beta \tau}$ for $t\in [0,\tau]$.

	Then, provided $\tau$ is sufficiently large, we have $x_i\in \cW^u_{loc}[x_{3-i}]$
 and  $d^u(x_1,x_2)\leq C(\delta)e^{-\beta\tau/2}$.
 \end{lemma}

\begin{proof}
	By the definition of $d^Q(\cdot, \cdot)$, $\cW^u[y]$, and $d^u(\cdot, \cdot)$ on $X$, it suffices to prove the analogous statements for  $\sigma(y), \sigma(x_1), \sigma(x_2)$.
	Indeed, the projection map $\cW^u[x]\to \cW^u[\sigma(x)]$ is a $d^u$-isometry, and the hypotheses already imply that $x_i\in \cW^u\left[x_{3-i}\right]$, so it only remains to establish the distance bounds.
	We thus assume for the remainder that $y, x_1, x_2\in Q$.

	Fix $\beta>0$.  There is $\lambda_0>0$ and a measurable function $C\colon Q\to [1,\infty)$ such that for $\nu$-a.e.\ $z$,  $d^Q(g_{-t} z, g_{-t}x)\le C(z) e^{-\lambda_0 t}$ for all $x\in \cW^u_{loc}[z]$ and $t\ge 0$ (see \autoref{prop:stableman} and \autoref{rem:WlocinLyapcharts}).

	Fix a choice of $0<\ve<\frac 1 {10}\min\{ \lambda_0, \beta\}$ and let $\cL$ be a choice of $(\ve,1)$-Lyapunov charts.  
	Given such a choice of Lyapunov charts $\cL$, we write $\cW^u_{\cL}[y]$ for the path component of $\cW^u[y]\cap \cL[y]$ containing $y$.  For $t\in [0,1]$ we have (see \autoref{prop:stableman})
	\begin{equation} \label{eq:overflowlyapcharts}g_{t}(\cW^u_{\cL}[y])\cap \cL[g_ty]= \cW^u_{\cL}[g_ty].\end{equation}

	Fix $L\gg 1$ sufficiently large.  Let $K_0$ be a compact subset of measure at least $1-\delta/5$ such that for $z\in K_0$
	\begin{enumerate}
	\item the identity map $( \cW^u_{loc}[z], d^Q)\to ( \cW^u_{loc}[z], d^u)$ is  $L$-Lipschitz,
	\item $d^Q(g_{-t} z, g_{-t}x)\le Le^{-\lambda_0 t}$ for all $x\in \cW^u_{loc}[z]$ and $t\ge 0$,
	\item $B^u(z, L\inv):= \{x\in \cW^u[z]: d^u(x,z)<L\inv\} \subset \cW^u_{\cL}[z]$, and 
	\item $r(z)\ge L\inv $ (where  $r(\cdot)$ is the function that  $\ve$-slowly decreases  along 2-sided orbits from   \autoref{rmk:on_lyapunov_radius}.)
	\end{enumerate}
	Fix $\tau_0>0$ so that for all $t\ge \tau_0$,  
	$$L^2e^{-t(\lambda_0 -\ve) }+ L e^{- t(\beta- \ve)}\le 1.$$
	Fix $L'\gg L$ sufficiently large and take $K_1\subset K_0$ to be a  compact subset of measure at least $1-2/5\delta$ such that for $z\in K_1$
	\begin{enumerate}[resume]
	\item $d^u( g_{\tau_0}z_1, g_{\tau_0}z_2)\le L'd^u(z_1,z_2)$ for all $z_1,z_2\in  \cW^u_{loc}[z].$
	\end{enumerate}
	Let $K\subset K_0\cap  g_{\tau_0}K_1$   be a compact subset of measure at least $1-\delta$.
	Finally, fix $\tau\ge \tau_0$ sufficiently large so that  $$L^4e^{-\lambda_0 \tau} + L^3 e^{-\tau\beta}\le 1.$$

	We claim the proposition holds with the above choices.
	Indeed, 
	suppose $y\in K$.  It follows for all $t\in [\tau_0,\tau]$ and $i= 1,2$ that $$d^Q(g_{-t} y, g_{-t}x_i)\le Le^{-t \lambda_0 } + e^{-t\beta}\le L\inv e^{-\ve t}\le  r(y) e^{-\ve t}\le r(g_{-t}(y)).$$
	It follows that $g_{-t}x_i\in  \cL[(g_{-t} y)]$ for both $i= 1,2$ and all $t\in [\tau_0,\tau]$. 
	Also, supposing $g_{-\tau} y\in K\subset K_0$ and $g_{-\tau} x_i\in \cW^u_{loc}[g_{-\tau} y]$,
	we have  
	$$d^u(g_{-\tau} y, g_{-\tau}x_i)\le L d^Q(g_{-\tau} y, g_{-\tau}x_i)\le  L^3e^{-\tau\lambda_0 } + L^2 e^{-\tau\beta}\le L^{-1}$$ 
	and so $g_{-\tau}x_i\in \cW^u_{\cL}[g_{-\tau}y]$ for both $i= 1,2$.  
	By  \eqref{eq:overflowlyapcharts}, 
	we have 
	$$g_{-t}x_i\in \cW^u_{\cL}[g_{-t} y]$$ for all  $\tau_0\le t\le \tau$ and $i=1,2$.  In particular, 
	$$g_{-\tau_0}x_i\in \cW^u_{\cL}[g_{-\tau_0} y]$$ for both $i=1,2$.  Also since $g_{-\tau_0} y\subset K_1 \subset K_0$, we have 
	$$d^u(g_{-\tau_0} y, g_{-\tau_0}x_i)\le Ld^Q(g_{-\tau_0} y, g_{-\tau_0}x_i)\le L^2e^{-\tau_0\lambda_0  } + Le^{-\tau_0\beta}\le 1$$
	whence $g_{-\tau_0}x_i\in \cW^u_{loc}[g_{-\tau_0} y]$ for both $i= 1,2$.  
	Since $g_{-\tau_0} y\in K_1\subset g_{-\tau_0}(K)$ we have 
	$$ d^u(x_1,x_2)\le L'd^u(g_{-\tau_0}x_1,g_{-\tau_0}x_2)\le (L'L) d^Q(g_{-\tau_0}x_1,g_{-\tau_0}x_2)\le( L'L) e^{-\beta \tau}$$
	and the proposition follows.
\end{proof}


\subsection{The skew product and the measurable cover for random dynamics}
\label{sec:subsec:skew_product_measurable_cover}

Let $\index{$X$@$\hat{X}$}\hat{X} $ be the measurable cover of $\hat Q$  
built as in \autoref{ssec:the_finite_measurable_cover}. 
Then, the flow
$\hat{T}^t$ on $\hat{Q}$ lifts to a flow (also denoted \index{$T$@$\hat{T}^t$}$\hat{T}^t$) on
$\hat{X}$. All the constructions of \autoref{ssec:the_finite_measurable_cover}
carry over with no modifications.



\section{Generalized u-Gibbs states and QNI}
	\label{sec:generalized_u_gibbs_states_and_qni}








\subsection{Subgroups compatible with the measure}
	\label{ssec:subgroups_compatible_with_the_measure}

\subsubsection{Setup}
	\label{sssec:setup_subgroups_compatible_with_the_measure}
Recall that unstable manifolds $\cW^{u}[q]$ admits a subresonant structure, see \autoref{appendix:subresonant_linear_algebra} for the formalism.
We also have the family of subresonant, and strictly subresonant, transformations of the unstable manifold denoted \index{$G$@$\bbG^{sr}(\cW^u[q])$}$\bbG^{sr}(\cW^u[q]), \bbG^{ssr}(\cW^u[q])$\index{$G$@$\bbG^{ssr}(\cW^u[q])$} (see \autoref{cor:strictly_subresonant_maps_on_subresonant_manifolds} for the definition of strictly subresonant maps of an unstable manifold).
Note that if $q'\in \cW^u[q]$ then we have a natural identification $\cW^u[q]=\cW^u[q']$ and hence an identification of the associated groups.
To study the dynamics on unstable manifolds in more detail, we introduce the following notation.

Recall that \index{$L\cW^u(q)$}$L\cW^u(q)$ denotes the fiber of the linearization cocycle, and we have $g_t$-equivariant embeddings $\cW^{u}[q]\into L\cW^u(q)$ such that the linear cocycle maps on fibers of $L\cW^u$ induce subresonant maps between unstable manifolds (see \autoref{thm:linearization_of_stable_dynamics_single_diffeo}).
For $q\in Q$, let \index{$G$@$\bbG^{ssr}(\cW^u[q])$}$\bbG^{sr}(q)\subset \GL(L \cW^u(q))$ denote the group of linear maps induced by subresonant maps of $\cW^u[q]$.


\subsubsection{Holonomies and relating different groups}
	\label{sssec:holonomies_and_relating_different_groups}
Recall that the cocycle $L\cW^u$ admits unstable holonomies $P^{L\cW^u}_{GM}$, which we will call the ``Gauss--Manin'' connection.
So for $q'\in \cW^u_{loc}[q]$ the isomorphism \index{$P^{L\cW^u}_{GM}(q,q')$}$P^{L\cW^u}_{GM}(q,q')$ identifies $\bbG^{sr}(q)$ with $\bbG^{sr}(q')$ via conjugation, and similarly for $\bbG^{ssr}$.  
Equivalently, all these natural isomorphisms can be seen via the identification of $\bbG^{sr}(q)$ with the subresonant automorphisms of $\cW^u[q]$, and then the set-theoretic identification of $\cW^u[q]$ with $\cW^u[q']$ inside $Q$.
We will write:
\[
	P_{GM}^{\bullet}(q,q')\colon \bbG^{\bullet}(q) \to \bbG^{\bullet}(q') \quad \text{ for }q'\in \cW^u[q]
\]
where $\bullet\in\{sr,ssr\}$ for the corresponding identifications (see also \autoref{eqn:defining_P_on_Gssr}).

\subsubsection{Lifting to $X$.}
All the constructions of 
\autoref{sssec:setup_subgroups_compatible_with_the_measure} and
\autoref{sssec:holonomies_and_relating_different_groups} naturally
lift from $Q$ to $X$. We will thus use notation such as
\index{$G$@$\bbG^{ssr}(x)$}$\bbG^{ssr}(x)$ for $x \in X$ to denote the lift. 

\subsubsection{Stabilizers and subgroups}
	\label{sssec:stabilizers_and_subgroups}
For $x \in X$, let also \index{$G$@$\bbG^{ssr}_x(x)$}$\bbG^{ssr}_x(x)$ denote the \emph{stabilizer} of $x\in \cW^u[x]$ inside $\bbG^{ssr}(x)$.
Note that this family of subgroups is \emph{not} preserved by the holonomy identifications.

Recall that the strictly subresonant maps $\bbG^{ssr}(x)$ are a unipotent algebraic group, and we will study unipotent algebraic subgroups $\index{$U^+(x)$}U^+(x)\subseteq \bbG^{ssr}(x)$, which we identify with their Lie algebras \index{$u$@$\fraku^+(x)$}$\fraku^+(x)$.
For such subgroups, we will also denote by \index{$U^+_x(x)$}$U^+_x(x)$ and \index{$u$@$\fraku^+_x(x)$}$\fraku^+_x(x)$ the group and Lie algebra of the stabilizer of $x\in \cW^u[x]$.

Additionally, for a subgroup $U^+(x)$ we will denote by $\index{$U^+[x]$}U^+[x]:=U^+(x)\cdot x\subset \cW^u[x]$ the corresponding orbit of $U^+$ on the unstable manifold.
As a $U^+(x)$-homogeneous space, it is isomorphic to $U^+(x)/U^+_x(x)$.
Note that since both groups are unipotent (as subgroups of the unipotent $\bbG^{ssr}$), they are unimodular and hence there is a unique up to scaling Haar measure on the quotient.

An important property of $U^+[x]$ is that it is an algebraic subset of $\cW^u[x]$ (and not just a connected component of an algebraic subset), when $\cW^u[x]$ is viewed with its subresonant structure.
Indeed, by \cite[III\S2, Prop.~6]{Serre2002_Galois-cohomology} the Galois cohomology group $H^1(\bR,\bbU)$ vanishes for any unipotent algebraic group $\bbU$, therefore $U^+[x]$ coincides with the real points of the algebraic variety corresponding to the orbit $U^+(x)\cdot x$ (see also \cite[Thm.~3.1.1]{Zimmer1984_Ergodic-theory-and-semisimple-groups}).

The following is one of the fundamental definitions we will use, following \cite[Def.~6.2]{EskinMirzakhani_Invariant-and-stationary-measures-for-the-rm-SL2Bbb-R-action-on-moduli-space}.
For notation on measurable partitions, see \autoref{ssec:measurable_partitions}.

\begin{definition}[Compatible family of subgroups]
\label{def:compatible_family_of_subgroups}
	Let $\nu$ be an ergodic, $g_t$-invariant probability measure
        on $X$ which is the lift of an ergodic, $g_t$-invariant
        probability measure on $Q.$
	We say that a measurable family of connected algebraic subgroups $x\mapsto U^+(x)\subset \bbG^{ssr}(x)$ is \emph{compatible with $\nu$} if the following hold:
	\begin{enumerate}
		\item The subgroups are $g_t$-equivariant.
		\item Let \index{$B$@$\gB_0$}$\gB_0$ be the partition of $X$ defined in
                  \autoref{ssec:measurable_partitions} and
                  \S\ref{sec:subsec:pulling:back:Markov}. 
		Then the sets of the form $U^+[x]\cap \gB_0[x]$ form a measurable partition of $X$.

		Furthermore, the conditional measures along $U^+[x]\cap \gB_0[x]$ are proportional to the restriction of the unique $U^+(x)$-invariant measure on $U^+[x]$, which will be called the Haar measure on $U^+[x]$.
		\item Given $x'\in \cW^u[x]$, we have natural identification of $\bbG^{ssr}(x) $ and $\bbG^{ssr}(x')$, viewing both as acting on $\cW^u[x]=\cW^u[x']$.
		By further restricting to a set of full measure, if $x',x\in X_{0}$ and $x'\in U^+[x]$, then  $U^+(x)=U^+(x')$.
		
		Equivalently, viewed inside $\bbG^{ssr}(x)$ and $\bbG^{ssr}(x')$, the groups $U^+(x)$ and $U^+(x')$ are identified by the holonomy $P_{GM}^{ssr}(x,x')$ as in \autoref{sssec:holonomies_and_relating_different_groups}.
		\item $U^+(x)$ is maximal in the sense that it is the largest connected subgroup that preserves the Haar measure on $U^+[x]$.
	\end{enumerate}
We say that the family $x \to U^+(x)$ is \emph{non-trivial} if $\dim
U^+[x] > 0$ for $\nu$-a.e.\ $x$. 
\end{definition}

\begin{remark}[On measures with local invariance]
	\label{rmk:on_measures_with_local_invariance}
	The above \autoref{def:compatible_family_of_subgroups} is meant to formalize the following situation.
	Suppose $\cW^u[x]$ is equipped with its leafwise (conditional) measure \index{$\nu^u_x$}$\nu^u_x$.
	We would like to consider the possibility that $\nu^u_x$-a.e. point $y$ has associated to it a group $U^+(y)\subset \bbG^{ssr}(\cW^u[x])$, but the groups $U^+(y)$ need not be the same.
	Nonetheless, the above definition expresses the possibility that the measure $\nu^u_x$ is ``Haar along $U^+$-orbits''.
\end{remark}

\begin{definition}[Generalized $u$-Gibbs state]
\label{def:generalized:u:Gibbs}
Let $\nu$ be an ergodic, $g_t$-invariant probability measure on $X$
which is the lift of an ergodic, $g_t$-invariant probability measure
on $Q$. Suppose there exists a non-trivial family $x \to U^+(x)$ of
subgroups compatible with $\nu$. Then, we call $\nu$ a ``generalized
$u$-Gibbs state''. 
\end{definition}



\subsection{The QNI condition for extra invariance}
\label{sec:subsec:conditions_for_extra_invariance}

Let $x\mapsto U^+(x)\subset \bbG^{ssr}(x)$ be   
a non-trivial measurable family of connected algebraic subgroups which is
compatible with $\nu$ in the sense of
\autoref{ssec:subgroups_compatible_with_the_measure}. 

Our strategy is to show that under some conditions on $U^+(\cdot)$,
there
exists a family of  subgroups 
$U_{new}^+(x)$ of $\bbG^{ssr}(x)$ compatible with $\nu$ in the sense of
\autoref{def:compatible_family_of_subgroups} such that for almost all
$x$ we have  $U^+(x) \subset U_{new}^+(x) $ and $U^+(x) \neq  U_{new}^+(x) $.
(See \autoref{thm:inductive:step} for the precise
statement).

\subsubsection{The QNI Condition}
We start with a technical condition which encapsulates what is
needed for our proof to go through. It was introduced in
\cite{Katz_Measure-rigidity-of-Anosov-flows-via-the-factorization-method}
in a somewhat different situation. 
\begin{definition}[QNI]
\label{def:QNI}
We say that $U^+(\cdot)$ satisfies QNI (quantitative
non-integrability) if the following holds. There exists
\begin{itemize}
\item $\index{$\alpha_0$}\alpha_0 > 0$ and,
\item For every $\delta > 0$ there exists a compact set $K$ with
  $\nu(K) > 1-\delta$, and constants $C = C(\delta) > 0$ and $\ell_0 =
  \ell_0(\delta) > 0$ such that:

  If $\ell > \ell_0$ and $q_{1/2} \in K$, $q_1 \equiv g_{\ell/2}
  q_{1/2} \in K$ and $q \equiv g_{-\ell/2} q_{1/2} \in K$ then

\item There is a subset $S = S(q_{1/2},\ell) \subset \cB_{\ell/2}[q_{1/2}]$ 	with $|S| \ge
  (1-\delta) |\cB_{\ell/2}[q_{1/2}|]$ (recall that \index{$|\bullet|$}$|\bullet|$ is the Haar measure on $U^+[q_1/2]$).
  \linebreak
  Furthermore for each $y_{1/2} \in S$:

\item
There is a subset $S'= S'(y_{1/2},q_{1/2},\ell) \subset \gB^-_{\ell/2}[q_{1/2}]$ with   $\nu^s_{q_{1/2}}(S') \ge (1-\delta) \nu^s_{q_{1/2}}(\gB^-_{\ell/2}[q_{1/2}])$ so that if $q'_{1/2} \in S'$
\begin{equation}
	\label{eq:def:QNI}
	d^Q_{min}\left(\cB_0[q'_{1/2}], \cW^{cs}_{loc}[y_{1/2}]\right) \ge C e^{-\alpha_0 \ell}, 
\end{equation}
where $d^Q_{min}(\cdot, \cdot)$ is defined in \S\ref{sec:subsec:distances}.
\end{itemize}
\end{definition}

\begin{remark}
The center-stable manifolds in \autoref{eq:def:QNI} are defined in
\autoref{sssec:center_stable_manifolds_factorization}. It is easy to
see that for $\ell$ large enough and on a set of $y_{1/2}$ of measure
arbitrarily close to $1$, the validity of \autoref{eq:def:QNI} depends
on the center-stable manifold $\cW^{cs}_{loc}[y_{1/2}]$ only via its 
its $k$-jet for some $k \gg \alpha_0$. Thus the arbitrary choice in
\autoref{sssec:center_stable_manifolds_factorization} does not affect
the validity of the QNI condition. 
\end{remark}

\subsubsection{Choice of constants}
	\label{sssec:choice_of_constants}
The QNI condition in \autoref{def:QNI} provides a constant $\alpha_0>0$.
We will then choose constants $\alpha_2$ in \autoref{lemma:lower:bound:on:3cA} and $\alpha_3$ in \autoref{sec:subsec:choice:of:alpha3}.
Finally \autoref{sec:subsec:choice:of:N} provides an integer $N>0$ which will give a degree of approximation.




\begin{remark}
\label{remark:place:QNI:used}
The QNI condition is used only in the proof of
\autoref{lemma:lower:bound:on:3cA}.   
\end{remark}



\subsection{Random Dynamics: subgroups compatible with the measure and the QNI Condition}
	\label{ssec:random:subgroups_compatible_with_the_measure}
In this section we restate the definitions of \autoref{ssec:subgroups_compatible_with_the_measure} and \autoref{sec:subsec:conditions_for_extra_invariance} for the case of random dynamics.
Let $\hat{T}^t: \hat{X} \to \hat{X}$ be as in
\S\ref{sec:subsec:skew_product_measurable_cover}. 
As in \autoref{sssec:stabilizers_and_subgroups}, let
\index{$G$@$\bbG^{ssr}_x(x,\omega,s)$}$\bbG^{ssr}_x(x,\omega,s)$ denote the stabilizer of $x \in
\cW^u[x,\omega.s]$ inside $\bbG^{ssr}(x,\omega,s)$.  

As above the strictly subresonant maps $\bbG^{ssr}(x,\omega,s)$ are a unipotent algebraic group, and we will study unipotent algebraic subgroups $\index{$U^+(x,\omega,s)$}U^+(x,\omega,s)\subseteq \bbG^{ssr}(x,\omega,s)$, which we identify with their Lie algebras \index{$u$@$\fraku^+(x,\omega,s)$}$\fraku^+(x,\omega,s)$.
For such subgroups, we will also denote by \index{$U^+_x(x,\omega,s)$}$U^+_x(x,\omega,s)$ and \index{$u$@$\fraku^+_x(x,\omega,s)$}$\fraku^+_x(x,\omega,s)$ the group and Lie algebra of the stabilizer of $x\in \cW^u[x,\omega,s]$.

Additionally, for a subgroup $U^+(x,\omega,s)$ we will denote by $\index{$U^+[q,\omega,s]$}U^+[q,\omega,s]:=U^+(x,\omega,s)\cdot q\subset \cW^u[x,\omega,s]$ the corresponding orbit of $U^+$ on the unstable manifold.
As a $U^+(x,\omega,s)$-homogeneous space, it is isomorphic to $U^+(x,\omega,s)/U^+_x(x,\omega,s)$.
Note that since both groups are unipotent (as subgroups of the
unipotent $\bbG^{ssr}$), they are unimodular and hence there is a
unique up to scaling Haar measure on the quotient. 

The following is the analogue of
\autoref{def:compatible_family_of_subgroups} (which in turn follows
 \cite[Def.~6.2]{EskinMirzakhani_Invariant-and-stationary-measures-for-the-rm-SL2Bbb-R-action-on-moduli-space}):

\begin{definition}[Random Compatible family of subgroups]
\label{def:random:compatible_family_of_subgroups}
Let $\nu$ be a ergodic $\mu$-stationary measure on $Q$, and let
\index{$\nu$@$\hat{\nu}$}$\hat{\nu}$ be  ergodic, $\hat{T}^t$-invariant probability measure on $\hat X$ 
given by (\ref{eq:def:hat:nu}). 
We say that a measurable family of connected algebraic subgroups
$(x,\omega,s) \mapsto U^+(x, \omega,s)\subset \bbG^{ssr}(x,\omega,s)$ is \emph{compatible with $\hat{\nu}$} if the following hold:
\begin{enumerate}
\item The subgroups are $\hat{T}^t$-equivariant.
\item Let $\gB_0$ be the partition of $X$ defined in
  \autoref{ssec:measurable_partitions} and
  \S\ref{sec:subsec:pulling:back:Markov}. 
  Then for almost all $\omega$,
  the sets of the form $U^+[x,\omega,s]\cap \gB_0[x,\omega,s]$ form a partition of $Q$.
  
  Furthermore, the conditional measures along $U^+[x,\omega,s]\cap \gB_0[x,\omega,s]$ are proportional to the restriction of the unique $U^+(x,\omega,s)$-invariant measure on $U^+[x,\omega,s]$, which will be called the Haar measure on $U^+[x,\omega,s]$.

\item Given $(x',\omega',s')\in \widehat{\cW}^u[x,\omega,s]$, we have natural
  identification of $\bbG^{ssr}(x,\omega,s)$ and $\bbG^{ssr}(x',\omega',s')$, viewing both as acting on $\cW^u[x,\omega,s]=\cW^u[x',\omega',s']$.
By further restricting to a set $\hat{X}_0 \subset \hat{X}$ of full measure, if $(x',\omega',s')$ and
$(x,\omega,s) \in \hat{X}_0$ and $(x',\omega',s') \in U^+[x,\omega,s]$, then
$U^+(x',\omega',s)=U^+(x,\omega,s)$.
		
\item $U^+(x,\omega,s)$ is maximal in the sense that it is the largest
  (connected) subgroup of $\bbG^{ssr}(x,\omega,s)$ that preserves the Haar measure on $U^+[x,\omega,s]$.
\end{enumerate}
\end{definition}

\begin{remark} In the random dynamics setting, $U^+(x,\omega,s) = \{ e
  \}$ is allowed. 
\end{remark}
  
As in the $SL(2,\reals)$ action setting, this is meant to formalize
the following:
\begin{remark}[On measures with local invariance]
	\label{rmk:random:on_measures_with_local_invariance}
The above \autoref{def:random:compatible_family_of_subgroups} is meant to formalize the following situation.
	Suppose $\cW^u[x,\omega,s]$ is equipped with its leafwise (conditional) measure $\hat{\nu}^u_{(x,\omega,s)}$.
	We would like to consider the possibility that $\hat{\nu}^u_{(x,\omega,s)}$-a.e. point $(y,\omega',s')$ has associated to it a group $U^+(y,\omega',s')\subset \bbG^{ssr}(\cW^u[x,\omega,s])$, but the groups $U^+(y,\omega',s')$ need not be the same.
	Nonetheless, the above definition expresses the possibility that the measure $\hat{\nu}^u_{(x,\omega,s)}$ is ``Haar along $U^+$-orbits''.
\end{remark}

Let $(x,\omega,s)\mapsto U^+(x,\omega,s)\subset \bbG^{ssr}(x,\omega,s)$ be   
a measurable family of connected algebraic subgroups which is
compatible with $\hat{\nu}$ in the sense of
\autoref{ssec:random:subgroups_compatible_with_the_measure}.

\begin{definition}[Random QNI]
\label{def:Random:QNI}
For $(q,\omega,s) \in \hat X$, let
\begin{displaymath}
\hat{U}^+[q,\omega,s] = W_1^+[\omega,s] \cross U^+[q,\omega,s] \subset \hat{\cW}^u[q,\omega,s].  
\end{displaymath}
We say that $U^+(\cdot,\cdot)$ satisfies Random QNI
if the following holds. There exists
\begin{itemize}
\item $\index{$\alpha_0$}\alpha_0 > 0$ and,
\item For every $\delta > 0$ there exists a compact set $K$ with
  $\hat{\nu}(K) > 1-\delta$, and constants $C = C(\delta) > 0$ and $\ell_0 =
  \ell_0(\delta) > 0$ such that:

  If $\ell > \ell_0$ and $\hat{q}_{1/2} = (q_{1/2},\omega_{1/2},s_{1/2}) \in K$, $\hat{q}_1 = \hat{T}^{\ell/2} \hat{q}_{1/2} \in K$ and $\hat{q} \equiv \hat{T}^{-\ell/2} \hat{q}_{1/2} \in K$ then

\item There is a subset $S = S(\hat{q}_{1/2},\ell) \subset \cB_{\ell/2}[\hat{q}_{1/2}]$ 
  with $|S| \ge
  (1-\delta) |\cB_{\ell/2}[\hat{q}_{1/2}]|$ (where $|\bullet|$ is Haar measure on $U^+[\hat{q}_{1/2}]$).
 And for each $\hat{y}_{1/2} \in S$:
\item There is a subset $S'= S'(\hat{y}_{1/2},\ell) \subset \gB^-_{\ell/2}[\hat{q}_{1/2}]$ with   $\hat{\nu}^s_{\hat{q}_{1/2}}(S') \ge (1-\delta) \hat{\nu}^s_{\hat{q}_{1/2}}(\gB^-_{\ell/2}[\hat{q}_{1/2}])$ so that if $\hat{q}'_{1/2} \in S'$
\begin{equation}
\label{eq:def:Random:QNI}
 d^Q_{min}(\cB_{\ell/2}[\hat{q}'_{1/2}], \cW^{cs}_{loc}[\hat{y}_{1/2}]) \ge C e^{-\alpha_0 \ell}, 
\end{equation}
where $d^Q_{min}(\cdot, \cdot)$ is defined in \S\ref{sec:subsec:distances}.
\end{itemize}
\end{definition}




\section{Main Theorems and Proofs}
	\label{sec:main_theorems_and_proofs}


\subsection{The Main Theorems}
\label{sec:subsec:the:main:theorems}

\subsubsection{Extra Invariance}
Our main result (which applies even if we have a single diffeomorphism or a flow) is the following:

\begin{theorem}
\label{thm:inductive:step}
Suppose $\nu$ is an ergodic, $g_t$-invariant measure on $Q$. 
Suppose $x \mapsto U^+(x)$ is a nontrivial 
family of subgroups of $\bbG^{ssr}(x)$ compatible
with $\nu$ in the sense of \autoref{def:compatible_family_of_subgroups}.  

Suppose the family $x \to
U^+(x)$ satisfies the QNI condition (see
\autoref{def:QNI}).  Then, there
exists a family of  subgroups 
$U_{new}^+(x)$ of $\bbG^{ssr}(x)$ compatible with $\nu$ in the sense of
\autoref{def:compatible_family_of_subgroups} such that for almost all
$x$ we have  $U^+(x) \subset U_{new}^+(x) $ and $U^+(x) \neq  U_{new}^+(x) $.
\end{theorem}

All of our other results are derived from \autoref{thm:inductive:step} (and its version for skew-products \autoref{thm:random:inductive:step}). 

We note the following immediate:
\begin{corollary}
\label{cor:end:induction}
Suppose $\nu$ is an ergodic, $g_t$-invariant measure on $Q$.
Suppose there exists 
a family of subgroups of $\bbG^{ssr}(x)$ compatible
with $\nu$ in the sense of \autoref{def:compatible_family_of_subgroups}.  

Then there exists (another) family $x \to
U^+(x)$ of subgroups compatible with $\nu$ such that the
the QNI condition \autoref{def:QNI} does not hold for $U^+$. 
\end{corollary}

\begin{proof}
  This is obtained by repeatedly applying \autoref{thm:inductive:step}.
\end{proof}

\autoref{thm:inductive:step} also applies to skew-products. More precisely, we have the following:
\begin{theorem}
\label{thm:random:inductive:step}
Suppose $\nu$ is an ergodic, $\mu$-stationary measure on $Q$, and let $\hat{\nu}$ be as in (\ref{eq:def:hat:nu}). 
Suppose $\hat{x} \to U^+(\hat{x})$ is a family of subgroups of $\bbG^{ssr}(x)$ compatible
with $\hat{\nu}$ in the sense of \autoref{def:random:compatible_family_of_subgroups}.
(Note that $U^+ = \{ e \}$ is allowed).  Suppose the family $q \to
U^+(q)$ satisfies the Random QNI condition (see
\autoref{def:Random:QNI}).  Then, there
exists a family of  subgroups 
\index{$U_{new}^+(x)$}$U_{new}^+(x)$ of $\bbG^{ssr}(x)$ compatible with $\hat{\nu}$ in the sense of
\autoref{def:random:compatible_family_of_subgroups} such that for almost all
$x$ we have  $U^+(x) \subset U_{new}^+(x) $ and $U^+(x) \neq  U_{new}^+(x) $.
\end{theorem}

As above, we have the following immediate corollary:
\begin{corollary}
\label{cor:random:end:induction}
Suppose $\nu$ is an ergodic, $\mu$-stationary measure on $Q$, and let $\hat{\nu}$ be as in (\ref{eq:def:hat:nu}). 
Then there exists a family $x \to
U^+(x)$ of subgroups compatible with $\nu$ such that the
the Random QNI condition 
\autoref{def:Random:QNI}) does not hold for $U^+$. 
\end{corollary}

\subsubsection{$\SL(2,\bR)$-actions}

Our main results for general $\SL(2,\R)$ actions are
\autoref{theorem:SL2R:P:invariant:is:SL2R:invariant} and
the following theorem:

\begin{theorem}
\label{theorem:SL2R:measure:classification}
Suppose $\SL(2,\reals)$ acts on $Q$, and suppose $\nu$ is an $\SL(2,\R)$-invariant measure on $Q$. Suppose there is no $P$-invariant $\nu$-measurable subbundle of $TQ$ on which the sum of the Lyapunov exponents is negative. Then
\begin{itemize}
\item[{\rm (i)}] For almost all $q \in Q$,
there exists a subgroup $U^+(q)$ of
$\bbG^{ssr}(\cW^u[q])$ such that the conditional of $\nu$ along $\cW^u[q]$ is the Haar measure along the orbit $U^+[q] \equiv U^+(q)q$.
\item[{\rm (ii)}] For almost all $q \in Q$ there exists a subgroup \index{$U^-(q)$}$U^-(q)$ of
$\bbG^{ssr}(\cW^s[q])$ such that the conditional of $\nu$ along $\cW^s[q]$ is the Haar measure along the orbit $\index{$U^-[q]$}U^-[q] \equiv U^-(q)q$. 
\item[{\rm (iii)}] The sum of the Lyapunov exponents of the bundle $T_qU^+[q] \oplus T_q U^-[q]$ is $0$. 
\item[{\rm (iv)}] The family of subgroups $q \to U^+(q)$ is compatible
  with the measure (see
  \autoref{def:compatible_family_of_subgroups}) and the 
  QNI condition \autoref{def:QNI} does not
  hold.
\end{itemize}
\end{theorem}

\subsubsection{Random Dynamics}

Our main results for random dynamics are
\autoref{theorem:random:stationary:is:invariant} and the following:
\begin{theorem}
\label{theorem:random:measure:classification}
Suppose $\mu$ is as in \autoref{theorem:random:stationary:is:invariant}. 
Suppose $\nu$ is an ergodic $\mu$-invariant measure on $Q$, and suppose that
there is no $\mu$-invariant $\nu$-measurable subbundle of $TQ$ (see \autoref{def:mu_invariant_subbundle}) on which the sum of the Lyapunov exponents is negative. Then,
\begin{itemize}
\item[{\rm (i)}] For almost all $(q,\omega,s) \in \hat{Q}$,
there exists a subgroup $U^+(q,\omega,s)$ of
$\bbG^{ssr}(\cW^u[q,\omega,s])$ such that the conditional of $\hat{\nu}$
along $\cW^u[q,\omega,s]$ is the Haar measure along the orbit
$U^+[q,\omega,s] \equiv U^+(q,\omega,s)q$.
\item[{\rm (ii)}] For almost all $(q,\omega,s) \in \hat{Q}$,
there exists a subgroup \index{$U^-(q,\omega,s)$}$U^-(q,\omega,s)$ of
$\bbG^{ssr}(\cW^s[q,\omega,s])$ such that the conditional of $\hat{\nu}$
along $\cW^s[q,\omega,s]$ is the Haar measure along the orbit
$\index{$U^-[q,\omega,s]$}U^-[q,\omega,s] \equiv U^-(q,\omega,s)q$.
\item[{\rm (iii)}] The sum of the Lyapunov exponents of the bundle $T_qU^+[q,\omega,s] \oplus T_q U^-[q,\omega,s]$ is $0$. 
\item[{\rm (iv)}] The family of subgroups $(q,\omega,s) \to U^+(q,\omega,s)$ is compatible
  with the measure (see
  \autoref{def:random:compatible_family_of_subgroups}) and the 
  Random QNI condition \autoref{def:Random:QNI} does not
  hold.

\end{itemize}
\end{theorem}

The following theorem summarizes our progress towards Conjecture~\ref{conj:symplectic}.

\begin{theorem}
\label{theorem:symplectic}
Suppose $Q$ is a symplectic manifold, and suppose $\mu$ is a finitely
supported probability measure on the group of symplectic diffeomorphisms of $Q$,
which satisfies uniform expansion on all the isotropic subspaces of $TQ$. Suppose $\nu$ is a $\mu$-stationary measure on $Q$. Then $\nu$ is $\mu$-invariant. Furthermore, $T_q U^+[q] \oplus T_q U^-[q]$ is symplectic and also
conclusions (i)-(iv) of \autoref{theorem:random:measure:classification} hold for $\nu$. 
\end{theorem}


\subsection{Proof of the main theorems for \texorpdfstring{$\SL(2,\R)$}{SL(2,R)} actions assuming \autoref{thm:inductive:step}.} 

\def \essspan{\mathrm{ess\, span}}
\subsubsection{Stable support}
	\label{sssec:stable_support_Ls}
Recall that $\nu^s_{q}$ denotes the conditional measure
of $\nu$ along $\cW^s[q]$.
For a.e.\ $q \in Q$,
 let \index{$L$@$\cL^s[q]$}$\cL^s[q] 
\subset \cW^s[q]$ denote the smallest real-algebraic subset
for which there exists $\epsilon(q) > 0$ such that the intersection of the 
ball of radius $\epsilon(q)$ with support of $\nu^s_q$ is contained in $\cL^s[q]$.
Then, $\cL^s[q]$ is $g_t$-equivariant.
Since
the action of $g_{-t}$ is expanding along $\cW^s[q]$, we see that for
almost all $q$ and any $\epsilon > 0$, $\cL^s[q]$ is the smallest
real-algebraic subset of $\cW^s[q]$ such that $\cL^s[q]$ contains
$\operatorname{support} (\nu^s_{q}) \cap B^{Q}(q,\epsilon)$.

\subsubsection{Smoothness and irreducibility}
	\label{sssec:smoothness_and_irreducibility}
Let $T_q\cL^s[q] \subset T_x\cW^s[q]$ be the tangent space at $q$ to $\cL^s[q]$.
Note that for $\nu$-a.e. $q$ the set $\cL^s[q]$ is smooth at $q$, since for an algebraic set, the subset of singular points is a proper algebraic subset.
Note that $T_q\cL^s[q]$ is $g_t$-invariant by construction.  
Furthermore, we note that $\cL^s[q]$ is irreducible.
Indeed, note that a.e. $q$ sits on a unique irreducible component of $\cL^s[q]$ (otherwise, we can replace $\cL^s[q]$ by the union of the intersection loci of any pair of irreducible components).
Therefore, $\cL^s[q]$ has to agree with that irreducible component for a.e. $q$.


\subsubsection{The bundle $\cZ$}
	\label{sssec:the_bundle_cz}
Given $q\in Q$, let $T_q Q = \index{$E^s(q)$}E^s(q) \oplus \index{$E^c( q)$}E^c( q) \oplus \index{$E^u( q)$}E^u( q) $ be the splitting into negative, zero, and positive Lyapunov exponents.
Let $\index{$\pi^{cs}_{q}$}\pi^{cs}_{q}\colon T_q Q \to E^u(q)$ be the canonical projection.

We let \index{$Z$@$\cZ(q)$}
\begin{equation}
\label{eq:def:cZ}
\cZ(q) = \index{$e$@$\essspan$}\essspan \{n\inv T_{n \cdot q}\cL^s[n\cdot q]  : {n\in N} \}.  
\end{equation}
(Here $\essspan$ denotes the essential span). 

Note that the bundle $\cZ$ is $P$-invariant. 
Let $\index{$n$@$\bar{\mathfrak{n}}(q)$}\bar{\mathfrak{n}}(q) \subset T_xQ$ denote the tangent space to the
orbit of $\bar{N} \subset \SL(2,\R)$, and let \index{$Z$@$\cZ'(q)$}
\begin{equation}
\label{eq:def:cZprime}
\cZ'(q) = \cZ(q) + \bar{\mathfrak{n}}(q).   
\end{equation}
Then the bundle $\cZ'$ is also $P$-invariant.

\subsubsection{Lifting to $X$}
	\label{sssec:lifting_to_x_main_proof}
Recall that the measure $\nu$ lifts to $X$ and we continue to denote it by $\nu$.
For $x \in X$ declare \index{$L$@$\cL^s[x]$}$\cL^s[x]$
to be $\cW^s[x] \cap \sigma^{-1}(\cL^s[\sigma(x)])$, where $\sigma\colon X
\to Q$ is the projection map.
Similarly, we let $\index{$T_x \cL^s[x]$}T_x \cL^s[x]:=T_{\sigma(x)}\cL^s[\sigma(x)]$ and similarly $\index{$Z$@$\cZ(x)$}\cZ(x):=\cZ(\sigma(x))$ and $\index{$Z$@$\cZ'(x)$}\cZ'(x):=\cZ'(\sigma(x))$.
The center-stable projection $\index{$\pi^{cs}_x$}\pi^{cs}_x$ is defined analogously.

Let $U^+( x)\subset \bbG^{ssr}(x)$ be as in \autoref{def:compatible_family_of_subgroups} and denote by $U^+[x] = U^+( x)\cdot x$ the $U^+(x)$-orbit of $x$. 

\begin{lemma}
\label{lemma:equivalent:stopping:conditions}
The following are equivalent:
\begin{itemize}
\item[{\rm (1)}] For a positive (and thus full) $\nu$-measure subset
  of $x\in X$ we have 
\begin{equation}
\label{eq:cZ:subset:TU}
  \pi^{cs}_{x} \left( \cZ(x) \right) \subset T_xU^+[x].
\end{equation}
\item[{\rm (2)}] For a positive (and thus full) $\nu$-measure subset of $x\in X$ there exists a positive (and thus full) Lebesgue measure subset $E_x \subset N$ such that for $n \in E_x$ we have
	\[
		\pi^{cs}_{nx}\left(n \cdot T_x\cL^s\left(x\right)\right) \subset T_{nx} U^+[nx].
	\]
\end{itemize}      
\end{lemma}

We remark that to go from a positive to a full measure subset of $x
\in X$ in (1) and (2) we use the ergodicity of $g_t$. To go from a
positive to a full measure subset of $n \in N$ in (2) we use
\cite[Lemma~8.3]{EskinMirzakhani_Invariant-and-stationary-measures-for-the-rm-SL2Bbb-R-action-on-moduli-space}
(with $U^+$ replaced by $N$).

\begin{proof}
Suppose (1) holds. Then, in view of the definition of $\cZ(x)$, for a.e.\ $x \in X$ and a.e.\ $n \in N$, $n^{-1}\cdot T_{nx}\cL^s[nx] \subset \cZ(x)$. Then,
letting $y = n x$, we have, for a.e.\ $y \in X$ and a.e.\ $n \in N$, 
$n^{-1} \cdot T_{y} \cL^s[y] \subset \cZ(n^{-1} y)$. Replacing $y$ by $x$ and $n$ by $n^{-1}$, we get for a.e.\ $x \in X$ and a.e.\ $n \in N$, 
\begin{displaymath}
n \cdot T_x \cL^s[x] \subset \cZ(n x).
\end{displaymath}
Thus, for a.e.\ $x \in X$ and a.e.\ $n \in N$, 
\begin{displaymath}
\pi^{cs}_{nx}(n \cdot T_x \cL^s[x]) \subset \pi_{nx}^{cs}(\cZ(nx)) \subset T_{nx} U^+[nx],
\end{displaymath}
where we have used the condition (1) for the final inclusion. Thus (2) holds.

Now suppose (2) holds. Then, 
for a.e.\ $y \in X$, for a.e.\ $n \in N$, 
\begin{displaymath}
\pi^{cs}_{ny}(n \cdot T_y \cL^s[y]) \subset T_{ny} U^+[ny].
\end{displaymath}
Using Fubini's theorem
we can let $ny = x$. Then, for a.e. $x \in X$, and a.e. $n \in N$,
\begin{displaymath}
\pi^{cs}_x(n \cdot T_{n^{-1} x} \cL^s[n^{-1}x]) \subset T_x U^+[x].
\end{displaymath}
Replacing $n$ by $n^{-1}$, we get for a.e.\ $x \in X$, for a.e. $n \in N$, 
\begin{displaymath}
\pi^{cs}_x(n^{-1} \cdot T_{n x} \cL^s[n x]) \subset T_x U^+[x].
\end{displaymath}
Since this is true for a fixed $x$ for a.e. $n \in N$, we get using the definition of $\cZ(x)$, 
\begin{displaymath}
\pi^{cs}_x(\cZ(x)) \subset T_x U^+[x].
\end{displaymath}
Thus (1) holds. 
\end{proof}

\begin{lemma}
\label{lemma:AWB:to:QNI}
Suppose one of the two equivalent conditions of \autoref{lemma:equivalent:stopping:conditions} is not satisfied.
Then the family $x \to	U^+(x)$ satisfies the QNI condition (see \autoref{def:QNI}).
\end{lemma}

\begin{proof}
See \cite{ElliotSmith}.
\end{proof}

\subsubsection{Facts from entropy theory for flows}
Let \index{$g_1$}$g_1$ denote $g_t$ for $t=1$; similarly, let
\index{$g_{-1}$}$g_{-1}$ denote $g_t$ for $t=-1$. Let
\index{$h_\nu(g_t)$}$h_\nu(g_t)$ 
denote the entropy of the transformation $g_t$ with respect to the
invariant measure $\nu$. For almost all $q \in Q$ we have the Lyapunov
splitting for $g_t$ (see \autoref{thm:tempered_oseledets})
\begin{displaymath}
T_q Q = \oplus_i \index{$E^{\lambda_i}(q)$}E^{\lambda_i}(q),
\end{displaymath}
where the Lyapunov exponents \index{$\lambda_i$}$\lambda_i$ satisfy $\lambda_1 >
\lambda_2 > \dots$. For $x \in \sigma^{-1}(q)$, we define
$\index{$E^{\lambda_i}(x)$}E^{\lambda_i}(x) = E^{\lambda_i}(q)$.

\begin{remark}[Entropy and finite covers]
	\label{rmk:entropy_and_finite_covers}
Recall that $X\to Q$ is a finite cover, with measures and flows denoted $\nu,g_t$ on both spaces.
Given a measurable partition $\xi$ of $ Q$ whose atoms are contained in $\cW^u$-manifolds, we write $\tilde \xi$ for the associated measurable partition of $ X$ such that $\tilde \xi[x] \to \xi[x]$ is a bijection for every $ x\in  Q$.
Then $\tilde \xi$ is a finite refinement of $\sigma\inv (\xi)$.  

Suppose now that $\xi$ is $g_t$-invariant.
Then it follows that $h_{ \nu}( g_t , \tilde  \xi) = h_{  \nu}( g_t ,  \xi) $.  Indeed, we have 
$h_{ \nu}( g_t , \tilde \xi) \ge h_{  \nu}( g_t ,  \xi) $ since $\tilde \xi$ is a finite refinement of $\sigma\inv \xi$; the reverse inequality follows from \cite[Lem.\ 6.1]{Katok-RH-2010}.

In particular, $h_{\nu}(g_t)$ is independent of whether we work on $X$ or $Q$.
\end{remark}

\begin{proposition}
	\label{prop:LYfactsgt}
For $\nu$-a.e.\ $x\in X$, the following hold:
\begin{enumerate}
\item \label{LYfact1gt} $h_{ \nu}(g_1)\ge 
 \sum_{\lambda_i>0} \lambda_i \dim E^{\lambda_i}(x)\cap T_x U^+[ {x}]$.  
\item \label{LYfact3gt} $h_{ \nu}(g_{-1})\le 
\sum_{\lambda_i<0}- \lambda_i \dim E^{\lambda_i}({x})\cap T_x \cL^s[ x].$


\item\label{LYfact5gt} We have $h_{ \nu}(g_{-1})= h_{ \nu} (g_{-1},  \gB_0^-)$.  

\end{enumerate}
\end{proposition}
\noindent Note that the last two assertions also make sense on $Q$, and follow on $X$ from the corresponding statement on $Q$, see \autoref{rmk:entropy_and_finite_covers}.
\begin{proof}
	We note that \autoref{LYfact5gt} follows from \autoref{prop:full_entropy} (applied to time-reversed flow $g_{-t}$).
	Next, \autoref{LYfact3gt} follows from \cite[Thm.~7.2]{Brown_RodriguezHertz_Wang} combined with \autoref{LYfact5gt}.
	
	Finally \autoref{LYfact1gt} follows from the same \cite[Thm.~7.2]{Brown_RodriguezHertz_Wang} by noting that $h_{\nu}(g_{1})\geq h_{\nu}(g_1;\frakU)$ where $\frakU$ is the measurable partition obtained by refining the $U^+[x]$-orbits by $\gB_0$.
	We note that the proof of \cite[Thm.~7.2]{Brown_RodriguezHertz_Wang} holds in our setting, for the lift of the $\cW^u$-foliation to the measurable finite cover $X$.
\end{proof}

\begin{proposition}\label{prop:LedrappierFactsgt}
Suppose for $  \nu$-a.e. $ q\in Q $, that  $\nu^u_{q}$ is supported on an embedded submanifold $\index{$M_{q}$}M_{q}\subset \cW^u[   q]$.  
Then 
$$h_{\nu}(g_1)\le \sum_{\lambda_i>0} \lambda_i \dim E^{\lambda_i}( {q})\cap T_q M_{ {q}}.$$
Moreover, equality holds if and only if 
$ \nu^u_{  q}$ equivalent to the Lebesgue measure on $M_{  q }$ for a.e.\ $ q\in Q$.  
\end{proposition}

\begin{proof}
	This follows again from \cite[Thm.~7.2]{Brown_RodriguezHertz_Wang}.
\end{proof}

\begin{proposition}\label{prop:Ledrappier_gt_newinvariance}

Suppose that $\nu$ is an ergodic $g_t$-invariant probability measure on $Q$.  
Suppose for a.e.\ $q\in Q$ that the conditional measure $\nu^{\overline N}_q$ along the $\overline N$-orbit $\overline N\cdot q$ is absolutely continuous with respect to the Haar measure along the orbit $\overline N\cdot q$.  Then $\nu$ is $\overline N$-invariant.


\end{proposition}
\begin{proof}
Indeed, it follows exactly as in the proof of \cite[Thm.\ 3.4]{Ledrappier-84} that  conditional measure $\nu^{\overline N}_q$ along a.e.\ $\overline N$-orbit is equivalent to the Haar measure with a constant density.  This implies $\overline N$-invariance.  \end{proof}

\begin{proposition}[Entropy balancing]
\label{prop:entropy:trick}
Suppose that
\begin{equation}
\label{eq:cZprime:non:negative:sum}
\text{the sum of the Lyapunov exponents along $\cZ'$ is non-negative,}
\end{equation}
where $\cZ'$ is as in \autoref{eq:def:cZprime}.
Furthermore, suppose that for $\nu$-a.e.\ $x \in X$,
\autoref{eq:cZ:subset:TU} holds. 

Then $\nu$ is $\SL(2,\R)$-invariant.
Furthermore,
\begin{itemize}
\item[{\rm (i)}] $U^+(x)$ depends only on $\sigma(x)$, where $\sigma: X \to Q$ is the covering map of \autoref{ssec:the_finite_measurable_cover}. 
\item[{\rm (ii)}] For almost all $q \in Q$, the conditional of $\nu$ along $\cW^u[q]$ is the Haar measure along the orbit $U^+[q]$. 
\item[{\rm (iii)}] For almost all $q \in Q$ there exists a subgroup \index{$U^-(q)$}$U^-(q)$ of $\bbG^{ssr}(\cW^s[q])$ such that the conditional of $\nu$ along $\cW^s[q]$ is the Haar measure along the orbit $\index{$U^-[q]$}U^-[q] \equiv U^-(q)q$. In particular, $\cL^s[q] = U^-[q]$. 
\item[{\rm (iv)}] The sum of the Lyapunov exponents of the bundle $\cZ'$ is $0$. 
\end{itemize}
\end{proposition}

\begin{proof}
Since for $t > 0$ the action of $g_t$ on $\overline{N}$ is contracting, $\pi^{cs}_{x} \left( \cZ(x) \right) = \pi^{cs}_{x} \left( \cZ'(x) \right)$. Therefore, in view of (\ref{eq:cZ:subset:TU})
\begin{equation}
\label{eq:cZprime:subset:TU}
\pi^{cs}_{x} \left( \cZ'(x) \right) = \cZ'(x)\cap T_x \cW^u[x]
\subset T_xU^+[x].
\end{equation}


Let $g_1$ denote $g_t$ for $t=1$.  
From \autoref{prop:LYfactsgt}\autoref{LYfact1gt}, we have a lower bound
\begin{align*}
h_\nu(g_1)
& \ge \sum_{\lambda_i>0} \lambda_i \dim E^{\lambda_i}(x)\cap T_x U^+[x]\\
& \ge \sum_{\lambda_i>0} \lambda_i \dim E^{\lambda_i}(x)\cap \cZ'(x).
\end{align*}
where the final inequality follows from \autoref{eq:cZprime:subset:TU}.    

Similarly, by the definition of $\cZ'(x)$
we have $T_x(\overline{N}\cL^s[x])\subset \cZ'(x)$ whence by \autoref{prop:LYfactsgt}\autoref{LYfact3gt},
\begin{align*}
h_\nu(g_{-1}) & \le \sum_{\lambda_i<0} - \lambda_i \dim E^{\lambda_i}(x)\cap T_x(\overline{N}\cL^s[x])
\\&\le \sum_{\lambda_i<0} - \lambda_i \dim E^{\lambda_i}(x)\cap \cZ'(x).
\end{align*}
In view of (\ref{eq:cZprime:non:negative:sum}), 
\begin{align*}
0
&\le  \sum_{\lambda_i} \lambda_i \dim E^{\lambda_i}(x)\cap \cZ'(x)\\
&=  \sum_{\lambda_i<0} \lambda_i \dim E^{\lambda_i}(x)\cap \cZ'(x)
  +\sum_{\lambda_i>0} \lambda_i \dim E^{\lambda_i}(x)\cap \cZ'(x)\\
&\le - h_\nu(g_{-1}) + h_\nu(g_1)
\end{align*}
Since $h_\nu(g_{-1}) =h_\nu(g_1)$, we conclude that the sum of the Lyapunov exponents along $\cZ'$ is $0$. Also, 
\begin{multline}
\label{eq:equal:entropy}
\sum_{\lambda_i<0} - \lambda_i \dim E^{\lambda_i}(x)\cap T_x(\overline{N}\cL^s[x])  = h_\nu(g_{-1}) = \\ =h_\nu(g_1)=  \sum_{\lambda_i>0} \lambda_i \dim E^{\lambda_i}(x)\cap T_xU^+[x].
\end{multline}
From  \autoref{prop:LedrappierFactsgt}  (for $g_{-1}$), we conclude that stable conditionals measure $\nu^s_x$ is absolutely continuous along  $\overline{N}\cL^s[x]$; since the foliation by $\overline N$-orbits is smooth, we conclude the conditionals measure $\nu^{\overline N}_x$ is absolutely continuous along the orbit $\overline N \cdot x$.  By \autoref{prop:Ledrappier_gt_newinvariance}, we conclude that $\nu$ is $\overline N$-invariant; this in particular, implies $\overline{N} \cdot \cL^s[x] = \cL^s[x]$.

Let $\index{$L$@$\cL^u[x]$}\cL^u[x] \subset \cW^u[x]$ denote the smallest algebraic subset of full $\nu^u_x$ measure. Note that since the measure on $X$ is a lift of the measure on $Q$, $\cL^u[x]$ depends only on $\sigma(x)$ where $\sigma: X \to Q$ is the projection map.

There exists $w \in \SL(2,\R)$ such that $w g_1 w^{-1} = g_{-1}$. Since $\nu$ is now $\SL(2,\R)$-invariant, we get, for a.e.\ $x \in X$, $\nu^u_x = w \nu^{s}_{w^{-1} x}$. In particular, $\cL^u[x] = w \cL^s[w^{-1}x]$. Thus, for a.e.\ $x \in X$, $\nu^u_x$ is absolutely continuous on $\cL^u[x]$. Note that $U^+[x] \subset \cL^u[x]$.
By \autoref{prop:LedrappierFactsgt}, we have $$h_\nu(g_1)=  \sum_{\lambda_i>0} \lambda_i \dim E^{\lambda_i}(x)\cap T_x\cL^u[x] 
\ge  \sum_{\lambda_i>0} \lambda_i \dim E^{\lambda_i}(x)\cap T_xU^+[x].$$
\autoref{eq:equal:entropy} then implies that 
$\dim U^+[x] = \dim \cL^u[x]$.   
Thus, $U^+[x] = \cL^u[x]$ since $U^+[x]$ is an algebraic subset of $\cW^u[x]$ (see \autoref{sssec:stabilizers_and_subgroups}) and since $\cL^u[x]$ is irreducible (see \autoref{sssec:smoothness_and_irreducibility}).
In particular $U^+[x] = U^+[\sigma(x)]$.

Let $U^-(x) = w U^+(w^{-1} x) w^{-1}$. 
Since $w\in \SL(2,\R)$, $\nu$ is $\SL(2,\R)$-invariant, and since $wg_tw^{-1} = g_{-t}$, the action by $w$ intertwines stable and unstable manifolds and thus  also intertwines the subresonant structures on stable and unstable manifolds.  In particular,  $U^-(x)$ is a subgroup of $\bbG^{ssr}(\cW^s[x])$ 
and $\cL^s[x] = U^-(x) x \equiv U^-[x]$. 
\end{proof}
  

\begin{proof}[Proof of \autoref{theorem:SL2R:P:invariant:is:SL2R:invariant} and \autoref{theorem:SL2R:measure:classification}.]
In view of \autoref{cor:end:induction}, we may assume that
$x \to U^+(x)$ is a family of subgroups compatible with the measure
(see \autoref{def:compatible_family_of_subgroups}) such that the QNI condition does not hold. Then, in view of 
\autoref{lemma:random:AWB:to:QNI},  (\ref{eq:cZ:subset:TU}) holds. 
Let $\cZ'$ be as in (\ref{eq:def:cZprime}). Then, $\cZ'$ is $P$-invariant.
Therefore, in view of the assumption
in \autoref{theorem:SL2R:P:invariant:is:SL2R:invariant} or \autoref{theorem:SL2R:measure:classification}, (\ref{eq:cZprime:non:negative:sum}) holds. 
Now all the conclusions follow from \autoref{prop:entropy:trick}. 
\end{proof}

\subsubsection{The non-negative Laplacian condition}

\begin{proposition}
\label{prop:SL2R:nonnegative:laplacian:implies:no:negative:subbundle}
Suppose $\SL(2,\R)$ acts smoothly on $Q$, let $\nu$ be a $P$-invariant measure, and suppose \autoref{def:SL2R:nonnegative:laplacian} holds.
Then,
there is no $P$-invariant $\nu$-measurable subbundle of $TQ$ on which the sum of the Lyapunov exponents is negative.
\end{proposition}

\begin{proof}
This is essentially \cite[Theorem A.3]{EskinMirzakhani_Invariant-and-stationary-measures-for-the-rm-SL2Bbb-R-action-on-moduli-space}.
\end{proof}

\begin{proof}[Proof of \autoref{theorem:intro:SL2R:P:invanriant:is:SL2R:invariant}] This follows immediately from \autoref{theorem:SL2R:P:invariant:is:SL2R:invariant} and \autoref{prop:SL2R:nonnegative:laplacian:implies:no:negative:subbundle}.
\end{proof}

\subsection{Proof of the main theorems for Random Dynamics assuming
  \autoref{thm:random:inductive:step}.}

\subsubsection{Auxiliary bundles in the endgame}
For a.e.\ $\hat{q} \in \hat{Q}$,
 let \index{$L$@$\cL^s[\hat{q}]$}$\cL^s[\hat{q}] 
\subset \cW^s[\hat{q}]$ denote the smallest real-algebraic subset
containing, for some $\epsilon > 0$, the intersection of the 
ball of radius $\epsilon$ with
the support of the measure \index{$\nu$@$\hat{\nu}^s_{\hat{q}}$}$\hat{\nu}^s_{\hat{q}}$, which is the conditional measure
of $\hat{\nu}$ along $W^s[\hat{q}]$. Then, $\cL^s[\hat{q}]$ is
$\hat{T}^t$-equivariant. Since 
the action of $\hat{T}^{-t}$ is expanding along $W^s[\hat{q}]$, we see that for
almost all $\hat{q}$ and any $\epsilon > 0$, $\cL^s[\hat{q}]$ is the smallest
real-algebraic subset of $\cW^s(\hat{q})$ such that $L^s[\hat{q}]$ contains
$\operatorname{support} (\hat{\nu}^s_{\hat{q}}) \cap B^{Q}(\hat{q},\epsilon)$. 

We lift the measure $\hat{\nu}$ to $\hat{X}$,
and for $\hat{x} \in \hat{X}$ declare \index{$L$@$\cL^s[\hat{x}]$}$\cL^s[\hat{x}]$
to be $\cW^s[\hat{x}] \cap \sigma^{-1}(\cL^s[\sigma(\hat{x})])$, where $\index{$\sigma$}\sigma: \hat{X}
\to \hat{Q}$ is the projection map given by $\sigma(x,\omega,s) =
(\sigma(x),\omega,s)$.  

For $\index{$x$@$\hat{x}$}\hat{x} = (x,\omega,s) \in \hat{X}$, 
let $T_x\cL^s[\hat{x}] \subset T_x\cW^s[\hat{x}]$ be the tangent space at $x$ to $\cL^s[\hat{x}]$. Note that $T_x\cL^s[\hat{x}]$ is $\hat{T}^t$-invariant by construction.  


Let $U^+(\hat{x})\subset \bbG^{ssr}(\hat{x})$ be as in \autoref{def:random:compatible_family_of_subgroups}.   
Let $U^+[\hat{x}] = U^+(\hat{x})\cdot x \subset \cW^u[\hat{x}]$ denote the
$U^+(\hat{x})$-orbit of $x$.  
Given $\hat{q}=(q,\omega,s)\in \hat{Q}$, let $T_qQ =
\index{$E^s(\hat{q})$}E^s(\hat{q}) \oplus
\index{$E^c(\hat{q})$}E^c(\hat{q}) \oplus
\index{$E^u(\hat{q})$}E^u(\hat{q})$ be the splitting into negative,
zero, and positive Lyapunov exponents. We then declare for $\hat{x}
\in \hat{X}$,
\begin{displaymath}
\index{$E^s(\hat{x})$}E^s(\hat{x}) = E^s(\sigma(\hat{x})), \qquad \index{$E^c(\hat{x})$}E^c(\hat{x}) =
E^c(\sigma(\hat{x})), \qquad \index{$E^u(\hat{x})$}E^u(\hat{x}) = E^u(\sigma(\hat{x})).
\end{displaymath}
Let $\index{$\pi^{cs}_{\hat{x}}$}\pi^{cs}_{\hat{x}}\colon T_xQ \to E^u(\hat{x})$ be the canonical projection.

For $x \in X$, let \index{$Z$@$\cZ(x)$}
\begin{equation}
\label{eq:def:random:cZ}
\cZ(x) = \essspan \{ T_x \cL^s[x,\omega,0]  \st {\omega\in \Diff_\infty(Q)^\Z} \},  
\end{equation}
where the measure on $\Diff_\infty(Q)^\Z$ is $\mu^\Z$.
Note that $\cZ(x)$ is  $\mu$-invariant (as in \autoref{def:mu_invariant_subbundle})
since $\cL^s$ is $\hat{T}^t$-invariant.

The following is the random dynamics version of \autoref{lemma:AWB:to:QNI}:
\begin{lemma}
\label{lemma:random:AWB:to:QNI}
Suppose for a positive (and thus full) $\nu$-measure subset of $\hat{x}\in \hat{X}$ we have:
\begin{equation}
\label{eq:random:continue:induction}
  \pi^{cs}_{\hat{x}} \left( \cZ(x) \right)\not \subset T_{\hat{x}}U^+[\hat{x}].
\end{equation}
Then the family $\hat{x} \to	U^+(\hat{x})$ satisfies the Random QNI condition (see \autoref{def:Random:QNI}).
\end{lemma}

\begin{proof}
See \cite{ElliotSmith}.
\end{proof}

The following is the analogue of
\autoref{lemma:equivalent:stopping:conditions}: 
\begin{lemma}
\label{lemma:random:equivalent:stopping:conditions}
The following are equivalent:
\begin{itemize}
\item[{\rm (1)}] For a positive (and thus full) $\hat{\nu}$-measure subset
  of $\hat{x} = (x,\omega,s) \in \hat{X}$ we have 
\begin{equation}
\label{eq:random:cZ:subset:TU}
  \pi^{cs}_{(x,\omega,0)} \left( \cZ(x) \right) \subset T_xU^+[x,\omega,0].
\end{equation}
\item[{\rm (2)}] For a positive (and thus full) $\hat{\nu}$-measure
  subset of $\hat{x}=(x,\omega,s) \in \hat{X}$ there exists a positive (and thus full) Lebesgue measure subset $E_x \subset \Diff_\infty(Q)^\N$ such that for $\omega'_- \in E_x$ we have
	\[
		\pi^{cs}_{(x,\omega,0)}\left(T_x\cL^s[x,\omega'_-,0)]\right) \subset T_x U^+[x,\omega,0].
	\]
\end{itemize}      
\end{lemma}

\subsubsection{Entropy theory for the skew extensions}
Recall that we assume $\mu$ is finitely supported.  

We work in the random dynamics setting.  We have $\hat T^t\colon (\hat X, \hat \nu)\to (\hat X, \hat \nu)$ which is a finite measurable extension of $ \hat T^t\colon (\hat  Q,   \nu)\to (\hat  Q,  \hat  \nu)$.  Moreover, $(\hat Q, \hat \nu)$ is an extension over $(S^\Z, \mu^\Z) \times [0,1]$ with fibers diffeomorphic to $Q$.


We note $U^+[{x},\omega,s]= U^+[ {x},\omega',s]$ for all $\omega'$ with $(\omega')^+ =\omega^+$ and for $\lambda_i>0$, the quantity $\dim E^{\lambda_i}(\hat{x})\cap T_x U^+[\hat {x}]$ is $\hat T$-invariant and thus constant a.s.
We also note that $\cL^s[x,\omega,s]=  \cL^s[x,\omega',s]$ for $\omega'$ with $(\omega')^- =\omega^-$.  

We recall that if $\mu= \sum {c_i} \delta_{f_i}$ is a discrete probability measure on $\Diff^\infty(Q)$, we define the entropy of $\mu$ to be $$\index{$H(\mu)$}H(\mu) = -\sum  c_i \log c_i.$$
We also write $\hat T = \hat T^1$.  
If $\zeta$ is a measurable partition of $\hat Q$ with $\hat T\zeta \prec \zeta$ we define 
\begin{displaymath}
\index{$h_{\hat \nu}(\hat T,\zeta)$}h_{\hat \nu}(\hat T,\zeta):= H_{\hat \nu}(T\inv \zeta\mid \zeta)= - \int \log \hat \nu|_{\zeta[\hat x]}(T\inv \zeta(\hat x))\, d \hat \nu(\hat x).
\end{displaymath}
Here, $\index{$\nu\mid_{\zeta[\hat x]}$}\nu|_{\zeta[\hat x]}$ is the conditional measure at $\hat x$ relative to the measurable partition $\zeta$. 

Given $\hat x=(x,\omega,s)\in \hat Q$, we write $\index{$W_1^+[\hat x]$}W_1^+[\hat x] = \{x\}\times \{\omega' : (\omega')^+= \omega^+\}\times \{s\}$.  
Let $\zeta$ denote the partition of $\hat Q$ whose atom through $\hat x$ is $\zeta(x)=W_1^+(\hat x) $. Then we define 
\begin{displaymath}
\index{$h_{\hat \nu}(\hat T,W_1^+)$}h_{\hat \nu}(\hat T,W_1^+):= h_{\hat \nu}(\hat T,\zeta).
\end{displaymath}
We similarly define $\index{$W_1^-[\hat x]$}W_1^-[\hat x] $ and $\index{$h_{\hat \nu}(\hat T^{-1},W_1^-)$}h_{\hat \nu}(\hat T^{-1},W_1^-)$.


\begin{proposition}\label{prop:LYfacts} For $\hat \nu$-a.e.\ $\hat x= (x,\omega,s)\in \hat X$, the following hold:
\begin{enumerate}
\item \label{LYfact1} $h_{\hat \nu}(\hat T)\ge h_{\hat \nu}(\hat T,W_1^+) + \sum_{\lambda_i>0} \lambda_i \dim E^{\lambda_i}(\hat{x})\cap T_x U^+[\hat {x}]$.
\item \label{LYfact2}$h_{\hat \nu}(\hat T,W_1^+) = h_{\mu^\Z\times ds} (T) = H(\mu)$.
\item \label{LYfact3} $h_{\hat \nu}(T^{-1})\le h_{\hat \nu}(\hat T^{-1},W_1^-) + \sum_{\lambda_i<0}- \lambda_i \dim E^{\lambda_i}(\hat{x})\cap T_x \cL^s[\hat x]$.  
\item \label{LYfact4} $h_{\hat \nu}(\hat T^{-1},W_1^-) \le H(\mu)$ with equality  
if and only if $\nu$ is $\mu$-invariant (invariance principle). 


\item\label{LYfact5} We have $h_{\hat \nu}(\hat T^{-1})= h_{\hat \nu} (\hat T^{-1},\widehat \gB_0^-)$.  Moreover, if $\nu$ is $\mu$-invariant then 
$$h_{\hat \nu}(\hat T^{-1})= H(\mu) + h_{\hat \nu} (\hat T^{-1},\gB_0^-).$$

\item\label{LYfact6} We have  $h_{\hat \nu}(\hat T) = h_{\hat \nu}(\hat T,W_1^+)$ if and only if  the leafwise measure $\hat \nu^u_x$ is supported at a single point for a.e.\ $x$.

\end{enumerate}
\end{proposition}
\begin{proof}
	The first three statements are analogous to \autoref{prop:LYfactsgt}.

	Next, \autoref{LYfact4} follows in the same way as the proof of \cite[\S2.1, Prop.~2]{Ledrappier_Positivity-of-the-exponent-for-stationary-sequences-of-matrices}.
	
	The first identity of \autoref{LYfact5} follows as in \autoref{prop:full_entropy} from \cite{LedrappierYoung1985_The-metric-entropy-of-diffeomorphisms.-II.-Relations-between-entropy-exponents}.
	The second identity follows by using \autoref{LYfact4}, which implies that $\hat \nu$ has a local product structure on $\hat{\frakB}_0^-$ and hence the entropy adds up from each factor.

	To prove \autoref{LYfact6}, we note that applying the previous two parts to $\hat{T}$ instead of $\hat{T}^{-1}$ we find that $h_{\hat \nu}(\hat{T},\gB_0)=0$.
	It follows that the relative entropy vanishes $H_{\hat{\nu}}(\gB_0|\hat{T}^{-1}\gB_0)=0$ from which, by iteration, it follows that $\hat\nu^u_x$ is supported in the intersection of a nested sequence of atoms of $\hat{T}^{-k}\gB_0$ and hence is atomic.
\end{proof}

\begin{proposition}\label{prop:LedrappierFacts}
  Suppose that for $\hat \nu$-a.e. $\hat q=(q,\omega,s) \in \hat Q$
  the measure $\hat \nu^s_{\hat q}$ is supported on an embedded submanifold $N_{\hat q}\subset \cW^q[ \hat q]$.  
Then 
$$h_{\hat\nu}(\hat T)\le H(\mu) + \sum_{\lambda_i<0} (-\lambda_i) \dim E^{\lambda_i}(\hat{q})\cap T_q N_{\hat{q}}.$$
Moreover, equality holds if and only if 
$\hat\nu^s_{\hat q}$ equivalent to the Lebesgue measure on $N_{\hat q}$ for a.e.\ $\hat q\in \hat Q$.  
\end{proposition}
\begin{proof}
Let $ \xi$ be a measurable partition of $\hat Q$ with
\begin{enumerate}
\item $\xi\prec \hat T(\xi) $
\item $\xi(\hat q)\subset N_{\hat q}$
\item $\xi(\hat q)$ contains an open neighborhood of $\hat q$ in $N_{\hat q}$.
\end{enumerate}
We may take $\xi(\hat q)= N_{\hat q}\cap \gB_0^-$.  
Then $h_{\hat\nu}(\hat T) = H(\mu) + h_{\hat\nu} (\hat T, \xi)$ and so we conclude that 
$h_{\hat\nu} (\hat T, \xi)\le  \sum_{\lambda_i<0}- \lambda_i \dim E^{\lambda_i}(\hat{q})\cap T_q N_{\hat{q}}$.  It follows exactly as in the proof of \cite[Thm.\ 3.4]{Ledrappier-84} that equality holds if and only if $\hat \nu^s_{\hat q}$ is equivalent to the Lebesgue measure on $N_{\hat q}$ for a.e.\ $\hat q$.
\end{proof}

\subsubsection{Entropy and the endgame}
The following is the random dynamics version of
\autoref{prop:entropy:trick}:

\begin{proposition}
\label{prop:random:entropy:trick}
Suppose $\nu$ is a $\mu$-stationary measure on $X$ which is a lift of
a $\mu$-stationary measure on $Q$ (also denoted by $\nu$). Suppose $\hat{x}
\to U^+(\hat{x})$ is a family of subgroups compatible with the measure
$\hat{\nu}$ on $\hat{X}$.  Let $\cZ$ be as in (\ref{eq:def:random:cZ}). 
Suppose also that
\begin{equation}
\label{eq:random:cZ:non:negative:sum}
\text{the sum of the Lyapunov exponents along $\cZ$ is non-negative. }
\end{equation}
Furthermore, suppose that for $\hat{\nu}$-a.e.\ $x \in X$, 
(\ref{eq:random:cZ:subset:TU}) holds. 
Then, $\nu$ is $\mu$-invariant (see \autoref{def:statandinvmsr}).  
 Also,
\begin{itemize}
\item[{\rm (i)}] $U^+(x,\omega,s)$ depends only on $(\sigma(x),\omega,s)$, where $\sigma: X \to Q$ is the covering map of Section~\ref{ssec:the_finite_measurable_cover}. 
\item[{\rm (ii)}] For almost all $\hat{q} \in \hat{Q}$, the
  conditional of $\hat{\nu}$ along $\cW^u[\hat{q}]$ is the Haar
  measure along the orbit $U^+[\hat{q}]$. 
\item[{\rm (iii)}] For almost all $\hat{q} \in \hat{Q}$ there exists a subgroup $U^-(\hat{q})$ of
$\bbG^{ssr}(\cW^s[\hat{q}])$ such that the conditional of $\hat{\nu}$ along $\cW^s[\hat{q}]$ is the Haar measure along the orbit $U^-[\hat{q}] \equiv U^-(\hat{q})q$. In particular, $\cL^s[\hat{q}] = U^-[\hat{q}]$. 
\item[{\rm (iv)}] The sum of the Lyapunov exponents of the bundle
  $\cZ$ is $0$, where $\cZ$ is as in (\ref{eq:def:random:cZ}).  
\end{itemize}
\end{proposition}

\begin{proof}
From (\ref{eq:random:cZ:subset:TU}), for a.e.\ $\hat{x} =
(x,\omega,s)$, 
we have
$$\cZ(x) \cap E^u(\hat{x}) \subset T_x U^+[\hat{x}].$$ 
By \autoref{prop:LYfacts}\ref{LYfact1},
\begin{align*}
h_{\hat{\nu}}(\hat{T})
& \ge H(\mu) + \sum_{\lambda_i>0} \lambda_i \dim E^{\lambda_i}(\hat{x})\cap T_x U^+[\hat{x}]\\
& \ge H(\mu) + \sum_{\lambda_i>0} \lambda_i \dim E^{\lambda_i}(\hat{x})\cap \cZ(x).
\end{align*}
In the above, the term $H(\mu)=h_{\hat{\nu}}(\hat{T}, W_1^+)$ is the contribution of the ``combinatorial unstable'' $W_1^+$.

Similarly, by the definition of $\cZ(x)$
we have $T_x\cL^s[\hat{x}]\subset \cZ(x)$ whence by \autoref{prop:LYfacts}\ref{LYfact3},
\begin{align*}
h_{\hat{\nu}}(\hat{T}^{-1}) & \le h_{\hat{\nu}}(\hat{T}^{-1},W_1^-) + \sum_{\lambda_i<0} - \lambda_i \dim E^{\lambda_i}(x)\cap T_x\cL^s[x]
\\&\le h_{\hat{\nu}}(\hat{T}^{-1},W_1^-) + \sum_{\lambda_i<0} - \lambda_i \dim E^{\lambda_i}(\hat{x})\cap \cZ(x).
\end{align*}
In the above, $h_{\hat{\nu}}(\hat{T}^{-1}, W_1^-)$ denotes the contribution of the ``combinatorial stable'' $W_1^-$ to the entropy. Note  by \autoref{prop:LYfacts}\ref{LYfact4} 
that $h_{\hat{\nu}}(\hat{T}^{-1}, W_1^-) \le H(\mu)$.   

In view of (\ref{eq:random:cZ:non:negative:sum}), 
\begin{align*}
0
  &\le H(\mu) - h_{\hat{\nu}}(\hat{T}^{-1}, W_1^-) + \sum_{\lambda_i} \lambda_i \dim E^{\lambda_i}(\hat{x})\cap \cZ(x)\\
&= H(\mu) - h_{\hat{\nu}}(\hat{T}^{-1}, W_1^-) + \sum_{\lambda_i<0} \lambda_i \dim E^{\lambda_i}(\hat{x})\cap \cZ(x)
  +\sum_{\lambda_i>0} \lambda_i \dim E^{\lambda_i}(\hat{x})\cap \cZ(x)\\
&\le - h_{\hat{\nu}}(\hat{T}^{-1}) + h_{\hat{\nu}}(\hat{T})
\end{align*}
Since $h_{\hat{\nu}}(\hat{T}) =h_{\hat{\nu}}(\hat{T}^{-1})$, we conclude that $h_{\hat{\nu}}(\hat{T}^{-1}, W_1^-) = H(\mu)$, hence \autoref{prop:LYfacts}\ref{LYfact4} implies that $\nu$ is $\mu$-invariant.   

We also get that
\begin{multline}
\label{eq:random:equal:entropy}  
H(\mu) + \sum_{\lambda_i<0} - \lambda_i \dim \left(E^{\lambda_i}(\hat{x})\cap T_x\cL^s[\hat{x}]\right)  = h_{\hat{\nu}}\left(\hat{T}^{-1}\right) =h_{\hat{\nu}}\left(\hat{T}\right) = \\ =  H(\mu) + \sum_{\lambda_i>0} \lambda_i \dim \left(E^{\lambda_i}(\hat{x})\cap T_x U^+[\hat{x}]\right).
\end{multline}
By \autoref{prop:LedrappierFacts}, for a.e.\ $\hat{x} \in \hat{X}$, $\hat{\nu}^s_{\hat{x}}$ is absolutely continuous on $\cL^s[\hat{x}]$.

Analogously to the definition of $\cL^s[\hat{x}]$, let
\index{$L$@$\cL^u[\hat{x}]$}$\cL^u[\hat{x}] \subset \cW^u[\hat{x}]$ denote the smallest algebraic subset of full $\hat{\nu}^u_{\hat{x}}$ measure. Note that since the measure on $\hat{X}$ is a lift of the measure on $\hat{Q}$, $\cL^u[\hat{x}]$ depends only on $\sigma(x)$ where $\sigma: \hat{X} \to \hat{Q}$ is the projection map.

Let \index{$\mu$@$\check{\mu}$}$\check{\mu}$ be the measure on $\Diff_\infty(Q)$ given by
$\check{\mu}(g) = \mu(g^{-1})$. Since $\nu$ is $\mu$-invariant, it is
also $\check{\mu}$-invariant, and in particular
$\check{\mu}$-stationary. We have 
\begin{displaymath}
\cW^u[x,\omega,s] = \cW^s[x, \check{\omega}, s] \qquad
\cW^s[x,\omega,s] = \cW^u[x, \check{\omega}, s], 
\end{displaymath}
where
\begin{displaymath}
\index{$\omega$@$\check{\omega}$}\check{\omega} = \dots,  \omega_2^{-1}, \omega_1^{-1}, \omega_0^{-1},
\omega_{-1}^{-1}, \omega_{-2}^{-1} \dots, 
\end{displaymath}
Note that for any $E$, $\mu^\Z( \{\omega \in E \}) =
\check{\mu}^\Z(\{ \check{\omega} \in E \})$. 
Then,
\begin{displaymath}
\hat{\nu}^u_{(x,\omega,s)} = \hat{\nu}^s_{(x,\check{\omega}, s)} \qquad
\hat{\nu}^s_{(x,\omega,s)} = \hat{\nu}^u_{(x,\check{\omega}, s)}. 
\end{displaymath}
In particular, $\cL^u[x,\omega,s] = \cL^s[x,\check{\omega},s]$ for a.e.\ $\hat{x} \in \hat{X}$
and since $\hat{\nu}^s_{\hat{x}}$ is absolutely continuous on $\cL^s[\hat{x}]$, 
it follows that $\hat{\nu}^u_{\hat{x}}$ is absolutely continuous on $\cL^u[\hat{x}]$.
Note that $U^+[\hat{x}] \subset \cL^u[\hat{x}]$, and we must have
$\dim U^+[\hat{x}] = \dim \cL^u[\hat{x}]$ (or else we contradict (\ref{eq:random:equal:entropy})).
Thus, $U^+[\hat{x}] = \cL^u[\hat{x}]$ since $\cL^u[\hat{x}]$ is irreducible (see \autoref{sssec:smoothness_and_irreducibility}), and in particular, $U^+[\hat{x}] = U^+[\sigma(\hat{x})]$.

Let $U^-(x,\omega,s) = U^+(x,\check{\omega},s)$. Then, $U^-(\hat{x})$ is a subgroup of $\bbG^{ssr}(\cW^s[\hat{x}])$ and $\cL^s[\hat{x}] = U^-(\hat{x}) x \equiv U^-[x]$. 
\end{proof}

\begin{proof}[Proof of \autoref{theorem:random:stationary:is:invariant} and \autoref{theorem:random:measure:classification}.]
In view of \autoref{cor:random:end:induction} and
\autoref{lemma:random:AWB:to:QNI},  (\ref{eq:random:cZ:subset:TU})
holds. 
Note that the bundle $\cZ$ is non-random. 
Therefore, in view of the assumption
in \autoref{theorem:random:stationary:is:invariant} or
\autoref{theorem:random:measure:classification},
(\ref{eq:random:cZ:non:negative:sum}) holds. 
The rest of the assertions follow from
\autoref{prop:random:entropy:trick}. 
\end{proof}


\subsubsection{Uniform Expansion}

\begin{proposition}
  \label{prop:uniform:expansion:implies:no:negative:subbundle}
Let $\mu$ be as in
\autoref{theorem:random:stationary:is:invariant}, and suppose
suppose for $1 \le d < \dim Q$, $\mu$ satisfies 
the uniform expansion condition \autoref{def:uniformly_expanding}.
Let $\nu$ be a $\mu$-stationary measure. Then, for every 
non-random $\nu$-measurable proper subbundle of $TQ$ the sum
of the Lyapunov exponents is positive (except for the zero bundle).

If $Q$ is a symplectic manifold and $\mu$ preserves the symplectic form on $Q$ and has uniform expansion on isotropic subspaces, then there exists no non-random $\nu$-measurable proper subbundle of $TQ$ on which the sum of the Lyapunov exponents is negative. 
\end{proposition}

\begin{proof}
This follows immediately from iterating
\autoref{def:uniformly_expanding}. For the last sentence, note
that for any symplectic subbundle, the sum of the Lyapunov exponents
is $0$. 
\end{proof}

\begin{proof}[Proof of \autoref{theorem:intro:stationary:is:invariant}]
  This follows immediately from \autoref{theorem:random:stationary:is:invariant} and \autoref{prop:uniform:expansion:implies:no:negative:subbundle}.
\end{proof}

\begin{proof}[Proof of \autoref{theorem:intro:ue:measure:classification}]
We apply \autoref{theorem:random:measure:classification}.
Let $\index{$T$@$\cT(x)$}\cT(x) = T_xU^+ \oplus T_x U^-$. Then, by
\autoref{theorem:random:measure:classification}(iii), the sum of the
Lyapunov exponents along $\cT$ is $0$. But, in view of
\autoref{prop:uniform:expansion:implies:no:negative:subbundle} and the
assumption in the statement of
\autoref{theorem:intro:ue:measure:classification}, for any subbundle
$\cT'$ of $TQ$ with $1 \le \dim \cT' < \dim Q$ the sum of the Lyapunov expoents along $\cT'$ is positive. Therefore, either $\dim \cT = 0$ or $\dim \cT = \dim Q$.

If $\dim \cT = 0$ then by \cite{RandomBarreiraPesin}
(and using the assumption that there are no $0$ Lyapunov exponents), the dimension of $\nu$ is $0$.\footnote{It may be possible to avoid the assumption that there are no $0$ exponents in this part of the argument.} Then, by \autoref{prop:Margulis:function:implies:positive:dimension}, $\nu$ has an atomic part. 
Then by the ergodicity of $\hat{T}$, $\nu$ is atomic, and all atoms
have the same measure. Therefore $\nu$ is finitely supported. 

Now suppose $\dim \cT = \dim Q$, which implies that the conditional measures on the (Pesin) stables and unstables are absolutely continuous. Furthermore, we are assuming that there are no zero Lyapunov exponents.
Then, $\nu$ is absolutely continuous with respect to Lebesgue measure by \cite[Corollary H]{LedrappierYoung1985_The-metric-entropy-of-diffeomorphisms.-II.-Relations-between-entropy-exponents}.

The proof of ergodicity of the volume measure is identical to that in 
\cite[\S10]{DolgopyatKrikorian2007_On-simultaneous-linearization-of-diffeomorphisms-of-the-sphere}.
Indeed, the uniform expansion gaps property guarantees that the conclusion of \cite[Cor.~4(c)]{DolgopyatKrikorian2007_On-simultaneous-linearization-of-diffeomorphisms-of-the-sphere} holds, with $d/2$ replaced by the dimension where Lyapunov exponents switch from positive to negative (see also \cite[\S3.2]{Zhang2019_On-stable-transitivity-of-finitely-generated-groups-of-volume-preserving}).
Then, the arguments in \cite[\S10.3-6]{DolgopyatKrikorian2007_On-simultaneous-linearization-of-diffeomorphisms-of-the-sphere} hold and yield the result.
\end{proof}


\subsection{\Teichmuller Dynamics}
	\label{ssec:teichmuller_dynamics_proofs}

We now proceed to the proofs of the results stated in \autoref{ssec:teichmuller_dynamics_statements}.
The notation follows loc. cit. and we do not repeat it here.

\subsubsection{Some lemmas on Lie groups}

Let $\index{$G$}G=\SL(2,\R)$, $\index{$G^m$}G^m = G \cross G \cross \dots G$ ($m$ times). Let $\index{$G_\Delta$}G_\Delta \subset G^m$ be the diagonally embedded copy of $G$ in $G^m$.

Let $\beta$ be a partition of $\{1, \dots, m\}$. Let $\index{$\norm$@$\sim_\beta$}\sim_\beta$ denote the associated equivalence relation, i.e.\ $k \sim_\beta l$ iff $k$ and $l$ are in the same atom of $\beta$. Let 
\begin{displaymath}
\index{$G_\beta$}G_\beta = \{ (g_1, \dots, g_m) \in G^m \st g_j = g_k \text{ if } j \sim_\beta k \}  
\end{displaymath}
Note that $G_\beta \cap G_\gamma = G_\delta$ where $\delta$ is the smallest partition such that both $\beta$ and $\gamma$ are refinements of $\delta$. Also, if $\beta$ is the partition of $\{1,\dots,m\}$ which has only one atom, then $G_\beta = G_\Delta$.

We will also need the following generalization. Suppose $S$ is a
(possibly empty) subset of $\{1,\dots, m\}$ and suppose $\beta$ is a
partition of $S$. We set\index{$G$@$\tilde{G}_\beta$}
\begin{multline*}
  \tilde{G}_\beta = \{ (g_1, \dots, g_m) \in G^m \st g_j = e \text{ if } j \not\in S, \\ \qquad g_j = g_k \text{ if  $j \in S$, $k \in S$, $j \sim_\beta k$}\}
\end{multline*}
Let \index{$Z^m$}$Z^m$ denote the center of $G^m$.

We recall the following standard lemma. 
\begin{lemma}
\label{lemma:normalized:by}
Suppose $F$ is a subgroup of $G^m$ which is normalized by $G_\Delta$. Then there exists a subset $S$ of $\{1,\dots,m\}$, a partition $\beta$ of $S$ and a subgroup $Z'$ of $Z^m$ such that 
$F= \tilde{G}_\beta Z'$.
\end{lemma}

As an immediate corollary we have the following:
\begin{corollary}
\label{cor:subgroups:containing:diagonal}
Suppose $F$ is a connected subgroup of $G^m$ which contains $G_\Delta$. Then, $F=G_\beta$ for some partition $\beta$ of $\{1,\dots, m\}$.   
\end{corollary}

  


\begin{proof}[Proof of \autoref{lemma:normalized:by}]
	Let $H:=\operatorname{PSL}(2,\bR):=G/Z$ be the quotient by the center.
	It suffices to prove the analogous statement for $H^m,H_\Delta$ and with $H_{\beta},\wtilde{H}_\beta$ as before.
	We recall the Goursat lemma: for any groups $H_1, H_2$, subgroups $S\subset H_1\times H_2$ that project surjectively to each factor are classified by normal subgroups $N_i\subset H_i$ an an isomorphism $f\colon H_1/N_1\toisom H_2/N_2$.
	The group $S$ is the preimage in $H_1\times H_2$ of the graph of $f$ in $H_1/N_1\times H_2/N_2$.

	We proceed by induction on $m$, the statement for $m=1$ being immediate.
	Let $F_{m-1}\subset H^{m-1}$ be the projection of $F$ to the first $m-1$ factors.
	It satisfies the inductive assumptions so $F_{m-1}=\wtilde{H}_{\beta'}$ with $\beta'$ a partition of a subset of $\{1,\dots,m-1\}$.
	We apply the Goursat lemma to $F\subset F_{m-1}\times H$, let $N_1\subset F_{m-1}, N_2\subset H$ be the corresponding groups.
	Since $H$ is simple, either $N_2=\id$ or $N_2=H$.

	If $N_2=H$ then $F\isom F_{\beta'}\times H$, clearly of the needed form.
	If $N_2=\id$ then $F_{m-1}/N_1\isom H$, and the normal subgroups of $F_{m-1}$ with this property are of the form $H_{\beta''}$, with $\beta''$ obtained from $\beta'$ by removing a full equivalence class $C$ from $\beta'$, with $C\subseteq \{1,\dots,m-1\}$.
	Note that the only automorphisms of $\operatorname{PSL}(2,\bR)$ are inner, or given by $g\mapsto \left(g^{-1}\right)^t$.
	Now the constraint that $F$ is normalized by $H_{\Delta}$ implies that $F\isom F_{\beta}$, where $\beta$ is obtained from $\beta'$ by adding the equivalence class $C\cup \{m\}$.
\end{proof}

\subsubsection{Stables, unstables and projections.}

Suppose $q = (q_1, \dots, q_m)$. Let \index{$p_j$}$p_j$ denote the natural map from $H^1(M_j, \Sigma_j, \R)$ to $H^1(M_j,\R)$, and let \index{$p$}$p: H^1(M,\Sigma,\R) \to H^1(M,\R)$ denote $\bigoplus_{j=1}^m p_j$.  
Let \index{$H^\perp(q_j)$}
\begin{multline*}
  H^\perp(q_j) = \{ v \in H^1(M_j, \Sigma_j, \R) \st 
  \text{$p(v_j)$ is symplectically orthogonal to} \\ \text{the plane spanned by $p_j(\Re q_j)$ and $p_j(\Im q_j)$.} \} 
\end{multline*}
Let \index{$H^\perp(q)$}
\begin{displaymath}
H^\perp(q) = \bigoplus_{j=1}^m  H^\perp(q_j)
\end{displaymath}
Let $\index{$H^+(q_j)$}H^+(q_j) = \R \Im q_j$, $\index{$H^-(q_j)$}H^-(q_j) = \R \Re q_j$. Let
\begin{displaymath}
\index{$H^+(q)$}H^+(q) = \bigoplus_{j=1}^m H^+(q_j) \qquad \index{$H^-(q)$}H^-(q) = \bigoplus_{j=1}^m H^-(q_j)
\end{displaymath}
Then 
if we identify $\cH(\vec{\alpha})$ with $H^1(M,\Sigma;\mathbb{C})$
using the period map, then
\begin{displaymath}
\cW^u[q] = \{ q + v \st v \in H^1(M,\Sigma,\R), \qquad p_j(v) \wedge p_j(\Im
q) = 0, \forall j\}
\end{displaymath}
The condition on the wedge product, or equivalently the symplectic pairing, is needed so that $q+v \in
\cH_1(\vec{\alpha})$ (and not $\cH(\vec{\alpha})$). Thus, 
\begin{equation}
\label{eq:Teich:Wu:form}
  \cW^u[q] = q+ H^+(q) + H^\perp(q) \equiv q + \cW^u(q)
\end{equation}
\begin{equation}
\label{eq:Teich:Ws:form}
\cW^s[q] = q + i H^-(q) +i H^\perp(q) \equiv q + \cW^s(q)
\end{equation}
Then $\cW^s[q]$ is naturally identified with $T_q\cW^s[q]$, and similarly,
$\cW^s[q]$ is naturally identified with $T_q\cW^s[q]$. 

Let \index{$\iota_{\tau}$}$\iota_{\tau}\colon W^s(q) \to W^u(q)$ denote the linear isomorphism given by for $c_j \in \R$, $v_j \in H^\perp_j$
\begin{multline}
\label{eq:def:iota}
\iota_{\tau}\big( i (c_1 \Re q_1 + v_1), \dots, i(c_m \Re q_m + v_m)\big) = \\ = \left(c_1 \tau^2\Im q_1 + \tau v_1, \dots, c_m \tau^2 \Im q_m + \tau v_m\right)
\end{multline}
Let \index{$\pi$@$\bar{\pi}^{cs}_q[q']$}$\bar{\pi}^{cs}_q[q']$ denote the unique intersection point of $\cW^{cs}[q']$ with $\cW^u[q]$.

Let $\index{$n_\tau$}n_\tau = \begin{pmatrix} 1 & \tau \\ 0 & 1 \end{pmatrix} \subset N
\subset SL(2,\reals)$.
Let $f_{q,\tau}: \cW^s[q] \to \cW^u[q]$ be defined as\index{$f_{q,\tau}(q')$}
\begin{displaymath}
f_{q,\tau}(q') = n_\tau^{-1} \bar{\pi}^{cs}_q(n_\tau q'). 
\end{displaymath}
Let \index{$f_{q,\tau,t}(q')$}$f_{q,\tau,t}: \cW^s[q] \to \cW^u[q]$ denote $g_{-t} \circ f_{g_t q,\tau} \circ g_t$.
\begin{lemma}
\label{lemma:explicit:projection}
Suppose $q \in \cH_1(\vec{\alpha})$, $q' \in \cW^s[q]$. 
Write $q=(q_1, \dots, q_n)$, $q' = (q_1', \dots, q_n')$. Then, in view of (\ref{eq:Teich:Ws:form}),
\begin{displaymath}
q'_j = q_j + i( c_j \Re q_j + v_j)
\end{displaymath}
for some $c_j \in \R$ and $v_j \in H^\perp_j$. Then, we have
\begin{itemize}
\item[{\rm (a)}] $f_{q,\tau}(q') = \big(\phi_{q_1,\tau}(c_1,v_1), \dots, \phi_{q_m,\tau}(c_m,v_m)\big)$, where 
\begin{displaymath}
\phi_{q,\tau}(c,v) = q-\frac{c\tau^2}{1+c\tau}\Im q + \frac{\tau}{1+c\tau} v.
\end{displaymath}
\item[{\rm (b)}] $f_{q,\tau,t}(q') = \left(\phi_{q_1,\tau,t}(c_1,v_1), \dots, \phi_{q_m,\tau,t}(c_m,v_m)\right)$, where 
\begin{displaymath}
\phi_{q,\tau,t}(c,v) = q- \frac{e^{-2t}c\tau^2}{1+e^{-2t}c\tau} \Im q +\frac{e^{-2t}\tau}{1+e^{-2t}c\tau} v.
\end{displaymath}
\item[{\rm (c)}] The derivative of $f_{q,\tau}(q')$ with respect to $q'$
  at $q'=q$ is $\iota_{\tau}$, where $\iota_{\tau}$ is as in (\ref{eq:def:iota}).
\end{itemize}
\end{lemma}

\begin{proof}

We begin with the proof of (a). It is enough to work one factor at a time.  
It is also convenient to express the coordinates as a $2\times \dim H^1(M,\Sigma)$ matrix:
\[
	q_j = \begin{bmatrix}
		1 & 0 & 0\\
		0 & 1 & 0
	\end{bmatrix}
	\quad
	q_j' = 
	\begin{bmatrix}
		1 & 0 & 0\\
		c & 1 & v
	\end{bmatrix}
	\text{ w.r.t. }
	\begin{bmatrix}
		\bR \Re q_j & \bR \Im q_j & H^{\perp}(q_j)\\
		i\bR \Re q_j & i\bR \Im q_j & iH^{\perp}(q_j)
	\end{bmatrix}
\]
Then $\SL_2(\bR)$ acts by ordinary matrix multiplication, whereas center-stable projection acts as:
\begin{align*}
	\pi_q^{cs}\left(
	\begin{bmatrix}
	a & b & v\\
	c & d & w					
	\end{bmatrix}
	\right)
	&
	=
	\begin{bmatrix}
	1 & \tfrac{b}{a} & \tfrac{v}{a}\\
	0 & 1 & 0
	\end{bmatrix}\\
	\pi_{g_t q}^{cs}\left(
	\begin{bmatrix}
	a & b & v\\
	c & d & w					
	\end{bmatrix}
	\right)
	&
	=
	\begin{bmatrix}
	e^{t} & e^{t}\tfrac{b}{a} & e^t\tfrac{v}{a}\\
	0 & e^{-t} & 0
	\end{bmatrix}
	\text{ since }
	g_tq = \begin{bmatrix}
		e^t & 0 & 0\\
		0 & e^{-t} & 0
	\end{bmatrix}.
\end{align*}
We now compute directly:
\begin{align*}
	n_{\tau}q_j' & = 
	\begin{bmatrix}
		1 + \tau c & \tau & \tau v\\
		c & 1 & v
	\end{bmatrix}\\
	\pi_{q}^{cs}n_{\tau}q_j' & = 
	\begin{bmatrix}
		1 & \tfrac{\tau}{1 + \tau c} & \tfrac{\tau}{1 + \tau c}v\\
		0 & 1 & 0
	\end{bmatrix}
	\\
	n_{\tau}^{-1}\pi_{q}^{cs}n_{\tau}q_j' & = 
	\begin{bmatrix}
		1 & \tfrac{-\tau^2 c}{1 + \tau c} & \tfrac{\tau}{1 + \tau c}v\\
		0 & 1 & 0
	\end{bmatrix}
\end{align*}
For the second part, the calculation is analogous:
\begin{align*}
	g_t q_j' & = 
	\begin{bmatrix}
		e^t  & 0  & 0\\
		e^{-t} c & e^{-t} & e^{-t}v
	\end{bmatrix}\\
	n_\tau g_t q_j' & = 
	\begin{bmatrix}
		e^t + \tau e^{-t}c  & \tau e^{-t}  & \tau e^{-t}v\\
		e^{-t} c & e^{-t} & e^{-t}v
	\end{bmatrix}\\
	\pi_{g_tq}^{cs} n_\tau g_t q_j' & = 
	\begin{bmatrix}
		e^t  & \frac{\tau}{e^t + \tau e^{-t}c}  & \frac{\tau v}{e^t + \tau e^{-t}c}\\
		0 & e^{-t} & 0
	\end{bmatrix}\\
	n_{-\tau}\pi_{q}^{cs} n_\tau g_t q_j' & = 
	\begin{bmatrix}
		e^t  & \frac{-\tau^2 e^{-2t}c}{e^t + \tau e^{-t}c}  & \frac{\tau v}{e^t + \tau e^{-t}c}\\
		0 & e^{-t} & 0
	\end{bmatrix}\\
	g_{-t}n_{-\tau}\pi_{q}^{cs} n_\tau g_t q_j' & = 
	\begin{bmatrix}
		1  & \frac{-\tau^2 e^{-t}c}{e^t + \tau e^{-t}c}  & \frac{\tau e^{-t} v}{e^t + \tau e^{-t}c}\\
		0 & 1 & 0
	\end{bmatrix}\\
	& = 
	\begin{bmatrix}
		1  & \frac{-\tau^2 e^{-2t}c}{1 + \tau e^{-2t}c}  & \frac{\tau e^{-2t} v}{1 + \tau e^{-2t}c}\\
		0 & 1 & 0
	\end{bmatrix}
\end{align*}
which precisely matches the claimed formulas.

Part (c) follows immediately from part (a) by taking the derivative at $q_j' = q_j$ (which corresponds to $c_j =0$, $v_j=0$).
\end{proof}

\begin{lemma}[No negative subbundles]
\label{lemma:no:negative:subbundles}
Suppose $V \subset H^1(M,\Sigma,\R)$ is a subbundle which is
equivariant under the action of the Kontsevich-Zorich cocycle. Then,
the sum of the Lyapunov exponents of $V$ is non-negative. Also if $V$
is an equivariant subbundle of $H^1(M,\Sigma;\mathbb{C})$, then the
sum of the Lyapunov exponents of $g_t$ on $V$ is non-negative.    
\end{lemma}

\begin{proof}
For a subbundle of a single factor $H^1(M_j,\Sigma_j)$ the result follows from \cite[Thm.~A.3]{EskinMirzakhani_Invariant-and-stationary-measures-for-the-rm-SL2Bbb-R-action-on-moduli-space}.
Given an arbitrary $V$, proceed by induction.
Indeed, its projection to the last factor satisfies the claim, and the kernel is contained in the preceding $n-1$ factors.
\end{proof}


\begin{lemma}
	\label{lemma:homogeneous:support}
For a.e.~$q \in Q$, $\nu_q^s$ gives measure $0$ to any proper algebraic subset of $\cL^s[q]$.
\end{lemma}
A closely related result (concerning only the foliation by fast stable manifolds) is established in \cite{LedrappierXie2011_Vanishing-transverse-entropy-in-smooth-ergodic-theory}; we follow the same method.
\begin{proof}
	Suppose, by contradiction, that there exists a minimal $d<\dim \cL^{s}$ such that for a.e. $q$ there exists at least one algebraic irreducible $\cD^s\subset \cW^s[q]$ with $\nu^s_q(\cD^s)>0$ and $\dim \cD^s=d$.
	Note that for any two such $\cD^s$, their intersection has $\nu^s_q$-measure zero by minimality of $d$.
	There are at most countably many such $\{\cD^s_i[q]\}_{i\in I(q)}$, and their union has full $\nu^s_q$-measure.
	Our goal is to show that $\nu^s_q$ is supported on a single one, implying that $\cD^s[q]=\cL^s[q]$.

	First, note that for a.e. $q$ there is a single $\cD^s_{i}$ denoted $\cD^s[q]$ such that $q\in \cD^s[q]$.
	Indeed, since $\nu_s^q$ is supported on the union of the $\cD_i^s[q]$ (by ergodicity), if there was a positive measure set inside the intersection of two (or more) $\cD_i^s[q]$ then we could consider algebraic sets of lower dimension.

	Consider now the partition $\frakD$ with atoms $\frakD[q]:=\cD^s[q]\cap \gB_0^-[q]$.
	It is a measurable partition, since for each pair $(k,m)$ we can construct a standard measure space (fibering over $Q$) whose fiber over $q$ consists of $m$ polynomials of degree at most $k$ on $\cW^s[q]$, whose vanishing locus intersects $\gB_0^-[q]$.
	On this space, the set of irreducible subvarieties of dimension $d$ with positive $\nu_q^s$-mass is measurable, and we can take the union over all $(k,m)\in \bN^2$.

	Note also that $\frakD$ is an \emph{invariant} partition, i.e. $g_{t}\left(\frakD[q]\right)\subseteq \frakD\left[g_t q\right]$ for $t\geq 0$.
	Furthermore, we have for $t\geq 0$ that $g_t(\frakD)\vee \frakD=g_t(\gB_0^-)\vee \frakD$, which will be useful when computing entropy.

	Next, we recall the \emph{mutual information function} of two measurable partitions $\frakA,\frakB$:
	\[
		I_{\nu}(q;\frakA|\frakB):= -\log \Big(\nu_{\frakB}[q]\big([q]_{\frakA}\big)\Big)
	\]
	and relative entropy:
	\[
		H_{\nu}(\frakA|\frakB):=\int_Q I_{\nu}(q;\frakA|\frakB)d\nu(q)
	\]
	Recall from \cite[\S7.1]{Rokhlin1967_Lectures-on-the-entropy-theory-of-transformations-with-invariant-measure} that for $g_t$-invariant measurable partition $\frakA$ we have
	\[
		h_{\nu}\left(g_t;\frakA\right):=H_{\nu}\left(\frakA|g_{t}^{-1}\frakA\right)
		\text{ and }h_{\nu}(g_t) :=\sup_{\frakA} h_{\nu}(g_t;\frakA)
	\]
	where the $\sup$ is over all $g_t$-invariant $\frakA$ (see \cite[\S9]{Rokhlin1967_Lectures-on-the-entropy-theory-of-transformations-with-invariant-measure}).
	Note that since $g_t$ is invertible we have $H_{\nu}(\frakA\vert g_t^{-1}\frakA)= H_{\nu}(g_t \frakA\vert \frakA)$.

	We claim first that, because $\frakD$ is finer and countably many of its atoms form an atom of $\gB_0^-$, we have (for fixed $t_0>0$):
	\begin{align}
		\label{eqn:entropy_inequality_hard_direction}
		h_{\nu}\left(g_{t_0};\frakD\right)\geq h_{\nu}\left(g_{t_0};\gB_0^-\right)
	\end{align}
	Indeed, for any partition $\frakA$ with $g_{t_0}(\frakA)$ finer than $\frakA$, we have
	\begin{align*}
		\tfrac{1}{n}I_{\nu}\left(q;g_{t_0}^n\frakA\vert\frakA\right) & = 
		\tfrac 1n
		\left[
		I_{\nu}\left(q; g_{t_0}^{n}\frakA \vert g_{t_0}^{n-1}\frakA  \right) 
		+ \cdots +
		I_{\nu}\left(q; g_{t_0}\frakA \vert \frakA  \right)
		\right]\\
		& = \tfrac{1}{n}\left[
		I_{\nu}\left(g_{t_0}^{-(n-1)}q; g_{t_0}\frakA \vert \frakA  \right)
		+
		\cdots
		+
		I_{\nu}\left(q; g_{t_0}\frakA \vert \frakA  \right)
		\right]\\
		& \xrightarrow{n\to +\infty} \int_Q I_{\nu}\left(q; g_{t_0}\frakA | \frakA \right) d\nu(q)
		= h_{\nu}(g_{t_0};\frakA)
	\end{align*}
	where in the last line we used the Birkhoff ergodic theorem.
	We next note that
	\[
		I_{\nu}\left(q;\frakD\vert\gB_0^-\right) = 
		- \log \left(\nu_{q}^s\left([q]_{\frakD}\right)\right) \in [0,+\infty)
	\]
	is a.e. finite.
	We next compute, using that $g_{t_0}^{n}\frakD = g_{t_0}^n\gB_0^-\vee \frakD$:
	\begin{align*}
		I_{\nu}\left(q;g_{t_0}^{n}\frakD \vert \frakD \right) & = 
		I_{\nu}\left(q;g_{t_0}^{n}\gB_0^- \vee \frakD\vert \frakD \right) =
		- \log \left(
		\frac{\nu_{q}^s([q]_{g_{t_0}^n\gB_0^-} \cap [q]_{\frakD} ) }
		{\nu_q^s\left([q]_\frakD\right)}
		\right)\\
		& = - 
		\log 
		\left({\nu_{q}^s([q]_{g_{t_0}^n\gB_0^-} \cap [q]_{\frakD})}\right)
		+ \log
		\left({\nu_q^s\left([q]_\frakD\right)} \right)\\
		& \geq 
		- 
		\log 
		\left({\nu_{q}^s([q]_{g_{t_0}^n\gB_0^-}}\right)
		+ \log
		\left({\nu_q^s\left([q]_\frakD\right)} \right)\\
		& = 
		I_{\nu}\left(q;g_{t_0}^n \gB_0^-\vert \gB_0^-\right)
		+ \log
		\left({\nu_q^s\left([q]_\frakD\right)} \right)
	\end{align*}
	Dividing by $n$ and sending it to $+\infty$ yields \autoref{eqn:entropy_inequality_hard_direction}.
	
	We now compute the difference.
	With $t_0>0$ fixed as above, let $B:=\gB_0[q]$ be an on atom of $\gB_0^-$, and let $\{B_j\}_{j\in J(q,t_0)}$ be the atoms of $g_{t_0}(\gB_0^-)$ contained in $B$, so we have $B=\coprod B_j$.
	Recall that $\nu_q^s$ is the conditional (probability) measure of $\nu$ on $B$.
	We then set
	\[
		m_j:=\nu_q^s(B_j) \quad d_{i,j}:=\nu_q^s\left(B_j \cap \cD_i^s[q]\right)
		\quad\text{so}\quad
		\sum_{j}m_j = 1 \quad \sum_{i}d_{i,j}=m_j.
	\]
	We also set $d_i:=\sum_{j}d_{i,j}=\nu_q^s(\cD_i^s[q]\cap \gB_0^-[q])$.
	The contribution of $B$ to $h_{\nu}(g_{t_0};\frakB)$ is (with the standard convention $0\log 0 :=0$):
	\[
		\sum_{j} - \log( m_j)\cdot m_j
	\]
	and to $h_{\nu}(g_{t_0};\frakD)$ it is
	\begin{multline*}
		\sum_{i} d_i \cdot \left(\sum_{j} -\log\left(\tfrac{d_{i,j}}{d_i}\right)\cdot \tfrac{d_{i,j}}{d_i}\right)=\\
		=\sum_{j} \sum_{i} -\log\left(\tfrac{d_{i,j}}{d_i}\right)\cdot \tfrac{d_{i,j}}{d_i}\cdot d_i
		\leq 
		\sum_{j} -\log\left(m_j\right)m_j
	\end{multline*}
	where in the last line we applied Jensen to the concave function $-x\log x$, with values $d_{i,j}/d_i$ and weights $d_i$.
	Integrating over $Q$ yields $h_{\nu}\left(g_t;\frakD\right)\leq h_{\nu}(g_t;\gB_0^-)$ so combined with \autoref{eqn:entropy_inequality_hard_direction}, it follows that the entropies agree, and furthermore we have a.e. equality in Jensen.
	From the equality case of Jensen we conclude that $d_{i,j}/d_i=m_j$, for all $i,j$, for a.e. $q$.

	Therefore, the conditional measures $\nu_{\frakD,q}$ and $\nu^s_q$, with respect to the partitions $\frakD,\gB_0^-$, give the same mass to each atom $B_j$ of $\gB_0^-\vee g_{t_0}(\gB_0^-)$.
	We now repeat the argument with $t_0\to +\infty$ and note that $\bigvee_{t\geq 0}g_t(\gB_0^-)$ is the partition of $Q$ into points.
	Therefore, for a.e. $q$ the conditional measures for $\frakD,\gB_0^-$ agree and hence there is a single $\cD[q]=\cL^s[q]$.
\end{proof}

\begin{lemma}
\label{lemma:homogeneous:dynamics}
Fix $\eta > 0$ arbitrary.
Fix also a linear identification of $W^s(q)$ with $\bR^m$ (measurably in $q$) such that $T_q \cL^s[q]$ is given by the vanishing of the last $m - \dim \cL^s$ coordinates.
There exists a function $c_{\eta,m,d}(q)$ with $c_{\eta,m,d}(q) > 0$ for a.e.\ $q$, such
that the following holds:  
Suppose $h$ is a polynomial of degree $d$ on $\gB_0^-[q]$ with
\begin{displaymath}
	|D^I h(q)| = 1, 
\end{displaymath}
where $D^I$ is a partial derivative of degree $m$, and $I$ is a multi-index involving only derivatives in the first $\dim \cL^s$ coordinates.

Then,
\begin{displaymath}
\nu^s_q\left(\{ q' \in \gB_0^-[q] \st \abs{h\left(q'\right)} < c_{\eta,m,d}(q) \}\right) \le \eta
\nu^s_q(\gB_0^-[q]). 
\end{displaymath}
\end{lemma}

\begin{proof}
Let \index{$P$@$\cP_q$}$\cP_q$ denote the space of polynomials $h$ of
degree $d$ on $\cL^s[q]$, consisting of all polynomials on $\cW^s[q]$ modulo those that vanish on $\cL^s[q]$.
Define next\index{$P$@$\cP_q^\ast$}
\[
 	\cP_q^\ast :=\{ h \in \cP_q \st \sup_{q' \in \gB_0^-[q]\cap \cL^s[q]} |h(q')|= 1\}\subset \cP_q
 \]
which is a compact set. 
Therefore, for a.e.\ $q$
there exists a constant $C =
C(q,d,m) > 0$ such that for $h \in \cP_q^\ast$ and all the (finitely many) multi-indices $I$ as in the statement, we have:
\begin{equation}
  \label{eq:tmp:DI:hq}
|D^I h(q)| < C.   
\end{equation}
Suppose the assertion of the lemma is not true.
Then, there exists a sequence $h_n \in \cP_q$  and $\delta_n \to 0$
such that $|D^I h_n(q)| = 1$ and 
\begin{displaymath}
\nu^s_q(\{ q' \in \gB_0^-[q] \st |h_n(q')| < \delta_n \}) \ge \eta
\nu^s_q(\gB_0^-[q]). 
\end{displaymath}
Let $s_n = \sup_{q' \in \gB_0[q]\cap \cL^s[q]} |h_n(q')|$, and 
let $\bar{h}_n = h_n/s_n$. Then, $\bar{h}_n \in \cP_q^\ast$.
We have, by (\ref{eq:tmp:DI:hq}),
\begin{displaymath}
\frac{1}{s_n} = \abs{D^I \bar{h}_n(q)} \le C,
\end{displaymath}
Then, if $|h_n(q')| < \delta_n$ then,
$|\bar{h}_n(q')| < C \delta_n$. Thus, 
\begin{displaymath}
\nu^s_q\left(\left\lbrace q' \in \gB_0^-[q] \st |\bar{h}_n(q')| < C \delta_n \right\rbrace\right) \ge \eta
\nu^s_q(\gB_0^-[q]). 
\end{displaymath}
Since $\cP_q^\ast$ is compact, the polynomials $\bar{h}_n$ converge to a
polynomial $\bar{h}_\infty$. Then,
\begin{displaymath}
\nu^s_q(\{ q' \in \gB_0^-[q] \st |\bar{h}_\infty(q')| =0 \}) \ge \eta
\nu^s_q(\gB_0^-[q]). 
\end{displaymath}
Since $\bar{h}_\infty$ is not the zero polynomial, 
this contradicts \autoref{lemma:homogeneous:support}.
\end{proof}

\subsubsection{The QNI condition on $\cH_1(\vec{\alpha})$.}
Let $X$ be the measurable cover of $\cH_1(\vec{\alpha})$ constructed
in \S\ref{ssec:the_finite_measurable_cover}. 

\begin{proposition}
\label{prop:Teich:to:QNI}
Suppose the QNI condition does not hold. Then,
\begin{multline}
\label{eq:pics:q:nqprime}
\bar{\pi}_x^{cs}(n x') \subset U^+[q], \\
\text{for a.e.\ $x \in X$ for all
  $x' \in \cL^s[x]$, and all $n \in N$} 
\end{multline}
\end{proposition}

\begin{proof}
	Define
	\begin{displaymath}
	 \index{$\psi_x(x',\tau)$}\psi_x(x',\tau) = d\left(n_\tau x', \bar{\pi}^{cs}_{x'}\left(U^+[x]\right)\right)^2,
	 \qquad\text{ where } n_\tau = \begin{pmatrix} 1 & \tau \\ 0 & 1 \end{pmatrix}.
	\end{displaymath}
	We break up the proof into two parts.

	\noindent Claim 1: There exists $\alpha_0>0$ such that for any $\delta>0$, there exists $c(\delta)>0$ such that for a $(1-\delta)$-fraction of $x'\in \frakB_\ell^-[x]$, there exists some $\tau\in \left[-e^{-\ell/2},e^{-\ell/2}\right] $ such that $\psi_{x}(x',\tau)>c(\delta)e^{-\alpha_0\ell}$.

	\noindent Claim 2:
	Suppose that for $x'_{1/2}\in \frakB_{\ell/2}^-[x_{1/2}]$ there exists $u'x'_{1/2}\in \cB_{\ell/2}[x'_{1/2}]$ with
	\[
		d^u\left(u'x'_{1/2}, \bar{\pi}^{cs}_{x'}\left(U^+[x]\right)\right)\geq \ve.
	\]
	Then for any $\delta>0$ there exists $a(\delta)>0$ such that, for a $(1-\delta)$-fraction of $ux_{1/2}\in \cB_{\ell/x}[x_{1/2}]$ we have
	\[
		d^u\left(U^+[x'_{1/2}], \bar{\pi}^{cs}_{x'}\left(ux_{1/2}\right)\right)\geq a(\delta )\ve.
	\]
	Note that $N\subset U^+(q)$ so Claim 1 implies the assumption of Claim 2 holds for a $(1-\delta)$ fraction of $x'_{1/2}\in \frakB_{\ell/2}^-[x_{1/2}]$.
	In turn, Claim 2 implies that the QNI condition as in \autoref{def:QNI} holds and completes the proof.
	Claim 2 however is a direct consequence of \cite[Lemma~6.18]{EskinMirzakhani_Invariant-and-stationary-measures-for-the-rm-SL2Bbb-R-action-on-moduli-space}, so it remains to show Claim 1.

	To proceed, note that for each $x$ the function $\psi_x$ is real-analytic in $x'$ and $\tau$.
	It also suffices to prove that there exists $k\in \bN$, with the same quantifies on $x',c(\delta)$ we have that
	\[
		\abs{
		\frac{d^k}{d\tau^k}\psi_{x}\left(x',\tau\right)\Big\vert_{\tau=0}
		}
		> c(\delta)e^{-\alpha_1\ell}
	\]
	and from this derive the claim, using the same argument as in \autoref{lemma:homogeneous:dynamics}.

	Suppose (\ref{eq:pics:q:nqprime}) does not hold. Then,
	$\psi_x$ is not identically $0$ for $x$ in a set of positive
	measure. Then, by equivariance $\psi_x$ is not identically $0$ for
	almost all $x$.
	Then, there exists $m$ and $k$ such that for almost
	all $x$, there exists a multi-index $I$ of length $m$ such that
	\begin{displaymath}
	\frac{\partial^I}{\partial^I x'} \left(\frac{d}{d\tau}\right)^k \psi_x(x',\tau) \Big\vert_{\stackrel{x'=x}{\tau=0}}\ne 0.
	\end{displaymath}
	Then for any $\delta > 0$ there exists a compact set $K$ with $\nu(K) > 1-\delta$ and a constant $c(\delta) > 0$ such that for $x_{1/2} \in K$,
	\begin{displaymath}
	\left|\frac{\partial^I}{\partial^I x'_{1/2}} \left(\frac{d}{d\tau_{1/2}}\right)^k \psi_{x_{1/2}}(x'_{1/2},\tau_{1/2}) \Big\vert_{\stackrel{x'_{1/2}=x_{1/2}}{\tau_{1/2}=0}}\right|\ge c(\delta).
	\end{displaymath}
	Write $x_{1/2} = g_{\ell/2} x$,
	$x'_{1/2} = g_{\ell/2} x'$ where $x' \in \gB^-_0[x]$. Let $\tau_{1/2} = e^{-\ell} \tau$, so that $n_{\tau_{1/2}} = g^{-\ell/2} n_\tau g^{\ell/2}$.
	Then, 
	\begin{displaymath}
	\left|\frac{\partial^I}{\partial^I x'} \left(\frac{d}{d\tau}\right)^k \psi_{x_{1/2}}\left(g_{\ell/2} x',e^{-\ell} \tau\right) \Big\vert_{\stackrel{x'=x}{\tau=0}}\right|\ge c(\delta) e^{-\alpha \ell}.
	\end{displaymath}
	where $\alpha > 0$ depends only on the Lyapunov spectrum.
	Choose $M > 2 \alpha$.
	Then, there exists an integer $d$ depending
	only on the Lyapunov spectrum, $m$, and $k$, and for every $\delta > 0$ a compact set
	$K_1$ with $\nu(K_1) > 1-\delta$ 
	such that for $x_{1/2} \in K_1$ there exists a polynomial
	$P_{x_{1/2}}(\cdot)$ of degree $d$ given by a Taylor expansion, such that for $x'_{1/2} \in \gB_{\ell/2}^-[x_{1/2}]$ we have:
	\begin{displaymath}
	\left|
	\frac{\partial^{I'}}{\partial^{I'} x'_{1/2}}
	\left[
	\left(\frac{d}{d\tau_{1/2}}\right)^k
	\psi_{1/2}(x'_{1/2},\tau_{1/2})\big\vert_{\tau_{1/2} = 0} 
	- 
	P_{x_{1/2}}(x'_{1/2})
	\right]
	\right|
	\le C(\delta) e^{-M \ell}, 
	\end{displaymath}
	for all multi-indices ${I'}$ with $|{I'}|\leq m$.
	Therefore,
	\begin{align}
		\label{eqn:P_derivative_lower_bound}
		\left|\frac{\partial^I}{\partial^I x'} P_{x_{1/2}}(g_{\ell/2} x') \Big\vert_{\stackrel{x'=x}{\tau=0}}\right|\ge c(\delta) e^{-\alpha \ell}.
	\end{align}
	Then, by \autoref{lemma:homogeneous:dynamics} applied to (a
	normalized version of) the polynomial $P_{x_{1/2}}(g_{\ell/2} x')$ for $x' \in \gB_0^-[x]$, 
	for at least $(1-\delta)$-fraction of $x'_{1/2} \in \gB_{\ell/2}^-[x_{1/2}]$,
	\begin{displaymath}
		\abs{P_{x_{1/2}}\left(x'_{1/2}\right)}\ge c(\delta) e^{-\alpha \ell}.
	\end{displaymath}
	Since $M > 2\alpha$, by \autoref{eqn:P_derivative_lower_bound}
	we get least $(1-\delta)$-fraction of $x'_{1/2} \in \gB_{\ell/2}^-[x_{1/2}]$,
	\begin{align}
		\label{eqn:lower_bound_psi_derivative}
			\left| 
		\left(\frac{d}{d\tau}\right)^k
		\psi_{x_{1/2}}\left(x'_{1/2}, e^{-\ell} \tau\right)
		\vert_{\tau=0}
		\right|
		\geq c(\delta)e^{-\alpha \ell}
	\end{align}
	and this completes the proof.
\end{proof}

\subsubsection{Proving linearity of the stable and unstable supports}
\begin{proposition}
	\label{prop:EM:theorem2.1}
Suppose (\ref{eq:pics:q:nqprime}) holds. Then $\nu$ is $\SL(2,\R)$
invariant.
In addition, there exists an $\SL(2,\R)$-equivariant map
$L$ from $\cH_1(\vec{\alpha})$ to $\mathbb{C}$-linear subspaces of
$H^1(M,\Sigma;\mathbb{C})$, such that $L$ is the complexification of an $\bR$-linear subspace $L_{\bR}\subset H^1\left(M,\Sigma;\bR\right)$, and 
such that for a.e.\ $q \in \cH^1(\vec{\alpha})$, we have that $q \in L(q)$, and the conditional measure $\nu^u_q$ of $\nu$ along $\cW^u[q]$ is supported on $\cW^u[q] \cap L(q)$ and of Lebesgue class on it, and the conditional measure $\nu^s_q$ of $\nu$ along $\cW^s[q]$ is supported on $\cW^s[q] \cap L(q)$ and of Lebesgue class on it.
\end{proposition}

\begin{proof}
Since $n x \in \cW^u[x]$, we have, by the definition of $\bar{\pi}^{cs}_x$, $\bar{\pi}^{cs}_{nx} = \bar{\pi}^{cs}_x$.   Suppose (\ref{eq:pics:q:nqprime}) holds. Then
\begin{multline*}
\bar{\pi}^{cs}_{nx}(n x') \subset U^+[n x], \\
\text{for a.e.\ $x \in X$ for all
  $x' \in \cL^s[x]$, for all $n \in N$} 
\end{multline*}
Then, taking the derivative with respect to $x'$ at $x'=x$, we get that for $x$ in a set of full measure and all $n \in N$, 
\begin{displaymath}
\pi^{cs}_{nx}(n \cdot T_{x} \cL^s[x])) \subset T_{nx} U^+[nx].
\end{displaymath}
Then, by \autoref{lemma:equivalent:stopping:conditions}, for a.e.\ $x \in X$,
(\ref{eq:cZ:subset:TU}) holds.
Let $\cZ(x)$ be as in (\ref{eq:def:cZ}), and let $\cZ'(x)$ be as
(\ref{eq:def:cZprime}). Then, by \autoref{lemma:no:negative:subbundles},
(\ref{eq:cZprime:non:negative:sum}) holds. Thus, we may apply
\autoref{prop:entropy:trick}. In particular, we get that
$U^+[x]$ depends only on $\sigma(x)$, 
$\nu$ is $\SL(2,\R)$-invariant, the stable and unstable supports are homogeneous, 
\begin{displaymath}
\cZ(q) = \cZ'(q) = T_q U^+[q]\oplus T_q\cL^s[q],
\end{displaymath}
and the sum of the Lyapunov exponents along $\cZ'$ is $0$.

Note that since $U^+$ is $N$-invariant, (\ref{eq:pics:q:nqprime}) implies that
for a.e.\ $q \in \cH_1(\vec{\alpha})$ for all $q' \in \cL^s[q]$, for all $\tau \in \R$,
\begin{equation}
\label{eq:conjugate:pics:q:nqprime}
n_\tau^{-1}\bar{\pi}^{cs}_q(n_\tau q') \in U^+[q], \\
\end{equation}
Thus, for all $q' \in \cL^s[q]$, and for all $\tau \in \R$,
\begin{displaymath}
f_{q,\tau}(q') \in U^+[q]. 
\end{displaymath}
Differentiating this at $q' = q$ and using \autoref{lemma:explicit:projection} (c), we get
\begin{displaymath}
\iota_{\tau}(T_q\cL^s[q]) \subset T_q U^+[q]. 
\end{displaymath}
Alternatively, 
\begin{equation}
\label{eq:iota:inverse}
  T_q\cL^s[q] \subset \iota_{\tau}^{-1}(T_q U^+[q]). 
\end{equation}

Also, since both $U^+$ and $\cL^s$ are $g_t$-equivariant, we have
for all $q' \in \cL^s[q]$, and for all $\tau \in \R$ and all $t \in \R$
\begin{equation}
\label{eq:tmp:f:q:tau:t}
  f_{q,\tau,t}(q') \in U^+[q].
\end{equation}
(Initially (\ref{eq:tmp:f:q:tau:t}) holds for $t$ small, so we can think of the Kontsevich-Zorich cocycle as constant. But then, by analyticity, it holds for all $t$).

\begin{claim}
We have $\forall \tau\in \bR$:
\begin{equation}
\label{eq:iota:TqcL:TqU}
\iota_{\tau}(T_q\cL^s[q]) = T_q U^+[q]. 
\end{equation}
In particular,   
\begin{displaymath}
\dim \cL^s[q] = \dim U^+[q]. 
\end{displaymath}
\end{claim}

\begin{proof}[Proof of claim]
We have the Lyapunov decomposition of the Kontsevich-Zorich cocycle:
\begin{displaymath}
H^1(M, \Sigma, \R) = \bigoplus_{k} E^{\lambda_k}.  
\end{displaymath}
Then,
\begin{displaymath}
T_q U^+[q] = \bigoplus_{k} E^{\lambda_k}(q) \cap T_q U^+[q], 
\end{displaymath}
and
\begin{displaymath}
T_q \cL^s[q] = \bigoplus_{k} i E^{\lambda_k}(q) \cap T_q \cL^s[q]. 
\end{displaymath}
Note that
\begin{equation}
\label{eq:iota:Elambda}
\iota_{\tau}(i E^{\lambda_k}(q)) = E^{\lambda_k}(q)
\end{equation}
where on the
right-hand side we have identified $H^1(M,\Sigma,\R)$ with $\{ x \in
H^1(M,\Sigma;\mathbb{C}) \st \Im x = 0 \}$. Under this identification,
$T_q U^+[q]$ is a $P$-invariant subbundle of $H^1(M,\Sigma,\R)$, and
thus, by \autoref{lemma:no:negative:subbundles}, 
\begin{equation}
\label{eq:sumk:lambdak}
  \sum_k \lambda_k \dim E^{\lambda_k}(q) \cap T_q U^+[q] \ge 0. 
\end{equation}
The Lyapunov exponent of $g_1$ on $E^{\lambda_k}$ is $1+\lambda_k$,
and the Lyapunov exponent of $g_1$ on $i E^{\lambda_k}$ is
$-1+\lambda_k$, and $0 \le \lambda_k \le 1$.  Then, using
(\ref{eq:iota:Elambda}), (\ref{eq:iota:inverse}) and (\ref{eq:sumk:lambdak}),
\begin{multline}
\label{eq:long:teich}
  \sum_k (1-\lambda_k) \dim (i E^{\lambda_k}(q) \cap T_q \cL^s[q]) \le \\
\le   \sum_{k}(1-\lambda_k) \dim (\iota_{\tau}^{-1} (E^{\lambda_k}(q)) \cap \iota_{\tau}^{-1}(T_q U^+[q])) \\ = 
  \sum_{k}(1-\lambda_k) \dim (E^{\lambda_k}(q) \cap T_q U^+[q])  = \\
  \sum_{k}(1+\lambda_k) \dim (E^{\lambda_k}(q) \cap T_q U^+[q]) - 2 \sum_k \lambda_k \dim E^{\lambda_k}(q) \cap T_q U^+[q] \le \\
  \le \sum_{k}(1+\lambda_k) \dim (E^{\lambda_k}(q) \cap T_q U^+[q]).
\end{multline}
But, since the sum of the Lyapunov exponents along $\cZ$ is $0$, we have
\begin{displaymath}
\sum_k(1-\lambda_k) \dim (i E^{\lambda_k}(q) \cap T_q \cL^s[q]) = \sum_k (1+\lambda_k) \dim(E^{\lambda_k}(q) \cap T_q U^+[q]).
\end{displaymath}
Thus, all the inequalities in (\ref{eq:long:teich}) are equalities, and thus
\begin{displaymath}
T_q \cL^s[q] = \iota_{\tau}^{-1}(T_q U^+[q]). 
\end{displaymath}
This concludes the proof of the claim.
\end{proof}



Let $\index{$H^{\cL}(q)$}H^{\cL}(q) = H^-(q) \cap T_q \cL^s[q]$. Then, $H^{\cL}(q)$ is the tangent space to the orbit through $q$ of a subgroup $\index{$N$@$\bar{N}^{\cL}(q)$}\bar{N}^{\cL}(q) \subset \bar{N}^m$.
By ergodicity of the measure, $\bar{N}^{\cL}(q)$ is independent of $q$.
Since the conditional measure of $\nu$ along $\cW^s[q]$ is of Lebesgue class along $\cL^s[q]$, and at $\nu^s$-a.e. point $\bar{N}^{\cL}$ is tangent to $\cL^s[q]$, it follows that $\cL^s[q]$ is $\bar{N}^{\cL}$-invariant.
It is now a standard argument that the conditional measures $\nu^s$, and hence $\nu$, are $\bar{N}^{\cL}$-invariant.


Let
\begin{equation}
\label{eq:HU:HL}
\index{$H^{U}(q)$}H^{U}(q) = H^+(q) \cap T_q U^+[q] = \iota_{\tau}(H^{\cL}(q)).  
\end{equation}
Then, $H^{U}(q)$ is the tangent space to the orbit through $q$ of a subgroup \index{$N^{U}(q)$}$N^{U}(q) \subset N^m$.

The same argument shows that $\nu$ is $N^{U}$-invariant. Also we have $N^U = \iota(N^{\cL})$, where $\index{$\iota$}\iota: \bar{N}^m \to N^m$
is defined by
\begin{displaymath}
  \iota\left(\begin{pmatrix} 1 & 0 \\ s_1 & 1 \end{pmatrix}, \dots,
    \begin{pmatrix} 1 & 0 \\ s_m & 1 \end{pmatrix}\right) =
  \left(\begin{pmatrix} 1 & s_1 \\ 0 & 1 \end{pmatrix}, \dots,
    \begin{pmatrix} 1 & s_m \\ 0 & 1 \end{pmatrix}\right)
\end{displaymath}
Then, by \autoref{cor:subgroups:containing:diagonal},
$\bar{N}^{\cL}$ and $N^U$ generate $\SL(2,\R)_\beta$ where $\beta$ is the smallest partition such that $\bar{N}^{\cL} \subset \SL(2,\R)_\beta$ and $N^U \subset \SL(2,\R)_\beta$. 
In particular, $\cL^s$ is $\bar{N}_\beta$-invariant, where
$\bar{N}_\beta = N^m \cap \SL(2,\R)_\beta$. 

\begin{claim}
\label{lemma:Nbeta:product}
For a.e.\ $q$, $\cL^s[q] = \bar{N}_\beta(\cL^s[q] \cap H^\perp[q])$.   
\end{claim}

\begin{proof}[Proof of claim]
We have, for all $q' \in \cL^s[q]$, 
$f_{q,\tau,t}(q') \in U^+[q]$.
From \autoref{lemma:explicit:projection}(b), we have, as $t \to \infty$,
\begin{displaymath}
f_{q,\tau,t}(q') = q + e^{-2t}\iota_\tau(q'-q) + O(e^{-4t})
\end{displaymath}
Thus, for all $q' \in \cL^s[q]$, $\iota_{\tau}(q'-q) \in T_qU^+[q]$.
In particular, if we write
\begin{displaymath}
q' = q+i w + i v
\end{displaymath}
with $w \in H^-(q)$, $v \in H^\perp(q)$, we have $\iota_{\tau}(w) \in
H^U(q)$, which implies, in view of (\ref{eq:HU:HL}), that $w \in H^\cL(q)$. Thus, there exists $\bar{n} \in \bar{N}^L$ such that $\bar{n}q = q+iw$. Then, $\bar{n}^{-1} q' \in \cL^s[q] \cap H^\perp[q]$.
Since $\bar{N}^L \subset \bar{N}_\beta$, the claim follows. 
\end{proof}

\begin{claim}
\label{lemma:L:star:shaped}
For a.e.\ $q$, $\cL^s[q]$ is star-shaped with respect to $q$. In other words, for any $q' \in \cL^s[q]$, for any $t \in \R$,
\begin{equation}
\label{eq:star:shaped}
q + t(q'-q) \in \cL^s[q].   
\end{equation}
\end{claim}

\begin{proof}[Proof of claim]
  In view of the $\bar{N}_\beta$-invariance of $\cL^s[q]$ and \autoref{lemma:Nbeta:product}, it is enough to prove (\ref{eq:star:shaped}) for $q' \in \cL^s[q] \cap i H^\perp[q]$. But then, \autoref{lemma:explicit:projection} (a) implies that
\begin{equation}
\label{eq:fq:tau:qprime}
f_{q,\tau}(q') = q+ \iota_{\tau}(q'-q) \in U^+[q] \cap H^\perp[q].
\end{equation}
Since for $\tau \ne 0$, $f_{q,\tau}$ is injective, and we have $\dim \cL^s[q] = \dim U^+[q]$, we get, for all $\tau \ne 0$,
$f_{q,\tau}(\cL^s[q]) = U^+[q]$.
Then, $(f_{q,\tau'}^{-1} \circ f_{q,\tau})(\cL^s[q]) = \cL^s[q]$.
Also, $(f_{q,\tau'}^{-1} \circ f_{q,\tau})(\cL^s[q] \cap i H^\perp[q]) = \cL^s[q] \cap iH^\perp[q]$. But, it follows from (\ref{eq:fq:tau:qprime}) that for $q' \in \cL^s[q] \cap i H^\perp[q])$, 
\begin{displaymath}
(f_{q,\tau'}^{-1} \circ f_{q,\tau})(q') = q + (\tau-\tau')(q'-q). 
\end{displaymath}
Since $\tau'$ is arbitrary, the claim follows.
\end{proof}

We now continue the proof of \autoref{prop:EM:theorem2.1}.
We have for $x \in \cL^s[q]$, $\cL^s[x] = \cL^s[q]$. Applying \autoref{lemma:L:star:shaped} to $x$ in place of $q$, we get that 
$\cL^s[q]$ is star-shaped with respect to any $x \in \cL^s[q]$.

Suppose $x$, $y \in \cL^s[q]$ are arbitrary.
Then, since $\cL^s[q]$ is star-shaped with respect to $x$, for all $t
\in \R$, $x + t(y-x) \in \cL^s[q]$. Then, $\cL^s[q]$ contains the line
through $x$ and $y$. This implies that $\cL^s[q]$ is an affine
subspace of $\cW^s[q]$. The same argument shows that $U^+[q]$ is an
affine subspace of $\cW^u[q]$, and then (\ref{eq:iota:TqcL:TqU})
shows that $\iota(\cL^s[q]) = U^+[q]$. Let
\begin{displaymath}
L(q) = \GL(2,\R)_\beta \left( (U^+[q]-q) + i (\cL^s[q]-q)\right). 
\end{displaymath}
Then $L(q)$ is a complex-linear subspace of $H^1(M,\Sigma;\mathbb{C})$, which is moreover the complexification of an $\bR$-linear subspace due to \autoref{eq:iota:TqcL:TqU}.
The proposition now follows since
$\nu$ is $\SL(2,\R)_\beta$-invariant, and the conditional measures of $\nu$
along $U^+[q]$ and $\cL^s[q]$ are Lebesgue.
\end{proof}

\subsubsection{Completion of the proof}

\begin{theorem}(cf. \cite[Theorem~2.1]{EskinMirzakhani_Invariant-and-stationary-measures-for-the-rm-SL2Bbb-R-action-on-moduli-space})
\label{theorem:EM:theorem2.1}
Suppose $\nu$ is an ergodic $P$-invariant measure on
$\cH_1(\vec{\alpha})$. Then, $\nu$ is $\SL(2,\R)$
invariant. In addition, there exists an $\SL(2,\R)$-equivariant map
$L$ from $\cH_1(\vec{\alpha})$ to $\mathbb{C}$-linear subspaces of
$H^1(M,\Sigma;\mathbb{C})$ 
such that for a.e.\ $q \in \cH^1(\vec{\alpha})$, $q \in L(q)$, the
conditional measure $\nu^u_q$ of $\nu$ along $\cW^u[q]$ is the
restriction of Lebegue measure to 
$\cW^u[q] \cap L(q)$ and the conditional measure $\nu^s_q$ of $\nu$
along $\cW^s[q]$ is the restriction of Lebesgue measure to
$\cW^s[q] \cap L(q)$. 
\end{theorem}
\begin{proof}
	The result follows from \autoref{prop:EM:theorem2.1}, since we can apply \autoref{thm:inductive:step} until the QNI condition no longer holds, in which case we apply \autoref{prop:Teich:to:QNI} and conclude that \autoref{eq:pics:q:nqprime} holds.
\end{proof}

The derivation of \autoref{theorem:Teich:P:measures} and
\autoref{theorem:Teich:closure:submanifold} are exactly the same
as the corresponding arguments in \cite[Sections 16-19]{EskinMirzakhani_Invariant-and-stationary-measures-for-the-rm-SL2Bbb-R-action-on-moduli-space}
and \cite{EskinMirzakhaniMohammadi_Isolation}.



\subsection{Outline of the proof of
  \autoref{thm:inductive:step}. }
  \label{ssec:outline_proof_inductive_step}

We begin by describing a set of general preliminary tools, and proceed to an outline of the proof in \autoref{sssec:choosing_the_eight_points_intro_outline}.

\subsubsection{Normal forms, Linearization, and the Space of Cosets}
	\label{sssec:normal_forms_linearization_and_the_space_of_cosets}
The unstable manifolds $\cW^u[q]$ admit normal form coordinates and a subresonant structure, see \autoref{appendix:subresonant_linear_algebra} and \autoref{appendix:normal_forms_on_manifolds_and_cocycles}.
In particular, there exist finite-dimensional Lie groups $\bbG^{sr}(q)\supset \bbG^{ssr}(q)$ acting on $\cW^u[q]$ by subresonant, resp. strictly subresonant maps.
Furthermore, we can ``linearize'' the unstable dynamics via a linear cocycle $L\cW^u(q)$ and equivariant embeddings \index{$L_q$}$L_q\colon \cW^u[q]\into L\cW^u(q)$.

Next, given a compatible family of subgroups $U^+(q)\subset \bbG^{ssr}(q)$ as in \autoref{def:compatible_family_of_subgroups}, we consider the space \index{$C$@$\cC(q)$}$\cC(q)$ of $\bbG^{ssr}(q)$-translates of the $U^+(q)$-orbit $U^+[q]:=U^+(q)\cdot q$.
To study its dynamics, we linearize it again via a linear cocycle \index{$L\cC$}$L\cC$ and equivariant embedding $\cC(q)\into L\cC(q)$ (see \autoref{sec:c_and_lc_and_boldH}).

Because the compatible family of subgroups $U^+$ is only defined on the finite cover $X$ and not on $Q$ (see \autoref{sssec:the_measurable_cover}) many of our constructions are carried out on $X$; see \autoref{ssec:the_finite_measurable_cover} for the relevant constructions.

We recall also the measurable partition $\gB_0[q]$ subordinated to the unstables, and we will also work with \index{$B$@$\cB_0[q]$}$\cB_0[q]:=U^+[q]\cap \gB_0[q]$ and the corresponding \index{$B$@$\cB_0(q)$}$\cB_0(q)\subset U^+(q)$.

\subsubsection{Factorization}
	\label{sssec:factorization_proof_outline}
Continuing with a compatible family of subgroups $U^+(q)\subset \bbG^{ssr}(q)$, we will construct in \autoref{sec:factorization} and \autoref{sec:the_curly_a_operators}:
\begin{itemize}
	\item A cocycle \index{$V^s$}$V^s$ which is smooth along stable manifolds, and has stable holonomy, together with a measurable $g_t$-equivariant section \index{$\iota_s$}$\iota_s\colon Q\to V^s$.
	The section yields a measurable $g_t$-equivariant embedding\index{$F^s_q$}
	\[
		F^s_q\colon \cW^s[q]\to V^s(q)
	\]
	by evaluating the section at $q'\in \cW^s[q]$ and moving the value to $V^s(q)$ by the holonomy.
	\item A $4$-variable map\index{$A$@$\cA_2(q_1,u,\ell,t)$}
	\begin{align}
		\label{eqn:curlyA_intro}
		V^s(q) \xrightarrow{\cA_2(q_1,u,\ell,t)} L\cC(g_t u q_1)
	\end{align}	
	where \index{$L\cC$}$L\cC$ is the linearization cocycle from \autoref{sssec:normal_forms_linearization_and_the_space_of_cosets} (see \autoref{sssec:construction_of_the_associated_bundles} for its construction and \autoref{prop:relation_to_hausdorff_distance} for its key property).
	The map is obtained as a composition
	\[
		V^s(q)\xrightarrow{\cA_2(q_1,u,\ell)}L\cC(q_1)\xrightarrow{u^*}L\cC(uq_1)
		\xrightarrow{g_t}L\cC(g_tuq_1)
	\]
	\item Let $U'\subset \cW^u[g_t u q_1]$ be the curve represented by $\cA_2(q_1,u,\ell,t)(F^s_q(q'))$.
	The Hausdorff distance between $U^+[g_t q_1']$ and $U'$ (measured at $g_t uq_1$) is exponentially small in $\ell$.
	\item It follows from the properties of $L\cC$ that the Hausdorff distance between $U^+[g_t u q_1]$ and $U'$ measured at $g_t u q_1$ is, up to a (slowly varying in the basepoint) multiplicative constant, equal to $\norm{\cA_2(q_1,u,\ell,t)(F^s(q'))}$.
	\item 
	The norm also satisfies bi-Lipschitz estimates, see \autoref{sec:new:bilipshitz}. 
	
	\item 
	An interpolation map, which is subresonant:\index{$\phi$@$\wtilde{\phi}$}
		\[
			\wtilde{\phi}\colon \cW^u[y_{1/2}] \to \cW^u[z_{1/2}]
			\text{ and satisfies }
			\wtilde{\phi}(y_{1/2})=z_{1/2}.
		\]
\end{itemize}
We note that what is needed for subsequent arguments are a number of properties of the above constructions (particularly $\cA_2$).
To construct the above objects however, we will make use of the ``half-way points'' \index{$y_{1/2}$}$y_{1/2}=g_{-\ell/2} u q_1$ and \index{$z_{1/2}$}$z_{1/2}=\cW^{cs}[y_{1/2}]\cap \cW^u_{loc}[g_{-\ell/2}q_1']$ (for two bottom-linked $Y$-configurations $Y,Y'$, see \autoref{sssec:y_configurations_outline_intro}).

\subsubsection{Inert subspaces}
	\label{sssec:inert_subspaces_intro_outline}
To study the divergence of two $U^+$-translates, we note that the cocycle $L\cC$ contains a tautological $1$-dimensional subspace \index{$L^{\tau}\cC$}$L^{\tau}\cC$ spanned by the point corresponding to $U^+(q)$.
We work with the cocycle \index{$H$@$\bbH$}$\bbH:=L\cC/L^\tau\cC$, and descend $\cA_2$ from \autoref{eqn:curlyA_intro} to take values $\bbH$.
We will also work with the natural map \index{$j$@$\bbj$}$\bbj\colon \cC\to \bbH$, obtain as the composition $\cC\into L\cC\onto \bbH$.

In \autoref{ssec:inert_subspaces} we introduce the notion of \emph{inert} subbundles \index{$E$@$\bbE^{\lambda_i}$}$\bbE^{\lambda_i}\subset \bbH$ and \index{$E$@$\bbE$}$\bbE:=\oplus \bbE^{\lambda_i}$.
Then (using Zimmer amenable reduction) we can refine the Lyapunov decomposition of $\bbE$ as \index{$E$@$\bbE^{\lambda_i}_j$}$\bbE^{\lambda_i}=\oplus \bbE^{\lambda_i}_j$ and introduce an indexing set \index{$\Lambda''$}$ij\in \Lambda''$.
Furthermore, we obtain scalar multiplicative cocycles \index{$\lambda_{ij}(q;t)$}$\lambda_{ij}(q;t)$ measuring growth in $\bbE^{\lambda_i}_j$.

The index set $\Lambda''$ carries a ``synchronization'' equivalence relation (see \autoref{prop:subbundles_e_ij_bdd}) with equivalence class denoted \index{$i$@$[ij]$}$[ij]\in \wtilde{\Lambda}$ and associated subcocycles \index{$E$@$\bbE_{[ij],bdd}(q)$}$\bbE_{[ij],bdd}$.
The properties of these decompositions are contained in \autoref{ssec:bounded_subspaces_synchronized_exponents_and_the_equivalence_relation}.

Associated to the cocycle $\bbE_{[ij],bdd}$ we obtain subsets \index{$E$@$\cE_{ij}(x)$}$\cE_{ij}(q)\subset \cC(q)$ defined as $\bbj^{-1}\left(\bbE_{[ij],bdd}\right)$ (see \autoref{sssec:the_subset_cE_ij_x}).
This leads to measurable equivalence relations
\index{$\norm$@$\sim_{ij}$}$\sim_{ij}$ (see
\autoref{sssec:equivalence_relations_ij}) with atoms
\index{$\Psi_{ij}[q]$}$\Psi_{ij}[q]\subset \gB_0[q]$ (see
\autoref{sssec:the_equivalence_classes_psi__ij}) and associated
leafwise measures \index{$f$@$\wtd f_{ij}[q]$}$\wtd f_{ij}[q]$ (coming from conditional measures \index{$f_{ij}(q)$}$f_{ij}(q)$).

It is these leafwise measure that will be ultimately shown to have extra invariance.

\subsubsection{Choosing the eight points}
	\label{sssec:choosing_the_eight_points_intro_outline}
Fix a constant $\alpha\gg 1$ and only take $t\leq \alpha \ell$.

A simplified scheme for choosing the eight points, analogously to
\cite{Eskin:Lindenstrauss:short}, is as follows: (see Figure~\ref{fig:outline})
\begin{itemize}
	\item[{\rm (i)}]
	Choose ${q}_1$ in some good set, such that in particular for
	most $t$ we have $g_t {q}_1, g_{-t}q_1 \in K_*$, 
	and for most $u$ and most $t$, we have $g_t u {q}_1 \in K_*$.

	\item[{\rm (ii)}]
	Let $\cA_2({q}_1, u, \ell, t)$ be as in \autoref{eqn:curlyA_intro} (see \autoref{sssec:the_4_variable_curly_a_2} for more details), and define ${q} = g_{-\ell} {q}_1$ and ${q}_2 = g_t u {q}_1$, where \index{$\tau_{(\epsilon)}(q_1,u,\ell)$}$t=\tau_{(\epsilon)}(q_1,u,\ell)$ is the solution to the equation $\|\cA_2({q}_1,u,\ell,t)\| = \epsilon$
	(see \autoref{sssec:defining_tau_epsilon} for $\tau_{(\epsilon)}$).


	By \autoref{lemma:monotonicity}, for fixed ${q}_1$, $u$, $\epsilon$, the function $\tau_{(\epsilon)}({q}_1,u,\ell)$ is bilipshitz in $\ell$, and by \autoref{prop:bilip:hattau:epsilon} for	most choices of $\ell$, we have ${q},q_2 \in K_*$.

	\item[{\rm (iii)}]
	For all $ij \in \tilde{\Lambda}$ (see \autoref{sssec:inert_subspaces_intro_outline}) let \index{$t_{ij}$}$t_{ij}= t_{ij}({q}_1,u,\ell)$ be defined by the equation $\lambda_{ij}(u {q}_1,t) = \lambda_{ij}({q}_1,t_{ij})$.
	Since $\lambda_{ij}(x,t)$ is bilipshitz in $t$, the same argument shows that for most choices of $\ell$, we have ${q}_{3,ij} :=	g_{t_{ij}} {q}_1 \in K_*$. 

	\item[{\rm (iv)}]
	Let \index{$M_u$}$M_u \subset V^s_{*}$ (see \autoref{sec:subsec:norm:Vs:star} for $V^s_{*}$) be the subspace of \autoref{lemma:bad:subspace} for the linear map $\cA_2({q}_1, u, \ell, t)$ restricted to $V^s_*({q})$. 
	By \autoref{prop:can:avoid:most:Mu}, we can choose ${q}' \in K_*$ with $q'\in \cW^s_{loc}[q]$, such that $\|F_{q}^s(q')\| \approx 1$ and 
	such that $F_{q}^s(q')$ avoids most of the subspaces $M_u$ as $u$ varies over $\cB_0(q)$ (recall that $\cB_0(q)\subset U^+(q)$ comes from intersecting with the Markov partition $\gB_0$).
	Then, for most $u$ we have:
	\begin{multline*}
		d^\cH_{loc}(U^+[{q}_2], \phi^{-1}_\tau(U^+[{q}_2')]) \approx \| \cA_2({q}_1,u,\ell,t)
		F_{q}^s(q') \| \approx\\
		\approx  \| \cA_2({q}_1,u,\ell,t)\| \|F_{q}^s(q') \|  \approx \epsilon,
	\end{multline*}
	as required, where $\phi_{\tau}$ is the evolved interpolation map, which satisfies estimates spelled out in \autoref{ssec:evolved_interp} and in particular does not move points too much.

	\item[{\rm (v)}]
	By an argument based on \autoref{prop:some:fraction:bounded},
	for most choices of $u$ the point ${q}_2'$ is close to $\Psi_{ij}[q_2]$ for some $ij \in \tilde{\Lambda}$ (see \autoref{sssec:inert_subspaces_intro_outline} for notation).

	\item[{\rm (vi)}]
	We can now proceed and let ${q}_1' = g_\ell {q}'$,	${q}_2' = g_t u {q}_1'$ where $t =	\tau_{(\epsilon)}({q}_1, u, \ell)$, and let ${q}_{3,ij}' =
	g_{t_{ij}} {q}_1'$.
	Recall that ${\nu}$ is $g_{t}$-invariant and $U^+$-invariant.
	Since $\lambda_{ij}({q}_1,t_{ij}) = \lambda_{ij}(u{q}_1,t)$, and one can show that $\lambda_{ij}(u q_1', t) \approx	\lambda_{ij}(q_1',t_{ij})$, we have
	\begin{align}
		\label{eqn:AA'_intro_outline}
		\begin{split}
		\tilde{f}_{ij}[{q}_2] & \propto A_* \tilde{f}_{ij}[{q}_{3,ij}]\\
		\tilde{f}_{ij}[{q}_2'] & \propto A'_* \tilde{f}_{ij}[{q}_{3,ij}']
		\end{split}
	\end{align}
	where $A,A'$ are almost the same (bounded) strictly subresonant map.
	Here and below, \index{$\norm$@$\propto$}$\propto$ means that the measures are proportional, whereas \index{$\norm$@$\approx$}$\approx$ will mean that the measures are close in the appropriately metrized space of leafwise measures up to scale (see \autoref{def:d:star}).
	
	Since ${q}_{3,ij}$ and ${q}_{3,ij}'$ are very close, we can ensure that 
	$\tilde{f}_{ij}[{q}_{3,ij}'] \approx \tilde{f}_{ij}[{q}_{3,ij}]$.
	One might expect to find $\tilde{f}_{ij}[q_2] \approx \tilde{f}_{ij}[q_2']$ but a more careful analysis (see \autoref{prop:nearby:linear:maps}) shows that the maps $A$ and $A'$ of \autoref{eqn:AA'_intro_outline} are not exactly the same.
 We get 
	\begin{displaymath}
		 \tilde{f}_{ij}[{q}_2'] \approx \psi(q_2,q_2')_*\tilde{f}_{ij}[{q}_2].  
	\end{displaymath}
	where $\psi(q_2,q_2')$ is defined in \autoref{eq:def:nontilde:psi}.
	Applying the argument with a sequence of $\ell$'s going to infinity,
	and passing to a limit along a subsequence, 
	we obtain points ${\tilde q_2}$ and ${\tilde q_2'}\in \cW^u_{loc}[\tilde q_2]$ 
	where $\tilde q_2'$ is not in $U^+[\tilde q_2]$.  
	We also have that (in some appropriate sense) 
	$\psi(q_2,q_2')$ converges to $\wp^{+}(\tilde q_2,\tilde q_2')$ 
	 satisfying
	\[
		\tilde{f}_{ij}[\tilde q_2'] \propto \wp^{+}(\tilde q_2,\tilde q_2')_* \tilde{f}_{ij}[\tilde q_2]
	\]
	Here $\wp^+(\tilde q_2,\tilde q_2')$ is a strictly subresonant map of $\cW^u[\tilde q_2]$ arising from the measurable connection $P^+$ on $L\cW^u$ (see \autoref{sssec:standard_measurable_connection} and \autoref{ssec:measurable_connections_and_subresonant_maps}).
	
	A standard argument \autoref{ssec:the_implication_have_friend_to_have_invariance}  then shows there exist a closed subgroup $U^+_{new}(\tilde q_2)$ compatible with the measure and with $U^+(\tilde q_2)\subsetneq U^+_{new}(\tilde q_2)$ such that $\tilde{f}_{ij}[\tilde q_2]$ is $U^+_{new}(\tilde q_2)$-invariant.  
	
\end{itemize}

\subsubsection{Y-configurations}
	\label{sssec:y_configurations_outline_intro}
The formal proof uses the same ideas, but we need to take a bit more care, mostly because we also need to make sure that ${q}_2'$ and ${q}_{3,ij}'$ belong to $K_*$.
We now give a slightly more precise outline of the strategy.

\begin{figure}[htbp!]
	\centering
	\includegraphics[width=0.49\linewidth]{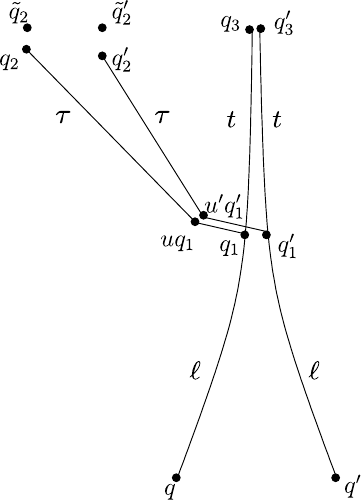}
	\includegraphics[width=0.49\linewidth]{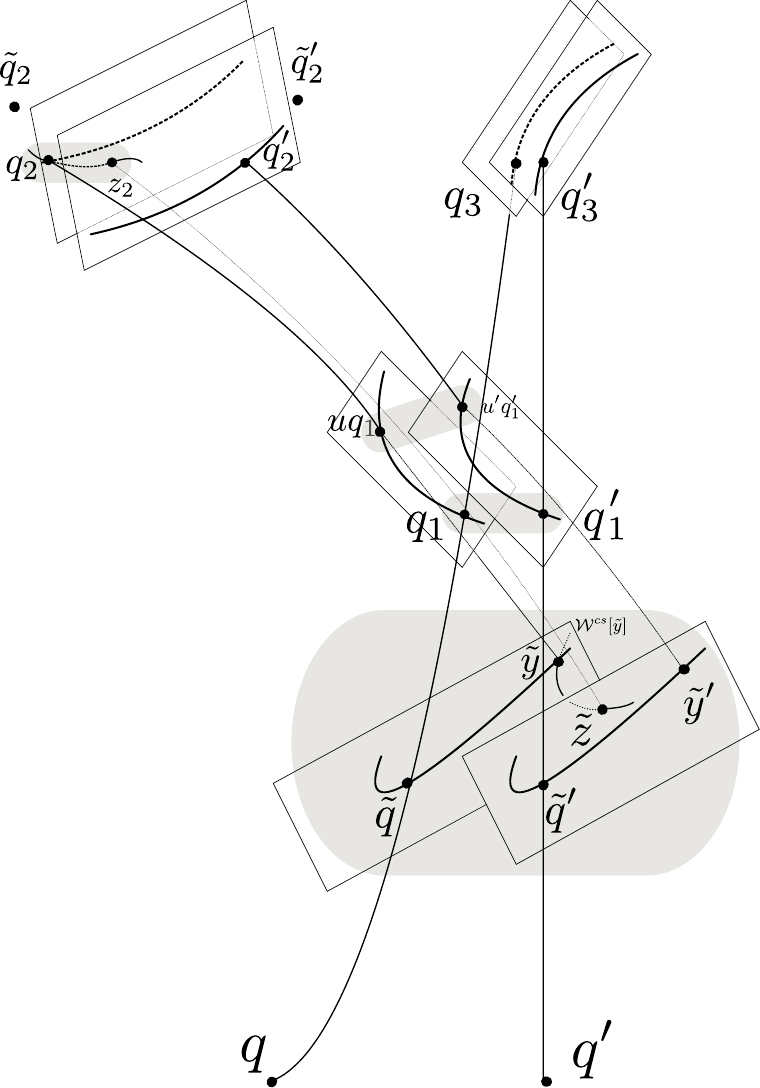}
	\caption{Left: The general arrangement of points in the argument.
	Right: More detailed notation.
	Points in a transparent gray blob are exponentially close in the parameter $\ell$, for some exponent that depends only on the Lyapunov spectrum, and constants that depend on the Lusin set.}
	\label{fig:outline}
\end{figure}

We will make use of the notion of $Y$-configurations \index{$Y(q_1,u,\ell,\tau,t_0)$}$Y(q_1,u,\ell,\tau,t_0)$ with $q_1\in X, u\in \cB_0(q_1)$ and $\ell,\tau,t_0\geq 0$ (see \autoref{ssec:definitions_y_configurations}).
The data in a $Y$-configuration determines the points $q,q_2,q_3$ (see \autoref{def:Y_configuration} and \autoref{fig:outline}), and for a compact set $K$ we say $Y\in K$ if $q,q_1,q_2,q_3\in K$.
Two $Y$-configurations $Y,Y'$ can be \emph{bottom-linked} (see \autoref{def:bottom_linked_Y_configurations}) which means that $\ell=\ell'$ and $q'\in \gB_0^-[q]$, and \emph{right-linked} (see \autoref{def:right_linked_Y_configurations}) if in addition $t_0=t_0'$.

\subsubsection{Left-balanced, $ij$-top-balanced, and $ij$-balanced}
	\label{sssec:left_balanced_ij_top_balanced_and_fully_balanced}
We will say (see \autoref{ssec:left_balanced_y_configurations}) that a $Y$-configuration is \emph{$(C,\ve)-$left-balanced} if $\norm{\cA_2(q_1,u,\ell,\tau)}\in \left[\tfrac 1C \ve,C\ve\right]$.
Note that this condition does not depend on the upper-right leg of the $Y$-diagram (i.e. the parameter $t$).
Furthermore, given $\ve$ there exists a unique $\tau_{(\ve)}:=\tau_{(\ve)}(q_1,u,\ell)$ such that $Y(q_1,u,\ell,\tau_{(\ve)},t)$ is $(1,\ve)$-left-balanced (by monotonicity of $\norm{\cA_2}$, see \autoref{eq:A:future:bilipshitz}).

Given $ij \in \tilde{\Lambda}$, we will say (see \autoref{def:top_balanced_y_configuration}) that a $Y$-diagram is \emph{$(C,ij)$-top-balanced} if $q_3=g_{t_3}^{ij}(q_1), q_2 = g_{t_2}^{ij}(uq_1)$ and $|t_2-t_3|\leq C$.
Again, note that given $q_1,u,\ell, \tau$ there exists a unique $t_0$ such that $Y(q_1,u,\ell,\tau,t_0)$ is $(0,ij)$-top-balanced.

Finally, we will (see \autoref{def:y_ij_fully_balanced}) denote by \index{$Y_{ij,balanced}(q_1,u,\ell)$}$Y_{ij,balanced}(q_1,u,\ell)$ the unique $Y$-configuration which is both $(1,\ve)$-left-balanced and $(0,ij)$-top-balanced.

The arguments of (i),(ii), (iii) of \autoref{sssec:choosing_the_eight_points_intro_outline} above and Fubini's theorem show that for an appropriate compact $K$ and an almost full density set of $\ell$, there are many $Y$-configurations with $Y_{ij,balanced}(q_1,u,\ell)\in K$, see \autoref{eq:Yij:balanced:in:K} below for the exact statement.

\subsubsection{Coupled $Y$-configurations}
	\label{sssec:coupled_y_configurations_intr_outline}
We say that two $Y$-configurations $Y = Y_{ij}({q}_1, u,\ell)$ and $Y' =
Y_{ij}({q}_1',u',\ell')$ with the same $ij$ are {\em coupled} if
$\ell = \ell'$, $u = u'$, ${q}(Y') \in \cW^s_{loc}[{q}(Y)]$, and also the conditions in (iv) of \autoref{sssec:choosing_the_eight_points_intro_outline} above are satisfied.

Then the argument described (iv), specifically \autoref{prop:can:avoid:most:Mu}, shows that we can (for most values
of $\ell$) choose points ${q}_1$, ${q}_1'$ such that for most
$u$ and all $ij$, the $Y$-configurations $Y_{ij}({q}_1,u,\ell)$ and
$Y_{ij}({q}_1',u,\ell)$ are both good and also are coupled
(see \autoref{sssec:choice_of_parameters_nr_2} ``Choice of parameters \#2'' below for the precise statement).

We then choose $u$ as in (v).
(See \autoref{claim:12:12}, \autoref{claim:large_set_of_u}, and \autoref{sssec:choice_of_parameters_5} ``Choice of parameters \#5'').

We are now almost
done, except for the fact that the lengths of the legs of
$Y_{ij} = Y_{ij}({q}_1,u,\ell)$ and $Y_{ij}' = Y_{ij}({q}_1',u,\ell)$ are not same. (The bottom
leg of $Y_{ij}$ has length $\ell$, and so does the bottom leg of $Y'_{ij}$, but
the two top legs of $Y_{ij}$ can potentially have different lengths than
the corresponding legs of $Y'_{ij}$).
We show that the lengths of the corresponding legs are close (see \autoref{eqn:tau_parameters_can_be_chosen_close} and \autoref{eq:t:ij:t:ij:prime:close}).
Then we proceed to part (vi).




\section{The space \texorpdfstring{$\cC$}{C} and the cocycles \texorpdfstring{$L\cC$}{LC} and \texorpdfstring{$\bbH$}{bold H}}
	\label{sec:c_and_lc_and_boldH}



\subsection{Some results on nilpotent Lie groups}
	\label{ssec:some_results_on_nilpotent_lie_groups}

In this subsection, we will take the group \index{$N$}$N$ to be the strictly subresonant maps $\bbG^{ssr}$ and $U^+$ to be given by groups as in \autoref{def:compatible_family_of_subgroups}.
In particular, the nilpotent Lie group $N$ will have nontrivial automorphisms given by the conjugation action of the subresonant maps $\bbG^{sr}$ (cf. \autoref{prop:group_structure_on_subresonant_maps}).

\subsubsection{Setup}
	\label{sssec:setup_some_results_on_nilpotent_lie_groups}
All our nilpotent Lie groups are connected and simply connected, so that most arguments can be translated to Lie algebras via the exponential map.
The Lie groups are denoted in Roman font, such as $N,U^+$ and the Lie algebras in corresponding Gothic script $\index{$n$@$\frakn$}\frakn,\index{$u$@$\fraku^+$}\fraku^+$ etc. without further mention.
Assume that we have a weight filtration \index{$N_{\leq \lambda_i}$}$N_{\leq \lambda_i}$ such that $[N_{\leq \lambda_i},N_{\leq \lambda_j}]\subseteq N_{\leq \lambda_i+\lambda_j}$ and all $\lambda_i$ are strictly negative.
It also induces a filtration on $U^+$ by $\index{$U^+_{\leq \lambda_i}$}U^+_{\leq \lambda_i}:=U^+\cap N_{\leq \lambda_i}$ and analogously on Lie algebras.

\begin{remark}[On filtrations]
	\label{rmk:on_filtrations}
	We could either take $N=N_{\leq -1}$ and then the induced filtration is the lower central series.
	Or, in cases of interest for dynamics we could take on $\frakn$ the Lyapunov filtration.
\end{remark}

We are grateful to Brian Chung for providing the proof of the next result.

\begin{lemma}[Transverse structure of quotients]
	\label{lem:transverse_structure_of_quotients}
	Let $N$ be a nilpotent Lie group and $U^+$ a Lie subgroup.
	Suppose that $\frakv\subset \frakn$ is a subspace satisfying $\frakv\oplus \fraku^+=\frakn$.
	Consider the maps
	\begin{align*}
		\frakv \oplus \fraku^+& \xrightarrow{\phantom{XX}\exp_2\phantom{XX}}N \xrightarrow{\phantom{XX}\log\phantom{XX}}\frakn\\
		v\oplus u & \mapsto \exp(v)\cdot \exp(u) =:\exp_2(u,v)
	\end{align*}
	and call their composite \index{$ch_{\frakv}$}$ch_{\frakv}$, which is a polynomial map.
	Assume that we additionally have:
		\[
			\frakn_{\leq \lambda_i} = (\fraku^+\cap \frakn_{\leq \lambda_i})\oplus (\frakn_{\leq \lambda_i} \cap \frakv).
		\]
	where $\fraku_{\leq \lambda_i}^+:=\fraku^+\cap \frakn_{\leq \lambda_i}$ is the induced filtration on $\fraku^+$.
	Then:
	\begin{enumerate}
		\item The two weight filtrations on $\frakv$ agree, one coming from the identification with $\frakn/\fraku^+$, the other from viewing it as a subspace of $\frakn$.
		Equip $\frakv$ with the corresponding subresonant structure.
		The linear projection map $\pi_{\frakv}\colon \frakn \to \frakv$ is subresonant, as is the inclusion $\frakv\into \frakn$.
		\item The map $ch_{\frakv}$ is a strictly subresonant polynomial map and has a polynomial inverse, which is also strictly subresonant.
		\item Furthermore, given two such $\frakv_1,\frakv_2$ as above, the composition $\pi_{\frakv_1}\circ ch_{\frakv_1}^{-1}\circ ch_{\frakv_2}$ viewed as a map $\frakv_2\to \frakv_1$ (where $\pi_{\frakv_1}$ is the linear projection in $\frakn$ to $\frakv_1$ along $\fraku^+$) is a subresonant map.
	\end{enumerate}
\end{lemma}
\begin{proof}
	Observe that the assumption is equivalent to the existence of a direct sum decomposition for each of $\fraku^+,\frakv,\frakn$ so as to have
	\[
		\frakn_{\lambda_i} = \frakv_{\lambda_i} \oplus \fraku^+_{\lambda_i} \quad \text{ and } \quad \frakg_{\leq \lambda_i}=\bigoplus_{\mu\leq \lambda_i} \frakg_{\mu} \quad \forall \frakg\in \{\frakn,\fraku^+,\frakv\}.
	\]
	The uniqueness of weight filtration on $\frakv$ follows, and the corresponding inclusions and projections are subresonant (linear) maps.

	Recall next the general Baker--Campbell--Hausdorff formula
	\[
		e^{v}e^{u} = e^{u+v + \psi_1(u,v) + \psi_2(u,v) + \cdots}
	\]
	where $\psi_i(u,v)$ involves $i$-fold commutators of $u$ and $v$.
	The map $ch_\frakv$ is precisely the exponent on the RHS.

	To check that $ch_{\frakv}$ is subresonant, observe that the first term is the identity map so it suffices to check that each subsequent term is a subresonant polynomial map from $\frakn$ to itself of strictly negative degree (see also \autoref{sssec:filtration_on_polynomial_maps}).
	Each of these maps however is a composition of maps of the form $\Lambda^2\frakn \xrightarrow{[,]}\frakn$ and $\frakn\otimes \frakn\to \Lambda^2\frakn$, which are clearly subresonant (and linear), as well as the polynomial map $\frakn\to \frakn\otimes \frakn$ given by $(v,u)\mapsto v\otimes u$ for the decomposition $\frakn=\frakv\oplus \fraku^+$.
	This is also a subresonant map, again since it is a composition of the subresonant linear map $\frakn\to \frakn\oplus \frakn$ given by $\pi_{\frakv}\oplus 1-\pi_{\frakv}$ and the subresonant polynomial map $\frakn\oplus\frakn\to \frakn\otimes\frakn$ given by multiplying coordinates.

	Inverses (when they exist) of subresonant maps are subresonant, and since $ch_{\frakv_i}$ has invertible derivative the last claim also follows.
\end{proof}

\begin{corollary}[Subresonant structure on the quotient space]
	\label{cor:subresonant_structure_on_the_quotient_space}
	The manifold $C:=N/U^+$ carries a natural subresonant structure, given by selecting any $\frakv\subset \frakn$ satisfying the conditions of \autoref{lem:transverse_structure_of_quotients} and using $ch_{\frakv}$ to obtain charts on it.

	Furthermore, let \index{$A$@$\Aut(N,U^+)$}$\Aut(N,U^+)$ be the Lie group of automorphisms of $N$ that preserve the subgroup $U^+$ (as a set), and additionally assume that the action of $\Aut(N,U^+)$ preserves the weight filtration on $N$.
	Then $\Aut(N,U^+)$ acts by subresonant maps on $C$.

	Additionally, there exists a linear representation $LC$ of $\Aut(N,U^+)$ and polynomial embedding $C\into LC$ such that subresonant automorphisms of $C$ are induced by linear maps on $LC$.
	Furthermore, $LC$ admits a filtration by subrepresentations, such that the associated quotients are tensor powers of the action on $\frakn/\fraku^+$.
\end{corollary}
\begin{proof}
	The compatibility of subresonant charts on $C$ is the content of \autoref{lem:transverse_structure_of_quotients}(iii).

	The linear action of $\Aut(N,U^+)$ on $\frakn$ by Lie algebra automorphisms is conjugated by the exponential map to the action at the Lie group level (by the Baker--Campbell--Hausdorff formula).
	Furthermore, given one $\frakv\subset \frakn$ satisfying the assumptions of \autoref{lem:transverse_structure_of_quotients}, its image under elements of $\Aut(N,U^+)$ still satisfies the assumptions, since $\Aut(N,U^+)$ is required to preserve the weight filtration of $\frakn$.
	Following through the identifications and induced maps, it follows that $\Aut(N,U^+)$ acts on $C$ by subresonant polynomial maps.

	To construct the linearizing representation $LC$, we apply the construction from \autoref{sssec:veronese_embedding_for_subresonant_spaces} and use \autoref{prop:linearization_of_subresonant_maps}.
\end{proof}

\subsubsection{Construction of the associated bundles}
	\label{sssec:construction_of_the_associated_bundles}
Working on the finite cover $X$, we apply now the above formalism to the setting of unstable manifolds $\cW^u$, assuming a compatible family of subgroups $U^+(x)\subset \bbG^{ssr}(x)=:N$ as in \autoref{def:compatible_family_of_subgroups}.
Define \index{$C$@$\cC(x)$}$\cC(x):=\bbG^{ssr}(x)/U^+(x)$, which is a measurable fiber bundle over a set of full measure.
The fibers are equipped with a subresonant structure, as per \autoref{cor:subresonant_structure_on_the_quotient_space}.
Because the $g_t$-dynamics preserves $U^+(x)$ as a subgroup of $\bbG^{ssr}(x)$, it also acts by subresonant maps on $\cC(x)$.
Let $L\cC\to X$ be the linearization cocycle provided by \autoref{cor:subresonant_structure_on_the_quotient_space}, so that we have an equivariant embedding \index{$j$}$j\colon \cC(x)\into L\cC(x)$.\index{$L\cC(x)$}

Recall next that we have the ``Gauss--Manin'' connection $P^{ssr}_{GM}(q,q')$ between $\bbG^{ssr}(q)$ and $\bbG^{ssr}(q')$, defined in \autoref{sssec:holonomies_and_relating_different_groups}.
\autoref{def:compatible_family_of_subgroups} requires that for $x'\in U^+[x]$ and in a set of full measure, we have that $P^{sr}_{GM}(x,x')(U^+(x))=U^+(x')$.
Therefore, we have induced identifications
\[
	\index{$P^{L\cC}_{GM}$}P^{L\cC}_{GM}\colon L\cC(x)\to L\cC(x') \quad \text{for }x'\in U^+[x]\cap X_0
\]
where $X_0$ is a fixed set of full measure.

\begin{proposition}
	\label{prop:LC_Jordan_normal_form}
The cocycle $g_t\colon L\cC\to L\cC$ over $g_t\colon X\to X$ can be
put in Jordan normal form.
\end{proposition}
\begin{proof}
By the definition of $X$, the cocycle $g_t\colon \liessr\to \liessr$ can be put in to Jordan normal form when lifted to the finite cover $X\to Q$.
Consider the cocycle $E(x):=\Lie (\bbG^{ssr}(x))/\Lie(U^+(x))$, which agrees with the derivative $D_{U^+(x)}g_t\colon \cC(x)\to \cC(g_tx)$ at $U^+(x)$.
As a quotient of the cocycle $ g_t\colon \liessr\to \liessr$, $E$ can be put in to Jordan normal form when lifted to $X$ by \autoref{prop:jordan_normal_form_and_standard_operations}(i).  
The linearization cocycle $g_t\colon L\cC\to L\cC$ has a block triangular structure, whose diagonal blocks are subcocycles of tensor powers of $E$, as in \autoref{cor:subresonant_structure_on_the_quotient_space}.
By \autoref{prop:jordan_normal_form_and_standard_operations}(iii) and \autoref{prop:jordan_normal_form_and_standard_operations}(ii), $L\cC$ can be put in to Jordan normal form.
\end{proof}

\subsubsection{The set $\cC[x]$}
\label{sec:subsubsec:cCbracketX}
Let \index{$C$@$\cC[x]$}$\cC[x] = \{ h U^+[x] \st h \in \cC(x) \}$. Then the natural map $\cC(x)
\to \cC[x]$ is a bijection. We will view the points of $\cC[x]$ as 
submanifolds inside $\cW^u[x]$ and speak about Hausdorff distances
(denoted by \index{$d$@$\dist_{Hausd}$}$\dist_{Hausd}$). 

\subsubsection{The tautological subspace}
\label{sec:subsec:def:LC:tau}
Define \index{$L\cC^{\tau}(x)$}$L\cC^{\tau}(x)\subset L\cC(x)$ to be the tautological line
spanned by $j(U^+(x))$ viewed as a point in $j(\cC(x))\subset
L\cC(x)$.
Quotienting by it has the effect of moving this basepoint to the origin of the vector space on which we have a linear cocycle.

\subsubsection{The cocycle $\bbH$}
	\label{ssec:the_bold_h_cocycle}
Define the quotient cocycle:
\begin{align}
	\label{eq:bold_H_definition}
	\bbH := L\cC / L\cC^{\tau}
\end{align}
Vectors in \index{$H$@$\bbH$}$\bbH$, and some of its subbundles $\bbE$ defined later, will be denoted $\bfv$.
When a distinction is necessary, we will write\index{$j$@$\bbj$}
\[
	\bbj \colon \cC \into \bbH, \quad  \bbj(\cU) = \bigl(L \circ j (\cU)\bigr)/L\cC^{\tau}
\]
for the natural inclusion.
By \autoref{thm:jordan_normal_form}\ref{Jordanextention:iiiii}, since $L\cC$ can be put in Jordan normal form on $X$, so can $\bbH$.

\subsubsection{Measurable connections on $\bbH$}
	\label{sssec:measurable_connections_on_bold_h}
Recall that the finite cover $X$ is constructed by bringing $\Lie \bbG^{ssr}$ to Jordan normal form.
The measurable connections $P^{\pm}$ lift to $X$ and preserve all $g_t$-invariant subcocycles of $\Lie \bbG^{ssr}$ on $X$ by the construction of lifts of stable and unstable manifolds on $X$ and the Ledrappier invariance principle (see \autoref{sssec:measurable_connections_on_finite_covers} and \autoref{thm:ledrappier_invariance_principle}).
The cocycle $\bbH$ is obtained from elementary linear-algebraic constructions with $\Lie \bbG^{ssr}$ and $\Lie U^+$, hence acquires induced measurable connections $P^{\pm}$.

\subsubsection{Metric structures, distances}
	\label{sssec:metric_structures_distances} 
Using the norm on $\bbH$ provided by \autoref{prop:good_norms}, together with the embedding $\bbj \colon \cC \into \bbH$, yields a distance function on $\cC(x)$ denoted \index{$d_{\cC}$}$d_{\cC}$.

Recall that the metric $d^u$ is defined
in \autoref{sec:subsec:distances}. 
Recall that the partition $\gB_0$ is chosen so that any atom of
$\gB_0$ has diameter at most $1$ in the $d^u(\cdot, \cdot)$ metric (see \autoref{rmk:choice_of_measurable_partition_and_du_smallness}).
We make a measurable choice of a distinguished point \index{$x^\ast$}$x^\ast$ in each
atom $\gB_0[x]$.
Let \index{$B^\ast[x]$}$B^\ast[x]$ be the ball in the $d^u(\cdot,
\cdot)$ metric of radius $1$ centered at $x^\ast$. Then, for all $x$,
\begin{displaymath}
\gB_0[x] \subset B^\ast[x] \subset \cW^u[x]. 
\end{displaymath}
Let also \index{$d^x_{\cH,loc}$}$d^x_{\cH,loc}$ denote the Hausdorff distance on $\cC[x]$
when its elements are viewed as subsets of $B^\ast[x]$. 
\begin{align}
	\label{eq:Hausdorff_distance_def}
	d^x_{\cH,loc}(c_1,c_2):=\dist_{Hausd}\left(c_1\cap B^\ast[x], c_2\cap B^\ast[x]\right).
\end{align}
We also have the reference basepoint $\index{$c_0(x)$}c_0(x):=U^+[x]\in \cC[x]$, which by assumption is $g_t$-equivariant.

\begin{proposition}[Relation to Hausdorff distance]
	\label{prop:relation_to_hausdorff_distance}
	There exist $g_t$-tempered functions $A_1(x), A_2(x), A_3(x)$ with the following properties.
	\begin{enumerate}
		\item Let \index{$v$@$\frakv(x)$}$\frakv(x)\subset \frakg^{ssr}(x)$ be the family of transversals obtained as orthogonal complements to $\fraku^+(x)$ for admissible metrics (see \autoref{sssec:admissible_metrics_on_admissible_cocycles}).
		Then these transversals satisfy the assumptions in \autoref{lem:transverse_structure_of_quotients}.
		Endow $\cC(x)$ with euclidean metrics \index{$d_{ch}$}$d_{ch}$ coming from $\frakv(x)$ via the charts $ch_{\frakv(x)}$.
		Then we have
		\[
		\frac {1}{A_1(x)}\cdot d_{\cC}(c_1,c_2)
		\leq
		d_{ch}(c_1,c_2) \leq 
		A_1(x)\cdot d_{\cC}(c_1,c_2)
		\]
		for $c_1,c_2\in \cC(x)$ in the ball of radius $1$
		for $d_{\cC}$ around the basepoint $c_0$.
		\item Additionally, for $c_1,c_2$ intersecting
                  $B^\ast[x]$ we have:
		\[
			\frac {1}{A_2(x)}\cdot d_{\cC}(c_1,c_2)
			\leq
			d^x_{\cH,loc}(c_1,c_2) \leq 
			A_2(x)\cdot d_{\cC}(c_1,c_2)
		\]
		where $d^x_{\cH,loc}$ is the Hausdorff distance from
                \autoref{eq:Hausdorff_distance_def}.
              \item Also, for $c_1, c_2$ as above we have
\begin{multline*}
\frac{1}{A_3(x)} d^x_{\cH,loc}(c_1,c_2) \le \min(r(x),
\dist_{Hausd}\left(c_1\cap B^u(x;r(x)), c_2\cap B^u(x;r(x)) \right)) \\
\le  A_3(x) d^x_{\cH,loc}(c_1,c_2)
\end{multline*}
where $r(x)$ is as in \autoref{rmk:on_lyapunov_radius}.  
	\end{enumerate}
\end{proposition}

\begin{proof}
	The family of transversals $\frakv(x)\subset\frakg^{ssr}$ satisfies the assumptions of \autoref{lem:transverse_structure_of_quotients} by the definitions of admissible metrics, since the Lyapunov subspaces will yield the required direct sum decompositions of $\frakg^{ssr}$.

	The derivative of the chart map $ch_{\frakv}$ at the origin is the identification between $\frakv$ and $\frakg^{ssr}/\fraku^+$ and is in fact an isometry for admissible metrics on both sides.
	Furthermore the (injective) chart maps $ch_{\frakv}\colon \frakv \to \cC(x)\subset L\cC(x)$ are given by polynomial expressions.
	We can identify all the objects in question with a model situation, independently of $x$, and compare metrics at the cost of introducing multiplicative factors that vary in a tempered way along the $g_t$-orbit of $x$.
	The first claim then follows.

	For part (ii) we continue to work in a fixed model situation (with fixed metric, groups $\bbG^{ssr},U^+$, etc.) independent of $x$.
	Let $c_0 \in \cC$ be the basepoint corresponding to $U^+$ viewed as an element of $\cC$.
	For the rest of the proof, we will write $U^+[x]\subset
        \cW^u[x]$ with the understanding that the argument applies to
        a full measure set of $x$, up to tempered multiplicative
        factors.

	Now for any $X\in \frakv$ with $\norm{X}\leq 1$, the claimed upper bound for $d_{\cH}(c_0,\exp(X)c_0)$ is clear.
	By applying further transformations of the form $\exp(Y)$ with $Y\in \frakv,\norm{Y}\leq 1$ the upper bound on $d_{\cH}(c_1,c_2)$ for arbitrary $c_1,c_2$ in the unit ball of $\cC$ follows.
	Both claims follow since $\exp(Y)$ has bounded derivatives on the ball of radius $1$ when $\norm{Y}\leq 1$.

	For the lower bound, using an $\exp(Y)$ as before, it suffices to show it for $d_{\cH}(c_0,\exp(X)c_0)$.
	Now, restrict first to $X$ in the unit sphere in $\frakv$, it yields a compact (and finite-dimensional) family of vector fields on $\cW^u$.
	Each vector field is not identically zero on $U^+[x]\subset \cW^u[x]$ (otherwise it would give a trivial transformation of $U^{+}[x]$).
	Since the vector fields are also real-analytic, there exist constants $\ve,c>0$ such that for any $X\in\frakv$ with $\norm{X}=1$, there exists $p\in U^+[x]\cap \cW^u_{loc}[x]$ such that $\norm{X(p)}\geq c$ and the angle between $X(p)$ and $T_pU^+[x]$ is at least $\ve$, where $X(p)$ denotes the vector field induced by $X$ at $p$.

	It follows that for any $X\in\frakv\setminus 0$ with $\norm{X}_{\frakv}\leq 1$ there exists $p\in U^+[x]\cap \cW^u_{loc}[x]$ such that the angle between $X(p)$ and $T_pU^+[x]$ is at least $\ve$, and $\norm{X(p)}\geq c\norm{X}_{\frakv}$.
	Furthermore there exists a constant $A>0$ such that $X$ has $C^1$-norm bounded by $A'$ on $\cW^u_{loc}[x]$.
	The desired lower bound now follows for $\norm{X}\leq \frac{1}{A''}\leq 1$ for some $A''=A''(A,c,\ve)$, by exponentiating the vector field $X$ and following along $p$.
	Note finally that scalings are among the subresonant maps of
        $\cW^u[x]$, so by rescaling (or equivalently redefining the
        metric) we can assume that the result holds for $\norm{X}\leq
        1$ and modified constants.

Part (iii) follows from Part (ii) and a compactness argument. 
\end{proof}

\subsubsection{Local Hausdorff distance on unstable}
	\label{sssec:local_hausdorff_distance_on_unstable}

\begin{proposition}[Hausdorff distance and norm of vector]
	\label{prop:hausdorff_distance_and_norm_of_vector}\label{lemma:hausdorff:distance:to:norm}

Consider $g\in \bbG^{ssr}(\cW^u[x])$ and  $\cU= g U^+[x] \in \cC[x]$.  Let $$\bbv= \bbj( g U^+(x)) \in \bbj(\cC(x))\subset \bbH(x).$$
Then there exists a measurable function $c_1(x)>0$ with the
following property: Suppose $\|\bfv\| \le 1$.
Then we have        
	\[
		\frac{1}{c_1(x)}\norm{\bbv} \leq d^x_{\cH,loc}(U^+[x], \cU) \leq c_1(x) \norm{\bbv}.
	\]
(We use the convention that $d^x_{\cH,loc}(U^+[x], \cU) = 1$ if $\cU$
does not intersect $B^\ast[x]$). Also, if $r(x)$ is as in
\autoref{rmk:on_lyapunov_radius},
\begin{displaymath}
  \frac{1}{c_1(x)}\norm{\bbv} \leq
\dist_{Hausd}\left(U^+[x] \cap B^u[x,r(x)], \cU \cap B^u[x,r(x)]\right) 
    \leq c_1(x) \norm{\bbv}.
\end{displaymath}
(In the above equation we use the convention that $\dist_{Hausd}(A,B)
= r(x)$ if $A$ or $B$ is empty). 

\end{proposition}

\begin{proof}
This follows from \autoref{prop:relation_to_hausdorff_distance} (ii)
and (iii), and the fact that for a.e. $x \in X$, the image of the
linearization map $j: \cC(x) \to L\cC(x)$ is transverse to the
tautological subspace $\cC^\tau(x)$. 
\end{proof}

Recall that \index{$L_y$}$L_y: \cW^u[y] \to L\cW^u(y)$ denotes the linearization map. 
Given $h \in \bbG^{ssr}(\cW^u[y])$, we denote by \index{$L$@$(L_y
  h)$}$(L_y h)$ the linear endomorphism of $L\cW^u(y)$ such that $(L_y
h) \circ L_y = L_y \circ h$. 

\begin{lemma}
	\label{lemma:L:and:C}
There exist measurable functions $C_1: X \to \reals_+$ and $c_2: X
\to \reals_+$ such that the following holds:
\begin{itemize}
\item[{\rm (a)}] Suppose $h \in \bbG^{ssr}(\cW^u[y])$ is such that
  $\|L_y h - I \| \le c_2(y)$. Then
\begin{displaymath}
d_{\cH,loc}^y(h U^+[y], U^+[y]) \le C_1(y) \|L_y h - I\|.  
\end{displaymath}
\item[{\rm (b)}] For any $\cU' \in \cC[y]$ with $d_{\cH,loc}^y(\cU',
  U^+[y]) < c_2(y)$ there exists $h \in \bbG^{ssr}(\cW^u[y])$ such that
\begin{displaymath}
 \|L_y h - I \| \le C(y) d_{\cH,loc}^y(\cU', U^+[y]).
\end{displaymath}
\end{itemize}
\end{lemma}

\begin{proof} 
We begin the proof of (a). 
Let $\tilde{h} = L_y h$. Then, $L_y \circ h = \tilde{h} \circ L_y$.
Note that by \autoref{lemma:L:is:bilipshitz}, 
$L_y$ is bilipshitz with constant
$\kappa(y)>1$. 
Then, for any $z \in \cW^u[y]$ with $d^u(y,z) \le 1$, we have
$\|L_y z \| \le \|L_y y\| + \kappa(y)$. Then, 
\begin{multline*}
d^u(h z, z) \le \kappa(y) \| L_y(h z) - L_y(z) \| = \| (\tilde{h} -I)
L_y z \| \le \\
\le \|\tilde{h} - I \| \|L_y z\| \le \kappa(y) (\|L_y y\|+\kappa(y))
\| \tilde{h}-I \|. 
\end{multline*}
This implies (a). 

We now begin the proof of (b). We can choose $h$ so that $h = \exp(H)$
where $H \in \frakv(y)$ where $\frakv(y)$ is as in
\autoref{prop:relation_to_hausdorff_distance}.
By \autoref{prop:relation_to_hausdorff_distance}, there exists
$C(y)$ such that
\begin{displaymath}
C(y)^{-1} \|H\|_y \le d^y_{\cH,loc}(U^+[y], h U^+[y]) \le C(y) \|H\|_y,
\end{displaymath}
where $\| \cdot \|_y$ is the norm on $\frakv(y)$ considered in
\autoref{prop:relation_to_hausdorff_distance}.

Let $\rho_y$ denote the homomorphism from $\bbG_{ssr}(W^u[y])$ to
  $GL(LW^u[y])$, and let $\phi_y$ denote the corresponding map of
the Lie algebras. Then, there exists $C_2(y) > 1$ such that
\begin{displaymath}
C_2(y)^{-1} \|H\|_y \le \|\phi_y(H)\|_y \le C_2(y) \|H\|_y, 
\end{displaymath}
where in the middle term $\| \cdot \|_y$ is a norm on the Lie algebra
of $GL(LW^u[y])$. 

Note that given $C_3(y)$ there exists $C_4(y)$ such that for any $M$
in the Lie algebra of 
$\GL(LW^u[y])$ with $\|M \|_y \le C_3(y)$, we have
\begin{displaymath}
\|\exp(M) - I \|_y \le C_4(y) \|M\|_y
\end{displaymath}
Then, 
\begin{displaymath}
\|L_y h - I \|_y = \|\rho_y(\exp H) - I \|_y = \| \exp( \phi_y(H)) - I \|_y
\le C \|H\|_y.  
\end{displaymath}
\end{proof}



\subsection{The change of basepoint map}
	\label{ssec:the_change_of_basepoint_map}

\subsubsection{Setup}
	\label{sssec:setup_the_change_of_basepoint_map}
Suppose $U^+(x)$ is a compatible family of subgroups in the sense of
\autoref{def:compatible_family_of_subgroups}, with associated full
measure set $X_0$ as in \autoref{sssec:construction_of_the_associated_bundles}.
Let $\cC\to X_0$ and $L\cC\to X_0$ be the associated fiber bundle and linearization bundle from \autoref{sssec:construction_of_the_associated_bundles}.
Recall that the constructions in \autoref{sssec:construction_of_the_associated_bundles} yield 
linear maps
\[
	\index{$P^{L\cC}_{GM}(x,ux)$}P^{L\cC}_{GM}(x,ux) \colon L\cC(x)\to L\cC(ux)
\]
assuming $x,ux\in X_0$ and $u\in U^+(x)$.

The properties of the next construction are the subject of Lemma 6.5 of \cite{EskinMirzakhani_Invariant-and-stationary-measures-for-the-rm-SL2Bbb-R-action-on-moduli-space}.
\begin{definition}[u-star]
	\label{def:u_star}
	With notation as above, for $x,ux\in X_0$ and $u\in U^+(x)$, define:
	\[
		\index{$u_*(x,ux)$}u_*(x,ux) \colon L\cC(x)\to L\cC(ux)
	\]
	as $P^{L\cC}_{GM}(x,ux)$.
\end{definition}

\begin{remark}[Equivariance and endpoint dependence of $u_*$]
	\label{rmk:equivariance_and_endpoint_dependence_of_u_star}
	From the definition of $u_*$ via the Gauss--Manin connection, it follows that whenever $x_2=ux_1$ the map $u_*(x_1,x_2)$ depends only on $x_1,x_2$ and not the particular choice of $u\in U^+(x_1)$ such that $x_1=ux_2$.

	Furthermore, it is also immediate that $u_*(x_1,x_2)$ intertwines the action of $g_t$, namely $u_*(g_t x_1, g_t x_2) \circ g_t = g_t \circ u_*(x_1,x_2)$.
\end{remark}
Note that for measure-generic points $x,x'\in \cW^u[x]$, we only have a map $P^+(x,x')\colon L\cC(x)\to L\cC(x')$, but in general no ``Gauss--Manin'' connection such as the one on $L\cW^u$.
Nonetheless, on cocycles and points for which $P^+$ and $u_*$ are both defined, they agree by \autoref{prop:u_star_preserves_flags} on \emph{inert subspaces}, which are defined in \autoref{ssec:inert_subspaces}.

\begin{proposition}[$u_*$ preserves flags]
	\label{prop:u_star_preserves_flags}
	Let $L\cC(x)^{\geq \bullet}$ be the backwards Lyapunov flag of the cocycle $L\cC$ at $x$ in the set of full measure $X_0$.

	Then, for $u\in U^+(x)$ with $ux\in X_0$, we have that $u_*(x,ux)$ takes the backwards Lyapunov flag at $x$ to that at $ux$, and its action on the associated graded agrees with the holonomy map from \autoref{sssec:standard_measurable_connection}.
\end{proposition}

\begin{proof}
	The collection of maps $P^{L\cC}_{GM}(x,ux)$ give holonomies in the sense of \autoref{sssec:holonomies_discussion}.
	They are also smooth when $u$ varies in $U^+(x)$.

	The statement of the proposition is then a consequence of standard properties of holonomies, see \autoref{sssec:compatibility_of_holonomies_and_standard_measurable_connection}.
\end{proof}

\subsubsection{On $\bbH$}
	\label{sssec:on_bold_h}
Recall that in  \autoref{eq:bold_H_definition} we are considering a quotient of $L\cC$ by the one-dimensional tautological subbundle corresponding to the choice of $U^+(x)$, viewed as the identity coset in $\cC(x)=\bbG^{ssr}(x)/U^+(x)$.
The statement of \autoref{prop:u_star_preserves_flags} holds in this
case as well, because both the measurable connection and
$P^{L\cC}_{GM}$ preserve this tautological section.

\begin{lemma}
	\label{lemma:ustar:bounded}
For every $\delta>0$ there is a compact set $K\subset X$ of measure $1-\delta$, and $C(\delta)<+\infty$ such that given $x\in K$ and $y\in\cB_0[x]\cap K$, we have $\|u_*(x,y)\|\leq C(\delta)$.
\end{lemma}
Remember $u_*(x,y)$ is defined in \autoref{def:u_star} and the norm is measured with respect to the good norms given by \ref{prop:good_norms}. Remember also that $u_*(y,x)=(u_*(x,y))^{-1}$.

\begin{proof}
Let $r(x,y)=\max\{\|u_*(x,y)\|, \|u_*(y,x)\|\}$.
  Let us define on $X\times X$ the measure $\hat \nu$ by $\int\phi(x,z)d\hat\nu=\int\int\phi(x,z)d\nu^{\cB_0}_x(z)d\nu(x)$, where $\nu_x^{\cB_0}$ are the conditional measures from \autoref{sssec:measures_on_leaves}.
  Notice that $\hat\nu$ is supported on the set of points $(x,z)$ s.t. $z\in\cB_0[x]$.
  Since $r(x,z)<+\infty$ for $\hat\nu$-a.e. $(x,z)$ and $u_*$ depends measurably on $(x,z)$, we have as $C\to +\infty$ that $\hat\nu\left(\{(x,z):r(x,z)> C\}\right)\to 0$.
 
  Let $C>0$ be such that $$\hat\nu(\{(x,z):r(x,z)>C\})< \frac{1}{4}\delta$$ and let $C(\delta)=C^2$. We get by Markov's inequality and since $$\hat\nu\left(\{(x,z):r(x,z)> C\}\right)=\int\nu^{\cB_0}_x(\{z:r(x,z)> C\}))d\nu(x)<\frac{1}{4}\delta.$$ that \begin{eqnarray*}
  &&\nu\left(\{x\in X:\nu^{\cB_0}_x\left(\{z:r(x,z)>C\}\right)\geq 1/4\}\right)\\
  &\leq&4\hat\nu\left(\{(x,z):r(x,z)> C\}\right)<\delta
  \end{eqnarray*}
 and hence $$\nu\left(\{x\in X:\nu^{\cB_0}_x\left(\{z:r(x,z)\leq C\}\right)\geq 3/4\}\right)>1-\delta.$$

Let $K\subset \{x\in X:\nu^{\cB_0}_x\left(\{z:r(x,z)\leq C\}\right)\geq 3/4\}$ be compact with $\nu(K)\geq 1-\delta$ and remember that $C(\delta)=C^2$ and that $u_*(y,x)=(u_*(x,y))^{-1}$ and $u_*(x,z)u_*(z,y)=u_*(x,y)$.

Finally, if $x,y\in K$ and $y\in\cB_0[x]$, then
$\nu^{\cB_0}_x\left(\{z:r(x,z)\leq C\}\right)\geq 3/4$ and $\nu^{\cB_0}_y\left(\{z:r(y,z)\leq C\}\right)\geq 3/4$. Since $y\in\cB_0[x]$,
$\nu^{\cB_0}_x=\nu^{\cB_0}_y$ and so there is $z\in\cB_0[x]$ such that simultaneously $r(x,z)\leq C$ and $r(y,z)\leq C$ which gives $$\|u_*(x,y)\|\leq\|u_*(x,z)\|\|u_*(z,y)\|\leq r(x,z)r(y,z)\leq C^2=C(\delta).$$
\end{proof}




\section{\texorpdfstring{$Y$}{Y}-configurations}
	\label{sec:Y:configs}


\subsection{Definitions}
	\label{ssec:definitions_y_configurations}

\subsubsection{Setup}
	\label{sssec:setup_definitions_y_configurations}
For the rest of this subsection, all the points considered are on the finite cover $X\to Q$.

\begin{definition}[$Y$-configuration]
	\label{def:Y_configuration}
Suppose \index{$\alpha_3$}$\alpha_3>1$.
(The constant $\alpha_3$ will be fixed in \autoref{sec:subsec:choice:of:alpha3}).
Fix $q_1$, $u\in \cB_0(q_1)$.
Also fix parameters $\ell,\tau,t_0>0$ and $0<\tau\le \alpha_3 \ell$.
A \emph{$Y$-configuration}  \index{$Y(q_1,u,\ell,\tau,t_0)$}$Y(q_1, u, \ell, \tau,t_0)$ 
is the collection of points  
\begin{itemize}\item $q_1$ \item $q:= g_{-\ell} q_1$  \item $q_2 :=
  g_\tau u q_1$ \item $q_3 := g_{t_0} q_1$. 
\end{itemize}
\end{definition}

\begin{definition}[pre-$Y$-configuration]
A \emph{pre-$Y$}-configuration \index{$Y_{pre}(q_1,u,\ell,\tau)$}$Y_{pre}(q_1,u,\ell,\tau)$ is the
collection of points
\begin{itemize}\item $q_1$ \item $q:= g_{-\ell} q_1$ \item $q_2 := g_\tau u q_1$. 
\end{itemize}
\end{definition}
\subsubsection{Notation.}
If $Y$ is a $Y$-configuration or a pre-$Y$-configuration, we denote
its points by \index{$q(Y)$}$q(Y)$, \index{$q_1(Y)$}$q_1(Y)$, \index{$q_2(Y)$}$q_2(Y)$ and if applicable
\index{$q_3(Y)$}$q_3(Y)$. Similarly, we use the notation \index{$l$@$\ell(Y)$}$\ell(Y)$, \index{$u(Y)$}$u(Y)$, \index{$\tau(Y)$}$\tau(Y)$, \index{$t_0(Y)$}$t_0(Y)$
to retrieve the parameters defining $Y$. We also use the notations
\index{$q_{1/2}(Y)$}$q_{1/2}(Y)$ to denote $g_{-\ell(Y)/2} q_1(Y) = g_{\ell(Y)/2} q(Y)$,
and \index{$y_{1/2}(Y)$}$y_{1/2}(Y)$ to denote $g_{-\ell(Y)/2} u(Y)q_1(Y)$. 

\subsubsection{Notation.}
\label{sec:subsec:Yconf:notation}
If $Y$ is a $Y$-configuration, and there is no potential for
confusion, we will occasionally write $\index{$q$}q,\index{$q_{1/2}$}q_{1/2}, \index{$y_{1/2}$}y_{1/2},\index{$q_1$}q_1,\index{$q_2$}q_2,\index{$q_3$}q_3,\index{$l$@$\ell$}\ell,\index{$u$}u,\index{$\tau$}\tau,\index{$t_0$}t_0$ for
$q(Y),q_{1/2}(Y), y_{1/2}(Y), q_1(Y),q_2(Y),q_3(Y),\ell(Y),u(Y),\tau(Y),t_0(Y)$
respectively. Also if $Y'$ is a $Y$-configuration we will write
$\index{$q'$}q',\index{$q'_{1/2}$}q'_{1/2}, \index{$y'_{1/2}$}y'_{1/2}$, $\index{$q_1'$}q_1',\index{$q_2'$}q_2',\index{$q_3'$}q_3',\index{$l$@$\ell'$}\ell',\index{$u'$}u',\index{$\tau'$}\tau',\index{$t_0'$}t_0'$ for 
$q(Y'),q_{1/2}(Y'), y_{1/2}(Y')$, $q_1(Y')$, $q_2(Y')$, $q_3(Y')$,
$\ell(Y')$, $u(Y')$, $\tau(Y')$, $t_0(Y')$ respectively.

We use the same notational convention for pre-$Y$-configurations. 

\subsubsection{Warning}
	\label{warning:upper_bound_on_time}
We emphasize that in order for a configuration of points to be called
a $Y$ or pre-$Y$ configuration, the condition $\tau < \alpha_3 \ell$
must be satisfied.

\begin{definition}($Y \in K$)
	\label{def:Y_in_K}
Let $K \subset X$ be a subset, and let $Y$ be a $Y$-configuration or a
pre-$Y$-configuration. Using the notation
\S\ref{sec:subsec:Yconf:notation}, 
we say that \index{$Y \in K$}$Y \in K$ if $q \in K$, $q_1 \in K$, $q_{1/2} \in K$, $u
q_1 \in K$, $g_{-\ell/2} u q_1 \in K$, $q_2 \in K$, and (provided $Y$ is a
$Y$-configuration), $q_3\in K$. 
\end{definition}

\begin{definition}[bottom-linked $Y$-configurations]
	\label{def:bottom_linked_Y_configurations}
If $Y$ and $Y'$ are two $Y$-configurations or pre-$Y$ configurations,
we say (using the notation \S\ref{sec:subsec:Yconf:notation}) that $Y$
and $Y'$ are \emph{bottom-linked} if $\ell = \ell'$ and $q' \in
\gB_0^-[q] \subset \cW^s[q]$.
\end{definition}
Note that $\gB_0^-[q]$ is a subset of the finite cover $X$, so not only are the projections of $q,q'$ to $Q$ in the same stable manifold, but their $U^+$-groups are identified by $P^-(q,q')$.

\begin{definition}[right-linked $Y$-configurations]
	\label{def:right_linked_Y_configurations}
If $Y$ and $Y'$ are two bottom-linked $Y$-configurations, we say that $Y$ and $Y'$
are \emph{right-linked} if $t_0 = t_0'$. In that case $q_3'\in \cW^s[q_3]$. 
\end{definition}

\begin{definition}[$z$-points from $Y$-configurations]
	\label{def:z_points_from_y_configurations}
	Suppose $Y,Y'$ are two bottom-linked $Y$-configurations.
	Define the associated $z$-points:
	\[
		\index{$z_{1/2}$}z_{1/2}:=\cW^{cs}[y_{1/2}] \cap \cW^{u}[q'_{1/2}]
	\]
	Note that the intersection of the two manifolds gives a point on $Q$, and it is lifted to the finite cover using the flag $\frakF(q'_{1/2})$ (see \autoref{sec:subsec:stable:manifold:measurable:cover}).
\end{definition}
In \autoref{sssec:double_prime_points} we will define further points \index{$z''$}$z''$.


\begin{definition}[left-shadowing $Y$-configurations up to time $t$]
	\label{def:left_shadowing_Y_configuration_up_to_time_t}
Suppose $Y$ and $Y'$ are two bottom linked $Y$ or pre-$Y$
configurations, and suppose $0 \leq t \leq \tau(Y)$.
Using the notation
\autoref{sec:subsec:Yconf:notation}, 
we say that $Y'$ is \emph{$(k,\ve)$-left-shadows $Y$ up to time $t$} if for all
$0 \le s \le t$, we have that:
\[
	g_{\ell/2 + s} z_{1/2} \in \cW^{u}_{loc}[g_{s}u' q_1']
\]
where $\cW^u_{loc}$ denotes the local unstable (see \autoref{sssec:unstable_lyapunov_distance_and_charts}) lifted to $X$ (see \autoref{sec:subsec:stable:manifold:measurable:cover}).

In \autoref{sssec:choice_of_k_epsilon}, we will choose once and for all values of $k,\ve$ and only refer to \emph{left-shadowing} without reference to $k,\ve$.
\end{definition}
Note that, in view of \autoref{warning:upper_bound_on_time}, assuming that for $\delta>0$ and an appropriate compact $K$ of measure at least $1-\delta$, if $y_{1/2}\in K$ then:
\begin{align}
	\label{eqn:z_stays_in_Lyapunov_chart}
	\begin{split}
	d^{\cL}_{y_{1/2}}\left(g_{\ell/2 + s} \left(\sigma(z_{1/2})\right), \sigma(g_s u q_1) \right)\leq C(\delta)e^{-\alpha'\ell}\\
	d^{Q}\left(g_{\ell/2 + s} \left(\sigma(z_{1/2})\right), \sigma(g_s u q_1) \right)\leq C(\delta)e^{-\alpha'\ell}\\
	\text{ and therefore }
	g_{\ell/2 + s} \left(\sigma(z_{1/2})\right) \in \cL_{k,\ve}[\sigma(g_s u q_1)].
	\end{split}
\end{align}
for some $\alpha'>0$ (independent of $\delta$), where $\cL_{k,\ve}$ denotes the Lyapunov charts (see \autoref{sssec:lyapunov_charts_from_introduction}).



\begin{definition}[left-shadowing $Y$-configurations]
	\label{def:left_shadowing_Y_configuration}
Suppose $Y$ and $Y'$ are two bottom linked $Y$ or pre-$Y$
configurations, and suppose $0 \leq t \leq \tau(Y)$.
Using the notation
\S\ref{sec:subsec:Yconf:notation}, we say that $Y'$  \emph{left-shadows $Y$} if $Y'$ left-shadows $Y$ up to $\tau(Y)$ \textbf{and} $\tau(Y') =
\tau(Y)$. 
\end{definition}

 





\subsection{Interpolation map}
	\label{ssec:interpolation_map}
In this section, most constructions occur on $Q$, except for \autoref{cor:interpolation_and_generalized_subspaces} which involves the finite cover $X$.

\subsubsection{Setup}
	\label{sssec:setup_interpolation_map}
Let \index{$E$}$E$ be a smooth and natural cocycle (see \autoref{def:smooth_and_natural_cocycles}).
Let \index{$F$@$\cF^sE$}$\cF^sE$ denote the bundle of flags in $E$, with the same dimensions as the (forward, aka stable) Lyapunov flags.
At forward-regular points, we have a measurable equivariant choice of such flags \index{$E^{\leq \bullet}$}$E^{\leq \bullet}$.

We will free use notation and terminology on center-stable manifolds from \autoref{ssec:center_stable_manifolds}.
In particular, recall that we fixed once and for all a $C^\infty$ realization of the center-stable manifolds.

\begin{lemma}[Subexponential separation]
	\label{lem:subexponential_separation}
	For any $\delta>0$ there exists a compact set $K$ of measure at least $1-\delta$ such that for every $y\in K$ the center-stable manifold $\cW^{cs}[y]$ is defined, and furthermore there exists a decreasing function $\index{$r(N')$}r(N')>0$ with $r(N')\leq e^{-N'/2}$, defined for $N'\geq 1$ with the following properties.

	For a constant $\Lambda>2$ that only depends on the Lyapunov spectrum, and $C(N')>0$ a constant, and for all $ s,t\in \left[0, \frac{N'}{2\Lambda}\log \frac{1}{d^{Q}(y,z)}\right]$, we have the following.
	Suppose that $y\in K$ and $z\in \cW^{cs}[y]$ with $d^{Q}(y,z)\leq r(N')$, then for any $\ov{z}$ such that $d^Q(z,\ov{z})\leq d(y,z)^{N'}$ and furthermore $\ov{z}\in g_{-s}\cW^{cs}[g_s y]$
	we have that, setting $\ov{z}_t:=g_t \ov{z}$:
	\[
		\begin{split}
		d^Q(\ov{z}_t, g_t z) & \leq C(N') e^{\Lambda t}d^Q(y,z)^{N'}\\
		d^Q(g_t z, g_t y) & \leq C(N')
			e^{\tfrac{1}{N'}t}d^Q(y,z)
		\end{split}
	\]
	Furthermore, let \index{$\chi_{E,i}$}$\chi_{E,i}$ be the chart
        for $g_t y$ (see \autoref{eqn:smooth_natural_charts_covering_exp_q} and the discussion following it).
	Then we have analogous estimates for flags:\index{$d$@$\dist$}
	\[
	\begin{split}
		\dist\bigl(\chi_{E,i}(E^{\le \bullet}_{g_ty}(\ov{z}_t)), 
		\chi_{E,i}(g_{t}E^{\le \bullet}_{y}(z)
		\bigr)
		& \leq C(N') 
			e^{\Lambda t}d^Q(y,z)^{N'}
			\\
		\dist\bigl(
		\chi_{E,i}(g_{t}E^{\le \bullet}_{y}(z)),
		\chi_{E,i}\left(E^{\leq \bullet}(g_t y)\right)
		\bigr)
		& \leq C(N') 
			e^{\tfrac{1}{N'}t}d^Q(y,z)
	\end{split}
	\]
	where $\dist$ is defined in \autoref{sssec:summary_lyapunov_charts_and_distances}.
\end{lemma}
\noindent Note that at the maximal allowed value of $t$, the right-hand side of the all the above inequalities is bounded by $C(N')d^Q(y,z)^{\tfrac 12}$.
\begin{proof}
Fix  $k\in \N$ sufficiently large; as $k\to +\infty$ we will achieve arbitrarily large $N'$ in the statement.
Let $\cW^{cs}_{fake, k}[y]$ and  $z'\in \cW^{cs}_{fake, k}[y]\mapsto E^{\le \bullet}_{fake,k,y}(z')$ be the fake objects constructed in  \autoref{sssec:center_stable_manifolds_factorization}.  
Recall also that we fixed a $C^\infty$ realization of the center-stable manifold \index{$W$@$\cW^{cs}[y]$}$\cW^{cs}[y]$ and the forwards flag of $z\in \cW^{cs}[y] \mapsto E^{\le \bullet}_{y}(z)$\index{$E^{\le \bullet}_{y}(z)$} along $\cW^{cs}[y]$. 
There is a measurable function $y\mapsto R_r(y)$ such that if $y\in V_i$ (the $i$th trivialization chart for $E$) then for any $z\in  \cW^{cs}[y]$ there is $\bar z^k\in \cW^{cs}_{fake, k}[y]$ with the following properties:
 \begin{enumerate}
	\item We have
	\[
		d^Q(z,\bar z^k)\le \index{$R_r(y)$}R_r(y) d^Q(y,z)^{k}
	\]
	and analogously for the flags:
	\[
		\dist\bigl(\chi_{E,i}(E^{\le \bullet}_{fake,k,y}(\bar z^k)),\chi_{E,i}(E^{\le \bullet}_{y}(z))\bigr)\le R_r(y) d^Q(y,z)^{k}
	\]
	\item For constants \index{$\epsilon$@$\ve(k)$}$\ve(k)\to 0$ as $k\to +\infty$ we have: 
	\[
		d^Q(g_t y, g_t \bar z^k)\leq C_1(y)e^{\ve(k) t} d^Q(y,\bar{z}^k)
	\]
	and analogously for the flags:
	\[
		\dist\bigl(\chi_{E,i}(E^{\le \bullet}_{fake,k,y}(g_t\bar z^k)),\chi_{E,i}(E^{\le \bullet}_{y}(g_t z))\bigr)\leq
		C_1(y)e^{\ve(k)t} d^Q(y,\bar{z}^k)
	\]
\end{enumerate}
We may assume $K$ was chosen so that $R_r(y), C_1(y)$ above are uniformly bounded from above by some $C_2(k)$ and the distance $d^Q$ is uniformly comparable with the Lyapunov distance at $y$.

Now we estimate:
\begin{align*}
	d^Q(g_t y, g_t z) &  \leq d^Q(g_ty, g_t \bar z^k)  + d^Q(g_t \bar{z}^k, g_t z) \\
	& \leq C_2(k)\left[
	e^{\ve(k) t} d^Q(y,\bar z^k) 
	+ e^{\Lambda t}d^Q(\bar z^k,z)\right]\\
	& \leq C_2(k)\left[
	e^{\ve(k) t} d^Q(y,\bar z^k) 
	+ e^{\Lambda t}d^Q(y,z)^k\right]
\end{align*}
and this is the first inequality that we claimed.
Note that the exponent $\Lambda$ can be taken to be the maximal Lyapunov exponent of the flow (or the fixed scalar $2$, whichever is larger).

The estimate for flags is analogous.
\end{proof}

\subsubsection{Double-prime points}
	\label{sssec:double_prime_points}
Recall that there exist $C_1,\alpha_1>0$ such that $d^Q(y_{1/2},z_{1/2})\leq C_1 e^{-\alpha_1 \ell}$.
Choose $N'$ in \autoref{lem:subexponential_separation} such that $r(N')=C_1 e^{-\alpha_1\ell}$ (note that $r(N')\to 0$ and $N'\to +\infty$ as $\ell\to +\infty$).
Let $K_0$ be a compact set of measure at least $1-\tfrac{\delta}{10}$ on which the Oseledets theorem on $E$, its flags, etc. are uniform and transverse with constant $C_\delta$.
Using the Birkhoff ergodic theorem, let $K\subset K_0$ be a compact set of measure at least $1-\delta$ such that for any $y\in K$, whenever $\ell>\ell(\delta)$, we have that for some \index{$t''$}$t''\in \left[\tfrac{N'}{2}\ell,N'\ell\right]$ that $g_{t''}y\in K_0$.
Fix such a measurable choice of \index{$t''$}$t''(y,\ell)$ and set:
\[
	y'':=g_{t''}y_{1/2} \text{ and }z'':=g_{t''}z_{1/2}.
	\index{$y''$} \index{$z''$}
\]
With this notation, we have:
\begin{lemma}[$Y$-configurations in uniform Pesin sets]
	\label{lem:y_configurations_in_uniform_pesin_sets}
	Let $E\to Q$ be a smooth and natural cocycle. Fix
	\index{$\ve$}$\ve>0$ and \index{$N$}$N$.  Let
	\index{$\Lambda_E$}$\Lambda_E$ 
	denote the Lyapunov spectrum of $E$, and let
	$Q^{\pm,\Lambda_E}_{N,\ve}$ denote the uniform Pesin sets
	defined in
	\autoref{sssec:forward_pesin_sets_of_smooth_and_natural_cocycles}. 

	If $y_{1/2}\in Q^{\pm,\Lambda_E}_{N,\ve}$,  $q'_{1/2}\in Q^{\pm,\Lambda_E}_{N,\ve}$, and $y''\in Q^{-,\Lambda_E}_{N,\ve}$,  and if $\ell>0$ is sufficiently large, then  $z_{1/2}$ and $z''$ are in the uniform backwards Pesin set $Q^{-,\Lambda_E}_{3N,3\ve}$.
\end{lemma}
\begin{proof}
The backwards regularity of $z_{1/2}$ follows from that of $q_{1/2}'$, which is assumed, combined with the fact that $z_{1/2}\in \cW^{u}_{loc}[q_{1/2}']$ and furthermore $d^{Q}(z_{1/2},q_{1/2}')\leq C(\delta)e^{-\sigma \ell}$ where $\sigma>0$ depends only on the Lyapunov spectrum while $C(\delta)$ is a constant that only depends on the choice of Lusin set $K$ of measure at least $1-\delta$.

To establish the claim for $z''$, it suffices to do so assuming that $M= \ell/2+t''$ is integral.
We consider the sequence of transformations defined by 
\[
A_n' = 
\begin{cases} 
g_{-1}\vert_{E(g_{-n}y'')}  		& 0\le n\le M\\
P^{-}(q_{1/2},q'_{1/2})\circ P^{+}(y_{1/2},q'_{1/2})\circ g_{-1}\vert_{E(g_{-M}y'')}  		& n=M\\
g_{-1}\vert_{E(g_{-n+M}q'_{1/2})}  	& M+1\le n
\end{cases}
\]
Note that $g_{-M}y''=y_{1/2}$.
Define also
\[
	B_n' = 
	g_{-1}\vert_{E(g_{-n}z'')}  		\quad \text{for } 0\leq n
\]
Note that $A_n', B_n'$ are linear transformations between fibers of cocycles at different points.
Let now $A_n,B_n$ be the sequence of matrices obtained from $A_n',B_n'$ by applying the local trivializations $\chi_{E,i}$ associated to the smooth and natural cocycle $E$, see \autoref{sssec:local_trivializations_of_smooth_and_natural_cocycles}.
Concretely
\[
	A_n = \chi_{E,i_{n+1}}(x_{n+1})\circ A_{n}'\circ \chi_{E,i_n}^{-1}(x_n)
\]
where $x_n$ is the sequence of points determined from the definition of $A_n'$, and analogously $B_n$ is defined using a sequence $y_n:=g_{-n}z''$.
We will show that the sequence of transformations $A_n$ is $(\epsilon,N)$-forwards regular, and furthermore that $\norm{A_n-B_n}$ satisfies the assumptions of \autoref{theorem:shadowing:flags}, which then gives the desired claim.
 
The $(3\ve,3N)$-regularity of $A_n$ follows from the corresponding statement for $A_n'$, defining the backwards flag for $A_n'$ by pushing the backwards flags of $y''$ along $g_{-t}$ and $P^{-}(q_{1/2},q'_{1/2})\circ P^{+}(y_{1/2},q'_{1/2})$.
Now use that the local trivializations satisfy tempered bounds, see \autoref{sssec:local_trivializations_of_smooth_and_natural_cocycles}, and note that $P^+(y_{1/2},q'_{1/2})$ and $P^{-}(q_{1/2},q'_{1/2})$ preserve the Oseledets decompositions, and are furthermore bounded (in fact exponentially close to the identity if using the local trivializations and the fact that $q_{1/2},q'_{1/2},y_{1/2}$ belong to the Lusin set $K$).
\end{proof}

Consider now two bottom-linked $Y$-configurations.
We will omit the dependence of the points on the $Y$-configurations and follow the notation of \autoref{sec:subsec:Yconf:notation} and simply write $q_\bullet,q'_\bullet$ etc. for their points and also $z_{1/2}, y_{1/2}$, etc. as in \autoref{def:z_points_from_y_configurations}.

For the next statement, recall from \autoref{lem:y_configurations_in_uniform_pesin_sets} that $z_{1/2}$ is backwards-regular so it has a backwards flag $E^{\geq \bullet}(z_{1/2})$, and that $z_{1/2}$ and $y_{1/2}$ are center-stably related.

We also set \index{$F$@$\cF_{\ell,C_0}^s E(z_{1/2})$}$\cF^{s}_{{\ell,C_0}}\subset \cF^{s}E(z'')$ to consist of flags that are $C_0$-transverse to the unstable flags at $z''$.

\begin{proposition}[Fake forward flags on fake center-stables]
	\label{prop:fake_forward_flags_on_fake_center_stables}
Let $K,K_0$ be the compact sets as in \autoref{sssec:double_prime_points}.
Then
\begin{enumerate}
\item for any two flags $F_1,F_2\in \cF_{\ell,C_0}^s E(z_{1/2})$, we
have
\begin{displaymath}
\dist(g_tF_1, g_t F_2)\leq C_0e^{-{N'}\alpha_E \ell} \quad \forall t\in[0,\alpha_3\ell].
\end{displaymath}
where $\alpha_E>0$ and $\ell$ is larger than some $\ell(E)>0$ depending on the Lyapunov spectrum of $E$.
\item The choice of  flag of $ E^{\le \bullet}_{y_{1/2}}(z_{1/2})$ as defined in \autoref{def:fixed_realization_center_stable} is contained in $\cF_{\ell,C_0}^s E(z_{1/2})$.

\item
\label{item:estimate_on_fdd_trivialization}
For $F\in \cF_{\ell,C_0}^s E(z_{1/2})$ denote by \index{$E^{\lambda_i}_{F}(z_{1/2})$}$E^{\lambda_i}_{F}(z_{1/2})$ the components of the splitting of $E(z_{1/2})$ induced by the Oseledets filtration $E^{\geq\bullet}(z_{1/2})$ and the flag $E_{y_{1/2}}^{\leq \bullet}(z_{1/2})$.
Then for any $t\in [0, \alpha_3 \ell]$ and any $i$:
\begin{displaymath}
	\dist\left(\chi_{E,j}\left(E^{\lambda_i}(g_t y_{1/2})\right),
	g_t E^{\lambda_i}_{F}(z_{1/2})\right) \le e^{-\alpha_E \ell},  
\end{displaymath}
where $\chi_{E,j}$ is an ambient chart that contains $g_t y_{1/2}$, and $\alpha_E>0$ depends only on $E$ and the Lyapunov spectrum.

\item 
\label{item:z_hat_fake_ball_of_flags}
Suppose that for $\hat{y}_{1/2}:=g_{s/2}y_{1/2}$ with $s\in[0,\ell]$ we have $\hat{y}_{1/2}\in K$.
Set $\hat{z}_{1/2,s}:=\cW^{cs}[\hat{y}_{1/2}]\cap \cW^u[g_{s/2} q'_{1/2}]$ and $\hat{z}_{1/2}:=g_{-s}\hat{z}_{1/2,s}$.
Let also $\cF_{\ell-s,C_0}^s E(\hat{z}_{1/2,s})$ be the set set of flags at $\hat{z}_{1/2,s}$ obtained by the construction above, and set
\[
	\cF_{\ell-s,C_0}^s E(\hat{z}_{1/2}):= g_{-s}
	\left(\cF_{\ell-s,C_0}^s E(\hat{z}_{1/2,s})\right).
\]
Then we have the estimate
\[
	\dist\Big(
	\chi_{E,j}\left(\cF_{\ell,C_0}^s E({z}_{1/2})\right),
	\chi_{E,j}\left(\cF_{\ell-s,C_0}^s E(\hat{z}_{1/2})\right)
	\Big) \le e^{-N_2'\ell},  
\]
where $N_2'\to +\infty$ as $N_1'\to +\infty$.
\end{enumerate}

\end{proposition}
\begin{proof}
	To establish (i), we apply
        \autoref{lem:y_configurations_in_uniform_pesin_sets} to
        conclude that $z''\in Q^{-,\Lambda_E}_{3N,3\ve}$ (i.e. it is
        backwards regular).
	Therefore, $g_{-t}$ is exponentially contracting on the space of flags so applying $g_{-t''}$ we find that for any two $F_1,F_2\in \cF_{\ell,C_0}^s E(z_{1/2})$ we have
	\[
		\dist(F_1, F_2)\leq C_0e^{-{N'}\alpha_E' \ell}
	\]
	for some $\alpha_E'>0$ that only depends on the Lyapunov spectrum of $E$.
	Now we can apply $g_t$ for $t\in [0,\alpha_3 \ell]$ to this estimate and conclude part (i) for an appropriate choice of $\alpha_E$ and $\ell>\ell(E)$ to guarantee that $N'$ is sufficiently large.

For conclusion (ii), since $y''\in K_0$ and $z''$ is exponentially close to $y''$, by \autoref{lem:subexponential_separation} applied to $y'',z''$ we have that
 $\chi_{E,j}\left(g_{\ell/2+t''}E^{\le \bullet}_{y_{1/2}}(z_{1/2})\right)$
 is exponentially close to the flag $\chi_{E,j}\left(E^{\le \bullet}_{y''}[z'']\right)$, which is itself exponentially close to $\chi_{E,j}(E^{\le \bullet}[y''])$ since $y'',z''$ are exponentially close.
 
 Since $y''$ is biregular (with specified bounds) the flags $E^{\le \bullet}[y'']$ and $E^{\ge \bullet}[y'']$ are transverse.  
Also $\chi_{E,j}(E^{\ge \bullet}[ z''])$ and $\chi_{E,j}(E^{\ge \bullet}[y''])$ are exponentially close by shadowing in \autoref{lem:y_configurations_in_uniform_pesin_sets} combined with the \Holder estimates of \autoref{prop:holder_properties_of_fdd_cocycles}\autoref{item:Holder_fdd_cocycles}.
We conclude that the evolved stable flag $\chi_{E,j}(g_{\ell/2+t''}E^{\le \bullet}_{y_{1/2}}(z_{1/2}))$ is sufficiently transverse to the unstable flag  $\chi_{E,j}(E^{\ge \bullet}(z''))$.

For part (iii), observe that we proved in the previous part that if $F=E_{y_{1/2}}^{\leq \bullet}(z_{1/2})$ then the filtrations determining each of the flags are exponentially close, therefore so are the flags.
The same estimate holds for any other $F\in \cF_{\ell,C_0}^s E(z_{1/2})$ so the claim follows.

The estimates in (iv) are analogous to \autoref{lem:subexponential_separation} with $z:=z_{1/2}$ and $\ov{z}_t:=\hat{z}_{1/2}$, noting that $d^Q(y_{1/2},z_{1/2})^{N'_1}\leq e^{-N_2'\ell}$ with $N_2'\to +\infty$ as $N_1'\to +\infty$.
Then, for any sufficiently large $k\in \bN$ we consider the part-time equivariant fake flags $E^{\leq\bullet}_{y_{1/2},k,fake}(z_{1/2,k})$, where $z_{1/2,k}:=\cW^{cs}_{k,fake}[y_{1/2}]\cap \cW^{u}[q'_{1/2}]$. and verify that for each of $z_{1/2},\hat{z}_{1/2}$ the distances to the $z_{1/2,k}$-flags satisfy the required bounds.
\end{proof}





For the next statement, recall that forward dynamically defined (fdd) cocycles are introduced in \autoref{def:forward_dynamically_defined_cocycle}; we will use the analogous backwards dynamically defined (bdd) cocycles and make use of the fact that $L\cW^u$ is bdd, see \autoref{prop:linearization_cocycle_is_fdd}.

\begin{proposition}[Fake future flag and subresonant map]
\label{prop:fake_future_flag_and_subresonant_map}
	Consider the cocycle \index{$L\cW^u$}$L\cW^u$ used for the linearization of subresonant normal forms for the unstable manifolds (see \autoref{thm:linearization_of_stable_dynamics_single_diffeo}), which is backwards dynamically defined (bdd) and in particular smooth along unstables.
	Let \index{$E$}$E$ be an associated smooth and natural cocycle, with bdd subcocycles $\index{$S_1$}S_1\subset \index{$S_2$}S_2\subseteq E$ such that $L \cW^u\isom S_2/S_1$ (see \autoref{prop:subquotient_representation_of_fdd_cocycles} for their existence).
	We can then choose $C_0$ sufficiently large and apply \autoref{prop:fake_forward_flags_on_fake_center_stables} to have the following properties.

	There exists a measurable choice of flag \index{$F$@$\tilde{F}$}$\tilde{F}\in
        \cF_{C_0,\ell}^s E(z_{1/2})$ (where $\cF_{C_0,\ell}^s
        E(z_{1/2})$ is provided by
        \autoref{prop:fake_forward_flags_on_fake_center_stables})
    that furthermore induces a flag \index{$F$@$\tilde{F}^{sr}$}$\tilde{F}^{sr}$ on $S_2/S_1(z_{1/2})\isom L\cW^u(z_{1/2})$, with the following properties:
	\begin{enumerate}
		\item There exists a linear map\index{$\Phi$@$\tilde{\Phi}_1$}
		\[
			\tilde{\Phi}_1\colon L\cW^u(y_{1/2})\to L\cW^u(z_{1/2})
		\]
		such that $\tilde{F}^{sr}=\tilde{\Phi}_1(L\cW^u(y_{1/2})^{\leq \bullet})$
		i.e. the flag $\tilde{F}^{sr}$ is the image of the stable flag $L\cW^u(y_{1/2})^{\leq \bullet}$ at $y_{1/2}$.
		\item There exists a subresonant map\index{$\phi$@$\tilde{\phi}_1$}
		\[
			\tilde{\phi}_1\colon \cW^u[y_{1/2}]\to \cW^u[z_{1/2}]
		\]
		such that its linearization is $\tilde{\Phi}_1$.
	\end{enumerate}
\end{proposition}

\begin{proof}[Proof of \autoref{prop:fake_future_flag_and_subresonant_map}]
  It suffices to perform the constructions at the points $y'',z''$, and then apply $g_{-t''}$ to obtain them at $y_{1/2},z_{1/2}$.
  As a preliminary, note that for some
  constant $C_1$ depending on the compact set $K$, and $\alpha''>0$
  depending on the Lyapunov spectrum, we have the bound
  $d^Q(y'',z'')\leq C_1 e^{-\alpha'' \ell}$.
  Furthermore, $z''$ is in a backwards-uniform Pesin set by \autoref{lem:y_configurations_in_uniform_pesin_sets}, i.e. there exist $N_0\in \bN,\ve_0>0$, which depend on the compact set $K_0$ in \autoref{prop:fake_forward_flags_on_fake_center_stables} such that $z''\in Q^{-,\Lambda_E}_{N_0,\ve_0}$ in the sense of \autoref{sssec:forward_pesin_sets_of_smooth_and_natural_cocycles}.

We will now construct a subresonant map
  $\phi''\colon \cW^{u}[y'']\to \cW^u[z'']$.
Let then $I_{y'',z''}^{L\cW^u}\colon L\cW^{u}(y'')\to L\cW^{u}(z'')$ be the
identification provided by \autoref{eqn:trivializations_fdd}.
By \autoref{lem:common_parametrization_of_linearizations}, it follows that 
there exists $h\in \GL\left(L\cW^u(y'')\right)$ such that
\begin{align*}
		h\cdot
        \left(I^{L\cW^u}_{y'',z''}\right)^{-1}
        \left(L_{z''}(\cW^u[z''])\right) & =
        L_{y''}\left(\cW^u[y'']\right) && \text{on linearized unstables}\\
        h\cdot
        \left(I^{L\cW^u}_{y'',z''}\right)^{-1}
        \left(L\cW^u(z'')^{\geq \bullet}\right)
        &=L\cW^u(y'')^{\geq \bullet} && \text{on flags}
\end{align*}
and
\[
\norm{h-\id}\leq C(\delta) e^{-\alpha_1'' \ell}
\]
where \index{$L_x$}$L_x\colon \cW^u[x]\to L\cW^u(x)$ is the linearization map,
and \index{$L\cW^u(x)^{\geq \bullet}$}$L\cW^u(x)^{\geq \bullet}$ is the unstable flag of $L\cW^u(x)$, and $\alpha_1''>0$ is some constant that only depends on the Lyapunov spectrum.
Let $\Phi'' = I_{y'',z''} \circ h^{-1}$.
Then, $\Phi''
(L_{y''}(\cW^u[y'']) = L_{z''}(\cW^u[z''])$,
$\Phi''(L\cW^u(y'')^{\geq\bullet}) = L\cW^u(z'')^{\geq\bullet}$,
and therefore, by
\autoref{cor:polynomial_maps_between_different_subresonant_manifolds},
$\Phi''$ is a linearization of a subresonant map $\phi''\colon \cW^u[y'']
\to \cW^u[z'']$.


Using the notation of \autoref{prop:fake_forward_flags_on_fake_center_stables}, it now suffices to construct $\phi''_E\colon E(y'')\to E(z'')$ with the following properties.
It preserves the filtration
  $S_1\subset S_2\subset E$ on both spaces, and it takes the stable
  flag $E(y'')^{\leq \bullet}$ into the the set of flags in
  $\cF^sE(z'')$ which are $C_{0}$-transverse to the unstable flag
  $E(z'')^{\geq \bullet}$, where $C_0$ is the same absolute constant
  as in \autoref{prop:fake_forward_flags_on_fake_center_stables}.
  Furthermore $\phi''_E$ is required to induce the subresonant
  polynomial map $\phi''\colon \cW^u[y'']\to \cW^u[z'']$ just constructed.
  Indeed,
  given such a $\phi''_E$, we can apply $g_{-t''}$ (where $t''$ is
  such that $y'' = g_{t''} y_{1/2}$) to the flag
  $\phi''_E\left(E(y'')^{\leq \bullet}\right)$ to obtain $\tilde{F}$
  at $z_{1/2}$.
  By
  \autoref{prop:fake_forward_flags_on_fake_center_stables} the flag
  $\tilde{F}$ is in $\cF^s_{C_0,\ell} E(z_{1/2})$.
  Because the
  $g_t$-dynamics is by subresonant maps on unstables, the resulting
  flag $\tilde{F}^{sr}$ on $L\cW^u(z_{1/2})$ is also induced by a
  subresonant map
  $\tilde{\phi} = g_{-t''} \circ \phi'' \circ g_{t''}$.
  So all the
  desired conclusions are satisfied.
 
  We now construct $\phi''_E$.
  First, recall that by \autoref{prop:holder_properties_of_fdd_cocycles}, the bundles $S_i$ are \Holder-continuous on the uniform backwards Pesin set $ Q^{-,\Lambda_E}_{N_0,\ve_0}$ to which $y'',z''$ belong.
  So there exists $C_2,\beta_2$ such that in the smooth local trivialization of $E$ at $y''$ we have
  \[
  	d_{Gr}(S_i(z''),S_i(y''))\leq 
  	C_2 d(z'',y'')^{\beta_2} \leq
  	C_1 C_2^{\beta_2} e^{-\beta_2\alpha_1''\ell}
  \]
   where $d_{Gr}$ is the distance in the Grassmannian.
   Let us also note that the identifications $I_{y'',z''}^{L\cW^u}$ used in the construction of $\phi''$ are built in \autoref{eqn:trivializations_fdd} using the same trivialization of $E$, the maps $S_i$ into the Grassmannian, and orthogonal projection.
   Let also $I_{y'',z''}^{E}$ be the associated identification between the fibers of $E$ at $y'',z''$.
   Then our desired $\phi''_{E}$ is $\wtilde{h}\circ I_{y'',z''}^E$, where $\wtilde{h}$ is the lift of $h\in \GL(L\cW^u(y''))$ to $E(y'')$ under the identification $L\cW^u(y'')\isom S_2(y'')\cap S_1^{\perp}(y'')$, and extension as the identity on the orthogonal complement.
   Then the desired $C_0$-transversality of the flags $\phi''_E\left(E^{\leq \bullet}(y'')\right)$ and $E^{\geq \bullet}(z'')$ follows from the boundedness of the trivialization of $E$ in the chart at $y''$, and the \Holder regularity of the unstable flags of bdd cocycles restricted to uniform Pesin sets, see \autoref{prop:holder_properties_of_fdd_cocycles}.
\end{proof}

\begin{theorem}[Fake holonomies and strictly subresonant maps]
\label{thm:fake_holonomies_and_ssr_map}
	\leavevmode
	\begin{enumerate}
		\item
		\label{item:fake_decomposition}
		There exists a decomposition $L\cW^u(z_{1/2})=\oplus L\cW^u(z_{1/2})^{\lambda_i}$ such that taking associated flags $L\cW^u(z_{1/2})^{\geq \lambda_i}$ gives the backward/unstable Lyapunov filtration of $L\cW^u$, and taking the flags $L\cW^u(z_{1/2})^{\leq \lambda_i}$ gives the filtration $\tilde{F}^{sr}$ from \autoref{prop:fake_future_flag_and_subresonant_map}.
		\item
		\label{item:fake_Pplus}

		There exists a map:\index{$P^+_{fake}(q'_{1/2},z_{1/2})$}
		\[
			P^+_{fake}(q'_{1/2},z_{1/2})\colon 
			L\cW^u(q_{1/2}')\to
			L\cW^u(z_{1/2})
		\]
		taking the Lyapunov decomposition on the left side to the decomposition on the right provided by \autoref{item:fake_decomposition}, and on the associated graded coincides with the standard unstable holonomy.
		Furthermore $P^+_{fake}$ is induced by a strictly subresonant map\index{$P$@$\wp_{fake}^+$}
		\[
			\cW^u[q'_{1/2}]\xrightarrow{\wp_{fake}^+} \cW^u[z_{1/2}]
		\]
		(note that these are the same submanifold) taking $q'_{1/2}$ to $z_{1/2}$.
		\item \label{item:P+_fake_z_zhat_comparison}
		Suppose that for some $N'_1>0$, we have a $\hat{z}_{1/2}\in \cW^u[q'_{1/2}]$ with $d^Q(z_{1/2},\hat{z}_{1/2})\leq e^{-N_1'\ell}$.
		Suppose also that we have an analogous decomposition $\oplus L\cW^u(\hat{z}_{1/2})^{\lambda_i}$ with
		\[
			d^{Gr}_{q'_{1/2}}\left(P_{GM}(z_{1/2},q'_{1/2})L\cW^{u}(z_{1/2})^{\lambda_i},
			P_{GM}(\hat{z}_{1/2},q'_{1/2}) L\cW^u(\hat{z}_{1/2})^{\lambda_i}
			\right)
			\leq e^{-N_1'\ell}\quad \forall i
		\]
		where $P_{GM}$ denotes the Gauss--Manin connection on $L\cW^u$ and $d^{Gr}_{q'_{1/2}}$ denotes the distance on Grassmannian at $q'_{1/2}$.
		Then the analogous map $P^+_{fake}(q'_{1/2},\hat{z}_{1/2})$, after post-composing with the Gauss--Manin connection, satisfies the bound:
		\[
			\norm{
			P^+_{fake}(q'_{1/2},{z}_{1/2})
			\circ
			P^+_{fake}(q'_{1/2},\hat{z}_{1/2})^{-1}
			}_{L\cW^u(q'_{1/2})}
			\leq e^{-N_2'\ell}
		\]
		where the map inside the norm is conjugated by the Gauss--Manin connection to be viewed as a map of $L\cW^u(q'_{1/2})$, and with $N_2'\to +\infty$ as $N_1'\to +\infty$.
	\end{enumerate}
\end{theorem}

\noindent Recall that we defined strictly subresonant maps of an unstable manifold in \autoref{sssec:setup_subgroups_compatible_with_the_measure}.
Recall also from \autoref{sssec:holonomies_and_relating_different_groups} that there exists a ``Gauss--Manin'' connection $P_{GM}$ on $L\cW^u$ for points on the same unstable by \autoref{thm:linearization_of_stable_dynamics_single_diffeo} (which lists its properties).
To show that $\wp_{fake}^+$ in the above proposition is strictly subresonant we must check that the linear map $P_{GM}(q'_{1/2},z_{1/2})\circ P^+_{fake}(q'_{1/2},z_{1/2})$ on $L\cW^u(q'_{1/2})$ is induced by a strictly subresonant map in the sense of \autoref{def:groups_of_resonant_and_strictly_subresonant_maps}.

\begin{proof} Note that the point $q'_{1/2}$ is biregular. 
	The map $P^+_{fake}(q'_{1/2},z_{1/2})$ is unique and well-defined given its requirements to map decompositions in \autoref{item:fake_Pplus} and to agree with the standard holonomy on graded pieces.
	It suffices to show that it is the linearization of \emph{some} subresonant map, since then the map $P_{GM}(q'_{1/2},z_{1/2})\circ P^+_{fake}(q'_{1/2},z_{1/2})$ will be strictly subresonant by \autoref{def:groups_of_resonant_and_strictly_subresonant_maps} and its action on the graded Lyapunov pieces.

	To exhibit the required subresonant map
	consider the composed subresonant maps
	\[
		\wp_{tmp}\colon 
		\cW^u[q_{1/2}']
		\xrightarrow{\wp^-(q_{1/2}',q_{1/2})}
		\cW^u[q_{1/2}]
		\xrightarrow{\wp^+(q_{1/2},y_{1/2})}
		\cW^u[y_{1/2}]
		\xrightarrow{\tilde{\phi}_1}
		\cW^u[z_{1/2}]
	\]
	where $\tilde{\phi}_1$ is provided by \autoref{prop:fake_future_flag_and_subresonant_map}, and \index{$P$@$\wp^{\pm}$}$\wp^{\pm}$ by \autoref{prop:subresonant_maps_from_holonomies}.
	By construction each subresonant maps admits a linearization which furthermore preserves the decompositions of $L\cW^u$ at the source and target.
	So the resulting linearization preserves the decompositions, but it might not act as the standard holonomy on each graded piece.

	Next we use the semidirect product structure of subresonant maps provided by \autoref{prop:group_structure_on_subresonant_maps}\autoref{item:sresn_gp_semidir}, using the Lyapunov splitting at $q_{1/2}'$ to define resonant maps.
	We thus have $\wp_{tmp} = \wp_{tmp}^{ssr} \cdot \wp_{tmp}^{r}$ where $\wp_{tmp}^{ssr}$ is strictly subresonant and $\wp_{tmp}^r$ is resonant.
	The linearization of $\wp_{tmp}^{r}$ preserves the Lyapunov splitting at $q_{1/2}'$, while $\wp_{tmp}^{ssr}$ acts as the identity (i.e. the Gauss--Manin connection) on the associated graded pieces since it is strictly subresonant.
	Since on associated graded pieces the Gauss--Manin connection and the standard holonomies agree, it follows that the map $\wp_{tmp}^{ssr}$ has all the required properties of $\wp^{+}_{fake}(q_{1/2}',z_{1/2})$.

	To prove the last part, we note that the corresponding bounds hold for the distances between the map $\wtilde{\phi}_1$ associated to $z_{1/2}$ and its analogue $\hat{\phi}_1$ associated to $\hat{z}_{1/2}$, from the assumption on the decompositions of $L\cW^u$.
	Indeed, the bounds on the decomposition follow from \autoref{prop:fake_forward_flags_on_fake_center_stables}\autoref{item:z_hat_fake_ball_of_flags}, and the construction of $\wtilde{\phi}_1$ from \autoref{prop:fake_future_flag_and_subresonant_map}.
	Following through the above construction of $\wp^+_{fake}$ yields the bound on the distance between the $P^+_{fake}$-maps associated to $z_{1/2}$ and $\hat{z}_{1/2}$.
\end{proof}

\begin{figure}[htbp!]
	\centering
	\includegraphics[width=0.8\linewidth]{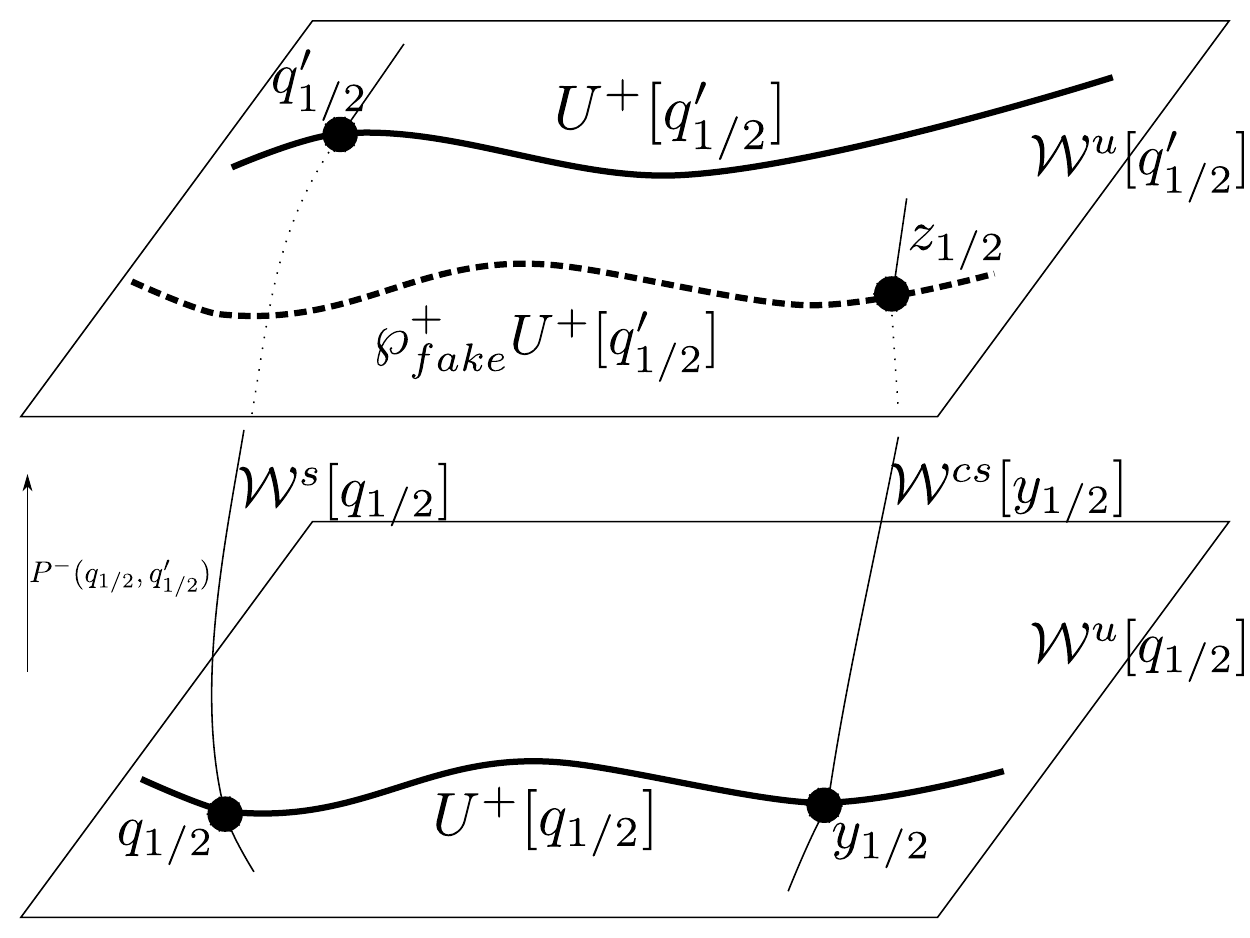}
	\caption{Interpolating the curves.}
	\label{fig:interpolating_curves}
\end{figure}

\subsubsection{Defining the interpolation map}
	\label{sssec:defining_the_interpolation_map}
We are now ready to define two key maps that appear in later arguments.
First define the map \index{$\phi$@$\tilde{\phi}$}$\tilde{\phi}$ by the composition:
\begin{equation}
	\label{eqn_cd:phi_tilde_floppy}
\begin{tikzcd}[column sep = large]
	\cW^u[y_{1/2}]
	\arrow[r, "{\wp^{+}(y_{1/2},q_{1/2})}"]
	\arrow[rrr, "{\tilde{\phi}}", bend right  = 15]
	& 
	\cW^u[q_{1/2}]
	\arrow[r, "{\wp^-(q_{1/2},q_{1/2}')}"]
	&
	\cW^u[q_{1/2}']
	\arrow[r, "{\wp^+_{fake}(q_{1/2}',z_{1/2})}"]
	&
	\cW^u[z_{1/2}]
\end{tikzcd}
\end{equation}
Here \index{$P$@$\wp^+(y_{1/2},q_{1/2})$}$\wp^+(y_{1/2},q_{1/2})$ and \index{$P$@$\wp^-(q_{1/2},q_{1/2}')$}$\wp^-(q_{1/2},q_{1/2}')$ are defined in \autoref{prop:subresonant_maps_from_holonomies}.

The linearization of each arrow then induces a linear map \index{$\Phi$@$\tilde{\Phi}$}$\tilde{\Phi}= L\td \phi$ by the composition:
\begin{equation}
	\label{eqn_cd:Phi_tilde_linearization}
\begin{tikzcd}[column sep = large]
	L\cW^u(y_{1/2})
	\arrow[r, "{P^{+}(y_{1/2},q_{1/2})}"]
	\arrow[rrr, "{\tilde{\Phi}}", bend right  = 15]
	& 
	L\cW^u(q_{1/2})
	\arrow[r, "{P^-(q_{1/2},q_{1/2}')}"]
	&
	L\cW^u(q_{1/2}')
	\arrow[r, "{P^+_{fake}(q_{1/2}',z_{1/2})}"]
	&
	L\cW^u(z_{1/2})	
\end{tikzcd}
\end{equation}
Note that $\tilde{\Phi}$ and $\tilde{\phi}$ are distinct from the maps $\tilde{\Phi}_1,\tilde{\phi}_1$ from \autoref{prop:fake_future_flag_and_subresonant_map}.

\subsubsection{Interpolation map and generalized subspaces}
	\label{sssec:interpolation_map_and_generalized_subspaces}
The construction of the interpolation maps $\tilde{\Phi},\tilde{\phi}$ can be performed without passing to the finite cover $X\to Q$, in other words the map only depends on the images of the points under the map $\sigma$.
Note that the map $\wp^{+}_{fake}$ is also defined only in terms of data on $Q$.

Suppose now that $U^+(q)\subset \bbG^{ssr}(q)$ is a family of generalized subspaces, defined for $q\in X$, as in \autoref{def:compatible_family_of_subgroups}.
Now the measurable connection maps $P^{\pm}(\bullet,\bullet)$ are
defined on $X$, and furthermore these maps have the property
$P^\pm(x,y)U^+(x)P^{\pm}(y,x)=U^+(y)$ for $x,y\in X$ in a full measure set, see
\autoref{sssec:measurable_connections_on_finite_covers}.

\begin{proposition}[Interpolation and generalized subspaces]
	\label{cor:interpolation_and_generalized_subspaces}
	Suppose that $U^+(q)\subset \bbG^{ssr}(q)$ is a family of groups compatible with the measures, and let $U^+[q]\subset \cW^u[q]$ be the corresponding family of orbits, for all $q$ in a full measure set.
	Assume that $y_{1/2}\in U^+[q_{1/2}]$.

	Then $U^+[y_{1/2}]=U^+[q_{1/2}]$ and the subresonant map $\tilde{\phi}\colon \cW^u[y_{1/2}]\to \cW^u[z_{1/2}]$ taking $y_{1/2}$ to $z_{1/2}$ defined by \eqref{eqn_cd:phi_tilde_floppy} is such that the image $\tilde{\phi}\left(U^+[y_{1/2}]\right)$ is the image of $U^+[q_{1/2}']$ under the strictly subresonant map $\wp_{fake}^+$.
\end{proposition}
\begin{proof}
	The map $\tilde{\phi}$ is obtained at the composition of 
	$\wp^{+}(y_{1/2},q_{1/2}),
	\wp^-(q_{1/2},q_{1/2}')$ and
	$\wp^+_{fake}(q_{1/2}',z_{1/2})$.
	The maps on cocycles $P^{+}(y_{1/2},q_{1/2})$ and $P^-(q_{1/2},q_{1/2}')$ induce maps on Lie algebras and preserve the family of $g_t$-equivariant subalgebras $\fraku^+$, by the Ledrappier invariance principle.
	It follows that
    $P^+(y_{1/2},q_{1/2})U^+(y_{1/2})P^{+}(q_{1/2},y_{1/2})=U^+(q_{1/2})$ and since by assumption $y_{1/2}\in U^+[q_{1/2}]$, it follows that in fact $U^+[y_{1/2}]= U^+[q_{1/2}]$.

	Similarly $\wp^-(q_{1/2},q_{1/2}')(U^+[q_{1/2}])=U^+[q_{1/2}']$, and the final conclusion follows since $\wp^+_{fake}(q_{1/2}',z_{1/2})$ is strictly subresonant by \autoref{thm:fake_holonomies_and_ssr_map}.
\end{proof}

\subsubsection{Properties needed for subsequent arguments}
	\label{sssec:properties_needed_for_subsequent_arguments}
We collect the properties that we derived in this section and that will be used in the next ones:
\begin{enumerate}
	\item There exists a decomposition $L\cW^u(z_{1/2})=\oplus L\cW^{u}(z_{1/2})^{\lambda_i}$, such that the associated unstable flag is compatible with the one from $q'_{1/2}$, while the ``fake'' stable flag is provided by \autoref{prop:fake_future_flag_and_subresonant_map}.
	\item The interpolation map $\wtilde{\phi}\colon \cW^{u}[y_{1/2}]\to \cW^{u}[z_{1/2}]$ is such that its linearization takes the Lyapunov decomposition at $y_{1/2}$ to the decomposition at $z_{1/2}$.

	\item Furthermore, $\wtilde{\phi}$ satisfies \autoref{cor:interpolation_and_generalized_subspaces}.
\end{enumerate}
Note that the interpolation map and its construction does not depend on the compatible family of subgroups $U^+$.


\subsection{Estimates on evolved interpolation map}
	\label{ssec:evolved_interp}

Let $Y$, $Y'$ be two bottom-linked $Y$ configurations. Using the
notation \S\ref{sec:subsec:Yconf:notation}, let $y_{1/2} = g_{-\ell/2}
u q_1$, and let $z_{1/2} = \cW^{cs}[y_{1/2}] \cap
\cW^{u}[q'_{1/2}]$. 

Recall the map $\tilde{\phi}\colon \cW^u[y_{1/2}]\to \cW^u[z_{1/2}]$  given in \eqref{eqn_cd:phi_tilde_floppy} and its linearization 
 $\tilde{\Phi}\colon L\cW^u(y_{1/2})\to 	L\cW^u(z_{1/2})$
 given by \eqref{eqn_cd:Phi_tilde_linearization} in
 \ref{sssec:defining_the_interpolation_map}.
 Let \index{$z$}$z = g_{\ell/2}
 z_{1/2}$, so that $z \in \cW^u[q_1']$ and it is also exponentially (in $\ell$) close to $\cW^{cs}[uq_1]$.
 
 Given $0<\tau<\alpha_3 \ell$, define \index{$\phi_\tau$}${\phi_\tau}\colon \cW^u[g_\tau uq_1]\to \cW^u[g_\tau z]$ by 
\begin{equation}
  {\phi_\tau}:= g_{\tau+\ell/2} \circ \tilde{\phi}\circ g_{-(\tau+\ell/2)}.
\end{equation}
We remark that $  {\phi_\tau}$ is defined using $Y$ and $q'\in \cW^s[q]$.  In the future, we may denote $\phi_\tau$ by\index{$\phi_\tau(Y,q')$} 
\begin{equation}\label{eq:phitau}
{\phi_\tau}(Y,q')\colon \cW^u[g_\tau uq_1]\to \cW^u[g_\tau z].\end{equation}

\begin{theorem}[Interpolation does not move points much]
\label{theorem:interpolation:does:not:move:points}
There exists $0<\alpha<1$ such that for all sufficiently large $k$ and sufficiently small $\ve>0$ the following holds.
For every $\delta > 0$ there exists a
compact set $K$ with $\nu(K) > 1-\delta$ such that the following
holds:
Suppose $Y = Y_{pre}(q_1,u,\ell,\tau)$ and $Y' =
Y_{pre}(q_1',u',\ell,\tau)$ are two bottom-linked
pre-$Y$-configurations, and suppose that $Y'$ $(k,\ve)$-left-shadows $Y$. 
Suppose that, using the notation of
\autoref{sec:subsec:Yconf:notation}, $q$, $q_{1/2}$, $q_1$, $u q_1$, $g_{-\ell/2} u q_1$,
$q'$, $q_{1/2}'$, $q_1'$ all belong to $K$.  
Note that $q_2$ is not necessarily assumed in $K$.

Then, for all $x\in \cW^u_{loc}[g_\tau uq_1]$, 
\begin{align*}
		d^{Q}(x, \phi_\tau(x))& \leq C(\delta) e^{-\alpha\ell}.
\end{align*}
\end{theorem}
\noindent Recall that $\cW^u_{loc}[q]\subset \cW^u[q]$ is the unit ball around $q$ in the $d^u$-distance.

The proof of \autoref{theorem:interpolation:does:not:move:points}
relies on \autoref{lemma:L:is:bilipshitz} and \autoref{lem:common_parametrization_of_linearizations} and will be given in \autoref{sssec:proof_of_theorem:interpolation:does:not:move:points}, after we establish the following preliminary results.

\begin{proposition}
	\label{prop:distortion_of_phi}
There exists $\beta_{2},\beta_{1/2}>0$ such that for all sufficiently large $k$ and sufficiently small $\ve>0$ the following holds.
For every $\delta > 0$ there exists a
compact set $K$ with $\nu(K) > 1-\delta$ such that the following
holds:
Suppose $Y = Y_{pre}(q_1,u,\ell,\tau)$ and $Y' =
Y_{pre}(q_1',u',\ell,\tau)$ are two bottom-linked
pre-$Y$-configurations, and suppose that $Y'$ $(k,\ve)$-left-shadows $Y$. 
Suppose, using the notation of
\autoref{sec:subsec:Yconf:notation} and \autoref{def:Y_in_K} that $Y\in K$
and also $q'$, $q_{1/2}'$, $q_1'$ all belong to $K$.

Let $z_{1/2}$ be the unique intersection point of
$\cW^{cs}_{loc}[y_{1/2}]$ and $\cW^u_{loc}[q'_{1/2}]$, and let \index{$z_2$}$z_2 =
g_{\tau + \ell/2} z_{1/2}$. 
Then we have the estimates:
\begin{align}
	\norm{ L\wtilde{\phi} -  I^{L \cW^u}_{q_{1/2},z_{1/2}}}
	& \leq H_{1/2} (\delta) e^{-\beta_{1/2} \ell} \nonumber \\
	\norm{\gr_\bullet L\wtilde{\phi} - \gr_\bullet I^{L \cW^u}_{q_{1/2},z_{1/2}}}
	& \leq H_{1/2} (\delta) e^{-\beta_{1/2} \ell} \nonumber \\
	\label{eqn:grLphi_ILWu_bound}
	\norm{\gr_\bullet L\phi_{\tau} (q_2,z_2)- \gr_\bullet I^{L \cW^u}_{q_2,z_2}}
	& \leq H_2 (\delta) e^{-\beta_{2} \ell}
\end{align}
\end{proposition}
\noindent Note that by \autoref{lem:y_configurations_in_uniform_pesin_sets} the point $z_{1/2}$ is in a uniform Pesin set, as is $q_{1/2}$ by assumption.

\begin{proof}
We recall some notation first: $d^{\cL}_q$ denotes the distance in the
$(k,\ve)$-Lyapunov chart centered at $q$, and $\alpha$ is the minimal rate of expansion along unstables.
We reduce first the second estimate to the first one.

We have $d^{\cL}_{y_{1/2}}(z_{1/2}, y_{1/2})\le e^{-\alpha \ell/2}$ whence for all
$0\le t\le \alpha_3\ell$, 
\[
d^{\cL}_{g_{t+\ell/2}y_{1/2}}(g_{t+\ell/2} z_{1/2}, g_{t+\ell/2} y_{1/2})\le e^{2\ve(t+\ell/2) - \alpha \ell/2}.
\]
Similarly, for all $0\le s \le t+\ell/2$ we have 
\[
	d^{\cL}_{g_{-s}q_2}(g_{-s}z_2,g_{-s}q_2)
	\le
	e^{2\ve(t+\ell/2-s) - \alpha\ell/2}\le e^{2\ve(\alpha_3 \ell+\ell/2) - \alpha\ell/2}\le 
e^{- \alpha\ell/4}
\]
by taking $\ve$ sufficiently small.
With $\gamma=\frac{\alpha}{8(\alpha_3+ \frac 1 2)}$ we have
\[
e^{- \alpha\ell/4}\le e^{- \alpha\ell/8}e^{-\gamma s}
\]
for all $0\le s \le t+\ell/2\le (\alpha_3 + 1/2) \ell.$
We then apply \autoref{distortionest} with the cocycle $g_{s}^{L \cW^u}$ starting at $q_{1/2}, z_{1/2}$ for time $s\in [0, t+\ell/2]$, we find that, for $\ve_1>0$ arbitrarily small, if $k\geq k_0(E),\ve\in (0,\ve_0(E,\ve_1))$ we have that
\[
	\norm{\gr_\bullet L\phi_{\tau} (q_2,z_2)- \gr_\bullet I^{L \cW^u}_{q_2,
  z_2}}\le 
  e^{\ve_1 (t+\ell/2)} 
  \norm{\gr_\bullet L\wtilde{\phi} - \gr_\bullet I^{L \cW^u}_{q_{1/2},z_{1/2}}}
\]
It now suffices to prove the second estimate and take $\ve_1>0$
sufficiently small.
Note that the second estimate follows from the first, so we have reduced to the first inequality.

Recall that $E$ is the smooth and natural cocycle from \autoref{prop:fake_future_flag_and_subresonant_map} associated to $L\cW^u$.
Let $\wtilde{F}\in \cF^{s}_{C_0,\ell}\left(z_{1/2}\right)$ be the flag constructed in that proposition, which is in particular also adapted to $L\cW^u$.
Then \autoref{prop:fake_forward_flags_on_fake_center_stables}\autoref{item:estimate_on_fdd_trivialization} with $t=0$ implies that (everything is for the ambient trivialization $\chi_{E}$ of $E$):
\begin{align}
	\label{eqn:y_12_z12_flag_estimate}
	\dist\left(\chi_{E,j}\left(E^{\lambda_i}\left(y_{1/2}\right),
	E^{\lambda_i}_{\wtilde{F}}\left(z_{1/2}\right)\right)\right)
	\leq e^{-\alpha_E \ell}.
\end{align}

We can ensure that for the compact set $K$, there exists $r(\delta,E)>0$ and $C(\delta,E,k)\geq 1$ such that for $q\in K$ the trivializations $\chi_{E,j}$ and $\exp_q^{E}$ are defined on $Q[q;r(\delta,E)]$ (see \autoref{sssec:notation_for_fdd_cocycles}) and the $C^k$-norm of $\chi_{E,j}^{-1}\circ \exp_q^{E}$ is bounded by $C_0:=C(\delta,E,k)$ on $B(q;r(\delta,E))$.

Recall now that $y_{1/2},q_{1/2},q'_{1/2}\in K$.
We can assume that $\ell$ is large enough such that the following hold:
\begin{itemize}
	\item $d^Q$ between any of the points involved is less than $r(\delta,E)$.
	\item The trivializations $\exp_q^E$ and $I^{E,L\cW^u}$ (see \autoref{eqn:trivializations_E_adapdated_F}) are uniformly \Holder-comparable with a constant $C_1\geq 1$ and exponent $\beta(E,L\cW^u)$.
\end{itemize}
We can therefore pass freely (at a cost of a uniform constant) between trivializations at $y_{1/2},q_{1/2},q'_{1/2}$ and $\chi_{E,j}$, so all estimates stated below are for $\chi_{E,j}$.
We will add the subscript $(E,j)$ to any map expressed in this trivialization.

It follows from \autoref{cor:holder_properties_of_measurable_holonomy} that exists $\gamma_0>0$ such that the measurable connection maps satisfy
\begin{align*}
	\norm{P^+_{E,j}(y_{1/2},q_{1/2}) - \id} & \leq C_0 e^{-\gamma_0 \ell}\\
	\norm{P^{-}_{E,j}(q_{1/2},q'_{1/2}) - \id} & \leq C_0 e^{-\gamma_0 \ell}
\end{align*}
and additionally, using \autoref{eqn:y_12_z12_flag_estimate}, we also have that
\[
	\norm{P^{+}_{fake,E,j}(q'_{1/2},z_{1/2}) - \id}  \leq C_0 e^{-\gamma_0 \ell}.
\]
Combining these estimates we find that
\[
	\norm{\wtilde{\Phi}_{E,j} - \id}\leq C_2 e^{-\gamma_0\ell}
\]
where $\wtilde{\Phi}_{E,j}$ is defined on $E$ analogously to the definition in \autoref{eqn_cd:Phi_tilde_linearization} of $\wtilde{\Phi}$ on $L\cW^u$.

Since the trivializations $\exp_q^E$ and $I^{E,L\cW^u}_{q,-}$ (see \autoref{eqn:trivializations_E_adapdated_F}) are \Holder-comparable on $K$ by assumption we conclude that
\[
	\norm{\wtilde{\Phi}_{E} - I^{E,L\cW^u}_{y_{1/2},z_{1/2}}}\leq C_3 e^{-\gamma_0\ell}
\]
where $\wtilde{\Phi}_E\colon E(y_{1/2})\to E(z_{1/2})$ is the map induced by $\wtilde{\Phi}_{E,j}$.
We finally note that $\wtilde{\Phi}_E$ induces $\wtilde{\Phi}=L\wtilde{\phi}$ on $L\cW^u$, and the induced map by $I^{E,L\cW^u}$ on $L\cW^u$ satisfies
\[
	\norm{I^{E,L\cW^u}_{y_{1/2},z_{1/2}}-  I^{L\cW^u}_{y_{1/2},z_{1/2}} }
	\leq C_4 e^{-\gamma_1\ell}
\]
as maps $L\cW^u(y_{1/2})\to L\cW^u\left(z_{1/2}\right)$ and hence the desired claim follows.
\end{proof}

\begin{lemma}[$L\phi_\tau$ and $I$ are close]
	\label{lemma:Lphi:close:to:I}
There exists $0<\alpha<1$ and for every $\delta > 0$ there exists a
compact set $K$ with $\nu(K) > 1-\delta$ such that the following
holds:
Suppose $Y = Y_{pre}(q_1,u,\ell,\tau)$ and $Y' =
Y_{pre}(q_1',u',\ell,\tau)$ are two bottom-linked
pre-$Y$-configurations, and suppose that $Y'$ $(k,\ve)$-left-shadows $Y$. 
Suppose, using the notation of
\autoref{sec:subsec:Yconf:notation} and \autoref{def:Y_in_K} that $Y\in K$
and also $q'$, $q_{1/2}'$, $q_1'$ all belong to $K$.

Then
\begin{displaymath}
\norm{\left(I^{L\cW^u}_{g_\tau u q_1, g_\tau z}\right)^{-1} \circ L \phi_\tau - \mathbf{1} } \le C(\delta)e^{-\alpha \ell}.
\end{displaymath}
\end{lemma}

\begin{proof}
The inequality on $\gr_{\bullet}L\cW^u$ follows from \autoref{eqn:grLphi_ILWu_bound}, so to conclude the same on $L\cW^u$ we need to estimate the distance between the corresponding decompositions.

Set \index{$y$}$y:=uq_1$ and consider the Lyapunov splitting $L\cW^u(g_\tau y)=\oplus
L\cW^u(g_\tau y)^{\lambda_i}$, and the decomposition 
$L\cW^u(g_\tau z)=\oplus
L\cW^u(g_\tau z)^{\lambda_i}$, where $L\cW^u(g_\tau z)^{\lambda_i} =
g_{\tau + \ell/2} L\cW^u(z_{1/2})^{\lambda_i}$ and
$L\cW^u(z_{1/2})^{\lambda_i}$ is as in 
\autoref{thm:fake_holonomies_and_ssr_map}. 
The map $L\phi_\tau$ intertwines these decompositions. 
The map $I^{L\cW^u}_{g_\tau u q_1, g_\tau z}$ intertwines the
associated backwards flags $L\cW^u(g_\tau y)^{\ge \lambda_i}$ and
$L\cW^u(g_\tau z)^{\ge \lambda_i}$ (since these give the backwards
Lyapunov filtration).

Recall that the map $I^{L\cW^u}_{x,y}$ is defined in \autoref{eqn:trivializations_fdd} using a smooth and natural cocycle $E$ with bdd subcocycles $S_1\subseteq S_2\subseteq E$ such that $L\cW^u\isom S_2/S_1$.
Recall also that the splitting $\oplus L\cW^{u}(z_{1/2})^{\lambda_i}$ is constructed in \autoref{prop:fake_future_flag_and_subresonant_map} from a flag $\wtilde{F}$ at $z_{1/2}$.
Now \autoref{prop:fake_forward_flags_on_fake_center_stables}\autoref{item:estimate_on_fdd_trivialization} shows that, in some ambient chart $\chi_{E,j}$ containing $y$, the splittings $\chi_{E,j}\left(E^{\lambda_i}(g_\tau y)\right)$ and $g_\tau E^{\lambda_i}_{\wtilde{F}}(z)$ are exponentially close, with exponent depending only on $E,L\cW^u$ and constant depending on $\delta$.

We also know from \autoref{prop:holder_properties_of_fdd_cocycles} that $S_1\left(g_\tau y\right)$ and $S_2\left(g_{\tau}z\right)$ are exponentially close when expressed in the chart $\chi_{E,j}$, with exponent depending only on $E,L\cW^u$ and constant depending on $\delta$.
We recall next that the trivializations $\chi_{E,j}$ and $I^{E,L\cW^u}$ of $E$ are \Holder-equivalent (see \autoref{eqn:trivializations_E_adapdated_F}) and the trivialization $I^{L\cW^u}$ and the one induced by $I^{E,L\cW^u}$ on $L\cW^u$ are \Holder-equivalent as well (see the compatibility claim in \autoref{prop:lyapunov_adapted_trivializations}).

We now conclude that $I_{g_\tau y,g_{\tau}z}^{L\cW^u}\left(L\cW^u(g_\tau y)^{{\lambda_i}}\right)$ and $g_\tau L\cW^u_{\wtilde{F}}(z)^{{\lambda_i}}$ are exponentially close and \autoref{eqn:grLphi_ILWu_bound} completes the proof.
\end{proof}

\subsubsection{Proof of \autoref{theorem:interpolation:does:not:move:points}}
	\label{sssec:proof_of_theorem:interpolation:does:not:move:points}
Let us introduce the notation \index{$x''$}$x'' = \phi_\tau(x)$, and for this proof only $y_\tau:=g_\tau u q_1$, and $z_{\tau}:=g_{\tau}z$.

Given any compact set $K_1$ of measure $1-\delta_1$ with $\delta_1>0$ sufficiently small (depending on $\delta$ and $k$), we first reduce to the case when $y_{\tau}\in K_1$.
Given an arbitrary $y_{\tau}$, and since $y_{0}\in K$ and we can assume by Birkhoff that a point in $K$ visits $K_1$ for a $(1-\tfrac 12 \delta_1)$-fraction of its time, we conclude that for some $\tau_1\in[\tau(1-\tfrac 35\delta_1),(1-\tfrac 1 {10}\delta_1)\tau]$ we have $y_{\tau_1}\in K_1$.
The estimates in \autoref{rmk:basic_c_k_bounds} then imply that the distances gets distorted by a factor of $A_ke^{\delta_1\kappa_k \tau}$ for some $A_k(\delta)$.
Now it suffices to pick $\delta_1$ sufficiently small depending on $k$.
We require, in addition, that for $y\in K_1$ the Lyapunov distance $d^{\cL}_y$ is bi-Lipschitz with $d^Q$ up to a factor of $C_1(\delta_1)$, that for $z\in \cW^u_{loc}[y]$ the quantity $d^u(y,z)$ is bi-Lipschitz with $d^Q(y,z)$ up to a factor of $C_1(\delta_1)$, and that $\ell\geq\ell_0$ such that $e^{-\tfrac 1 {10}\delta_1 (\lambda_+-\ve)\ell_0}C_1(\delta_1)^2\leq 1$ (where $\lambda_+$ is the smallest strictly positive Lyapunov exponents).
With this assumption we have that $g_{-\delta_1 \ell/10}\left(\cW^{u}_{loc}[g_{-\delta_1 \ell/10}y]\right)\subset \cL_{k,\ve}[y]$.
Note that we can always increase the constant $C(\delta)$ to account for the case when $\ell,\tau$ are smaller than the required lower bounds.

To proceed to the proof, recall that $L_q\colon \cW^u[q]\into L\cW^u(q)$ denotes the linearization map based at a point $q$, and note that 
\begin{displaymath}
	L_{z_{\tau}} x'' =  (L\phi_\tau) (L_{y_{\tau}} x).
\end{displaymath}
Using \autoref{lem:common_parametrization_of_linearizations} to obtain $h$ below, let $x'\in \cW^{u}_{loc}[y_{\tau}]$ be the unique point such that
\[
	L_{y_{\tau}}x' = h\cdot \left(I^{L\cW^u}_{y_{\tau},z_{\tau}}\right)^{-1}(L_{z_{\tau}}x'') \text{ inside }L\cW^{u}(y_{\tau})
\]
By the triangle inequality
\[
	d^{\cL}_{y_{\tau}}(x,x'')\leq d^\cL_{y_\tau}(x,x') + d^\cL_{y_{\tau}}(x',x'')
\]
it suffices to estimate each summand on the right.

For the first summand, note that by \autoref{lem:common_parametrization_of_linearizations} there exists $h\in \GL(L\cW^u[y_{\tau}])$ and that we have
\begin{align*}
	\tfrac{1}{\kappa_0} d^\cL_{y_{\tau}}(x,x') & \leq \norm{L_{y_{\tau}}x - L_{y_{\tau}}x'}
	= \norm{
	L_{y_{\tau}}x - h\cdot \left(I^{L\cW^u}_{y_{\tau},z_{\tau}}\right)^{-1}(L_{z_{\tau}}x'')
	}\\
	& = \norm{
	L_{y_{\tau}}x - h\cdot \left(I^{L\cW^u}_{y_{\tau},z_{\tau}}\right)^{-1}(L\phi_\tau) (L_{y_{\tau}} x)
	}
\end{align*}
For the first inequality we used \autoref{lemma:L:is:bilipshitz}.
But now $L_{y_{\tau}}x$ is a vector of uniformly bounded norm in $L\cW^u(y_{\tau})$ since $x\in \cW^u_{loc}[y_{\tau}]$, the norm of $\id-h$ is bounded using \autoref{lem:common_parametrization_of_linearizations} by $\kappa_1e^{-\alpha_1 \ell}$, the norm of $\left(I^{L\cW^u}_{y_{\tau},z_{\tau}}\right)^{-1}(L\phi_\tau)-\id$ is bounded using \autoref{lemma:Lphi:close:to:I} by $e^{-\alpha_2 \ell}$.
It follows that the last expression above is bounded by $\kappa_3e^{-\alpha_2' \ell}$ for appropriate $\kappa_3,\alpha_2'$ where $\alpha_2'$ depends only on the Lyapunov spectrum and $\kappa_3$ on $\delta$.

To estimate the second term, i.e. $d^\cL_{y_{\tau}}(x',x'')$, we will use \autoref{prop:measuring_distances_in_linearization} with $y=y_{\tau}$ and $z=g_\tau z$:
\[
	d^{\cL}_{y_{\tau}}(x'',x') 
	\leq
	\kappa_4\left(
	\norm{I_{y_{\tau},g_\tau z}^{-1}\circ L_{g_\tau z}(x'') - L_{y_{\tau}}(x')}
	+ d^{\cL}_{y_{\tau}}(g_\tau z,y_{\tau})
	\right)^{\beta_4}
\]
Now the second term in the parentheses is bounded by $\kappa_5 e^{-\alpha_5 l}$ for appropriate $\kappa_5,\alpha_5$.
Using the definition of $x'$ at the start of the proof, we can rewrite the first term in the parentheses as
\[
	\norm{I_{y_{\tau},g_\tau z}^{-1}\circ L_{g_\tau z}(x'') - h\cdot I^{-1}_{y_{\tau},z_{\tau}}\circ L_{g_\tau z}(x'')}
\]
Now the norm of $L_{z_{\tau}}(x'')$ is bounded by some $\kappa_6$
since $x''\in \cW^u_{loc}[z_{\tau}]$, the norm of
$I^{-1}_{y_{\tau},z_{\tau}}$ is bounded by $\kappa_7$, and finally the
norm of $h-\id$ is bounded by $\kappa_1e^{-\alpha_1 \ell}$ as before.
The conclusion follows.
\hfill \qed

\subsubsection{Choice of $(k,\ve)$}
	\label{sssec:choice_of_k_epsilon}
We now fix \index{$k$}$k\in \bN$ sufficiently large with the following properties.
First, we require that \autoref{theorem:interpolation:does:not:move:points} and \autoref{prop:distortion_of_phi} holds.
We additionally require that for the constant $\alpha_3$ from \autoref{sec:subsec:choice:of:alpha3}, we have that \index{$\epsilon$@$\ve$}$\ve\alpha_3$ is less than any of the positive constants $\alpha$ arising in the preceding sections, as well as \autoref{sec:the_curly_a_operators} through \autoref{sec:extra_invariance}.

From now on, left-shadowing of $Y$-diagrams means $(k,\ve)$-left shadowing in the sense of \autoref{def:left_shadowing_Y_configuration_up_to_time_t}.

\section{Factorization}
	\label{sec:factorization}
The main goal of this section is \autoref{theorem:3varA:haltime:factorizable}. 
Recall the $Y$-configuration notation from \autoref{sec:Y:configs}.
Let $Y$ and $Y'$ be two bottom linked pre-$Y$ configurations. We use
the notation \autoref{sec:subsec:Yconf:notation}. A function $Z(Y,Y')$
is called admissible if it depends only on $q$, $u$, $\ell$, $q'$, and
not on $u'$, $\tau$ and $\tau'$. Thus, if $Z$ is an admissible
function, we may write $Z(q,q',u,\ell)$ for $Z(Y,Y')$.








\subsection{Basic constructions}
	\label{ssec:basic_constructions_factorization}

\begin{definition}[half-time factorizable]
	\label{def:half_time_factorizable}
	Let \index{$T$}$T\to Q$ be a measurable a.e. defined vector bundle.
	Suppose $q \in Q$, $\ell >
        0$, $u \in \cB_0(g_\ell q)$,  and $q' \in \cW^s_{loc}[q]$.  
	A measurable admissible function $Z(q,q',u,\ell)$ with values in $T(q_{1/2})$ is \emph{factorizable of order \index{$N$}$N>0$} if there exist:
	\begin{itemize}
		\item Cocycles $\index{$V^s$}V^s,\index{$V^u$}V^u$, smooth along stables, resp. unstables, and admitting stable resp. unstable holonomies, equipped with $g_t$-equivariant measurable sections $\index{$\iota_s$}\iota_s,\index{$\iota_u$}\iota_u$.

		Define $g_t$-equivariant embeddings \index{$F_q^{s}$}$F_q^{s}\colon \cW^{s}[q]\into V^{s}(q)$ (depending measurably on $q$) by 
		\[
			F_q^{s}(p)=H^{s}(p,q)\iota_{s}(p) \text{ where }\index{$H^s$}H^s \text{ is the stable holonomy of }V^s.
		\]
		Define the analogous construction for points on $\cW^u[q]$ and \index{$F_q^u$}$F_q^u$.
		\item A measurable family of bilinear maps\index{$A$@${\crA}(q_{1/2})$}
		\[
			{\crA}(q_{1/2})\colon V^s(q_{1/2})\times V^u(q_{1/2})\to T(q_{1/2}).
		\]
	\end{itemize}
such that for any $\delta>0$ there exists a compact set $K$ of measure at least $1-\delta$, as well as constants $\ell_0=\ell_0(K,N)$ and $C=C(K,N)$, such that provided $q,q',q_{1/2},q_{1/2}',y_{1/2},q_1,uq_1$ all belong to $K$, we have:
	\[
		\norm{
		Z(q,q',y,\ell) - 
		{\crA}(q_{1/2}) 
		\left( g_{\ell/2} F_q^s(q'), 
		g_{-\ell/2} F_{q_1}^u(uq_1) \right) 
		}
		\leq
		C e^{-N\ell}
	\]
	provided that $\ell>\ell_0$.

	When $Z(q,q',u, \ell)$ is factorizable to any order $N>0$, we will say that $Z$ is factorizable.
\end{definition}
With our assumptions and notation as above we have
\begin{align*}
	g_{\ell/2} F_q^s(q') & = F_{q_{1/2}}^s(q_{1/2}')\\
	g_{-\ell/2} F_{q_1}^u(uq_1) & = F_{q_{1/2}}^u(y_{1/2})
\end{align*}
so we can rewrite the main inequality above as
\[
	\norm{
		Z(q,q',y,\ell) - 
		{\crA}(q_{1/2}) 
		\left( 
		F_{q_{1/2}}^s(q_{1/2}'), 
		F_{q_{1/2}}^u(y_{1/2}
		\right) 
		}
		\leq
		C e^{-N\ell}
\]

\begin{figure}[htbp!]
	\centering
	\includegraphics[width=0.5\linewidth]{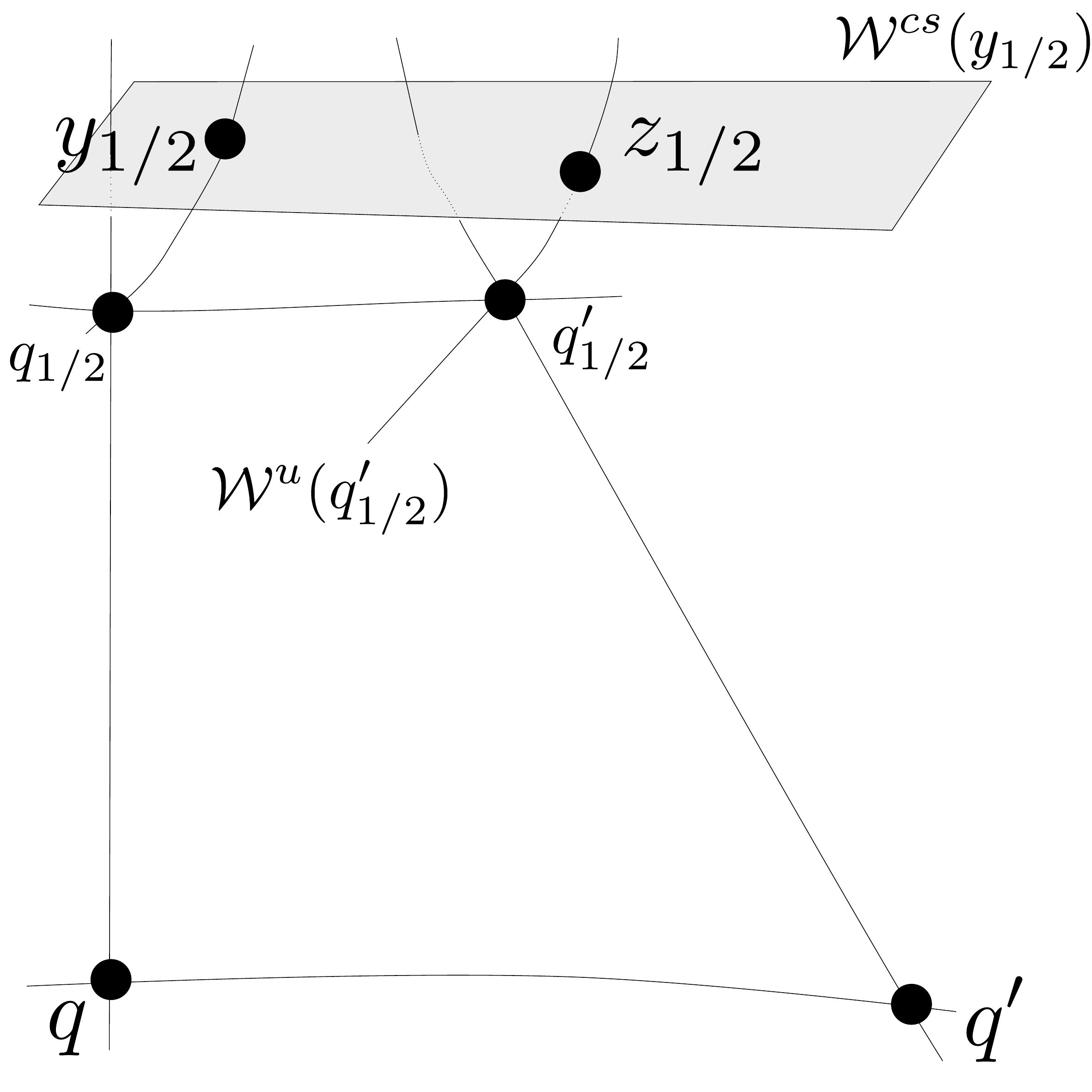}
	\caption{The initial, and halfway points.}
	\label{fig:halfway_diagram}
\end{figure}

\begin{remark}[On factorizability]
	\label{rmk:on_factorizability}
	\leavevmode
	\begin{enumerate}
		\item The target bundle $T\to Q$ is not necessarily a cocycle over $g_t$, although in our intended application (\autoref{theorem:3varA:haltime:factorizable}) it will be.
		However, at intermediate stages we will also work with some $T$ that don't have a cocycle structure.
		\item We will call \index{$N$}$N$ the order of approximation of the factorizable map $Z$ and $V^{s/u}, T$ the associated bundles.
		When the dependence on $N$ is relevant, we write \index{$V^{s/u,(N)}$}$V^{s/u,(N)}$ for the bundle; note that $T$ is independent of $N$.
		Similarly, the maps $F$ and $\crA$ can depend on $N$, and when this dependence is relevant we write it as \index{$F^{(N)}$}$F^{(N)}$ and \index{$A$@$\crA^{(N)}$}$\crA^{(N)}$.
		\item When we can choose $V^{s/u},\iota_s,\iota_u$ and
                  $\crA$ independently of $N$, to make $Z=\crA
                  \circ F$ equal on the nose, we will say that $Z$ is
                  \emph{exactly factorizable}.

		\item In most statements and proofs below, $N$ will be fixed and the discussion will hold for every choice of $N$.
		In other situations,
		we will apply the assumptions of factorizability of $Z$ with a larger $N$ to conclude the factorizability of another function $Z'$ for a smaller $N'$, but with $N'$ going to infinity as $N$ goes to infinity.

		\item We will also make use of the following notation:
		\begin{align}
			\label{eq:Zqtilde}
			\index{$Z$@$\wtilde{Z}(q_{1/2})$}\wtilde{Z}(q_{1/2}):=\lim_{\ell\to +\infty}Z(g_{-\ell/2}q_{1/2},g_{-\ell/2}q_{1/2},0,\ell) \in T(q_{1/2})
		\end{align}
		Let us see that this in fact makes sense.
		Because of the equivariance of the sections $\iota_s$ and $\iota_u$ we have
		\[
			\crA(q_{1/2})
			\big(g_{\ell/2}\iota_s(q),g_{-\ell/2}\iota_u(q_1)\big)=
			\crA(q_{1/2})\big(\iota_s(q_{1/2}),\iota_u(q_{1/2})\big)
		\]
		so the expression that approximates $Z(q,q,0,\ell)$ is
                in fact independent of $\ell$.
		We can send $\ell$ to $+\infty$ and define alternatively
		\[
			\wtilde{Z}(q_{1/2}):=
			\crA(q_{1/2})\big(\iota_s(q_{1/2}),\iota_u(q_{1/2})\big)
		\]
		Note that the expression on the right-hand side, although a priori dependent on the bundles $V^{s/u}$ and degree of approximation $N$, must in fact be independent of these choices since it is the limit of $Z(q,q,0,\ell)$ as $l\to +\infty$.
	\end{enumerate}
\end{remark}

\begin{proposition}[Contraction in $V^s$]
\label{prop:contraction_in_v}
	Let $Z$ be factorizable and fix a degree of approximation $N$.
	Let $V^s$ be a bundle associated to a factorizable function as in \autoref{def:half_time_factorizable} and let $F_q^{s}\colon \cW^{s}[q]\into V^{s}(q)$ be the $g_t$-equivariant embedding given by the equivariant section $\iota_{s}$ of $V^{s}$ (analogous statements will hold for the unstable objects).
	Then for any $\delta>0$ we can choose the compact set $K$ of measure at least $1-\delta$ such that there exist constants with the following additional properties.
	\begin{enumerate}
		\item \label{item:holonomy_invariance_of_iota_s}
		For $q,q'\in K$ and $q'\in \cW^s[q]$, the difference $\iota_s(q)-H^s(q',q)\iota_s(q')\in V^s(q)$ is in the negative Lyapunov subspaces of $V^s$, while the section $\iota_s(q)$ is in the zero Lyapunov subspace of $V^s$.
		
		The analogous property holds for $\cW^u[q]$, $\iota_u$ and $V^u$ except that the difference is in the positive Lyapunov subspaces.
		\item There exists $\beta>0$ depending on $Z$ and $N$, and a constant $C$ depending on $Z,N$ and $K$, such that we have the estimates
		\[
			\norm{F^{s}_{q_{1/2}}(q_{1/2})-F^{s}_{q_{1/2}}\left(q_{1/2}'\right)}\leq Ce^{-\beta \ell}
		\]
		where the norm is in the fiber $V^s\left(q_{1/2}\right)$.

		The analogous estimate for $F^u_{q_{1/2}}(q_{1/2})-F^u_{q_{1/2}}(y_{1/2})$ holds.
		\item For $\wtilde{Z}(q_{1/2})$ defined in \autoref{eq:Zqtilde} we have the estimate
		\[
			\norm{\wtilde{Z}(q_{1/2})-Z(q,q',u,\ell)}\leq C e^{-\delta_0 \ell}
		\]
		for some $\delta_0>0$ that only depends on $Z$.
	\end{enumerate}
\end{proposition}
Note that the degree of approximation $N$ does not appear in part (iii) of the proposition above.
\begin{proof}
	Part (ii) follows from (i) and the equivariance of the sections.
	The proofs for the stable and unstable are analogous, so we
        only cover the stable case. In this proof only, we write $V$
        for $V^s$. 

	For (i) we can take $-\delta$ to be half the largest strictly negative Lyapunov exponent (so, closest to $0$) of $V$ and further restrict to a compact set where the Oseledets theorem holds uniformly.
	It is clear that $\iota_s(q)$ is in the zero Lyapunov subspace since it is equivariant, and since the stable holonomy $H^s_V$ preserves the Lyapunov filtration, it suffices to show that the projection of $\iota_s(q)-H^s_V(q',q)\iota_s(q')$ to $\gr_0^sV:=V^{\leq 0}/V^{<0}$ vanishes.
	But if we consider the projection $\iota_0(q)$ of $\iota_s(q)$ to $\gr^s_0 V$, the induced holonomy map $H^s_V$, which agrees with the holonomy on graded bundles, will preserve $\iota_0$ by the Ledrappier Invariance Principle \autoref{thm:ledrappier_invariance_principle}.
	The conclusion then follows.

	For (iii), pick \emph{some} degree of approximation $N>0$ and apply (ii) on a compact set of almost full measure where $\norm{\crA}$ is bounded by some constant.

\end{proof}

\begin{proposition}[Polynomials and factorizability]
	\label{prop:polynomials_and_factorizability}
	\leavevmode
	\begin{enumerate}
	\item 
	Suppose that $Z(q,q',u,\ell)$ satisfies the definition of factorizability in \autoref{def:half_time_factorizable} with bundles $V^{s/u,(N)},T$, but $\crA^{(N)}$ is a polynomial map (and such data is available for every $N$).
	Then $Z$ is factorizable, i.e. we can replace the bundles $V^{s/u,(N)}$ by different ones to satisfy the definition of factorizability with $\crA^{(N')}$ bilinear.
	\item 
	If a function $Z(q,q',u,\ell)$ is factorizable and $f\colon T\to T'$ is a measurable map of bundles, which is fiberwise polynomial, then $f(Z)$ is also factorizable.
	\end{enumerate}
\end{proposition}

\begin{proof}
	A polynomial map between finite-dimensional vector spaces $W_1\xrightarrow{p} W_2$ can be factored as
	\[
		W_1 \xrightarrow{\Sym^{\leq d}}\Sym^{\leq d}W_1 \xrightarrow{Lp} W_2
	\]
	where the degree of $p$ is at most $d$, and $Lp$ is a linear map, while the (Veronese) embedding $\Sym^{\leq d}$ of $W_1$ is not linear.
	
	For (i), we replace $F_q$ by $\Sym^{\leq d}\circ F_q$ and thus replace $V^{(N)}$ by the bundle $\Sym^{\leq d}V^{(N)}$.
	The compact sets $K$ are taken so that the degree of $\crA^{(N)}$ is bounded by $d$.

	Part (ii) follows from (i) since $\crA_f^{(N)}:=f\circ \crA^{(N)}$ is now a polynomial map and can be used to factorize $f(Z)$.
\end{proof}

Let us derive a useful consequence.

\begin{corollary}[Assembly of factorizable maps]
	\label{cor:assembly_of_factorizable_maps}
	Suppose that $Z_1,\ldots, Z_k$ are factorizable, with target bundles $T_1,\ldots, T_k$ and set $Z=Z_1\oplus\cdots Z_k$.
	\begin{enumerate}
		\item Suppose $f\colon T_1\times \cdots \times T_k\to T'$ is a measurable, fiberwise polynomial map.
		Then $f(Z)$ is also factorizable.
		\item Suppose $f\colon T_1\times \cdots \times T_k\setminus S\to T'$ is a measurable, fiberwise smooth map, where $S$ is an equivariant measurable subset that restricts to a closed set in a.e. fiber, and such that $\wtilde{Z}(q_{1/2})$ is almost surely not in $S$.
		Then the $f(Z)$ is also factorizable.

		\item 
		\label{item:smooth_postcomposition_assembly}
		Suppose $F\subset T_1\times\cdots \times T_k$ is a measurable subset that restricts to a smooth submanifold in a.e. fiber, and that furthermore almost surely $Z(q,q',u,l)\in F\left(q_{1/2}\right)$.
		Suppose also that $f\colon F\to T$ is a measurable, fiberwise smooth function.
		Then $f(Z)$ is factorizable.
	\end{enumerate}
\end{corollary}
An instance of (i) in the above corollary is when $Z_1$ and $Z_2$ are valued in $T_1=\Hom(E,T')$ and $T_2=\Hom(T',T'')$ respectively, and $f(Z_1,Z_2)\in \Hom(E,E'')$ is their composition.
\begin{proof}
	Part (i) is an immediate consequence of \autoref{prop:polynomials_and_factorizability}.
	For part (ii) we will use the estimates in \autoref{prop:contraction_in_v}.
	Fix a desired degree of approximation $N$ as in \autoref{def:half_time_factorizable} for $f(Z)$.
	Let also $V^s,V^u$ etc. be provided by the definition of factorizability of $Z$ to degree $M$.

	We restrict to a compact set $K$ of almost full measure such that there exists $\ve_0>0$ such that $\dist\left(\wtilde{Z}(q_{1/2}),S(q_{1/2})\right)\geq \ve_0$.
	For any given $D\in \bN$, we can further choose $K$ (still of almost full measure) to assume that the $C^D$-norm of $f$ on the fibers over $q_{1/2}\in K$ are bounded by some $c(D,K)$ at $\wtilde{Z}(q_{1/2})$ on a ball of radius $\ve_0$.
	We can now Taylor-expand $f(Z)$ along the fibers as
	\[
		f(Z)=\sum_{|I|< D}
		a_I\left(\wtilde{Z}(q_{1/2})\right)
		\left[\wtilde{Z}(q_{1/2}) - Z\right]^{I}
		+ O_{D,K}\left(\norm{\wtilde{Z}(q_{1/2}) - Z}^D\right)
	\]
	By \autoref{prop:contraction_in_v} we can arrange the last term to be $O\left(e^{-D\cdot \delta_0 \ell}\right)$ for some fixed $\delta_0$.
	We can take $D$ sufficiently large so that $D\cdot \delta_0 >N$.
	For the remaining polynomial expression in $Z$ we can use part (i) to deduce that it is factorizable, and hence we can choose appropriate bundles to approximate it to degree $N$.

	To prove the last part, we restrict in each fiber to a neighborhood of $F$, such that the fiberwise nearest-point projection $\pi_F$ to $F$ is smooth and well-defined.
	Then claim then follows by applying the preceding part to $\wtilde{f}:=f\circ \pi_F$.
\end{proof}


\subsection{Some Factorizable Functions}
\label{sec:subsec:some:factorizable:functions}

\subsubsection{On normal forms and holonomies}
	\label{sssec:on_normal_forms_and_holonomies}
Recall that the existence of holonomies is equivalent to the existence of trivializations of the bundle (smooth along stables, say) such that the transition maps, and the dynamics, are constant maps.
The property of a cocycle $E$ to be smooth along stables is defined in \autoref{def:regularity_along_stables} (as well as \autoref{sssec:cocycles_smooth_along_stable_manifolds}), its holonomies $H^{-}_{\gr^{\bullet \leq} E}$ (on the associated graded) is constructed in \autoref{cor:holonomies_on_graded}, and the measurable connection \index{$P^-_E$}$P^{-}_E$ (for any admissible cocycle) is introduced in \autoref{sssec:standard_measurable_connection}.
When the cocycle $E$ is clear from context, we will write \index{$H^-$}\index{$P^-$}$H^-,P^-$ for the corresponding holonomy and measurable connection.

Factorizing holonomies is accomplished by the next result:

\begin{proposition}[Sections of bundles with holonomy are factorizable]
	\label{prop:sections_of_bundles_with_holonomy_are_factorizable}
	Suppose that \index{$E$}$E$ is a cocycle smooth along stables, and equipped with a stable holonomy \index{$H^-$}$H^-$ which is itself smooth along stables.
	\begin{enumerate}
		\item \label{item:holonomy_section_factorizable}
		If $s$ is a measurable $g_t$-equivariant section of $E$ then $H^-(q_{1/2}',q_{1/2})\circ s(q_{1/2}')\in E(q_{1/2})$ is precisely factorizable.
		\item The map
		\[
			H^-(q_{1/2}',q_{1/2})\circ P^{-}(q_{1/2},q_{1/2}')\in \End\left(E(q_{1/2})\right)
		\]
		is also precisely factorizable.
	\end{enumerate}
	The analogous statements along unstables hold.
\end{proposition}
\begin{proof}
	The equivalent definition of smooth holonomies is that there exist smooth local trivializations of the bundle along the stable manifolds such that the transition maps, and the dynamics, act by constant maps.
	In other words, in the right choice of coordinates the maps $H^-$ are the identity.

	For the section $s$, observe that by the definition of holonomy, since $s(q_{1/2}')=g_{\ell/2}s(q')$ we have that
	\[
		H^-(q_{1/2}',q_{1/2})\circ s(q_{1/2}')
		=
		A(q,\ell/2)\circ H^-(q',q)\circ s(q')
	\]
	where $A(q,\ell/2)$ denotes the cocycle action on the bundle $E$ from $q$ to $q_{1/2}=g_{\ell/2}q$.
	We can then take $F_q(q'):=H^-(q',q)\circ s(q')$ and $V^s=E$, $V^u$ again trivial, and $\crA(q_{1/2})$ to be the identity map.

	For the measurable connection $P^-$, recall that its definition is via the measurable isomorphism $\iota_E\colon \gr_{\bullet}E\to E$ at the Oseledets-biregular points.
	The bundle $\gr_{\bullet}E$ carries a natural holonomy $H^{-}$ (since the irreducible factors have a single Lyapunov exponent) and so does $E$, thus $\iota_E(q)\in \Hom\left(\gr_{\bullet}E,E\right)(q)$ is an equivariant section and by part (i) it is factorizable.

	The map that we explicitly want to factorize is (rewriting it
        using $\iota_E$):
	\begin{align}
		\label{eq:HPgames}	
		H^-_E(q_{1/2}',q_{1/2})\circ P^-_E(q_{1/2},q_{1/2}')
		=
		\underbrace{H^-_E(q_{1/2}',q_{1/2})\circ 
		\iota_E(q_{1/2}')
		\circ
		H^-_{\gr_{\bullet E}}\left(q_{1/2},q_{1/2}'\right)}_{H^-_{\Hom}\left(q_{1/2}',q_{1/2}\right)\circ \iota_E(\wtilde{q'})}
		\circ 
		\iota_E^{-1}(q_{1/2})
	\end{align}
	where the underbraced expression is the transport of $\iota_E(\wtilde{q'})$ via the holonomy of the bundle $\Hom(\gr_{\bullet}E,E)$.
	By part (i) this expression is factorizable, and since $\iota_{E}^{-1}(q_{1/2})$ is also clearly factorizable, the desired conclusion follows from \autoref{cor:assembly_of_factorizable_maps}.
\end{proof}

Recall that $\index{$L_{q}^{\pm}$}L_{q}^{\pm}\colon \cW^{u/s}[q]\into
L\cW^{u/s}(q)$ denotes the linearization maps from
\autoref{thm:linearization_of_stable_dynamics_single_diffeo}.
We begin with a basic construction:
\begin{proposition}[$y_{1/2}$ and $q_{1/2}'$ are factorizable]
	\label{prop:wtilde_y_and_wtildeq_is_factorizable}
	\leavevmode
	\begin{enumerate}
		\item The point $L_{q_{1/2}}^-\left(q'_{1/2}\right)\in L\cW^{s}(q_{1/2})$ is precisely factorizable.
		\item The point $L_{q_{1/2}}^+\left(y_{1/2}\right)\in L\cW^{u}(q_{1/2})$ is precisely factorizable.
	\end{enumerate}
\end{proposition}

\begin{proof}
	The discussion is the same for both $y_{1/2}$ and $q_{1/2}'$ so we only treat the latter.
	The cocycle $V^s$ is taken to be $L\cW^s$ which has holonomies $P^{GM}_{L\cW^s}$ by construction and the measurable section $\iota_s(q)\in V^s(q)$ such that $L_{q_{1/2}}^-(q'_{1/2}) = P^{GM}_{L\cW^s}\left(q'_{1/2},q_{1/2}\right)\iota_s(q_{1/2}')$.
	The claim then follows directly from \autoref{prop:sections_of_bundles_with_holonomy_are_factorizable}\ref{item:holonomy_section_factorizable}.
\end{proof}




We now combine the previous propositions with \autoref{cor:assembly_of_factorizable_maps}.
For the manifold $Q$, recall that we have a locally finite collection of charts $\chi_i$.
For a smooth and natural bundle $E$ on $Q$, fix locally finite trivializing charts $\chi_{E,i}$ which are compatible with the charts $\chi_i$, see \autoref{sssec:local_trivializations_of_smooth_and_natural_cocycles}.
On compact sets $K$, we can assume that only finitely many charts are used and there is a Lebesgue number of the charts, so all points that are sufficiently close for our considerations are in the same chart.

\subsubsection{Connection from smooth charts}
	\label{sssec:connection_from_smooth_charts}
Suppose that $E$ is a smooth and natural cocycle, which is therefore equipped with a family of global charts $\chi_{E,i}$ on open sets $U_i$.
Then, for every $x,y\in U_i$ we obtain an identification
\[
	\index{$P^{\chi,i}_{E}(x,y)$}P^{\chi,i}_{E}(x,y)\colon E(x)\toisom E(y) 	
\]
which varies smoothly as $x,y$ vary in $U_i$.

\begin{proposition}[Factorizability in smooth global charts]
\label{prop:factorizability_in_smooth_global_charts}
Suppose $E$ is a smooth and natural cocycle.
\begin{enumerate}
		\item\label{item:prop:factorizability_in_smooth_global_charts:item1}
		Consider the trivial bundle $Q\times \bR^{\dim Q}$ and let $\chi_i$ denote smooth charts on $Q$ containing $q_{1/2},q'_{1/2},y_{1/2}$.
		Then the functions $\chi_{i}^{-1}(y_{1/2}),\chi_{i}^{-1}(y_{1/2})$ and $\chi_i^{-1}(q_{1/2}')$, viewed as valued in the trivial bundle $\bR^{\dim Q}$, are factorizable.
		\item\label{item:prop:factorizability_in_smooth_global_charts:item2}
		Suppose $s$ is a $g_t$-equivariant measurable section of $E$.
		Then with notation as in \autoref{sssec:connection_from_smooth_charts} we have that
		\[
			P_{E}^{\chi,i}(q'_{1/2},q_{1/2})s(q_{1/2}) \in E\left(q_{1/2}\right)
		\]
		is factorizable.


		\item
		\label{item:prop:factorizability_in_smooth_global_charts:item3}
		Recall that \index{$P^-_E$}$P^-_E$ denotes the measurable connection of $E$, see \autoref{sssec:standard_measurable_connection}.
		Then the map
		\[
			P_{E}^{\chi,i}\left(q_{1/2}',q_{1/2}\right)
			\circ
			P^{-}_E\left(q_{1/2},q_{1/2}'\right)\in \End E\left(q_{1/2}\right)
		\]
		is factorizable.

	\end{enumerate}
\end{proposition}
Note that in applying this proposition we loose precise factorizability and only have factorizability.
\begin{proof}
For (i), note that $L_{q_{1/2}}^-\left(q'_{1/2}\right)\in L\cW^s(q_{1/2})$ is factorizable by \autoref{prop:wtilde_y_and_wtildeq_is_factorizable}.
Then \autoref{cor:assembly_of_factorizable_maps}\ref{item:smooth_postcomposition_assembly} applied to post-composition with the smooth maps $\chi_i^{-1}\circ \left(L_{q_{1/2}}^{-}\right)^{-1}\colon L\cW^s(q_{1/2})\to \bR^{\dim Q}$ then gives the claim.
Here we denote by $\left(L_{q_{1/2}}^-\right)^{-1}\colon L_{q_{1/2}}^-\left(\cW^s[q_{1/2}]\right)\to \cW^s[q]$ the inverse of the linearization map, restricted to its image; it is a smooth bijection.

The factorizability of $\chi_i^{-1}(y_{1/2})$ and $\chi_{i}^{-1}\left(q_{1/2}\right)$ is analogous.

Recall from \autoref{cor:holonomy_linearization_of_arbitrary_cocycles} that associated to $E$ there is a cocycle $L^-E$ and equivariant linear embedding $L^-_q\colon E(q)\to L^-E(q)$, such that furthermore $L^-E$ admits a holonomy denoted $P^{GM}_{L^-E}$; note that $E$ is not invariant under this holonomy.

Let now $\Psi_{q'_{1/2}}\colon E(q_{1/2}) \to L^-E(q_{1/2})$ be given
by
\begin{displaymath}
\Psi_{q'_{1/2}} = P^{GM}_{L^-E}(q'_{1/2}, q_{1/2}) \circ
L^-_{q'_{1/2}} \circ P_{E}^{\chi,i}(q'_{1/2},q_{1/2}).
\end{displaymath}
Then $\Psi_{q'_{1/2}}$ is a smooth function of $q'_{1/2}$. Therefore,
the factorizability of 
$\Psi_{q'_{1/2}}$ follows from combining \ref{item:prop:factorizability_in_smooth_global_charts:item1} with 
\autoref{cor:assembly_of_factorizable_maps}\ref{item:smooth_postcomposition_assembly}. Let
$\Phi_{q'_{1/2}}:L^-E(q_{1/2}) \to E(q_{1/2})$
denote $\Psi_{q'_{1/2}}^{-1}$ precomposed with the
nearest point projection from $L^-E(q_{1/2})$ to the image of
$\Psi_{q'_{1/2}}$. Then, $\Phi_{q'_{1/2}}$ is also factorizable by 
\autoref{cor:assembly_of_factorizable_maps}. 
Now, \ref{item:prop:factorizability_in_smooth_global_charts:item2}
follows since
\begin{displaymath}
P_{E}^{\chi,i}(q'_{1/2},q_{1/2}) \circ s(q'_{1/2}) =
\Phi_{q'_{1/2}} \circ P^{GM}_{L^-E}(q'_{1/2},q_{1/2}) \circ (L^-_{q'_{1/2}} \circ s)(q'_{1/2}) 
\end{displaymath}
and $P^{GM}_{L^-E}(q'_{1/2},q_{1/2}) \circ (L^-_{q'_{1/2}} \circ s)(q'_{1/2})$
is factorizable by  
\autoref{prop:sections_of_bundles_with_holonomy_are_factorizable}.
Also, \ref{item:prop:factorizability_in_smooth_global_charts:item3}
follows since 
\begin{multline*}
P_{E}^{\chi,i}(q'_{1/2},q_{1/2}) \circ P^-_E(q_{1/2}, q'_{1/2}) = \\ = 
P_{E}^{\chi,i}(q'_{1/2},q_{1/2})  \circ (L^-_{q'_{1/2}})^{-1} \circ
P^-_{L^-E}(q_{1/2}, q'_{1/2}) \circ L^-_{q_{1/2}} = \\
= \Phi_{q'_{1/2}} \circ (P^{GM}_{L^-E}(q'_{1/2},q_{1/2}) \circ
P^-_{L^-E}(q_{1/2}, q'_{1/2}) ) \circ L^-_{q_{1/2}},
\end{multline*}
the expression in parenthesis is factorizable by
\autoref{prop:sections_of_bundles_with_holonomy_are_factorizable}, and
$L^-_{q_{1/2}}$ is available as part of the assembly map.

\end{proof}

\begin{remark}[Factorizability of equivariant subspaces]
	\label{rmk:factorizability_of_equivariant_subspaces}
	Suppose that $E$ is smooth and natural, and $V\subset E$ is a measurable $g_t$-equivariant subbundle.

	Note that $V(q'_{1/2}) = P_E^-(q_{1/2},q'_{1/2}) V(q_{1/2})$, and
	$V(q_{1/2})$ is available in the assembly map, while
	$P_E^-(q_{1/2},q'_{1/2})$ is factorizable in the sense of
	\autoref{prop:factorizability_in_smooth_global_charts}\ref{item:prop:factorizability_in_smooth_global_charts:item3}.
	We will thus say that $P^{\chi,i}_E(q'_{1/2},q_{1/2}) V\left(q'_{1/2}\right)$ is factorizable, although we have not formally defined factorizability of a subspace.

	The construction and factorizability of	$V(y_{1/2})$ is analogous.
\end{remark}

\subsubsection{Jets of manifolds and ideals of functions}
	\label{sssec:jets_of_manifolds_and_ideals_of_functions}
Suppose $\phi_1,\phi_2\colon \bR^d\to \bR^{n}$ are two germs of smooth maps, with $\phi_i(\bold0)=\bold 0$, and $\rk D\phi_i(\bold 0)=d\leq n$.
We will say that they have \emph{image tangent to order $k$} if there exists a germ of diffeomorphism $\chi$ of $(\bR^d,\bold 0)$ such that $\phi_1$ and $\phi_2\circ \chi$ have the same Taylor expansion up to order $k$.
We let \index{$M$@$\cM_k$}$\cM_k$ denote such an equivalence class.
Let also \index{$T$@$\cT_k(\bR^n,\bold0)$}$\cT_k(\bR^n,\bold0)$ denote the finite-dimensional ring of Taylor coefficients of functions at $\bold 0$ up to order $k$.
Then $\cM_k$ determines an ideal \index{$I$@$\cI\cM_k$}$\cI\cM_k\subset \cT_k(\bR^n,\bold 0)$ given by those functions that vanish up to order $k$ when pulled back by some representative $\phi\colon \bR^d\to \bR^n$ from $\cM_k$.

Let us explain how to recover $\cM_k$ from $\cI\cM_k$ (assuming $k\geq 1$):
\begin{proposition}[Taylor coefficients from ideals of functions]
	\label{prop:taylor_coefficients_from_ideals_of_functions}
	Suppose that we have the decomposition $\bR^n=\bR^d\times \bR^{n-d}$ and we have a germ of manifold $\cM_k$ such that the tangent space of $\cM_k$ projects isomorphically to the $\bR^d$-factor.
	\begin{enumerate}
		\item 
		Then there exists a unique degree $k$ polynomial map $\phi\colon \bR^d\to\bR^d\times \bR^{n-d}$ in the equivalence class determining $\cM_k$, such that $\phi(\bbx)=(\bbx,\psi(\bbx))$ with $\psi\colon \bR^d\to \bR^{n-d}$.

		\item 
		Furthermore, if we fix a basis of equations $\xi_{1},\dots \xi_{n-d}$ cutting out the tangent space of $\cM_k$, as well as lifts $\wtilde{\xi}_i\in \cI\cM_k$, then these lifts generate the ideal $\cI\cM_k$ and the Taylor coefficients of $\psi$ are rational functions of the coefficients of $\wtilde{\xi}_i$.
	\end{enumerate}
\end{proposition}
\begin{proof}
	The first claim is an immediate consequence of the constant rank theorem, which in this case requires only a finite-dimensional calculation.
	Let therefore $\psi=(\psi_1,\dots \psi_{n-d})$ be the corresponding function such that $\cM_k$ is its graph.

	For the second claim, let $\bby=(y_1,\dots y_{n-d})$ be coordinates on $\bR^{n-d}$.
	Then one set of generators for the ideal is $y_i-\psi_i(\bbx)$ and the second claim for them is clear.
	For any other set of generators $\wtilde{\xi}_i$, there exist uniquely determined polynomials $p_{i}^j$ such that $y_i-\psi_i(\bbx)\equiv\sum p_{i}^j(\bbx,\bby) \wtilde{\xi}_j(\bbx,\bby)$ up to order $k$.
	Note that the constant terms of $p_{i}^j$ are a linear matrix, necessarily inverse to the first order term matrix in $\bby$ of the $\wtilde{\xi}_j$.
	Subsequent coefficients of $p_{i}^j$ amount to solving further linear equations in the coefficients of $\wtilde{\xi}_j$, guaranteed to have solutions since $\wtilde{\xi}_j$ are generators.
\end{proof}
In our intended applications, given a basis $\xi_i$ of linear functions cutting out the tangent space, we will take the unique lifts $\wtilde{\xi}_i$ that minimize the natural background norm on the space of $k$-jets of functions.


\subsubsection{Factorizing Taylor Coefficients}
\label{sssec:factorizing_taylor_coefficients}
Let \index{$J$@$\cJ^kQ$}$\cJ^kQ$ denote the smooth and natural cocycle of $k$-jets of
real-valued smooth functions, and furthermore the fibers have a ring
structure. We fix locally trivializing charts \index{$\chi_{\cJ^kQ,i}$}$\chi_{\cJ^kQ,i}$ for
$\cJ^kQ$ which are compatible with the charts $\chi_i$, see \autoref{sssec:local_trivializations_of_smooth_and_natural_cocycles}.
Let \index{$I$@$\cI^u\cJ^kQ$}$\cI^u\cJ^kQ\subset \cJ^k Q$ be the subcocycle of functions that
vanish to order $k$ along the unstable manifolds. This is an 
ideal for the ring structure and uniquely determines the $k$-jet of
the unstable manifold.

Similarly, let $\cI^{cs}\cJ^kQ\subset \cJ^k Q$ denote the ideal of
functions that vanish to order $k$ along the center-stable
manifold. (The center stable manifold may not be well-defined, but
it's $k$-jet, hence the cocycle $\cI^{cs}\cJ^kQ$ is well-defined, see
\autoref{prop:properties:fake:flags} \autoref{item:fakeflag:smoothness}.)

Now suppose $E$ is a smooth and natural cocycle. We fix locally
trivializing charts \index{$\chi_{E,i}$}$\chi_{E,i}$ for $E$ which are
compatible with the charts \index{$\chi_i$}$\chi_i$. 
Let \index{$E^{\ge\bullet}$}$E^{\ge\bullet}$ denote the unstable Lyapunov
flag of $E$.
Then, for almost all $q$, $E^{\ge\bullet}(q)$ is a
well-defined point in the flag variety $\cF^{\ge\bullet} E (q)$.
Similarly,
if $E^{\le\bullet}$ denotes the stable Lyapunov flag of $E$, then for 
almost all $q$, $E^{\le\bullet}(q)$ is a
well-defined point in the flag variety $\cF^{\le\bullet} E (q)$. 
Also,
for $z \in \cW^{cs}[y_{1/2}]$ (not necessarily forwards regular), we
denote by $E_{y_{1/2}}^{\le\bullet}(z)$ the forwards flag defined in
\autoref{def:fixed_realization_center_stable}.

\subsubsection{Trivializing maps to bundles, along submanifolds}
	\label{sssec:trivializing_maps_to_bundles_along_submanifolds}
Suppose $E\to Q$ is a smooth and natural cocycle.
Let $\cF_E\to Q$ be an associated bundle with fiber the set of flags in $E$ of fixed dimensions (omitted from the notation), e.g. the Grassmannian of $k$-planes for some $k$.
Using the charts $\chi_{i}$ we will trivialize the bundle \index{$F$@$\cF_E$}$\cF_E$ with the maps $P^{\chi,i}_E(x,q_{1/2})$ to identify each fiber with the fiber at $q_{1/2}$, for $x$ in a sufficiently small neighborhood of $q_{1/2}$.
Next, we assume given at $q_{1/2}$ a flag $Z(q_{1/2})\in \cF_E(q_{1/2})$.
Then, using an orthogonal decomposition to turn the flag $Z$ into a decomposition $E(q_{1/2})=\oplus Z_i(q_{1/2})$, we introduce the (standard) charts on $\cF_{E}(q_{1/2})$ near $Z(q_{1/2})$ into the vector space $\oplus_{j>i}\Hom\left(Z_i,Z_j\right)$.

Suppose now that we have a smooth section $\psi\colon \cW^u[q'_{1/2}]\to \cF_{E}$ over an unstable manifold.
Using the trivialization $\cF_E\vert_{U_i}\isom U_i\times \cF_E(q_{1/2})$ and the previous charts, we obtain a function $\cW^u[q'_{1/2}]\xrightarrow{\psi_Z} \oplus_{j>i}\Hom(Z_i,Z_j)$.
Again in the charts $\chi_i$ starting from the decomposition $T_{q_{1/2}}Q=W^u(q_{1/2})\oplus W^{cs}(q_{1/2})$, let $W^{\bullet}_{\chi,i}(q_{1/2})\subset \bR^n$ be the affine linear subspaces that pass through $\chi_i(q_{1/2})$ and are tangent to the respective unstable/center-stable submanifolds.
Then we know that $\chi_i\left(\cW^{u}[q'_{1/2}]\right)$ is the graph of a smooth function
\begin{align}
	\label{eqn:unstable_manifold_in_chi_i_charts_graph}
	\index{$\eta^u_{q'_{1/2}}$}\eta^u_{q'_{1/2}}\colon W^u_{\chi,i}(q_{1/2})\to W^{cs}_{\chi,i}(q_{1/2})
\end{align}
Extend now $\psi$ to a scalar-valued function $\psi_{\chi,i}\colon \chi_{i}(U_i)\to \oplus_{j>i}\Hom(Z_i,Z_j)$ defined on all of $\chi_i(U_i)\subset \bR^n$ where
\begin{align}
 	\label{eqn:psi_to_psi_chi_i}
 	\psi_{\chi,i}=\psi_Z\circ \chi_i^{-1}\circ (\id,\eta^u_{q'_{1/2}})\circ \bar{\pi}^{cs}_u
 \end{align}
 where \index{$\pi$@$\bar{\pi}^{cs}_u$}$\bar{\pi}^{cs}_u$ is the projection in $\bR^n$ to $\cW^{u}_{\chi,i}(q_{1/2})$ along $\cW^{cs}_{\chi,i}(q_{1/2})$.

We will make use of analogous constructions along the center-stable manifold with $y_{1/2}$ instead of $q'_{1/2}$.
Here is a first application:

\begin{proposition}[Factorizing the unstable along the stable]
	\label{prop:factorizing_the_unstable_along_the_stable}
	The subspaces
	\[
		P^{\chi,i}_{\cJ^k Q}(q'_{1/2},q_{1/2}) \cI^u\cJ^kQ(q'_{1/2}) \text{ and } P^{\chi,i}_{\cJ^k Q}(y_{1/2},q_{1/2}) \cI^{cs}\cJ^kQ(y_{1/2})
	\]
	are factorizable in the sense of \autoref{rmk:factorizability_of_equivariant_subspaces}.

	Furthermore, for any order $k$ the Taylor coefficients of $\eta^u_{q'_{1/2}}$ as defined in \autoref{eqn:unstable_manifold_in_chi_i_charts_graph} are factorizable, and the same holds for its analogue \index{$\eta^{cs}_{y_{1/2}}$}$\eta^{cs}_{y_{1/2}}$.
\end{proposition}

\begin{proof}
The fist assertion is an immediate consequence of
\autoref{prop:factorizability_in_smooth_global_charts}
\autoref{item:prop:factorizability_in_smooth_global_charts:item3}.
The second one follows from the first, combined with \autoref{prop:taylor_coefficients_from_ideals_of_functions}.
\end{proof}

\subsubsection{Realized $C^\infty$ jets of center-stable manifolds and flags of cocycles}
For the next statement, we recall and will use the notation from the appendix \autoref{def:fixed_realization_center_stable}.

\begin{proposition}[Factorizing the jet of the unstable flag]
	\label{prop:factorizing_the_jet_of_the_unstable_flag}
	Let $E\to Q$ be a smooth and natural cocycle.
\begin{enumerate}
	\item
	Let $S\subset E$ be a backwards dynamically defined subbundle.
	Let $\Gr E$ denote the Grassmannian of $E$, and let \index{$\psi^{S}(q'_{1/2})$}$\psi^{S}(q'_{1/2})\colon \cW^u[q'_{1/2}] \to \Gr E$ denote the map taking $z \in \cW^u[q'_{1/2}]$ to $S(z)$.
	Let \index{$\psi^S_{\chi,i}(q'_{1/2})$}$\psi^S_{\chi,i}(q'_{1/2})$ be the associated (vector-space valued) map constructed in \autoref{sssec:trivializing_maps_to_bundles_along_submanifolds}, using the ambient charts.
	Then, for any order $k$ the Taylor coefficients at $\chi_i(q_{1/2}')$ of $\psi^{S}_{\chi,i}(q'_{1/2})$ are factorizable.

	\item
	Let $\psi^{u}(q'_{1/2})\colon \cW^u[q'_{1/2}] \to \cF^{\ge\bullet} E$ denote the map taking $z \in \cW^u[q'_{1/2}]$ to the Oseledets flag $E^{\ge\bullet}(z)$.
	Let $\psi^{u}_{\chi,i}(q'_{1/2})$ be the associated map constructed in \autoref{sssec:trivializing_maps_to_bundles_along_submanifolds}.
	Then, for any order $k$ the Taylor coefficients at $\chi_i(q'_{1/2})$ of $\psi^{u}_{q'_{1/2}}$ are factorizable.

	\item Let \index{$\psi^{cs}(y_{1/2})$}$\psi^{cs}(y_{1/2})\colon \cW^{cs}[y_{1/2}] \to \cF^{\le\bullet} E$ denote the map taking $z \in \cW^{cs}[y_{1/2}]$ to the flag
	$E_{y_{1/2}}^{\le\bullet}(z)$.
	Let \index{$\psi^{cs}_{\chi,i}(y_{1/2})$}$\psi^{cs}_{\chi,i}(y_{1/2})$ be the associated map constructed in \autoref{sssec:trivializing_maps_to_bundles_along_submanifolds}.
	Then, for any order $k$, the Taylor coefficients at $\chi_i(y_{1/2})$ of $\psi^{cs}_{\chi,i}(y_{1/2})$ is factorizable.
\end{enumerate}
\end{proposition}

\begin{proof}
Note that since the components $E^{\ge \lambda_i}$ of the flag $E^{\ge \bullet}$
are backwards dynamically defined, (ii) follows immediately from (i). 

Let now $\cJ^k E$ denote the smooth and natural cocycle of $k$-jets of
real-valued smooth functions on the total space of the bundle $E$.
(This is naturally a cocycle over the total space of $E$, but we
restrict it to a cocycle over $Q$ which is the zero section of the
total space of $E$).

For the proof of (i), let
\begin{displaymath}
	I^{S}(q) = \{ (x,v) \st x \in \cW^u[q], v \in S(x) \}. 
\end{displaymath}
Then $I^{S}(q)$ is a smooth (immersed) submanifold of the total space of $E$.
Let $I^{S} \cJ^k E(q) \subset \cJ^k E(q)$ denote the ideal of jets of functions which vanish to order $k$ on $I^S(q)$. 
Then, $I^{S} \cJ^k E$ is a measurable equivariant subbundle of $\cJ^k E$ and thus by \autoref{prop:factorizability_in_smooth_global_charts}
\ref{item:prop:factorizability_in_smooth_global_charts:item3}, $P^{\chi,i}_{\cJ^k E}(I^{S} \cJ^k E(q'_{1/2}))$ is factorizable.
Next, by \autoref{prop:taylor_coefficients_from_ideals_of_functions} the ideal $I^{S} \cJ^{k'} E(q'_{1/2})$ determines the Taylor coefficients of $\psi^S$ at $q'_{1/2}$.
Using that all the other functions need to pass from $\psi^S$ to $\psi^S_{\chi,i}$ appearing in \autoref{eqn:psi_to_psi_chi_i} are defined at $q_{1/2}$, combined with \autoref{cor:assembly_of_factorizable_maps}, yields the claim.


The proof of (iii) is virtually identical using the equivariance of $k$-jets in \autoref{prop:properties:fake:flags}.
\end{proof}

Recall from \autoref{sssec:standard_measurable_connection} that associated to the unstable filtration of a smooth cocycle we have the associated graded bundle $\gr_{\bullet}^u E$ which is smooth along unstables and admits a smooth unstable holonomy $H^+$.
Using a chart $\chi_i$ for $E$ that identifies each fiber with some $\bR^{m_E}$, we have an identification $\gr^u_{\bullet}E\to E\to \bR^{m_E}$ by using the orthogonal complement with respect to the standard Euclidean structure on $\bR^{m_E}$ to split the filtration defining $\gr^{u}_{\bullet}E$, followed by the trivializing chart of $E$.
Let \index{$H^{+}_{\chi,i}(x,y)$}$H^{+}_{\chi,i}(x,y)\colon \bR^{m_E}\to \bR^{m_E}$ be the conjugation of the holonomy map by these trivializing charts, with $y\in \cW^u[x]$.


\begin{proposition}[Factorizing the jet of holonomy of the associated graded]
	\label{prop:factorizing_the_jet_of_the_holonomy_of_graded}
	Recall from \autoref{sssec:trivializing_maps_to_bundles_along_submanifolds} that $\bar{\pi}_u^{cs}$ is the projection in the chart $\chi_i$ along the (linear) center-stable manifold to the unstable, and $\left(\id,\eta^u_{q'_{1/2}}\right)$ is the map which gives the unstable through $q'_{1/2}$ in the same chart.

	Then for any $k$ the Taylor coefficients of order $k$ of
	\[
	H^+_{\chi,i}\left(q'_{1/2}, \chi_{i}^{-1}\circ (\id,\eta^u_{q'_{1/2}}) \circ \bar{\pi}^{cs}_u(-)\right)
	\colon \bR^{\dim Q}\to \End(\bR^{m_E}) 	
	\]
	are factorizable.

\end{proposition}

\begin{proof}
	We proceed analogously to the proof of \autoref{prop:factorizing_the_jet_of_the_unstable_flag}.
	Namely, let
	\begin{align*}
		IH^+(q):= & \{ (x,H) \colon x \in \cW^u[q], H\in \Hom\left(E(q),E(x)\right) \text{ such that }\\
		&  H\left(E^{\geq \bullet}(q)\right)=E^{\geq \bullet}(x) \text{ and } 
		gr^u_\bullet H = H^+(x,q)\}
	\end{align*}
	where $\gr^u_\bullet H$ denotes the map induced by $H$ on the associated graded.
	Then $IH^+(q)$ is a smooth submanifold of the linear bundle $\Hom(E(q),E)$ where we view $E(q)$ as a fixed vector space, and so the fiber above $x\in \cW^u[q]$ is $\Hom(E(q),E(x))$.
	The intersection of $IH^+(q)$ with a fiber is naturally a smooth algebraic subset (acted upon freely and transitively by the group of linear transformations that act as the identity on the associated graded).

	To obtain the relevant cocycle on our space $Q$, consider first $Q\times Q$, with the cocycle $\Hom(E_1,E_2)$ where $E_i$ corresponds to the pullback of $E$ under projection to the $i$-th factor.
	Then we have the cocycle of jets $\cJ^k\Hom(E_1,E_2)$, which we can restrict to the diagonal $Q\into Q\times Q$ and denote by $\cJ^k_{\Delta}\Hom(E_1,E_2)$.
	Then $IH^+(q)$ determines an ideal of functions $\cI\cH^+(q)\subset \cJ^k_{\Delta}\Hom(E_1,E_2)$ vanishing on it.

	It now follows from \autoref{prop:taylor_coefficients_from_ideals_of_functions} that the Taylor coefficients of the desired function are factorizable.
\end{proof}

\begin{remark}[Factorizing holonomy]
	\label{rmk:factorizing_holonomy}
	Suppose that $E$ is a bdd cocycle, smooth along unstables, and equipped with a holonomy map $H^+_E$.
	Then the analogous $H^+_{\chi,i}$ from \autoref{prop:factorizing_the_jet_of_the_holonomy_of_graded} is factorizable, with the same proof.
	We will apply this remark to $E=L\cW^u$.
\end{remark}

\subsection{More factorizable functions}
	\label{ssec:more_factorizable_functions}

\subsubsection{Setup}
	\label{sssec:setup_more_factorizable_functions}
In this section, we will continue to work extensively with the charts
$\chi_i\colon Q\to \bR^{\dim Q}$ and assume that all points with subscript
${1/2}$ are in a single such chart. 

Recall that center-stable manifolds are defined in
\autoref{def:fixed_realization_center_stable}. 


\begin{proposition}[$z_{1/2}$ is factorizable]
	\label{prop:wtildez_is_factorizable}
	The point
	\[
		\chi_i\left( 
		\cW^{cs}_{loc}\left[y_{1/2}\right]
	\bigcap  
	{\cW^u_{loc}\left[q_{1/2}'\right]}
	 \right)\in \bR^{\dim Q}
	 \text{ see \autoref{fig:halfway_diagram}}
	\]
	is factorizable.
\end{proposition}
\begin{proof}
	For $q_{1/2}$, we can apply a smooth coordinate transformation $\tilde{\cA}_1(q_{1/2})$ after $\chi_i$ so as to arrange that, with coordinates $(x,y)\in \bR^n\times \bR^m$:
	\begin{itemize}
		\item $q_{1/2}$ lands at $(0,0)\in \bR^{n}\times \bR^m$.
		\item The center-stable manifold $\cW^{cs}[q_{1/2}]$ becomes $y=0$.
		\item The unstable manifold $\cW^u[q_{1/2}]$ becomes $x=0$.
	\end{itemize}
	We will write $\tilde{q}_{1/2}'\in \bR^n$ and $\tilde{y}_{1/2}\in \bR^m$ for the images of the corresponding points on the coordinate planes.
	The coordinates of these points are factorizable, by \autoref{prop:factorizability_in_smooth_global_charts}.

	In these coordinates, the unstable manifold at $q_{1/2}'$ and center-stable manifold at $y_{1/2}$ can be written using smooth functions $\xi^u_{\tilde{q}_{1/2}'}\colon \bR^m \to \bR^n$ and $\xi^{cs}_{\tilde{y}_{1/2}}\colon \bR^n\to \bR^m$ as
	\[
		\cW^u[q_{1/2}'] = \text{image of} \begin{bmatrix}
			\tilde{q}_{1/2}' + \xi^u_{\tilde{q}_{1/2}'}(y)\\
			y
		\end{bmatrix}
		\quad
		\cW^{cs}[y_{1/2}] = \text{image of} \begin{bmatrix}
			x\\
			\tilde{y}_{1/2} + \xi^{cs}_{\tilde{y}_{1/2}}(x)
		\end{bmatrix}
	\]
	In view of
  \autoref{prop:factorizing_the_unstable_along_the_stable}, 
          for any $k\in \bN$, the Taylor coefficients of
          $\xi^{\bullet}_{\bullet}$ up to degree $k$ are
          factorizable.
          Also, there exist $c_k$ and $\delta_k>0$
        such that the $C^k$-norm of these functions are bounded by
        $c_ke^{-\delta_k \ell}$, both in their dependence on $x/y$ and
        $q_{1/2}'/y_{1/2}$.

	To find the point of intersection via the implicit function theorem, we can write
	\[
		F\left(x,y,\tilde{q}_{1/2}',\tilde{y}_{1/2}\right) =
		\begin{bmatrix}
		 	\tilde{q}_{1/2}' - x + \xi^u_{\tilde{q}_{1/2}'}(y)\\
		 	y - \tilde{y}_{1/2} - \xi^{cs}_{\tilde{y}_{1/2}}(x)
		 \end{bmatrix} 
	\]
	The differential of this function is
	\[
		DF = 
		\begin{bmatrix}
			-1 & O(e^{-\delta_k \ell}) & 1 + O(e^{-\delta_k\ell}) &  0\\
			 O(e^{-\delta_k \ell}) & 1 & 0 & 1 + O(e^{-\delta_k\ell})
		\end{bmatrix}
	\]
	Note that $F(0,0,0,0)=(0,0)$ and the first block is invertible, so we can write $F^{-1}(0,0)$ as
	\[
		(x,y) = G_0(\tilde{q}_{1/2}',\tilde{y}_{1/2}) = \left(\tilde{q}_{1/2}',\tilde{y}_{1/2}\right) +
		G_1(\tilde{q}_{1/2}',\tilde{y}_{1/2})
	\]
	where furthermore we have an estimate $\norm{D^kG}\leq O\left(e^{-\delta_k l}\right)$ in the range where $q_{1/2}'$ and $y_{1/2}$ are themselves $O\left(e^{-\delta_0\ell}\right)$.

	By successive applications of \autoref{cor:assembly_of_factorizable_maps}, it follows that the point $\tilde{z}_{1/2}=(\tilde{q}_{1/2}',\tilde{y}_{1/2})+G_1(\tilde{q}_{1/2}',\tilde{y}_{1/2})$ is itself factorizable.
\end{proof}

\subsubsection{Trivializing bdd bundles}
	\label{sssec:trivializing_bdd_bundles}
Suppose $L$ is a bdd cocycle as defined in \autoref{def:forward_dynamically_defined_cocycle}, with $E$ an associated smooth and natural cocycle, and two bdd subcocycles $S_1\subset S_2\subset E$ such that $L\isom S_2/S_1$.
Suppose $\chi_{E,i}$ are charts trivializing $E$.
Denote by $L_{\chi,i}\subset \bR^{\dim E}$ the subspace obtained by taking the orthogonal complement of $\chi_{E,i}(S_1)$ inside $\chi_{E,i}(S_2)$, with respect to the standard Euclidean inner product.

More generally, given a filtration $L^{\geq \bullet}$ of $L$ by bdd cocycles, we again let $\gr^{k}_{\chi,i} L$ denote the orthogonal complement of $\chi_{E,i}\left(S_1 + L^{\geq k-1}\right)$ inside $\chi_{E,i}\left(S_1 + L^{\geq k}\right)$.
In the case when the filtration is the Oseledets filtration of $L$, using the splitting provided by this construction, we denote by $H^{+}_{\chi,i}$ the holonomy map with respect to the decomposition $L_{\chi,i}\isom \oplus_k \gr^k_{\chi,i}L$.

Recall now that $L\cW^u$ is a subquotient of the smooth and natural cocycle
$E$ and is itself bdd, by \autoref{prop:linearization_cocycle_is_fdd}.
In what follows, we are fixing a trivialization of the bundle $E$,
and assuming that all the points with subscript $1/2$ belong to a
single chart.

\begin{lemma}
\label{lemma:unstable:flag:factorizable}
The backwards Lyapunov filtration $(L\cW^u)^{\ge \bullet}(z_{1/2})$ is a bdd cocycle, and in the notation of \autoref{sssec:trivializing_bdd_bundles} let the subspace $\left(L\cW^u\right)^{\geq \bullet}_{\chi,i}(z_{1/2})$ is factorizable.

Furthermore, the holonomy on the associated graded
$H^+_{\chi,i}(q'_{1/2},z_{1/2}): \gr^u_{\chi,i} L\cW^u(q'_{1/2}) \to \gr^u_{\chi,i}
L\cW^u(z_{1/2})$ is factorizable. 
\end{lemma}

\begin{proof} 
The cocycle $L\cW^u$ is backward dynamically defined (bdd) by \autoref{prop:linearization_cocycle_is_fdd} and thus there exists a smooth and natural cocycle $E$ and equivariant subbundles $S_2$ and $S_1$ of $E$, smooth along
unstables, such that $L\cW^u = S_2/S_1$.
We will work exclusively in the trivialization $\chi_{E,i}$ of $E$ and omit it from the notation below, as well as from the subscripts for the objects, but these are always implicitly present.

The unstable flag $E^{\ge\bullet}(z_{1/2})$ is smooth as a
function of $z_{1/2} \in \cW^u[q'_{1/2}]$. We can then approximate
$E^{\ge\bullet}(z_{1/2})$ using a Taylor series based at $q'_{1/2}$.
In view of \autoref{prop:factorizing_the_jet_of_the_unstable_flag},
the Taylor coefficients of this Taylor series are factorizable. 
Thus, the flag $E^{\ge\bullet}(z_{1/2})$ is factorizable.

Since $S_2$ and $S_1$ are bdd and hence also smooth as a function of $z_{1/2} \in
\cW^u[q'_{1/2}]$, we can also approximate them using a Taylor series
based at $q'_{1/2}$. In view of
\autoref{prop:factorizing_the_jet_of_the_unstable_flag}, the Taylor
coefficients of these Taylor series are also factorizable. Therefore,
we can factorize $S_2(z_{1/2})$ and $S_1(z_{1/2})$. Since the desired flag
$(L\cW^u)^{\ge\bullet}(z_{1/2})$ can be recovered from $E^{\ge\bullet}(z_{1/2})$,
$S_2(z_{1/2})$ and $S_1(z_{1/2})$, we have shown the factorizability
of $(L\cW^u)^{\ge\bullet}(z_{1/2})$.

Let $H^+_E(q'_{1/2},z_{1/2}): \gr^u E(q'_{1/2})\to \gr^u E(z_{1/2})$ denote the
holonomy of $E$ on the associated graded $\gr^u E$. Then, $H^+_E(q'_{1/2},z_{1/2})$ is smooth as a function of $z_{1/2} \in
\cW^u[q'_{1/2}]$. We approximate it using a Taylor series based at
$q'_{1/2}$. By
\autoref{prop:factorizing_the_jet_of_the_holonomy_of_graded}, the
coefficients of this Taylor series are factorizable. Therefore
$H^+_E(q'_{1/2},z_{1/2})$ is factorizable.



Since at every bireguar point $q$, 
$S_2(q) = \oplus_i S_2(q) \cap E^{\lambda_i}(q)$, and $S_2$ is
backwards dynamically defined, we have 
$S_2^{\ge \lambda_i} = S_2 \cap E^{\ge \lambda_i}$. Hence,
\begin{displaymath}
\gr^u S_2 \equiv \bigoplus_i S_2^{\ge \lambda_i}/S_2^{\ge
  \lambda_{i+1}} = \bigoplus_i (S_2 \cap E^{\ge \lambda_i})/(S_2 \cap E^{\ge \lambda_{i+1}}). 
\end{displaymath}
Let $h_2: \gr^u S_2 \to \gr^u E$ denote $\bigoplus_i h_{2,i}$, where
\begin{displaymath}
h_{2,i}: (S_2 \cap E^{\ge \lambda_i})/(S_2 \cap E^{\ge
  \lambda_{i+1}}) \to E^{\ge \lambda_i}/E^{\ge \lambda_{i+1}}
\end{displaymath}
is defined by $h_{2,i}(v) = v + E^{\ge \lambda_{i+1}}$.
Then, $h_2$ is an injective
equivariant map. Let $h_1: \gr^u S_1 \to \gr^u S_2$
be the analoguous map for $S_1$. Note that $\gr^u L\cW^u = \gr^u(S_2/S_1)$ is
naturally isomorphic to $\gr^u S_2/h_1(S_1)$, and also that $h_2$
induces an isomorphism between $\gr^u S_2/h_1(S_1)$ and $h_2(S_2)/(h_2
\circ h_1)(S_1)$. Let $S_2' = h_2(S_2)$, $S_1' = (h_2 \circ
h_1)(S_1)$. Then, $S_1' \subset S_2' \subset \gr^u E$, and we have an
equivariant map  $j_q: \gr^u L\cW^u(q) \to \gr^u E(q)$.


By the Ledrappier
invariance principle, $S_2'$ and $S_1'$ are invariant by 
$H^+_E$.
Therefore, $j_{x}^{-1} \circ H^+_E(q,x) \circ j_q: L\cW^u(q)
\to L\cW^u(x)$ is a smooth holonomy map.
Since such maps are unique (see \cite[Prop.~4.2(c)]{KalininSadovskaya2013} for the case of a single Lyapunov block, and we are working on a direct sum of such blocks), this map has to coincide with the
standard holonomy
$H^+(q,x)$.  Thus, $H^+(q'_{1/2}, z_{1/2})$, which is the holonomy on the associated graded of $L\cW^u =
S_2/S_1$, can be recovered from $H^+_E(q'_{1/2}, z_{1/2})$,
$S_2(q'_{1/2})$, $S_1(q'_{1/2})$, $S_2(z_{1/2})$ and $S_1(z_{1/2})$
and the flags $E^{\ge \bullet}(q'_{1/2})$, $E^{\ge
  \bullet}(z_{1/2})$. Since all of these
quantities were shown to be factorizable, $H^+(q'_{1/2},
z_{1/2})$ is also factorizable.  
\end{proof}

\begin{proposition}[Factorizability of realized flags]
\label{prop:factorizability_of_realized_fake_flags}
The forward flag $\chi_{E,i}\left(E^{\le\bullet}_{y_{1/2}}(z_{1/2})\right)$ as defined in 
\autoref{def:fixed_realization_center_stable} is factorizable. 
%
\end{proposition}
\begin{proof}
The forward flag $E^{\le\bullet}_{y_{1/2}}(z_{1/2})$ is smooth as a
function of $z_{1/2} \in W^{cs} [y_{1/2}]$.
We can then approximate
$\chi_{E,i}\left(E^{\le\bullet}_{y_{1/2}}(z_{1/2})\right)$ using a Taylor series based at $\chi_i(y_{1/2})$.
In view of \autoref{prop:factorizing_the_jet_of_the_unstable_flag},
the Taylor coefficients of this Taylor series are factorizable. 
Thus, the flag $\chi_{E,i}(E^{\le\bullet}_{y_{1/2}}(z_{1/2}))$ is factorizable.
\end{proof}

\begin{corollary}
\label{cor:good:flag:factorizable}
	The flags $\chi_{E,i}\left(\tilde{F}\right)$ and $\chi_{E,i}\left(\tilde{F}^{sr}\right)$ of \autoref{prop:fake_future_flag_and_subresonant_map} are
	factorizable.   
\end{corollary}

\begin{proof}
The factorizability of $\tilde{F}$ follows immediately from
\autoref{prop:factorizability_of_realized_fake_flags}, the fact that
$\tilde{F} \in \cF^{s}_{\ell}E(z_{1/2})$ and
\autoref{prop:fake_forward_flags_on_fake_center_stables}.
The factorizability of $\tilde{F}^{sr}$ follows from the
factorizability of $\tilde{F}$ and of the subspaces $S_1(z_{1/2})$,
$S_2(z_{1/2})$ (established in the course of the proof of
\autoref{lemma:unstable:flag:factorizable}).
\end{proof}

\begin{lemma}
\label{lemma:Pfake:factorizable}
The map $P^+_{fake}(q'_{1/2},z_{1/2}): L\cW^u(q'_{1/2}) \to L\cW^u(z_{1/2})$ of
\autoref{thm:fake_holonomies_and_ssr_map}, after conjugation by $\chi_{E,i}$ is factorizable.   
\end{lemma}

\begin{proof}
In view of \autoref{cor:good:flag:factorizable} and
\autoref{lemma:unstable:flag:factorizable} the decomposition
$L\cW^u(z_{1/2}) = \oplus L\cW^u(z_{1/2})^{\lambda_i}$
of \autoref{thm:fake_holonomies_and_ssr_map}
is factorizable, after applying $\chi_{E,i}$, is factorizable.
Also, in view of \autoref{lemma:unstable:flag:factorizable},
the holonomy maps from
$q'_{1/2}$ to $z_{1/2}$ on the associated graded are factorizable.
Now the assertion of the lemma follows from
\autoref{cor:assembly_of_factorizable_maps}.   
\end{proof}




\subsection{Factorizability of the interpolation map.}
\label{sec:subsec:factizability:of:interpolation}

\subsubsection{}
	\label{sec:subsec:def:S:LC}
Recall that we defined in \autoref{sssec:defining_the_interpolation_map} the interpolation map
\[
	\wtilde{\phi}\colon \cW^{u}[y_{1/2}]\to \cW^{u}[z_{1/2}]
\]
For the next construction, we use the (canonical) identifications $\cW^{u}[y_{1/2}]=\cW^{u}[q_{1/2}]$ and $\cW^{u}[z_{1/2}]=\cW^{u}[q'_{1/2}]$.

Define then \index{$Z_{1/2}(q,q',u,\ell)$}$Z_{1/2}(q,q',u,\ell)$ to be the element in $L\cC(q_{1/2})$ defined by
\begin{align}
	\label{eqn:Z_12_def}
	Z_{1/2}\left(q,q',u,\ell\right)
	=  L\left(\td \phi\inv\left(  U^+[q'_{1/2}]\right)\right).
\end{align}
In this formula $L$ refers to the linearization map
$\cC[q_{1/2}] \to L\cC[q_{1/2}]$.
For the next statement, we equip $L\cC$ with the norm $\| \cdot \|$ constructed in \autoref{prop:good_norms}.

\begin{theorem}[Half-time factorizability of 3-variable map]
	\label{theorem:3varA:haltime:factorizable}
The function $Z_{1/2}(\cdot,\cdot,\cdot,\cdot)$ is half-time factorizable in the sense of \autoref{def:half_time_factorizable}.
Therefore, for every $N \in \mathbb{N},\delta>0$, there exists a compact set $K$ of measure at least $1-\delta$, a cocycle \index{$V^s$}$V^s$, a $g_t$-equivariant embedding \index{$F_q^{s}$}$F_q^{s}\colon \cW^{s}[q]\into V^{s}(q)$ and a linear map 
\index{$A$@$\cA_{1/2}(q,u, \ell)$}$\cA_{1/2}(q,u, \ell) \colon V^s(q) \to L\cC(q_{1/2})$, such that provided $q,q',q_{1/2},q'_{1/2}\in K$ we have for any $q' \in \gB^-_0[q]\cap K$ that: 
\begin{equation}
\label{eq:thm:3varA:halftime:factorizable}
\norm{\cA_{1/2}(q,u,\ell) F^s_q(q') - L\left(\td \phi\inv(  U^+[q'_{1/2}])\right) } \le
e^{-N \ell}.
\end{equation}
Furthermore, we can choose $\cA_{1/2}(\cdot, \cdot, \cdot)$ so that
\begin{equation}
\label{eq:3varcA:q:in:LC}
\cA_{1/2}(q,u,\ell) F^s_q(q) \in L\cC^{\tau}(q_{1/2}), 
\end{equation}
where $L\cC^{\tau}$ is the tautological subspace defined in
\S\ref{sec:subsec:def:LC:tau}. 
\end{theorem}

In the above statement, we absorbed the cocycle $V^u\left(g_{-\ell/2}uq\right)$ and its equivariant section into the map $\cA_{1/2}(q,u,\ell)$.

\begin{proof}
Recall from \autoref{sssec:defining_the_interpolation_map} that $\tilde{\phi} = \wp^+_{fake}(q'_{1/2},z_{1/2}) \circ
\wp^-(q_{1/2},q'_{1/2}) \circ \wp^+(y_{1/2},q_{1/2})$.
Let
\begin{align*}
h & = \wp^+(q_{1/2},y_{1/2}) \circ \wp^-(q'_{1/2},q_{1/2}) \circ
\wp^+_{fake}(z_{1/2}, q'_{1/2})\circ \id(q'_{1/2},z_{1/2}) \circ \wp^-(q_{1/2},q'_{1/2})
\\
& = \wtilde{\phi}^{-1}\circ \id(q'_{1/2},z_{1/2}) \circ \wp^-(q_{1/2},q'_{1/2})
\colon \cW^{u}[q_{1/2}] \to \cW^{u}[q_{1/2}]
\end{align*}
where $\id(q'_{1/2},z_{1/2})\colon \cW^u[q'_{1/2}]\to \cW^u[z_{1/2}]$ is the tautological identification.
Then we have
\begin{align*}
\tilde{\phi}^{-1}(U^+[q'_{1/2}]) 
= h \left(U^+[q_{1/2}]\right)
\end{align*}
since $\wp^-(q_{1/2},q'_{1/2})(U^+[q_{1/2}])=U^{+}[q'_{1/2}]$.

Let $H: L\cW^u(q_{1/2}) \to L\cW^u(q_{1/2})$ denote $L_{q_{1/2}} \circ
h \circ L_{q_{1/2}}^{-1}$. Recall that for $a,b \in \cW^u[c]$,
\begin{displaymath}
L_b \circ \wp^+(a,b) \circ L_a^{-1} = P^+(a,b), 
\end{displaymath}
and
\begin{displaymath}
L_b \circ L_a^{-1} = P^{L\cW^u}_{GM}(a,b). 
\end{displaymath}
Then we can rewrite the linearization of $h$ as:
\begin{displaymath}
H = P_{GM}^{L\cW^u}( y_{1/2},q_{1/2}) \circ P^+(q_{1/2},y_{1/2}) \circ
P^-(q_{1/2},q'_{1/2}) \circ P^*(q'_{1/2},q'_{1/2}) \circ
P^-(q_{1/2},q'_{1/2}),
\end{displaymath}
where $P^*(q'_{1/2},q'_{1/2}) =
P^+_{fake}(z_{1/2},q'_{1/2})
P_{GM}^{L\cW^u}(q'_{1/2}, z_{1/2})
$.
We have that $P^*(q'_{1/2},q'_{1/2})$ is
factorizable by \autoref{rmk:factorizing_holonomy} and \autoref{lemma:Pfake:factorizable}. 
Therefore $H$ is factorizable by
\autoref{prop:factorizability_in_smooth_global_charts}\autoref{item:prop:factorizability_in_smooth_global_charts:item3} and
\autoref{lemma:unstable:flag:factorizable}.

Let \index{$G$@$\cG(x)$}$\cG(x)$ denote the image of
$\bbG^{ssr}(\cW^u[x])$ under the linearization map $L_x$. Then $\cG(x)$
is a subgroup of the general linear group of $L\cW^u(x)$. Let $\Psi_x:
\cG(x) \to L\cC(x)$ be the unique map such that $\Psi_x(L_x f) = L (f U^+[x])$, for any $f\in \bbG^{ssr}(\cW^u[x])$.
Then, 
\begin{displaymath}
Z_{1/2}(q,q',u,\ell) = \Psi_{q_{1/2}} (H). 
\end{displaymath}
The map $\Psi_{q_{1/2}}$ is clearly factorizable (as it only depends
on data at $q_{1/2}$). Therefore, since $H$ is factorizable, and in
view of \autoref{cor:assembly_of_factorizable_maps}, $Z$ is
factorizable. 
 
Therefore there exist $F^s$, $F^u$, $\crA$ such that
for any $\delta>0$ there exists a compact set $K$ of measure at least
$1-\delta$, constants $\ell_0=\ell_0(K,N)$ and $C=C(K,N)$, such that
provided $q,q',q_{1/2},q_{1/2}',y_{1/2},q_1,uq_1$ all belong to $K$,
and $\ell > \ell_0$, we have:
\begin{equation}
\label{eq:Z:factorizable}
		\norm{
		Z_{1/2}(q,q',y,\ell) - 
		{\crA}(q_{1/2})( 
		g_{\ell/2} F_q^s(q'), 
		g_{-\ell/2} F_{q_1}^u(uq_1) ) 
		}
		\leq
		C e^{-N\ell}
\end{equation}
Note that $Z_{1/2}(q,q,u,\ell) \in L\cC^\tau(q_{1/2})$. Therefore, there
exists $v_0 \in L\cC^\tau(q_{1/2})$ such that
\begin{displaymath}
  \| \crA(q_{1/2})( g_{\ell/2} F_q^s(q), g_{-\ell/2} F_{q_1}^u(uq_1)) 
  -v_0\| \le e^{-N \ell}.  
\end{displaymath}
We may find a linear map $L_0:
L\cC(q_{1/2}) \to L\cC(q_{1/2})$ with  $\|L_0 - I \| =
O(e^{-N \ell})$ such that
\begin{displaymath}
L_0 \circ \crA(q_{1/2})( g_{\ell/2} F_q^s(q), g_{-\ell/2} F_{q_1}^u(uq_1)) 
   = v_0. 
\end{displaymath}
Then, (\ref{eq:Z:factorizable}) still holds with $\crA$
replaced by $L_0 \circ \crA$. We now define
\begin{displaymath}
\cA_{1/2}(q,u, \ell) \colon V^s(q) \to 
L\cC(q_{1/2})
\end{displaymath}
as
\begin{displaymath}
		\cA_{1/2}(q,u,\ell) v = L_0 \circ \crA(q_{1/2})\left(
		g_{\ell/2}(v), g_{-\ell/2}F^u_q(uq_1)
		\right)
\end{displaymath}
for $v\in V^s(q)$ where $u_*$ is as in \autoref{def:u_star}.
Then $\cA_{1/2}(\cdot, \cdot, \cdot)$ is a linear
map. The estimate (\ref{eq:thm:3varA:halftime:factorizable}) (for a
different $N$) follows
from (\ref{eq:Z:factorizable}) and \autoref{lemma:ustar:bounded}.
Also (\ref{eq:3varcA:q:in:LC}) holds by construction. 
\end{proof}



\section{The \texorpdfstring{$\cA$}{curly A} operators}
	\label{sec:the_curly_a_operators}


\subsection{A lower bound for \texorpdfstring{$\cA_{1/2}$}{curly A 1/2}}
	\label{ssec:a_lower_bound_for_curly_a}

For now on, we work with the cocycles $V^s$ and map $\cA_{1/2}(-,-,-)$ from \autoref{theorem:3varA:haltime:factorizable}, for some $N$ sufficiently large to be fixed in \autoref{eq:choice:of:N}.

\subsubsection{The subspace $V^s_*(q)$}
For a.e. $x \in X$, let \index{$\nu_x^s$}$\nu_x^s$ denote the leaf-wise measure of $\nu$ along $\cW^s[x]$,
normalized so that $\nu^s_x(\gB_0^-[x]) = 1$. 

For a.e $q \in X_0$, let \index{$V^s_*$}$V^s_*(q) 
\subset V^s(q)$ denote the smallest subspace
containing the support of the measure $(F^s_q)_*\hat{\nu}^s_q$, where
$\hat{\nu}^s_q$ is the restriction of $\nu^s_q$ to $\gB_0^-[q]$.
By the Markov property of the partition $\gB_0^-$, we see that
for $t > 0$ we have that
$g_{-t} V^s_*(q) \subset V^s_*(g_{-t} q)$.
Therefore by ergodicity of
$g_{-t}$, the dimension of $V^s_*(q)$ is almost everywhere constant,
and $g_t V^s_*(q) = V^s_*(g_t q)$ almost everywhere. 

\subsubsection{The norm $\| \cdot \|_{V^s_*}$}
\label{sec:subsec:norm:Vs:star}
Let \index{$V^s_-(q)$}$V^s_-(q) \subset V^s_*(q)$ denote the direct sum of the
Lyapunov subspaces of $V^s_*(q)$ with negative Lyapunov
exponent. Then, by \autoref{prop:contraction_in_v},
\begin{displaymath}
V^s_*(q) = \reals F^s_q(q) \oplus V^s_-(q).
\end{displaymath}
By a standard Lyapunov norm construction, we can
measurably pick an inner product \index{$\norm$@$\langle \cdot, \cdot \rangle_{V^s_*(q)}$}$\langle \cdot, \cdot \rangle_{V^s_*(q)}$ on
$V^s_*(q)$, such that the following hold:
\begin{itemize}
\item[{\rm (a)}] $\reals F^s_q(q)$ and $V^s_-(q)$ are orthogonal.
\item[{\rm (b)}] $\|F^s_q(q) \|_{V^s_*(q)} = 1$, where $\| v
  \|_{V^s_*(q)} = \langle v,v \rangle_{V^s_*(q)}^{1/2}$.
\item[{\rm (c)}] There exists \index{$\kappa_1(N)$}$\kappa_1(N) > 1$ depending only on $N$
  and the Lyapunov spectrum such that  for $F \in V^s_-(q)$ and $t > 0$,
\begin{align}
\label{eqn:def_kappa_1_Vs}
e^{-\kappa_1(N) t} \| F \|_{V^s_*(q)} \le \| g_t F
\|_{V^s_*(g_t q)} \le e^{- \kappa_1(N)^{-1} t} \| F \|_{V^s_*(q)}. 
\end{align}
\end{itemize}
Note that $F^s_q(q)\neq 0$ on a set of full measure because $\cA_{1/2}(q,u,\ell)F^s_q(q)\neq 0$.

\begin{lemma}
\label{lemma:away:from:subspace}
Let \index{$\tilde{\nu}_q$}$\tilde{\nu}_q$ denote the measure on $V^s_-(q)$ which is the
pushforward of restriction of the conditional measure $\nu^s_q$ to 
$\cW^s_{loc}[q]$ by the map $q' \to F_q(q') - F_q(q)$. Then for every
$\delta > 0$ there exists a compact set $K$ with $\nu(K) \ge 1-
\delta$ and constants $C(\delta)$ and $\rho(\delta)> 0$ such that 
for any $q \in K$ and any subspace $V' \subset V^s_-(q)$,
\begin{displaymath}
\tilde{\nu}_q(B(q,C(\delta)\setminus Nbhd_{\rho(\delta)}(V')) \ge 1-\delta.
\end{displaymath}
\end{lemma}

\begin{proof}
This follows immediately from the definitions (cf. \cite[Lemma 5.4,
Lemma 5.5]{EskinMirzakhani_Invariant-and-stationary-measures-for-the-rm-SL2Bbb-R-action-on-moduli-space}).
\end{proof}

\begin{corollary}
\label{cor:can:choose:basis}
For every $\delta > 0$ there exists 
a compact subset $K \subset X$
with $\nu(K) > 1-\delta$ such that for any $q \in K$ and any subset
$E$ with $\tilde{\nu}_x(E) > (1-\delta)\tilde{\nu}_q(V_-^s)$ there
exists a basis $\{v_1, \dots, v_n\}$ of $V^s_-$ such that for all $1 \le j
\le n$ we have $v_j \in E$, $\|v_j\| \le C(\delta)$ and the distance
between $v_j$ and the span of $\{v_1, \dots, v_{j-1}$ is at least
$\rho(\delta)$.    
\end{corollary}

\begin{proof}
This follows by repeated application of
\autoref{lemma:away:from:subspace}. 
\end{proof}

\subsubsection{Bounds on cocycles}
	\label{sssec:bounds_on_cocycles}
Let \index{$\kappa$}$\kappa > 1$ be such that we have the following bounds for $t\geq 0$:
\begin{align}
	\label{eqn:du_lower_growth_bound}
	\begin{split}
	e^{\kappa^{-1}t}d^u(x,y) & \le  d^u(g_t x, g_t y)\\ 
	e^{\kappa^{-1}t}d^\cL_z(x,y) & \le  d^\cL_{g_tz}(g_t x, g_t y) \le e^{\kappa t} d^\cL_{z}(x,y)\\
	\end{split}
\end{align}
where $x,y\in \cW^u[z]$ in the first line, and $x,y\in \cW^u_{loc}[z]\cap \cL_{k,\ve}[z]$ in the second.
Additionally, we require that for any $c\in L\cC(x)$, or any $c\in \bbH(x)$ (where $\bbH(x):=L\cC(x)/L\cC^{\tau}(x)$ is defined in \autoref{eq:bold_H_definition}) we have for $t\geq 0$:
\begin{align}
	\label{eqn:kappa_growth_good_norms}
	\begin{split}
	e^{\kappa^{-1}t}\norm{v} \le & \norm{g_t v} \le e^{\kappa t}\norm{v}\\
	e^{\kappa^{-1}t}\norm{c} \le & \norm{g_t c} \le e^{\kappa t}\norm{c}.
	\end{split}
\end{align}
For the bounds on linear cocycles in \autoref{eqn:kappa_growth_good_norms}, the corresponding norms are defined in \autoref{prop:good_norms}.
We choose $\kappa$ large enough so that the same bounds hold in the stable directions.
Note also that all of the cocycles listed above depend only on the Lyapunov spectrum of the measure and not the degree of approximation $N$.

\subsubsection{Defining $\cA_1$}
	\label{sssec:defining_curly_a_1}
Let $\cA_{1/2}(\cdot,\cdot,\cdot)$ be as in \autoref{theorem:3varA:haltime:factorizable}, recall that
\begin{align*}
	\cA_{1/2}(q,u,\ell)& \colon V^s(q)\to L\cC(q_{1/2})
\intertext{and set\index{$A$@${\cA}_{1}(q_1,u,\ell)$}}
	{\cA}_{1}(q_1,u,\ell) = g_{\ell/2} \circ \cA_{1/2}(g_{-\ell} q_1,u,\ell)
	& \colon
	V^s(q) \to L\cC(q_1)
\end{align*}
Note that with our conventions (also later on for $\cA_2$) we have that $\cA_i$ takes values in $L\cC(q_i)$.
We also define 
\begin{equation}
	\label{eq:Z_1_Z12_def}
	\index{$Z_1(q,q',u,\ell)$}Z_1(q,q',u,\ell) = g_{\ell/2} \circ Z_{1/2}(q,q',u,\ell).
\end{equation}
%
%

\begin{lemma}[Almost-equivariance of $\cA_{1/2},\cA_{1}$]
\label{lemma:3varA:almost:equivariant}
For every $N >0$ there exists $N' >0$ with $N' \to
\infty$ as $N \to \infty$ such that the following holds:
for every $\delta > 0$ there exists a compact set $K$ with $\nu(K) >
1-\delta$ such that if 
$q \in K$, $g_{\ell/2} q \in K$, $g_{-\ell/2} u g_\ell q \in K$, 
$0.9\ell < s < \ell$, $\hat{q} :=
g_{\ell-s} q \in K$, $q_1:=g_{\ell}q\in K$, and $g_{s/2} \hat{q} \in K$, $g_{-s/2} u g_s \hat{q} \in K$
we have, for $\ell > \ell_0(\delta)$, 
\begin{align}
\label{eq:lemma:3varZ:almost:equivariant}
\| Z_1(q,q',u,\ell)- Z_1(\hat q,\hat q',u,s)\| &\le e^{-N' \ell}
\\
\label{eq:lemma:3varA_12:almost:equivariant}
\|\cA_{1/2}(q,u,\ell) - g_{-(\ell-s)/2} \circ \cA_{1/2}(\hat{q}, u, s) \circ
g_{\ell-s} \| &\le e^{-N' \ell}\\
\label{eq:lemma:3varA_1:almost:equivariant}
\|\cA_{1}(q_1,u,\ell) - \cA_{1}(q_1, u, s) \circ
g_{\ell-s} \| &\le e^{-N' \ell}.
\end{align}
In the last two bounds the norm $\| \cdot \|$
refers to the operator norm with respect to the norms $\| \cdot
\|_{V^s_*(q)}$ in the domain and $\| \cdot \|_{L\cC(g_{-\ell/2} u
  g_\ell q)}$ in the codomain. 
\end{lemma}

\begin{proof}
We treat each of the above bounds separately.

\noindent \textbf{Proof of \autoref{eq:lemma:3varZ:almost:equivariant}}
We first argue that $$\| Z_{1/2}(q,q',u,\ell)-g_{-(\ell-s)/2} Z_{1/2}(\hat q,\hat q',u,s)\| \le e^{N'_1\ell}
$$
with $N'_1\to +\infty$ as $N\to +\infty$.
Assuming this bound, note that since $Z_{1}$ differs from $Z_{1/2}$ by an application of $g_{\ell/2}$ (by \autoref{eq:Z_1_Z12_def}), and the cocycle in which $Z$ takes values, namely $L\cC$ is independent of $N$ and hence has fixed Lyapunov spectrum, the desired bound follows with $N'=N''-\Lambda_{L\cC}$ where $\Lambda_{L\cC}$ is the largest Lyapunov exponent of $L\cC$.

Next, recall from \autoref{eqn:Z_12_def} that
\[
	Z_{1/2}(q,q',u, \ell) = L\left(\td \phi\inv\left(  U^+[q'_{1/2}]\right)\right)
\]
where $L$ is a linearization map and $\wtilde{\phi}$ is the interpolation map from \autoref{sssec:defining_the_interpolation_map}.
Note that both the linearization $L$ and $U^+[q']$ are $g_t$-equivariant, therefore in order to prove
\begin{align}
	\label{eqn:Z_12_almost_equivariance_ineq}
	\norm{Z_{1/2}(q,q',u, \ell) - 
	g_{-(\ell-s)/2} \circ Z_{1/2}(\hat{q},\hat{q}',u,s)}
	\leq e^{-N_2'\ell}
\end{align}
it suffices to establish the bound (where we now explicitly include the dependence of $\wtilde{\phi}$ on $q,q',u,\ell$):
\[
	\norm{L\wtilde{\phi}^{-1}(q,q',u,\ell)-g_{-(\ell-s)/2} L\wtilde{\phi}^{-1}(\hat{q},\hat{q}',u,s)g_{(\ell-s)/2}}
	\leq e^{-N_2''\ell}
\]
with $N_2''\to +\infty$ as $N\to +\infty$, and then use again that the cocycle in which the functions take values is independent of $N$.
For notational simplicity, we will prove the same bound but for $\wtilde{\phi}$ instead of $\wtilde{\phi}^{-1}$.
Recall again from \autoref{sssec:defining_the_interpolation_map} that
\[
	\wtilde{\phi}(q,q',u,\ell) =
	\wp^+_{fake}(q_{1/2}',z_{1/2})\circ
	\wp^-(q_{1/2},q_{1/2}')\circ
	\wp^{+}(y_{1/2},q_{1/2})
\]
where $\wp^+_{fake}(q_{1/2}',z_{1/2})$ is the map from \autoref{thm:fake_holonomies_and_ssr_map}.
Again $\wp^{\pm}$ are $g_t$-equivariant, and the norms $\norm{L\wp^{-}(q_{1/2},q_{1/2}')}$ and $\norm{L\wp^+\left(y_{1/2},q_{1/2}\right)}$ are uniformly bounded (since $q_{1/2},q'_{1/2},y_{1/2}$ are biregular), so it suffices to prove that
\[
	\norm{L\wp^{+}_{fake}\left(q_{1/2}',z_{1/2}\right) - g_{-(\ell-s)/2}\circ L\wp^{+}_{fake}\left(\hat{q}_{1/2}',\hat{z}_{1/2}\right)g_{(\ell-s)/2}}
	\leq e^{-N'_3\ell}.
\]
Since $g_{-(\ell-s)/2}\hat{q}_{1/2}'=q_{1/2}'$ the second term above can be rewritten as $L\wp^{+}_{fake,\ell-s}\left(q_{1/2}',g_{-(\ell-s)/2}\hat{z}_{1/2}\right)$
where $\wp^{+}_{fake,\ell-s}$ is defined by taking the fake forward flag at $\hat{z}_{1/2}$ from \autoref{thm:fake_holonomies_and_ssr_map} and applying $g_{-(\ell-s)/2}$ to it.
We can thus rewrite the difference above as:
\[
	L\wp^{+}_{fake}\left(q_{1/2}',z_{1/2}\right)\left(
	\id - L\wp^{+}_{fake,fake,\ell-s}(z_{1/2},g_{-(\ell-s)/2}\hat{z}_{1/2})\right)
\]
where $\wp^{+}_{fake,fake,\ell-s}$ is defined using the fake forward stable flags at $z_{1/2}$ and $g_{-(\ell-s)/2}\hat{z}_{1/2}$, as well as the backwards flags (defined and equivariant at both points) and the holonomy on the associated graded associated to the backwards flag.

However, \autoref{thm:fake_holonomies_and_ssr_map}\autoref{item:P+_fake_z_zhat_comparison}
implies that $\norm{\id-L\wp^+_{fake,fake,\ell-s}}\leq e^{-N_5'\ell}$ and concludes the proof.

\noindent \textbf{Proof of \autoref{eq:lemma:3varA_1:almost:equivariant} and \autoref{eq:lemma:3varA_12:almost:equivariant}}
Note that \autoref{eq:lemma:3varA_1:almost:equivariant} is implied by \autoref{eq:lemma:3varA_12:almost:equivariant}, with exponent $N'-\kappa$ going to infinity as $N$ goes to infinity.
Indeed the terms inside the norm differ only by an application of $g_{\ell/2}$, which has norm bounded by $e^{\kappa \ell}$ and $\kappa$ is independent of $N$ (see \autoref{eqn:kappa_growth_good_norms}).

To prove \autoref{eq:lemma:3varA_12:almost:equivariant}, let $K'$ be
the compact set from \autoref{theorem:3varA:haltime:factorizable} but
of measure at least $1-\delta^2$. We may also assume that
\autoref{cor:can:choose:basis} holds for $K'$ and $\delta$
replaced by $\delta^2$.


Let now $K\subseteq K'$ be its compact subset of measure at least $1-2\delta$ of points of density for the stable conditionals $\nu^s_q\vert_{\gB_0^-}$ provided by \autoref{lemma:gB:vitali}, and such that $C_{\delta}$ is uniformly bounded on $K$ by some constant $C(\delta)$.
In particular if $q\in K$ then $\nu^s_q(K'\cap \gB_t^-[q])\geq (1-\delta)\nu^s_q(\gB_t^-[q])$ for $t\geq 0$.

By assumption $q,\hat{q}=g_{\ell-s}q\in K$ so let $q' \in \gB_0^-[q]\cap K'$ and $\hat{q}' = g_{\ell-s} q'\in \cW^s[\hat{q}]$ with $\hat{q}'\in K$.
This can be arranged for $q'\in \gB_0^-[q]$ of density at least $1-2\delta$.

Recall that $q_1 = g_{\ell} q = g_{s} \hat{q}$, and 
let $Z_{1/2}(\cdot, \cdot, \cdot, \cdot)$ be as in
\autoref{sec:subsec:def:S:LC}.
Then, we proved earlier in \autoref{eqn:Z_12_almost_equivariance_ineq} that:
\begin{displaymath}
\norm{Z_{1/2}(q,q',u,\ell) - 
g_{-(\ell-s)/2} \circ Z_{1/2}(\hat{q},\hat{q}',u,s)} \leq e^{-N_2''\ell}
\end{displaymath}
Therefore, by \autoref{theorem:3varA:haltime:factorizable} and
\S\ref{sec:subsec:norm:Vs:star}, since $q'\in K'$ and $g_{\ell-s}q'\in K'$ we have:
\begin{displaymath}
\|\cA_{1/2}(q,u,\ell) F^s_q(q') - g_{-(\ell-s)/2} \circ \cA_{1/2}(\hat{q},u,s)
F^s_{\hat{q}} (g_{\ell-s} q') \| \le e^{-N' \ell}, 
\end{displaymath}
where $N' \to \infty$ as $N \to \infty$. 
By the equivariance of $F^s$ we have $F^s_{\hat{q}}(g_{\ell-s} q') =
g_{\ell-s} F^s_q(q')$. Therefore,
\begin{multline}
\label{eq:var:AF:almost:equivariant}
  \|\cA_{1/2}(q,u,\ell) F^s_q(q') - g_{-(\ell-s)/2} \circ \cA_{1/2}(\hat{q},u,s)
\circ g_{\ell-s} F^s_q(q') \| \\ \le e^{-N' \ell}. 
\end{multline}
We now choose a basis
for $V^s_-(q)$ as in \autoref{cor:can:choose:basis} consisting
of vectors for which (\ref{eq:var:AF:almost:equivariant}) holds. Then,
\begin{displaymath}
\|\cA_{1/2}(q,u,\ell) - g_{-(\ell-s)/2} \circ \cA_{1/2}(\hat{q}, u, s) \circ
g_{\ell-s} \| \le C(\delta) e^{-N' \ell}.
\end{displaymath}
After replacing $N'$ by $N'+1$ and choosing $\ell$ large enough, we
see that \autoref{eq:lemma:3varA_12:almost:equivariant} holds. 
\end{proof}

\subsubsection{A lower bound on $\cA_{1/2}, \cA_1$}
\label{sec:subsec:lower:bound:on:3cA}
Let now $N$ be as in \autoref{theorem:3varA:haltime:factorizable}.
Let $\cB^X$ denote the set of pairs $(q_1,u)$ where $q_1 \in X$ and $u
\in \cB_0(q_1)$.  
\begin{lemma}[Choice of \index{$\alpha_2$}$\alpha_2$]
	\label{lemma:lower:bound:on:3cA}
Suppose the family $q \to
U^+(q)$ satisfies the QNI condition (see
\autoref{def:QNI}).
Then, there exists $\alpha_2 > 0$ and $N_0$ (depending only on the family $q \to
U^+(q)$) such that for any $N > N_0$ there exists a measurable
function $\ell_0: \cB^X\cross \reals
\to \reals^+$ such that for almost every $q_1 \in X$
and $u \in \cB_0(q_1)$ and any $\ell > \ell_0(q_1,u,N)$ there exists a nonzero $F_- \in  
V^s_-(g_{-\ell} q_1)$ such that
\begin{align}
\label{eq:lower:bound:on:3cA_12}
\|\cA_{1/2}(g_{-\ell} q_1, u, \ell) F_- \|_{L\cC(q_{1/2})} &\ge e^{-\alpha_2 \ell}
\|F_-\|_{V^s_-(g_{-\ell} q_1)}\\
\label{eq:lower:bound:on:3cA_1}
\|\cA_{1}(q_1, u, \ell) F_- \|_{L\cC(q_{1})} &\ge e^{-\alpha_2 \ell}
\|F_-\|_{V^s_-(g_{-\ell} q_1)}
\end{align}
In this lemma, $N$ is the parameter of
\autoref{def:half_time_factorizable}. 
\end{lemma}


\begin{proof}
As in the proof of \autoref{lemma:3varA:almost:equivariant}, the bound in \autoref{eq:lower:bound:on:3cA_1} follows from that in \autoref{eq:lower:bound:on:3cA_12}, with the same $\alpha_2$, by applying $g_{\ell/2}$ and noting again that on the cocycle $L\cC$ the forward flow is expanding.

Let $K$, $\ell_0$, $\alpha_0$ be as in \autoref{def:QNI}. 
Suppose first that $\ell > \ell_0$ is such that $q_{1/2} \equiv
g_{-\ell/2} q_1 \in K$, and also $q \equiv g_{-\ell} q_1 \in
K$.
Suppose also that $u \in \cB_0(q_1)$ is such that $y_{1/2} \equiv g_{-\ell/2} u q_1 \in S$,
where $S = S(q_{1/2},\ell)$ is as in \autoref{def:QNI}. 

Let $\alpha$ be as in
\autoref{theorem:interpolation:does:not:move:points}.
Note that if \autoref{theorem:interpolation:does:not:move:points} holds with some $\alpha>0$ then it also holds for any smaller $\alpha>0$.
Similarly, if the QNI condition \autoref{def:QNI} holds for some $\alpha_0$, then it holds for any larger $\alpha_0$ as well.

Suppose then that $\alpha_0\gg \alpha$, by a uniform factor that only depends on the Lyapunov spectrum.
Then since $z_{1/2}\in \cW^{cs}[y_{1/2}]$, by the QNI condition we have that
\[
	d^Q\left(U^+[q'_{1/2}],z_{1/2}\right)\geq e^{-\alpha_0\ell}.
\]
By selecting the compact set $K$ accordingly, we have (since $q'_{1/2}\in K$)
\[
	d^\cL_{q'_{1/2}}\left(U^+[q'_{1/2}],z_{1/2}\right) \geq \tfrac{1}{C_1(\delta)} d^Q\left(U^+[q'_{1/2}],z_{1/2}\right)\geq 
	\tfrac{1}{C_1(\delta)}e^{-\alpha_0\ell}.
\]
Note that the two sets are in the same unstable.
We would like to obtain an estimate of this form, but after applying $\wtilde{\phi}^{-1}$.
To do so, consider $q'_t:=g_t q'_{1/2}, z_t:=g_t z_{1/2}, y_t:=g_t y_{1/2}$ and in the time interval $(\alpha_0-\tfrac \alpha{10}) [\tfrac{\ell}\kappa,\kappa\ell]$ choose $t$ such that $q'_t,y_t\in K_1$ for a compact set $K_1$ (with further standard dynamical assumptions) and such that we also have:
\[
	e^{-\tfrac \alpha{100} \alpha\ell} \geq d^{\cL}_{q'_t}\left(U^+\left[q'_{t}\right], z_{t}\right) \geq e^{-\tfrac\alpha{10}\ell}
\]
using \autoref{eqn:du_lower_growth_bound}.
Using \autoref{theorem:interpolation:does:not:move:points}, that $\phi_t\left(y_{t}\right)=z_t$, and the triangle inequality we have that
\[
	d^\cL_{y_t}\left(y_t, \phi_t^{-1}\left(U^+[q'_t]\right)\right)
	\geq C_3(\delta)
	e^{-\tfrac\alpha{10} \ell}
	-
	2 C_2(\delta) e^{-\alpha \ell}
	\geq e^{-\tfrac\alpha{5}\ell}
\]
Explicitly, setting $m_t\in\phi_t^{-1} \left(U^+[q'_t]\right)$ to be a minimizer of the distance to $y_t$, we have that $d^Q(m_t,\phi_t(m_t))\leq C_2(\delta)e^{-\alpha \ell}$ and $d^Q(y_t,\phi_t(y_t))\leq C_2(\delta)e^{-\alpha \ell}$ (and $\phi_t(y_t)=z_t$) by \autoref{theorem:interpolation:does:not:move:points}.
We finally pass between $d^\cL_{y_t}$ and $d^Q$ at the cost of another $C_3(\delta)$ to obtain the above inequality.
Applying now $g_{-t}$ we find that
\[
	d^\cL_{y_{1/2}}\left(y_{1/2},\wtilde{\phi}^{-1}\left(U^+[q'_{1/2}]\right)\right)
	\geq
	e^{-\kappa t} e^{\tfrac\alpha{5}\ell}
	\geq e^{-\left[\kappa^2(\alpha_0-\tfrac \alpha{10})+\tfrac\alpha{5}\right] \ell}.
\]
Finally, again using that $y_{1/2}\in K$ we obtain that
\[
	d^Q\left(y_{1/2},\wtilde{\phi}^{-1}\left(U^+[q'_{1/2}]\right)\right)
	\geq
	\tfrac{1}{C_1(\delta)}
	 e^{-\left[\kappa^2(\alpha_0-\tfrac \alpha{10})+\tfrac\alpha{5}\right] \ell}
	 \geq
	 \tfrac{1}{C_1(\delta)}
	 {e^{-\alpha_2'\ell}}
\]
where $\alpha_2':=\kappa^2(\alpha_0-\tfrac \alpha{10})+\tfrac\alpha{5}$ only depends on the Lyapunov spectrum and the QNI condition.
Since $y_{1/2}\in U^+[q_{1/2}]$, it follows that
\[
	d^{q_{1/2}}_{\cH,loc}(U^+[q_{1/2}], \wtilde{\phi}^{-1}\left(U^+[q'_{1/2}]\right))
	\geq 
	\tfrac{1}{C_1'(\delta)}
	 {e^{-\alpha_2'\ell}}
\]
where $d^q_{\cH,loc}$ is the local Hausdorff distance defined in \autoref{eq:Hausdorff_distance_def}.
Let $\gamma^+(q_{1/2})\in L\cC(q_{1/2})$ be the vector corresponding to $U^+[q_{1/2}]$ and similarly for $\gamma^+(z_{1/2})$ corresponding to $\tilde{\phi}^{-1}\left(U^+[q'_{1/2}]\right)$.
Using now \autoref{prop:relation_to_hausdorff_distance}, we have that
\begin{displaymath}
\norm{\gamma^{+}(z_{1/2})-\gamma^{+}(q_{1/2})}
\geq
\frac{1}{C_2(\delta)} e^{-\alpha_2' \ell}. 
\end{displaymath}


Let us now write $F_q(q'):=F_q(q) + F_{q,-}(q')$ with $F_{q,-}(q')\in V^s_-(q)$ and the coefficient of $F_q(q)$ is $1$ by \autoref{prop:contraction_in_v}\ref{item:holonomy_invariance_of_iota_s}.
By \autoref{theorem:3varA:haltime:factorizable} we have
\[
	\norm{\cA_{1/2}(q,u,l)F_q(q') - \gamma^+(z_{1/2})} \leq C_3(\delta) e^{-N\ell}
	\text{ and }
	\cA_{1/2}(q,u,l)F_q(q) = \gamma^+(q_{1/2}).
\]
We conclude that for all $q' \in S'\cap K$ we have:
\begin{multline*}
\norm{\cA_{1/2}(q,u,\ell) F_{q,-}(q')}
=
\norm{\cA_{1/2}(q,u,\ell) \big(F_q(q)-F_q(q')\big)}\geq\\
\geq \norm{\gamma^+(q_{1/2})-\gamma^+(z_{1/2})} - 
C_1(\delta) e^{-N \ell}
\ge \tfrac{1}{C_2(\delta)} e^{-\alpha_2' \ell} -
C_1(\delta) e^{-N \ell}
\end{multline*}
If we choose $N > 10 \alpha_2'$ then for $\ell$ large enough we get
\begin{displaymath}
\|\cA_{1/2}(q,u,\ell) F_{q,-}(q')\| \ge 
\tfrac{1}{C_2(\delta)} e^{-2 \alpha_2' \ell}.
\end{displaymath}
Next, note that $\norm{F_{q,-}(q')}\leq C_3(\delta)e^{-\beta(N)\ell}$ by \autoref{prop:contraction_in_v}(ii), for some $\beta$ that depends on $N$.
Note that in fact and upper bound of $O_\delta(1)$ suffices.
Since $S'\cap K$ contains most of the measure of $\gB^-_0[q]$, it follows
that for $u$ such that $g_{-\ell/2} u q_1 \in S(q_{1/2},\ell)$
\begin{equation}
\label{eq:3cA:tmp:norm:lower:bound}
\norm{\cA_{1/2}(q,u,\ell)\vert_{V^s_-(q)}}
\ge 
\tfrac{1}{C(\delta)} e^{-2 \alpha_2' \ell}.
\end{equation}
By the Birkhoff ergodic theorem, we may assume that for $q_1 \in K$
and $L > L_0(\delta)$,
\begin{displaymath}
|\{ \ell \in [0,L] \st g_{-\ell} q_1 \in K \}| \ge (1 - c_1(\delta))L, 
\end{displaymath}
where $c_1(\delta) \to 0$ as $\delta \to
0$. Then,
\begin{displaymath}
|\{(u,\ell) \in \cB_0(q_1) \cross [0,L] \st g_{-\ell q_1} \in K \} \text{
    and } g_{-\ell/2} u q_1 \in S(q_1,\ell) \}| (1 - c_1(\delta))L, 
\end{displaymath}
where $c_2(\delta) \to 0$ as $\delta \to
0$. (The measure $| \cdot |$ on the left hand side of the above
equation is the product of the measures $| \cdot |$ on $\cB_0(q_1)$
and on $[0,L]$). 
Then, by Fubini's theorem (applied to $\cB_0[q_1] \cross [0,L]$) there
exists $Q_1'(q_1) \subset \cB_0(q_1)$ with $|Q_1'(q_1)| \ge
(1-c_2(\delta)) |\cB_0(q_1)|$ such that for $u \in Q_1'(q_1)$ and $L$
sufficiently large depending on $\delta$,  
\begin{displaymath}
|\{ \ell \in [0,L] \st g_{-\ell} q_1 \in K \text{
    and } g_{-\ell/2} u q_1 \in S(q_1,\ell) \}| \ge (1 - c_3(\delta))L, 
\end{displaymath}
where $c_2(\delta)$ and $c_3(\delta)$ tend to $0$ as $\delta \to 0$. 

Now suppose $u \in Q_1'(q_1)$, and suppose $\ell$ is sufficiently
large depending on $\delta$. We can choose $0 < s <  \ell$ such that
$\ell - s < c_4(\delta) \ell$ where $c_4(\delta) \to 0$ as $\delta
\to 0$ such that $g_{-s} q_1 \in K$ and $g_{-s/2} u q_1 \in
S(q_1,s)$. 
Then, by (\ref{eq:3cA:tmp:norm:lower:bound}),
\begin{displaymath}
\|\cA_{1/2}(g_{-s} q_1, u, s) \| 
\ge 
\tfrac{1}{C(\delta)} e^{-2 \alpha_2' s}.
\end{displaymath}
We choose $\delta > 0$ sufficiently small so that
$c_4(\delta) (\kappa+\kappa_1(N)) < \kappa \alpha_2'$.
By \S\ref{sec:subsec:norm:Vs:star} (c) and
\autoref{prop:good_norms} (iv) applied to $L\cC$
\begin{multline*}
\|g_{-(\ell-s)/2} \cA_{1/2}(g_{-s} q_1, u, s) \circ g_{\ell-s} \| \ge\\
\geq 
\tfrac{1}{C(\delta)} e^{-(2 \kappa\alpha_2' + c_4(\delta) (\kappa +
  \kappa_1(N))) \ell} 
  \ge \tfrac{1}{C(\delta)} e^{-3 \kappa\alpha_2' \ell}.
\end{multline*}
Then, by \autoref{lemma:3varA:almost:equivariant},
\begin{displaymath}
\|\cA_{1/2}(g_{-\ell} q_1, u, \ell) \| \ge \tfrac{1}{C(\delta)} e^{-3 \kappa\alpha_2' \ell} - e^{-N' \ell}.
\end{displaymath}
Since $\alpha_1$, $\kappa$ and $\kappa_1$ are independent of $\delta$,
$c_4(\delta) \to 0$ as $\delta \to 0$ 
and $N' \to \infty$ as $N \to \infty$, we can choose $\delta$ small
enough (dependent on $N$) and $N_0$ large enough so that for $\ell > \ell_0'(\delta,N)$,
(\ref{eq:lower:bound:on:3cA_12}) holds with $\alpha_2 = 4 \kappa\alpha_2'$. 

As $\delta$ decreases to $0$, the sets $\{ (q_1,u) \st q_1 \in K =
K(\delta),  u \in Q_1'(q_1) \}$ exhaust $\cB^X$. Therefore, we can
define $\ell_0(q_1,u,N)$ to be $\ell_0'(\delta,N)$ where $\delta$ is the
largest number such that $q_1 \in K(\delta)$, $u \in Q_1'(\delta)$. 
\end{proof}

\subsubsection{Choice of $\alpha_3$.}
\label{sec:subsec:choice:of:alpha3}
We choose \index{$\alpha_3$}$\alpha_3 > 0$ such that
\begin{displaymath}
\alpha_3 > 4 \kappa \alpha_2,
\end{displaymath}
where $\alpha_2$ is as in \autoref{lemma:lower:bound:on:3cA} and $\kappa$ is as in \autoref{eqn:kappa_growth_good_norms}.
Then, for any $\ell > 0$, 
\begin{equation}
\label{eq:choice:of:alpha3}
  e^{\kappa^{-1}(\ell/2+\alpha_3 \ell)} e^{-\alpha_2 \ell} > 1, 
\end{equation}



\subsection{Modifying \texorpdfstring{$\cA_{1}(\cdot,\cdot,\cdot)$}{curly A1}}
	\label{ssec:modifying_curly_a1}

\label{sec:subsec:modifying:3cA_1}

In the rest of this subsection, we assume that the family $q \to U^+(q)$ satisfies the QNI condition (see \autoref{def:QNI}).
\medskip

For the next statement, the constant $\kappa_1(N)$ is defined in \autoref{eqn:def_kappa_1_Vs} to control the contraction on $V^s_-$.

\begin{proposition}[Modifying $\cA_1$ to be Lipschitz in time]
	\label{prop:modifying:3cA}
We can modify the map $\cA_{1}(\cdot, \cdot, \cdot)$, denoted $\cA_1^{new}$ so that
\autoref{theorem:3varA:haltime:factorizable} still holds (with $\cA_{1/2}^{new}:=g_{-\ell/2}\circ \cA_1^{new}$), and for any $\delta>0$ there exists a compact $K$ of measure at least $1-\delta$ such that provided $q,q_{1/2},q_1,q',q'_{1/2}\in K$ we have
\[
	\norm{\cA_1(q_1,u,\ell)-\cA_1^{new}(q_1,u\ell)}\leq e^{-N'\ell}
\]
where $N'$ is as in \autoref{lemma:3varA:almost:equivariant}

Furthermore, there exists a measurable function $\ell_0''(\cdot,\cdot): Q
\cross \reals \to \reals_+$
and there exists $\rho < \kappa_1(N)^{-1}/2$ so that for any $\ell >
\ell_0''(q_1,N) + \ell_0''(uq_1,N)$, 
any $F \in V^s_-(g_{-\ell} q_1)$ with
\begin{equation}
\label{eq:bar:cA:F:big}
\norm{{\cA}_{1}^{new}(q_1,u,\ell) F } \ge e^{-(N'/4)\ell} \|F\|
\end{equation} 
and for any {$0 < s < \ell/(2 \kappa_1(N))$,} we have:
\begin{displaymath}
e^{-\rho s} \norm{{\cA}_{1}^{new}(q_1,u,\ell) F} \le \norm{{\cA}_{1}^{new}(q_1,u,\ell+s) g_{-s} F\norm{\le e^{\rho s} }{\cA}_{1}^{new}(q_1,u,\ell) F}.
\end{displaymath}
\end{proposition}

After the proof of this result, we will write $\cA_i$ instead of $\cA_{i}^{new}$ and simply assume that $\cA_i$ satisfies the conclusions.

\begin{proof} Choose $\delta > 0$ so that $100 \delta \kappa_1(N) <
N'/4$, where $\kappa_1(N)$ is defined in \autoref{eqn:def_kappa_1_Vs}.
Let $K$ be as in \autoref{lemma:3varA:almost:equivariant}.
Let
\begin{equation}
\label{eq:tmp:def:Egood}
E_{good}(q_1) = \{ \ell \in \reals_+ \st g_{-\ell} q_1
\in K \text{ and } g_{-\ell/2} q_1 \in K \}.  
\end{equation}
Let $\ell_0''(\cdot,\cdot): Q \cross \reals \to \reals_+$ be such that for $\ell > \ell_0''(x,N)$,
\begin{equation}
\label{eq:tmp:most:s:in:lusin}
|\{ s \in [0,\ell] \st g_{-s} x \in K \}| > (1-10\delta) \ell. 
\end{equation}
Fix $D \gg 0$. 
For a fixed $q_1$, we inductively pick $\ell_j \in E_{good}(q_1)$ to
be as small as possible provided that $ \ell_j - \ell_{j-1} \ge D$.
Now suppose $\ell > \ell_0''(q_1,N) + \ell_0''(uq_1,N)$, and suppose $j$
is such that $\ell_{j-1} \le \ell \le \ell_j$. 
Then, by (\ref{eq:tmp:def:Egood}) and (\ref{eq:tmp:most:s:in:lusin}),
$\ell_j - \ell_{j-1} < 100 \delta \ell_j$. Write
$\bar{\cA}_j = {\cA}_{1}(q_1,u,\ell_j)$,
$\bar{\cA}_{j-1} = {\cA}_{1}(q_1,u,\ell_{j-1})$.

For $\ell_{j-1} \le \ell \le \ell_j$, we define
\begin{displaymath}
\bar{\cA}_1^{new}(q_1,u,\ell) = \frac{\ell-\ell_{j-1}}{\ell_j -
  \ell_{j-1}} \bar{\cA}_j \circ g_{\ell_j-\ell} + 
\frac{\ell_j-\ell}{\ell_j -   \ell_{j-1}} \bar{\cA}_{j-1}
\circ g_{\ell-\ell_{j-1}}
\end{displaymath}
Note that by \autoref{lemma:3varA:almost:equivariant}, \autoref{eq:lemma:3varA_1:almost:equivariant}, we have:
\begin{displaymath}
\norm{\bar{\cA}_j - \bar{\cA}_{j-1} \circ
g_{-(\ell_j - \ell_{j-1})} } \le e^{-N' \ell_j},
\end{displaymath}
where $N' \to \infty$ as $N \to \infty$ and so in particular
\begin{align}
	\label{eqn:A1new_A1_comparison}
	\norm{\cA_{1}^{new}(q_1,u,\ell)-\cA_1(q_1,u,\ell)}\leq e^{-N'\ell}
\end{align}
Suppose first that $\ell_{j-1} \le \ell \le \ell+s \le \ell_j$.
Then,
\begin{displaymath}
\bar{\cA}_1^{new}(q_1,u,\ell+s) \circ g_{-s} - \bar{\cA}_1^{new}(q_1,u,\ell)
= \frac{s}{\ell_j-\ell_{j-1}}  \left( \bar{\cA}_j - \bar{\cA}_{j-1}
  \circ g_{\ell_j - \ell_{j-1}} \right) \circ g_{\ell-\ell_j}.
\end{displaymath}
Thus, for any $F \in V^s_-(g_{-\ell} q_1)$,
\begin{multline*}
\norm{\bar{\cA}_1^{new}(q_1,u,\ell+s) \circ g_{-s} F - \bar{\cA}_1^{new}(q_1,u,\ell) F }
\le \\
\frac{s}{|\ell_j-\ell_{j-1}|}
\norm{\bar{\cA}_j - \bar{\cA}_{j-1} \circ g_{\ell_j - \ell_{j-1}}} \norm{g_{\ell-\ell_j} F} \leq 
\\
\le \frac{s}{D} e^{\kappa_1(N)(\ell-\ell_{j-1}) -N' \ell_j} \|F\|\le \frac{s}{D}
e^{-(3N'/4) \ell}\|F\|.
\end{multline*}
Now suppose
\begin{equation}
\label{eq:tmp:cond:AF:big}
\norm{\bar{\cA}_1^{new}(q_1,u,\ell) F} \ge e^{-(N'/2)\ell}  \|F\|.
\end{equation}
Then,
\begin{displaymath}
\left(1-\frac{s}{D}
e^{-(N'/4) \ell}\right) \le \frac{\|\bar{\cA}_1^{new}(q_1,u,\ell+s) \circ g_{-s}
  F\|}{\|\bar{\cA}_1^{new}(q_1,u,\ell) F\|} \le \left(1+\frac{s}{D}
e^{-(N'/4) \ell} \right).
\end{displaymath}
Recall that we are assuming that
$\ell_{j-1} \le \ell < \ell + s \le \ell_j$, and $\ell_j-\ell_{j-1} <
100 \delta \ell_j$. Therefore $0 < s < 100 \delta \ell$.
Then, for any $\rho > 0$,
provided that $\ell$ is large enough,
\begin{equation}
\label{eq:tmp:cAnew:bilip}
e^{-\rho s} \le \frac{\|\bar{\cA}_1^{new}(q_1,u,\ell+s) \circ g_{-s}
F\|}{\|\bar{\cA}_1^{new}(q_1,u,\ell) F\|} \le e^{\rho s}. 
\end{equation}
Indeed, we have
\[
	e^{-s\rho} \leq \max\left(1-\frac{s\rho}{100},1-\frac{1}{100}\right)\text{assuming }s\geq 0, \rho\leq \tfrac 1{10}
\]
If $s\rho\leq 1$ then we have
\[
	1-\frac{s\rho}{100}\leq 1 - \frac{s}{D}e^{-(N'/4)\ell}
\]
as soon as $\ell$ is sufficiently large (depending on $\rho$), and if $s\rho\geq 1$ then we have
\[
	1-\frac{1}{100}\leq 1 - \frac{s}{D}e^{-(N'/4)\ell}
\]
again as soon as $\ell$ is sufficiently large (recall that $s\leq 100\delta \ell$).
Analogously $1+\frac{s}{D} e^{-(N'/4) \ell} \leq e^{\rho s}$ for $\ell$ sufficiently large depending on $\rho$ and thus \autoref{eq:tmp:cAnew:bilip} follows.

It remains to relax the assumption that $\ell+s < \ell_j$. Suppose
(\ref{eq:bar:cA:F:big}) holds. Then, in view of
\autoref{lemma:3varA:almost:equivariant}, for
{$0 < s < \ell/(2 \kappa_1(N))$},
\begin{multline*}
\norm{{\cA}_1(q_1,u,\ell+s) g_{-s} F } 
\geq \norm{{\cA}_1(q_1,u,\ell) F } - e^{-N'\ell}\norm{F}\\
\ge \left(e^{-(N'/4)\ell} -
e^{-N' \ell}\right) \|F\| \ge \\ \ge
\tfrac{1}{2} e^{-(N'/4)\ell} 
e^{-\kappa_1(N) s}  \|g_{-s}F\| \ge \tfrac{1}{2} e^{-(1/2 + N'/4)\ell}\norm{g_{-s}F}.
\end{multline*}
Using \autoref{eqn:A1new_A1_comparison} we conclude that
\[
	\norm{{\cA}_1^{new}(q_1,u,\ell+s) g_{-s} F } \geq
	\tfrac{1}{4} e^{-(1/2 + N'/4)\ell}\norm{g_{-s}F}.
\]
Thus, provided $N'$ is large enough and $\ell$ is large enough and {$0 < s < \ell/(2 \kappa_1(N))$}, 
\begin{displaymath}
\|\bar{\cA}_1^{new}(q_1,u,\ell+s) g_{-s} F \| \ge e^{-(N'/2)\ell}
\|g_{-s} F\|.
\end{displaymath}
Then, the estimate (\ref{eq:tmp:cAnew:bilip}) can be iterated, 
and the condition (\ref{eq:tmp:cond:AF:big}) persists under iteration.
Thus (\ref{eq:tmp:cAnew:bilip}) holds for any {$0 < s < \ell/(2 \kappa_1(N))$}.
\end{proof}
From now on, we assume that $\cA_{1/2}(\cdot, \cdot, \cdot)$ and
${\cA}_{1}(\cdot,\cdot,\cdot)$ are such that \autoref{prop:modifying:3cA} holds. 



\subsection{The 4-variable \texorpdfstring{$\cA_2$}{curly A2}}
	\label{sec:subsec:4var:curlyA}
	\label{sssec:the_4_variable_curly_a_2}


Recall that the cocycle $\bH$ is defined in
\autoref{ssec:the_bold_h_cocycle}. 
We define
\index{$A$@$\cA_2(q_1,u, \ell,t)$}$$\cA_2(q_1,u, \ell,t) \colon V^s_*(q) \to \bbH(q_2)$$
as
$$\cA_2(q_1,u, \ell,t ) =  g_{t+\ell/2} \circ u_*(q_{1/2}, y_{1/2}) \circ   \cA_{1/2}(q,u, \ell) \mod L\cC^{\tau}(q_2)$$
where $u_*$ is defined in \autoref{def:u_star}.
Note that because $u_{*}$ is $g_t$-equivariant, we have an equivalent definition of $\cA_2$ as:
\begin{align}
	\label{eqn:curly_A2_in_terms_of_curly_A1}
	\cA_2(q_1,u, \ell,t ) =  g_{t} \circ u_*(q_{1}, uq_1) \circ   \cA_{1}(q_1,u, \ell) \mod L\cC^{\tau}(q_2)
\end{align}

using the definition of $\cA_1$ from \autoref{sssec:defining_curly_a_1}.
Note that we have
\begin{equation}
\label{eq:4varA:equivariant}
\cA_2(q_1, u,\ell,\tau+t) = g_t \circ \cA_2(q_1,u,\ell,\tau). 
\end{equation}
Also, in view of \autoref{eq:3varcA:q:in:LC},
\begin{equation}
\label{eq:4varcA:q:0}
\cA_2(q_1, u,\ell,\tau) F^s_q(q)= 0. 
\end{equation}

\subsection{Bilipshitz estimates}
	\label{sec:new:bilipshitz}
We denote by
\index{$A$@$\|\cA_2(q_1,u,\ell,\tau)\|$}$\|\cA_2(q_1,u,\ell,\tau)\|$ 
the operator norm of $\cA_2(q_1,u,\ell,\tau)$ with respect to the
dynamical norms $\| \cdot \|_{V^s_*}$  on $V^s_*(g_{-\ell}q_1)$ 
and $\| \cdot \|$ on $\bH(g_t u q_1)$.


\begin{lemma}[Future bilipschitz]
\label{lemma:A:future:bilipshitz}
We have, for $t > s$, 
\begin{multline}
\label{eq:A:future:bilipshitz}
e^{\kappa^{-1}(t-s)} \|\cA_2(q_1,u,\ell,s)\| \le
\|\cA_2(q_1,u,\ell,t)\|\le \\ \le e^{\kappa(t-s)} \|\cA_2(q_1,u,\ell,s)\|, 
\end{multline}
where $\kappa$ is as in \autoref{eqn:kappa_growth_good_norms}.
\end{lemma}

\begin{proof}
Let $q = g_{-\ell} q_1$. Let $F \in V^s_*(q)$ be such that $\| F \|_{V^s_*(q)} =
1$ and 
\begin{displaymath}
\|\cA_2(q_1,u,\ell,t)\| = \|\cA_2(q_1,u,\ell,s) F \|. 
\end{displaymath}
Then, using \autoref{prop:good_norms}\autoref{good_norms_bilipschitz},
\begin{multline*}
\|\cA_2(q_1,u,\ell,t)\| =  \|\cA_2(q_1,u,\ell,t) F \| = \|g_{t-s} \circ \cA_2(q_1,u,\ell,s) F \| \\
   \le e^{\kappa(t-s)} \|\cA_2(q_1,u,\ell,s) F \|
   \le e^{\kappa(t-s)} \|\cA_2(q_1,u,\ell,s) \|.
\end{multline*}
Thus proves the upper bound in (\ref{eq:A:future:bilipshitz}). 
To prove the lower bound,
let $F \in V^s_*(q)$ be such that $\| F \|_{V^s_*(q)} =
1$ and 
\begin{displaymath}
\|\cA_2(q_1,u,\ell,s)\| = \|\cA_2(q_1,u,\ell,s) F \|. 
\end{displaymath}
Then, using \autoref{prop:good_norms}\autoref{good_norms_bilipschitz},
\begin{multline*}
\|\cA_2(q_1,u,\ell,t)\| \ge  \|\cA_2(q_1,u,\ell,t) F \| = \|g_{t-s} \circ \cA_2(q_1,u,\ell,s) F \| \\
  \ge e^{\kappa^{-1}(t-s)} \|\cA_2(q_1,u,\ell,s) F \|
  = e^{\kappa^{-1}(t-s)} \|\cA_2(q_1,u,\ell,s) \|.
\end{multline*}
\end{proof}

\begin{lemma}[Past bilipschitz]
\label{lemma:A:past:bilipshitz}
There exists $N_0 > 1$ such that if the parameter
$N$ in \autoref{theorem:3varA:haltime:factorizable} is greater
than $N>N_0$, then there exist $\index{$\kappa_3$}\kappa_3=\kappa_3(N) > 1$ and 
$\beta=\beta(N) > 0$ and 
a function $\ell_0: X \cross \reals \to \reals^+$
such that for $\ell > \ell_0(q_1,N) + \ell_0(u q_1,N)$, $\tau < \alpha_3 \ell$
and {$0 < s < \ell/(2 \kappa_1(N))$},
\begin{multline}
\label{eq:A:past:bilipshitz}
e^{-\kappa_3 s} \|{\cA}_{2}(q_1,u,\ell,\tau)\| \le
\|{\cA}_{2}(q_1,u,\ell+s,\tau)\|\le \\ \le e^{-\kappa_3^{-1}s}
\|{\cA}_{2}(q_1,u,\ell,\tau)\|,
\end{multline}
provided $\|\cA_2(q_1,u,\ell,\tau)\| > e^{-\beta \ell}$.
\end{lemma}

\begin{proof}
Let $q = g_{-\ell} q_1$.
Let $N'$ be as in \autoref{lemma:3varA:almost:equivariant} and
\autoref{prop:modifying:3cA}, so that $N' \to \infty$ as $N
\to \infty$.  {Let $\beta = N'/8$.}

\noindent\textbf{Preliminary bound.} We first prove that
\begin{equation}
\label{eq:cA:ell:plus:s:lower:bound}
\|{\cA}_{2}(q_1,u,\ell+s,\tau)\| \ge e^{-(N'/8) \ell} -
e^{-(N'-\kappa \alpha_3) \ell}.
\end{equation}
Let $F' \in V^s_*(q)$ be such that
$\|{\cA}_{2}(q_1,u,\ell,\tau) F' \| = \|{\cA}_{2}(q_1,u,\ell,\tau)\| \|F' \|$. In
view of (\ref{eq:4varcA:q:0}) and \S\ref{sec:subsec:norm:Vs:star} (a), we have
$F' \in V^s_-(q)$. We have
\begin{equation}
\label{eq:tmp:cAF:one}
  \|{\cA}_{2}(q_1,u,\ell,\tau) F' \| = \|{\cA}_{2}(q_1,u,\ell,\tau)\| \|F'\| \ge
e^{-(N'/8)\ell}\|F'\|. 
\end{equation}
By \autoref{lemma:3varA:almost:equivariant} and \autoref{eqn:curly_A2_in_terms_of_curly_A1}:
\begin{multline*}
\|{\cA}_{2}(q_1,u,\ell+s,\tau) g_s F' - {\cA}_{2}(q_1,u,\ell,\tau) F'\| =\\
\norm{g_{\tau}\circ u_*(q_1,uq_1)\circ \left[{\cA}_{1}(q_1,u,\ell+s) g_s - {\cA}_{1}(q_1,u,\ell)\right] F'}\leq \\
\leq 
C(q_1,uq_1)e^{\kappa \tau - N' \ell} \|F'\|\le C(q_1,uq_1)e^{-(N'-\kappa \alpha_3) \ell}\|F'\|.  
\end{multline*}
for a constant $C(q_1,uq_1)$ depending measurably on the points.
Therefore using \autoref{eq:tmp:cAF:one},
\begin{multline*}
\|{\cA}_{2}(q_1,u,\ell+s,\tau) g_s F'\| \ge \left(e^{-(N'/8) \ell} -
C(q_1,uq_1)e^{-(N'-\kappa \alpha_3) \ell}\right) \|F'\| \ge \\ \ge 
\left(e^{-(N'/8) \ell} - C(q_1,uq_1)e^{-(N'-\kappa \alpha_3) \ell} \right) \|g_s F'\|
\end{multline*}
since $s > 0$ and $\norm{F'}\geq \norm{g_s F'}$.
This completes the proof of \autoref{eq:cA:ell:plus:s:lower:bound}.

In view of \autoref{prop:modifying:3cA}, for   $0<s < \ell/(2\kappa_1(N))$,
restricted to $V^s_-$ we have
\begin{align}
	\label{eqn:A_1_H_ell_s}
	{\cA}_{1}(q_1,u,\ell+s) \circ g_{-s} &  = {\cA}_{1}(q_1,u,\ell) \circ H_{\ell,s}
\end{align}
where $H_{\ell,s}: V^s_{-}(q) \to V^s_-(q)$ is a linear map
satisfying
\begin{equation}
\label{eq:norm:Hls}
\max(\|H_{\ell,s}F'\|, \|H_{\ell,s}^{-1}F'\|) \le e^{\rho s} \|F'\|, 
\end{equation}
for any $F' \in V^s_{-}(q)$ satisfying
\begin{equation}
\label{eq:F:lower:bound:one}
\|{\cA}_{1}(q_1,u,\ell)F'\|_{V^s_*(q)} \ge e^{-(N'/4)\ell}\|F'\|.   
\end{equation}
Therefore, by applying $g_{\tau}\circ u_*(q_1,uq_1)$ to \autoref{eqn:A_1_H_ell_s} we find:
\[
	{\cA}_{2}(q_1,u,\ell+s,\tau) \circ g_{-s}  = {\cA}_{2}(q_1,u,\ell,\tau) \circ H_{\ell,s},
\]
and so
\begin{displaymath}
{\cA}_{2}(q_1,u,\ell+s,\tau) = {\cA}_{2}(q_1,u,\ell,\tau) \circ
H_{\ell,s} \circ g_s.  
\end{displaymath}
when restricting the domain to $V^s_-$.

\noindent \textbf{Lower bound.} We now prove the lower bound in \autoref{eq:A:past:bilipshitz}.
Let $F'$ be as in \autoref{eq:tmp:cAF:one}.
We have
\begin{multline}
\|{\cA}_{1}(q_1,u,\ell) F'\| \ge e^{-\kappa \tau} \norm{u_*(uq_1,q_1)}
\|{\cA}_{2}(q_1,u,\ell,\tau) F'\| \geq \\
\ge C(q_1,uq_1) e^{-((N'/8)+\kappa \alpha_3)\ell}
\|F'\| \ge e^{-(N'/4)\ell} \|F'\|
\end{multline}
once $\ell_0(q_1),\ell_0(uq_1)$ are sufficiently large.
Therefore \autoref{eq:F:lower:bound:one} holds as long as $N$ and
therefore $N'$ are large enough.
Let $F = g_{-s} \circ H_{\ell,s}^{-1} F'$ so that $F' = H_{\ell,s} \circ g_s
F$. 
Then, $F \in V^s_-(g_{-s} q)$. We may normalize $F$ so that $\|F\|_{V^s_-(g_{-s} q)} = 1$.
Using \S\ref{sec:subsec:norm:Vs:star}(c) and (\ref{eq:norm:Hls}),
\begin{multline*}
  \|{\cA}_{2}(q_1,u,\ell+s,\tau)\| \ge  \|{\cA}_{2}(q_1,u,\ell+s,\tau) F \| =
  \|{\cA}_{2}(q_1,u,\ell,\tau) \circ H_{\ell,s} \circ g_s F \| = \\
  = \|{\cA}_{2}(q_1,u,\ell,\tau)\| \|g_s \circ H_{\ell,s} F \|
  \ge e^{-(\kappa_1(N)+\rho)s} \|{\cA}_{2}(q_1,u,\ell,\tau)\|.
\end{multline*}
Thus the lower bound in \autoref{eq:A:past:bilipshitz} holds. 

\noindent \textbf{Upper bound.} We now prove the upper bound in (\ref{eq:A:past:bilipshitz}).
Let $F \in V^s_*(g_{-s}q)$ be such that $\| F
\|_{V^s_*(g_{-s} q)} =
1$ and 
\begin{displaymath}
\|{\cA}_{2}(q_1,u,\ell+s,\tau)\| = \|{\cA}_{2}(q_1,u,\ell+s,\tau) F \|. 
\end{displaymath}
In view of \autoref{eq:4varcA:q:0} and \S\ref{sec:subsec:norm:Vs:star}(a), we have $F \in V^s_-(g_{-s} q)$.
Let $F' = g_s F$, then we claim \autoref{eq:F:lower:bound:one} holds for it.
Indeed, if this is not the case, then by
\autoref{lemma:3varA:almost:equivariant},
\begin{multline*}
  \|{\cA}_{1}(q_1,u,\ell+s) F \| \le \|{\cA}_{1}(q_1,u,\ell) g_s F
    \| + e^{-N' \ell} \|F\| = \\
    = \norm{{\cA}_{1}(q_1,u,\ell) F'}  + e^{-N'  \ell}
    \leq e^{-(N'/4)\ell} \norm{F'}  + e^{-N'  \ell}\leq\\
    \leq e^{-(N'/4)\ell} e^{-\kappa_1(N)^{-1} s} + e^{-N'  \ell}
\end{multline*}
where we used that $\norm{F'}=\norm{g_s F}\leq e^{-\kappa_1(N)^{-1}s}$ (and $\norm{F}=1$).
Then,
\begin{multline*}
\|{\cA}_{2}(q_1,u,\ell+s,\tau) \| = 
  \|{\cA}_{2}(q_1,u,\ell+s,\tau) F\| = \\ = \| g_{\tau} \circ u_*(q_1,uq_1) \circ
  {\cA}_{1}(q_1,u,\ell+s) F\| \le 
  C(q_1,uq_1)\cdot e^{\kappa \tau}\left(e^{-(N'/4)\ell} e^{-\kappa_1(N)^{-1} s} + e^{-N'  \ell}\right).
\end{multline*}
However, since $\tau < \alpha_3 \ell$ and $s>0$, 
this contradicts
\autoref{eq:cA:ell:plus:s:lower:bound} as long as $N$ (and thus $N'$)
are large enough. 
Therefore, (\ref{eq:F:lower:bound:one}) holds. 
Then, using \S\ref{sec:subsec:norm:Vs:star}(c) and (\ref{eq:norm:Hls}),  
\begin{multline*}
  \|{\cA}_{2}(q_1,u,\ell+s,\tau)\| =  \|{\cA}_{2}(q_1,u,\ell+s,\tau) F \|
 = \|{\cA}_{2}(q_1,u,\ell,\tau) \circ H_{\ell,s}\circ g_s F\| \\
\le \|{\cA}_{2}(q_1,u,\ell,\tau)\| \| H_{\ell,s} \circ g_s F \| \le 
e^{\left(-\kappa_1(N)^{-1} +\rho\right) s} \|{\cA}_{2}(q_1,u,\ell,\tau)\|.
\end{multline*}
By \autoref{prop:modifying:3cA} the exponent $\rho$ can be made arbitrarily small and positive and this proves the upper bound in (\ref{eq:A:past:bilipshitz}). 
\end{proof}

\subsubsection{Choice of $N$.}
\label{sec:subsec:choice:of:N}
We now choose the parameter \index{$N$}$N$ in \autoref{theorem:3varA:haltime:factorizable} so that 
\begin{equation}
\label{eq:choice:of:N}
e^{\kappa(4\alpha_3+1/2) \ell} e^{-N \ell} < e^{-\ell}, 
\end{equation}
and also so that \autoref{lemma:lower:bound:on:3cA} and
\autoref{lemma:A:past:bilipshitz} hold.

\subsubsection{Defining $\tau_\epsilon$}
	\label{sssec:defining_tau_epsilon}
For $\epsilon > 0$, almost all $q_1 \in X_0$, almost all $u
\in \cB_0(q_1)$ and $\ell > 0$, let
\begin{displaymath}
\index{$\tau$@$\tau_{(\epsilon)}(q_1,u, \ell)$}\tau_{(\epsilon)}(q_1,u, \ell) = \sup \{\tau \st \tau > 0 \text{ and }
\|\cA_2(q_1,u,\ell,\tau)\| \le \epsilon  \}.
\end{displaymath}

\begin{lemma}
	\label{lemma:monotonicity}
There exists \index{$\kappa_4$}$\kappa_4 > 1$ (depending only on the
Lyapunov spectrum) with the following property: for almost all $q_1
\in X_0$, $u \in \cB_0(q_1)$, for all $\ell > \ell_0(q_1) +
\ell_0(u q_1)$ and $s > 0$, 
\begin{displaymath}
\tau_{(\epsilon)}(q_1,u,\ell) + \kappa_4^{-1} s < \tau_{(\epsilon)}(q_1,u, \ell + s)< \tau_{(\epsilon)}(q_1,u,\ell) + \kappa_4 s. 
\end{displaymath}
\end{lemma}

\begin{proof}
  It is enough to prove the lemma assuming {$0 < s < \ell/(2 \kappa_1(N))$}.
Let $\tau = \tau_{(\epsilon)}(q_1,u,\ell)$, $\tau_s = \tau_{(\epsilon)}(q_1,u,\ell+s)$. We have
\begin{displaymath}
\epsilon = \|\cA_2(q_1,u,\ell,\tau)\| = \|\cA_2(q_1,u,\ell+s,\tau_s)\|.
\end{displaymath}
Therefore,
\begin{displaymath}
  \frac{\|\cA_2(q_1,u,\ell+s,\tau_s)\|}{\|\cA_2(q_1,u,\ell+s,\tau)\|}
  = \frac{\|\cA_2(q_1,u,\ell,\tau)\|}{\|\cA_2(q_1,u,\ell+s,\tau)\|} 
\end{displaymath}
In view of \autoref{lemma:A:past:bilipshitz}, this implies
\begin{displaymath}
  e^{\kappa_3^{-1} s} \le \frac{\|\cA_2(q_1,u,\ell+s,\tau_s)\|}{\|\cA_2(q_1,u,\ell+s,\tau)\|} \le e^{\kappa_3 s}. 
\end{displaymath}
Now, \autoref{lemma:A:future:bilipshitz} implies
\begin{displaymath}
\kappa^{-1} \kappa_3^{-1} s \le \tau_s - \tau < \kappa \kappa_3 s. 
\end{displaymath}
so $\kappa_4:=\kappa\kappa_3$ suffices.
\end{proof}

\begin{proposition}
\label{prop:bilip:hattau:epsilon}
There exists $\kappa > 1$ depending only on the Lyapunov spectrum, and
such that for almost all $q_1  \in X_0$, almost all $u \in
\cB_0(q_1)$, any $\ell > \ell_0(q_1) + \ell_0(u q_1)$ and any measurable 
subset $E_{bad} \subset \reals^+$, 
\begin{displaymath}
|\tau_{(\epsilon)}(q_1,u,E_{bad}) \cap [\tau_{(\epsilon)}(q_1,u,0),\tau_{(\epsilon)}(q_1,u,\ell)]|  \le \kappa
|E_{bad} \cap [0,\ell]|
\end{displaymath}
and
\begin{displaymath}
|\{ t \in [0,\ell]| \st \tau_{(\epsilon)}(q_1,u,t) \in E_{bad}
\}|  \le \kappa |E_{bad} \cap
      [\tau_{(\epsilon)}(q_1,u,0),\tau_{(\epsilon)}(q_1,u,\ell)]|.
\end{displaymath}
\end{proposition}
\begin{proof}
This follows immediately from \autoref{lemma:monotonicity}. 
\end{proof}

\noindent For future use we record the following restatement of \autoref{theorem:3varA:haltime:factorizable}, in terms of $\cA_2$ rather than $\cA_{1/2}$:
\begin{proposition}
\label{prop:factorization:main}
There exists a constant $\alpha > 0$ and  
for every $\delta > 0$ there is a compact set
$K$ with $\nu(K) > 1-\delta$ such that the following holds: Suppose
$Y$ and $Y'$ are two bottom-linked pre-$Y$ configurations, such that
$Y'$ left-shadows $Y$. Suppose $Y \in K$ and also (using the notation
\S\ref{sec:subsec:Yconf:notation}) $q' \in K$, $q'_{1/2} \in K$, $q_1'
\in K$. Then, 
the linear map $\cA_2(q_1,u, \ell,\tau) \colon V^s_*(q) \to \bbH(q_2)$ satisfies 
\begin{displaymath}
\cA_2(q_1,u,\ell,\tau) F^s_q(q') =
\bfj\big(\phi_\tau^{-1}\bigl(U^+[g_\tau q_1']
\bigr)\big) + O(e^{-\alpha \ell}). 
\end{displaymath}
\end{proposition}
\noindent This is a simplified formulation of \autoref{theorem:3varA:haltime:factorizable} since the constant $N$ has now been fixed, and hence also the cocycle $V^s_*$.





\section{Some technical lemmas}
	\label{sec:some_technical_lemmas}


\subsection{Left-balanced \texorpdfstring{$Y$}{Y}-configurations}
	\label{ssec:left_balanced_y_configurations}

\begin{definition}[Left-balanced $Y$-configurations]
 \label{def:left:balanced:Y}
Fix $C\ge1$.  We say a pre-$Y$-configuration or a $Y$-configuration $Y$ is
\emph{$(C,\epsilon)$-left-balanced} if using the notation of
\autoref{sec:subsec:Yconf:notation} and \autoref{eqn:curly_A2_in_terms_of_curly_A1} we have
\[
	\frac 1 C  \epsilon \le \| \cA_2(q_1, u, \ell,\tau)\| \le  C \epsilon.
\]
Equivalently, $Y$ is $(C,\epsilon)$-left-balanced if using the
notation of \autoref{sec:subsec:Yconf:notation} and \autoref{sssec:defining_tau_epsilon} we have
\begin{displaymath}
|\tau - \tau_{(\epsilon)}(q_1,u,\ell)| < C'
\end{displaymath}
for some $C'$ depending on $C$. 
\end{definition}

We denote by \index{$Y^{(\epsilon)}_{pre,lb}(q_1,u,\ell)$}$Y^{(\epsilon)}_{pre,lb}(q_1,u,\ell)$ the unique
$(1,\epsilon)$-left-balanced pre-$Y$ configuration $Y$ with $q_1(Y) =
q_1$, $\ell(Y) = \ell$, $u(Y) = u$.  


\subsection{A selection lemma}

\begin{lemma}
\label{lemma:bad:subspace}
For any $\rho > 0$ there is a constant $c(\rho)$ with the following
property:  
Let $A: \cV \to \cW$ be a linear map between Euclidean spaces.
Then there exists a proper subspace
$\cM \subset \cV$ such that for any $v$ with
$\|v\|=1$ and $d(v,\cM) > \rho$, we have 
\begin{displaymath}
\|A\| \ge \|A v\| \ge c (\rho) \|A\|. 
\end{displaymath}
Furthermore, if $A,\cV,\cW$ vary measurably on some parameter, then $\cM$ can be chosen to vary measurably as well.
\end{lemma}
\begin{proof}
The matrix $A^t A$ is symmetric, so it has a complete
orthogonal set of eigenspaces $W_1, \dots, W_m$ corresponding to
eigenvalues $\mu_1 > \mu_2 > \dots \mu_m$. Let $\cM = W_1^{\perp}$. 
\end{proof}

\noindent In this subsection, let $(\cB, | \cdot |)$ be a standard Lebesgue space equipped with a finite measure space.
\medskip


The following is a slight modification of 
\cite[Proposition~5.3]{EskinMirzakhani_Invariant-and-stationary-measures-for-the-rm-SL2Bbb-R-action-on-moduli-space}
since we are using partition atoms and not balls. However, the proof
is virtually unaffected. 

\begin{proposition}
\label{prop:can:avoid:most:Mu}
For every $\delta > 0$ there exist constants 
$c_1(\delta) > 0$,  $\epsilon_1(\delta) > 0$
with $c_1(\delta) \to 0$ and $\epsilon_1(\delta) \to 0$ 
as $\delta \to 0$, and also constants
$\rho(\delta) > 0$, $\rho'(\delta) > 0$, and
$C(\delta) < \infty$ such that the following holds: 

For any subset $K' \subset X_0$ with $\nu(K') > 1-\delta$, 
there exists a subset $K \subset K'$ 
with $\nu(K) > 1-c_1(\delta)$
such that the following holds: suppose for each $x\in X_0$ we have a
measurable map from $\cB$ to proper subspaces of
$V^s_*(x) $, written as 
$u \to \cM_u(x)$, where $\cM_u(x)$ is
a proper subspace of $V^s_*(x) $. Then, for 
any $q \in K$ there exists $q' \in K' \cap \gB_0^-[q]$ with
\begin{equation}
\label{eq:rho:delta:le:Fq:minus:Fqprime}
C(\delta)^{-1} \le \|F_{q}^s(q') \|_{V^s(q)} \le C(\delta)
\end{equation}
and such that
\begin{multline}
\label{eq:qprime:avoids:Mu}
d_{V^s(q)}(F_{q}^s(q'), 
\cM_u(q)) > \rho(\delta) \\ \text{\rm  for at least
  $(1-\epsilon_1(\delta))$-fraction of $u \in \cB$. } 
\end{multline}
\end{proposition}


\subsection{A switching lemma}
	\label{ssec:a_switching_lemma}

Recall that left-shadowing $Y$-configurations are introduced in \autoref{def:left_shadowing_Y_configuration}.

\begin{lemma}
\label{lemma:swithching}
There exists $\alpha > 0$ and  
for every $\delta > 0$ there exists a compact set $K$ and for every
$\ell > 0$ a subset $K_\ell \subset K$ 
with $\nu(K_\ell) >
1 - c(\delta)$ with $c(\delta) \to 0$ as $\delta \to 0$ and
$c_1(\delta) > 0$ such that the
following holds. Suppose $Y$, $Y'$ are bottom linked 
pre-$Y$-configurations, with $Y \in K$ and $Y' \in K$. Suppose also
that $Y'$  left-shadows $Y$ and $Y$ left-shadows $Y'$.  Also
suppose (using the notation \S\ref{sec:subsec:Yconf:notation}) that $q
\in K_\ell$, $q'\in K_\ell$. Then, assuming $\|\cA_2(q_1,u,\ell,\tau)\| \le
c_1(\delta)$ and $\|\cA_2(q_1',u',\ell,\tau)\| \le c_1(\delta)$, we have 
\begin{equation}
\label{eq:cA:cAprime}
\|\cA_2(q_1', u', \ell, \tau) \| \le C(\delta) \|\cA_2(q_1, u, \ell, \tau)\| +
O(e^{-\alpha \ell}),
\end{equation}
and
\begin{equation}
\label{eq:cAprime:cA}
\|\cA_2(q_1, u, \ell, \tau) \| \le C(\delta) \|\cA_2(q_1', u', \ell, \tau)\| +
O(e^{-\alpha \ell}).
\end{equation}
\end{lemma}

\subsubsection{Proof .}
Let $K = K(\delta)$ be the set where
\autoref{prop:factorization:main} and
\autoref{theorem:interpolation:does:not:move:points}
apply. We also assume that the functions
$C_1(\cdot)$, $c_1(x)$, $C(x)$
of \autoref{lemma:L:and:C} and
\autoref{prop:hausdorff_distance_and_norm_of_vector}
are bounded by $C(\delta)$ for $x \in K$, and that the function $r(x)$
of \autoref{lemma:du:vs:dQ} is bounded from below by
$C(\delta)^{-1}$ for $x \in K$. 

Recall that
\begin{displaymath}
F_q^s(q') = H^s_{V^s}(q',q) \iota_s(q'),
\end{displaymath}
where $\iota_s(q') \in V^s(q')$ is an equivariant section, and
$H^s_{V^s}(q',q)$ is a holonomy map. For any $\delta > 0$ we may
choose a compact set $K''$ with $\nu(K'') > 1-\delta$ such that
for all $q \in K''$,
\begin{equation}
\label{eq:bound:iota}
\|\iota_s(q)\|_{V^S(q)} \le C(\delta)
\end{equation}
and also that
\begin{equation}
\label{eq:bound:holonomy}
\|H^s_{V^s}(q',q)\| \le C(\delta)  
\end{equation}
for all $q \in K''$, $q' \in K'' \cap \gB_0^-[q]$. 
Let $K' = K'' \cap K \cap g_{-\ell/2} K \cap g_{-\ell} K$. 
Then, let $K_\ell$ be the subset denoted by $K$ in
\autoref{prop:can:avoid:most:Mu}.

By symmetry, it is enough to prove
(\ref{eq:cA:cAprime}). Since $q,q' \in K''$, we have
\begin{equation}
\label{eq:Fs:qprime:q:norm}
\|F^s_{q}(q') \|_{V^s(q)} \le C(\delta). 
\end{equation}
Let $\cM'$ denote the subspace
as in \autoref{lemma:bad:subspace} for $\cA_2(q_1',u',\ell,\tau)$.
Note that by (\ref{eq:3varcA:q:in:LC}), 
$F^s_{q'}(q') \in \cM'$. 
By \autoref{prop:can:avoid:most:Mu} (with the function $u\to
\cM_u$ the constant function $\cM'$) we can choose $q'' \in \gB_0^-[q]
\cap K'$ 
with $d_{V^s(q')}(F^s_{q'}(q''), \cM') > \rho(\delta)$.
Then, by (\ref{eq:bound:iota}) and (\ref{eq:bound:holonomy})
\begin{equation}
\label{eq:Fs:qtwoprime:qprime:norm}
  \|F^s_{q'}(q'') \|_{V^s(q')} \le C(\delta). 
\end{equation}
Then, again by (\ref{eq:bound:holonomy}), 
\begin{equation}
\label{eq:Fs:qtwoprime:q:norm}
\|F^s_{q}(q'') \|_{V^s(q)} \le C(\delta).
\end{equation}
By \autoref{lemma:bad:subspace}, we have
\begin{equation}
\label{eq:Aprime:norm:upper:bound}
\|\cA_2(q_1', u', \ell,\tau) \| \le C(\delta) \| \cA_2(q_1',u',\ell,\tau)
F^s_{q'}(q'') \|.  
\end{equation}
Recall that using our notational conventions
\S\ref{sec:subsec:Yconf:notation}, $q_2 = g_\tau u q_1$, $q_2' =
g_\tau u' q_1'$. 
Using the notation as in \eqref{eq:phitau}, let $\phi_\tau=\phi_\tau(Y,q'): \cW^u[q_2] \to \cW^u[q_2']$ and $\phi_\tau''=\phi_\tau(Y,q''): \cW^u[q_2] \to \cW^u[g_\tau q_1'']$ be the interpolation maps relative to $u q_1$, and let
$\phi'_\tau=\phi_\tau(Y',q''): \cW^u[q_2'] \to \cW^u[g_\tau q_1'']$ be the interpolation
map relative to $u' q_1'$. Then,
\begin{displaymath}
\cU \equiv \phi_\tau^{-1}(U^+[q_2']) \in \cC[q_2]
\end{displaymath}
\begin{displaymath}
\cU'' \equiv (\phi''_\tau)^{-1}(U^+[g_\tau q_1'']) \in \cC[q_2],
\end{displaymath}
and
\begin{displaymath}
\cU' \equiv  (\phi'_\tau)^{-1}(U^+[g_\tau q_1'']) \in \cC[q_2'].
\end{displaymath}
By \autoref{prop:factorization:main}, 
\begin{equation}
\label{eq:AF:q:qprime}
\cA_2(q_1, u, \ell, \tau) F^s_{q}(q') = \bfj (\cU) +
O(e^{-\alpha \ell})
\end{equation}
\begin{equation}
\label{eq:AF:q:qtwoprime}
\cA_2(q_1, u, \ell, \tau) F^s_{q}(q'') = \bfj (\cU'') +
O(e^{-\alpha \ell})
\end{equation}
\begin{equation}
\label{eq:AF:qprime:qtwoprime}
\cA_2(q_1', u', \ell, \tau) F^s_{q'}(q'') = \bfj (\cU') +
O(e^{-\alpha \ell}).
\end{equation}
Note that by our assumptions, $\|\bfj(\cU)\|$, $\|\bfj(\cU')\|$ and
$\|\bfj(\cU'')\|$ are all bounded by $2c_1(\delta)C(\delta)$. We can choose
$c_1(\delta)$ so that for $x \in K$, $c_1(\delta) C(\delta) \ll r(x)$, where
$r(x)$ is as in \S\ref{rmk:on_lyapunov_radius}.
By \autoref{prop:hausdorff_distance_and_norm_of_vector},
and since $q_2 \in K$, $q_2' \in K$,  
\begin{displaymath}
C(\delta)^{-1} \|\bfj (\cU)\| \le d_{\cH,loc}^{q_2}(\cU, U^+[q_2]) \le
C(\delta) \|\bfj (\cU)\|,
\end{displaymath}
\begin{displaymath}
C(\delta)^{-1} \|\bfj (\cU'')\| \le d_{\cH,loc}^{q_2}(\cU'',U^+[q_2]) \le
C(\delta) \|\bfj (\cU'')\|,
\end{displaymath}
\begin{equation}\label{eq:iphonecharger}
  C(\delta)^{-1} \|\bfj (\cU')\| \le d_{\cH,loc}^{q_2'}(\cU',U^+[q_2']) \le
C(\delta) \|\bfj (\cU')\|.
\end{equation}
By \autoref{lemma:L:and:C} there exist $h,h'' \in
\bbG^{ssr}(\cW^u[q_2])$
such that
\begin{displaymath}
  \cU = h U^+[q_2], \qquad \cU'' = h'' U^+[q_2], 
\end{displaymath}
and
\begin{displaymath}
\|L(h) - I \|_{L\cW^u(q_2)} \le C(\delta) \|\bfj(\cU)\|, \qquad \|L(h'') - I \|_{L\cW^u(q_2)} \le C(\delta) \|\bfj(\cU'')\|. 
\end{displaymath}
Then,
\begin{displaymath}
\|L(h'' \circ h^{-1}) - I \|_{L\cW^u(q_2)} \le C(\delta)(\|\bfj(\cU)\| + \|\bfj(\cU'')\|). 
\end{displaymath}
It then follows from \autoref{lemma:L:is:bilipshitz} 
that for $x \in \cW^u_{loc}[q_2]$,
\begin{equation}
\label{eq:h:twoprime:circ:h:inverse}
d^u(x, h'' \circ h^{-1} x) \le C(\delta) (\|\bfj(\cU)\| + \|\bfj(\cU'')\|).
\end{equation}
Let
\begin{displaymath}
  h' = (\phi'_\tau)^{-1} \circ \phi_\tau'' \circ h'' \circ h^{-1}
  \circ \phi_\tau^{-1}. 
\end{displaymath}
Then, 
\begin{displaymath}
h' U^+[q_2'] = (\phi'_\tau)^{-1} \circ \phi_\tau'' \circ h'' \circ
h^{-1} (h U^+[q_2]) =  (\phi'_\tau)^{-1} \circ \phi_\tau''
\circ (\phi_\tau'')^{-1}(U^+[g_\tau q_1'']) = \cU'.
\end{displaymath}
By (\ref{eq:h:twoprime:circ:h:inverse}) and
\autoref{theorem:interpolation:does:not:move:points},
for $x \in \cW^u_{loc}[q_2']$,
\begin{displaymath}
d^Q(\sigma(x),\sigma(h'x)) \le C(\delta) (\|\bfj(\cU)\| + \|\bfj(\cU'')\|).
\end{displaymath}
Then, using \autoref{lemma:du:vs:dQ}, 
\begin{displaymath}
d^u(x, h'x) \le C(\delta) (\|\bfj(\cU)\| + \|\bfj(\cU'')\|).
\end{displaymath}
Therefore, since $h' U^+[q_2'] = \cU'$, 
\begin{displaymath}
d_{\cH,loc}^{q_2'}(U^+[q_2'], \cU') \le C(\delta) (\|\bfj(\cU)\| +
\|\bfj(\cU'')\|).
\end{displaymath}
Then, by \eqref{eq:iphonecharger}, 
\begin{displaymath}
\|\bfj(\cU')\| \le C(\delta) (\|\bfj(\cU)\| +\|\bfj(\cU'')\|). 
\end{displaymath}
Therefore, using (\ref{eq:AF:q:qprime}), (\ref{eq:AF:q:qtwoprime}) and
(\ref{eq:AF:qprime:qtwoprime}), we get 
\begin{multline*}
\|\cA_2(q_1', u', \ell, \tau) F^s_{q'}(q'')\| \le C(\delta) (\|\cA_2(q_1, u,
\ell, \tau) F^s_{q}(q')\| + \|\cA_2(q_1, u, \ell, \tau) F^s_{q}(q'')\|) +
O(e^{-\alpha \ell}) \\
\le C_1(\delta) \|\cA_2(q_1, u, \ell, \tau) \| + O(e^{-\alpha \ell}), 
\end{multline*}
where for the last estimate we used (\ref{eq:Fs:qprime:q:norm}) and
(\ref{eq:Fs:qtwoprime:q:norm}). Now, (\ref{eq:cA:cAprime}) follows from
(\ref{eq:Aprime:norm:upper:bound}). 
\qed\medskip

For the next statement, see \autoref{def:left:balanced:Y} for left-balanced $Y$-configurations.

\begin{corollary}
\label{cor:tau:u:bilipshitz}
There exists $\alpha > 0$ and 
for every $\delta > 0$ there exists a compact set $K'$ with $\nu(K') >
1 - c(\delta)$ with $c(\delta) \to 0$ as $\delta \to 0$ and for every
$\ell > 0$ there exists a compact set $K_\ell$ with $\nu(K_\ell) >
1-c'(\delta)$ with $c'(\delta) \to 0$ as $\delta \to 0$ 
such that the
following holds. Suppose $Y = Y_{pre}(q_1,u,\ell,\tau)$ and $Y_1 =
Y_{pre}(q_1,u_1,\ell,\tau_1)$ are $(1,\epsilon)$-left-balanced pre-$Y$
configurations. (Note that $Y$ and $Y_1$ share the entire bottom
leg). Suppose also that $Y \in K'$, $Y_1 \in K'$,  $q \in K_\ell$, 
and also that $0 \le s \leq \max(\tau,\tau_1)$
is such that $u_1q_1 \in \cB_{s}[u q_1]$. Then,   for $\ell$
sufficiently large depending on $\delta$ and $\epsilon$, 
\begin{equation}
\label{eq:cor:tau:u:biliphitz}
(1/4) \kappa^{-2} (\tau-s) - C(\delta)  <
\tau_1-s < 4 \kappa^2 (\tau-s) + C(\delta),
\end{equation}
where $\kappa$ is as in
\autoref{prop:good_norms}\autoref{good_norms_bilipschitz}. 
\end{corollary}
In particular, note that if $s\in[\tau,\tau_1]$ then this yields a bound on $|\tau-\tau_1|$ depending on $\delta$.
\begin{proof}
The conclusion is symmetric in the exchange of $\tau$ and
$\tau_1$ (up to modifying $C(\delta)$), so we may assume that $\tau_1 \le \tau$. 
Let $K$ be as in \autoref{lemma:swithching}. By the
  ergodic theorem, we can
choose $K' \subset K$ and $T_0 = T_0(\delta)$ such that for all $x' \in
K'$ and any $L > T_0(\delta)$, 
\begin{displaymath}
\{ t' \in [L,2L] \st g_{-t'} x' \in K\} \ge (1-\kappa^{-2}/4) L.
\end{displaymath}
Let $s' = \min(\tau_1,s)$.
Then, we can find $0 < t \leq s'$ such that setting $t':=\tau-t$ we have $g_{-t'}(g_{\tau}u q_1) = g_t u q_1 \in K$, and similarly $g_t u_1 q_1 \in K$ and furthermore with $L:=T_0(\delta)+(\tau-s')$ we have:
\[
	\tau - t = t' \leq (1+\kappa^{-2}/2)L
\]
and thus:
\begin{equation}
\label{eq:sprime:minus:t}
s'-t < (1+\kappa^{-2}/2)(\tau-s')+T_1(\delta).
\end{equation}
Since $\epsilon < 1$, we may
also ensure that $\|\cA_2(q_1,u,\ell,t)\| < c_1(\delta)$, and
$\|\cA_2(q_1',u',\ell,t)\| < c_1(\delta)$ where $c_1(\delta)$ is as
in \autoref{lemma:swithching}. 
By \autoref{rmk:choice_of_measurable_partition_and_du_smallness} and
\autoref{lemma:swithching} with $q_1' = q_1$, $\tau = t$ and
$u' = u_1$, we get
\begin{multline*}
C(\delta)^{-1} \|\cA_2(q_1, u, \ell, t)\| -
O(e^{-\alpha \ell}) \le 
  \|\cA_2(q_1, u_1, \ell, t) \| \\ \le C(\delta) \|\cA_2(q_1, u, \ell, t)\| +
O(e^{-\alpha \ell}),
\end{multline*}
By \autoref{lemma:A:future:bilipshitz}, and since $\|\cA_2(q_1,u,\ell,
\tau)\| = \epsilon$, we have
\begin{displaymath}
e^{-\kappa (\tau-t)} \epsilon \le
  \|\cA_2(q_1, u, \ell, t)\| \le e^{-\kappa^{-1} (\tau-t)} \epsilon.
\end{displaymath}
Similarly, since $\|\cA_2(q_1,u_1,\ell,
\tau_1)\| = \epsilon$, by \autoref{lemma:A:future:bilipshitz}, 
\begin{displaymath}
e^{-\kappa (\tau_1-t)} \epsilon \le
  \|\cA_2(q_1, u_1, \ell, t)\| \le e^{-\kappa^{-1} (\tau_1-t)} \epsilon.
\end{displaymath}
Combining the last three displayed equations, we get, for $\ell$
  sufficiently large depending on $\delta$ and $\epsilon$, 
\begin{displaymath}
\kappa^{-2} (\tau-t) - C(\delta)  < \tau_1 -t \le \kappa^2 (\tau-t) + C(\delta).
\end{displaymath}
In view of (\ref{eq:sprime:minus:t}), this implies
\begin{equation}
 \label{eq:tau:sprime:tau1:sprime}
  (1/2) \kappa^{-2}(\tau-s') - C(\delta) \le \tau_1 - s' \le 2
  \kappa^2 (\tau-s') + C(\delta). 
\end{equation}
Recall that $s' = \min(\tau_1,s)$. 
If $s' = s$, then it follows immediately from
(\ref{eq:tau:sprime:tau1:sprime}) that (\ref{eq:cor:tau:u:biliphitz})
holds. If $s' = \tau_1$, then it follows from
(\ref{eq:tau:sprime:tau1:sprime}) that
\begin{displaymath}
0 < \tau - s < \tau - s' < C(\delta),
\end{displaymath}
and thus (\ref{eq:cor:tau:u:biliphitz}) still holds.
\end{proof}



\subsection{Matching and volume distortion}
	\label{ssec:matching_and_volume_distortion}

\subsubsection{Partitions and $U^+$}
	\label{sssec:partitions_and_u_plus}
Recall that the measurable partition $\gB_0[\cdot]$ is defined in \autoref{ssec:measurable_partitions}.
For $\epsilon_0' > 0$, let\index{$B$@$\gB_0^{(\epsilon_0')}[x]$}
\begin{displaymath}
\gB_0^{(\epsilon_0')}[x] = \{ y \in \gB_0[x] \st d^u(y,\partial \gB_0[x]) \ge
\epsilon_0'\},
\end{displaymath}
and let\index{$B$@$\cB_0^{(\epsilon_0')}[x]$}
\begin{displaymath}
\cB_0^{(\epsilon_0')}[x] = \gB_0^{(\epsilon_0')}[x]  \cap U^+[x].
\end{displaymath}
Note that $\cB_0^{(\epsilon_0')}$ admits a ``Haar measure'' $| \cdot |$ that is defined in \autoref{sssec:measurable_partitions_and_subgroups_compatible_with_the_measure}.

\begin{lemma}
\label{lemma:parts:staying:close}
For every $\delta > 0$ and every $\epsilon_0' > 0$ sufficiently small
depending on $\delta$ there exists \index{$\epsilon_0$}$\epsilon_0 =
\epsilon_0(\delta,\epsilon_0')>0$ and a compact set
$K \subset X$ with $\nu(K) >
1-\delta$ so that the following holds:

Suppose $\epsilon < \epsilon_0$. 
Suppose $Y=Y_{pre}(q_1,u,\ell,\tau)$ and
$Y'=Y_{pre}(q_1',u',\ell,\tau)$ are bottom-linked
pre-$Y$-configurations such that $Y'$ left-shadows $Y$. 
Suppose $Y \in K$, and using the notation
\S\ref{sec:subsec:Yconf:notation}, $q' \in K$, $q'_{1/2} \in K$, $q_1'
\in K$. Suppose also that
\begin{displaymath}
d^{q_2}_{\cH,loc}(U^+[q_2], \phi_{\tau}^{-1}(U^+[q_2'])) \le \epsilon, 
\end{displaymath}
and
\begin{equation}
\label{eq:gB0:minus:close}
d^Q_{\cH}(\cB_0[q_1], \cB_0[q_1']) < \epsilon,  
\end{equation}
where \index{$d^Q_{\cH}$}$d^Q_{\cH}$ denotes the Hausdorff distance of the projection to
$Q$, with respect to the ambient distance $d^Q(\cdot, \cdot)$.

Let
\begin{displaymath}
A_{\tau} = U^+[ q_2] \cap \gB_0^{(2\epsilon_0')}[q_2],
\end{displaymath}
\begin{displaymath}
  A_{\tau}' = U^+[ q_2'] \cap \phi_{\tau}(\gB_0^{(\epsilon_0')}[q_2]).
\end{displaymath}
Then, 
\begin{displaymath}
\frac{|g_{-\tau}A_{\tau} |}{|\cB_0[q_1]|} \le
\kappa(\delta)\frac{ |g_{-\tau}A_{\tau}'|}{|\cB_0[q_1']|},
\end{displaymath}
for an appropriate constant $\kappa(\delta)$.
\end{lemma}

\subsubsection{Proof of \autoref{lemma:parts:staying:close}.}
In the proof of
this lemma we normalize the measure $|\cdot|$ on $U^+[q_1]$ so
that $|\cB_0[q_1]| = 1$ and similarly we normalize the measure
$|\cdot|$ on $U^+[q_1']$ so that $|\cB_0[q_1']|=1$.

Let $A_0 = g_{-\tau}A_{\tau}$, $A_0' = g_{-\tau}A_{\tau}'$. 
Let $\phi_{\tau}$ be the interpolation map. 
Let $\tilde{A}_{\tau} = \phi_{\tau}^{-1}(A'_{\tau})$. Then,
\begin{displaymath}
\tilde{A}_0 \equiv g_{-\tau}\tilde{A}_{\tau} = \phi_0^{-1}(A'_0).
\end{displaymath}
Since the Hausdorff distance between $U^+[q_2]$ and $\phi_{\tau}^{-1}(U^+[q_2'])$ is bounded by $\epsilon$, by \autoref{lemma:L:and:C}
there exists $h_{\tau} \in \bbG^{ssr}(\cW^u[q_2])$
with $\|L_{q_2}h_{\tau} - \id \|_{L\cW^u(q_2)} \le C(\delta) \epsilon$ 
and such that
$\phi_{\tau}^{-1}(U^+[q_2'])   = h_{\tau} U^+[q_2]$. 
There exits $\epsilon_0 = \epsilon_0(\delta,\epsilon_0') >0$
so that if $q\in K$ then for any $h\in \bbG^{ssr}(\cW^u[q])$ with
$\|L_q h - \id \|_{L\cW^u(q)} <C(\delta) \epsilon_0$, we have that
$$h\gB_0^{(2\epsilon_0')}[q] \subset \gB_0^{(\epsilon_0')}[q].$$
Then, since $\epsilon < \epsilon_0$, we have that
\begin{displaymath}
h_{\tau} A_{\tau} \subset \tilde{A_{\tau}}. 
\end{displaymath}
Let $h_0 = g_{-\tau} \circ h_{\tau} \circ g_{\tau} \in \bbG^{ssr}(\cW^u[uq_1])$.
Then, $h_0$ is exponentially close to the identity, and $h_0 A_0
\subset \tilde{A}_0 = \phi_0^{-1}(A_0')$. Thus,
\begin{displaymath}
  (\phi_0 \circ h_0) A_0 \subset A_0'.
\end{displaymath}
Recall that $\phi_0=\wp^{+}(q_1',uq_1')\circ\wp^-(q_1,q_1')\circ\wp^+(uq_1,q_1)$ where $\wp^\pm$ denotes the subresonant maps induced by the measurable connections $P^{\pm}$ on $L\cW^u$.
In particular, since all the points are in the compact set $K$ we can assume that $\phi_0$ has uniformly bounded Jacobian.
Therefore
$|A_0| \le \kappa(\delta) |A_0'|$ for an appropriate constant $\kappa(\delta)$.
\qed\medskip

The following is a trivial modification of 
\cite[Lemma~12.8]{EskinMirzakhani_Invariant-and-stationary-measures-for-the-rm-SL2Bbb-R-action-on-moduli-space}.

\begin{lemma}
\label{lemma:segments:equal:length}
Suppose $P$ and $P'$ are finite measure subsets of $\reals^n$ with
$|P|=|P'|$, admitting partitions 
\begin{displaymath}
P = \bigsqcup_{j=1}^N P_j, \quad P' = \bigsqcup_{j=1}^N P'_j, 
\end{displaymath}
Let  
$Q \subset P$ and $Q' \subset P'$ be subsets with $|Q| > (1-\delta)|P|$,
$|Q'| > (1-\delta)|P'|$. 

Suppose there exists $\kappa > 1$ such that for all $1 \le j \le N$, 
$$|P_j{\cap Q} | \le \kappa|P_j'|.$$  Then there exists $\hat{Q} \subset Q$ with $|\hat{Q}|
\ge (1-2 \kappa \delta)|P|$ such that for all   $1\le j\le N$, if 
$\hat{Q} \cap P_j \ne \emptyset$, then $Q' \cap P_j'  \ne \emptyset$. 
\end{lemma}
\subsubsection{Proof.} Let $J = \{ j \st P_j \cap Q \ne \emptyset\}$, and let 
$J' = \{ j \st Q' \cap P_j' \ne \emptyset \}$, and let
\begin{displaymath}
\hat{Q} = \{ x \in Q \st \text{ for all $j$ with $x \in P_j$, we have
  $j \in J'$.} \}
\end{displaymath}
Thus, if $x \in Q \setminus \hat{Q}$, then there exists $j \in J$ with
$x \in Q \cap P_j$ but $j \not\in J'$. 
Then, 
\begin{displaymath}
|Q\setminus \hat{Q}| \le  \sum_{j \in J\setminus J'} |Q \cap P_j| 
 \le \kappa  \sum_{j \not\in J'}
|P_j'| \le  \kappa  |(Q')^c|,   
\end{displaymath}
since if $j \not\in J'$ then $P_j' \subset (Q')^c$. Thus, $|Q\setminus
\hat{Q}| \le \kappa  \delta |P|$, and so $|\hat{Q}| \ge (1-2
\kappa  \delta) |P|$. 
\hfill \qed 



\subsection{Transport Lemmas}
	\label{ssec:transport_lemmas}

\begin{lemma}[$\cA_2$ and $u_*$ are compatible]
	\label{lemma:cA:ustar}
There exists $\alpha > 0$ and 
for every $\delta > 0$ there exists a subset $K \subset Q$ with
$\nu(K) > 1-\delta$ such that the following holds:
Suppose $Y = Y_{pre}(q_1, u, \ell,\tau)$, $Y_1 = Y_{pre}(q_1,u_1,\ell,\tau)$
and $Y' = Y_{pre}(q',u',\ell,\tau)$ are pre-$Y$-configurations. (Note
that $Y$ and $Y_1$ share the entire bottom leg). Suppose $Y'$  and $Y$
are bottom-linked, and also that $Y'$  left-shadows
$Y$. Suppose also that $u_1 q_1 \in \cB_{\tau}[u q_1]$, so in particular,
$Y_1$ also left-shadows $Y$. 

Furthermore, suppose $Y \in K$, $Y_1 \in K$ and also (using the
notation \S\ref{sec:subsec:Yconf:notation}) $q' \in K$, $q'_{1/2} \in
K$. 
Then,
\begin{displaymath}
\cA_2(q_1,u_1,\ell,\tau) F^s_q(q') = (g_{\tau} \circ u_*(uq_1,u_1q_1) \circ g_{-\tau}) \cA_2(q_1,u,\ell,\tau) F^s_q(q') +
O(e^{-\alpha \ell}). 
\end{displaymath}
\end{lemma}

\begin{proof}
By \autoref{prop:factorization:main},
\begin{displaymath}
\cA_2(q_1,u,\ell,\tau) F^s_q(q') = \bfj\big(\phi_{\tau}^{-1}\bigl(U^+[g_{\tau}q_1']
\bigr)\big) + O(e^{-\alpha \ell}). 
\end{displaymath}
Let $h \in  \bbG^{ssr}(\cW^u[q_{\tau} u q_1])$ be such that
\begin{displaymath}
\phi_{\tau}^{-1}\bigl(U^+[g_{\tau}q_1'] \bigr) = h U^+[g_{\tau}q_1]
\end{displaymath}
Then, also by \autoref{prop:factorization:main}, 
\begin{displaymath}
\cA_2(q_1,u_1,\ell,\tau) F^s_q(q') = \bfj\big((\phi'_{\tau})^{-1}\bigl(U^+[g_{\tau}q_1']
\bigr)\big) + O(e^{-\alpha \ell}),
\end{displaymath}
where $\phi'_{\tau}$ is the interpolation map with $u$ replaced by $u_1$.
Note that $g_{\tau} u_1 q_1 \in \cB_0(g_{\tau} u q_1)$.
By
\autoref{theorem:interpolation:does:not:move:points}, for $x\in  \cW^u_{loc}[g_\tau uq_1]$ we have 
$$d^{Q}(x, \phi_\tau(x)) \leq C(\delta) e^{-\alpha\ell}$$
and
$$d^{Q}(x, \phi_\tau'(x)) \leq C(\delta) e^{-\alpha\ell}$$
and so 
$$d^{Q}(x,(\phi_{\tau}')\inv\circ \phi_\tau(x))\leq C(\delta) e^{-\alpha\ell}$$
for $x\in   \cW^u_{loc}[g_\tau uq_1]$.  
Then 
by \autoref{lemma:L:and:C}
there exists $h_1 \in \bbG^{ssr}(\cW^u[q_{\tau} u q_1])$ with $\|L_x h_1 - I \|= O(e^{-\alpha \ell})$ such that 
\begin{displaymath}
(\phi_{\tau}')^{-1}\bigl(U^+[g_{\tau}q_1'] \bigr) 
=(\phi_{\tau}')^{-1}\circ \phi_\tau \bigl( h U^+[g_{\tau}q_1] \bigr) 
= h_1 h U^+[g_{\tau}q_1].
\end{displaymath}
Therefore,
\begin{displaymath}
\bfj(h_1 h U^+[g_{\tau}q_1]) = \bfj(h U^+[g_{\tau}q_1]) + O(e^{-\alpha \ell}).
\end{displaymath}
By the equivariance of $u_*$,  
\begin{displaymath}
g_{\tau} \circ u_*(uq_1,u_1q_1) \circ g_{-\tau} = u_*(g_{\tau} u q_1, g_{\tau} u_1
q_1). 
\end{displaymath}
By the definition of $u_*$, the map $u_*(g_{\tau} u q_1,g_{\tau} u_1 q_1):
\bH(g_{\tau} u q_1) \to \bH(g_{\tau} u_1 q_1)$ satisfies
\begin{displaymath}
u_*(g_{\tau} u q_1,g_{\tau} u_1 q_1)\, \bfj(h U^+[g_{\tau}q_1]) = \bfj(h
U^+[g_{\tau}q_1]) \in \bH(g_{\tau} u_1 q_1),
\end{displaymath}
where the vector $\bfj(h U^+[g_{\tau}q_1])$ on the left hand side is
considered to be an element of $\bH(g_{\tau} u q_1)$. Now the statement
follows.  
\end{proof}

For the next statement, see \autoref{def:left:balanced:Y} for left-balanced $Y$-configurations.
It is an enhancement of the previous result by allowing $u_1q_1$ to be in the larger partition atom $\cB_{\tau-L}[uq_1]$ (at the cost of requiring the $Y$-configurations to be left-balanced):
\begin{lemma}
\label{lemma:cA:transport}
For every $\delta > 0$ there exists a subset $K' \subset Q$ with
$\nu(K') > 1-c_1(\delta)$ with $c_1(\delta) \to 0$ as $\delta \to 0$
and for every $\ell > 0$ there exists a subset $K_\ell \subset Q$ with
$\nu(K_\ell) > 1-c_2(\delta)$ with $c_2(\delta) \to 0$ as $\delta \to 0$
such that if $\epsilon$ is sufficiently small depending on $\delta$
and $\ell$ is sufficiently large depending on $\delta$ and $\epsilon$, 
the following holds:

Suppose $Y = Y^{(\epsilon)}_{pre,lb}(q_1, u, \ell)$ and $Y_1 = Y^{(\epsilon)}_{pre,lb}(q_1,u_1,\ell)$
are $(1,\epsilon)$-left-balanced pre-$Y$-configurations. (Note
that $Y$ and $Y_1$ share the entire bottom leg). Suppose $Y'$ is a
pre-$Y$ configuration which is bottom-linked to $Y$ and left-shadows
$Y$. Furthermore, suppose $Y \in K'$, $Y_1 \in K'$ and also (using the
notation \S\ref{sec:subsec:Yconf:notation}) $q \in K_\ell$, $q' \in
K_\ell$, $q'_{1/2} \in K'$. 

Let $\tau = \tau(Y)$, $\tau_1 = \tau(Y_1)$. 
Suppose $0 \le L \ll \tau$ is such that $u_1q_1 \in \cB_{\tau-L}[u
q_1]$. (In particular, this implies that $Y_1$ shadows $Y$ up to time
$\tau-L$). Finally, suppose that  
\begin{displaymath}
\cA_2(q_1,u,\ell,\tau) F^s_q(q') = \bfv + \bfw,
\end{displaymath}
with $\bfw = O(e^{-\alpha \ell})$. Then,
\begin{multline}
\label{eq:lemma:cA:transport}
  \cA_2(q_1,u_1,\ell,\tau_1) F^s_q(q') = (g_{\tau_1} \circ u_*(uq_1,u_1q_1) \circ
g_{-\tau})  \bfv + \\
O(e^{-\alpha \ell+2\kappa^3 L}).  
\end{multline}
\end{lemma}

\begin{proof}
Let $K_\ell$ be as in \autoref{lemma:swithching} and
\autoref{cor:tau:u:bilipshitz}, and let
$K$ be the intersection of the sets $K$ of 
\autoref{lemma:cA:ustar} and \autoref{lemma:ustar:bounded} 
and the set $K'$ of \autoref{cor:tau:u:bilipshitz}. 

By the ergodic theorem
we can choose $K' \subset K$ with $\nu(K') > 1-c_1(\delta)$ with
$c_1(\delta) \to 0$ as $\delta \to 0$ and $L_0 = L_0(\delta) > 0$
such that for all $x \in K'$ and any $L > L_0(\delta)$,
\begin{displaymath}
\{ t \in [L,2L] \st g_{-t} x \in K\} \ge (1-\kappa^{-4}/4) L, 
\end{displaymath}
where $\kappa$ is as in
\autoref{prop:good_norms}\autoref{good_norms_bilipschitz}.
By \autoref{cor:tau:u:bilipshitz} applied with $s=\tau-L$ (whose conditions $u_1q_1\in \cB_s[uq_1]$ and $s\in[0,\max(\tau,\tau_1)]$ are satisfied by construction), we have
\begin{equation}
\label{eq:tau1:tau:tmp}
  \tfrac 16 \kappa^{-2} L < \tau_1 - (\tau-L) < 6 \kappa^2 L
\end{equation}
Then, we can find $t \in [L,2L]$ such that (taking $x=g_{\tau}uq_1, x_1=g_{\tau u_1q_1}\in K'$) we have 
$g_{\tau-t} u q_1 \in K$, $g_{\tau-t} u_1 q_1 \in K$.
By \autoref{lemma:ustar:bounded}, 
there exists $C = C(\delta)$ such that
\begin{displaymath}
\|u_*(g_{\tau-t} u q_1,g_{\tau-t} u_1 q_1) \| \le C(\delta).  
\end{displaymath}
Since $g_{-t}: \bH(g_\tau u q_1) \to \bH(g_{\tau-t} u q_1)$ is a
contraction, we have
\begin{equation}
\label{eq:u:star:tmp}
\| u_*(g_{\tau-t} u q_1,g_{\tau-t} u_1 q_1) \circ g_{-t} \| \le C(\delta).  
\end{equation}
Using \autoref{lemma:cA:ustar} for the first line
\begin{align}
\label{eq:cA:q1:u1:ell:tau:t}
\cA_2(q_1,u_1,\ell,\tau-t) F^s_q(q')
  & = u_*(g_{\tau-t} u q_1,g_{\tau-t} u_1 q_1) \circ
    \cA_2(q_1,u,\ell,\tau-t) F^s_q(q') + O(e^{-\alpha \ell}) \notag\\
  & = u_*(g_{\tau-t} u q_1,g_{\tau-t} u_1 q_1) \circ g_{-t}
    \cA_2(q_1,u,\ell,\tau) F^s_q(q') + O(e^{-\alpha \ell}) \notag \\
  & =  u_*(g_{\tau-t} u q_1,g_{\tau-t} u_1 q_1) \circ g_{-t} 
    \bfv  + O(e^{-\alpha \ell}), \notag \\ 
  & = (g_{\tau-t} \circ u_*(u_1 q_1, uq_1) \circ g_{-\tau}) 
    \bfv  + O(e^{-\alpha \ell}),  
\end{align}
where to pass from the second to the third line we have used
(\ref{eq:u:star:tmp}), and to pass to the last line we used the
equivariance of $u_*$. 
We have
\begin{equation}
\label{eq:cA:q1:u1:ell:tau1}
\cA_2(q_1,u_1,\ell,\tau_1) F^s_q(q') = g_{\tau_1 -(\tau-t)} \circ
\cA_2(q_1,u_1,\ell,\tau-t) F^s_q(q')
\end{equation}
In view of (\ref{eq:tau1:tau:tmp}), $0 < \tau_1 - (\tau-t) < 6 \kappa^2
L$. Thus, (\ref{eq:lemma:cA:transport})
follows from (\ref{eq:cA:q1:u1:ell:tau:t}),
(\ref{eq:cA:q1:u1:ell:tau1}) and \autoref{prop:good_norms}\autoref{good_norms_bilipschitz}.
\end{proof}

\section{Divergence Estimates, Time Changes, and Foliations}
	\label{sec:divergence_estimates_time_changes_and_foliations}

\paragraph{Outline of section}
The goal of this section is to adapt the constructions from \S8 and \S9 of \cite{EskinMirzakhani_Invariant-and-stationary-measures-for-the-rm-SL2Bbb-R-action-on-moduli-space}.


\subsection{Inert subspaces}
	\label{ssec:inert_subspaces}

\subsubsection{Setup}
	\label{sssec:setup_inert_subspaces}
Here we follow \S8 of \cite{EskinMirzakhani_Invariant-and-stationary-measures-for-the-rm-SL2Bbb-R-action-on-moduli-space}.
Recall that $\bold H$ was defined in \autoref{eq:bold_H_definition}.
For ease of notation, we denote the Lyapunov spectrum of $\bold H$ by $\{\lambda_i\}$ and use the indexing convention $\lambda_i>\lambda_{i+1}$.
Recall that we have the Oseledets decomposition\index{$H$@$\bold H^{\lambda_i}$}
\[
	\bold H = \oplus_{\lambda_i} \bold H^{\lambda_i}
\]
and associated Oseledets filtrations $\bold H^{\geq \lambda_i}, \bold H^{\leq \lambda_i}$.

\subsubsection{Definition of $U^+$-inert subspaces}
	\label{sssec:definition_of_u_inert_subspaces}
Recall that $\cB_0(x)$ is defined in \autoref{eq:U_plus_frakb_partition}.
Define now the following coarsening of the Oseledets filtration:\index{$H$@$\bold H^{\leq \lambda_i}_{in}(x)$}
\begin{multline*}
	\bold H^{\leq \lambda_i}_{in}(x):= 
	\{
	\bfv\in \bold H(x) \colon \text{ for Haar-a.e. }u\in \cB_0(x),
        \text{ we have } \\
	u_*(x,ux) \bfv \in \bold H^{\leq \lambda_i}(ux)
	\}
\end{multline*}
where $u_*(x,y)$ was defined in \autoref{ssec:the_change_of_basepoint_map}.
In other words, these are vectors which keep growing at rate at most $\lambda_i$ into the future, when pushed to any nearby point along $U^+[x]$.
From the definition it follows that $\bold H^{\leq \lambda_i}\subseteq \bold H^{\leq \lambda_{i-1}}$ (with $\lambda_i<\lambda_{i-1}$) but this filtration might well be trivial.

Define now the \emph{inert subbundles}:\index{$E$@$\bbE^{\lambda_i}$}
\[
	\bold E^{\lambda_i} := \bold H^{\leq \lambda_i}_{in} \cap \bold H^{\geq \lambda_i}
\]
Here are some basic properties, immediate from the definitions:
\begin{proposition}[Properties of inert subbundles]
\label{prop:properties:of:inert:subbundles}
	On a subset of $X$ of full measure:
	\begin{enumerate}
		\item The bundles $\bold E^{\lambda_i}$ are $g_t$-equivariant and so yield $g_t$-cocycles.
		\item We have that $\bold E^{\lambda_i}\subset \bold H^{\lambda_i}$.
		\item For Haar-a.e. $u\in U^+(x)$ we have that
		\[
			\bold E^{\lambda_i}(ux) = u_*(x,ux) \bold E^{\lambda_i}(x).
		\]
	\end{enumerate}
\end{proposition}
See \cite[\S8.1]{EskinMirzakhani_Invariant-and-stationary-measures-for-the-rm-SL2Bbb-R-action-on-moduli-space} for the proofs.
Define now the sum of inert subbundles:\index{$E$@$\bbE$}
\[
	\bbE = \bigoplus_{\lambda_i} \bbE^{\lambda_i} \subset \bbH.
\]
Note that the subspace corresponding to the top Lyapunov exponent is in $\bbE$.
The following is a direct translation of Prop.~8.5(a) of \cite{EskinMirzakhani_Invariant-and-stationary-measures-for-the-rm-SL2Bbb-R-action-on-moduli-space}.
\begin{proposition}[Convergence to $\bbE$]
	\label{prop:most:inert}
	For every $\delta>0$ there exists a compact set $K\subset Q$ of measure at least $1-\delta$, and a constant $L_2(\delta)$ with the following properties.
	Take $q\in K$ and $\bbv\in \bbH(q)$.
	For any $L'>L_2(\delta)$ there exists $t\in(L',2L')$ such that for a set of $u\in \cB_0(g_{-t}q)$ of Haar measure at least $(1-\delta)|\cB_0(g_{-t}q)|$ we have
	\[
		d\left(
		\frac{g_s \circ u_*(g_{-t}q, u g_{-t}q) \circ g_{-t}\bbv}{g_s \circ u_*(g_{-t}q, u g_{-t}q) \circ g_{-t}\bbv}
		,
		\bbE(g_s u g_{-t}q)
		\right)
		\leq C(\delta)e^{-\alpha\cdot t}
	\]
	where $s>0$ satisfies
	\[
		\norm{\bbv} = \norm{g_s \circ u_*(g_{-t}q, u g_{-t}q) \circ g_{-t}\bbv}
	\]
	and the exponent $\alpha$ depends only on the Lyapunov spectrum.
\end{proposition}
The statement says that by applying $g_{-t}$ to a vector $\bbv$, followed by a typical $u\in \cB_0(g_{-t}q)$, and then by $g_s$, the resulting vector will be close to the inert subspace.
We refer for the detailed proof to \cite[Prop.~8.5(a)]{EskinMirzakhani_Invariant-and-stationary-measures-for-the-rm-SL2Bbb-R-action-on-moduli-space}, which is based on the lemmas preceding it.

\begin{proposition}[Agreement of $u_*$ and $P^+$]
\label{prop:agreement:of:u:star:and:p}
	For a set of full measure of $q,q'$, if $q'=uq\in \cB_0[q]$, the restriction to inert subspaces of the map $u_*(q,q')\colon \bbE(q)\to \bbE(q')$ agrees with the standard measurable connection $P^+(q,q')$ as defined in \autoref{sssec:standard_measurable_connection}, when $P^+$ is restricted to inert subspaces as well.
\end{proposition}
The above result is Lemma~9.1 of \cite{EskinMirzakhani_Invariant-and-stationary-measures-for-the-rm-SL2Bbb-R-action-on-moduli-space}.
\begin{proof}
	By \autoref{prop:u_star_preserves_flags} the map $u_*$ preserves the backwards Lyapunov filtration and agrees with the standard measurable connection on the associated graded.
	On the other hand, for the inert subbundles the map $u_*$ has the property that $u_*(q,q') \bbE^{\lambda_i}(q) = \bbE^{\lambda_i}(q')$, by \autoref{prop:properties:of:inert:subbundles}.
	Therefore it has no off-diagonal components relative to the Oseledets decomposition.
	It follows that it must agree with $P^+$ on the inert subspaces $\bbE^{\lambda_i}$.
\end{proof}



\subsection{Jordan canonical form of a cocycle}
	\label{ssec:jordan_canonical_form_of_a_cocycle}

\subsubsection{Setup}
	\label{sssec:setup_jordan_canonical_form_of_a_cocycle}
Recall, from \autoref{thm:jordan_normal_form} 
that for any measurable cocycle $\bbE$, with Lyapunov decomposition $\bbE=\oplus \bbE^{\lambda_i}$, we have a further filtration in each Lyapunov block \index{$E$@$\bbE^{\lambda_i}_{\leq j}$}$\bbE^{\lambda_i}_{\leq j}\supsetneq \bbE^{\lambda_i}_{\leq j-1}$.
Furthermore, the cocycle is block-diagonal under the action of $g_t$, and for a metric which will be denoted $\norm{}$ we have that
\begin{align}
	\label{eq:growth_v_in_Elambda_i_j}
	\norm{g_t \bfv} = e^{\lambda_{ij}(q;t)}\norm{\bfv} \quad \forall \bfv\in \bbE^{\lambda_i}_{\leq j}/\bbE^{\lambda_i}_{\leq j-1}(q).
\end{align}
Note that by \autoref{prop:good_norms}(v), \index{$\lambda_{ij}(-;t)$}$\lambda_{ij}(-;t)$ is constant on the atoms of the partition $\gB_0$.
Furthermore, they are monotonic in $t$ by
\autoref{prop:good_norms}(iii) and the cocycle property of
$\lambda_{ij}$. Furthermore, \autoref{prop:good_norms}(iii) and the
cocycle property of $\lambda_{ij}$ imply that
\begin{align}
	 \label{eq:bilip:lambda:ij}
	\frac{1}{\kappa}|t-t'| \leq 
	| \lambda_{ij}(q;t) - \lambda_{ij}(q;t')|
	\leq
	{\kappa}|t-t'|
\end{align}
where $\kappa$ depends only on the Lyapunov spectrum.

Define also\index{$E$@$\bbE^{\lambda_i}_j$}

\[
	\bbE^{\lambda_i}_{j}\text{ as the orthogonal complement of }\bbE^{\lambda_i}_{\leq j-1}
	\text{ inside }
	\bbE^{\lambda_i}_{\leq j}
\]
with respect to the metric $\norm{-}$ introduced in \autoref{prop:good_norms}.
In \cite[\S9.1]{EskinMirzakhani_Invariant-and-stationary-measures-for-the-rm-SL2Bbb-R-action-on-moduli-space} this is denoted by $\bbE_{ij}'$.
The notation $\bbE^{\lambda_i}_j$ will only be needed for this section.

We will denote by \index{$\Lambda''$}$\Lambda''$ the ``fine'' Lyapunov spectrum consisting of pairs $(i,j)\in \Lambda''$ such that $\lambda_i$ is in the Lyapunov spectrum of the cocycle and $j=1\ldots n_i$ denotes the index of the sub-block corresponding to $\lambda_i $ in the Jordan normal form of the cocycle.



\subsubsection{Time changes}
	\label{sssec:time_changes}
The functions $\lambda_{ij}(q;t)$ are monotonic in $t$ by \autoref{prop:good_norms}\autoref{item:lambda_ij_bilipschitz} combined with the cocycle property of $\lambda_{ij}(q;t)$.
Therefore we can define the time-changed flow \index{$g^{ij}_t$}$g^{ij}_t$ with the property that the expansion factor in $\bbE^{\lambda_i}_{\leq j}/\bbE^{\lambda_i}_{\leq j-1}$ is constantly $\lambda_i$, i.e.
\[
	\norm{g_{t}^{ij}v} = e^{\lambda_i t}\norm{v}.
\]

We record the following elementary fact:

\begin{proposition}[Backwards nesting of partitions]
	\label{prop:backwards_nesting_of_partitions}
	Suppose that $q'\in \gB_0[q]$.
	Then for any $\lambda_i$ and $j\in[1,n_i]$ we have that
	\[
		g_{-t}^{ij}q'\in \gB_0[g_{-t}^{ij}q].
	\]
\end{proposition}
This is immediate from the constructions, see also \cite[Lemma 9.2]{EskinMirzakhani_Invariant-and-stationary-measures-for-the-rm-SL2Bbb-R-action-on-moduli-space}.



\subsection{The foliations \texorpdfstring{$\cF_{ij}$}{Fij} and parallel transport \texorpdfstring{$\bbR(q,q')$}{\bbR(q,q')}}
	\label{ssec:the_foliations_cf_ij}

In this section, we define and establish basic properties of certain partitions $\cF_{ij}$ of the center-unstable manifolds that are adapted to later arguments.

\subsubsection{The map $\bold R$}
	\label{sssec:the_map_bold_r}
From now on, $X_0\subset X$ is a set of full measure on which the objects used below are defined.
For $q \in X_0$, let
\begin{displaymath}
\index{$G[q]$}G[q] = \{ g_s u g_{-t} q 
\st t \ge 0, s \ge 0, u \in \cB_0(g_{-t} q) \}
\subset X.
\end{displaymath}
For $q' =g_s u g_{-t} q \in G[q]$, let \index{$R$@$\bbR(q,q')$}
\begin{align}
  \label{eq:def:R}
\bbR(q,q') =  g_s\circ u_*(g_{-t} q, g_{-s} q') \circ g_{-t} \colon \bbH(q)\to \bbH(q')
.
\end{align}
Here $g_s$ is the cocycle dynamics on $\bbH$ and $u_*(g_{-t} q,
g_{-s} q'): 
\bbH(g_{-t} q) \to \bbH(g_{-s} q')$ is as in
\autoref{def:u_star}.
The map $\bbR(q,q')$ does not depend on the choice of $s,t,u$ by the properties listed in \autoref{rmk:equivariance_and_endpoint_dependence_of_u_star}.
We will typically consider $\bbR(q,q')$ as a map from $\bbE(q) \to \bbE(q')$.  

\subsubsection{Growth of vectors}
	\label{sssec:growth_of_vectors}
By \autoref{eq:growth_v_in_Elambda_i_j}
for any $\bfv \in \bbE^{\lambda_i}_{j}(q)$, and any
$q' = g_s u g_{-t} q \in G[q]$, there exist $\bfv' \in \bbE^{\lambda_i}_{j}(q')$ with $\|\bfv' \| = \|\bfv\|$ and also $\bfv'' \in \bbE^{\lambda_i}_{\leq j-1}(q')$ such that:
\begin{equation}
	\label{eq:Rxy:bfv}
	\bbR(q,q') \bfv = e^{\lambda_{ij}(q,q')} \bfv' + \bfv''
\end{equation}
where:
\begin{equation}
\label{eq:def:lambdaij:xy}
\index{$\lambda_{ij}(q,q')$}\lambda_{ij}(q,q') = \lambda_{ij}(q;-t) + \lambda_{ij}(u g_{-t}q; s).
\end{equation}

\subsubsection{The foliations $\cF_{ij}$}
	\label{sssec:the_foliations_cf_ij}
For $q \in X_0$ and $\lambda_i,j$ as in \autoref{sssec:time_changes}, let 
\index{$F$@$\cF_{ij}[q]$}$\cF_{ij}[q]$ denote the set of $q' \in G[q]$ such that there exists
$\ell \ge 0$ with the property that
\begin{equation}
\label{eq:def:Fv}
	g^{ij}_{-\ell}q' \in \cB_0[g^{ij}_{-\ell} q]. 
\end{equation}
By \autoref{prop:backwards_nesting_of_partitions}, if \autoref{eq:def:Fv} holds
for some $\ell$ then it also holds for any larger $\ell$. 
Alternatively,
\begin{displaymath}
\cF_{ij}[q] = \left\lbrace g_{\ell}^{ij} u g_{-\ell}^{ij} q \colon \ell \ge 0, \quad u \in
\cB(g_{-\ell}^{ij} q) \right \rbrace \subset X. 
\end{displaymath}
By \autoref{eq:def:lambdaij:xy} it follows that 
\begin{equation}
\label{eq:lambdaij:zero:on:Fij}
\lambda_{ij}(q,q') = 0 \qquad \text{ if $q' \in \cF_{ij}[q]$.}
\end{equation}
We will refer to the sets $\cF_{ij}[q]$ as {\em leaves}.
Locally at $q$, the leaf $\cF_{ij}[q]$ through $q$ is a piece of $U^+[q]$.
More precisely, for $q'\in \cF_{ij}[q]$:
\begin{align}
	\label{eqn:F_ij_pieces_of_Uplus}
	\left(\cF_{ij}[q] \cap \gB_0[q']\right) =  U^+[q'] \cap \gB_0[q'].
\end{align}
Note that for any compact subset $A \subset \cF_{ij}[q]$ there exists a large enough $\ell$ such that $g^{ij}_{-\ell}(A)$ is contained in a set of
the form $\cB_0[z] \subset U^+[z]$ for some $z\in Q_0$.
Then the same holds for $g^{ij}_{-t}(A)$ and any $t > \ell$. 

\subsubsection{Measures on leaves}
	\label{sssec:measures_on_leaves}
By \autoref{sssec:the:measure:on:Uplus:orbits}
the sets $\cB_0[q]$ support a ``Lebesgue measure'' $| \cdot |$, 
i.e. the pushforward of the Haar measure on $U^+(q)/U^+_q(q)$ to $\cB_0[q]$ under the action map $u \to u q$ 
(recall that $U^+_q(q)$ is the stabilizer of $q$ in $U^+(q)$).

It follows that the leaves $\cF_{ij}[q]$ also support a Lebesgue
measure (defined up to normalization), which will also be denoted by 
\index{$\norm$@$\abs{\bullet}$}$| \cdot |$.
Specifically, if $A \subset \cF_{ij}[q]$ and $B \subset \cF_{ij}[q]$
are compact subsets, define
\begin{equation}
\label{eq:def:lebesgue:leaf}
	\frac{|A|}{|B|} \equiv \frac{|g^{ij}_{-\ell}(A)|}{|g^{ij}_{-\ell}(B)|}, 
\end{equation}
where $\ell$ is chosen large enough such that both
$g^{ij}_{-\ell}(A)$ and $g^{ij}_{-\ell}(B)$ are contained in a
set of the form $\cB_0[z]$, for some $z \in X_0$.
It is immediate to check that replacing $\ell$ by a larger number, the right-hand-side of \autoref{eq:def:lebesgue:leaf} is not changed.

\subsubsection{Balls in $\cF_{ij}$}
	\label{sssec:balls_in_cf_ij}
Define the ``balls'' $\cF_{ij}[q,\ell] \subset \cF_{ij}[q]$ by 
\begin{equation}
\label{eq:def:Fij:ball}
	\index{$F$@$\cF_{ij}[q,\ell]$}\cF_{ij}[q,\ell] = \{ q' \in \cF_{ij}[q] \colon g^{ij}_{-\ell} q' \in \cB_0[
g^{ij}_{-\ell} q] \}.
\end{equation}
These have the property that, for $q\in X_0$ and $q'\in \cF_{ij}[q]$, and sufficiently large $\ell$:
\begin{displaymath}
	\label{lemma:symmetric:difference:balls}
	\cF_{ij}[q,\ell] = \cF_{ij}[q',\ell]. 
\end{displaymath}
Indeed, for $\ell$ large enough we have that $g^{ij}_{-\ell} q' \in \cB_0[g^{ij}_{-\ell} q]$, 
and then $\cB_0[g^{ij}_{-\ell} q'] = \cB_0[g^{ij}_{-\ell} q]$.

\subsubsection{The flows $g^\bfv_t$}
	\label{sssec:the_flows_g_v_t}
Let $q \in X_0$ and $\bfv \in \bbE(q)$.
Define now a time change \index{$\tau$@$\tau_\bfv(q;t)$}$\tau_\bfv(q;t)$ and reparametrized flow by:
\[
	\index{$g^\bfv_t$}g^\bfv_t q = g_{\tau_\bfv(q;t)} q\quad 
	\text{ such that }
	\quad
	\|(g^\bfv_t)_* \bfv\|_{g^{\bfv}_t q} = e^t \|\bfv\|_q.
\]
(Note that we are not defining $g_t^{\bfv} q'$ for $q' \ne q$).
The composition rule for the flow is then:
\begin{displaymath}
g_{t+s}^{\bfv} q = g_s^{(g_t)_* \bfv} g_t^{\bfv} q.
\end{displaymath}
For $q' \in G[q]$ and $\ell \in \bR$, define
\begin{equation}
\label{eq:def:tilde:g}
\index{$g$@$\tilde{g}^{\bfv,q}_{-\ell}$}\tilde{g}^{\bfv,q}_{-\ell}
:= g_{-\ell}^{\bfw} q', \quad \text{ where $\bfw =
  \bbR(q,q')\bfv$. }
\end{equation}
(When there is no potential for confusion about the point $q$ and the
vector $\bfv$ used, we
denote $\tilde{g}^{\bfv,q}_{-\ell}$ by $\index{$g$@$\tilde{g}_{-\ell}$}\tilde{g}_{-\ell}$.)
Note that \autoref{prop:backwards_nesting_of_partitions} still holds if
$g^{ij}_{-t}$ is replaced by $\tilde{g}^{\bfv,q}_{-t}$.

\subsubsection{The foliations $\cF_\bfv$}
	\label{sssec:the_foliations_cf_bfv}
For $\bfv \in \bbE(q)$ we
can define the foliations 
\index{$F$@$\cF_{\bfv}[q]$}$\cF_{\bfv}[q]$
and the ``balls'' 
\index{$F$@$\cF_{\bfv}[q,\ell]$}$\cF_{\bfv}[q,\ell]$ as in
\autoref{eq:def:Fv} and
\autoref{eq:def:Fij:ball}, with $\tilde{g}^{\bfv,q}_{-t}$
playing the role of $g^{ij}_{-t}$.

For $q' \in \cF_{\bfv}[q]$, we have
\begin{displaymath}
	\cF_{\bfv}[q] = \cF_{\bfw}[q'], \qquad \text{ where } \bfw = \bbR(q,q') \bfv. 
\end{displaymath}
We can 
define the measure (up to normalization) \index{$\norm$@$\abs{\bullet}$}$| \cdot |$ on
$\cF_{\bfv}[q,\ell]$ as in \autoref{eq:def:lebesgue:leaf}.
Then \autoref{lemma:symmetric:difference:balls} holds for $\cF_{\bfv}[q]$
without modifications.
The following follows immediately from the
construction: 
\begin{proposition}
	\label{lemma:Fv:norm}
For $q\in X_0$ and any $\bfv \in \bbE(q)$, for a.e.\ $q' \in
\cF_\bfv[q]$, we have that
	\begin{displaymath}
		\|\bbR(q,y) \bfv \|_{q'} = \|\bfv\|_q. 
	\end{displaymath}
\end{proposition}


\subsection{A maximal inequality}
\begin{lemma}
\label{lemma:fake:ergodicity:Fij}
Suppose $K \subset X$ with $\nu(K) > 1-\delta$. Then, for any
$\theta' > 0$ there exists a
subset $K^* \subset X$ with $\nu(K^*) > 1-2\kappa^2
\delta/\theta'$ such that 
for any $x \in K^*$ and any $\ell > 0$, 
\begin{equation}
\label{eq:fake:ergodicity:Fij}
|\cF_{ij}[x,\ell] \cap K| > (1-\theta') |\cF_{ij}[x,\ell]|. 
\end{equation}
Furthermore, if $x\in K^*$ then $\cB_0[x]:=\left(U^+[x]\cap \gB_0[x]\right)\subset K^*$.
\end{lemma}

\begin{proof}
Note that if $x$ satisfies \autoref{eq:fake:ergodicity:Fij} and $x'\in \cB_0[x]$ then $x'$ also satisfies the same equation, since $\cF_{ij}[x,\ell]$ consists of pieces of $U^+$-orbits as in \autoref{eqn:F_ij_pieces_of_Uplus}, so the last asserted property of $K^*$ can be ensured by saturating it accordingly.

For $t > 0$ let now
\begin{displaymath}
  \cB_t^{ij}[x] = g_{-t}^{ij}
  (\gB_0[g_t^{ij} x] \cap U^+[g_t^{ij} x]) = \cB_\tau[x], 
\end{displaymath}
where $\tau$ is such that $g_\tau x = g_t^{ij} x$. 
Let $s > 0$ be arbitrary. Let $K_s = g_{-s}^{ij} K$. Then
$\nu(K_s) > 1-\kappa \delta$.
Then, by 
\autoref{lemma:cB:vitali:substitute}, there exists a subset $K_s'$ with
$\nu(K_s') \ge (1-2 \kappa \delta/\theta')$ such that
for $x \in K_s'$ and all $t > 0$, 
\begin{displaymath}
| K_s \cap \cB_t^{ij}[x] | \ge (1-\theta'/2) |\cB_t^{ij}[x]|. 
\end{displaymath}
Let $K_s^* = g_s^{ij} K_s'$, and 
note that $g_s^{ij}\cB_t^{ij}[x] = \cF_{ij}[g_s^{ij} x, s-t]$. Then, for all 
$x \in K_s^*$ and all $0 < s - t < s$, 
\begin{displaymath}
|\cF_{ij}[x,s-t] \cap K| \ge (1-\theta'/2) |\cF_{ij}[x,s-t]|. 
\end{displaymath}
We have $\nu(K_s^*) \ge (1-2\kappa^2 \delta/\theta')$. Now take a
sequence $s_n \to\infty$, and let $K^*$
be the set of points which are in infinitely many $K_{s_n}^*$. 
\end{proof}


\subsection{Bounded subspaces, Synchronized Exponents, and the Equivalence Relation}
	\label{ssec:bounded_subspaces_synchronized_exponents_and_the_equivalence_relation}

\paragraph{Outline of section}

In this section, we recall from \cite{EskinMirzakhani_Invariant-and-stationary-measures-for-the-rm-SL2Bbb-R-action-on-moduli-space} Prop.~10.1, Prop.~10.2, and Prop.~10.3.
The rest of \S10 of loc.cit. contains the proofs of these statements.

Recall from \autoref{sssec:setup_jordan_canonical_form_of_a_cocycle} that $\Lambda''$ indexes the ``fine Lyapunov spectrum'' on $\bbE$.
\begin{proposition}[Subbundles $\bbE_{[ij],bdd}$]
	\label{prop:subbundles_e_ij_bdd}
	\index{$E$@$\bbE_{[ij],bdd}(q)$}
	There exists a family of subcocycles $\bbE_{[ij],bdd}\subset \bbE$, indexed by equivalence classes $[ij]\subset \Lambda''$ that partition $\Lambda''$, such that their direct sum
	\[
		\bigoplus_{[ij]}\bbE_{[ij],bdd}\text{ maps injectively to }\bbE
	\]
	and furthermore the bundles satisfy \autoref{prop:some:fraction:bounded}, \autoref{prop:ej:bdd:transport:bounded}, and \autoref{prop:Rxy:v:small:implied:sync:bounded}.
\end{proposition}
\noindent In fact, the next three propositions characterize $\bbE_{[ij],bdd}$.
Let \index{$\Lambda$@$\wtilde{\Lambda}$}$\wtilde{\Lambda}$ denote the set of equivalence classes $[ij]\subset \Lambda''$ as above.


\begin{proposition}
	\label{prop:some:fraction:bounded}
There exists $\theta > 0$ depending only on $\nu$ and $n \in \bN$ depending only on the dimension of $Q$ such that the
following holds:
for every $\delta > 0$, {every $T > 0$} and every $\eta > 0$, 
there exists a subset $K=K(\delta,T,\eta)$ of measure at least
$1-\delta$ and $L_0 = L_0(\delta,T,\eta) > 0$ 
such that the following holds:

Suppose $q \in X$, 
$\bfv \in \bbE(q)$, $L \ge L_0$, 
and 
\begin{displaymath}
| g_{{[-T,T]}} K \cap \cF_\bfv[q,L]| \ge (1-(\theta/2)^{n+1})
|\cF_\bfv[q,L]|. 
\end{displaymath}
Then, for at least $(\theta/2)^n$-fraction of $q' \in \cF_\bfv[q,L]$,
\begin{displaymath}
d\left(\frac{\bbR(q,q') \bfv}{\|\bbR(q,q') \bfv\|},
    \bigcup_{[ij] \in \tilde{\Lambda}} \bbE_{[ij],bdd}(q') \right) < \eta. 
\end{displaymath}
\end{proposition}

In
\cite{EskinMirzakhani_Invariant-and-stationary-measures-for-the-rm-SL2Bbb-R-action-on-moduli-space}
Prop.~10.1, it is assumed that $T=1$. However, the same proof also
shows this more general version (in fact, in the paragraph labeled
``Proof of Proposition~10.1'', one should replace $L_{n-1}+1$ by
$L_{n-1}+T$ and $L_1+1$ by $L_1+T$; no other changes are needed).


\begin{proposition}
\label{prop:ej:bdd:transport:bounded}
There exists a measurable function $C_3\colon X \to \bR^+$ finite almost
everywhere such that for all $q \in Q$, for all $q' \in
\cF_{ij}[q]$, for all $\bfv \in
\bbE_{[ij],bdd}(q)$, we have
	\begin{displaymath}
		C_3(q)^{-1}C_3(q')^{-1} \|\bfv\| \le \|\bbR(q,q') \bfv \| \le 
		C_3(q)C_3(q') \|\bfv\|.  
	\end{displaymath}
\end{proposition}

\begin{proposition}
\label{prop:Rxy:v:small:implied:sync:bounded}
There exists \index{$\theta$}$\theta > 0$ (depending only on $\nu$) and a subset
$X_0 \subset X$ with $\nu(X_0) = 1$ such that the following holds:

Suppose $q \in X_0$, $\bfv \in \bbH(q)$, 
and there exists $C > 0$ such that for all $\ell
> 0$, and at least
$(1-\theta)$-fraction of $q' \in \cF_{ij}[q,\ell]$, 
\begin{displaymath}
	\| \bbR(q,q') \bfv \| \le C \|\bfv \|.  
\end{displaymath}
Then, $\bfv \in \bbE_{[ij],bdd}(q)$. 
\end{proposition}

We refer to the discussion after Prop.~10.3 in \cite{EskinMirzakhani_Invariant-and-stationary-measures-for-the-rm-SL2Bbb-R-action-on-moduli-space} for the technical issues that the above propositions address.
All the proofs of the above propositions are contained in \S10 of loc.cit. and, since they take place in the $g_t$-orbit of a single unstable, carry over to our setting with virtually no modifications.

\subsection{A variant of \autoref{prop:Rxy:v:small:implied:sync:bounded}.}

We will need the following:
\begin{proposition}
\label{prop:variant:Rxy:v:small:implied:sync:bounded}
There exists \index{$\theta$}$\theta > 0$ (depending only on $\nu$) $\epsilon > 0$ 
and a subset
\index{$X_0$}$X_0 \subset X$ with $\nu(X_0) = 1$ such that the following holds:

Suppose $q \in X_0$, $\bfv \in \bbH(q)$, 
and there exists $C > 0$ such that for all $L > 0$, for at least
$(1-\epsilon)$-fraction of $\ell \in [0,L]$, and at least
$(1-\theta)$-fraction of $q' \in \cF_{ij}[q,\ell]$, 
\begin{displaymath}
	\| \bbR(q,q') \bfv \| \le C \|\bfv \|.  
\end{displaymath}
Then, $\bfv \in \bbE_{[ij],bdd}(q)$. 
\end{proposition}

\begin{definition}
\label{def:bounded:subspace}
Suppose $x \in X$, $\epsilon > 0$, $\theta > 0$, $ij \in \Lambda'$. 
A vector $\bfv\in \bH(x)$ is called {\em $(\epsilon,\theta,ij)$-bounded} if there
exists $C_x < \infty$ such that for a set $\cD_x \subset \reals^+$ of
asymptotic density at least $(1-\epsilon)$, 
for all $\ell \in \cD_x$ and for ($1-\theta$)-fraction
of $y \in \cF_{ij}[x,\ell]$,
\begin{equation}
\label{eq:def:thetabounded}
\|\bbR(x,y)\bfv \| \le C_x\|\bfv\|.
\end{equation}
\end{definition}

\begin{lemma}
\label{lemma:bounded:subspace}
Let $n = \dim \bH(x)$ (for a.e $x$). If there exists no non-zero
$(\epsilon/n^2,\theta/n^2,ij)$-bounded vectors in $\bH(x)$,
we set \index{$V$@$\bbV_{ij}(x)$}$\bbV_{ij}(x) = \{0\}$. Otherwise, we define 
$\index{$V$@$\bbV_{ij}(x)$}\bbV_{ij}(x) \subset \bH(x)$ to
be the linear span of the $(\epsilon/n^2,\theta/n^2,ij)$-bounded vectors in $\bH(x)$.
This is a subspace of $\bH(x)$, and any vector in this
subspace is $(\epsilon,\theta,ij)$-bounded. Also, 
\begin{itemize}
\item[{\rm (a)}] $\bbV_{ij}(x)$ is $g_t$-equivariant, i.e.\ $g_t
  \bbV_{ij}(x) = \bbV_{ij}(g_t x)$. 
\item[{\rm (b)}] For almost all $u \in \cB(x)$, $\bbV_{ij}(ux) =
  (u)_* \bbV_{ij}(x)$. 
\end{itemize}
\end{lemma}

\begin{proof}
This is virtually the same as the proof of
\cite[Lemma~10.6]{EskinMirzakhani_Invariant-and-stationary-measures-for-the-rm-SL2Bbb-R-action-on-moduli-space}.
The assertion (b) follows immediately from the definition since for $u
\in \cB(x)$, $\cF_{ij}[ux,\ell] = \cF_{ij}[x,\ell]$. 
\end{proof}

\begin{remark}
\label{remark:def:theta}
Formally, from its definition, the subspace
$\bbV_{ij}(x)$ depends on the choice of $\epsilon$ and $\theta$.
It is clear that
as we decrease $\theta$ or $\epsilon$,
the subspace $\bbV_{ij}(x)$ decreases. 
In view of \autoref{lemma:bounded:subspace}, there
exist $\epsilon_0 > 0$,  $\theta_0 > 0$ and $m \ge 0$ 
such that for all $\epsilon < \epsilon_0$, all $\theta < \theta_0$ and almost
all $x \in X$, the dimension of $\bbV_{ij}(x)$ is $m$. We will
always choose $\epsilon$ and $\theta$ so that $\epsilon < \epsilon_0$
and $\theta < \theta_0$. 
\end{remark}
\medskip

\begin{lemma}
\label{lemma:switching:Rxy:Ryx}
Suppose $\bbV$ is a subbundle of
$\bH$, which is $g_t$ and $(u)_\ast$-equivariant. 
Suppose that $\theta > 0$, $\delta > 0$ and $\ell > 0$, and
suppose there exists a set $K_\ell \subset X$
with $\nu(K_\ell) > 1-\delta$ and a constant $C_1 < \infty$ (depending
on $\delta$ and $\theta$ but not on $\ell$, such that
such that for all $x \in K_\ell$, 
and at least $(1-\theta)$-fraction of $y \in
\cF_{ij}[x,\ell]$, 
\begin{equation}
\label{eq:tmp:Rxy:v:small}
\|\bbR(x,y) \bfv\| \le C_1 \|\bfv \| \quad\text{ for all $\bfv \in \bbV(x)$.}
\end{equation}
Then there exists a subset $K''(\ell) \subset X$
with $\nu(K''(\ell)) > 1-c(\delta)$ where $c(\delta) \to 0$ as
$\delta \to 0$, and there exists $\theta'' = \theta''(\theta,\delta)$ with
$\theta'' \to 0$ as $\theta \to 0$ and $\delta \to 0$
such that for all $x \in
K''(\ell)$, for at least $(1-\theta'')$-fraction of $y \in
\cF_{ij}[x,\ell]$, 
\begin{equation}
\label{eq:tmp:Rxy:v:twosided:small}
C_1^{-1} \|\bfv\| \le \|\bbR(x,y) \bfv\| \le C_1 \|\bfv \| \quad\text{
  for all $\bfv \in \bbV(x)$}.
\end{equation}
\end{lemma}

\begin{proof}
This is the same as the proof of
  \cite[Lemma~10.16]{EskinMirzakhani_Invariant-and-stationary-measures-for-the-rm-SL2Bbb-R-action-on-moduli-space}
\end{proof}

Recall that any bundle is measurably trivial. 
As in
\cite[Lemma~10.13]{EskinMirzakhani_Invariant-and-stationary-measures-for-the-rm-SL2Bbb-R-action-on-moduli-space},
we make let \index{$\mu$@$\tilde{\mu}_\ell$}$\tilde{\mu}_\ell$ be the measure on $X \cross \bP(\bbV)$ defined by
\begin{equation}
\label{eq:def:tilde:mu}
\tilde{\mu}_\ell(f) = \int_X \int_{\bP(\bbV)}
\frac{1}{|\cF_{ij}[x,\ell]|} \int_{\cF_{ij}[x,\ell]}
f(x,\bbR(y,x)\bfv) \, dy \, d\rho_0(\bfv) \, d\nu(x)
\end{equation}
where \index{$\rho_0$}$\rho_0$ is the ``round'' measure on $\bP(\bbV)$. 
(In fact, $\rho_0$ can be any measure on $\bP(\bbV)$ in the measure class of
Lebesgue measure, independent of x and fixed once and for all).

The following is analogous to \cite[Lemma~10.17]{EskinMirzakhani_Invariant-and-stationary-measures-for-the-rm-SL2Bbb-R-action-on-moduli-space}:
\begin{lemma}
\label{lemma:Vbdd:undependent:u}
Suppose $\bbV(x) = \bbV_{ij}(x)$. Then
there exists a function $C: X 
\to \reals^+$ finite almost
everywhere such that for all $x \in X$, 
all $\bfv \in \bbV_{ij}(x)$, 
and all $y \in \cF_{ij}[x]$, 
\begin{displaymath}
C(x)^{-1}C(y)^{-1}\|\bfv\| \le \|\bbR(x,y) \bfv\| \le C(x)C(y)\|\bfv\|,
\end{displaymath}
\end{lemma}

\begin{proof}
Let $\epsilon_0, \theta_0$ be as in Remark~\ref{remark:def:theta}. 
Choose $\epsilon < \epsilon_0^2$, $\theta < \theta_0$. 
Let $\cD_x$ be as in the definition of $\bbV_{ij}(x)$,
so that for almost all $x$
By Fubini's theorem and the definition of $\bbV_{ij}$ there exists
$\cD \subset \reals^+$ of asymptotic density at least $(1-\epsilon^{1/2})$
such that for $\ell
\in \cD$, the set $K'_\ell = \{ x \in X \st \ell \in \cD_x \}$
satisfies $\nu(K'_\ell) > 1-\epsilon^{1/2}$.
There exists a constant $C_1 > 0$ and a compact set $\hat{K}\subset X$ with $\nu(\hat{K}) <
1-\epsilon^{1/2}$ such that for $x \in K$, $C_x < C_1$, where $C_x$ is
as in \autoref{def:bounded:subspace}. Let $\delta = 2
\epsilon^{1/2}$, and let $K_\ell = \hat{K} \cap K'_\ell$. Then,
$\nu(K_\ell) > 1-\delta$, and for $x \in K_\ell$ and $\ell \in \cD$,
(\ref{eq:tmp:Rxy:v:small}). Then, we can apply
\autoref{lemma:switching:Rxy:Ryx} to get $K''(\ell)$ and $\theta'' >
0$ such that for $\ell \in \cD$, $x \in K''(\ell)$ and at least
$(1-\theta'')$-fraction of $y \in \cF_{ij}[x,\ell]$,
(\ref{eq:tmp:Rxy:v:twosided:small}) holds. 

Choose a sequence $\ell_k \to \infty$ such that $\ell_k \in \cD$, and 
$\tilde{\mu}_{\ell_k}$ converges to a measure $\tilde{\mu}_\infty$. The rest
of the proof of
\cite[Lemma~10.17]{EskinMirzakhani_Invariant-and-stationary-measures-for-the-rm-SL2Bbb-R-action-on-moduli-space}
goes through without modifications. 
\end{proof}

\begin{lemma}
\label{lemma:Vij:is:Eij}
We have for a.e. $x \in X$, $\bbV_{ij}(x) = \bbE_{[ij],bdd}(x)$.   
\end{lemma}

\begin{proof}
Let $\theta_0$ be as in
\autoref{prop:Rxy:v:small:implied:sync:bounded}. Choose
$\theta < \theta_0$, and choose $0 < \delta < 2\kappa^2 \theta^2$. 
Let $C(\cdot)$ be as in \autoref{lemma:Vbdd:undependent:u}. Since
$C(\cdot)$ is finite a.e.\, 
there exists $C_1 > 0$ and a set $K \subset X$ with
$\nu(K) > 1-\delta$ such that for $x \in K$, $C(x) < C_1$. Let
$K^\ast$ be as in \autoref{lemma:fake:ergodicity:Fij}, so in
particular $\nu(K^\ast) > 0$. Let $\Psi =
\{ x \in X \st g_{-t} x \in \Psi \text{ for some $t > 0$ } \}$.  
  
Suppose $x \in \Psi$ and $\bfw \in \bbV_{ij}(x)$. Let $t > 0$ be such
that $x' \equiv g_{-t} x \in K^\ast$. Then, $\bfw' \equiv g_{-t} \bfw \in
\bbV_{ij}(x')$. By \autoref{lemma:fake:ergodicity:Fij}, for any
$\ell > 0$, and at least $(1-\theta)$-fraction of $y \in
\cF_{ij}[x,\ell]$, we have $y \in K$, so $C(y) < C_1$. Thus, in view
of \autoref{lemma:Vbdd:undependent:u}, 
\begin{displaymath}
\|\bbR(x,y) \bfw'\| \le C_1 \|\bfw'\|.
\end{displaymath}
Then, by \autoref{prop:Rxy:v:small:implied:sync:bounded}, we get $\bfw'
\in \bbE_{[ij],bdd}(x')$. Then, $\bfw \in \bbE_{[ij],bdd}(x)$. 
\end{proof}

\begin{proof}[Proof of
  \autoref{prop:variant:Rxy:v:small:implied:sync:bounded}]
This follows immediately from \autoref{def:bounded:subspace},
\autoref{lemma:bounded:subspace} and \autoref{lemma:Vij:is:Eij}. 
\end{proof}





\section{The equivalence relation on \texorpdfstring{$\cW^u$}{W+}}
	\label{sec:the_equivalence_relation}

\subsection{The Equivalence Relation}
	\label{ssec:the_equivalence_relation}

Going forward, our space $\cC$ corresponds in our arguments to the space denoted $\operatorname{GSpc}$ in \cite{EskinMirzakhani_Invariant-and-stationary-measures-for-the-rm-SL2Bbb-R-action-on-moduli-space}.

\subsubsection{On local invariance of measures}
	\label{sssec:on_local_invariance_of_measures}
Recall that in \autoref{def:compatible_family_of_subgroups} we defined what it means to have a (measurable) compatible family of subgroups $U^+(x)\subset \bbG^{ssr}(x)$.
This notion formalizes the situation where a leafwise measure on an unstable has Haar conditionals in some directions, but the directions themselves are allowed to vary, even for points within the same unstable manifold, see \autoref{rmk:on_measures_with_local_invariance}.

With this in mind, starting from the compatible family $U^+(x)$ we will construct a measurable partition with atoms denoted \index{$\Psi_{ij}[x]$}$\Psi_{ij}[x]$ subordinated to the unstable (and even to $\gB_0[x]$).
The conditional measures \index{$f_{ij}[x]$}$f_{ij}[x]$ (supported on
$\Psi_{ij}[x]$) will lead to leafwise measures \index{$f$@$\wtd f_{ij}[x]$}$\wtd f_{ij}[x]$, and these measures will be invariant under $U^+(x)$.
Furthermore, the leafwise measures $\wtd f_{ij}[x]$ will be shown in \autoref{sec:extra_invariance}, if the assumptions of \autoref{thm:inductive:step} hold, to be invariant under a group $U^+_{new}(x)$ of dimension strictly larger than $U^+(x)$.

\subsubsection{The map $I(x,y)$.}
We will use the notation from \autoref{prop:compatibility_of_measurable_connection_and_the_family_of_subgroups}, with $\wp^+(x,y)$ for the strictly subresonant map induced by the measurable connection $P^+(x,y)$ on the cocycle $L\cW^u$ linearizing the unstable manifold, for $y\in \cW^u[x]$.
We then have that $U^+[y] = \wp^+(x,y) U^+[x]$ and $U^+(y) = \wp^+(x,y) U^+(x) \wp^+(x,y)^{-1}$, where $\wp^+(x,y) \in \bbG^{ssr}(\cW^u[x])$.
Also recall that $\cC(x) =\bbG^{ssr}(\cW^u[x])/U^+(x)$ parametrizes generalized subspaces in $\cW^u[x]$, with $g U^+(x) \in \cC(x)$ corresponding to the generalized subspace $g U^+[x] \subset \cW^u[x]$.

We now define \index{$I(x,y)$}$I(x,y): \cC(x) \to \cC(y)$ by
\begin{displaymath}
I(x,y) (g U^+(x)) = g \wp^+(y,x) U^+(y).  
\end{displaymath}
We claim that $g U^+(x) \in \cC(x)$ and $I(x,y) (g U^+(x)) \equiv g' U^+(y) \in \cC(y)$ parametrize the same generalize subspace in $\cW^u[x]$.
Indeed,
\begin{displaymath}
g' U^+[y] = g \wp^+(y,x) U^+[y] = g U^+[x], 
\end{displaymath}  
so the claim follows. 

\subsubsection{The subset $\cE_{ij}(x)$.}
	\label{sssec:the_subset_cE_ij_x}
Recall that we constructed an equivalence relation on the fine Lyapunov spectrum in \autoref{prop:subbundles_e_ij_bdd}, and denoted by $\wtilde{\Lambda}$ the corresponding set of equivalence classes.
Going forward, to ease the notation we will refer to the entire equivalence class by either $\index{$ij$}ij\in \wtilde{\Lambda}$ or $[ij]$.

For $ij \in \wtilde{\Lambda}$, let 
\begin{displaymath}
\index{$E$@$\cE_{ij}(x)$}\cE_{ij}(x) =
 \{ gU^+(x) \in \cC(x) \colon \bbj(gU^+(x)) \in
\bbE_{[ij],bdd}(x) \}.
\end{displaymath}
Here, $\bbj\colon \cC(x) \to \bbH(x)$ is the composition of the linearization map $\cC(x) \to L\cC(x)$ and the natural projection $L\cC(x) \to \bbH(x)$, see \autoref{eq:bold_H_definition}.

\subsubsection{A partition of $\cW^u[x]$.}
Recall that $\gB_0$ denotes the measurable partition
described in \autoref{sssec:measurable_partitions_and_subgroups_compatible_with_the_measure}.
The atom of $\gB_0$ containing $x$ is denoted by $\gB_0[x]$.
In this section, the only properties of $\gB_0$ we
will use is that it is subordinate to $\cW^u$, and that the atoms
$\gB_0[x]$ are relatively open in $\cW^u[x]$. 

\subsubsection{Equivalence relations.}
	\label{sssec:equivalence_relations_ij}
Fix $x_0 \in X$. 
For $x$, $x' \in \cW^u[x_0]$ we say that
\begin{displaymath}
x' \index{$\norm$@$\sim_{ij}$}\sim_{ij} x \text{ if $x' \in \gB_0[x]$ and $I(x',x)(U^+(x')) \in \cE_{ij}(x)$. }
\end{displaymath}
In the above equation, we think of $U^+(x')$ as the identity coset of $\cC(x') = \bbG^{ssr}(\cW^u[x])/U^+(x')$.

\begin{proposition}
\label{prop:sim:ij:is:equivalence:relation}
The relation $\sim_{ij}$ is a measurable equivalence relation.
In fact, $x\sim_{ij}x'$ if and only if $x'\in \gB_0[x]$ and $\cE_{ij}(x')=I(x,x')\cE_{ij}(x)$.
\end{proposition}
The main part of the proof of
\autoref{prop:sim:ij:is:equivalence:relation} is the following:
\begin{lemma}
	\label{lemma:Eij:bdd:locally:constant}
There exists a subset $\index{$\Psi$}\Psi \subset X$ 
with $\nu(\Psi) = 1$ such
that for any $ij \in \tilde{\Lambda}$, 
if $x_0 \in \Psi$, $x_1 \in \Psi$ and $x_1 \sim_{ij} x_0$ then 
$\cE_{ij}(x_1) = I(x_0,x_1)\cE_{ij}(x_0)$. 
\end{lemma}

\medskip\noindent
\textbf{Warning.}
In what follows
we will consider the condition $x' \sim_{ij} x$ to be undefined unless
$x$ and $x'$ both belong to the set $\Psi$ of
\autoref{lemma:Eij:bdd:locally:constant}. 

\subsubsection{Notation}
As in \S\ref{sec:subsubsec:cCbracketX}, 
let
$\cC[x] = \bbG^{ssr}(\cW^u[x])\cdot U^+[x]$; we view $\cC[x]$ as a collection of subsets of $\cW^u[x]$.
Then, for a.e.\; $y \in \cW^u[x]$, we have $U^+[y] = h U^+[x]$ for
some $h \in \bbG^{ssr}(\cW^u[x])$ by \autoref{prop:compatibility_of_measurable_connection_and_the_family_of_subgroups} and therefore $\cC[y] = \cC[x]$.

We think of a point $\cQ \in \cC[x_0]$ as a subset of $\cW^u[x_0]$, and call it a generalized subspace.
For $x \in \cW^u[x_0]$ let \index{$U$@$\cU_x$}$\cU_x: \cC(x) \to \cC[x]$ denote the natural map. More precisely, for $g U^+(x) \in \cC(x)$ we have
\begin{displaymath}
\cU_x( g U^+(x)) = g U^+[x] \subset \cC[x] = \cC[x_0].
\end{displaymath}  
We then have the inverse map $\cU_x^{-1}: \cC[x_0] \to \cC(x)$. 
We also use the notation \index{$E$@$\cE_{ij}[x]$}$\cE_{ij}[x] = \cU_x(\cE_{ij}(x))$. Clearly, $\cE_{ij}[x] \subset \cC[x]$. Note also that $I(x,y) = \cU_y^{-1} \circ \cU_x$. 

We will need the following lemma:
\begin{lemma}
\label{lemma:simple:matching}
For every $\delta > 0$ there exists  constants $\kappa(\delta)>1$ and
$\epsilon(\delta) > 0$ and a
compact set $K$ with $\nu(K) >
1-\delta$ such that the following holds:
Suppose $x_0 \in K,x_1 \in K, \ell > 0$, $\theta>0$ are such
that $x_1 \in \gB_0[x_0]$ and the following hold:
\begin{itemize}
\item For at least $(1-\theta)$ fraction of $y \in
  \cF_{ij}[x_0,\ell]$, $y \in K$ and
\begin{equation}
\label{eq:dy0:cH:loc}
  d^y_{\cH,loc}(U^+[y],U^+[z]) < \epsilon(\delta) \qquad
  \text{for some $z \in \cF_{ij}[x_1,\ell]$}. 
\end{equation}
\item $d^u(x_0,x_1) < \epsilon(\delta)$. 
\item For at least $(1-\theta)$ fraction of $y \in
  \cF_{ij}[x_1,\ell]$, $y \in K$.
\item For $k=0,1$,
\begin{equation}
\label{eq:most:fraction:close}
|\cB_0[g^{ij}_{-\ell} x_k]  \cap K| \ge 0.99
  |\cB_0[g^{ij}_{-\ell} x_k]|. 
\end{equation}
\end{itemize}
Then, for at least $(1-6\theta)$-fraction of $y \in
\cF_{ij}[x_0,\ell]$ we have $y \in K$, (\ref{eq:dy0:cH:loc}) holds,
and $\cB_0[z] \cap K$ is non-empty.  
\end{lemma}

\begin{proof}
For $\epsilon_0' > 0$, let $\gB_0^{(\epsilon_0')}[x]$ and
$\cB_0^{(\epsilon_0')}[x]$ be defined as in \autoref{sssec:partitions_and_u_plus}.
Suppose $\delta > 0$. Then there exists
$\epsilon_0'>0$ constants $C(\delta)$ and a compact set $K$ with
$\nu(K) > 1-\delta$ so that
\begin{itemize}
\item[{\rm (K1)}] For $x \in K$, $d^u(x,\partial \gB_0[x]) > 2 \epsilon_0'$.
\item[{\rm (K2)}] For $x \in K$, $C(x) < C(\delta)$, where $C(x)$ is as in
  \autoref{lemma:L:and:C}. 
\item[{\rm (K3)}] The map $x \to \cB_0[x]$ is uniformly continuous on $K$ in the
  Hausdorff topology. 
\item[{\rm (K4)}] There exits $\epsilon_0 = \epsilon_0(\delta,\epsilon_0') >0$
so that if $x\in K$ then for any $h\in \bbG^{ssr}(\cW^u[x])$ with
$\|L_x h - \id \|_{L\cW^u(x)} <C(\delta) \epsilon_0$, we have that
$$h\gB_0^{(2\epsilon_0')}[x] \subset \gB_0^{(\epsilon_0')}[x].$$
\end{itemize}

Let $x_k' = g^{ij}_{-\ell} x_k$, so that
$g^{ij}_{-\ell}(\cF_{ij}[x_k,\ell]) = \cB_0[x_k']$. 
Note that $\gB_0[x_0'] = \gB_0[x_1']$, and also $x_k' \in K$.

Let $\gB_\ell^{ij}[x] = g^{ij}_{-\ell}(\gB_0[g^{ij}_\ell x])$. Then,
the sets of the form $\gB_\ell^{ij}[y']$, where $y' \in \gB_0[x']$
partition $\gB_0[x_0'] = \gB_0[x_1']$. Thus, we can choose $y_n' \in
U^+[x_0']$ such that if we define
\begin{displaymath}
P^0_n = \gB^{ij}_\ell[y_n'] \cap U^+[x_0']
\end{displaymath}
\begin{displaymath}
P^1_n = \gB^{ij}_\ell[y_n'] \cap U^+[x_1'],
\end{displaymath}
then the sets $P_n^0$ partition $\cB_0[x_0']$ and the sets $P_n^1$
are contained in $\cB_0[x_1']$. Let
\begin{displaymath}
\cF^0_n = g^{ij}_\ell(P^0_n), \qquad \cF^1_n = g^{ij}_\ell(P^1_n).
\end{displaymath}
Then, if $y_n = g^{ij}_\ell y_n'$, we have
\begin{equation}
\label{eq:cFn:subset:gB}
\cF^0_n \subset
\gB_0[y_n], \qquad \cF^1_n \subset \gB_0[y_n],
\end{equation}
and
\begin{displaymath}
\cF_{ij}[x_0,\ell] = \bigcup_n \cF^0_n, \qquad 
\bigcup_n \cF^1_n \subset \cF_{ij}[x_1,\ell].
\end{displaymath}
Let
\begin{displaymath}
S_0 = \{ y \in \cF_{ij}[x_0',\ell] \st y \in K \text{ and
  (\ref{eq:dy0:cH:loc}) holds} \},
\end{displaymath}
\begin{displaymath}
S_1 = \{ z \in \cF_{ij}[x_1',\ell] \st z \in K \}
\end{displaymath}
and for $k=0,1$, let $S'_k = g^{ij}_{-\ell} S_k$. By assumption, for
$k=0,1$, $|S'_k| \ge (1-\theta) |\cB_0[x_k']|$.

Suppose $y'_n \in P^0_n \cap S'_0$. Then,
there exists $z_n' \in P_n^1$ so that $y_n = g^{ij}_\ell y_n'$ and
$z_n = g^{ij}_\ell z_n'$ satisfy $y_n \in K$ and
\begin{equation}
\label{eq:tmp:dy0:cH:loc}
  d^{y_n}_{\cH,loc}(U^+[y_n],U^+[z_n]) < \epsilon(\delta).
\end{equation}
Then, by \autoref{lemma:L:and:C}
there exists $h_n \in \bbG^{ssr}(\cW^u[y_n])$
with $\|L_{y_n}h_n - \id \|_{L\cW^u(y_n)} \le C(\delta) \epsilon(\delta)$ 
and such that
$U^+[z_n]   = h_n U^+[y_n]$.

We choose $\epsilon(\delta) < \epsilon_0$, where $\epsilon_0$ is as in
(K4). Then,  $h_n \gB_0^{(2\epsilon_0')}[y_n] \subset
\gB_0^{(\epsilon_0')}[y_n]$. Then, in view of (\ref{eq:cFn:subset:gB})
and (K1), 
\begin{displaymath}
h_n(\cF^0_n \cap K) \subset \cF^1_n. 
\end{displaymath}
Let $\tilde{h}_n = g^{ij}_{-\ell} \circ h_n \circ g^{ij}_\ell$. Then,
$\tilde{h}_n$ is even closer to the identity and 
\begin{equation}
\label{eq:tilde:hn:close}
\tilde{h}_n (P^0_n \cap S'_0) \subset P^1_n
\end{equation}
Note that $x_0' \in K$, $x_1' \in K$, $d^u(x_0',x_1') < \epsilon(\delta)$. 
In view of the condition (K3) and (\ref{eq:most:fraction:close}),
we can choose $\epsilon(\delta)$ small enough so that 
\begin{equation}
\label{eq:tmp:cB0:minus:close}
d^Q_{\cH}(\cB_0[x_0'], \cB_0[x_1']) < \epsilon_0',  
\end{equation}
where $d^Q_{\cH}$ is the Hausdorff distance of the projection to $Q$. 
We normalize the measures $| \cdot |$ on $U^+[x_0']$ and $U^+[x_1']$
so that $|\cB_0[x_0']| = |\cB_0[x_1']| = 1$. Then, it follows from
(\ref{eq:tilde:hn:close}) and (\ref{eq:tmp:cB0:minus:close})
that 
\begin{displaymath}
|P^0_n \cap S'| \le 3 | P^1_n|.
\end{displaymath}
We now apply \autoref{lemma:segments:equal:length} (with $P_n^0$ in
place of $P_n$, $P_n^1$ in place of $P_n'$, $S'_0$ in place of $Q$ and
$S'_1$ in place of $Q'$) to get a set
$\hat{S}' \subset \cB_0[x_0']$ with $|\hat{S}'| \ge
(1-6\theta)$ such that if $y_n' \in \hat{S}' \cap
\cP_n^0$ then there exists $z_n' \in S'_1 \cap \cP_n^1$. Let
$\hat{S} =
g^{ij}_\ell(\hat{S}')$. Then, $|\hat{S}| \ge (1-\kappa(\delta)\theta)
|\cF_{ij}[x_0,\ell]|$, and if $y_n = g^{ij}_\ell y_n' \in \hat{S}$
then $z_n = g^{ij}_\ell z_n'$ belongs to 
$\cF_{ij}[x_1,\ell] \cap \gB_0[y_n] \cap K$. 
\end{proof}

\medskip\noindent
\textbf{Proof of \autoref{lemma:Eij:bdd:locally:constant}.}
This proof follows \cite[\S11.1*]{EskinMirzakhani_Invariant-and-stationary-measures-for-the-rm-SL2Bbb-R-action-on-moduli-space}.

Let $\delta > 0$ be a small constant to be chosen later.
Let $K \subset X$ and $C > 0$ be such
that \autoref{lemma:simple:matching} holds with $\epsilon(\delta) <
C^{-1}$,  
for $x \in K$ \autoref{prop:hausdorff_distance_and_norm_of_vector} holds with $c_1(x) >
C^{-1}$, and for all $x \in K$, all $\bfv \in \bbE_{[ij],bdd}(x)$ and all
$\ell > 0$, for all $y \in \cF_{ij}[x,\ell] \cap K$, we have
\begin{equation}
\label{eq:tmp:redef:bounded}
	\|\bbR(x,y) \bfv\| < C \|\bfv\|. 
\end{equation}

By \autoref{lemma:fake:ergodicity:Fij} there exists a subset $K^*
\subset K$ with $\nu(K^*) \ge (1-4\kappa^2 \delta^{1/2})$ such that for
$x \in K^*$, $|\cB_0[x] \cap K| \ge 0.99 |\cB_0[x]|$, 
\autoref{eq:fake:ergodicity:Fij} holds with $\theta' =
\delta^{1/2}$, and $x\in K^*\implies \cB_0[x]\subset K^*$.
Furthermore, we may ensure that for $x \in K^*$, $K^* \cap
\cF_{ij}[x]$ is relatively open in $\cF_{ij}[x]$.
(Indeed, suppose $z \in \cF_{ij}[x] \cap \gB_0[x]$. Then,
for all $\ell$, $\cF_{ij}[x,\ell] = \cF_{ij}[z,\ell]$, 
and thus \autoref{eq:fake:ergodicity:Fij} holds for $z$.)

Let
\begin{displaymath}
\Psi  = \left\lbrace x \in X \colon \lim_{T \to \infty} \tfrac 1 T|\{ t \in [0,T] \colon
  g_{-t} x \in K^* \}| \ge (1-8 \kappa^2 \delta^{1/2}) \right\rbrace
\end{displaymath}
Then $\nu(\Psi) = 1$. From its definition, $\Psi$ is invariant
under $g_t$. Since $K^* \cap \cF_{ij}[x]$ is
relatively open in $\cF_{ij}[x]$, $\Psi$ is saturated by the leaves
of $\cF_{ij}$. This implies that $\Psi$ is (locally) invariant under
$U^+$. We also assume that for $x \in \Psi$, $\dim \cE_{ij}(x)$ is
constant. This condition is also locally invariant under $U^+$. 




Suppose $x_0 \in \Psi$, $x_1 \in \gB_0[x_0] \cap \Psi$.
We will use notation for Hausdorff distances from
\autoref{sssec:local_hausdorff_distance_on_unstable}.
By
\autoref{prop:compatibility_of_measurable_connection_and_the_family_of_subgroups},
$U^+[x_1] \in \bbG^{ssr}(\cW^u[x_0]) U^+[x_0]$. Therefore,
$U^+[x_1] \in \cC[x_0]$, and thus, $\cU_{x_0}^{-1}(U^+[x_1])$
is defined. Let
\begin{displaymath}
\bbu = \bbj(\cU_{x_1}^{-1}(U^+[x_1])) \in \bH(x_0). 
\end{displaymath}
Suppose $\cQ \in \cE_{ij}[x_1]$. Let
\begin{displaymath}
\bfv = \bbj(\cU_{x_1}^{-1}(\cQ)).
\end{displaymath}
Assume
\begin{equation}
\label{eq:bfk:bbuk}
\|\bfv\| \le 2 \| \bbu \|.
\end{equation}
Since, for $x_k$ fixed, the map $\bbj \circ \cU_{x_k}^{-1}: \cC[x_k]
\to \bH(x_k)$ is analytic, (\ref{eq:bfk:bbuk}) will hold assuming
$\cQ$ is sufficiently close to $U^+[x_1]$.

Note that
$\cQ \in \bbG^{ssr}(\cW^u[x_1]) U^+[x_1]$. Since 
$\bbG^{ssr}(\cW^u[x_0]) = \bbG^{ssr}(\cW^u[x_1])$ is a group, we get
$$\cQ \in \bbG^{ssr}(\cW^u[x_0]) U^+[x_0].$$ 
Therefore $\cU_{x_0}^{-1}(\cQ)$ is defined. 
Let
\begin{displaymath}
\bfw = \bbj(\cU_{x_0}^{-1}(\cQ)).
\end{displaymath}



For any $C_1>0$ sufficiently large (depending on $x_0$, $x_1$), we can
find $t \in \R$ with 
$C_1 < t < 2C_1$ such that
$x_0' \equiv g_{-t}^{ij} x_0 \in K^*$, 
$x_1' \equiv g^{ij}_{-t} x_1 \in K^*$.  By
\autoref{prop:backwards_nesting_of_partitions},
$x_1' \in \gB_0[x_0']$. 
Let $\bbu' = g_{-t}^{ij} \bbu$,
$\bfv' = g_{-t}^{ij} \bfv$,
$\bfw' = g_{-t}^{ij} \bfw$,
$\cQ' = g_{-t}^{ij} \cQ$. 
By choosing $C_1$ sufficiently large (depending on $x_0,x_1$), we can
ensure that
\begin{equation}
\label{eq:estimates:jplus:xk:prime}
  \|\bbu' \| \le C^{-3}, \qquad \|\bbv' \| \le C^{-3},
\end{equation}
where $C$ is defined in \autoref{eq:tmp:redef:bounded}.


Suppose $\ell' > 0$. Then, for $k=0,1$, since $x_k' \in K^*$, 
\begin{displaymath}
|\{ y_k' \in \cF_{ij}[x_k',\ell'] \colon y_k' \in K\}| \ge (1-\delta^{1/2}) 
|\cF_{ij}[x',\ell']|,
\end{displaymath}
Since $U^+[x_1] \in \cE_{ij}[x_0]$, we have $U^+[x_1'] \in
\cE_{ij}[x_0']$, and thus
\begin{displaymath}
\bbu' \in \bbE_{[ij],bdd}(x_0').
\end{displaymath}
Since $x_0' \in K^* \subset K$, we have by
(\autoref{eq:tmp:redef:bounded}), 
for all $y_0' \in
\cF_{ij}[x_0', \ell'] \cap K$,  
\begin{equation}
\label{eq:j:Uplus:y0:close:j:Uplus:y1}
\|\bbR(x_0',y_0') \bbu' \| \le C
\|\bbu' \| \le C^{-2}, 
\end{equation}
where we have used (\autoref{eq:estimates:jplus:xk:prime}) 
for the last estimate. Note that for a suitable $y_1'$, 
\begin{displaymath}
\bbR(x_0',y_0') \bbu_0' = \bbj(\cU_{y_0'}^{-1}(U^+[y_1'])).
\end{displaymath}
Therefore, by \autoref{prop:hausdorff_distance_and_norm_of_vector},
for all $y_0' \in \cF_{ij}[x_0',\ell'] \cap K$, for a suitable $y_1'
\in \cF_{ij}[x_1',\ell']$,
\begin{equation}
\label{eq:hd:U:y0:prime:U:y1:prime:bounded}
d^{y_0'}_{\cH,loc}(U^+[y_0'],U^+[y_1']) \le C^{-1}. 
\end{equation}

Let
\begin{displaymath}
\cD' = \{ \ell' \in \reals^+ \st g^{ij}_{-\ell'} x_0' \in K^\ast \text{
  and } g^{ij}_{-\ell'} x_1' \in K^\ast \}. 
\end{displaymath}
By the definition of $\Psi$,  the density of $\cD'$ is at least $1-4
\kappa^2 \delta^{1/2}$. Now suppose $\ell' \in \cD$. Then,
in view of \autoref{lemma:simple:matching}, for at least
$(1-6\delta^{1/2})$-fraction of $y_0' \in \cF_{ij}[x_0',
\ell']$,
\begin{multline}
\label{eq:conditions:on:y0:prime}
\text{ $y_0' \in K$  and there exists }\\ \text{$y_1' \in \cF_{ij}[x_1',\ell']
  \cap K \cap \gB_0[y_1']$ such that
(\ref{eq:hd:U:y0:prime:U:y1:prime:bounded}) holds. }
\end{multline}
Now suppose $y_0'$, $y_1'$ are such that
(\ref{eq:conditions:on:y0:prime}) holds. 
Since $\cQ \in \cE_{ij}[x_1]$, $\bfv \in
\bbE_{[ij],bdd}(x_1)$, and thus, $\bfv' \in
\bbE_{[ij],bdd}(x_1')$. 
Hence, by (\autoref{eq:tmp:redef:bounded}), since $y_1' \in K$, 
\begin{equation}
\label{eq:j:Uplus:y0:close:image:Qk:prime}
\|\bbR(x_1',y_1') \bfv'\| \le C
\| \bfv'\| \le C^{-2}. 
\end{equation}
where we used (\autoref{eq:estimates:jplus:xk:prime}) for the last
estimate. Note that
\begin{displaymath}
\bbR(x_1',y_1') \bfv' = \bbj(\cU_{y_1'}^{-1}(\bbR(x_1',y_1')
\cQ')). 
\end{displaymath}
Therefore, by
\autoref{prop:hausdorff_distance_and_norm_of_vector},
\begin{displaymath}
d^{y_1'}_{\cH,loc}(U^+[y_1'],\bbR(x_1',y_1')\cQ') \le C^{-1}. 
\end{displaymath}
Since $y_1' \in \gB_0[y_0']$, this implies
\begin{displaymath}
d^{y_0'}_{\cH,loc}(U^+[y_1'],\bbR(x_1',y_1')\cQ') \le C^{-1}. 
\end{displaymath}
Then, by the triangle inequality, 
\begin{equation}
\label{eq:final:hd:bound}
d^{y_0'}_{\cH,loc}(U^+[y_0'],\bbR(x_1',y_1')
\cQ') \le 2 C^{-1}.
\end{equation}
Note that
\begin{displaymath}
\bbR(x_0',y_0') \bfw' = \bbj(\cU_{y_0'}^{-1}(\bbR(x_1',y_1')
\cQ')).
\end{displaymath}
Then,
by \autoref{prop:hausdorff_distance_and_norm_of_vector}
\begin{equation}
\label{eq:tmp:bbR:two}
\|\bbR(x_0',y_0')\bfw'\| \le 2. 
\end{equation}

Thus, for $\ell' \in \cD'$ and at least $(1-6\delta^{1/2})$-fraction of $y_0'
\in \cF_{ij}[x_0',\ell']$,   (\ref{eq:tmp:bbR:two}) holds. Let
$\epsilon, \theta$ be as in
\autoref{prop:variant:Rxy:v:small:implied:sync:bounded}. By
choosing $\delta > 0$ small enough, we can make sure that $1-6\delta^{1/2} <
\theta$ and the density of $\cD'$ is at least $1-\epsilon$. Then,
by
\autoref{prop:variant:Rxy:v:small:implied:sync:bounded} we get
$\bfw' \in \bbE_{[ij],bdd}(x_0')$. Then, $\bfw \in
\bbE_{[ij],bdd}(x_0)$.


Thus, for all $\cQ \in \cE_{ij}[x_1]$ sufficiently close to
$U^+[x_1]$ so that (\ref{eq:bfk:bbuk}) holds,
we have
$\bbj(\cU_{x_0}^{-1}(\cQ)) \in \bbE_{[ij],bdd}(x_0)$, which
implies $\cQ \in \cE_{ij}[x_0]$. Note that for fixed $x_0$,
$x_1$, $\cE_{ij}[x_0]$ and $\cE_{ij}[x_1]$ are both analytic subsets
of $\cC[x_0] = \cC[x_1]$. Thus, we get, 
$\cE_{ij}[x_1] \subset \cE_{ij}[x_0]$. Since $x_0, x_1 \in \Psi$ we
have $\dim  \cE_{ij}[x_1] = \dim \cE_{ij}[x_0]$. 
This implies that
$\cE_{ij}[x_0] = \cE_{ij}[x_1]$. This is equivalent to $I(x_0,x_1)
\cE_{ij}(x_0) = \cE_{ij}(x_1)$. 
\qed\medskip

\subsubsection{Proof of the equivalence relation property in \autoref{prop:sim:ij:is:equivalence:relation}}
	\label{sssec:proof_of_the_equivalence_relation_property_in_prop:sim:ij:is:equivalence:relation}
We assume \autoref{lemma:Eij:bdd:locally:constant}.
Note that the map $\bbj:\cC(x) \to \bbH(x)$ takes the identity coset of $\cC(x)$ to $0$. We have $0 \in \bbE_{[ij],bdd}(x)$,
therefore,  
\begin{equation}
\label{eq:Uplus:x:in:cEij:x}
U^+(x) \in \cE_{ij}(x). 
\end{equation}
Thus $x \sim_{ij} x$. 

Suppose $x' \sim_{ij} x$. Then, $x' \in \gB_0[x]$, and so $x \in \gB_0[x']$. 
By (\autoref{eq:Uplus:x:in:cEij:x}), $U^+(x) \in \cE_{ij}(x)$, and by
Lemma~\autoref{lemma:Eij:bdd:locally:constant}, $\cE_{ij}(x') =
I(x,x')\cE_{ij}(x)$. Therefore, $I(x,x') (U^+(x)) \in \cE_{ij}(x')$, and thus $x
\sim_{ij} x'$. 

Now suppose $x' \sim_{ij} x$ and $x'' \sim_{ij} x'$. Then, $x'' \in
\gB_0[x]$. Also, $I(x'',x')(U^+(x'')) \in \cE_{ij}(x') = I(x,x')\cE_{ij}(x)$.
Since $I(x,x')^{-1} I(x'',x') = I(x'',x)$, we have $I(x'',x)(U^+(x'')) \in \cE_{ij}(x)$. Therefore $x'' \sim_{ij} x$. 
\qed\medskip



\subsubsection{The measurability property}
	\label{sssec:the_measurability_property}
For a.e. $x\in X$, the set $\cE_{ij}(x)\subset \cC(x)$ is an algebraic subset, and furthermore its $\bbG^{ssr}(x)$-orbits form an algebraic space \index{$Y$@$\cY(x)$}$\cY(x)$.
We can therefore form a measurable bundle $\cY\to X$, defined $\nu$-a.e., parametrizing these sets, with the property that $I(y,x)\cE_{ij}(y)\in \cY(x)$.

Let $X_{\gB_0}$ denote the space of atoms of the measurable partition $\gB_0$, which is by construction a standard measurable space and there is a measurable quotient map $\psi \colon X \to X_{\gB_0}$.
By a standard measurable selection theorem
there exists a measurable section $\sigma_0\colon X_{\gB_0}\to X$.
Let $\cY_0:=\sigma_0^*\cY$ be the pulled back bundle, which is a measurable bundle over $X_{\gB_0}$.
The collection of maps $I(x,y)$, measurable in $(x,y)\in X\times X$ restricted to the subset with $x\in \gB_0[y]$, defines a measurable map
\[
	\psi_{ij}\colon X\to \cY_0 \text{ via }\psi_{ij}(y):=I(y,\sigma_0(\psi(y)))\cE_{ij}(y) \in \cY(\sigma_0(\psi(y)))\isom \cY_0(\psi(y)).
\]
The atoms of the equivalence relation $\sim_{ij}$ are, by construction, the fibers of the measurable map $\psi_{ij}$.
Indeed by \autoref{lemma:Eij:bdd:locally:constant}, we have
\begin{displaymath}
x \sim_{ij} y \qquad \text{if and only if $y \in \gB_0[x]$ and }
\cE_{ij}(y) = I(x,y)\cE_{ij}(x).
\end{displaymath}
Together with the cocycle property of $I(x,y)$, this establishes the equality of fibers of $\psi_{ij}$ and equivalence classes of $\sim_{ij}$, and thus the measurability of $\sim_{ij}$.

\subsubsection{The equivalence classes $\Psi_{ij}[x]$}
	\label{sssec:the_equivalence_classes_psi__ij}
For a full measure subset of $X$, denote the atoms of the equivalence class by
\begin{displaymath}
\index{$\Psi_{ij}[x]$}\Psi_{ij}[x] = \{ x' \in \gB_0[x] \colon x' \sim_{ij} x \}.
\end{displaymath}
For correspondence of notation, note that the above sets in
\cite[Lemma~11.3]{EskinMirzakhani_Invariant-and-stationary-measures-for-the-rm-SL2Bbb-R-action-on-moduli-space}
are denoted by $\cC_{ij}[x]$.

Let \index{$\Psi_{ij}$}$\Psi_{ij}$ denote the $\sigma$-algebra
of $\nu$-measurable sets which are unions of the equivalence classes
$\Psi_{ij}[x]$.
By \autoref{sssec:the_measurability_property}, the atoms of $\Psi_{ij}$ are the fibers of the measurable maps $\psi_{ij}$, hence $\Psi_{ij}$ is a countably generated $\sigma$-algebra.

\begin{lemma}
\label{lemma:Psiij:equivariant}
There exists an invariant set of full measure $X_0\subset X$ such that, if $x\in X_0$, $t \in \bR$, and $u \in U^+(x)$ then we have:
\begin{itemize}
\item[{\rm (a)}] 
$\left(g_t \Psi_{ij}[x]\right) \cap \gB_0[g_t x] \cap g_t \gB_0[x] = \Psi_{ij}[g_t x] \cap
\gB_0[g_t x] \cap g_t \gB_0[x]$.
\item[{\rm (b)}]
$\left(u\Psi_{ij}[x]\right) \cap \gB_0[ux] \cap u\gB_0[x] = \Psi_{ij}[ux] \cap \gB_0[ux] \cap
u \gB_0[x]$. 
\end{itemize}
\end{lemma}

\begin{proof}
Note that the sets $U^+[x]$ and $\bbE_{[ij],bdd}(x)$ are
$g_t$-equivariant. Therefore, so are the $\cE_{ij}[x]$, which implies
(a). Part (b) is also clear since, in view of 
\autoref{prop:agreement:of:u:star:and:p} and the Ledrappier invariance principle \autoref{thm:ledrappier_invariance_principle},
 $u_*(x,ux) \bbE_{[ij],bdd}(x) =
\bbE_{[ij],bdd}(ux)$. 
\end{proof}


\begin{definition}[{The measures $f_{ij}[x]$}]
	\label{def:measures_fij}
We now define \index{$f_{ij}[x]$}$f_{ij}[x]$ 
to be the conditional measure
of $\nu$ along the $\Psi_{ij}[x]$. In other words, $f_{ij}[x]$
is defined so that for any measurable $\phi: X \to \bR$,
\begin{displaymath}
\mathbb{E}(\phi \mid \Psi_{ij})(x) = \int_X \phi \, df_{ij}[x]. 
\end{displaymath}
We view $f_{ij}[x]$ as a measure on $\cW^u[x]$ 
which is supported on $\Psi_{ij}[x]$. 
\end{definition}

\begin{definition}[{The measures $\wtd f_{ij}[x]$}]
	\label{def:measures_fij_tilde}
By \autoref{lemma:Psiij:equivariant},  for almost every $x$ and all $t\ge 0$ we have 
$$ \Psi_{ij}[ x]  \subset g_t \Psi_{ij}[g_{-t} x].$$
We define a locally finite Radon measure 
\index{$f$@$\wtd f_{ij}[x]$}$\wtd f_{ij}[x]$ 
on $\cW^u[x]$  as follows: we declare $\wtd f_{ij}[x](\Psi_{ij}[x]) = 1 $ and for any compact $K\subset \cW^u[x]$, we take $t>0$ so that $g_{-t}K\subset \gB_0[g_{-t}x]$ and declare 
$$\wtd f_{ij}[x](K) := \frac{
(g_t)_*  f_{ij}[g_{-t}x] (K)
}{(g_t)_*  f_{ij}[g_{-t}x] (\Psi_{ij}[x])}.$$
The $g_t$-invariance of the measure $\nu$ ensures that $\wtd f_{ij}[x](K)$ is well-defined independent of the choice of 
$t>0$ making  $g_{-t}K\subset \gB_0[x].$
\end{definition}





\subsubsection{The measures $\overline{f}_{ij}[x]$}
	\label{sssec:the_measures_overline_f_ij}
Let \index{$\pi_{U^+}$}$\pi_{U^+}\colon \cW^u[x] \to \cC[x]$ denote the map $x \mapsto
U^+[x]$.
We would like to define a measure \index{$f$@$\overline f_{ij}[x]$}$\overline f_{ij}[x]$ on
$\cC[x]$ by $\overline f_{ij}[x] = (\pi_{U^+})_\ast \wtd f_{ij}[x]$ but
this does not make sense since for $\cV \subset \cC[x]$
the $\wtd f_{ij}[x]$-measure of $\pi_{U^+}^{-1}(\cV)$ is typically
infinite.  Instead we proceed as follows. Fix a compact $K_U \subset
U^+[x]$ of positive Haar measure.
Suppose $\cV \subset \cC[x]$
is a compact set. For a.e.\ $c \in \cV$ choose $g_c \in \bbG^{ssr}(x)$ such
that $g_c U^+[x] = c$ and the map $c \mapsto g_c$ is measurable. 
(If $g_c' \in \bbG^{ssr}(x)$ is another such
choice, then $g_c' = g_c u_c$ for some $u_c \in U^+(x)$). 
Then, let $\cV' \subset \cW^u[x]$ be $\bigcup_{c
  \in \cV} g_c K_U$. Since the conditionals of $\wtd f_{ij}[x]$ along the
sets $U^+[g_c x]$ are Haar under the action of $g_c U^+(x) g_c^{-1}$, 
$\wtd f_{ij}[x](\cV')$ does not depend on the choice of $g_c$. Thus, we may
define the measure $\overline f_{ij}[x]$ on $\cC[x]$
by $\overline f_{ij}[x](\cV) = \wtilde{f}_{ij}(\cV')$.
It is easy to see that
up to normalization, $\overline f_{ij}[x]$ does not depend on the
choice of $K_U$. 

We normalize $\overline f_{ij}[x]$ so that $\overline
f_{ij}[x](\gB_\cC[x]) = 1$, where\index{$B$@$\gB_\cC[x]$}
\begin{align}
	\label{eqn:frakB_calC}
	\gB_\cC[x] := \{ U^+[y] \st y \in \gB_0[x]\}
\end{align}



\subsubsection{Recovering $\wtilde{f}_{ij}[x]$ from $\bar{f}_{ij}[x]$}
	\label{sssec:recovering_wtilde_f_ij_x_from_bar_f_ij}
Let us now explain how to compute $\wtilde{f}_{ij}[x](E)$ for a (bounded) measurable set $E\subset \cW^u[x]$, given $\bar{f}_{ij}[x]$.
Recall that we have a measurable function $g(c)\in \bbG^{ssr}(x)$ over $c\in \cC[x]$ such that $c=g(c)U^+[x]$ and let $\nu_{U^+[x]}$ denote a fixed Haar measure on $U^+[x]$, so that $g(c)_*\nu_{U^+[x]}$ is a Haar measure on $c$. 
Then
\[
	\wtilde{f}_{ij}[x](E):=
	\fint_{\cC[x]} 
	g(c)_{*}\nu_{U^+[x]}
	\left(E\right) d\bar{f}_{ij}(c)
\]
where $\fint$ denotes the integral normalized by the scalar obtained by using $\gB_0[x]$ instead of $E$ in the domain of integration.

\subsubsection{The measures $\bbf_{ij}(x)$}
	\label{sssec:the_measures_bbf_ij}
Recall the maps $\cU_x\colon \cC(x)\to \cC[x]$
 and $\bbj\colon \cC(x)\to \bbH(x)$.
 We let \index{$f$@$\bbf_{ij}(x)$}$\bbf_{ij}(x) = (\bbj \circ \cU_x^{-1})_* \overline f_{ij}[x].$
Again, $\bbf_{ij}(x)$ is a locally finite Radon measure on  $\bbH(x).$
We note that the support of each $\bbf_{ij}(x)$ is contained in $\bbE_{[ij],bdd}(x).$
We record the following elementary properties that follow from the constructions:
\begin{lemma}
\label{lemma:transformation:rule:bff}
\leavevmode
For $x,y$ in a set of full measure:
\begin{enumerate}
	\item 
	 If $y \in \gB_0[x]$ then we have $$I(x,y)_*\overline f_{ij}[x]\propto \overline f_{ij}[y].$$
	\item 
	If $y \in \cF_{ij}[x]$ then we have
	\begin{displaymath}
		\bbf_{ij}(y) \propto \bbR(x,y)_* \bbf_{ij}(x). 
	\end{displaymath}
	Moreover, if $y \in G[x]$ then we have
	\begin{displaymath}
		\bbf_{ij}(y) \propto \bbR(x,y)_* \bbf_{ij}(x). 
	\end{displaymath}
	Here $\propto$ means that the two measures agree up to a scalar factor.
\end{enumerate}
\end{lemma}

\begin{definition}[{The ``distance'' $d_*(\cdot, \cdot)$}]
\label{def:d:star}
We consider the set of Radon measures on a manifold, with the vague
topology (so $\mu_n \to \mu$ if and only if for all continuously
compactly supported $f$, $\mu_n(f) \to \mu(f)$). This topology is
metrizable, and let \index{$d$@$\hat{d}^x_*(\cdot, \cdot)$}$\hat{d}^x_*(\cdot, \cdot)$ be a metric inducing this
topology on $\cW^u[x]$.
Then, for measures $\nu_1,\nu_2$ on $\cW^u[x]$, let
\begin{displaymath}
\index{$d_*(\nu_1,\nu_2)$}d_*(\nu_1,\nu_2) = \hat{d}^x_*\left(
\frac{\nu_1}{\nu_1(\gB_0[x])}, \frac{\nu_2}{\nu_2(\gB_0[x])}\right).
\end{displaymath}
For measures $\nu_1, \nu_2$ on $\bH(x)$, define
\begin{displaymath}
\index{$d_*(\nu_1,\nu_2)$}d_*(\nu_1,\nu_2) = \hat{d}^x_*\left(
\frac{\nu_1}{\nu_1(D(x))}, \frac{\nu_2}{\nu_2(D(x))}
\right),
\end{displaymath}
where $D(x) = \bfj(\cU_x^{-1}(\frakB_\cC[x]))$ and $\frakB_\cC[x]$ is defined by \autoref{eqn:frakB_calC}.

\end{definition}

\subsection{Proposition 11.4 from EM}
	\label{ssec:proposition_11_4_from_em}

\subsubsection{Top-balanced $Y$-configurations}

\begin{definition}[top-balanced $Y$-configuration]
	\label{def:top_balanced_y_configuration}
Fix $C>0$.  Using the notation \S\ref{sec:subsec:Yconf:notation}, we
say a $Y$-configuration $Y$ is \emph{$(C, ij)$-top-balanced} if $q_3= g^{ij}_{t_3}(q_1)$ and $q_2= g^{ij}_{t_2}(uq_1)$ where $$|t_2-t_3|\le C.$$
\end{definition}

We denote by \index{$Y_{ij,tb}(q_1,u,\ell,\tau)$}$Y_{ij,tb}(q_1,u,\ell,\tau)$ the unique
$(0,ij)$-top-balanced $Y$ conguration with $q_1(Y) = q_1$, $u(Y) = u$,
$\ell(Y) = \ell$ and $\tau(Y) = \tau$.

The main goal of this section is to prove the next result, analogous to \cite[Prop.~11.4]{EskinMirzakhani_Invariant-and-stationary-measures-for-the-rm-SL2Bbb-R-action-on-moduli-space}:

\begin{proposition}
\label{prop:nearby:linear:maps}
For any $\delta>0$, there exists a compact set $K_0\subset X$ with $\nu(K_0)>1-\delta$, with the following properties.

Suppose $Y$, $Y'$ are bottom-linked, right-linked,
left-shadowing and $(C(\delta),ij)$-top-balanced $Y$-configurations with $Y
\in K_0$ and $Y' \in K_0$. 

Suppose (using the notation \S\ref{sec:subsec:Yconf:notation})
there exist $\tilde{q}_2 \in K_0$ and $\tilde{q}_2' \in K_0$ such that $\sigma(\tilde{q}_2') \in \cW^u[\sigma(\tilde{q}_2)]$,
	and also $d^Q(\tilde{q}_2, q_2) < \xi$ and $d^Q(\tilde{q}_2', q_2') < \xi$, for some $\xi>0$.


Then, provided $\xi$ is small enough and $\tau=\tau(Y)$ is large enough (depending on $K_0$), we have that
\begin{equation}
\label{eq:goal:same:sheet}
\tilde{q}_2' \in \cW^u[\tilde{q}_2]. 
\end{equation}
Furthermore, there exists $\xi' > 0$ (depending on $\xi$, $K_0$ and $C$
and $t_2$)
with $\xi' \to 0$ as $\xi \to 0$ and $\tau \to \infty$
such that 
\begin{equation}
	\label{eq:nearby:maps:final}
d_*( \wp^+(\tilde{q}_2,\tilde{q}_2')_* \wtd f_{ij}[\tilde{q}_2],
\wtd f_{ij}[\tilde{q}_2']) \le \xi', 
\end{equation}
where $d_*(\cdot, \cdot)$ is as in \autoref{def:d:star}. 

\end{proposition}

In the course of the proof, we will also prove
the following analogue of \cite[Lemma 11.6]{EskinMirzakhani_Invariant-and-stationary-measures-for-the-rm-SL2Bbb-R-action-on-moduli-space}.
\begin{lemma}
\label{lemma:tau:ij:nearby:close} 
For every $\delta > 0$ there exists a subset $K_0 \subset X$ with
$\nu(K_0) > 1-\delta$ and $C(\delta) < \infty$ 
such that the following holds: Suppose $Y \in K_0$ and $Y' \in K_0$
are bottom-linked, left-shadowing and $(C,ij)$-top-balanced $Y$-configurations. Then,
\begin{displaymath}
\label{lem:almost_agreement_of_tau_and_tau_prime}
  |t_0(Y) - t_0(Y')| \le C_1(\delta), 
\end{displaymath}
where $C_1(\delta)$ depends on $C(\delta)$ and $\delta$. 
\end{lemma}

The proof of \autoref{prop:nearby:linear:maps} and
\autoref{lemma:tau:ij:nearby:close} will rely on the interpolation
map from \autoref{ssec:interpolation_map},
instead of the more ad hoc construction from
\cite{EskinMirzakhani_Invariant-and-stationary-measures-for-the-rm-SL2Bbb-R-action-on-moduli-space}.

\subsubsection{Preliminary constructions at the halfway points}
Fix a compact set $K_1$ and suppose $Y \in K_1$, $Y' \in K_1$ are
bottom-linked, right-linked,
left-shadowing and $(C,ij)$-top-balanced $Y$-configurations.
We use the notation \S\ref{sec:subsec:Yconf:notation}. Suppose also
$\tilde{q}_2 \in K_1$ and $\tilde{q}_2' \in K_1$ are as in
\autoref{prop:nearby:linear:maps}.

Recall we write  \index{$y_{1/2}$}$y_{1/2}= g_{-\ell/2} uq_1$ and \index{$z_{1/2}$}$z_{1/2}= \cW^{cs}_{fake, loc}[y_{1/2}]  \cap \cW^{u}_{loc}[q'_{1/2}] $ as in  \autoref{fig:interpolating_curves}.  
Recall also \index{$y_{1/2}'$}$y_{1/2}' := g_{-(\tau+l/2)}q_2' = g_{-\ell/2} u'q_1'$.  
By assumption $y_{1/2}$, $y_{1/2}'$, and $\wtd q'$ are in $K_1$.  

	\label{sssec:preliminary_constructions_at_the_halfway_points}
Recall that in \autoref{ssec:interpolation_map} we constructed (see \autoref{fig:interpolating_curves}):
\begin{itemize}
	\item a subresonant ``interpolation'' map
	\[
		\wtilde{\phi}(y_{1/2},z_{1/2})\colon \cW^u[y_{1/2}]\to \cW^u[z_{1/2}]
	\]
	\item an 	for any biregular point $x_{1/2}$ on $\cW^u[z_{1/2}]=\cW^u[q_{1/2}']$ a strictly subresonant map
	\[
		\wp_{fake}^+(x_{1/2},z_{1/2}) \colon\cW^u[x_{1/2}] 
		\to \cW^u[z_{1/2}]
	\]
	We will denote the  inverse of $\wp_{fake}^+(x_{1/2},z_{1/2})$ by   $\wp_{fake}^+(z_{1/2},x_{1/2})$.
\end{itemize}
For the constructions, see \autoref{thm:fake_holonomies_and_ssr_map} and \autoref{cor:interpolation_and_generalized_subspaces}.

We now define\index{$\psi$@$\wtd \psi(y_{1/2},y_{1/2}')$}
\begin{equation}
\label{eqn:psi_tilde_12_defn}
\wtd \psi(y_{1/2},y_{1/2}') =\wp_{fake}^+(z_{1/2},y_{1/2}') \circ \wtd \phi(y_{1/2} , z_{1/2}).\end{equation}
Note from the construction of $\wtd \phi(y_{1/2} , z_{1/2})$ we have alternatively
\begin{equation}
\label{eq:tilde:psi:alt}
  \wtd \psi(y_{1/2},y_{1/2}') =
	{\wp^+\left(q_{1/2}',y_{1/2}'\right)}\circ 
	{\wp^-\left(q_{1/2},q_{1/2}'\right)}\circ 
	{\wp^{+}\left(y_{1/2},q_{1/2}\right)}.\end{equation}


\subsubsection{Pushing to the upper-left corner}
	\label{sssec:pushing_to_the_upper_left_corner}
Using the dynamics, we can now push these maps to the upper-left corner of the $Y$-diagram.

Define \index{$z_2$}
\begin{align*}
	z_2&:= g_{(\tau+\ell/2)}z_{1/2} && \text{the pushforward of a halfway point}.
\end{align*}
We construct analogous objects in the upper-left corner of \autoref{fig:outline}:
\begin{itemize}
	\item An interpolation map
			\index{$\phi_\tau(q_2,z_2)$}$\phi_\tau(q_2,z_2)\colon \cW^u[q_2]\to \cW^u[z_2]$,
			$$\phi_\tau(q_2,z_2) =g_{\tau+\ell/2} \circ \wtd \phi(y_{1/2}, z_{1/2} ) \circ g_{-\tau-\ell/2}$$
		\item Strictly subresonant polynomial maps $
		\index{$P$@$\wp_{fake}^{+}(x_2,z_2)$}\wp_{fake}^{+}(x_2,z_2)\colon \cW^u[x_2]\to \cW^u[z_2]
$,
			$$\wp_{fake}^{+}(x_2,z_2) =g_{\tau+\ell/2} \circ \wp_{fake}^+(x_{1/2},z_{1/2})  \circ g_{-\tau-\ell/2}$$
	for any biregular point $x_2\in \cW^u[z_2]$.
	\item The pushforward of $\tilde{\psi}(y_{1/2},y_{1/2}')$:
$
		\index{$\psi(q_2,q_2')$}\psi(q_2,q_2')\colon \cW^u[q_2]\to \cW^u[q_2'].
$,
\begin{equation}
\label{eq:def:nontilde:psi}
\psi(q_2,q_2')=g_{\tau+\ell/2} \circ\wtd \psi(y_{1/2},y_{1/2}')  \circ g_{-\tau-\ell/2}.
\end{equation}
\end{itemize}
Note that in view of their definitions in \autoref{eq:def:nontilde:psi} and \autoref{eq:tilde:psi:alt}, we have that
\begin{align}
\label{eq:psi_tilde_psi_preserve_Zimmer_flag}
	\text{both }\psi\text{ and }\tilde{\psi} \text{ preserve the Zimmer flags}.
\end{align}
since all maps involved do, by the Ledrappier invariance principle.

\subsubsection{Convergence at tilde-points}
Let \index{$\tau(q)$}$\tau(q)$ be the measurable family of maps from \autoref{prop:agreement_of_linearizations_and_compatible_family_of_subgroups}, and set
\[
	\index{$\eta(q,q')$}\eta(q,q'):=\tau(q')^{-1}\circ \tau(q)\colon \cW^u[q] \to  \cW^u[q'].	
\]
Note that the map is subresonant, and in fact resonant if $\cW^u[q],\cW^{u}[q']$ are equipped with resonant structures induced by the Oseledets decomposition at $q,q'$.


\begin{proposition}[Boundedness and Convergence]
\label{prop:convergence_at_tilde_points}
	For any $\delta>0$ there exists a compact set $K\subset X$
	with $\nu(K)>1-\delta$ with the following property. Suppose $Y$, $Y'$ are bottom-linked, right-linked,
	left-shadowing and $(C(\delta),ij)$-top-balanced $Y$-configurations with $Y
	\in K$ and $Y' \in K$. 
	
	Then we have that:
	\begin{enumerate}
		\item The norm of $\psi(q_2,q_2')$ is bounded.
		\label{psi_2_is_bounded}
		\item 
		Suppose (using the notation \S\ref{sec:subsec:Yconf:notation})
		there exist $\tilde{q}_2 \in K$ and $\tilde{q}_2' \in K$ such that $\sigma_0(\tilde{q}_2') \in \cW^u[\sigma_0(\tilde{q}_2)]$,
		and also $d^X(\tilde{q}_2, q_2) < \xi$ and $d^X(\tilde{q}_2', q_2') < \xi$, for some $\xi>0$.
		The transformation
		 $$\eta  ({q_2', \wtilde {q}_2'})\circ
		\psi(q_2,q_2')\circ\eta({  \wtilde {q}_2, q_2}) \colon 
		\cW^u[\wtilde {q}_2]\to  	\cW^u[\wtilde {q}_2']
		$$ converges to $\wp^+(\tilde{q}_2,\tilde{q}_2')$ as
        $\tau,\ell\to +\infty$ and $\xi \to 0$. 
	\end{enumerate}
\end{proposition}



\begin{proof}[Proof of \autoref{prop:convergence_at_tilde_points}]
We recall the following notation.  Suppose $\phi$ is a subresonant map
between two unstable manifolds. Then, its linearization when acting on
the linearized unstables $L\cW^u$ is denoted \index{$L\phi$}$L\phi$.

We first prove (ii).
Note that the map \index{$\Upsilon({\wtilde {q}_2, \wtilde {q}'})$}$$\Upsilon({\wtilde {q}_2, \wtilde {q}'})=\eta(q_2',\wtd q_2')\circ
\psi(q_2,q_2')\circ\eta  (  \wtilde {q}_2, q_2) \colon 	\cW^u[\wtilde {q}_2]\to  	\cW^u[\wtilde {q}'_2]$$ is a subresonant diffeomorphism taking $\wtilde {q}_2$ to  $\wtilde {q}'_2$.
We have $\psi(q_2,q_2') =\wp_{fake}^{+}(z_2,q_2') \circ
\phi(q_2,z_2)$ and in view of
(\ref{eq:psi_tilde_psi_preserve_Zimmer_flag}), $L\psi(q_2,q_2')$ intertwines the splittings of $L\cW^u[q_2]$ and $L\cW^u[q_2']$.

Assuming $\wtd q_2',\wtd q_2, q_2'$, and $q_2$ all belong to a fixed
Lusin set $K_2$, given $\epsilon>0$ we may assume 
\begin{enumerate}
\item $\norm{\gr_\bullet L(\eta  (  \wtilde {q}_2, q_2))-\gr_\bullet I^{L \cW^u}_{\wtd q_2, q_2}}< \epsilon $ \hfill
\item $\norm{\gr_\bullet L(\eta(q_2',\wtd q_2')) - \gr_\bullet I^{L \cW^u}_{q_2',\wtd  q_2'}} < \epsilon $ 
\item $\norm{\gr_\bullet L\phi (q_2,z_2) - \gr_\bullet I^{L \cW^u}_{q_2, z_2}}<\epsilon $ \hfill
\item $\norm{\gr_\bullet (I^{L \cW^u}_{q_2',\wtd q_2'}  \circ P_{fake}^{+}(z_2,q_2') \circ  I^{L \cW^u}_{\wtd q_2,z_2}  ) -  \gr_\bullet P^{+}(\wtd q_2,\wtd q_2') }<\epsilon.$  \hfill 
\end{enumerate}
The first two estimates follow by applying \autoref{prop:closeness_to_the_identity} to the trivialization of $L\cW^u$ provided by \autoref{prop:agreement_of_linearizations_and_compatible_family_of_subgroups}, which gives the analogous estimates but without $\gr_{\bullet}$.
But the estimate on $L\cW^u$ implies the one on $\gr_{\bullet}L\cW^u$.

For the last two estimates, note that $z_2$ is not assumed biregular, it is however backwards-regular and we can obtain the needed bounds.
Specifically, the third estimate follows from \autoref{prop:distortion_of_phi}.

For the fourth estimate, recall that $\gr_\bullet P_{fake}^{+}(z_2,q_2') = H_{z_2,q_2'}^{L\cW^u}$ where \index{$H_{z_2,q_2'}^{L\cW^u}$}$H_{z_2,q_2'}^{L\cW^u}$ is the holonomy on $\gr_\bullet L\cW^u$ as in \autoref{cor:holonomies_on_graded}.
The estimate follows from the \Holder continuity of holonomies on Pesin sets, as explicated in \autoref{cor:holonomies_on_graded}, \autoref{eqn:holonomy_graded_estimate}.

It now follows that
\begin{align*}
\|\gr_\bullet L\Upsilon_{\wtilde {q}_2, \wtilde {q_2}'} - & \gr_\bullet  P^{+}(\wtd q_2,\wtd q_2') \| \\
&\le  \| \gr_\bullet I^{L \cW^u}_{q_2',\wtd  q_2'} \gr_\bullet  P_{fake}^{+}(z_2,q_2') \gr_\bullet I^{L \cW^u}_{q_2, z_2}    \gr_\bullet  I^{L \cW^u}_{\wtd q_2, q_2}- \gr_\bullet  P^{+}(\wtd q_2,\wtd q_2')\| + C\epsilon \\
&=  \| \gr_\bullet I^{L \cW^u}_{q_2',\wtd  q_2'} \gr_\bullet  P_{fake}^{+}(z_2,q_2')  \gr_\bullet  I^{L \cW^u}_{\wtd q_2, z_2}- \gr_\bullet  P^{+}(\wtd q_2,\wtd q_2')\| + C\epsilon \\
&\le(C+1)\epsilon.
\end{align*}
It follows that 
$\gr_\bullet L\Upsilon_{\wtilde {q}_2, \wtilde {q_2}'} $ converges to $\gr_\bullet  P^{+}(\wtd q_2,\wtd q_2') $.  Since $L\Upsilon_{\wtilde {q}_2, \wtilde {q_2}'}$ intertwines the splittings, it follows that 
$ L\Upsilon_{\wtilde {q}_2, \wtilde {q_2}'} $ converges to $ P^{+}(\wtd q_2,\wtd q_2') $.
This implies  $\Upsilon_{\wtilde {q}_2, \wtilde {q_2}'}$ converges to  $\wp^{+}(\wtd q_2,\wtd q_2') $ and finishes the proof of (ii).

To prove (i), recall that
\[
	\gr_{\bullet}L\psi(q_2,q_2') = H^{L\cW^u}_{z_2,q_2'}\circ \gr_{\bullet}L\phi(q_2,z_2)
\]
where $H^{L\cW^u}$ is the holonomy on $\gr_{\bullet}L\cW^u$.
By \autoref{cor:holonomies_on_graded} the holonomy $H^{L\cW^u}_{z_2,q_2'}$ and $I^{L\cW^u}_{z_2,q_2'}$ are bounded distance apart (in fact arbitrarily close), which together with the third estimate above implies that $\gr_{\bullet}L\psi(q_2,q_2')$ is bounded on $K_2$.
Note that $L\psi$ preserves the Oseledets decomposition by construction, and the isomorphism induced by the Oseledets splitting between $\gr_{\bullet}L\cW^u$ and $L\cW^u$ is bounded on compact sets of arbitrarily large measure.
It follows that $L\psi(q_2,q_2')$ itself is bounded on such sets.
\end{proof}

%


\begin{remark}
Consider two points $a,b\in Q$ and a subresonant map $\psi \colon \cW^u[a]\to \cW^u[b].$  Suppose moreover that $\psi$ intertwines a compatible family of subgroups as in \autoref{def:compatible_family_of_subgroups} at $a$ and $b$; that is, suppose that 
$\psi(U^+[a]) = U^+[b]$ so that $\psi$ induces a map $\cC[a]\to \cC[b]$.  If moreover $\psi(a) = b$ then this map descends to a well-defined map \index{$L$@$\bbL\psi$}$$\bbL\psi\colon \bbH(a) \to \bbH(b).$$
\end{remark}
The $\psi$ from \autoref{eq:def:nontilde:psi} has this property, and this is the map we use going forward.
Recall that $\bbR(\cdot,\cdot)$ was defined in (\ref{eq:def:R}). Then
$\bbR(\cdot, \cdot)$ acts on $\bbH$, and thus on
$\bbE \subset \bbH$.

To single out the measurable connections on $\bbH$, we will write \index{$P$@$\bbP^\pm$}$\bbP^{\pm}$ for $P^\pm_\bbH$.

\begin{proposition}[Agreement of transports]
	\label{prop:agreement_of_transports}
Suppose $Y$, $Y'$ are bottom-linked, left-shadowing and right-linked
$Y$-configurations. We use the notation from
\autoref{sec:subsec:Yconf:notation}. 
	Consider the two maps
	\begin{align*}
		\bb\bbR(q_3',q_2')\circ \bbP^-(q_3,q_3')
		\colon &
		\bbE(q_3) \to \bbE(q_2')\\
		\bbL \psi(q_2,q_2')
 \circ \bb\bbR(q_3,q_2)
		\colon &
		\bbE(q_3) \to \bbE(q_2')
	\end{align*}
	\begin{enumerate}
		\item The two maps agree, when all points belong to a set of full measure on which the maps are defined.
		\item \label{item2:agreement_of_transports}
		For any $\delta>0$ there exists a compact set
		$K_1\subset X$ with $\nu(K)>1-\delta$, such that
		if $Y \in K$, $Y' \in K$ and both $Y$ and $Y'$
		are $(C,ij)$-top-balanced,
		then each of the maps above is of bounded norm on
		$\bbE_{[ij],bdd}$, with bound on the norm depending
		only on $\delta$ and $C$.
		\end{enumerate}
\end{proposition}
\begin{proof}
We first prove (i).
Note that by
\autoref{prop:agreement:of:u:star:and:p}, we have (acting on $\bbE$)
\begin{equation}
\label{eq:Rq3q2:alt}
\bbR(q_3,q_2) = g_{(\tau+\ell/2)} \circ \bbP^+(q_{1/2}, y_{1/2})
\circ g_{-(\tau+\ell/2)}
\end{equation}
Similarly,
\begin{equation}
\label{eq:Rq3q2:prime:alt}
\bbR(q_3',q_2') = g_{(\tau+\ell/2)} \circ \bbP^+(q_{1/2}', y_{1/2}')
\circ g_{-(\tau+\ell/2)}
\end{equation}
Also from (\ref{eq:tilde:psi:alt}), we have
\begin{equation}
\label{eq:Pplus:Pminus:psi}
\bbP^+(q_{1/2}', y_{1/2}') \circ \bbP^-(q_{1/2},q_{1/2}') = \bbL
\tilde{\psi}(y_{1/2},y_{1/2}') \circ \bbP^+(q_{1/2},y_{1/2}).   
\end{equation}
Therefore
\begin{align*}
  \bbR(q_3',q_2') & \circ \bbP^-(q_3,q_3')
  = g_{(\tau+\ell/2)} \circ \bbP^+(q_{1/2}',y_{1/2}') \circ
    \bbP^-(q_{1/2},q_{1/2}') \circ g_{-(\tau+\ell/2)}
                           & \text{by~(\ref{eq:Rq3q2:prime:alt})} \\
  & = g_{(\tau+\ell/2)} \circ
  \bbL\tilde{\psi}(y,y_{1/2}') \circ \bbP^+(q_{1/2},y_{1/2})
    \circ g_{-(\tau+\ell/2)} \qquad & \text{by~(\ref{eq:Pplus:Pminus:psi})} \\
    & = \bbL\psi(q_2,q_2') \circ g_{(\tau+\ell/2)} \circ
      \bbP^+(q_{1/2},y_{1/2}) 
    \circ g_{-(\tau+\ell/2)} \qquad & \text{by~(\ref{eq:def:nontilde:psi})} \\
   & = \bbL\psi(q_2,q_2') \circ \bbR(q_3,q_2) \qquad & \text{by~(\ref{eq:Rq3q2:alt})}
\end{align*}

To prove (ii), note that by the Ledrappier invariance principle (\autoref{thm:ledrappier_invariance_principle}), 
$\bbP^-(q_3,q_3') \bbE_{[ij],bdd}(q_3) = \bbE_{[ij],bdd}(q_3')$. Also,
by the choice of the Lusin set $K$, we can make sure that
$\|\bbP^-(q_3,q_3') \| \le C(\delta)$. We may also choose $K$ so
that the the function $C(\cdot)$ of
\autoref{prop:ej:bdd:transport:bounded} is bounded by
$C(\delta)$ on $K$. Then, by
\autoref{prop:ej:bdd:transport:bounded}, for $\bfv' \in
\bbE_{[ij],bdd}(q_3')$, $\|\bbR(q_3',q_2') \bfv' \| \le C(\delta) \|
\bfv'\|$. This completes the proof of (ii). 
\end{proof}

\begin{figure}[ht!]
\includegraphics{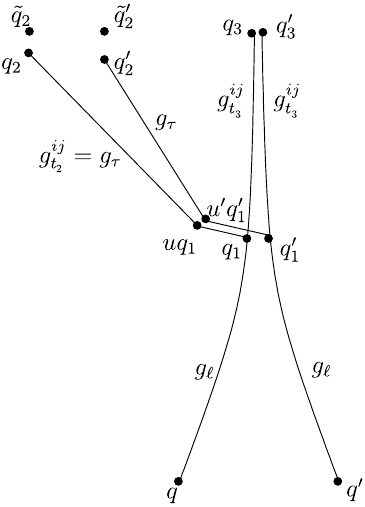}
\caption{}\label{fig:outline2}
\end{figure}

\subsubsection{Proof of \autoref{prop:nearby:linear:maps}}
	\label{sssec:proof_of_prop:nearby:linear:maps}
We first have 
\[
 	 d_*(\bbP^-(q_3,q_3') \bbf_{ij} (q_3) , \bbf_{ij}(q_3'))<\xi^{(2)}
 	 \text{ with $\xi^{(2)}\to 0$ as $\xi\to 0$.}
\]
We also have (by \autoref{prop:agreement_of_transports}\autoref{item2:agreement_of_transports}) that $\bbR(q_3',q_2')$  is a bounded map. 
By \autoref{lemma:transformation:rule:bff}, we have $d_*\left(\bbR(q_3',q_2')  \bbf_{ij}(q_3'),  \bbf_{ij}(q_2')\right)=0$.  
Assembling, we have $$ d_*\left(\bbR(q_3',q_2') \circ \bbP^-(q_3,q_3') \bbf_{ij} (q_3) , \bbf_{ij}(q_2')\right)<\xi^{(3)}
 	 \text{with $\xi^{(3)}\to 0$ as $\xi\to 0$.}$$

We again have $d_*\left(\bbR(q_3,q_2)\bbf_{ij}(q_3),\bbf_{ij}(q_2)\right) =0$
so 
\begin{align*}
d_*\big(\bbL \psi(q_2,q_2')
 \bbf_{ij} (  q_2), \bbf_{ij} (  q_2')\big) 
&=  d_*\left(\bbR(q_3',q_2') \circ \bbP^-(q_3,q_3')
 \bbf_{ij} (  q_3), \bbf_{ij} (  q_2')\right) \end{align*}
is bounded above by $\xi^{(4)}\to 0$ as $\xi\to 0$.
 
 By the above and \autoref{prop:convergence_at_tilde_points}, we have 
\begin{align*}
d_*(\bbP^+ (\wtd q_2,\wtd q_2') \bbf_{ij} (\wtd q_2), \bbf_{ij} (\wtd q_2')) 
&\le 
d_*(\bbP^+ (\wtd q_2,\wtd q_2') \bbf_{ij} (\wtd q_2), 
I^{\bbH}_{ q_2', \wtd q_2'}
\bbL\psi (  q_2,  q_2') I^{\bbH}_{\wtd q_2, q_2} \bbf_{ij} (\wtd q_2))
\\
&\quad + d_*( I^{\bbH}_{ q_2', \wtd q_2'} \bbL\psi (  q_2,  q_2') I^{\bbH}_{\wtd q_2, q_2} \bbf_{ij} (\wtd q_2)),
I^{\bbH}_{ q_2', \wtd q_2'} \bbL\psi (  q_2,  q_2') \bbf_{ij} (q_2)))
\\
&\quad + d_*( I^{\bbH}_{ q_2', \wtd q_2'} \bbL\psi (  q_2,  q_2') \bbf_{ij} (q_2), I^{\bbH}_{ q_2', \wtd q_2'}  \bbf_{ij} (q_2'))\\
&\quad + d_*( I^{\bbH}_{ q_2', \wtd q_2'}  \bbf_{ij} (q_2'), \bbf_{ij} (\wtd q_2')) 
\end{align*}
and all four terms on the right go to zero as $\xi\to 0$.

We have now shown $$d_*(\bbP^+ (\wtd q_2,\wtd q_2') \bbf_{ij} (\wtd q_2), \bbf_{ij} (\wtd q_2')) <\xi^{(5)}$$ for some $\xi^{(5)}$.
Pushing forward both measures by $(\bbj_{\wtd q_2'} \circ \cU_{\wtd q_2'}^{-1})^{-1}$ we have that 
$$d_*(\wp^+ (\wtd q_2,\wtd q_2') \overline f_{ij} [\wtd q_2], \overline f_{ij} [\wtd q_2']) <\xi^{(6)}.$$
By the construction in \autoref{sssec:recovering_wtilde_f_ij_x_from_bar_f_ij}, it also follows that:
$$d_*(\wp^+ (\wtd q_2,\wtd q_2') \wtd f_{ij} [\wtd q_2], \wtd f_{ij} [\wtd q_2']) <\xi^{(7)}$$
and the inequality asserted in \autoref{eq:nearby:maps:final} holds.

The convergence in \eqref{eq:goal:same:sheet}  now follows since the
$\psi$ maps converge to the $P^+$ maps and a point in the lifted
unstable is a point downstairs plus a flag and $P^+$  intertwines the
flags.  (See discussion in
\S\ref{sec:subsec:stable:manifold:measurable:cover}).  
\hfill \qed

\subsubsection{Proof of \autoref{lem:almost_agreement_of_tau_and_tau_prime}}
\label{sssec:proof_of_almost_agreement_of_tau_and_tau_prime}

We will need the following:
\begin{lemma}
\label{lemma:lambdaij:stable:nearby:close}
For every $\delta > 0$ there exists a compact $K \subset X$ with
$\nu(K) > 1-\delta$ such that the following holds: Suppose $x
\in K$, $x' \in \gB_0^-[x] \subset \cW^s[x]$ 
and suppose $t > 0$ is such that $g_t x \in K$, $g_t x' \in
K$. Then,
\begin{equation}
\label{eq:lambda:ij:x:t:minus}
|\lambda_{ij}(x,t) - \lambda_{ij}(x',t)| \le C(\delta). 
\end{equation}
\end{lemma}
\begin{proof}
For any $\delta > 0$ we can choose a subset $K_0 = K_0(\delta)$ such
that $\nu(K_0) > 1-\delta$ and also the Lyapunov splitting of $\bbH$ and the
dynamical norm $\|
\cdot \|$ vary continuously with respect to a trivialization and
ambient norm provided by \autoref{thm:ledrappier_invariance_principle}.
Then, provided $x,x', g_t x,
g_t x' \in K_0$ and $x' \in \gB_0^-[x]$, 
we have
\begin{equation}
\label{eq:bbPminus:bounded}
\|\bbP^-(x,x')\| \le C(\delta), \qquad \|\bbP^-(g_t x, g_t x')\| \le
C(\delta). 
\end{equation}
From \autoref{eq:growth_v_in_Elambda_i_j} which defined
$\lambda_{ij}$, let $\bbv \in \bbE^{\lambda_i}_{\leq
  j}(x)/\bbE^{\lambda_i}_{\leq j-1}(x)$ and let
$\bbv':=\bbP^-(x,x')\bbv$. In view of (\ref{eq:bbPminus:bounded}), 
$\norm{\bbv}$ and $\norm{\bbv'}$ are within a bounded multiplicative constant.

Next we have
\begin{align*}
	\bbv_2:=g_{t}\bbv &\text{ and } \bbv_2':=g_{t}\bbv'\\
	 \norm{\bbv_2}:=e^{\lambda_{ij}(x,t)}\norm{\bbv} & \text{ and }  \norm{\bbv_2'}:=e^{\lambda_{ij}(x',t)}\norm{\bbv'}
\end{align*}
We also have $\bbP^-(g_t x,g_t x') \bfv_2 = \bfv_2'$. Thus, in view of
(\ref{eq:bbPminus:bounded}), $\norm{\bbv_2}$ and $\norm{\bbv'_2}$ are also
within a bounded multiplicative constant depending on $\delta$, and so 
(\ref{eq:lambda:ij:x:t:minus}) follows. 
\end{proof}

\begin{proof}[Proof of \autoref{lem:almost_agreement_of_tau_and_tau_prime}.]
We use the same notation as in the proof of
\autoref{prop:nearby:linear:maps}.
Let $\tilde{\psi}$ and
$\psi$ be as in the proof of
\autoref{prop:nearby:linear:maps}.

First, we check that
$|\lambda_{ij}(y_{1/2};\tau+\ell/2)-\lambda_{ij}(y_{1/2}';\tau+\ell/2)|<c_1$
for a constant that depends only on $\delta$. This argument is similar
to the proof of \autoref{lemma:lambdaij:stable:nearby:close}. 
From \autoref{eq:growth_v_in_Elambda_i_j} which defined $\lambda_{ij}$, let $\bbv \in \bbE^{\lambda_i}_{\leq j}/\bbE^{\lambda_i}_{\leq j-1}(y_{1/2})$ and let $\bbv':=\bbL\wtilde{\psi}(y_{1/2},y_{1/2}')\bbv$.
Note that $\norm{\wtilde{\psi}}$ is uniformly bounded in terms of
$\delta$ by \autoref{prop:convergence_at_tilde_points} (i), so $\norm{\bbv}$ and $\norm{\bbv'}$ are within a bounded multiplicative constant.

Next we have
\begin{align*}
	\bbv_2:=g_{\tau+l/2}\bbv &\text{ and } \bbv_2':=g_{\tau+l/2}\bbv'\\
	 \norm{\bbv_2}:=e^{\lambda_{ij}\left(y_{1/2},\tau+l/2\right)}\norm{\bbv} & \text{ and }  \norm{\bbv_2'}:=e^{\lambda_{ij}\left(y_{1/2},\tau+l/2\right)}\norm{\bbv'}
\end{align*}
Recall that by \eqref{eq:psi_tilde_psi_preserve_Zimmer_flag}, the map $\bbL\psi$ preserves the Zimmer flags and intertwines the dynamics.
So $\bbL\psi(q_2,q_2')\bbv_2=\bbv_2'$.
Furthermore, by \autoref{prop:convergence_at_tilde_points}\autoref{psi_2_is_bounded}, $\bbL\psi(q_2,q_2')$ is bounded in norm, with bound depending on $\delta$.
So $\norm{\bbv_2}$ and $\norm{\bbv_2'}$ are within a bounded multiplicative constant depending on $\delta$, and so the bound
\[
	|\lambda_{ij}(y_{1/2};\tau+\ell/2)-\lambda_{ij}(y_{1/2}';\tau+\ell/2)|<C_1(\delta) \text{ follows.}
\]
Since $Y$ is $(C,ij)$-top-balanced, 
\begin{displaymath}
|\lambda_{ij}(y_{1/2};\tau+\ell/2)-\lambda_{ij}(q_{1/2};t_0(Y)+\ell/2)|
\le C.
\end{displaymath}
Similarly, since $Y'$ is $(C,ij)$-top-balanced,
\begin{displaymath}
|\lambda_{ij}(y_{1/2}';\tau+\ell/2)-\lambda_{ij}(q_{1/2}';t_0(Y')+\ell/2)|
\le C.
\end{displaymath}
Combining the above three displayed equations, we get
\begin{displaymath}
\lambda_{ij}(q_{1/2}';t_0(Y')+\ell/2)-\lambda_{ij}(q_{1/2};t_0(Y)+\ell/2)|
\le C(\delta). 
\end{displaymath}
By \autoref{lemma:lambdaij:stable:nearby:close},
\begin{displaymath}
\lambda_{ij}(q_{1/2}';t_0(Y')+\ell/2)-\lambda_{ij}(q_{1/2};t_0(Y')+\ell/2)|
\le C(\delta). 
\end{displaymath}
Therefore, 
\begin{displaymath}
|\lambda_{ij}(q_{1/2};t_0(Y')+\ell/2)-\lambda_{ij}(q_{1/2};t_0(Y)+\ell/2)|
\le C(\delta), 
\end{displaymath}
and the bound $|t_0(Y)-t_0(Y')|<C(\delta)$ follows because
$\lambda_{ij}$ is monotonic and with Lipschitz constant bounded above
and below by \autoref{prop:good_norms}.
\end{proof}

\section{Extra Invariance}
	\label{sec:extra_invariance}

\paragraph{Outline of section}
In this section we prove \autoref{thm:inductive:step}.
The proof of \autoref{thm:random:inductive:step} is analogous and follows the same steps, with only straightforward notational changes.

Recall from \autoref{def:measures_fij} and \autoref{def:measures_fij_tilde} that $f_{ij}[x]$ denotes conditional measures supported on the atoms $\Psi_{ij}[x]$, while $\tilde{f}_{ij}[x]$ denotes (conditional) leafwise measures, which locally on $\Psi_{ij}[x]$ agree with $f_{ij}[x]$, and are extended equivariantly (up to scaling) for the dynamics.

\subsubsection{Measurable trivializations of measures}
	\label{sssec:measurable_trivializations_of_measures}
In order to describe the continuity properties of the measures $\tilde{f}_{ij}[x]$ on $\cW^u[x]$, we make the following construction.
Let $T(x)\colon L\cW^u(x)\to \bR^{n_u}$ be a measurable trivialization of the bundle $L\cW^u(x)$, compatible with the appropriate subresonant structure on $\bR^{n_u}$.
Recall that $L_x\colon \cW^u[x]\to L\cW^u(x)$ denotes the linearization map, and set \index{$f$@$\frakf_{ij}(x)$}$\frakf_{ij}(x):=T(x)_* (L_x)_* \wtilde{f}_{ij}[x]$ to be a measure on $\bR^{n_u}$.

Let $\wp^+(x,y)$ be the map introduced in \autoref{prop:subresonant_maps_from_holonomies}. 
\autoref{thm:inductive:step} will be derived from the
following:
\begin{proposition}
\label{prop:half:inductive:step}
Suppose $U^+$
are as in \autoref{thm:inductive:step} and that the QNI
condition of \autoref{def:QNI} holds.
Then there exists $0 < \delta_0 < 0.1$, a subset $K_* \subset X$ with
$\nu(K_*) > 1-\delta_0$ such that
all the functions $\frakf_{ij}$, $\forall ij \in \tilde{\Lambda}$ 
are uniformly continuous on $K_*$, and $C > 1$ (depending on $K_*$) 
such that for every $0 < \epsilon < C^{-1}/100$ there exists a 
subset $E \subset K_*$ with 
$\nu(E) > \delta_0$, such that for every $x \in E$ there exists $ij
\in \tilde{\Lambda}$ and $y \in \Psi_{ij}[x] \cap K_*$ with 
\begin{equation}
\label{eq:dxy:within:C:epsilon}
	C^{-1} \epsilon \le d_{\cH,loc}^x(U^+[x],U^+[y]) \le C \epsilon
\end{equation}
and furthermore:
\begin{equation}
	\label{eq:half:inductive:step:two}
	\tilde{f}_{ij}[y] \propto \wp^+(x,y)_* \tilde{f}_{ij}[x]. 
\end{equation}
\end{proposition}
\noindent To establish that $\nu(E)>\delta_0$ we will show that any compact set $K_{00}$ with $\nu(K_{00})>1-2\delta_0$ contains a point $x\in K_{00}\cap K_*$ satisfying \autoref{eq:dxy:within:C:epsilon} and \autoref{eq:half:inductive:step:two}.

A broad outline of the proof of \autoref{prop:half:inductive:step} is contained in \autoref{ssec:outline_proof_inductive_step}.


\subsection{Having friends, proof of \autoref{prop:half:inductive:step}}

\subsubsection{Choice of parameters \#1.}
Fix $\theta > 0$ as in \autoref{prop:some:fraction:bounded} 
and \autoref{prop:ej:bdd:transport:bounded}.
We then choose
$\delta > 0$ sufficiently small; the exact value of $\delta$ will we
chosen at the end of this section. All subsequent constants
will depend on $\delta$. (In particular, $\delta \ll \theta$; we will
make this more precise below). Let $\epsilon > 0$ be arbitrary and
$\eta > 0$ be arbitrary;  however we will always assume that
$\epsilon$ and $\eta$ are sufficiently small depending on $\delta$.

We will show that \autoref{prop:half:inductive:step} holds
with $\delta_0 = \delta/10$. Let $K_* \subset X$ be any subset with
$\nu(K_*) > 1 - \delta_0$ on which all the functions $\frakf_{ij}$
are uniformly continuous. 
It is enough to show that there exists $C
= C(\delta)$ such that for any $\epsilon > 0$ and 
for an arbitrary compact set $K_{00} \subset X$ with $\nu(K_{00}) \ge (1-2
\delta_0)$, there exists $x \in K_{00} \cap K_*$, $ij \in
\tilde{\Lambda}$  and $y \in \Psi_{ij}[x] \cap K_*$ 
satisfying (\ref{eq:dxy:within:C:epsilon}) and
(\ref{eq:half:inductive:step:two}). Thus, let $K_{00} \subset X$ be an
arbitrary compact set with $\nu(K_{00}) > 1- 2\delta_0$. 

We can choose a compact set $K_0 \subset K_{00} \cap K_*$ with
$\nu(K_0) > 1-5 \delta_0 = 1-\delta/2$ so that 
\autoref{prop:nearby:linear:maps} and
\autoref{lemma:tau:ij:nearby:close} hold, and also
\autoref{lemma:swithching} holds for $K_0$ in place of $K$.    
In addition, 
there exists $\epsilon_0'(\delta) > 0$ such that for all $x \in K_0$, 
\begin{equation}
\label{eq:def:epsilon0:prime}
d^u(x,\partial \gB_0[x]) > 2 \epsilon_0'(\delta).
\end{equation}
Here, $d^u(\cdot, \cdot)$ is as in \S\ref{sec:subsec:distances} and
  by $\partial \gB_0[x]$ we mean the boundary of 
  $\gB_0[x]$ as a subset of  $\cW^u[x]$. 


We now apply \autoref{lemma:parts:staying:close} with $\epsilon_0'(\delta)$ and assume that the set appearing $K$ in the statement agrees with $K_0$.
We thus obtain $\epsilon_0(\delta)$ with the properties listed in the Lemma.
We will always assume that $\epsilon < \epsilon_0(\delta)<0.01 \epsilon'(\delta)$. 
  

Let $\kappa > 1$ be as in
\autoref{prop:bilip:hattau:epsilon}, and so that
(\ref{eq:bilip:lambda:ij}) and \autoref{prop:good_norms}
hold. 
Without loss of generality, assume $\delta < 0.01$. We now choose a
subset $K \subset K_0 \subset X$ with $\nu(K) > 1-\delta$ such that the
following hold:
\begin{itemize}
\item There exists a number $T_0(\delta)$ such that for any $x \in K$
  and any $T > T_0(\delta)$, 
\begin{displaymath}
\{ t \in [-T/2,T/2] \st g_t x \in K_0\} \ge (1-\kappa^{-2}/4) T, 
\end{displaymath}
where $\kappa$ is as in
\autoref{prop:good_norms}\autoref{good_norms_bilipschitz}. 
(This can be done by the Birkhoff ergodic theorem).
\item For each $q_1 \in K$ there exists a subset $Q(q_1) \subset
  \cB_0(q_1)$ with $|Q(q_1)| \ge (1-\delta)| \cB_0(q_1)|$
  such that

\begin{equation}
\label{eq:continue:indction:quantified}
\text{ for $u \in Q(q_1)$, $\ell_0(q_1,u) < C(\delta)$,}
\end{equation}
where  $\ell_0(\cdot, \cdot)$ is as in
\autoref{lemma:lower:bound:on:3cA}. (This is the only place where
the QNI condition is used). 
\item \autoref{prop:most:inert} holds.
\item
Let
\begin{equation}
\label{eq:def:T1:delta}
  T_1(\delta) = \frac{1}{\kappa} \log 4 c(\rho(\delta))^{-1} C(\delta),
\end{equation}
where $C(\delta)$ is as in (\ref{eq:rho:delta:le:Fq:minus:Fqprime}),
$\rho(\delta)$ is as in Proposition~\ref{prop:can:avoid:most:Mu}, and
$c(\cdot)$ is as in Lemma~\ref{lemma:bad:subspace}. Then, 
\autoref{prop:some:fraction:bounded} holds for $K$, with $T=T_1(\delta)$.
\item \autoref{theorem:interpolation:does:not:move:points} and
  \autoref{prop:factorization:main} hold. 
\item There exists a constant $C = C(\delta)$ such that for $x \in K$,
  $C_3(x)^2 < C(\delta)$ where $C_3$ is as in 
  \autoref{prop:ej:bdd:transport:bounded}. 

\item For $x \in K$,
  the function $c_1(x)$ of
  \autoref{lemma:hausdorff:distance:to:norm} is bounded from below by
  $C(\delta)^{-1}$. 
\item \autoref{lemma:hodge:norm:vs:dynamical:norm} holds for
  $K=K(\delta)$ and $C_1 = C_1(\delta)$.
\item \autoref{lemma:cA:ustar}
  holds for $K$. 
\item \autoref{cor:tau:u:bilipshitz} and
  \autoref{lemma:cA:transport} hold for $K$ in place of
  $K'$.
\item For $x \in K$, $r(x) \ge \delta$, where $r(\cdot)$ is as in
  \autoref{rmk:on_lyapunov_radius} and \autoref{lemma:du:vs:dQ}, and the constant appearing in the conclusion from \autoref{lem:lemma_stay_in_chart_stay_in_loc_new} are bounded by $C(\delta)$.
\item The map $q \to \gB_0[q]$ is uniformly continuous on $K$, see
  \autoref{ssec:measurable_partitions}. 
\item The map $q \to \cB_0[q]$ is uniformly continuous on $K$.
\end{itemize}


\begin{definition}[$Y_{ij,balanced}$]
	\label{def:y_ij_fully_balanced}
	Let \index{$Y_{ij,balanced}(q_1,u,\ell)$}$Y_{ij,balanced}(q_1,u,\ell)$ be the unique $(1,\epsilon)$-left
	balanced, $(0,ij)$-top-balanced $Y$ configuration with $q_1(Y) = q_1$,
	$u(Y) = u$, and $\ell(Y) = \ell$.	
\end{definition}

\begin{claim}
\label{claim:12.3}
There exists $\ell_3 = \ell_3(K_{00},\delta,\epsilon,\eta) > 0$, a  
set $K_3 = K_3(K_{00},\delta,\epsilon,\eta)$ of measure at least  
$1-c_3(\delta)$ and for each $q_1 \in K_3$ a subset $Q_3 =
Q_3(q_1,K_{00},\ell,\delta,\epsilon,\eta) \subset 
\cB_0(q_1)$  of
measure at least $(1-c_3'(\delta))|\cB_0(q_1)|$ such that for each
$q_1 \in K_3$ and $u \in Q_3$ there exists a subset $E_{good}(q_1,u)
\subset \reals^+$ such that the density of $E_{good}(q_1,u)$ (for
$\ell > \ell_3$) is at least $1-c_3''(\delta)$ such that if $q_1 \in
K_3$, $u \in Q_3$ and $\ell \in E_{good}(q_1,u)$, we have for each
$ij$, 
\begin{displaymath}
Y_{ij,balanced}(q_1,u,\ell) \in K.
\end{displaymath}
Also, we have
$c_3(\delta)$, $c_3'(\delta)$ and
$c_3''(\delta) \to 0$ as $\delta \to 0$.
\end{claim}

\noindent The proof is a routine application of the biLipschitz estimates from \autoref{prop:bilip:hattau:epsilon}:

\begin{proof}[Proof of claim]
Let 
\begin{displaymath}
\tilde{\cD}_{00}(q_1) =
\tilde{\cD}_{00}(q_1,K_{00},\delta,\epsilon,\eta) = \{ t > 0 \st g_t
q_1 \in K \text{ and } g_{-t} q_1 \in K\}.   
\end{displaymath}
For $ij \in \tilde{\Lambda}$, let
\begin{displaymath}
\tilde{\cD}_{ij}(q_1) =
\tilde{\cD}_{ij}(q_1,K_{00},\delta,\epsilon,\eta) 
= \{ \lambda_{ij}(q_1,t) \st g_t q_1
\in K, \quad t > 0 \}.
\end{displaymath}
Then by the ergodic theorem and (\ref{eq:bilip:lambda:ij}),
there exists a set $K_{\cD} = K_{\cD}(K_{00},\delta,\epsilon,\eta)$ 
with $\nu(K_{\cD})
\ge 1-\delta$ and 
$\ell_\cD = \ell_\cD(K_{00},\delta,\epsilon,\eta) > 0$ 
such that for $q_1 \in K_{\cD}$ and all $ij
\in \{00\} \cup \tilde{\Lambda}$,  $\tilde{\cD}_{ij}(q_1)$ has
density at least $1-2\kappa\delta$ for $\ell > \ell_{\cD}$. Let
\begin{displaymath}
E_2(q_1,u) = E_2(q_1,u,K_{00},\delta,\epsilon,\eta) = 
\{ \ell \st g_{\tau_{(\epsilon)}(q_1,  u, \ell)}u q_1 \in K \}, 
\end{displaymath}
From the definition of $E_2(q_1,u)$, it follows that if $\ell \in
E_2(q_1,u)$ then for any $ij$,
\begin{equation}
\label{eq:q2:Yijbalanced:in:K}
q_2(Y_{ij,balanced}(q_1,u,\ell)) \in K. 
\end{equation}
Let
\begin{multline*}
E_3(q_1,u) = E_3(q_1,u,K_{00},\delta,\epsilon,\eta) = \\ = \{ \ell \in
E_2(q_1,u) \st \forall ij \in \tilde{\Lambda}, 
\quad \lambda_{ij}(u q_1,\tau_{(\epsilon)}(q_1,u,\ell)) \in
\tilde{\cD}_{ij}(q_1) \}. 
\end{multline*}
Note that  for any $ij$ we have $\lambda_{ij}(u q_1,\tau_{(\epsilon)}(q_1,u,\ell)) \in
\tilde{\cD}_{ij}(q_1)$
if and only if
\begin{displaymath}
\lambda_{ij}(uq_1, \tau_{(\epsilon)}(q_1,u,\ell)) =
\lambda_{ij}(q_1, s) \text{ and } g_s q_1 \in K. 
\end{displaymath}
Thus, if $\ell \in E_3(q_1,u)$, then for any $ij$,
(\ref{eq:q2:Yijbalanced:in:K}) holds and also
\begin{equation}
\label{eq:q3:Yijbalanced:in:K}
q_3(Y_{ij,balanced}(q_1,u,\ell)) \in K.
\end{equation}

We choose $K_2 = K \cap K_\cD$, and
\begin{displaymath}
K_3 = K_2 \cap \{ x \in X \st |\{ u \in \cB_0(x) \st ux \in K_2
\}| > (1-\delta)|\cB_0(x)|\}.
\end{displaymath}
Suppose $q_1 \in K_3$,
and $u q_1 \in K_2$. 
Let 
\begin{displaymath}
E_{bad} = \{ t \st g_t u q_1 \notin K \}. 
\end{displaymath}
Then, since $u q_1 \in K_\cD$, for $\ell > \ell_{\cD}$,
 the density of $E_{bad}$ is at most $2 \kappa \delta$. We have
\begin{displaymath}
E_2(q_1,u)^c = \{ \ell \st
\tau_{(\epsilon)}(q_1,u,\ell) \in E_{bad} \}. 
\end{displaymath}
Then, by \autoref{prop:bilip:hattau:epsilon}, for $\ell >
\kappa \ell_{\cD}$,
the density of $E_2(q_1,u)$ is at least $1-4\kappa^2\delta$.

Let
\begin{displaymath}
\hat{\cD}(q_1,u) = \hat{\cD}(q_1,u,K_{00},\delta,\epsilon,\eta) = \{ t
\st \forall ij \in \tilde{\Lambda}, \quad 
\lambda_{ij}(u q_1, t) \in \tilde{\cD}_{ij}(q_1) \}. 
\end{displaymath}
Since $q_1 \in K_\cD$, for each $j$, for $\ell > \ell_{\cD}$, the
density of $\tilde{\cD}_{ij}(q_1)$ is at least
$1-2\kappa \delta$. 
Then, by (\ref{eq:bilip:lambda:ij}), for $\ell > \kappa\ell_{\cD}$,
the density of 
$\hat{\cD}(q_1,u)$ is at least
$(1-4|\tilde{\Lambda}|\kappa^2 \delta)$.
Now from the definitions we have
\begin{displaymath}
E_3(q_1,u) = E_2(q_1,u) \cap \{ \ell \st \tau_{(\epsilon)}(q_1,u,\ell)
\in \hat{\cD}(q_1,u) \}.  
\end{displaymath}
Therefore, in view of \autoref{prop:bilip:hattau:epsilon}, the
density of $E_3(q_1,u)$ is at least $1-c(\delta)$ where $c(\delta) \to
0$ as $\delta \to 0$.

Let $E_4(q_1,u) = \{ \ell > 0 \st g_{-\ell/2} u q_1 \in K\}$. Let
$E_5(q_1) = \{ \ell > 0 \st g_{-\ell/2} q_1 \in K \}$ and let
$E_6(q_1) = \{ \ell > 0 \st g_{-\ell} q_1 \in K \}$. Since $q_1 \in
K_\cD$ and $uq_1 \in K_\cD$ 
the density of $E_4(q_1,u) \cap E_5(q_1) \cap E_6(q_1)$ is at least
$1-c_1(\delta)$ where $c_1(\delta) \to 0$ as $\delta \to 0$. Let
$E_{good}(q_1,u) = E_3(q_1,u) \cap E_4(q_1,u) \cap E_5(q_1) \cap
E_6(q_1)$. Then, if $\ell \in E_{good}(q_1,u)$ then
(\ref{eq:q2:Yijbalanced:in:K}) and (\ref{eq:q3:Yijbalanced:in:K})
hold, and also $q_1\in K$, $g_{-\ell/2}q_1 \in K$, $g_{-\ell} q_1 \in
K$, and $g_{-\ell/2} u q_1 \in K$. Thus, $Y_{ij,balanced}(q_1,u,\ell) \in
K$. 
\end{proof}

\begin{claim}
There exists  a set $\cD_4 = \cD_4(K_{00},\delta,\epsilon,\eta) 
\subset \reals^+$ and 
a number $\ell_4 = \ell_4(K_{00},\delta,\epsilon,\eta) > 0$ so
that $\cD_4$ has density at least $1-c_4(\delta)$ 
for $\ell > \ell_4$, and for $\ell \in \cD_4$ 
a subset $K_4(\ell) = K_4(\ell,K_{00},\delta,\epsilon,\eta)
\subset X$ with $\nu(K_4(\ell)) > 1- c_4'(\delta)$,  
such that for any $q_1 \in K_4(\ell)$ there exists a subset
$Q_4(q_1,\ell) \subset Q_3(q_1, \ell) \subset \cB_0(q_1)$ 
with density at least $1-c_4''(\delta)$, so
that for all $\ell \in \cD_4$, for all $q_1 \in
K_4(\ell)$ and all $u \in Q_4(q_1,\ell)$, 
\begin{equation}
\label{eq:Yij:balanced:in:K}
Y_{ij,balanced}(q_1,u,\ell) \in K.
\end{equation}
(We have $c_4(\delta)$, $c_4'(\delta)$ and
$c_4''(\delta) \to 0$ as $\delta \to 0$). 
\end{claim}

\begin{proof}[Proof of Claim.] This follows from \autoref{claim:12.3}
by applying Fubini's theorem to
$X_\cB \cross \reals$, where $X_\cB = \{ (x,u) \st x \in X,\ u \in
\cB_0(x) \}$.   
\end{proof}

Suppose $\ell \in \cD_4$. Let $K_\ell$ be the intersection of the sets
$K_\ell$ of 
\autoref{lemma:swithching}, \autoref{cor:tau:u:bilipshitz} and
\autoref{lemma:cA:transport}. 
We now apply \autoref{prop:can:avoid:most:Mu} 
with $K'=g_{-\ell}K_4(\ell) \cap K_\ell$. We denote
the resulting set by $K_5(\ell)= K_5(\ell,K_{00},\delta,\epsilon,\eta)$.
In view of the choice of $\epsilon_1$,  
we have $\nu(K_5(\ell)) \ge 1 - c_5(\delta)$, where $c_5(\delta) \to
0$ as $\delta \to 0$.


\medskip

Let $\cD_5 = \cD_4$ and let $K_6(\ell) = g_{\ell}K_5(\ell)$.

\subsubsection{Choice of parameters \#2: Choice of $q$, $q'$, $q_1'$ (depending
  on $\delta$, $\epsilon$, $q_1$, $\ell$).} 
  \label{sssec:choice_of_parameters_nr_2}
Suppose $\ell \in \cD_5$ and $q_1 \in K_6(\ell)$.  Let 
$q = g_{-\ell} q_1$.
Then $q \in K_5(\ell)$. 
Let $\cA_2(q,u,\ell,t)$ be as in
\autoref{sssec:the_4_variable_curly_a_2}
and for $u \in Q_4(q_1,\ell)$ let \index{$M$@$\cM_u$}$\cM_u$ be the subspace of
\autoref{lemma:bad:subspace} applied to the restriction of  linear map $\cA_2(q_1,u, \ell,
\tau_{(\epsilon)}(q_1,u,\ell))$ to $V^s_*(q)$. 
By \autoref{prop:can:avoid:most:Mu} and the
definition of $K_5(\ell)$, we
can choose $q' \in g_{-\ell} K_4(\ell) \cap K_\ell$ with $q' \in \gB_0^-[q]$
and so that (\ref{eq:rho:delta:le:Fq:minus:Fqprime}) and
(\ref{eq:qprime:avoids:Mu}) hold with
  $\epsilon_1(\delta) \to 0$ as
  $\delta \to 0$.  
Let $q_1' = g_\ell q'$. Then $q_1' \in K_4(\ell)$.


\subsubsection{Standing Assumption.} We assume $\ell \in \cD_5$, $q_1 \in K_6(\ell)$
and $q$, $q'$, $q_1'$ are as in Choice of parameters \#2 (\autoref{sssec:choice_of_parameters_nr_2}). 

\subsubsection{Notation.} For $u \in \cB_0(q_1)$, $u' \in \cB_0(q_1')$, let
\begin{displaymath}
\tau(u) = \tau_{(\epsilon)}(q_1, u, \ell), \qquad \tau'(u') =
\tau_{(\epsilon)}(q_1', u', \ell).
\end{displaymath}
Then, for any $ij$,\index{$\tau(u)$}\index{$\tau'(u')$}
\begin{displaymath}
\tau(u) = \tau(Y_{ij,balanced}(q_1,u,\ell)), \qquad \tau'(u') = \tau(Y_{ij,balanced}(q_1',u',\ell)).
\end{displaymath}





\subsubsection{The numbers $t_{ij}(u)$ and $t_{ij}'(u')$.}
Suppose $u \in Q_4(q_1,\ell)$, and suppose $ij \in \tilde{\Lambda}$. 
Let \index{$t_{ij}(u)$}$t_{ij}(u)$ be defined by the equation 
\begin{equation}
\label{eq:def:tj}
\lambda_{ij}(u q_1, \tau_{(\epsilon)}(q_1,u,\ell)) =
\lambda_{ij}(q_1, t_{ij}(u)).  
\end{equation}
Alternatively,
\begin{equation}
\label{eq:alt:def:tj}
t_{ij}(u) = t_0(Y_{ij,balanced}(q_1,u,\ell)).   
\end{equation}
Then, since $\ell \in \cD_4$ and in view of
(\ref{eq:Yij:balanced:in:K}),  it
follows that 
\begin{equation}
\label{eq:gtj:q1:inK}
g_{t_{ij}(u)} q_1 \in K.
\end{equation}
Similarly, suppose $u' \in Q_4(q_1',\ell)$ and $ij \in
\tilde{\Lambda}$. 
Let \index{$t_{ij}'(u')$}$t_{ij}'(u')$ be defined by the equation
\begin{equation}
\label{eq:def:tj:prime}
\lambda_{ij}\left(u' q_1', \tau_{(\epsilon)}(q_1',u',\ell)\right) =
\lambda_{ij}(q_1', t_{ij}'). 
\end{equation}
Alternatively,
\begin{equation}
\label{eq:alt:def:tj:prime}
t_{ij}'(u') = t_0\left(Y_{ij,balanced}(q_1',u',\ell)\right).   
\end{equation}
Then, by the same argument, 
\begin{equation}
\label{eq:gtjprime:inK}
g_{t_{ij}'(u')} q_1' \in K.
\end{equation}

\subsubsection{The map $\bfv(u)$ and the generalized subspace $\cU(u)$.}
For $u \in \cB_0(q_1)$, let \index{$v$@$\bfv(u)$}
\begin{equation}
\label{eq:def:bfv:u}
\bfv(u) = \bfv(q,q',u,\ell,\tau) = \cA_2(q_1, u, \ell, \tau) F^s_q(q'),
\end{equation}
where $\tau = \tau_{(\epsilon)}(q_1,u,\ell)$, $F^s_q$ is as in
\autoref{theorem:3varA:haltime:factorizable}
and $\cA_2(\cdot, \cdot, \cdot, \cdot)$ is as in \autoref{sec:subsec:4var:curlyA}.

Let \index{$U$@$\cU(u)$}$\cU(u) = \phi_{\tau(u)}^{-1}(U^+[g_{\tau}q_1']) \subset
\cC[g_{\tau(u)} u q_1]$. Then, by
\autoref{prop:factorization:main},
\begin{displaymath}
\bfj(\cU(u)) = \bfv(u) + O(e^{-\alpha \ell}). 
\end{displaymath}


\subsubsection{Standing Assumption.}
\label{sec:subsubsec:standing:assumption:epsilons}
Let $\epsilon_0(\delta) =
\epsilon_0(\delta,\epsilon_0'(\delta)) < \epsilon_0'(\delta)$ be as in
\autoref{lemma:parts:staying:close}. We assume that $\epsilon < 0.01
\epsilon_0(\delta)$ and furthermore
$C(\delta) \epsilon <
0.01 \epsilon_0(\delta)$ for
any constant $C(\delta)$ arising in the course of the proof. In
particular, this applies to $C_2(\delta)$ and $C_2'(\delta)$ in the
next claim.

\begin{claim} There exists a subset $Q_5 =
  Q_5(q_1,\ell, K_{00},\delta,\epsilon,\eta) \subset Q_4(q_1,\ell)$ 
with $|Q_5| \ge (1-c_5''(\delta))|\cB_0(q_1)|$  
(with $c_5''(\delta) \to 0$ as
$\delta \to 0$), and a number
$\ell_5 = \ell_5(\delta,\epsilon)$ such that for all $u \in Q_5$
and $\ell > \ell_5$, 
\begin{equation}
\label{eq:tau:u:alpha3:ell}
\tau(u) < \frac{1}{2} \alpha_3 \ell,
\end{equation}
where $\alpha_3 > 0$ is as in \S\ref{sec:subsec:choice:of:alpha3}. 
In addition,
\begin{equation}
\label{eq:vu:has:norm:near:epsilon}
e^{-\kappa T_1(\delta)/4} \epsilon \le \|\bfv(u)\| \le e^{\kappa
  T_1(\delta)/4} \epsilon, 
\end{equation}
where $T_1(\delta)$ is as in (\ref{eq:def:T1:delta}). 
\begin{equation}
\label{eq:hdloc:Uplus:phit:inverse:Uplusprime}
C_1'(\delta) \epsilon  \le d^u_{\cH,loc}(U^+[g_{\tau(u)} u q_1), \phi_{\tau(u)}^{-1}(U^+[q_1'])) \le
C_1(\delta) \epsilon.
\end{equation}
Furthermore, given $u' \in U^+[q_1']$ set $z_{1/2}:=\cW^{cs}[g_{-\ell/2}uq_1]\cap \cW^u_{loc}[g_{-\ell/2}q_1']$ (see \autoref{def:z_points_from_y_configurations}).
If $u'$ is such that
\begin{equation}
\label{eqn:left_shadowing_condition_main_proof}
	g_{\ell/2 + \tau(u)} z_{1/2} \in \cW^{u}_{loc}[g_{\tau(u)}u' q_1']
\end{equation}
then $u' \in \cB_0(q_1')$. 
\end{claim}

\begin{proof}[Proof of claim.] Let $\cM_u$ be the subspace of
\autoref{lemma:bad:subspace} applied to the linear map $\cA_2(q_1,u, \ell,
\tau_{(\epsilon)}(q_1,u,\ell))$, where $\cA_2$ is as in
\autoref{sec:subsec:4var:curlyA}. Let $Q(q_1)$ be as in
(\ref{eq:continue:indction:quantified}), so $|Q(q_1)| \ge
(1-\delta)|\cB_0(q_1)|$.  
Let $Q_5' \subset Q_4 \cap Q(q_1)$ be such that for all $u \in
Q_5'$, 
\begin{displaymath}
d_{V^s(q)}(F^s_q(q'),\cM_u) \ge \beta(\delta)
\end{displaymath}
where $F^s_q(\cdot)$ is as in \autoref{theorem:3varA:haltime:factorizable}. Then,
(\ref{eq:tau:u:alpha3:ell}) follows from
(\ref{eq:continue:indction:quantified}),
\autoref{lemma:lower:bound:on:3cA}, the choice of $\alpha_3$ in
\S\ref{sec:subsec:choice:of:alpha3}  
and the fact that
$Q_5 \subset Q_1$. 
Also, by (\ref{eq:qprime:avoids:Mu}),
\begin{displaymath}
|Q_5'| \ge |Q_4| - (\delta + \epsilon_1(\delta))|\cB_0(q_1)| \ge (1-
\delta- \epsilon_1(\delta) - c_4''(\delta))|\cB_0(q_1)|.
\end{displaymath}
Then, let 
$Q_5 = \{ u \in Q_5' \st d(u, \partial \cB_0(q_1)) > \delta \}$, hence 
$$|Q_5| \ge
(1-c_5'(\delta) - c_4'(\delta)-c_n \delta)|\cB_0(q_1)|,$$ 
where $c_n$ depends only on the dimension.

We have $C(\delta)^{-1} \epsilon \le \|\cA_2(q_1,u,\ell,t)\| \le C(\delta)
\epsilon$ by the definition of $t =
\tau_{(\epsilon)}(q_1,u,\ell)$. 
We now apply \autoref{lemma:bad:subspace} to the linear map
$\cA_2(q_1,u,\ell,t)$. Then, for all $u \in Q_5$, 
\begin{displaymath}
c(\delta) \|\cA_2(q_1,u,\ell,t)\| \le \|\cA_2(q_1,u,\ell, t) F_q(q')\| \le \|\cA_2(q_1,u,\ell,t)\|. 
\end{displaymath}
Therefore, 
\begin{displaymath}
C'(\delta)^{-1} \epsilon \le 
\|\cA_2(q_1,u,\ell, t) F_q(q')\| \le C'(\delta) \epsilon
\end{displaymath}
This immediately implies (\ref{eq:vu:has:norm:near:epsilon}), in view
of the definition of $\bfv(u)$.  Then,
(\ref{eq:hdloc:Uplus:phit:inverse:Uplusprime}) follows from
\autoref{prop:hausdorff_distance_and_norm_of_vector}. 

Finally, suppose $u \in Q_5$, and $u' \in U^+(q_1')$ is such that
\autoref{eqn:left_shadowing_condition_main_proof} holds.
Let $z_t := g_{\ell/2 + t}z_{1/2}$.
Recall that the interpolation map satisfies $\phi\left(y_{1/2}\right)=z_{1/2}$ where $y_{1/2}:=g_{-\ell/2}uq_1$; therefore $\phi_{t}^{-1}\left(z_{t}\right) = g_{t} u q_1$ for $t\in [0,\tau(u)]$.

The assumption in \autoref{eqn:left_shadowing_condition_main_proof} states that
\[
	d^u\left(z_{\tau(u)},g_{\tau(u)}u'q_1'\right)\leq 1 
	\text{ therefore }
	d^u\left(z_0,u'q_1'\right)\leq e^{-\kappa \cdot \tau(u)}.
\]
Since $q_1'\in K$ (and $z_0=\phi_0(uq_1)$) this implies that there exists a $C_1(\delta)>0$ such that
\[
	d^Q\left(\phi_0(uq_1),u'q_1'\right)\leq C_1(\delta) e^{-\kappa \cdot \tau(u)}.
\]
Now \autoref{theorem:interpolation:does:not:move:points} implies that
\[
	d^Q\left(uq_1,\phi_0^{-1}(u'q_1')\right)\leq C_2(\delta) e^{-\kappa \cdot \tau(u)}
\]
for some $C_2(\delta)>0$.
Since $uq_1\in K$ we have
\[
	d^u\left(uq_1,\phi_0^{-1}(u'q_1')\right)\leq C_3(\delta) e^{-\kappa \cdot \tau(u)}.
\]



Again since $uq_1 \in K$, we have that $u q_1 \in
\cB_0^{(2\epsilon')}[q_1]$. Therefore, for $\ell$ sufficiently large,
\begin{displaymath}
\phi_0^{-1}(u' q_1') \in \cB_0^{(\epsilon')}[q_1].
\end{displaymath}
However, by \autoref{theorem:3varA:haltime:factorizable},
\begin{displaymath}
d^Q(\phi_0^{-1}(u' q_1'), u' q_1') = O(e^{-\alpha \ell}), 
\end{displaymath}
and since $K$ was chosen so that the map $q \to \gB_0[q]$ is
uniformly continuous with respect to the Hausdorff topology, $q_1 \in
K$, and $q_1' \in K$, 
we have, for $\ell$ sufficiently large,
\begin{displaymath}
d_{\cH}^Q(\gB_0[q_1], \gB_0[q_1']) \le 0.1 \epsilon'. 
\end{displaymath}
Therefore, $u' q_1' \in \gB_0[q_1']$, which implies $u' \in
\cB_0(q_1')$. 
\end{proof}

\subsubsection{Standing Assumption.} We assume $\ell > \ell_5$.

\subsubsection{Recap}
Under our standing assumptions, for $u \in Q_5(q_1,\ell)$, 
the $Y$-configuration $Y=Y_{ij,balanced}(q_1,u,\ell)$ satisfies $Y \in
K$. If $u' \in U^+(q_1')$ is such that
(\ref{eqn:left_shadowing_condition_main_proof}) holds, then the $Y$-configuration
$Y'= Y(q_1',u',\ell,\tau(u))$ is bottom-linked to $Y$ and left-shadows
$Y$. Also, (\ref{eq:vu:has:norm:near:epsilon}) and (\ref{eq:hdloc:Uplus:phit:inverse:Uplusprime}) hold. 

However, we do not know that $Y' \in K$ (we only know that
$q' \in K$, $q_1' \in K$ and $q_{1/2}' \in K$). 
The next two claims aim to deal with this issue. 

\begin{claim} Suppose $u \in Q_5(q_1,\ell)$, $u' \in
  Q_4(q_1',\ell)$ and (\ref{eqn:left_shadowing_condition_main_proof}) holds.
 Then, there exists $C_0 = C_0(\delta)$ such that
\begin{equation}
\label{eq:tau:epsilon:same:as:tau:epsilon:prime}
|\tau_{(\epsilon)}(q_1, u, \ell) -
\tau_{(\epsilon)}(q_1',u',\ell)| \le C_0(\delta). 
\end{equation}
\end{claim}

\begin{proof}[Proof of claim.] Let $\tau = \tau_{(\epsilon)}(q_1, u,
\ell)$, $\tau' = \tau_{(\epsilon)}(q_1',u',\ell)$. In view of
(\ref{eq:Yij:balanced:in:K}), 
we have $g_\tau u q_1 \in K$, $g_{\tau'} u' q_1' \in K$. 
Suppose $\tau \le
\tau'$. 
By the definition of $K$, there exists $s<\tau$ with $(\tau-s) < (1/2)\kappa^{-2}
(\tau'-\tau)+T_0(\delta)$ such that $g_s u q_1 \in K_0$, $g_s u' q_1' \in
K_0$. By \autoref{lemma:A:future:bilipshitz} and \autoref{lemma:swithching}, we have (as long as $\ell$ is
sufficiently large compared with $\epsilon$),
\begin{align*}
\epsilon = \|\cA_2(q_1',\ell, u', \tau')\| & \ge e^{\kappa^{-1}(\tau'-s)}
\|\cA_2(q_1',\ell, u', s)\| \\
& \ge e^{\kappa^{-1}(\tau'-s)} C(\delta)^{-1} \|\cA_2(q_1,\ell, u, s)\| \\
& \ge e^{\kappa^{-1}(\tau'-s)} C(\delta)^{-1} e^{-\kappa(\tau-s)}
  \|\cA_2(q_1,\ell, u, \tau)\| = \epsilon    
\end{align*}
Thus,
\begin{displaymath}
1 \ge e^{\kappa^{-1}(\tau'-s)} C(\delta)^{-1} e^{-\kappa(\tau-s)}
\end{displaymath}
hence 
\begin{displaymath}
\kappa^{-1}(\tau'-\tau) \le \kappa^{-1}(\tau'-s) \le \kappa(\tau-s) + C(\delta)
\le (1/2) \kappa^{-1}(\tau'-\tau) + C(\delta). 
\end{displaymath}
Thus, $(\tau'-\tau) < 2 \kappa C(\delta)$ as required.

The proof in the case $\tau' < \tau$ is identical.
\end{proof}



\begin{claim} There exists a subset $Q_6(q_1,\ell) =
  Q_6(q_1,\ell,K_{00},\delta,\epsilon,\eta) \subset Q_5(q_1,\ell)$ with 
$|Q_6(q_1,\ell)| > (1-c_6'(\delta)) |\cB_0(q_1)|$ 
and with $c_6'(\delta) \to 0$ as
$\delta \to 0$ and a number $\ell_6 = \ell_6(\epsilon,\delta)$ 
such that for all $\ell > \ell_6$ for all $u \in Q_6(q_1,\ell)$ there exists
$u' \in Q_4(q_1',\ell)$  such that
\begin{equation}
\label{eq:alt:d:psi:u:g:tau:u:uprime}
d^Q\left(g_{\tau(u)} u q_1, g_{\tau(u)} u' q_1'\right) < C(\delta) \epsilon.
\end{equation}
\end{claim}

\subsubsection{Remark.} The claim implies that for all $u \in
Q_6(q_1,\ell)$ there exists $u' \in \cB_0(q_1')$ such that for
any $ij$, if we write
$Y = Y_{ij,balanced}(q_1,u,\ell)$ and $Y' =
Y_{ij,balanced}(q_1',u',\ell)$ then we have $Y \in K$, $Y' \in K$ and
also that $Y$ and $Y'$ are bottom-linked and also $Y'$ shadows $Y$ up to the
time $\tau(Y)$.

\begin{proof}[Proof of Claim.]
Note that the sets $$\{ \cB_{\tau(u)}[u q_1] \st u \in Q_5(q_1,\ell) \}$$
are a cover of $Q_5(q_1,\ell)q_1$.
Then, since these sets satisfy the condition of
\autoref{lemma:gB:properties} (b), we can
find a pairwise disjoint subcover, i.e. 
find $u_j \in Q_5(q_1,\ell)$, $1 \le j \le N$,
with 
$Q_5(q_1,\ell) q_1 \subset \bigsqcup_{j=1}^N \cB_{\tau(u_j)}[u_j q_1]$; in particular, 
$\cB_{\tau(u_j)}[u_j q_1]$ and $\cB_{\tau(u_k)}[u_k q_1]$ are disjoint for
$j \ne k$. Let
\begin{displaymath}
P_j = g_{-\tau(u_j)} \cB^{(\epsilon_0')}_0[g_{\tau(u_j)} u_j q_1] \subset
\cB_{\tau(u_j)}[u_j q_1] \subset U^+[q_1]. 
\end{displaymath}
Then $P_j$ and $P_k$ are disjoint for $j \ne k$.

Using notation (\ref{eq:phitau}), let $\phi^j_s=\phi_s(Y(q_1, u_j,\ell),q_1')$, that is  $\phi^j_s$ denote the interpolation operator $\phi_s$ where $u$ is replaced by
$u_j$. Let 
\begin{multline*}
P_j' = \phi_0^{j} \left(g_{-\tau(u_j)} \gB^{(\epsilon_0')}_0[g_{\tau(u_j)}
u_j q_1]\right)  \cap U^+\left[q_1'\right] = \\
g_{-\tau(u_j)} \left(\phi^j_{\tau(u_j)}\left(\gB^{(\epsilon_0')}_0\left[g_{\tau(u_j)}
u_j q_1\right] \right) \cap U^+\left[g_{\tau(u_j)} q_1'\right]\right) \subset U^+\left[q_1'\right]. 
\end{multline*}
In view of \autoref{theorem:3varA:haltime:factorizable},
\begin{displaymath}
d_{\cH}^Q\left(P_j,P_j'\right) = O\left(e^{-\alpha \ell}\right).
\end{displaymath}
Also, since $K$ was chosen so that the maps $q \mapsto \gB_0[q], q\mapsto \cB_{0}[q]$ is
uniformly continuous with respect to the Hausdorff topology,
we have, for $\ell$ sufficiently large,
\begin{align}
	\label{eqn:uniform_continuity_gB_cB}
	d_{\cH}^Q(\gB_0[q_1], \gB_0[q_1']) & \le 0.1 \epsilon'(\delta)\\
	d_{\cH}^Q(\cB_0[q_1], \cB_0[q_1']) & \le 0.1 \epsilon'(\delta)
\end{align}
Thus, since $P_j \subset \gB_0^{(\epsilon')}[q_1]$, we have $P_j' \subset
\gB_0[q_1']$. Since by construction, $P_j' \subset U^+[q_1']$, we have
$P_j' \subset \cB_0[q_1']$.

Note that by (\ref{eq:hdloc:Uplus:phit:inverse:Uplusprime}) and
\autoref{theorem:interpolation:does:not:move:points}
for any $u \in P_j$, $u' \in P_j'$,
\begin{equation}
\label{eq:PjPjprime:staying:close}
d^Q(g_{\tau(u_j)} u q_1, g_{\tau(u_j)} u' q_1') < C(\delta)
\epsilon + O(e^{-\alpha \ell}). 
\end{equation}
We now claim that $P_j'$ and $P_k'$ are disjoint for $j \ne k$. 

Suppose $u' \in P_j' \cap P_k'$. Without loss of generality, we may assume $\tau(u_j) \le \tau(u_k)$.
Then $w \equiv g_{\tau(u_j)} u' q_1' =
\phi^j_{\tau(u_j)} (x_j) = \phi^k_{\tau(u_j)}(x_k)$, where
$x_j \in {\gB}^{(\epsilon_0')}_0[g_{\tau(u_j)} u_j q_1]$, $x_k \in
{\gB}^{(\epsilon_0')}_0[g_{\tau(u_j)} u_k q_1]$.
By \autoref{theorem:interpolation:does:not:move:points},
for all $0 \le t \le \tau(u_j)$,
\begin{displaymath}
d^Q(g_{-t} x_j, g_{-t} x_k) = O(e^{-\alpha \ell}). 
\end{displaymath}
Let $y_j = g_{\tau(u_j)} u_j q_1$, $y_k = g_{\tau(u_j)} u_k q_1$. 




By \autoref{lem:lemma_stay_in_chart_stay_in_loc_new} applied with $y=y_j, x_1=x_j, x_2=x_k$, we have:
\begin{displaymath}
d^u(x_j,x_k) \le C(\delta) e^{-\alpha
  \ell}. 
\end{displaymath}

For $\ell$ sufficiently large  this implies
that $d^u(x_j,x_k) < \epsilon_0'$. However, since
$x_j \in \gB_0^{(\epsilon_0')}[y_j]$, this forces
$x_k \in \gB_0[y_j]$. But, by construction,
$x_k \in \gB_0[y_k]$, hence $\gB_0[y_j] =
\gB_0[y_k]$.  Now, since $\tau(u_j) \le \tau(u_k)$, 
\begin{multline*}
\gB_{\tau(u_k)}[u_kq_1]  
\subseteq \gB_{\tau(u_j)}[u_k q_1] = \\ = g_{-\tau(u_j)} \gB_0[y_k]
= g_{-\tau(u_j)} \gB_0[y_j] = \gB_{\tau(u_j)}[u_j q_1], 
\end{multline*}
which is a contradiction for $j \ne k$. 
This proves the required
disjointness. 

By \autoref{lemma:parts:staying:close} 
provided $P_j \cap Q_5(q_1, \ell) \ne \emptyset$, we have
$\kappa^{-1}|P_{j,m}| \le |P_{j,m}'| \le
\kappa|P_j|$, where $\kappa$ depends only on the Lyapunov
spectrum, and we have normalized the measures $| \cdot |$ so that
$|\cB_0(q_1)| = |\cB_0(q_1')|=1$. (Note that
\autoref{lemma:parts:staying:close} 
applies since in view of \autoref{eq:def:epsilon0:prime}, 
$x \in K$ implies $x \in \gB_0^{(2\epsilon_0')}[x]$, and because \autoref{eqn:uniform_continuity_gB_cB} holds). 
We now apply \autoref{lemma:segments:equal:length}
to obtain the set $Q_6$. 

It remains to replace $\tau(u_j)$ by $\tau(u)$ in
(\ref{eq:PjPjprime:staying:close}). In view of
\autoref{cor:tau:u:bilipshitz}, 
\begin{equation}
\label{eq:tmp:tau:u:tau:uj}
|\tau(u) - \tau(u_j)| \le C_1(\delta).
\end{equation}
Then, provided $\epsilon$ is small enough depending on $\delta$ and
$\ell$ is large enough depending on $\epsilon$ and $\delta$, 
(\ref{eq:alt:d:psi:u:g:tau:u:uprime})  follows from
(\ref{eq:PjPjprime:staying:close}), 
(\ref{eq:tmp:tau:u:tau:uj}), as well as \autoref{eqn:elementary_compactness_2_point_estimate}.
\end{proof}

\begin{claim}
\label{claim:12:12}
There exist a number $\ell_8 =
  \ell_8(K_{00},\delta,\epsilon,\eta)$ and a constant $c_8(\delta)$ with
$c_8(\delta) \to 0$ as $\delta \to 0$ and for every $\ell > \ell_8$ 
a subset $Q_8(q_1,\ell) = Q_8(q_1,\ell,K_{00},\delta,\epsilon,\eta) 
\subset \cB(q_1)$ with $|Q_8(q_1,\ell)| \ge
(1-c_8(\delta)) |\cB(q_1)|$ so that 
for $u \in Q_8(q_1,\ell)$  we have
\begin{equation}
\label{eq:Ru:R2prime:u:close:Eplus}
d\left(\frac{\bfv(u)}{\|\bfv(u)\|}, \bbE(g_{\tau(u)} u q_1)\right) \le
C_8(\delta) e^{-\alpha' \ell},  
\end{equation}
where $\bfv(u)$ is defined in (\ref{eq:def:bfv:u})
and $\alpha'$ depends only on the Lyapunov spectrum.  
\end{claim}

\begin{proof}[Proof of claim.] 
Let $L' > L_2(\delta)$ be a constant to be chosen later, where
$L_2(\delta)$ is as in \autoref{prop:most:inert}. Also let
$\ell_8 = \ell_8(\delta,\epsilon,K_{00},\eta)$ be a constant to be
chosen later. Suppose $\ell > \ell_8$, 
and suppose $u \in Q_7^*(q_1,\ell)$, so in particular $g_{\tau(u)} u
q_1 \in K$. Let $t \in [L',2L']$ 
be such that \autoref{prop:most:inert} holds for $\bfv =
\bfv(u)$ and $x = g_{\tau(u)} u q_1$. 

Let  
$B_u \subset \cB(q_1)$ denote
$\cB_{\tau_{(\epsilon)}(q_1,u,\ell)-t}(u q_1)u$,
(where $\cB_t(x)$ is defined in
\autoref{sssec:measurable_partitions_and_subgroups_compatible_with_the_measure}). 
Suppose $u_1 \in B_u \cap Q_7(q_1,\ell)$, and write
\begin{displaymath}
g_{\tau(u_1)} u_1 q_1 = g_s u_2 g_t^{-1} g_{\tau(u)} u q_1. 
\end{displaymath}
Then, $u_2 \in \cB(g_t^{-1} g_{\tau(u)} u q_1)$ and $L' \le t \le  2L'$.
We have $s = \tau(u_1) - (\tau(u) - t)$. By
\autoref{cor:tau:u:bilipshitz},
\begin{multline*}
(1/4) \kappa^{-2} L - C(\delta) - O(e^{-\alpha' \ell}) < s < 4 \kappa^2 L + C(\delta) + O(e^{-\alpha' \ell}).
\end{multline*}
In the rest of this proof, we will use the notation $(g_s u_2
g_t^{-1})_*$ to denote the map from $\bH(g_{\tau(u)} u q_1) \to
\bH(g_{\tau(u_1)} u_1 q_1)$ given by
\begin{displaymath}
(g_s u_2 g_t^{-1})_* = g_s \circ u_*(g_{\tau(u)-t} uq_1, g_{\tau(u)-t}
u_1q_1) \circ g_{-t}.
\end{displaymath}
By \autoref{lemma:cA:transport}, 
\begin{equation}
\label{eq:tracking:one}
\| (g_s u_2 g_t^{-1})_* \bfv(u) - \bfv(u_1) \| = O(e^{\kappa' L'}
e^{-\alpha \ell}). 
\end{equation}
In view of
(\ref{eq:vu:has:norm:near:epsilon}), $\|\bfv(u_1) \| \approx
\epsilon$. Thus, $\|(g_s u_2 g_t^{-1})_* \bfv(u)\| \approx
  \epsilon$, and
\begin{displaymath}
\left\| \frac{(g_s u_2 g_t^{-1})_* \bfv(u)}
{\|(g_s u_2 g_t^{-1})_* \bfv(u)\|} - \frac{\bfv(u_1)}{\|\bfv(u_1)\|}
\right\| = O_\epsilon( e^{\kappa'L'-\alpha\ell}). 
\end{displaymath}
But, by \autoref{prop:most:inert}, 
for $1-\delta$ fraction
of $u_2 \in \cB(g_t^{-1}g_{\tau(u)} u q_1)$, 
\begin{displaymath}
d\left( \frac{(g_s u_2 g_{-t})_* \bfv(u)}{\|(g_s u_2 
    g_{-t})_* \bfv(u)\|} , 
     \bbE(g_{\tau(u_1)} u_1 q_1) \right) \le C(\delta) e^{-\alpha L'}, 
\end{displaymath}
Note that
\begin{displaymath}
\cB(g_t^{-1}g_{\tau(u)} u q_1) = g_{\tau_{(\epsilon)}(q_1, u, \ell) - t}
B_u. 
\end{displaymath}
Therefore, for $1-\delta$ fraction of $u_1 \in B_u$, 
\begin{equation}
\label{eq:tmp:claim:cover}
d\left(\frac{\bfv(u_1)}{\|\bfv(u_1)\|}, \bbE(g_{\tau(u_1)} u_1 q_1)\right) 
\le C(\epsilon,\delta) [e^{\kappa'L'-\alpha \ell} + e^{-\alpha L'}]
\end{equation}
We can now choose $L' > 0$ to be $\alpha' \ell$ where $\alpha' > 0$ is
a small constant depending only on the Lyapunov  spectrum, 
and $\ell_8 > 0$ so that for $\ell > \ell_8$ 
the right-hand-side of the above equation is at most $e^{-\alpha'
  \ell}$. 

The collection of balls $\{B_u\}_{u \in Q_7^*(q_1,\ell)}$ 
are a cover of $Q_7^*(q_1,\ell)$. These balls satisfy
the condition of \autoref{lemma:gB:properties} (b); hence we may
choose a pairwise disjoint subcollection which still covers 
$Q_7^*(q_1,\ell)$. Then, by summing (\ref{eq:tmp:claim:cover}), we see
that 
(\ref{eq:Ru:R2prime:u:close:Eplus}) holds for $u$ in a subset
$Q_8 \subset \cB_0[q_1]$ of measure at least $(1-c_8(\delta))|\cB_0[q_1]| =
(1-\delta)(1-c_7^*(\delta))|\cB_0[q_1]| $. 
\end{proof}

\begin{claim}
There exists a subset $Q_8^*(q_1,\ell) =
Q_8^*(q_1,\ell,K_{00},\delta,\epsilon,\eta) \subset 
Q_8(q_1,\ell)$ with $|Q_8^*| \ge (1-c_8^*(\delta)) |\cB(q_1)|$
such that for $u \in Q_8^*$ and any
$t > \ell_8(\delta)$
we have
\begin{equation}
\label{eq:Bt:Q8:density}
  |\cB_t(u q_1)u \cap Q_8(q_1,\ell)| \ge
  (1-c_8^*(\delta)) |\cB_t(u q_1)u|, 
\end{equation}
where $c_8^*(\delta) \to 0$ as $\delta \to 0$. 
\end{claim}
\begin{proof}[Proof of Claim.] This follows immediately from
\autoref{lemma:cB:vitali:substitute}. 
\end{proof}

\subsubsection{Choice of parameters \#3: Choice of $\delta$.}
Let $\theta' = (\theta/2)^n$, where $\theta$
and $n$ are as in \autoref{prop:some:fraction:bounded}. We can
choose $\delta > 0$ so that 
\begin{equation}
\label{eq:choice:of:delta}
c_8^*(\delta) < \theta'/2. 
\end{equation}

\begin{claim}
\label{claim:large_set_of_u}
	There exist sets $Q_9(q_1,\ell) = Q_9(q_1,\ell,K_{00},
  \delta,\epsilon, \eta) \subset Q_8^*(q_1,\ell)$ 
with $|Q_9(q_1,\ell)| \ge (\theta'/2)(1-\theta'/2) |\cB(q_1)|$ and $\ell_9 =
\ell_9(K_{00}, \delta,\epsilon, \eta)$, 
such that for $\ell > \ell_9$ and $u \in Q_9(q_1,\ell)$, 
\begin{equation}
\label{eq:psi:psiprime:close:to:bigcup}
d\left(\frac{\bfv(u)}{\|\bfv(u)\|},
  \bigcup_{ij \in \tilde{\Lambda}}   \bbE_{[ij],bdd}(g_{\tau(u)} u q_1)\right) < C(\delta) \eta. 
\end{equation}
\end{claim}

\begin{proof}[Proof of claim.]  
Suppose $u \in Q_8^*(q_1,\ell)$. Then, by
(\ref{eq:Ru:R2prime:u:close:Eplus}) and (\ref{eq:vu:has:norm:near:epsilon}),
we may write
\begin{displaymath}
\bfv(u) = \bfv'(u) + \bfv''(u),
\end{displaymath}
where $\bfv'(u) \in \bbE(g_{\tau(u)} u q_1)$ and $\|\bfv''(u) \| \le
C(\delta) e^{-\alpha' \ell}$.

Let $L= L_0(\delta,\eta)$ where $L_0(\delta,\eta)$ is as in
\autoref{prop:some:fraction:bounded}.
Suppose $u_1 \in
\cB_{\tau(u)-L}(uq_1)u \cap Q_8(q_1,\ell)$. 
Let $u_2$ be such that $u_1 = u_2 u$. Then,
by \autoref{lemma:cA:transport}, 
\begin{displaymath}
\bfv(u_1) = (g_{\tau(u_1)} \circ (u_2)_* \circ g_{-\tau(u)}) \bfv'(u) +
  O(e^{-\alpha' \ell + \kappa^3 L})
\end{displaymath}
Since $\bfv' \in \bbE(g_{\tau(u)} u q_1)$ and in view of
\autoref{prop:agreement:of:u:star:and:p}, we can rewrite the
above equation as
\begin{equation}
\label{eq:bfv:u1:R:g:tau}
\bfv(u_1) = \bbR(g_{\tau(u_1)} u_1 q_1, g_{\tau(u)} u q_1) \bfv'(u) +
  O(e^{-\alpha' \ell + \kappa^3 L}).
\end{equation}
In view of (\ref{eq:vu:has:norm:near:epsilon}), we have, for
sufficiently large $\ell$, 
\begin{displaymath}
e^{-\kappa T_1(\delta)/2} \epsilon \le \|\bfv'(u)\| \le e^{\kappa
  T_1(\delta)/2} \epsilon,   
\end{displaymath}
and
\begin{displaymath}
e^{-\kappa T_1(\delta)/4} \epsilon \le \|\bfv(u_1)\| \le e^{\kappa
  T_1(\delta)/4} \epsilon.    
\end{displaymath}
Let $s > 0$ be such that
$y \equiv g_s u_1 q_1 \in \cF_{\bfv'(u)}[g_{\tau(u)} u q_1,L]$. From
the definition of $\cF_{\bfv'(u)}$ we have
$\|\bbR(g_{\tau(u)} u q_1,y) \bfv'(u) \| = \|\bfv'(u)\|$. Thus,
\begin{equation}
\label{eq:s:tau:u1}
|s-\tau(u_1)| \le T_1(\delta),
\end{equation}
and hence $y \in g_{[-T_1(\delta),T_1(\delta)]}
K$.
Then, 
by \autoref{prop:some:fraction:bounded} applied with
$T=T_1(\delta)$, $L =
L_0(\delta,\eta)$ and $\bfv = \bfv'(u)$, we get that for a
at least $\theta'$-fraction of $y \in
\cF_{\bfv'(u)}[g_{\tau(u)} u q_1,L]$, 
\begin{displaymath}
d\left(\frac{\bbR(g_{\tau(u)} u q_1,y) \bfv'(u)}{\|\bbR(g_{\tau(u)} u q_1,y) \bfv'(u)\|}, 
\bigcup_{ij \in\tilde{\Lambda}} 
  \bbE_{[ij],bdd}(y)\right) < 2 \eta. 
\end{displaymath}
From (\ref{eq:bfv:u1:R:g:tau}) and the definition of $y$ we have
\begin{displaymath}
\bfv(u_1) = g_{\tau(u_1)-s} \bbR(g_{\tau(u)} u q_1,y) \bfv'(u) + O(e^{-\alpha' \ell + \kappa^3 L}).
\end{displaymath}
Therefore, in view of (\ref{eq:s:tau:u1}), 
there a subset $Q_{tmp}$ of $\cB_{\tau(u)-L}(uq_1)u$ with $|Q_{tmp}|
\ge \theta' |\cB_{\tau(u)-L}(uq_1)u|$ such that for $u_1 \in Q_{tmp} \cap Q_8(q_1,\ell)$,
\begin{displaymath}
d\left(\frac{\bfv(u_1)}{\|\bfv(u_1)\|}, 
\bigcup_{ij \in\tilde{\Lambda}} 
  \bbE_{[ij],bdd}(g_{\tau(u_1)} u_1 q_1)\right) < C_2(\delta) [2 \eta +
\epsilon^{-1} e^{2\kappa^3 L}e^{-\alpha'\ell} ].
\end{displaymath}
In view of (\ref{eq:Bt:Q8:density}) and
(\ref{eq:choice:of:delta}), the above equation holds for at
least $\theta'/2$-fraction of $u_1 \in \cB_{\tau(u)-L}(uq_1)u$. 
We may choose $\ell_9 = \ell_9(K_{00},\epsilon,\delta,\eta)$ so that
for $\ell > \ell_9$ the right-hand side of the above equation is at
most $C(\delta)\eta$. Thus, (\ref{eq:psi:psiprime:close:to:bigcup}) holds for
$(\theta'/2)$-fraction of $u_1 \in \cB_{\tau(u)-L}(uq_1)u$. 

The collection of balls $\{\cB_{\tau(u)-L}(uq_1)u\}_{u \in Q_8^*(q_1,\ell)}$ 
are a cover of $Q_8^*(q_1,\ell)$. These balls satisfy
the condition of \autoref{lemma:gB:properties} (b); hence we may
choose a pairwise disjoint subcollection which still covers 
$Q_8^*(q_1,\ell)$. Then, by summing over the disjoint subcollection, we see
that the claim holds on a set $E$ of measure at least
$(\theta'/2)|Q_8^*| \ge (\theta'/2)(1-c_8^*(\delta)) \ge
(\theta'/2)(1-\theta'/2)$. 
\end{proof}


\subsubsection{Choice of parameters \#4: Choosing $\ell$, $q_1$, $q$, $q'$, $q_1'$.} 
Choose $\ell > \ell_9(K_{00},\epsilon,\delta,\eta)$. 
Now choose $q_1 \in K_7(\ell)$, and let $q$, $q'$, $q_1'$ be as in Choice of
Parameters \#2.

\subsubsection{Choice of parameters \#5: Choosing $u$, $u'$, $q_2$, $q_2'$,
  $ij$, $q_{3,ij}$,
  $q_{3,ij}'$ (depending on $q_1$, $q_1'$, $u$, $\ell$).}
	\label{sssec:choice_of_parameters_5}
Choose $u \in Q_9(q_1,\ell)$, $u' \in
Q_4(q_1',\ell)$  so that (\ref{eq:vu:has:norm:near:epsilon}), 
(\ref{eq:hdloc:Uplus:phit:inverse:Uplusprime}) and
(\ref{eq:alt:d:psi:u:g:tau:u:uprime}) hold.
Let $Y_{bal}=Y_{ij,balanced}(q_1,u,\ell)$, so that $\tau(Y_{bal}) = \tau(u)$. Let
$Y_{bal}'=Y_{ij,balanced}(q_1',u',\ell)$, so that $\tau(Y'_*) = \tau(u')$. 
We have $Y_{bal} \in K$, $Y'_* \in K$.
By (\ref{eq:tau:epsilon:same:as:tau:epsilon:prime}), 
\begin{align}
	\label{eqn:tau_parameters_can_be_chosen_close}
|\tau_{(\epsilon)}(q_1, u, \ell) -
\tau_{(\epsilon)}(q_1',u',\ell)| \le C_0(\delta), 
\end{align}
therefore,
\begin{displaymath}
g_{\tau(u)} u' q_1' \in g_{[-C,C]} K,
\end{displaymath}
where $C = C(\delta)$.

By the definition of $K$ we can
find $C_4(\delta)$ and 
$s \in [0, C_4(\delta)]$ such that 
\begin{displaymath}
q_2 \equiv g_s g_{\tau(u)} u q_1 \in K_0, \qquad q_2' \equiv g_s g_{\tau(u)} u' q_1' \in K_0. 
\end{displaymath}
Let $Y_{sh}$ be the $Y$-configuration with the same points as $Y_{bal}$,
except that $q_2(Y_{sh}) = q_2$. Similarly, let $Y_{sh}'$ be the $Y$-configuration with the same points as $Y_{bal}$,
except that $q_2(Y_{sh}') = q_2'$. Then, $Y_{sh} \in K_0$, $Y_{sh}' \in
K_0$. 

In view of  (\ref{eq:vu:has:norm:near:epsilon}),
the fact that $s \in
[0,C_4(\delta)]$, \autoref{prop:good_norms}\autoref{good_norms_bilipschitz}. 
and \autoref{lemma:hausdorff:distance:to:norm}, we
get 
\begin{equation}
\label{eq:prelimit:hd:close}
\frac{1}{C(\delta)} \epsilon \le d_{\cH,loc}^{q_2}(U^+[q_2],
\phi_{\tau(u)+s}^{-1}(U^+[q_2']))  \le C(\delta) \epsilon. 
\end{equation}
Additionally, we have the inequality in \autoref{eq:alt:d:psi:u:g:tau:u:uprime}, and since $s \in [0,C_4(\delta)]$ by \autoref{eqn:elementary_compactness_2_point_estimate} we conclude that:
\begin{multline}
\label{eq:prelimit:points:close}
	d^{Q}(q_2,q_2') \leq C_5(\delta) d^Q(g_{\tau(u)uq_1,g_{\tau(u')}u'q_1'}) \\
	\leq
	 C_5(\delta) C(\delta) \epsilon \le 0.01 \epsilon_0'(\delta).  
\end{multline}
where we have used the standing assumption
\S\ref{sec:subsubsec:standing:assumption:epsilons}. 
Therefore, the $Y$-configuration $Y_{sh}'$ left-shadows the
$Y$-configuration $Y_{sh}$.

In view of (\ref{eq:psi:psiprime:close:to:bigcup}), we can choose
$ij \in \tilde{\Lambda}$ so that
\begin{equation}
\label{eq:bfvu:close:to:bEij:bdd}
d\left(\frac{\bfv(u)}{\|\bfv(u)\|}, \bbE_{[ij],bdd}(g_{\tau(u)} u q_1) \right) \le
C(\delta) \eta.
\end{equation}

Since $Y_{bal}$ and $Y_{bal}'$ are $(0,ij)$-top-balanced, and since $s \in
[0,C_4(\delta)]$, $Y_{sh}$ and $Y_{sh}'$ are $(C,ij)$-top-balanced,
for some $C = C(\delta)$. Then, by  \autoref{lem:almost_agreement_of_tau_and_tau_prime}, 
we have $|t_0(Y_{sh})-t_0(Y_{sh}')| \leq C(\delta)$. However the
top-right points of $Y_{sh}$ (resp. $Y_{sh}'$) agree with the
top-right points of $Y_{bal}$ (resp. $Y_{bal}'$). Therefore, 
\begin{equation}
\label{eq:t:ij:t:ij:prime:close}
|t_{ij}(u) - t_{ij}'(u')| \le C_5(\delta). 
\end{equation}
Therefore, by (\ref{eq:gtj:q1:inK}) and (\ref{eq:gtjprime:inK}), we
have
\begin{displaymath}
g_{t_{ij}(u)} q_1 \in K, \quad \text{ and } \quad g_{t_{ij}(u)} q_1' \in
g_{[-C_5(\delta),C_5(\delta)]} K. 
\end{displaymath}
By the definition of $K$, we can find $s'' \in [0, C_5''(\delta)]$
such that 
\begin{displaymath}
q_{3,ij} \equiv g_{s''+t_{ij}(u)} q_1 \in K_0, \quad \text{ and } \quad
q_{3,ij}' \equiv g_{s''+t_{ij}(u)} q_1' \in K_0.  
\end{displaymath}
Let $Y$ be the $Y$-configuration with the same points as $Y_{sh}$,
except that $q_3(Y) = q_{3,ij}$. Similarly, let $Y'$ be the
$Y$-configuration with the same points as $Y_{sh}'$,
except that $q_3(Y') = q_{3,ij}'$. Then, $Y\in K_0$, $Y' \in
K_0$.

By construction, $Y$ and $Y'$ are right-linked. Also, since $s'' \in
[0,C_5''(\delta)]$ since $Y_{sh}$ and $Y_{sh}'$ were $(C(\delta),ij)$
top-balanced, the same is true (with a different constant depending on
$\delta)$ of $Y$ and $Y'$. In  addition, $Y$ and $Y'$ are
left-shadowing and bottom-linked, since the same was true of $Y_{sh}$
and $Y_{sh}'$. We will apply
\autoref{prop:nearby:linear:maps} to $Y$ and $Y'$ in the next
subsection. 



\subsubsection{Taking the limit as $\eta \to 0$.} For fixed
  $\delta$ and $\epsilon$, 
we now take a sequence of $\eta_k \to 0$ (this forces $\ell_k \to
\infty$) and  pass to
limits along a subsequence. Let $\tilde{q}_2 \in K_0$ be the limit of the
$q_2$, and $\tilde{q}_2'\in K_0$ be the limit of the $q_2'$.
We may also
assume that along the subsequence $ij \in \tilde{\Lambda}$ is fixed, 
where $ij$ is as in (\ref{eq:bfvu:close:to:bEij:bdd}). By passing to
the limit in (\ref{eq:prelimit:hd:close}), we get
\begin{equation}
\label{eq:postlimit:hd}
\frac{1}{C(\delta)} \epsilon \le d_{\cH,loc}^{\tilde{q}_2}(U^+[\tilde{q}_2],
U^+[\tilde{q}_2'])  \le C(\delta) \epsilon.
\end{equation}
We now apply \autoref{prop:nearby:linear:maps} (with $\xi \to
0$ as $\eta_k \to 0$). By (\ref{eq:goal:same:sheet}), $\tilde{q}_2'
\in \cW^+[\tilde{q}_2]$. By applying $g_s$ to
(\ref{eq:bfvu:close:to:bEij:bdd}) and then passing to the limit, we get
$U^+[\tilde{q}_2'] \in \cE_{ij}(\tilde{q}_2)$. Finally, 
it follows from passing to the limit in
(\ref{eq:prelimit:points:close}) that
$d^u(\tilde{q}_2,\tilde{q}_2') \le \epsilon_0'(\delta)$, and thus,
since $\tilde{q}_2 \in K_0$ and $\tilde{q}_2' \in K_0$,
it follows from
(\ref{eq:def:epsilon0:prime}) that
$\tilde{q}_2' \in \gB_0[\tilde{q}_2]$. Hence, 
\begin{displaymath}
\tilde{q}_2' \in \Psi_{ij}[\tilde{q}_2].
\end{displaymath}
Now, by (\ref{eq:nearby:maps:final}), we have
\begin{displaymath}
f_{ij}(\tilde{q}_2) \propto P^+(\tilde{q}_2, \tilde{q}_2')_* f_{ij}(\tilde{q}_2'). 
\end{displaymath}
We have $\tilde{q}_2 \in K_0 \subset K_{00} \cap
K_*$, and $\tilde{q}_2' \in K_0 \subset K_*$. 
This concludes the proof of \autoref{prop:half:inductive:step}. 
\qed\medskip

\subsection{The implication have friend implies have invariance}
	\label{ssec:the_implication_have_friend_to_have_invariance}

\subsubsection{Setup}
	\label{sssec:setup_the_implication_have_friend_to_have_invariance}

We now show that \autoref{prop:half:inductive:step} implies \autoref{thm:inductive:step}.
Take a sequence $\epsilon_n \to 0$ and apply
\autoref{prop:half:inductive:step} with $\epsilon =
\epsilon_n$; after passing to a subsequence, assume $ij$ is
constant.
We get, for each $n$ a set $E_n \subset K_*$ 
with $\nu(E_n) > \delta_0$ and with the property that for every $x \in E_n$ there exists $y \in \Psi_{ij}[x] \cap K_*$ such that \autoref{eq:dxy:within:C:epsilon} and
\autoref{eq:half:inductive:step:two} hold for $\epsilon = \epsilon_n$.  
Let 
\begin{displaymath}
F = \bigcap_{k=1}^\infty \bigcup_{n=k}^\infty E_n \subset K_*,
\end{displaymath}
(so $F$ consists of the points which are in infinitely many
$E_n$).
Suppose $x \in F$.
Then
there exists a sequence $y_n \to x$ such that $y_n \in \Psi_{ij}[x]$,
$y_n \not\in U^+[x]$, and such that $\wtd f_{ij}[y_n] \propto \wp^+(x,y_n)_*\wtd  f_{ij}[x]$.
Because $y_n\in \Psi_{ij}[x]$ we have that $\tilde{f}_{ij}[y_n]\propto \tilde{f}_{ij}[x]$ (and this holds for any $x,y_n$ that are in the same atom $\Psi_{ij}[x]$).
We can therefore conclude that
\begin{displaymath}
	\wtd  f_{ij}[x] \propto \wp^+(x,y_n)_* \wtd  f_{ij}[x].
\end{displaymath}
Since $y_n\notin U^+(x)\cdot x$, we have that $\wp^+(x,y_n)\not \in U^+(x)$.

 

\subsubsection{Extra invariance up to scale}
	\label{sssec:extra_invariance_up_to_scale}
For $x \in X$, let $U_{new}^+(x)$ denote the maximal connected
subgroup of $\bbG^{ssr}(x)$ 
such that for $u \in U_{new}^+(x)$ we have that
\begin{equation}
	\label{eq:u:star:tilde:fij}
	u_* \wtd f_{ij}[x] \propto \wtd f_{ij}[x].
\end{equation}
Since $\wp^+(x,y_n)\not \in U^+(x)$ for every $n$ and since $\wp^+(x,y_n)\to \id$ as $y_n\to x$ we have for
$x\in F$, that $U_{new}^+(x)$ strictly contains $U^+(x)$.
We will now check that $U_{new}^+(x)$ yields a family of compatible subgroups as in \autoref{def:compatible_family_of_subgroups}.

Suppose $x \in F$, $t < 0$ and $g_t x \in F$.
Then, since the leafwise measures $\tilde{f}_{ij}[x]$ are $g_t$-equivariant (up to proportionality) we have that \autoref{eq:u:star:tilde:fij} holds for $u\in g_{-t} U^+(g_tx)$.
Therefore, by the maximality of $U_{new}^+(x)$, for $x \in F$ and $t < 0$ with $g_t x \in F$ we have that
\begin{equation}
	\label{eq:g:minus:t:U:new}
	g_{-t} U_{new}^+(g_t x) = U_{new}^+(x),
\end{equation}
and \autoref{eq:u:star:tilde:fij}
still holds.

We now verify that $U^+_{new}(x)$ can be defined for $x$ in a full measure set.
Since $\nu(F) > \delta_0 > 0$ and $g_t$ is
ergodic, for almost all $x \in X$ there exist arbitrarily large $t > 0$ such that $g_{-t} x \in F$.
Then, we define $U_{new}^+(x)$ to be $g_t U_{new}^+(g_{-t}x)$.
This is consistent in view of \autoref{eq:g:minus:t:U:new}.

\subsubsection{Extra invariance on the nose}
	\label{sssec:extra_invariance_on_the_nose}
We next claim for $u \in
U_{new}^+(x)$ that 
\begin{equation}
	\label{eq:extra_invariance_on_the_nose}
	u_* \wtd f_{ij}[x] = \wtd f_{ij}[x]
\end{equation}
and not just that the measures are proportional.
Indeed, from \autoref{eq:u:star:tilde:fij} we get that for a.e. $x$ and $u \in
U_{new}^+(x)$, 
\begin{equation}
	\label{eq:scaling:cocycle}
	u_* \wtd f_{ij}[x] = e^{\beta_x(u)}  \wtd f_{ij}[x],
\end{equation}
where $\beta_x\colon U_{new}^+(x) \to \reals$ is a homomorphism. 

It follows from (\ref{eq:scaling:cocycle}) that for a.e.\ $x \in X$,
$u \in U_{new}^+(x)$ and $t > 0$, 
\begin{equation}
\label{eq:beta:x:equivariant}
	\beta_{g_{-t} x}( g_{-t} u g_t) = \beta_x(u).
\end{equation}
We can write
\begin{displaymath}
\beta_x(u) = L_x(\log u),
\end{displaymath}
where $L_x\colon \fraku^+(x) \to \reals$ is a Lie algebra homomorphism
(which is in particular a linear map). Let $K \subset X$ be a positive
measure set for which there exists a constant $C$ with $\|L_x\| \le C$
for all $x \in K$. Now for almost all $x \in X$ and $u \in
U^+_{new}(x)$ there exists a
sequence $t_j \to \infty$ so that $g_{-t_j} x\in K$ and $g_{-t_j} u
g_{t_j} \to \id$, where $\id$ is the identity element of
$U^+_{new}$.
Then, \autoref{eq:beta:x:equivariant} applied to the
sequence $t_j$ implies that $\beta_x(u) = 0$ almost everywhere.
Therefore, for almost all $x \in X$, 
the conditional measure of $\nu$ along the orbit $U_{new}^+[x]$ is the
push-forward of the Haar measure on $U_{new}^+(x)$, and hence \autoref{eq:extra_invariance_on_the_nose} holds.


\subsubsection{Checking compatibility with admissible partitions}
	\label{sssec:checking_compatibility_with_admissible_partitions}
For checking part (ii) of \autoref{def:compatible_family_of_subgroups}, we proceed with any measurable partition $\frakB$, which is subordinate to the unstable foliation $\cW^u$.
Let $f_\frakB\colon X\to X_\frakB$ be the measurable map to the space of atoms of $\frakB$.
Let also $\sigma_\frakB\colon X_{\frakB}\to X$ be a section of the projection, i.e. a measurable choice of ``center'' in each atom.

Now, since the family of leafwise measures $\tilde{f}_{ij}[x]$ is measurable in $x$, the Lie algebra of its stabilizer $\fraku^+_{new}(x)$ is a measurable family of Lie algebras in $\frakg^{ssr}(x)$.
Recall that $U^+_{new}(x)$ was the corresponding (nilpotent) Lie group and $U^+_{new}[x]\subset \cW^u[x]$ are the corresponding orbits.
We need to check that the partition whose atoms are $U^+_{new}[x]\cap \frakB[x]$ is itself measurable.
We will follow the technique of \autoref{sssec:the_measurability_property}.

Form the measurable bundle $\cO\to X$ whose fiber is $\cO(x):=\cW^u[x]/U'(x)$, where $U'(x)$ ranges over all connected nilpotent Lie subgroups of $\bbG^{ssr}(x)$, of the same dimension as $U^+_{new}(x)$.
The Lie algebras $\fraku'(x)$ determine $U'(x)$ uniquely and the space of Lie algebras is an algebraic subset of a Grassmannian of $\frakg^{ssr}(x)$.

Let next $s_{new}\colon X\to \cO$ be the (tautological) measurable section which maps $x$ to the orbit $U^+_{new}[x]$.
The space $\cO$ carries on the unstables a ``Gauss--Manin'' connection $I(x,y)\colon \cO(x)\to \cO(y)$ obtained by the identification $\cW^u[x]=\cW^u[y]$ for $y\in \cW^u[x]$.

On the subset of $X\times X$ consisting of pairs $(x,y)$ with $x\in \frakB[y]$ the map $I(x,y)$ is measurable (recall that $\frakB$ is subordinate to $\cW^u$).
Recall also that $\sigma\colon X_\frakB\to X$ is a measurable choice
of center, and let $\cO_{\frakB}:=\sigma^* \cO$ be the pullback of the
space of orbits to the ``center'' in each atom; to ease notation we
will write $x_c:=\sigma_{\frakB}(f_{\frakB}(x))\in \frakB[x]$
as the ``projection to the center''.
We can now construct a measurable map $c_{\frakB}\colon X\to\cO_{\frakB}$ that assigns to $x$ the point $I(x,x_c)s_{new}(x)\in \cO(x_c)=\cO_{\frakB}(f_{\frakB}(x))$.
By construction, we have that $c_{\frakB}(x)=c_{\frakB}(y)$ if and only if $x,y$ are in the same atom of $\frakB$ and $y\in U^+_{new}[x]$, or equivalently $x\in U^+_{new}[y]$.
\hfill  \qed

\def \H{\mathcal H}
\def \P{\mathbb P}
\def\vol{\mathrm{vol}}
\def \calL{\mathcal L}
\def \calG{\mathcal G}
\def \calU{\mathcal U}
\def \calS{\mathcal S}
\def \calZ{\mathcal Z}
\def \essspan{\mathrm{ess\, span}}

\appendix


\section{Nonuniform partially hyperbolic background}
	\label{appendix:nonuniform_partially_hyperbolic_background}

We consider a continuous-time dynamical system $g_t\colon Q\to Q$ preserving an ergodic probability measure $\nu$.
The space $Q$ is assumed to be a manifold, though this is not always needed.
All statements are understood to hold $\nu$-a.e. unless further specified.

\subsubsection{Summary of constructions with charts and distances}
	\label{sssec:summary_lyapunov_charts_and_distances}
In addition to standard constructions of stable and unstable manifolds, we will make use of the following distances adapted to the dynamics:
\begin{description}
	\item[Ambient distance] \index{$d^Q$}$d^Q$ coming from the Riemannian metric on $Q$.
	\item[Lyapunov distance] \index{$d^{\cL}_q(q',q'')$}$d^{\cL}_q(q',q'')$ will only be defined for $q',q''\in \cL_{k,\ve}[q]$ (see \autoref{sssec:lyapunov_charts_and_distance}); the parameters $k,\ve$ will be suppressed from the notation.
	\item[Unstable distance] \index{$d^u(q,q')$}$d^u(q,q')$ will only be defined for $q'\in \cW^u[q]$ and will have a uniform contraction property under $g_{-t}$ for $t\geq 0$ (see \autoref{sssec:unstable_lyapunov_distance_and_charts}); the parameter $\ve>0$ is suppressed from the notation.
	\item[Induced ambient distance] \index{$d^{Q,u}$}$d^{Q,u}$ on $\cW^u[q]$ induced by the path metric from the ambient distance.
	\item[Distance on flags] For two flags $F_1, F_2$ of a fixed vector space $V$ with a Riemannian positive-definite bilinear form, we denote by \index{$d$@$\dist$} $\dist(F_1,F_2)$ the corresponding distance in the flag manifold induced by the unique metric which is invariant under $\Orthog(V)$.
\end{description}
We will also construct the following types of charts adapted to the dynamics:
\begin{description}
	\item[Lyapunov charts] \index{$L$@$\cL_{k,\ve}[q]$}$\cL_{k,\ve}[q]\subset Q$ will denote the $(k,\ve)$-Lyapunov chart at $q$ (see \autoref{sssec:lyapunov_charts_and_distance}).
	\item[Local unstable manifold] (see \autoref{sssec:unstable_lyapunov_distance_and_charts})\index{$W$@$\cW^{u}_{loc}[q]$}
	\[
		\cW^{u}_{loc}[q]:=\{q'\colon d^u(q,q')<1\} \subset \cW^u[q].
	\]
\end{description}


\subsection{General dynamics estimates}
	\label{ssec:general_dynamics_estimates}

We begin with some basic notions which are relevant to nonuniformly hyperbolic systems.

\begin{proposition}[Growth dichotomy]
	\label{prop:growth_dichotomy}
	Suppose that $a\colon Q\to \bR$ is measurable.
	Then for a set of full measure $Q_0$:
	\[
		\limsup_{t\to +\infty}\frac{1}{|t|}a(g_t q)
		=
		\limsup_{t\to -\infty}\frac{1}{|t|}a(g_t q)
		\in \{0,\infty\}
	\]
	The analogous result holds for $\liminf$ and $\{0,-\infty\}$.
\end{proposition}
The proof is contained in \cite[4.1.3]{Arnold1998_Random-dynamical-systems}.

\begin{definition}[Temperedness and slow variation]
	\label{def:temperedness_slow_variation}
	Consider a measurable function $A\colon Q\to \bR_{>0}$.
	Then $A$ is called \emph{tempered} if in \autoref{prop:growth_dichotomy}, the function $a:=\log A$ has both $\limsup$ and $\liminf$ equal to $0$.

	The function $A$ is called \emph{$\ve$-slowly varying} if $e^{-\ve|t|}\leq A(q)/A(g_t q)\leq e^{\ve|t|}$ for $\nu$-a.e. $q\in Q$ and $\forall t\in \bR$.
\end{definition}

The following properties are a consequence of Birkhoff ergodic theorem as well as \cite[4.3.2-4.3.3 and Prop.~4.1.3]{Arnold1998_Random-dynamical-systems}.
\begin{proposition}[On temperedness and slow variation]
	\label{prop:on_temperedness_and_slow_variation}
	\leavevmode
	\begin{enumerate}
		\item Suppose $A$ is $\ve$-slowly varying, or more generally satisfies
		\[
			\sup_{t\in [-1,1]}|\log \left(A(q)/A(g_t q)\right)|\leq \phi(q)\quad \text{ for some }  \phi\in L^1(Q,\nu).
		\]
		Then $A$ is tempered.
		\item Suppose $A$ is tempered.
		Then for any $\ve>0$ there exists an $\ve$-slowly varying $B$ such that $A\leq B$.
	\end{enumerate}
\end{proposition}

	


Suppose $E\to Q$ is a cocycle over $g_t$, with a metric on fibers satisfying the integrability assumption
\begin{align}
	\label{eqn:cocycle_integrability_condition}
	\sup_{t\in [-1,1]} \log \norm{g_t^E(q)} \in L^1(Q,\nu).
\end{align}
We state the Oseledets multiplicative ergodic theorem in a form that will be convenient for us (see \cite[Thm.~4.3.4]{Arnold1998_Random-dynamical-systems}).
\begin{theorem}[Tempered Oseledets]
	\label{thm:tempered_oseledets}
	There exists a measurable decomposition into subcocycles $E=\oplus \index{$E^{\lambda_i}$}E^{\lambda_i}$ with distinct $\lambda_i\in \bR$, and a tempered function $C(q)$ such that
	\begin{align*}
		C(q)^{-1} e^{-\ve|t|}\norm{v}\leq \norm{g_t v} e^{\lambda_i t} \leq C(q)e^{\ve |t|}\norm{v} && \forall v\in E^{\lambda_i}(q)\\
		\measuredangle\left(E^{\lambda_i}(q),E^{\lambda_j}(q)\right) > C(q)^{-1} && \forall i\neq j
	\end{align*}
	Note that by \autoref{prop:on_temperedness_and_slow_variation}, for any $\ve>0$ there exists an $\ve$-slowly varying $C_\ve(q)$ such that the above bounds hold with $C_{\ve}$ instead of $C$.
\end{theorem}

\subsubsection{Lyapunov norms}
	\label{sssec:lyapunov_norms}
We recall next a construction from \cite[\S3.5.1-3]{BarreiraPesin2007_Nonuniform-hyperbolicity}, see also \cite[Thm.~4.3.6]{Arnold1998_Random-dynamical-systems}.
For any $\ve>0$, there exists a measurable norm \index{$\norm{-}_{\ve}$}$\norm{-}_{\ve}$ on the cocycle such that
\[
	e^{-|t|\ve}\norm{v_i}_{\ve} \leq \norm{g_t v_i}_{\ve} e^{-t\lambda_i} \leq e^{t|\ve|}\norm{v_i}_{\ve} \quad \forall v_i \in E(q)^{\lambda_i}.
\]
Furthermore, the original norm \index{$\norm{-}$}$\norm{-}$ on $E$ and the $\ve$-Lyapunov norm $\norm{-}_{\ve}$ are related by
\[
	A_{\ve}(q)^{-1}\norm{v}\leq \norm{v}_{\ve}\leq A_{\ve}(q)\norm{v} \quad \forall v\in E(q)
\]
for an $\ve$-slowly varying measurable function $A_{\ve}(q)$.

\subsubsection{One-sided Lyapunov norms}
	\label{sssec:one_sided_lyapunov_norms}

We next describe a variation of the constructions in \cite[\S3.5.1-3]{BarreiraPesin2007_Nonuniform-hyperbolicity}, but adapted to unstable manifolds and backwards dynamics.

\begin{proposition}[Unstable Lyapunov norms]
	\label{prop:unstable_lyapunov_norms}
	Suppose $E\to Q$ is an admissible cocycle which varies smoothly along unstables, equipped with its norm which varies smoothly along unstables.
	Write \index{$E^u$}$E^u\oplus \index{$E^{cs}$}E^{cs}$ for the decomposition into strictly positive and nonpositive Lyapunov exponents, and let $\lambda>0$ be its smallest strictly positive Lyapunov exponent.

	There exists a set of full measure $Q_0\subset Q$, as well as its saturation by the unstable foliation $Q_{0}^u$, such that for any $\ve>0$ there exists a norm \index{$\norm{-}_{u,\ve}$}$\norm{-}_{u,\ve}$ on $E\vert_{Q_0^u}$ which varies smoothly along unstables, and a measurable $\ve$-slowly varying $C_{\ve}$, such that $\forall t\geq 0$:
	\begin{align*}
		& \norm{g_{-t}v}_{u,\ve}\leq e^{-(\lambda-\ve)t}\norm{v}_{u,\ve} && \forall v \in E^u\\
		C_{\ve}^{-1}(q)e^{-\ve t} \norm{v}_{u,\ve} \leq &\norm{g_{-t}v}_{u,\ve} 
		&& \forall v\in E^{cs}(q)
	\end{align*}
\end{proposition}
\begin{proof}
Write $E= E^{u} \oplus E^{u,\perp}$ where $E^u$ is the subcocycle with strictly positive Lyapunov exponents and $E^{u,\perp}$ is its orthogonal complement (both bundles vary smoothly along unstables).
Let $\Pi\colon E\to E^{u,\perp}$ be the projection map.
Note that $\Pi g_t \Pi v = \Pi g_t v$ for any vector $v$.

Let $\lambda$ be the smallest strictly positive Lyapunov exponent of $E$.
For $v\in E^u$ define now the norm (the integral converges by the Oseledets theorem):
\begin{align*}
	\norm{v}_{u,\ve}^2 & :=\int_{0}^{\infty}\norm{g_{-s}v}^2 e^{2(\lambda-\ve)s}ds
\end{align*}
and keep the original norm on $E^{u,\perp}$ and keep the decomposition orthogonal.
It follows by construction that $\norm{g_{-t}v}_{u,\ve}\leq e^{-(\lambda-\ve)t}\norm{v}_{u,\ve}$ for $t\geq 0$ and any $v \in E^u$.
Note also that
\[
	\norm{v} \leq \norm{v}_{u,\ve} \leq A(q) \norm{v} \quad \forall v\in E(q)
\]
for a tempered function $A$.

Now the second bound for the original norm follows from the Oseledets theorem again, and so it follows for the modified norm by adjusting the slowly varying functions.
\end{proof}

The next result allows us to relate the hyperbolic behavior of points along the same stable manifold.

\begin{theorem}[Shadowing Flags]
	\label{theorem:shadowing:flags}
	Let $A_n$ and $B_n$ be sequence of matrices.  Fix $D$ and $\eta>0$ with $\|A_n\|, \|B_n\|, \|A_n^{-1}\|, \|B_n^{-1}\|\le D e^{\eta n }$.
	Fix $\delta>0$ and $\kappa>\eta>0$. 
	Suppose $\{A_n\}$ is forwards regular with Lyapunov exponents $\Lambda$.  Suppose also that 
	$$\|A_n-B_n\|\le \delta e^{-\kappa n}, \quad \quad \|A_n^{-1}-B_n^{-1} \|\le \delta e^{-\kappa n}.$$
	Then
	\begin{enumerate}
	\item  $\{B_n\}$ is forwards regular with Lyapunov exponents $\Lambda$.  
	\item If $\{A_n\}\in \Omega^{+,\Lambda}_{N,\epsilon }$ and $\delta>0$ is sufficiently small, then $\{B_n\}\in \Omega^{+,\Lambda}_{2N,\epsilon }$
	\item The Lyapunov flags for $\{A_n\}$ and $\{B_n\}$ are $\delta^\theta$-close for some $\theta=\theta(\Lambda)$.
	\end{enumerate}
\end{theorem}
\begin{proof}
	This is in different notation \cite[Thm.~4.1]{Ruelle1979_Ergodic-theory-of-differentiable-dynamical-systems}, with the modification that we extend our finite sequence to an infinite one with $B_{n+1}=A_{n+1}$.
\end{proof}



\subsection{Lyapunov charts and distances}
	\label{ssec:lyapunov_charts_and_distances}

\subsubsection{Setup}
	\label{sssec:setup_lyapunov_charts_and_center_stable_manifolds}
We assume that $(Q,g_t)$ has an ambient systole function satisfying the assumptions from \autoref{sssec:systole_function}, and that the ergodic invariant measure $\nu$ satisfies the assumptions preceding \autoref{prop:rescaling_metrics_to_size_1}.

\subsubsection{Start of proof of \autoref{prop:rescaling_metrics_to_size_1}}
	\label{sssec:proof_of_prop:rescaling_metrics_to_size_1}
Recall that \index{$\norm$@$\tnorm{\bullet}$}$\tnorm{\bullet}$ denotes the standard norm on the tangent space $TQ(q)$ and $\exp_q\colon TQ(r_a(q))\to Q$ are our exponential charts.
Recall that in case of Finsler norms, by \autoref{rmk:on_distances_and_norms} that, up to a uniformly bounded constant, we may assume each norm $\tnorm{\bullet}_q$ is induced by an inner product \index{$\norm$@$\tIP{\bullet,\bullet}_q$}$\tIP{\bullet,\bullet}_q$ at every $q$.

\begin{lemma}[Standard hyperbolicity estimates]\label{lem:std}
For every $\ve>0$ there exists a measurable function \index{$A_\ve$}$A_\ve\colon Q\to [1,\infty)$ such that 
\begin{enumerate}
\item for all $t\in \R$ and $v\in TQ^{\lambda_i}(q)$,  $$\frac 1 {A_\ve(q)} e ^{-\ve/8 |t|} \tnorm{ v}\le e^{-\lambda_i t} \tnorm {D g_tv} \le A_\ve(q) e ^{\ve/8 |t|} \tnorm{ v}.$$
\item for all $\lambda_i\neq \lambda_j$, $\sin \angle( TQ^{\lambda_i}(q), TQ^{\lambda_j}(q))\ge \frac 1 {A_\ve(q)}$.  
\item  $A_\ve(g_t q)\le e^{\ve/4 |t|}A_\ve (q)$.
\end{enumerate}
\end{lemma}

\begin{proof}[Proof of \autoref{prop:rescaling_metrics_to_size_1}]
Let $ \norm{\bullet}_{\ve}$ be the Lyapunov norm and let $A_{\ve}\colon Q\to [1,\infty)$ the $\ve/2$-slowly growing measurable function in \autoref{sssec:lyapunov_norms} comparing $\norm{\bullet}_{\ve}$ and $\tnorm{\bullet}$.
By \autoref{prop:on_temperedness_and_slow_variation}, given  $\ve>0$ and $k\in \{1,2,3,\dots\}$ we may take  $\psi_{k,\ve}\colon Q\to [1,\infty)$ an $\ve/2$-slowly growing measurable function such that 
$$\psi_{k,\ve}(q)\ge[A_{\ve}(q)]^2  \max\left\{1,  \phi_a(q)^{e_0+e},  2\ve\inv \phi_a(q) ^{a_k} \right\} .$$
Set \[\phi_{k,\ve}(q) := A_{\ve}(q) \psi_{k,\ve}(q).\]   Then \autoref{prop:rescaling_metrics_to_size_1}\autoref{Lnorm1} holds.  
Define the norm $\norm{\bullet}_{k,\epsilon}$ as 
 \begin{equation}
 \norm{v}_{k, \ve} := \psi_{k,\ve}(q)\norm{v}_{\ve} \quad \text{ for }v\in TQ(q).
\end{equation}
For $v\in TQ(q)$, we have 
\begin{equation} \label{LyapNormBd1}\tnorm{v}\le A_{\ve}(q)  \norm{v}_{\ve}\le 
 \psi_{k,\ve}(q) \norm{v}_{\ve} = \norm{v}_{\ve,k}
\end{equation}
and
\begin{equation} \label{LyapNormBd2}
\norm{v}_{\ve}= \psi_{k,\ve}(q)  \norm{v}_{\ve}\le \psi_{k,\ve}(q)  A_{\ve}(q) \tnorm{v}.
\end{equation}
and so \autoref{prop:rescaling_metrics_to_size_1}\autoref{Lnorm2} holds.
Also, for $v\in TQ(q)$,
\[
\frac{\norm{D_q g_t v}_{k,\epsilon}}{\norm {v}_{k,\epsilon}}\frac{\norm{v}_{\epsilon}}{\norm{D_q g_t v}_{\epsilon}}= \psi_{k,\epsilon}(g_t(q))\psi_{k,\epsilon}(q)\inv.\]
Since $$
e^{-|t|\epsilon/2}\le\psi_{k,\epsilon}(g_t(q))\psi_{k,\epsilon}(q)\inv\le e^{|t|\epsilon/2},$$
\autoref{prop:rescaling_metrics_to_size_1}\autoref{Lnorm3} holds by the properties of the function $A_{\ve}$ in \autoref{sssec:lyapunov_norms}.

Finally, recall that relative to the norms $\tnorm{\bullet}$, for all $t\in [-1,1]$  the map 
\index{$G_{t,q}$}$G_{t,q}\colon B_q\big(q;r_d(q)\big)\to B\big({g_t(q)};r_a(g_t(q))\big)$,
$$G_{t,q}(v) = \exp_{g_t q}\inv \circ g_t\circ \exp_q (v)$$
is well-defined. 
If $\norm{v}_{k, \ve}\le 1$ then by the choice of $ \psi_{k,\ve}(q)$, \[ \tnorm{v}\le A_{\ve}(q) \norm{v}_{\ve} 
= A_{\ve}(q) \psi_{k,\ve}(q) \inv \norm{v } _{k,\ve}\le  \phi_a(q) ^{- e-e_0} \le   \phi_a(q) ^{- e}   r_a(q) = r_d(q)\]
and thus, for $t\in [-1,1]$,  $ G_{t,q} (\bullet)$ is defined on the $\norm{\bullet}_{k,\ve}$-ball $B_{k,\ve}(q;1)$.   
Furthermore, relative to the $\norm{\bullet}_{k,\ve}$-norms, for $v\in B_{k,\ve}(q;1)$ and   $2\le j\le k$ we have 
\begin{align*}
\norm{D^j _vG_{t,q} }_{ \norm{\bullet}_{k,\ve}} 
& \le \phi_{k,\ve}(g_t q)\norm{ D^j_v G_{t,q}  }_{\tnorm{\bullet} }
\left(\frac{A_{\ve}(q)}{\psi_{k,\ve}(q)}\right)^{j}\\
& \le e^{\ve t}  A_{\ve}(q) \psi_{k,\ve}(q) \norm{ D^j_v G_{t,q}  }_{\tnorm{\bullet} }
\left(\frac{A_{\ve}(q)}{\psi_{k,\ve}(q)}\right)^{j}\\
& \le e^{\ve t}  A_{\ve}(q)^{j+1} \norm{ D^j_v G_{t,q}  }_{\tnorm{\bullet} }
{\psi_{k,\ve}(q)}^{-(j-1)}\\
& \le e^{\ve t}  A_{\ve}(q)^{-j+3} \norm{ D^j_v G_{t,q}  }_{\tnorm{\bullet} }
 (2\ve\inv \phi_a(q) ^{a_k}) ^{-(j-1)}\\
& \le   \ve 
\end{align*}
and thus on the  $\norm{\bullet}_{k,\ve}$-ball $B_{k,\ve}(q;1)$,   \begin{equation}\label{eq:Lipschitzderivative}\norm{ D _vG_{t,q}  -Dg_t(q)}_{\norm{\bullet}_{k,\ve}} \le \ve \|v\|_{k,\ve}\end{equation}
and  \begin{equation}\label{eq:Ck in norm}\norm{ G_{t,q}  -Dg_t(q)}_{C^k,\norm{\bullet}_{k,\ve}} \le \ve. \end{equation}

Set $\lambda_0 =\max\{|\lambda_i|,1\}.$ 
Take $r_0 = e^{-2\lambda_0}$.   
For all $0<\ve<1/2$ sufficiently small, $e^{\lambda_0+\ve} +\ve\le e^{\lambda_0+ 2\ve}$.  
Thus for all $t\in [-1,1]$ and $v\in TQ(q)$ with $\norm{v}_{k, \ve}\le r_0$, by \eqref{eq:Lipschitzderivative},
$$\norm{G_{t,q} (v)}_{k,\ve}\le 
(e^{\lambda_0+\ve} + \ve  r_0) \norm{v}_{k, \ve}\le
e^{\lambda_0+2\ve}   r_0<e^{-\lambda_0+2\ve} <1.$$
Thus relative to $\norm{\bullet}_{k, \ve} $-norms, the image of ball $B_{k,\ve}(q;r_0)
$ is contained in the ball $B_{k,\ve}(g_tq;1)$; combined with \eqref{eq:Ck in norm}, \autoref{prop:rescaling_metrics_to_size_1}\autoref{Lnorm4} holds.  
\end{proof}

\subsubsection{Lyapunov charts and distance}
	\label{sssec:lyapunov_charts_and_distance}
Let $y\in Q$ be such that \autoref{prop:rescaling_metrics_to_size_1} holds.  
Let  \index{$L$@$\cL[q]$}$\cL[y] =  \exp_y(  B_{k,\ve}(y;r_0) )$ be the image of the ball of radius $r_0$ in the metric $\norm{\bullet}_{k,\ve}$.  
For $x_1, x_2\in \cL[y]$ introduce the Lyapunov distance based at $y$ via:
\[
	\index{$d^{\cL}_y(x_1,x_2)$}d^{\cL}_y(x_1,x_2):= \norm{\exp_{y}^{-1}(x_1) - \exp_{y}^{-1}(x_2)}_{k,\ve}
\]
Note also that for any $\delta>0$, there exists a compact set $K$ of measure at least $1-\delta$ and $\kappa=\kappa(\delta)$ such that if $y\in K$ and then $d^Q$ and $d^{\cL}_y$ are $\kappa$-biLipschitz in the Lyapunov chart at $y$.

\subsubsection{Unstable distance}
	\label{sssec:unstable_lyapunov_distance_and_charts}
We can apply \autoref{prop:unstable_lyapunov_norms} to $E:=TQ=TQ^{u}\oplus TQ^{cs}$ and obtain a set $Q_0$ of full measure, as well as its saturation by unstables $Q^u_{0}$, together with a norm $\norm{\bullet}_{u,\ve}$ on $TQ\vert_{Q^u_0}$.

For $q\in Q_0^u$ and $q'\in \cW^u[q]$, denote by \index{$d^u(q,q')$}$d^u(q,q')$ the length of the shortest path measured with the metric $\norm{\bullet}_{u,\ve}$.
By construction we have that
\begin{align}
	\label{eqn:uniform_contraction_of_du}
		d^{u}(g_{-t}q,g_{-t}q') \leq e^{-(\lambda_u-\ve)t} d^u(q,q') 
	\quad\forall t\geq 0.
\end{align}
where $\lambda_u>0$ is the smallest strictly positive exponent.
Set $\cW^u_{loc}[q]:=\{q'\colon d^u(q,q')<1\}$ and we drop the $\ve$-dependence from the notation.

\subsection{Stable manifolds}
	\label{ssec:lyapunov_charts_and_stable_manifolds}

\subsubsection{Local and global stable manifolds}
	\label{sssec:stable_manifolds}
For $q\in Q$ we define the \emph{global stable set of $q$} as \index{$W$@$\cW^s[q]$}
\[
	\cW^s[q] = \{q'\colon \exists \ve=\ve(q')>0 \text{ such that }
	e^{\ve t}d(g_tq,g_tq')\xrightarrow{t\to +\infty}0
	\}
\]
Assume now that $(Q,g_t,\nu)$ is measurably good smooth dynamical
system, i.e. satisfying \autoref{def:measurably_good_smooth_dynamics}.

We write $$\lambda_1>\lambda_2> \dots >\lambda_r>0\ge \lambda_{r+1}>\dots>\lambda_\ell$$ for the distinct Lyapunov exponents of $g_t$ with respect to $\nu$, listed without multiplicity.  
Let \index{$\lambda_s:$}$$\lambda_s:= \max\{|\lambda_i|: \lambda_i<0\}.$$


\begin{proposition}[Local stable manifolds relative to Lyapunov charts]
	\label{prop:stableman}
Fix $k\in \N$ and $\ve>0$ sufficiently small. 
Let \index{$W$@$\wtilde{\cW}^s_{k,\ve}(q)$}$\wtilde{\cW}^s_{k,\ve}(q)$ denote the path connected component of 
$$\exp_q^{-1}\left(\cW^s[q]\right) \cap B_{k,\ve}(q;r_0) $$
containing $0\in TQ(q)$.  
Then there exists $C_{\ve,k}>1$ and for $\nu$-a.e.\ $q$, a $C^\infty$ function $h_q\colon W^s(q)\to (W^c\oplus W^u)(q)$ such that following hold:

\begin{enumerate}
	\item $h_q(0) = 0$.
	\item $\|D_0h_q\|=0$ and $\|D_vh_q\|\le \frac 1 3 $ for all $v$.
	\item $\|h_q\|_{C^k} \le C_{\ve, k}$ 
	\item $\wtilde{\cW}^s_{k,\ve}(q)=\mathrm{graph}\left(h_q\right)\cap B_{k,\ve}(q;r_0) $
\end{enumerate}
Moreover, relative to any measurable frame identifying $W^s(q)$ and $ (W^0\oplus W^u)(q)$ with $\R^d$ and $\R^{\dim (Q)-d}$, 
\begin{enumerate}[resume]
\item the map $(Q,\nu)\to C^k(\R^d, \R^{\dim (Q)-d})$,
$$q\mapsto h_q$$ 
is measurable.
\end{enumerate}
Set \index{$W$@$\cW^s_{k,\ve}[q]$}$\cW^s_{k,\ve}[q]:=\exp_q\left(\wtilde{\cW}^s_{k,\ve}(q)\right)$.
Then, following notation in \autoref{prop:rescaling_metrics_to_size_1}
\begin{enumerate}[resume]
\item We have
$$
\wtilde{\cW}^s_{k,\ve}(q) = \left\lbrace v\in  B_{k,\ve}(q;r_0):  \limsup_{t\to\infty} \tfrac 1 t \log  \norm{\exp_{g_tq}^{-1} \circ g_t \circ \exp_q(v)}<0\right \rbrace.
 $$
 \item For all $v_1,v_2\in\wtilde{\cW}^s_{k,\ve}(q)$ and $t\ge 0$,
$$\norm{\exp_{g_tq}^{-1} \circ g_t \circ \exp_q(v_1)- \exp_{g_tq}^{-1} \circ g_t \circ \exp_q(v_2)}\le e^{-\lambda_s t + 2\ve t}$$
 \item For all $x,y\in \cW^s_{k,\ve}[q]$,
 $$d^Q(g_tx,g_ty)\le \phi_{k,\ve} (q) e^{-\lambda_s t + 2\ve t}d^Q(x,y).$$
  \item For all $t\geq 0$, we have
 \[
 	 g_{-t}\left(\cW^{s}_{k,\ve}[q]\right) \supset \cW^s_{k,\ve}[g_{-t}q].
 \]

\end{enumerate}
\end{proposition}

\noindent For $\nu$-a.e. $q$ we have  $$\cW^s[q]= \bigcup_ {t\ge 0} g_{-t} \left( \cW^s_{k,\ve}[g_tq]\right).$$
In particular, for $\nu$-a.e.\ $q$, the set $\cW^s[q]$ is an injectively immersed $C^\infty$ submanifold tangent to $W^s(q)$.

\begin{proof}[Proof of \autoref{prop:stableman}]
	This is standard and can be found, for instance, in \cite[Ch.~7]{BarreiraPesin2007_Nonuniform-hyperbolicity}.
\end{proof}

\begin{remark}\label{rem:WlocinLyapcharts}
We remark that given a.e.\ $z\in Q$, $\cW^u_{loc}[z]$ is a pre-compact subset of $\cW^u[z]$.
It follows from the construction of global (un)stable manifolds that there is a tempered $T_0(z)>0$ such that for all $t\ge T_0(z)$, $$g_{-t} \cW^u_{loc}[z]\subset \cL_{k,\ve}[g_{-t}z]$$
where $\cW^{u}_{loc}[z]$ is defined in \autoref{sssec:unstable_lyapunov_distance_and_charts}.
\end{remark}

\subsection{Center-stable manifolds}
	\label{ssec:center_stable_manifolds}

\subsubsection{Summary of constructions}
	\label{sssec:summary_of_constructions_center_stable_manifolds}
We recall that the formal jet of the center-stable manifold is well-defined and equivariant under the dynamics, but there does not exist a uniquely-defined ``center-stable'' manifold (or set) which is also equivariant for the dynamics.
Analogously, for a smooth and natural cocycle $E$, or more generally for a forward-dynamically defined one (see \autoref{def:forward_dynamically_defined_cocycle}), the formal jets of the forward Oseledets filtrations are also well-defined and equivariant.
For our applications, we will use two workarounds for this issue:
\begin{description}
	\item[Realized Center-Stables] We will fix, once and for all, a $C^\infty$ manifold realizing the jet of the center-stable manifold denoted \index{$W$@$\cW^{cs}_{loc}[q]$}$\cW^{cs}_{loc}[q]$.
	For cocycles we will also fix a $C^\infty$ realization \index{$E^{\leq \lambda_i}_q(q')$}$E^{\leq \lambda_i}_q(q')$ for $q'\in \cW^{cs}_{loc}[q]$ which varies smoothly in $q'$.

	\item[Fake Center-Stables] 
	For each degree of regularity $k$ we will construct a family \index{$W$@$\cW^{cs}_{fake,k}[q]$}$\cW^{cs}_{fake,k}[q]$ with better equivariance properties and controlled norms.
	For cocycles we will also construct a $C^k$ realization \index{$E^{\leq \lambda_i}_{fake, k, q}(q')$}$E^{\leq \lambda_i}_{fake, k, q}(q')$ for $q'\in \cW^{cs}_{fake,k}[q]$ which varies $C^k$ in $q'$.

	\item[Coherence] 
	The manifold $\cW^{cs}_{fake,k}[q]$ is tangent to order $k$ to $\cW^{cs}_{loc}[q]$, and both contain the stable manifold $\cW^{s}_{loc}[q]$ (in the neighborhood where all are defined).	
\end{description}

For the following constructions, recall that \index{$B_{k,\ve}(q;r)$}$B_{k,\ve}(q;r)$ denotes the ball of radius $r$ in $TQ(q)$ in the $(k,\ve)$-Lyapunov norm, see \autoref{prop:rescaling_metrics_to_size_1}.

\subsubsection{Fake center-stable manifolds}
	\label{sssec:center_stable_manifolds}
Fix an integer $k\geq 2$ and $\ve>0$ sufficiently small depending on $k$ and the Lyapunov spectrum.
Let $\cL_{k,\ve}[q]\subset Q$ denote the Lyapunov chart at $q$.
Using a standard construction with bump functions, define \index{$g$@$\wtilde{g}_t^{(k)}$}$\wtilde{g}_t^{(k)}\colon TQ\to TQ$ such that for $t\in [-1,1]$ we have:
\begin{align*}
	\wtilde{g}_t^{(k)}\colon T_q Q &\to T_{g_t q} Q\\
	\wtilde{g}_t^{(k)}(v) & = 
	\begin{cases}
		Dg_t(v) &\text{ if }v\notin B_{k,\ve}(q;r_0/2)\\ 
		\exp_{g_t q}^{-1}\circ g_t\circ \exp_q(v) &\text{ if }
		v \in B_{k,\ve}(q;r_0/4)
	\end{cases}
\end{align*}
where $r_0$	is defined in \autoref{prop:rescaling_metrics_to_size_1}.
With this fixed globalization, it follows from \cite[Thm.~6.6]{PughShub1989_Ergodic-attractors} that there exists a unique family of maps and manifolds
\begin{align*}
	\index{$\gamma^{cs}_{fake,k}$}\gamma^{cs}_{fake,k}\colon W^{cs}(q) & \to W^{u}(q)\\
	\index{$W$@$\wtilde{W}^{cs}_{fake,k}(q)$}\wtilde{W}^{cs}_{fake,k}(q) & 
	:= \left\lbrace(v,\gamma^{cs}_{fake,k}(v)) \colon v\in W^{cs}(q)\right\rbrace\\
	\index{$W$@$\cW^{cs}_{fake,k}[q]$}\cW^{cs}_{fake,k}[q] & 
	:= \exp_q  \left(\wtilde{W}^{cs}_{fake,k}(q) \cap B_{k,\ve}(q;r_0)\right)
\end{align*}
such that $w \in \wtilde{W}^{cs}_{fake,k}$ if and only if $e^{-\ve t}\wtilde{g}_t^{(k)} w  \xrightarrow{t\to +\infty} 0$.
In particular, it follows that
\begin{enumerate}
	\item We have
	\[
		\wtilde{g}_t^{(k)}\left(\wtilde{W}^{cs}_{fake,k}(q)\right) = \wtilde{W}^{cs}_{fake,k}\left(g_t q\right) \quad \forall t\in \bR
	\]
	\item We have
	\[
		g_t\left(\cW^{cs}_{fake,k}[q] \cap \exp_q\left(B_{k,\ve}(q;e^{-\ve t})\right) \right) \subset \cW^{cs}_{fake,k}[g_t q] \quad
		\forall t \geq 0.
	\]
\end{enumerate}

\subsubsection{Fake cocycle flags along center-stables}
	\label{sssec:face_cocycle_flags_along_center_stables}
%
Let $E\to Q$ be a smooth and natural cocycle (see \autoref{def:smooth_and_natural_cocycles}), and fix an integer $k\geq 2$ and $\ve>0$ sufficiently small depending on $k$ and the Lyapunov spectrum of $TQ$ and $E$.
Consider the vector bundle\index{$E$@$\wtilde{E}$}
\begin{align*}
	\wtilde{E}:=TQ\times_{Q} E :=\{(q,v,e)\colon q\in Q, v\in TQ(q), e\in E(q)\} \to TQ
\end{align*}
We defined in \autoref{sssec:center_stable_manifolds} a flow $\wtilde{g}_t^{(k)}$ on $TQ$, and we now extend this to a cocycle on $\wtilde{E}$.
To avoid confusion, we will write \index{$g_t^E$}$g_t^E$ for the cocycle on $E$ and \index{$g$@$\wtilde{g}_{t}^{k,\wtilde{E}}$}$\wtilde{g}_{t}^{k,\wtilde{E}}$ for the cocycle on $\wtilde{E}$.

First, recall that we have trivializations
\[
	E(q)\times B_{k,\ve}(q;r_0)
	\xrightarrow{\chi^{E}_{k,\ve}(q)}
	E\vert_{\cL_{k,\ve}[q]}
\]
which cover the map $\exp_q \colon B_{k,\ve}(q;r_0)\to \cL_{k,\ve}[q]$.
Using a bump function, define for $e\in \wtilde{E}(q,v)$ with $(q,v)\in TQ$ and $t\in [-1,1]$:
\begin{align*}
	\wtilde{g}_t^{k,\wtilde{E}}(e) & = 
	\begin{cases}
		g_t^E(q)(e) &\text{ if }v\notin B_{k,\ve}(q;r_0/2)\\
		g_t^E(\exp_q(v))\left(\chi^{E}_{k,\ve}(q)(v,e)\right) &\text{ if }
		v \in B_{k,\ve}(q;r_0/4)
	\end{cases}
\end{align*}
Using \cite[Thm.~6.3]{PughShub1989_Ergodic-attractors} or \cite[Rmk.~5.2(b)]{Ruelle1979_Ergodic-theory-of-differentiable-dynamical-systems} we obtain a $\wtilde{g}_{t}^{k,\wtilde{E}}$-equivariant filtration by order of growth $\wtilde{E}^{\leq \bullet}_{fake,k}(q,v)$ which is $C^k$ along $\wtilde{W}^{cs}_{fake,k}$ and we set
\[
	E^{\leq \bullet}_{fake, k, q}(q'):=
	\chi_{k,\ve}^{E}\left(q,\exp_{q}^{-1}(q')\right)
	\left(\wtilde{E}^{\leq \bullet}_{fake,k}\left(q,\exp_{q}^{-1}(q')\right)\right)
\]
for $q'\in \cL_{k,\ve}[q]$.

\begin{proposition}[Properties of center-stable jets]
	\label{prop:properties:fake:flags}
	There exists a set of full measure $Q_0$ such that all $q,q_1,q_2\in Q_0$  we have:
	\begin{enumerate}
		\item\label{item:fakeflag:equivariance}
		The $k$-jet of $\wtilde{W}^{cs}_{fake,k+1}(q)$ agrees with that of $\wtilde{W}^{cs}_{fake,k}(q)$ and is $g_t$-equivariant.
		\item\label{item:fakeflag:smoothness}
		Suppose $q'\in \cW^{s}_{loc}[q_1]=\cW^{s}_{loc}[q_2]$.
		Then the formal $k$-jets of $\cW^{cs}_{fake,k}[q_1]$ and $\cW^{cs}_{fake,k}[q_2]$ agree at $q'$.
		\item The analogous statements hold for $\wtilde{E}^{\leq \bullet}_{fake,\bullet}$.	
	\end{enumerate}
\end{proposition}
\begin{proof}
	By construction $\wtilde{g}_{t}^{(k')}$ agrees with $\wtilde{g}_{t}^{(k)}$ in a neighborhood of the zero section.
	Recall also that for fixed $(k,\ve)$ the center-stable manifold $\wtilde{W}^{cs}_{fake,k}(q)$ is unique and can be obtained by iterating the dynamics.
	Therefore, one can obtain $\wtilde{W}^{cs}_{fake,k}(q)$ as the limit of $\wtilde{g}_{-t}^{(k)}\wtilde{W}^{cs}_{fake,k+1}\left(g_{t}q\right)$ as $t\to +\infty$ and the first claim follows.
	Note that equivariance of the flags follows from their construction and uniqueness.

	For the second claim, it suffices to prove it if $q_2$ is sufficiently close to $q_1$, depending on the $(k,\ve)$-Lyapunov radius at $q_1$ and then use equivariance under the dynamics.
	But if $q_2,q_1$ are sufficiently close, the maps $\wtilde{g}_t^{(k)}(q_1)$ and $\wtilde{g}_{t}^{(k)}(q_2)$ are conjugated in a sufficiently small neighborhood containing the origin (using $\exp_{q_1}\circ \exp_{q_2}^{-1}$) and again it follows that $\wtilde{W}^{cs}_{fake,k}(q_2)$ can be obtained by iterating $\exp_{q_2}^{-1}\left(\wtilde{W}^{cs}_{fake,k}[q_1]\right)$.

	The proof of the last claim is analogous to the preceding two.
\end{proof}

\subsubsection{Realized Center-Stable Manifolds}
	\label{sssec:center_stable_manifolds_factorization}
We fix, once and for all, a measurable family $q\mapsto \index{$W$@$\wtilde{W}^{cs}(q)$}\wtilde{W}^{cs}(q)\subset TQ(q)$ of embedded, $C^\infty$ submanifolds that have the same formal jet at the origin as the family $\wtilde{W}^{cs}_{fake,k}$ as $k\to +\infty$.
We realize these manifolds as graphs of smooth functions $\wtilde{\phi}_{q}\colon W^{cs}(q)\to W^u(q)$.

Denote by \index{$W$@$\cW^{cs}_{loc}[q]$}$\cW^{cs}_{loc}[q]:=\exp_{q}\left(\wtilde{W}^{cs}(q)\cap TQ(q;r_{k,\ve}(q))\right)$ where the $(k,\ve)$ are fixed in \autoref{sssec:choice_of_k_epsilon}.
Note that for any $k\geq 2$ and sufficiently small $\ve,\delta>0$ there exist a compact set $K$ of measure at least $1-\delta$ and a constant $A(\delta,k,\ve)$ such that $\norm{\wtilde{\phi}_q}_{C^k(B_{k,\ve}(q;r_0))}\leq A(\delta,k,\ve)$.

\subsubsection{Realized Center-Stable Oseledets filtrations}
	\label{def:fixed_realization_center_stable}
Let now $E\to Q$ be a smooth and natural cocycle.
We fix once and for all a measurable family of $C^\infty$ sections $q'\mapsto \index{$E^{\le \bullet}_q(q')$}E^{\le \bullet}_q(q')$ from $\cW^{cs}_{loc}[q]$ into the bundle of all flags in $E$
such that given $k\in \N$ the map $q'\mapsto E^{\le \bullet}_q(q')$ is $C^k$ tangent at $q$ to $q'\mapsto E_{fake,k, q}^{\le\bullet}(q')$ as above.
The analogous remark about $C^k$ bounds on these maps apply.

\subsection{Measurable partitions}
	\label{ssec:measurable_partitions}

In this section we discuss some general constructions concerning measurable partitions.

\subsubsection{Some generalities on partitions}
	\label{sssec:some_generalities_on_partitions}
We refer to \cite[\S5-7]{EinsiedlerLindenstrauss2010_Diagonal-actions-on-locally-homogeneous-spaces} for background on the following constructions.
We mainly consider partitions that are measurable (i.e. generated by countably many Borel sets) and denote the atom of a partition $\gB_0$ containing a point $x$ by \index{$B$@$\gB_0[x]$}$\gB_0[x]$.
A probability measure $\nu$ then has conditional measures on each atom, denoted \index{$\nu_{\frakB}[x]$}$\nu_{\frakB}[x]$.
Recall also that a partition is subordinated to another if every atom of the first is contained in an atom of the second.

We work in the setting of a flow $g_t$ on a manifold $Q$, preserving an ergodic probability measure $\nu$, with $\cW^u[q]$ denoting the unstable manifold through $q$ (when it exists).
For a measurable partition $\frakB$ we will use the shorthand:\index{$B$@$\frakB_t[q]$}
\[
	\frakB_t[q]:=g_{-t}\frakB[g_t q]
\]
so in particular $\gB_0=\frakB$.

\begin{definition}[Unstable Markov partition]
	\label{def:unstable_markov_partition}
	A measurable partition $\gB_0$ of $Q$ is called a \emph{unstable Markov partition} if it is subordinated to the partition into unstable manifolds, \index{$B$@$\gB_0[q]$}$\gB_0[q]$ is relatively open and relatively compact in $\cW^u[q]$, and furthermore for a set of $q$ of full measure we have:
	\begin{itemize}
		\item We have $\frakB_t[q]\subset \gB_0[q]$ for all $t\geq 0$.  
		\item We have that $\gB_0[q]\subset \cW^u_{loc}[q]$, where $\cW^u_{loc}[q]$ is the ball of radius $1$ around $q$ in the $d^u$-distance (see \autoref{sssec:unstable_lyapunov_distance_and_charts}).
	\end{itemize}	
\end{definition}

We note that since $h_\nu(g_1)<\infty$, we have $h(g_1,\frakB) = h(g_{-1} \frakB\vert \frakB)\le h_\nu(g_1)<\infty$ for any measurable partition $\frakB$ with the Markov property $\frakB_t= g_{-t} \frakB \prec \frakB$ for $t\ge 0$.
Unstable Markov partitions $\frakB$ have the further property that they cary all the entropy of the flow $g_t$. 

\begin{proposition}\label{prop:full_entropy}
For any unstable Markov partition $\frakB$ as in \autoref{def:unstable_markov_partition},
$$h_\nu(g_1)= h(g_1,\frakB) := h(\frakB_1\vert \frakB)= h(g_{-1} \frakB\vert \frakB).$$
\end{proposition}

The proof of \autoref{prop:full_entropy} follows exactly as in the classical results \cite{LedrappierYoung_The-metric-entropy-of-diffeomorphisms.-I.-Characterization-of-measures-satisfying-Pesins} after showing the dynamics admits good Lyapunov charts.

\subsubsection{Hausdorff distance properties of the partitions}
	\label{sssec:hausdorff_distance_properties_of_the_partitions}
Let $\cK$ denote the space of compact subsets of $Q$.
The Hausdorff distance coming from the ambient metric $d^Q(\cdot, \cdot)$ on $Q$ turns $\cK$ into a metric space. 

The map $q \to \ov{\gB_0[q]}$ is measurable with respect to the Borel $\sigma$-algebra on the metric space $\cK$.
Then, by Lusin's theorem, for every $\delta > 0$ there
exists a compact set $K \subset Q$ such that the map $q \to \ov{\gB_0[q]}$
is continuous on $K$.

\subsubsection{Measurable partitions and subgroups compatible with the measure}
	\label{sssec:measurable_partitions_and_subgroups_compatible_with_the_measure}
Assume that we have a family $U^+(q)$ as defined in \autoref{def:compatible_family_of_subgroups}, relative to an unstable Markov partition $\gB_0$.
We then define for $t\in \bR$:\index{$B$@$\cB_t[q]$}\index{$B$@$\cB_t(q)$}
\begin{align}
	\label{eq:U_plus_frakb_partition}
	\begin{split}
	\cB_t[q]& :=\frakB_t[q]\cap U^+[q]\\
	\cB_t(q)& :=\{u\in U^+(q)\colon uq\in \cB_t[q]\}
	\end{split}
\end{align}
Thus $\cB_0[q]$ is a subset of $U^+[q]$, while $\gB_0[q]$ is a subset of the unstable manifold $\cW^u[q]$.
Note that $U^+[q]$ admits a natural Haar measure up to a normalizing scalar, and therefore $\cB_0[q]$ admits a natural Haar probability measure.

We also use a Markov partition \index{$B$@$\gB_0^-$}$\gB_0^-$ of $Q$ subordinate to the stables,
and pull it back to $X$ as in \S\ref{sec:subsec:pulling:back:Markov}. 
For $x \in X$, let
\begin{align*}
\index{$B$@$\gB_t[x]$}\gB_t[x] & = g_{-t}(\gB_0[g_t x]) 
\index{$B$@$\gB_t^-[x]$}\gB_t^-[x] & = g_{t}(\gB_0^-[g_{-t} x]) 
\end{align*}

We note the following lemma for future reference:
\begin{lemma} $ $
\label{lemma:gB:properties}
\begin{itemize}
\item[{\rm (a)}] For $t' > t \ge 0$, $\gB_{t'}[x] \subset \gB_t[x]$. 
\item[{\rm (b)}] Suppose $t \ge 0, t' \ge 0$, $x \in X_0$ and $x' \in X_0$ are such that 
$\gB_t[x] \cap \gB_{t'}[x'] \ne \emptyset$. Then either $\gB_t[x]
\supseteq \gB_{t'}[x']$ or $\gB_{t'}[x'] \supseteq \gB_t[x]$ (or
both). 
\end{itemize}
\end{lemma}

For the next statement, recall (see \autoref{sssec:unstable_lyapunov_distance_and_charts}) that $\cW^{u}_{loc}[q]$ is the ball of radius $1$ in the $d^u$-distance.
\begin{lemma}[Partition adapted to a measurable function]
	\label{lem:partition_adapted_to_a_measurable_function}
	Suppose that $K$ is a measurable set with $\nu(K)>0$ and $T_0\colon K\to \bR_{>0}$ is a measurable function.
	
	Then there exists an unstable Markov partition $\gB_0$ with $\gB_0[q]\subset \cW^u_{loc}[q]$ for a.e. $q$, and with the following additional properties.
	We can pick a measurable subset $K_0\subset K$,
	with the following properties:
	\begin{itemize}
		\item For any two $q\neq q'\in K_0$ we have $\gB_0[q]\cap \gB_0[q']=\emptyset$.
		Write $\gB_0[K_0]:=\bigsqcup_{q\in K_0}\gB_0[q]$.
		\item For a.e. $q\in Q$ the set $\{t\colon g_{-t}q\in \gB_0[K_0], t>0 \}$ is nonempty and has a smallest element $T_r(q)$ (called the recurrence time of $q$).
		\item For a.e $q$ we have $\gB_0[g_{-t} q]=g_{-t} \gB_0[q]$ for $t\in[0,T_r(q))$.
		\item For all $q\in K_0$ the ``recurrence time'' satisfies $T_r(q)>T_0(q)$.
	\end{itemize}
\end{lemma}
\subsubsection{Centers of atoms}
 	\label{sssec:centers_of_atoms}
Note that the function $T_r(q)$ is constant on atoms of $\gB_0$.
It also follows that there exists $T_r^-(q)\geq 0$ such that $g_{T_r^-(q)}q\in K_0$ and $g_tq\notin K_0$ for all $t\in [0,T_r^-(q))$.
It follows that for a full measure set of atoms $\gB_0[q]$, we can define its ``center'' \index{$q^*$}$q^*\in\gB_0[q]$ as $q^*:=g_{-T_r^-(q)}q_0$ where $q_0\in K_0$ and $g_{T_r^-(q)}q\in\gB_0[q_0]$.
Note that for an atom intersecting $K_0$, its center is the corresponding point of $K_0$.
\begin{proof}[Proof of \autoref{lem:partition_adapted_to_a_measurable_function}]
	The proof is entirely parallel to \cite[Prop.~3.7]{EskinMirzakhani_Invariant-and-stationary-measures-for-the-rm-SL2Bbb-R-action-on-moduli-space}, see also \cite[\S9]{MargulisTomanov1994_Invariant-measures-for-actions-of-unipotent-groups-over-local}.
\end{proof}



\begin{remark}[Choice of measurable partition]
	\label{rmk:choice_of_measurable_partition_and_du_smallness}
	We will construct \index{$B$@$\gB_0$}$\gB_0$ such that \autoref{prop:good_norms} holds, with the function $T_0$ chosen there.
	We record here as well that by \autoref{lem:partition_adapted_to_a_measurable_function} the atoms satisfy $\gB_0[q]\subset \cW^{u}_{loc}[q]$.
\end{remark}



\subsection{Conditional measures along (un)stable manifolds}
	\label{ssec:conditional_measures}

\begin{definition}[{The measures $\nu_{\gB_0}[x]$}]
	\label{def:measures_nu_gB0}
We now define \index{$\nu_{\gB_0}[x]$}$\nu_{\gB_0}[x]$ 
to be the conditional measure
of $\nu$ along the measurable partition
$\gB_0[x]$. In other words, $\nu_{\gB_0}[x]$
is the measure on $\gB_0[x]$ defined so that for any measurable
$\phi: X \to \bR$, 
\begin{displaymath}
\mathbb{E}(\phi \mid \gB_0)(x) = \int_X \phi \, d\nu_{\gB_0}[x]. 
\end{displaymath}
\end{definition}

\begin{definition}[{The measures $\nu^u_x$}]
By \autoref{lemma:Psiij:equivariant},  for almost every $x$ and all $t\ge 0$ we have 
$$ \gB_0[ x]  \subset g_t \gB_0[g_{-t} x].$$
We define a locally finite Radon measure 
\index{$\nu^u_x$}$\nu^u_x$ 
on $\cW^u[x]$  as follows: we declare $\nu^u_x(\gB_0[x]) = 1 $ and for any compact $K\subset \cW^u[x]$, we take $t>0$ so that $g_{-t}K\subset \gB_0[g_{-t}x]$ and declare 
$$\nu^u_x(K) := \frac{
(g_t)_*  \nu_{\gB_0[g_{-t}x]} (K)
}{(g_t)_*  \nu_{\gB_0[g_{-t}x]} (\gB_0[x])}.$$
The $g_t$-invariance of the measure $\nu$ ensures that $\nu^u_x$ is well-defined independent of the choice of 
$t>0$ making  $g_{-t}K\subset \gB_0[x].$
\end{definition}

Similarly, we get locally finite radon measures \index{$\nu^s_x$}$\nu^s_x$ along $W^s[x]$.

The following is used repeatedly:
\begin{lemma}[Points of density]
\label{lemma:gB:vitali}
Suppose $\delta > 0$ and $K \subset X_0$ is such that
$\nu(K) > 1-\delta$. Then there exists a subset $K^* \subset
K$ with $\nu(K^*) > 1 - \delta^{1/2}$ such that for any $x \in K^*$,
and any $t > 0$, 
\begin{displaymath}
\nu^u_x(K \cap \gB_t[x]) \ge (1-\delta^{1/2}) \nu^u_x(\gB_t[x]). 
\end{displaymath}
\end{lemma}

\begin{proof} Let $E = K^c$, so $\nu(E) \le \delta$. Let $E^*$
denote the set of $x \in X_0$ such that there exists some $\tau \ge 0$ with
\begin{equation}
\label{eq:gB:lots:of:bad:set}
\nu^u_x(E \cap \gB_\tau[x]) \ge \delta^{1/2} \nu^u_x(\gB_\tau[x]). 
\end{equation}
It is enough to show that $\nu(E^*) \le \delta^{1/2}$. 
Let $\tau(x)$ be the smallest $\tau>0$ so that (\ref{eq:gB:lots:of:bad:set})
holds for $x$. Then the (distinct) sets $\{\gB_{\tau(x)}[x]\}_{x \in E^*}$
cover $E^*$ and are pairwise disjoint by
\autoref{lemma:gB:properties} (b). Let 
\begin{displaymath}
F = \bigcup_{x \in E^*} \gB_{\tau(x)}[x].
\end{displaymath}
Then $E^* \subset F$. For every set of the form $\gB_0[y]$, let
$\Delta(y)$ denote the set of distinct sets $\gB_{\tau(x)}[x]$ where
$x$ varies over $\gB_0[y]$.  Then, by (\ref{eq:gB:lots:of:bad:set})
\begin{multline*}
\nu^u_y(F \cap \gB_0[y]) = \sum_{\Delta(y)} \nu^u_y(\gB_{\tau(x)}) \le \\
\le \delta^{-1/2} \sum_{\Delta(y)} \nu^u_y(E \cap \gB_{\tau(x)}[x]) \le
\delta^{-1/2} \nu^u_y(E \cap \gB_0[y]). 
\end{multline*}
Integrating over $y$, we get $\nu(F) \le \delta^{-1/2}
\nu(E)$. Hence, 
\begin{displaymath}
\nu(E^*) \le \nu(F) \le \delta^{-1/2} \nu(E) \le \delta^{1/2}.  
\end{displaymath}
\end{proof}

\subsubsection{The measure $| \cdot |$.}
\label{sssec:the:measure:on:Uplus:orbits}
Recall that $\cB_0[x] = \gB_0[x] \cap U^+[x]$. 
Since $U^+(x) \subset \bbG^{ssr}(x)$ is unipotent, it is
unimodular, and the same holds for $U^+_x(x)$ which is the stabilizer
of $x$ in $U^+(x)$). Therefore, we can define
``Lebesgue measure'' \index{$| \cdot |$}$| \cdot |$ on $\cB_0[x]$ as 
the pushforward of the Haar measure on $U^+(q)/U^+_q(q)$ to $\cB_0[q]$
under the action map $u \to u q$. (The measure $| \cdot |$ is only
defined up to normalization).

The same argument as \autoref{lemma:gB:vitali} 
also proves the following:
\begin{lemma}
\label{lemma:cB:vitali:substitute}
Suppose $\delta > 0$, $\theta' > 0$ and $K \subset X$, with 
$\nu(K) > 1-\delta$. 
Then there exists a subset $K^* \subset
K$ with $\nu(K^*) > 1 - \delta/\theta'$ 
such that for any $x \in K^*$, and any $t > 0$, 
\begin{displaymath}
|K \cap \cB_t[x]| \ge (1-\theta') |\cB_t[x]|,  
\end{displaymath}
and thus
\begin{displaymath}
|\{ u \in \cB_t(x) \st u x \in K \}| \ge (1-\theta') 
|\cB_t(x)|. 
\end{displaymath}
\end{lemma}




\section{Generalities on cocycles}
	\label{appendix:generalities_on_cocycles}

We will work with a number of flavors of cocycles.
We enumerate them here and outline the basic logical structure in their definitions, although some of the flavors below are not introduced until later sections.
We will consider:
\begin{description}
	\item[Smooth and Natural] Start with the tangent, and more generally higher-jet cocycles on $Q$, and their natural morphisms.
	Then standard linear-algebraic operations (quotients, duals, tensors) lead to the class of smooth and natural cocycles.
	These admit smooth trivializations and are defined everywhere on $Q$.
	See \autoref{def:smooth_and_natural_cocycles}.
	\item[Fdd and Bdd] Forward dynamically defined (fdd) cocycles are those cocycles that arise once we introduce the subcocycles of the smooth and natural that are given by stable Lyapunov filtrations, and again allow for standard linear-algebraic operations.
	These cocycles are defined only along a subset of $Q$ saturated by stable manifolds, and vary smoothly along stable manifolds.
	See \autoref{def:forward_dynamically_defined_cocycle}.
	The bdd (backward dynamically defined) cocycles are the analogues along unstables.
	\item[Admissible] These are the weakest notion, with only a measurable isomorphism to a measurable subcocycle of a smooth and natural one.
	See \autoref{def:admissiblecocycle}.
\end{description}


\subsection{Algebraic hulls}
	\label{ssec:algebraic_hulls}

\subsubsection{Setup}
	\label{sssec:setup_algebraic_hulls}
We will make use of the following basic principles in this section.
Given a fixed vector space $V$, an algebraic subgroup $H\subset \GL(V)$ can be defined by specifying a list of subspaces (or vectors) that $H$ fixes in some tensor constructions applied to $V$.
Conversely, a theorem of Chevalley says that any algebraic subgroup $H\subset \GL(V)$ can be defined as the stabilizer of some line in some tensor construction on $V$.

Given two vector spaces $V_1,V_2$, it is not possible to say that a linear isomorphism $V_1\to V_2$ is in a specific group.
It is possible, however, to say that the linear isomorphism is a bijection on the specific list of subspaces or tensors that determine the corresponding groups for $V_1$ and $V_2$.
Using this point of view, it is frequently possible and convenient to avoid measurable trivializations of cocycles and keep more accurate track of the regularity of various maps involved.

\subsubsection{Trivializations and principal bundles}
	\label{sssec:trivializations_and_principal_bundles}
Suppose that $E\to Q$ is a vector bundle and $V$ is a fixed vector space (of the same dimension as the fibers of $E$).
Instead of globally trivializing the bundle $E$ (in some measurable way, for instance) it is useful to consider the principal bundle of frames \index{$P_E$}$P_E$ with fiber given by $P_E(q)=\Isom(V, E(q))$, where $\Isom$ denotes linear isomorphisms.
The group $\GL(V)$ acts on the right by precomposing the framings and the action is transitive on the fibers.
If $E$ had a symplectic structure or symmetric pairing, we could consider isomorphism that preserve this extra structure (and consider a group smaller than $\GL(V)$).
If the bundle $E$ is a cocycle for the (left) group action of $A$ on $Q$, then $P_E$ inherits this left group action, which commutes with the right action of $\GL(V)$.

Suppose now that we have a subgroup $H\subset \GL(V)$.
Then a \emph{reduction of the structure group} of the bundle from $\GL(V)$ to $H$ is a choice of a principal $H$-subbundle $P_H\subset P_E$, which again can be measurable, continuous, smooth, etc.
Equivalently, a reduction of the structure group to $H$ is a section of the bundle $P_E/H$.
This is a \emph{reduction of the cocycle} if the principal bundle is $A$-invariant, i.e. $a(P_H(q))=P_H(aq),\forall a\in A$ and appropriate conditions on $q$, such as $\mu$-a.e. or all.
For instance, if $E$ carried a symplectic form and $V$ was symplectic, then a reduction to $H=\Sp(V)\subset \GL(V)$ would be given by those linear isomorphisms $V \to E(q)$ which respect the symplectic pairing.

Finally, suppose we have a principal $H$-bundle $P\to Q$, equipped with a left action of $A$ lifted from the base $Q$ (and commuting with the right action of $H$).
For any space $F$ with a left $H$-action we can form
\[
	\index{$F_P$}F_P = P\times_H F := \rightquot{P\times F}{(p,f)\sim (ph,h^{-1}f)}
\]
which is now an $F$-bundle over $Q$, with a natural left $A$-action.
The bundles of flags in the fibers of a linear bundle, spaces of measures, etc. can be obtained this way (though the more direct standard constructions are also straightforward).

Note that the fibers of $F_P$ are not canonically identified with $F$, the identification depending on the choice of point in the fiber of $P$.
Any two identifications differ by an element of $H$, and after such identifications the action maps $F_P(q)\to F_P{(aq)}$ are in $H$ (a statement which is independent of the choice of identifications).
Essential for later use will be the natural map\index{$\pi_F$}
\begin{align}
	\label{eq:principal_bundle_to_orbit_space}
	F_P \xrightarrow{\pi_F} \leftquot{H}{F}
\end{align}
which takes a point in the fiber to its equivalence class under the action of $H$.
By construction the map is ``constant on $H$-orbits'', i.e. choosing an identification of the fiber $F_P(q)\toisom F$ (or equivalently a basepoint in the fiber of $P$) gives a (right) action of $H$ on the fiber and for any $h\in H, f\in F_P(q)$ we have $\pi_{F}(hf)=\pi_{F}(f)$.

An important consequence is that the map is constant on the orbits of the dynamics, i.e.
\begin{align}
	\label{eq:map_to_orbit_space_dynamics_invariant}
	\pi_{F}(af)=\pi_{F}(f)  \text{ for any }a\in A.
\end{align}
Note that to define the map $\pi_F$ constructed in \autoref{eq:principal_bundle_to_orbit_space}, we could have avoided the introduction of principal bundles, etc. and simply use local trivializations and check that the map is independent of such trivializations.

\subsubsection{Separability of orbit space}
	\label{sssec:separability_of_orbit_space}
The fundamental fact that will be used frequently is the following (see \cite[Thm.~3.1.3]{Zimmer1984_Ergodic-theory-and-semisimple-groups}).
When $H$ is the real points of an algebraic group, and $F$ is the real points of an algebraic variety with an algebraic action of the group corresponding to $H$, the orbits of $H$ on $F$ are locally closed in the usual (analytic) topology.
In particular the quotient space $\leftquot{H}{F}$ satisfies the separability axiom $T_0$, i.e. for any two points, there exists an open set containing one but not the other.

\subsubsection{Algebraic hull of cocycle}
	\label{sssec:algebraic_hull_of_cocycle}
For more details on this see \cite[\S2]{EskinFilipWright_The-algebraic-hull-of-the-Kontsevich-Zorich-cocycle};
\cite[\S4.2]{Zimmer1984_Ergodic-theory-and-semisimple-groups} contains another discussion of the algebraic hull and \cite[Appendix B]{EskinFilipWright_The-algebraic-hull-of-the-Kontsevich-Zorich-cocycle} shows the notion below coincides with Zimmer's.

Continuing with the cocycle $E\to Q$ over the group action of $A$, assume that there is an $A$-invariant and ergodic measure $\mu$ on $Q$.
We will define now a subgroup \index{$Hull(E,A)(q)$}$Hull(E,A)(q)\subset \GL(E(q))$ for $\mu$-a.e. $q$, the algebraic hull of $E$ for the action of $A$.
The subgroups agree up to conjugation after an identification of the fibers $E(q)\toisom E(q')$.
Given the fixed reference vector space $V$, the algebraic hull can be viewed as a subgroup $H\subset \GL(V)$, well-defined only up to inner conjugation in $\GL(V)$.

Concretely, let \index{$T$@$\bbT(E)$}$\bbT(E)$ denote some tensor construction on $E$, e.g. $E^{\otimes n}$, or its dual.
An $A$-invariant section $s\colon Q\to \bbT(E)$ is a measurable map such that $s(aq)=as(q)$ for any $a\in A$ and $\mu$-a.e. $q\in Q$ (to obtain a set $Q'\subset Q$ of full $\mu$-measure on which the identity holds for every $a\in A$ see \cite[Prop~B.5]{Zimmer1984_Ergodic-theory-and-semisimple-groups}).
An $A$-invariant subbundle is defined analogously.
The \emph{algebraic hull} is the subgroup of $\GL(E(q))$ which fixes every $A$-invariant subbundle in every tensor representation $\bbT(E)$.
Note that the Noetherian property of algebraic sets implies that only finitely many tensors or subbundles need to be considered, and a theorem of Chevalley implies a single one suffices (again, see \cite[\S2]{EskinFilipWright_The-algebraic-hull-of-the-Kontsevich-Zorich-cocycle} for details and motivation).


It follows from the definitions that if $A_1\subset A$ is a subgroup then $Hull(E,A_1)\subset Hull(E,A)$.
We will be interested in algebraic hulls for the action of $g_t$, the group $P$, and all of $\SL_2\left(\bR\right)$.
We can also speak of measurable, continuous, real analytic, etc. algebraic hulls, depending on the regularity of the required tensors or subspaces.

\begin{proposition}[Algebraic hull and reduction of structure]
	\label{prop:algebraic_hull_and_reduction_of_structure}
	Let $E\to (Q,\mu)$ be an $A$-cocycle, with associated principal bundle $P$ with fiber $P(q)=\Isom(V,E(q))$.
	\begin{enumerate}
		\item Suppose $H\subset \GL(V)$ is an algebraic subgroup such that there exists an $A$-equivariant reduction of the structure group of $E$ to $H$, i.e. there exists an $A$-equivariant section $\sigma\colon Q\to P/H$.
		Then there exists an $A$-equivariant family of subgroups $H(q)\subset \GL(E(q))$, such that for any $H$-invariant subspace in a tensor construction $S\subset \bbT(V)$ there is an associated $A$-equivariant subbundle $S_E\subset \bbT(E)$ such that the fibers $S_E(q)$ are $H(q)$-invariant.

		The dependence of $H(q)$ and $S_E(q)$ on $q$ are at least as good as the dependence of $\sigma(q)$ on $q$.

		\item Conversely, suppose $S_E\subset \bbT(E)$ is an $A$-invariant subbundle in some tensor construction, with stabilizers $H(q)\subset \GL(E(q))$.
		Then there exists a subspace $S\subset \bbT(V)$ with stabilizer $H\subset \GL(V)$, and an $A$-equivariant reduction of the structure group of the cocycle from $\GL(V)$ to $H$.

		If $S_E$ is a measurable subbundle, then the reduction can be chosen measurable and $\mu$-a.e.
		If $S_E$ is continuous or better, then the reduction can be chosen on an $A$-invariant open set of full $\mu$-measure with the same regularity as $S_E$.

		\item
		\label{adapted_trivializations}
		There exists a measurable isomorphism $V\times Q\to E$ of the trivial vector bundle with $E$, such that the action of the cocycle is $\mu$-a.e. valued in $H\subset \GL(V)$ and $H$ is the algebraic hull of the trivialized cocycle.

		Equivalently, if we pull back under the trivialization any subspace in a tensor construction on $V$ which is invariant under $H$, the corresponding subbundle in the tensor construction on $E$ is invariant under the $A$-dynamics.
	\end{enumerate}
\end{proposition}

\begin{remark}[On trivializations of bundles]
	\label{rmk:on_trivializations_of_bundles}
	 A trivialization as in \autoref{prop:algebraic_hull_and_reduction_of_structure}\ref{adapted_trivializations} will be called a \emph{trivialization adapted to the algebraic hull}.

	 The separability of the space of orbits of algebraic group actions, as stated in \autoref{sssec:separability_of_orbit_space}, is used only in part (ii) of \autoref{prop:algebraic_hull_and_reduction_of_structure}.
	 It is the mechanism that allows us to go from invariant subbundles to reductions of the structure group of the cocycle.
\end{remark}

\begin{proof}[Proof of \autoref{prop:algebraic_hull_and_reduction_of_structure}]
	For (i), suppose given the $A$-equivariant section $\sigma\colon Q\to P/H$.
	Then the groups $H(q)$ are defined as the stabilizers of $\sigma(q)$ inside $\GL(E(q))$ acting on $P(q)/H$.
	Let also $P_H\subset P$ denote the principal $H$-bundle.
	Then given any tensor construction $\bbT(V)$ with $H$-invariant subspace $S\subset \bbT(V)$, take the associated vector bundles $S_E:=P_H \times_H S\subset P_H\times_H \bbT(V)=\bbT(E)$ to obtain an $A$-invariant subbundle.
	The fibers $S_E(q)\subset \bbT(E(q))$ will be stabilized by $H(q)\subset \GL(E(q))$.

	For the converse construction (ii), given an $A$-invariant subbundle $S_E\subset \bbT(E)$, by the construction in \autoref{eq:principal_bundle_to_orbit_space} the map sending the fiber $S_E(q)$ to its equivalence class in $\leftquot{\GL(V)}{\Gr(\bbT(V))}$ (where $\Gr$ denotes the appropriate Grassmannian) takes the same value $\mu$-a.e. or, in the case $S_E(q)$ depends continuously on $q$, to the same value on an $A$-invariant open set of full measure.
	Let in either case the set be denoted $U\subset Q$.

	By construction there exists a corresponding subspace $S\subset \bbT(V)$ such that the isomorphisms $V\to E(q)$ taking $S_E(q)$ to $S$ are nonempty for every $q\in U$.
	The left action of $H(q)$ (the stabilizer of $S_E(q)$ inside $\GL(E(q))$) on these isomorphisms is transitive, as is the right action of $H$ (the stabilizer of $S$ inside $\GL(V)$).
	This family of isomorphisms $P_H(q)\subset P(q)$ gives the desired reduction of the structure group of the cocycle to $H$.

	Continuing to (iii), let $H(q)$ be the algebraic hull of the cocycle $E$, and let $L\subset \bbT(E)$ be the $A$-equivariant line bundle that defined $H$ as its stabilizer (the existence of $L$ is guaranteed by Chevalley's theorem).
	Let $l\subset V$ be the line provided by the construction in part (ii).
	Then the set of isomorphisms $V\to E(q)$ taking $l$ to $L(q)$ is non-empty $\mu$-a.e., hence it has a measurable section (by standard measurable selection theorems).
	This measurable section is the required choice of trivialization.
\end{proof}

\begin{corollary}[Oseledets decompositions are conjugate]
	\label{cor:oseledets_decompositions_are_conjugate}
	Suppose $E\to (Q,\mu)$ is a cocycle over $g_t$ and $E=\oplus E^{\lambda_i}$ is its decomposition into Lyapunov subspaces.
	Then the algebraic hull $Hull(E,g_t)$ of $E$ respects this decomposition.

	On a set of full measure $Q_0$, for any $q,q'\in Q_0$ there exists a linear isomorphism $E(q)\to E(q')$ respecting the decomposition and also respecting any other $g_t$-invariant subspace in any other tensor construction on $E$.
	Furthermore the isomorphisms can be chosen to vary measurably in both $q$ and $q'$.
\end{corollary}
\begin{proof}
	The first part follows from the definition of algebraic hull.
	The second part follows by choosing a trivialization adapted to the algebraic hull to identify $E(q),E(q')$ to $V$, respecting all subspaces in all tensor powers that are invariant under the algebraic hull.
\end{proof}



\subsection{Holonomies and algebraic hulls}
	\label{ssec:holonomies_and_algebraic_hulls}

\subsubsection{Setup}
	\label{sssec:setup_holonomies_and_algebraic_hulls}
The proof of existence of holonomies in the uniformly hyperbolic setting is in \cite[\S3]{AvilaSantamariaViana2013_Holonomy-invariance:-rough-regularity-and-applications-to-Lyapunov-exponents}.
We develop some additional properties and estimates in \autoref{prop:holonomy_estimates_and_trivialization}.

\subsubsection{Holonomies}
	\label{sssec:holonomies_discussion}
Continuing with the earlier notation, suppose in addition that we have stable manifolds $\cW^{s}[q]\subset Q$ for the action of a $1$-parameter subgroup $g_t$.
The essential property is that for $\mu$-a.e. $q$ and $q'\in \cW^s[q]$ we have $\limsup_{t\to +\infty}\frac{1}{t}\log \dist(g_tq,g_tq')<0$.

We will say that the cocycle $E$ has \emph{stable holonomies} if there exists a set of full $\mu$-measure $Q_0$ such that if $q,q'\in Q_0$ and $q'\in \cW^s[q]$ we have operators\index{$P$@$\cP^{-}(q,q')$}
\[
	\cP^{-}(q,q')\colon E(q)\to E(q') \text{ with } \cP^{-}(q,q'')=\cP^{-}(q'',q')\circ \cP^{-}(q,q') \text{ whenever defined}
\]
and additionally
\[
	\cP^-(g_tq,g_tq')\circ g_t= g_t\circ \cP^-(q,q') \quad \text{ as maps }E(q)\to E(q').
\]
We assume that $\cP^-(q,q)$ is the identity.
The stable holonomies $\cP^-$ will be called measurable, smooth, etc. according to the dependence of $\cP(q,q')$ on $q'\in \cW^s[q]$ (provided the bundle $E$ has at least the regularity needed to speak of this).

\subsubsection{Standard measurable connection, smooth case}
	\label{sssec:standard_measurable_connection_smooth_along_stables}
If the cocycle $E$ is smooth along the stable foliation $\cW^s$, and it has a single Lyapunov exponent, then it has a standard stable holonomy which will be denoted \index{$H^-$}$H^-$ (see \autoref{ssec:cocycle_normal_forms}).

For a cocycle $E$ which is smooth along the stable foliation, but does not necessarily have a single Lyapunov exponent, let $E=\oplus E^{\lambda_i}$ be its decomposition into Lyapunov subspaces.
Then the bundles \index{$E^{\leq \lambda_i}$}$E^{\leq \lambda_i}:=\oplus_{\lambda\leq \lambda_i}E^{\lambda}$ vary smoothly on the stable \cite[\S5.2(b)]{Ruelle1979_Ergodic-theory-of-differentiable-dynamical-systems} and so do their quotients $E^{\leq \lambda_i}/E^{\leq \lambda_{i+1}}$ which have a single Lyapunov exponent $\lambda_i$ and therefore standard holonomies $H^{-}_i$.
These can be assembled into a \emph{standard holonomy}
\index{$H^-$}$H^-$ on
$\gr_{\bullet}E$, where $\gr_{\bullet} E$ is defined in \autoref{sssec:holonomy_construction}.
Using the splitting of $E$ given by the Lyapunov subspaces, which gives us an isomorphism
\[
	E = \bigoplus_{\lambda_i} E^{\lambda_i} 
	\toisom 
	\gr_\bullet^- E
	= \bigoplus_{\lambda_i} E^{\leq \lambda_i}/E^{\leq \lambda_{i+1}}
\]
we define the \emph{standard measurable connection of $E$} to be \index{$P^-$}$P^-:=\oplus P^{-}_i$, where each $P^-_i$ are obtained from $H^-_i$ by the above measurable isomorphism.
Note that the map \index{$P^-(q,q')$}$P^-(q,q')$ is only defined for $q,q'$ in a set of full measure $Q_0$ such that moreover $q'\in \cW^s[q]$.

The unstable measurable connection \index{$P^+(q,q')$}$P^+(q,q')$ is defined analogously, for a cocycle smooth along unstables and for $q,q'$ in a set of full measure $Q_0$ such that $q'\in \cW^u[q]$.

\begin{theorem}[Ledrappier invariance principle]
	\label{thm:ledrappier_invariance_principle_smooth_along_stables}
	Let $E\to Q$ be a cocycle over $g_t$, smooth along stables, with standard measurable connection \index{$P^-_E$}$P^-_E$.
	\begin{enumerate}
		\item \label{item:p_minus_invariance_of_bundles}
		If $S\subset E$ is a measurable $g_t$-invariant subbundle, then $S$ is preserved by $P^-_E$.
		\item 
		If $\bbT(E)$ denotes a tensor construction applied to $E$, then its standard measurable connection $P^-_{\bbT(E)}$ agrees with the measurable connection induced from $P^-_E$ by the tensor construction.
		\item 
		If we pick a measurable trivialization $V\times Q\to E$ adapted to the algebraic hull, then $P^-$ is valued in the algebraic hull $H\subset \GL(V)$.
	\end{enumerate}
\end{theorem}
\begin{proof}
	Part (i) is the original Ledrappier invariance principle, see 
	\cite{Ledrappier_Positivity-of-the-exponent-for-stationary-sequences-of-matrices,AvilaViana_Extremal-Lyapunov-exponents:-an-invariance-principle-and-applications,EskinMirzakhani_Invariant-and-stationary-measures-for-the-rm-SL2Bbb-R-action-on-moduli-space}.
	
	Part (ii) follows because tensor constructions are compatible with the Oseledets decomposition, so in particular the cocycles $\gr_{\bullet}\bbT(E)$ and $\bbT\left(\gr_{\bullet}E\right)$ are isomorphic, via identifications with $\bbT(E)$.
	It then follows that the induced $P^-$ maps also agree.
	
	Part (iii) follows from parts (i) and (ii) applied to the $g_t$-invariant subbundle $L$ in a tensor construction $\bbT(E)$ defining the algebraic hull.	
\end{proof}

We can enlarge to the class of cocycles to which the above constructions and results apply.
For the next definition, we recall that smooth and natural cocycles are defined in \autoref{def:smooth_and_natural_cocycles}; these cocycles are derived from the tangent and higher jet cocycles by elementary linear algebra operations.

\begin{definition}[{Admissible cocycles}]\label{def:admissiblecocycle}
A measurable cocycle $E\to Q$ is \emph{admissible} if there is a smooth and natural cocycle  $Big(E)$, called the \emph{smooth  envelope of $E$}, two $g_t$-equivariant measurable subcocycles $$S_{1}\subset S_2\subset Big(E)$$ and a measurable isometry $E\simeq S_2/S_1$.
\end{definition}



\subsubsection{Standard measurable connection}
	\label{sssec:standard_measurable_connection}
Suppose now that $E\to Q$ is admissible and $E\isom S_2/S_1$ with $S_i\subset Big(E)$.
By \autoref{thm:ledrappier_invariance_principle_smooth_along_stables}, we have that $P^-_{Big(E)}$ preserves $S_1,S_2$ hence induces a measurable connection \index{$P^-_{E}$}$P^-_{E}$ on $E$.

With this construction, we observe that \autoref{thm:ledrappier_invariance_principle_smooth_along_stables} extends to an arbitrary admissible cocycle:

\begin{corollary}[Ledrappier Invariance for admissible cocycles]
	\label{thm:ledrappier_invariance_principle}
	For any admissible cocycle $E$, the conclusions of \autoref{thm:ledrappier_invariance_principle_smooth_along_stables} hold.
\end{corollary}
The result follows from the original one using \autoref{thm:ledrappier_invariance_principle_smooth_along_stables}\autoref{item:p_minus_invariance_of_bundles}.

\subsubsection{Admissible metrics on admissible cocycles}
	\label{sssec:admissible_metrics_on_admissible_cocycles}
Given any cocycle $F$ with a metric to which the Oseledets theorem applies, we can define a new metric on $F$ which agrees with the old metric on each Lyapunov subspace, and for which the Lyapunov decomposition is orthogonal.
As a consequence of the Oseledets theorem, there exist tempered
(see \autoref{def:temperedness_slow_variation}) functions such that the two metrics are comparable up to multiplicative scalars provided by the functions (in one direction the constant $\dim F$ suffices).


We call metrics $F$ obtained on our cocycles by the above procedure \emph{admissible metrics}.
Note that even if $F$ is smooth, the admissible metrics are not in general smooth.
For the class of admissible metrics, given an admissible cocycle $F$ and its smooth envelope $E$, there exists an isometric map $F\into E$ which is not in general $g_t$-equivariant but which induces the identity on $\gr_{\lambda_i}F$ and $S_{2}^{\lambda_i}/S_1^{\lambda_i}$ where $\lambda_i$ are Lyapunov exponents.




\subsubsection{Compatibility of holonomies and standard measurable connection}
	\label{sssec:compatibility_of_holonomies_and_standard_measurable_connection}
Suppose that $E$ is an admissible cocycle, in particular it admits the standard measurable connection $P^{-}_E$ constructed in \autoref{sssec:standard_measurable_connection} above.
Suppose also that $E$ carries holonomies $\cP^{-}_E$ on the stable manifolds, or some equivariant subsets of the stable manifolds.
Suppose that $\cP^-_E$ is \Holder-continuous on sets of almost full measure, i.e. there exists an exponent $\alpha>0$ such that for any $\ve>0$ there exists a compact set $K_\ve$ of measure at least $1-\ve$, and a constant $C_\ve<+\infty$, such that $E$ admits \Holder trivializations on charts covering $K$, and such that $\cP^{-}_E(q,q')$ is \Holder-continuous on $K_\ve$, with constant $C_\ve$ and exponent $\alpha$.

With these assumptions, it is a standard fact that $\cP^{-}_E$ preserves the Lyapunov flags on $E$ along $\cW^s[q]$.
The measurable connection $P^{-}_E$ does so as well, and the induced maps on the associated graded of $E$ agree.
This last assertion is equivalent to the uniqueness of holonomies for cocycles with one Lyapunov exponent, and a \Holder structure on sets of almost full measure, which is a standard fact.

\subsubsection{Summary on holonomies, notation}
	\label{sssec:summary_on_holonomies_notation}
\begin{description}
	\item[Holonomy on graded] If $E$ is smooth along stables, then we have the holonomy maps \index{$H^-(q,q')$}$H^-(q,q')$ on the associated graded $\gr_\bullet E$.
	These depend smoothly on $q'\in \cW^s[q]$.
	\item[Measurable connection] If $E$ is any admissible cocycle, the measurable connections \index{$P^-$}$P^-$ are provided by Ledrappier's \autoref{thm:ledrappier_invariance_principle}, and analogously for \index{$P^+$}$P^+$.
	\item[Gauss--Manin] If $E$ is smooth along stables, then we have the flat connections \index{$P^{-}_{GM}$}$P^{-}_{GM}$ on the linearization \index{$L^-E$}$L^-E$, provided by normal forms as in \autoref{cor:holonomy_linearization_of_arbitrary_cocycles}.
	Analogous statements hold for unstables and \index{$P^+_{GM}$}$P^+_{GM}$.
\end{description}



\subsection{Jordan Normal Forms for Cocycles}
	\label{ssec:jordan_normal_forms_for_cocycles}


\subsubsection{Zimmer's amenable reduction}
\def\reals{\R}
The following general facts about linear cocycles over an action of $\reals$ are often called ``Zimmer's amenable reduction''.
We state a version adapted to our setup.

Suppose $E\to (Q,\nu)$ is a measurable linear cocycle over the $g_t$-action.
\begin{definition}[Jordan Normal Form of Cocycle]
	\label{def:Jordan_normal_form}
	The measurable linear cocycle $E\to (Q,\nu)$  \emph{can be put in Jordan normal form} if there exists a filtration by subcocycles:
	\begin{equation}
	\label{eq:jordanNF:flag}
	\{0\} = E_{0}(q) \subset E_{1}(q) \subset \dots \subset
	E_{p}(q) = E(q)
	\end{equation}
	and a measurable family of metrics on the quotient cocycles $E_{i}/E_{i-1}$ such that the induced action is by conformal transformations.

	We call a filtration $\gF= E_\bullet$ as in \eqref{eq:jordanNF:flag}, with its metric, a \emph{Jordan flag} and say the Jordan flag  \eqref{eq:jordanNF:flag} \emph{puts the cocycle in Jordan normal form}.
	The \emph{length} of a Jordan flag \eqref{eq:jordanNF:flag} is the number of subspaces $p$ appearing the filtration in \eqref{eq:jordanNF:flag}.   
\end{definition}

\begin{proposition}[Jordan Normal Form and Standard Operations]
	\label{prop:jordan_normal_form_and_standard_operations}
	Suppose that $E\to (Q,\nu)$ and $F\to (Q,\nu)$ are measurable linear cocycles. 
	\begin{enumerate}
		\item If $F$ is a subcocycle of $E$ and if $E$ can be put in Jordan normal form, then $F$ and $E/F$ can be put in Jordan normal form.  
		\item If $F$ is a subcocycle of $E$ and if both the quotient $E/F$ and $F$ can be put in Jordan normal form, then $E$ can be put in Jordan normal form.
		\item If $E$ and $F$ can be put in Jordan normal form then the cocycle $E\otimes F\to (Q,\nu)$ can be put in Jordan normal form.  
	\end{enumerate}
\end{proposition}
\begin{proof}
	For (i), suppose $E$ can be put in Jordan normal form, with subcocycles $E_i$.
	Then certainly the subcocycles $E_i\cap F$ put $F$ in Jordan normal form, and their images in the quotient $E/F$ will put $E/F$ in Jordan normal form.

	For (ii), suppose $F$ and $E/F$ can be put in Jordan normal form by $F_i$ and $(E/F)_j$ respectively.
	Then we can use the sequence $F_i$, followed by the preimages of $(E/F)_j$, to put $E$ in Jordan normal form.

	Finally, for (iii), suppose that $E_i,F_j$ put $E,F$ in Jordan Normal Form.
	Choose any total ordering of the indexes $(i,j)$ compatible with the partial order $(i,j)\prec (i',j')$ if $i\leq i'$ and $j\leq j'$.
	Define then $(E\otimes F)_{k}$ to be the span of all the $E_i\otimes F_j$ for the first $k$ pairs of $(i,j)$ in the selected total order (the compatibility with the partial order ensures that $(E\otimes F)_{k-1}\subset (E\otimes F)_{k}$).
	Then the successive quotients of this filtration are of the form $(E_i/E_{i-1})\otimes (F_j/F_{j-1})$ and so the filtration puts the cocycle in Jordan normal form.
\end{proof}

Given a cocycle  $E\to (Q,\nu)$ be a cocycle over a measurable flow and an equivariant finite measurable cover $(\hat Q,\hat \nu) \to (Q,\nu)$, we can always consider the pulled-back cocycle $   E\to (\hat Q,\hat \nu)$ on $\hat Q$.

\begin{theorem}[Jordan Normal Form]
\label{thm:jordan_normal_form}
	Given any cocycle $E\to (Q,\nu)$, there exists an equivariant finite measurable cover $\sigma\colon (\hat Q,\hat \nu)\to (Q,\nu)$ and a filtration \index{$F$@$\gF_{\bullet}$}$\frakF_{\bullet}$ of $E$ that puts the cocycle in Jordan Normal Form, with the following additional properties:
	\begin{enumerate}
		\item The associated graded cocycle $\gr^\frakF E$, and the metric that puts it in Jordan normal form, is pulled back from $(Q,\nu)$.
		\item The norm of the linear cocycle for the new metric is bounded above by the norm for the original metric.
\item \label{Jordanextention:iii} Suppose that $E$ is an admissible cocycle.   
		Then for $\nu$-a.e.\ $q$ and every biregular $q'\in \cW^s[q]$, the (finite) set of Jordan flags $\frakF_{\bullet}E(q)$ and $\frakF_{\bullet}E(q')$ 
		is intertwined by $P^-(q,q')$.

For $\nu$-a.e.\ $q$ and a Jordan flag $\frakF$ on  $E(q)$ we write $x=(q,\frakF)\in \hat{Q}$  and 
define the measurable set\index{$W$@$\cW^s[x]$}\index{$W$@$\cW^s[q,\frakF]$}
		\[
			\cW^s[x]=\cW^s[q,\frakF]:=\{(q',\frakF')\colon q'\in \cW^s[q]\text{ and }\frakF'=P^{-}(q,q')\frakF\}\subset\hat{Q}.
		\]
		Then the projection map identifies $\cW^s[q,\frakF]\toisom \cW^s[q]$.  

\item\label{Jordanextention:iiii} Suppose $E$ is smooth along stable manifolds.  Then for $\nu$-a.e.\ $q$ and every $q'\in \cW^s[q]$, there exists a finite collection of flags $\frakF_\bullet E(q')$ with the same combinatorics as \index{$F$@$\frakF_\bullet E(q)$}$\frakF_\bullet E(q)$  and  the collection $q'\mapsto  \frakF_\bullet E(q')$ varies smoothly along $\cW^s[q]$.  Moreover, there exists a smooth family of linear maps \index{$P$@$\widetilde P(q, q')$}$q'\mapsto \widetilde P(q, q')$ such that $\widetilde P(q,q')\frakF_\bullet E(q)= \frakF_\bullet E(q')$ and $\widetilde P(q,q')$ restricts to $P^-(q,q')\frakF_\bullet E(q)= \frakF_\bullet E(q')$ if $q$ and $q'$ are biregular.  


		Furthermore,  for $\nu$-a.e.\ $q\in Q$ and $\frakF$ a Jordan flag on $E(q)$, we define  $x=(q,\frakF)\in \hat{Q}$ and 
		\[
			\cW^s[x]=\cW^s[q,\frakF] :=\{(q',\frakF')\colon q'\in \cW^s[q]\text{ and }\frakF'=\wtd P^{}(q,q')\frakF\}\subset\hat{Q}.
		\]
		Then the projection map identifies $\cW^s[q,\frakF]\toisom \cW^s[q]$ and if two elements belong to the same $\cW^s[q,\frakF]$, then they are exponentially asymptotic in forward time.  Moreover, the restriction of $\cW^s[q,\frakF]$ to fibers over biregular points $q'\in \cW^s[q]$ coincides with the measurable set in part \ref{Jordanextention:iii}. 

		\item \label{Jordanextention:iiiii}Any other finite cover $(Q',\nu')\to (Q,\nu)$ on which $E$ can be put in Jordan Normal Form by a filtration $\frakF'_{\bullet}$ factors via $(Q',\nu')\to (\hat Q,\hat \nu)$, and the associated graded $\gr^{\frakF'}_{\bullet}E$ refines the pullback of $\gr^{\frakF}_\bullet E$.

		This last property is universal and characterizes $(\hat Q,\hat\nu)$ and the filtration $\frakF_\bullet$.
	\end{enumerate}
\end{theorem}

 



\subsubsection{Preliminary constructions for the proof of \autoref{thm:jordan_normal_form}}
	\label{sssec:preliminary_constructions_for_the_proof_of_lemmas_on_zimmer_flags}
Let \index{$G$@$\bbG_E$}$\bbG_E$ denote the algebraic hull of the cocycle $E$ for the $g_t$-action on $(Q,\nu)$.
Recall that this is an a.e. defined family of algebraic subgroups of $\GL(E(q))$, for $q$ in a set of full $\nu$-measure.
Since the acting group, i.e. $\bR$, is amenable, the algebraic hull is also amenable by Zimmer's amenable reduction theorem \cite[9.2.3]{Zimmer1984_Ergodic-theory-and-semisimple-groups}.
In particular, it has a unipotent radical \index{$G$@$\bbG_E^u$}$\bbG_E^u\subset \bbG_E$ such that the quotient is a product of compact groups and split tori.
The last statement follows from \cite[4.1.9]{Zimmer1984_Ergodic-theory-and-semisimple-groups} which characterizes amenable Lie groups as compact extensions of solvable Lie groups.

\subsubsection{Engel filtration}
	\label{sssec:engel_filtration}
Any unipotent subgroup $\bbG_E^u(q)\subset \GL(E(q))$ induces a canonical filtration of $E(q)$: the lowest term consists of vectors fixed by the full group (nonempty by Engel's theorem), and subsequent terms of the filtration are defined inductively by taking the quotient by the preceding term.
Let us call this the \emph{Engel filtration} which we denote by $$ \{0\} =V_0(q)\subsetneq V_1(q)\subsetneq \dots \subsetneq V_\ell(q)= E(q)$$




Let $\bbG_E^{\circ}$ denote the connected component of the identity of $\bbG_E$, the algebraic hull of $E$. 
We consider the image of $\bbG_E$ in $\GL(V_{j+1}(q)/V_{j}(q))$ of the Engel filtration and denote it by $\bbG_j$.
This is reductive, since it has no unipotent radical, and we denote by $\bbG_j^{\circ}$ its Zariski-connected component of the identity.
The solvable radical $\bbT_j^{sol}\subset \bbG_j^{\circ}$ is an $\bR$-torus, which has a canonical\footnote{Tori over $\bR$ are equivalent to the data of their cocharacter lattice (obtained by extending scalars to $\bC$) and an involution on it coming from complex conjugation; the split part corresponds to the sublattice fixed by the involution.} split subtorus $\bbT_j\subseteq \bbT_j^{sol}$.
Note that $\bbT_j$ is a normal subgroup of $\bbG_j$ and the quotient is compact.

\subsubsection{Zimmer flags}
	\label{sssec:zimmer_flags}
We now decompose $V_{j+1}/V_j$ according to the characters $\bbT^{*}_j$ of $\bbT_j$:
\begin{equation}
	\label{eq:split_torus_character_splitting}
	V_{j+1}/V_{j}=\bigoplus_{\alpha\in \bbT^*_j} V_{j,\alpha}(q)
\end{equation}
A \emph{Zimmer flag} at $q$ is any flag in $E(q)$ comprised of subspace of the form
\[
	V_{j-1}\oplus \bigoplus _{i=1}^{k} V_{j,\alpha_{i}}(q) \quad k=0,1,\ldots, k_j
\]
where $\alpha_1,\ldots,\alpha_{k_j}$ is allowed to be any enumeration of the nonvacuous characters.


\subsubsection{Constructing the conformal metrics}
	\label{sssec:constructing_the_conformal_metrics}
The group $\bbG_j$ has finitely many Zariski-connected components denoted $F_j$, and a Zariski-connected component of the identity denoted $\bbG_j^{\circ}$.
Recall also that $\bbT_j$ is normal and the quotient $\bbG_j/\bbT_j$ is compact.

We now return to the decomposition in \autoref{eq:split_torus_character_splitting}.
Note that the finite group $F_j$ can permute the $V_{j,\alpha}$'s with different $\alpha$'s.
Let $M_j$ denote the space of metrics on $V_{j+1}/V_j$ for which the decomposition $\oplus V_{j,\alpha}$ is orthogonal.
This is equipped with a natural metric making it a $\CAT(0)$-space.
Then the quotient $M_j/\bbT_j$ is still such a space, isomorphic to $\prod \PSL(V_{j,\alpha})/K(V_{j,\alpha})$, where $K(V_{j,\alpha})$ is a maximal compact of $\PSL(V_{j,\alpha})$.
Since the quotient group $\bbG_j/\bbT_j$ is compact and acts on $M_j/\bbT_j$, so its fixed point set is a nonempty geodesically convex subset $M_{j,f}\subset M_j/\bbT_j$.

\subsubsection{Equivariance properties}
	\label{sssec:equivariance_properties}
The above constructions were canonical and apply to a set of full measure of the fibers of $E\to(Q,\nu)$.
We now consider what happens over the entire dynamical system, and indicate the dependence on $q\in Q$.

The algebraic hull acts trivially on $M_{j,f}(q)$ by construction, so in fact there exists a single ``model'' space $M_{j,f}^{mod}$ and a $g_t$-equivariant map $M_{j,f}(q)\to M_{j,f}^{mod}$, as in \autoref{eq:principal_bundle_to_orbit_space} and as explained in the discussion preceding it.
Select one point $m_{j,f}\in M_{j,f}^{mod}$, let $m_{j,f}(q)\in M_{j,f}(q)$ be the corresponding equivariant family of spaces, and its preimage $M_{j,c}(q)\subset M(q)$ is also a geodesically convex subset on which the algebraic hull acts via $\bbT_{j}$.
The space $M_j(q)$ is a geodesically convex subset of $M^{full}(V_{j+1}/V_j(q))$ of all metrics on the fiber $V_{j+1}/V_j(q)$, and so is $M_{j,c}(q)$.
Let $\pi_0(q)\colon M^{full}(V_{j+1}/V_j(q))\to M_{j,c}(q)$ be the nearest point projection.
Note that $\pi_{0}(q)$ is distance nonincreasing, which will imply that the norms of the cocycle for the new metrics are not larger than for the original metrics. 
All objects and maps defined so far are equivariant for the dynamics.

Now $V_{j+1}/V_j(q)$ has its original cocycle metric $m_V(q)\in M^{full}(V(q))$, and we set $m_c(q):=\pi_0(q)(m_V(q))\in M_{j,c}(q)$ to be its projection.
Then the action of the dynamics $g_t$ takes $m_c(g_{-t}q)$ to the same $\bbT^s$-orbit as $m_c(q)$, and hence the action is conformal on each $V_{\alpha}$.

\subsubsection{The proof of \autoref{thm:jordan_normal_form}}
	\label{sssec:the_proof_of_thm:jordan_normal_form}
The collection of flags \index{$F$@$\frakF_\bullet E$}$\frakF_\bullet E$ has been defined in \autoref{sssec:zimmer_flags}.
This is a collection of flags on each $E(q)$ with the same associates graded pieces, up to permutation.

We endow the associated graded with the metric constructed in \autoref{sssec:constructing_the_conformal_metrics}.
The property that the dynamics is by conformal transformations, that the metric is induced already on $(Q,\nu)$, and the boundedness properties follow from the discussion in \autoref{sssec:constructing_the_conformal_metrics} and \autoref{sssec:equivariance_properties}.

The universality of the construction follows since by passing to a finite cover, the algebraic hull can only change by finite index, so any filtration of the cocycle on $(Q',\nu')$ will have to refine the Engel filtration and character pieces defined in \autoref{eq:split_torus_character_splitting}.

The $P^-$-equivariance of the Zimmer flags follows from the Ledrappier Invariance Principle, \autoref{thm:ledrappier_invariance_principle}, applied to an appropriate tensor construction on $E$ that has invariant subbundles corresponding to finite collections of Zimmer flags.
The smoothness of the flags $\frakF_\bullet E$ along stables follows from the fact that they refine the stable flags (which are smooth along stables), the fact that the holonomy on $\gr_{\bullet} E$ (which is smooth along stables) agrees with $P^-$, and that by Ledrappier again $P^-$ preserves $\frakF_\bullet E$.
\hfill \qed

\begin{remark}[Points in the cover and flags]
	\label{rmk:points_in_the_cover_and_flags}
	Let $E=\oplus E^{\lambda_i}$ be the Lyapunov splitting of the cocycle $E\to Q$.
	A point in the finite cover $\hat{Q} \to Q$ can be interpreted (equivalently) as the following additional data:
	\begin{itemize}
		\item A flag in each Lyapunov subspace $E^{\lambda_i}$.
		\item A flag of $E$, refining the stable flags $E^{\leq \lambda_i}$.
		\item A flag of $E$, refining the unstable flags $E^{\geq \lambda_i}$.
	\end{itemize}
	Their associated graded pieces are naturally identified.
	Each of these flags has additional properties listed in \autoref{thm:jordan_normal_form}, and the flag on the stable has smoothness properties on the stable, and analogously for the unstable. 
\end{remark}


\subsection{Good Norms}
\label{ssec:good:norms}

The next result is an analogue of \cite[Prop.~4.15]{EskinMirzakhani_Invariant-and-stationary-measures-for-the-rm-SL2Bbb-R-action-on-moduli-space} and its proof is analogous.
Recall that the bundle \index{$H$@$\bbH$}$\bbH:=L\cC/L\cC^{\tau}$ is introduced in \autoref{eq:bold_H_definition}.
We denote by \index{$H$@$\bbH^{\lambda_i}_{j}$}$\bbH^{\lambda_i}_{j}$ the flags on $\bbH^{\lambda_i}$ obtained from \autoref{thm:jordan_normal_form} applied to $\bbH$.

\begin{proposition}[Good norms]
\label{prop:good_norms}
	On a set of full measure $X_0\subset X$ there exists an inner product \index{$\norm$@$\ip{-,-}_q$}$\ip{-,-}_q$ on $\bbH(q)$ (or on any bundle for which the conclusions of \cite[Lemma~4.14]{EskinMirzakhani_Invariant-and-stationary-measures-for-the-rm-SL2Bbb-R-action-on-moduli-space} hold) with the following properties:
	\begin{enumerate}
		\item For $q\in X_0$, the distinct Lyapunov spaces $\bbH^{\lambda_i}(q)$ are orthogonal.
		\item Let $\bbH^{\lambda_i}_{j'}$ denote the orthogonal complement relative to the inner product $\ip{-,-}$ of $\bbH^{\lambda_i}_{j-1}$ in $\bbH^{\lambda_i}_{j}$.
		Then for $q\in X_0$, any $t\in \bR$ and any $v\in \bbH^{\lambda_i}_{j'}(q)$ we have
		\[
			g_t \bfv = e^{\lambda_{ij}(q;t)}\bfv' + \bfv''
		\]
		where $\lambda_{ij}(q;t)\in \bR$, $\bfv'\in \bbH^{\lambda_i}_{j'}(g_t q)$, $\bfv'' \in \bbH^{\lambda_i}_{j-1}(g_t q)$, and \index{$\norm{\bfv}$}$\norm{\bfv'}=\norm{\bfv}$.
		In particular
		\[
			\norm{g_t \bfv}\geq e^{\lambda_{ij}(q;t)}\norm{\bfv}.
		\]
		\item 
		\label{item:lambda_ij_bilipschitz}
		There exists $\kappa>1$ such that for $q\in X_0$ and any $t>0$:
		\[
			\frac{1}{\kappa}\cdot t \leq \lambda_{ij}(q;t)\leq \kappa \cdot t.
		\]
		\item \label{good_norms_bilipschitz}
		There exists $\kappa>1$ such that for $q\in X_0$, and $\bfv\in \bbH(q)$, and any $t>0$:
		\[
			e^{\frac{1}{\kappa}\cdot t}\norm {v} \leq 
			\norm{g_t \bfv}
			\leq e^{\kappa \cdot t} \norm{\bfv}.
		\]
		\item For $q\in X_0$ and $q'\in X_0\cap \gB_0[q]$ and any $t\leq 0$ we have
		\[
			\lambda_{ij}(q;t) = \lambda_{ij}(q';t)
		\]
		\item For $q\in X_0$ and $q'\in X_0\cap \gB_0[q]$, and any $v,w\in \bbH(q)$
		 we have
		\[
			\ip{P^+(q,q')\bfv,P^+(q,q')\bfw}_{q'} = \ip{\bfv,\bfw}_{q}.
		\]
	\end{enumerate}
\end{proposition}

\begin{proof}
	The proof of this result is analogous to \cite[Prop.~4.15, \S4.10]{EskinMirzakhani_Invariant-and-stationary-measures-for-the-rm-SL2Bbb-R-action-on-moduli-space}, with the following modifications.

	The measurable set $K$ and function $T_0\colon Q\to \bR_{> 0}$ are chosen such that, on $K$ we have uniformity in the Oseledets theorem for $\bbH$ with constants $C(K),\ve(K)>0$ (see \autoref{sssec:lyapunov_norms}), and the measurable connection $P^+$ satisfies \cite[Eqn.~(4.17)]{EskinMirzakhani_Invariant-and-stationary-measures-for-the-rm-SL2Bbb-R-action-on-moduli-space}.
	Then $T_0$ is chosen such that the conditions from \cite[Lemma~4.14]{EskinMirzakhani_Invariant-and-stationary-measures-for-the-rm-SL2Bbb-R-action-on-moduli-space} hold.
	Then, we apply \autoref{lem:partition_adapted_to_a_measurable_function} with this $T_0$ to obtain the measurable partition $\gB_0$ and follow the same constructions as in \cite[\S4.10]{EskinMirzakhani_Invariant-and-stationary-measures-for-the-rm-SL2Bbb-R-action-on-moduli-space}.
\end{proof}

The next result is an elementary consequence of Lusin's theorem:
\begin{lemma}
\label{lemma:hodge:norm:vs:dynamical:norm}
For every $\delta > 0$ there exists a compact subset $K(\delta)
\subset X_0$ with
$\nu(K(\delta)) > 1-\delta$ and a number $C_1(\delta) < \infty$ such
that for all $x \in K(\delta)$ and all $\bfv \in \bH(x)$, 
\begin{displaymath}
C_1(\delta)^{-1} \le \frac{\|\bfv\|_{1}}{\|\bfv\|_{2}} \le C_1(\delta),
\end{displaymath}
where $\|\cdot \|_i$ can be either the Lyapunov norm (see \autoref{sssec:lyapunov_norms}), or the natural norm arising on $\bbH$ as an admissible cocycle, or the the dynamical norm defined in
\autoref{prop:good_norms}.
\end{lemma}



\section{Subresonant algebra and geometry}
	\label{appendix:subresonant_linear_algebra}

\paragraph{Outline of section}
In \autoref{ssec:subresonant_linear_algebra} we start by recalling some basic facts regarding linear maps between vector spaces equipped with (weighted) filtrations.
We proceed to study polynomial maps between such vector spaces, called \emph{subresonant maps}.
In \autoref{ssec:resonant_maps_and_} we consider in more detail decompositions of such maps into \emph{resonant} and \emph{strictly subresonant} maps.

One key construction that we introduce is that of \emph{linearization}, which polynomially embeds one filtered vector space into another, such that the action of polynomial subresonant maps on the source is conjugated to a linear action on the target.
An example that illustrates most of our constructions is in \autoref{sssec:example_of_subresonant_maps}.



\subsubsection*{Ordering of Lyapunov exponents}
	\label{sssec:ordering_of_lyapunov_exponents}
We recall that our convention for indexing Lyapunov exponents is such that
\[
	\lambda_1>\lambda_2>\cdots
\]
with the implicit assumption that $\lambda_1>0$.
When discussing negative Lyapunov exponents, we will write them as
\[
	-\lambda_1 < - \lambda_2 <\cdots
\]
again with the assumption that $\lambda_1>0$.


\subsection{Subresonant linear algebra}
	\label{ssec:subresonant_linear_algebra}

All our vector spaces are finite dimensional, over $\bR$ or $\bC$.

\subsubsection{Filtered vector spaces}
	\label{sssec:filtered_vector_spaces}
By convention, a filtered vector space $V$ will be written as\index{$V^{\leq -\lambda_i}$}
\begin{align}
	\label{eq:filtered_vector_space_notation}
	\{0\}\subsetneq V^{\leq -\lambda_1} \subsetneq V^{\leq -\lambda_2}\subsetneq
	\cdots \subsetneq V^{\leq -\lambda_l} = V
\end{align}
where the ``Lyapunov exponents'' satisfy $-\lambda_1<-\lambda_2<\cdots < -\lambda_l$.
This indexing is convenient when assuming that $\lambda_i>0$; later we will also label the elements of the filtration by arbitrary $\lambda\in \bR$, but the convention is always that $\lambda<\mu$ implies $V^{\leq\lambda}\subseteq V^{\leq \mu}$.
We will also call the $\lambda$ that index the filtration the ``weights of the filtration''.
The following notation will be useful:
\[
	wt(v)\leq \mu \text{ if and only if }v\in V^{\leq \mu}
\]
and more generally\index{$wt(v)$}
\[
	wt(v)=\mu \text{ if and only if }v\in V^{\leq \mu}\setminus V^{<\mu}
\]
where $V^{<\mu}:=\bigcup_{\alpha<\mu}V^{\leq \alpha}$.

\subsubsection{Standard linear algebra operations}
	\label{sssec:standard_linear_algebra_operations}
Given two filtered vector spaces $(V,V^{\leq \bullet})$ and $(W,W^{\leq \bullet})$ a linear map
\[
	T\colon V\to W
\]
is said to have \emph{weight $\leq \mu$} if $T\left(V^{\leq \lambda}\right)\subseteq W^{\leq \lambda + \mu}$.
This induces a filtration on $\Hom(V,W)$.
In particular, giving the scalars the trivial filtration in degree $0$ imposes the filtration on the dual vector space
\[
	V^{\leq \lambda} = \left\lbrace \xi \in V^{\dual}
	\colon \xi\left(V^{<-\lambda}\right)=0
	 \right\rbrace.
\]
With these conventions, if we have a filtered vector space
\begin{align*}
	0\subsetneq V^{\leq -\lambda_1}\subsetneq & \cdots \subsetneq V^{\leq -\lambda_l}=V\\
	\intertext{then its dual is}
	0\subsetneq \left(V^{\dual}\right)^{\leq \lambda_l}
	\subsetneq 
	& \cdots
	\subsetneq \left(V^\dual\right)^{\leq \lambda_1}=V^{\dual}\\
\end{align*}
and in particular the dual of $V^{\leq -\lambda_i}$ is $V^\dual/\left(V^{\dual}\right)^{\leq \lambda_{i+1}}$.

The filtration on $V\otimes W$ is induced in the standard way from that on $V$ and $W$, therefore the tensor product spaces $V^{\otimes n}$ and $\Sym^n(V)$ acquire natural filtrations.
Note also that for a linear map $T\colon V\to W$ and its transpose \index{$T^t$}$T^t\colon W^{\dual}\to V^{\dual}$, the conditions $wt(T)\leq \alpha$ and $wt\left(T^t\right)\leq \alpha$ are equivalent.

\paragraph{Standing assumption:}
From now on, vector spaces $V,W$, etc. are assumed to have strictly negative weights, so that their duals $V^\dual$ and more generally spaces of polynomial functions have positive weights.

\subsubsection{Filtration on polynomial functions}
	\label{sssec:filtration_on_polynomial_functions}
Suppose that $V$ is filtered as in \autoref{eq:filtered_vector_space_notation}.
The ring of polynomial functions on $V$ is naturally\index{$S$@$\Sym^\bullet (V^\dual)$}
\[
	\Sym^\bullet (V^\dual) = \bigoplus_{i\geq 0} \Sym^i\left(V^\dual\right)
\]
and the weights of its filtration are of the form $\sum k_i \lambda_i$ with $k_i\in \bZ_{\geq 0}$ and $-\lambda_i$ are the weights of $V$.

Given a filtered subspace $S\subset V$, with the filtration on $S$ compatible with that on $V$ (for instance, but not necessarily, $S=V^{\leq -\lambda_j}$ for some $j$), we have the restriction of linear functions $V^\dual \onto S^{\dual}$ which is a map of dynamical weight $\leq 0$.
Analogously, the restriction of polynomial functions $\Sym^{\bullet}(V^{\dual})\onto \Sym^{\bullet}(S^\dual)$ is also of dynamical weight $\leq 0$.

We also have the useful property of polynomial functions
\begin{align}
	\label{eq:vanishing_sres_polynomials}
	p\in \Sym^{\bullet}(V^{\dual})^{\leq \lambda}
	\implies
	p(v+v')=p(v)\text{ if } v'\in V^{<-\lambda}
\end{align}
which follows from the analogous statement for linear functions (and is in fact a characterization for linear functions).
Indeed $p$ is a linear combination of products of linear functions of positive weights adding up to $\lambda$, and the property is true of each such linear function.

\subsubsection{Filtration on polynomial maps}
	\label{sssec:filtration_on_polynomial_maps}
Given two vector spaces $V,W$, a polynomial map \index{$P$@$\Poly(V,W)$}$\Poly(V,W)$ is an element of $\Sym^\bullet\left(V^\dual\right)\otimes W$.
In particular, the polynomial maps between filtered vector spaces are themselves filtered.
Note that $\Poly(V,W)$ includes affine maps (that do not take the origin of $V$ to that of $W$), since $\Sym^\bullet(V^{\dual})$ contains $\Sym^0(V^\dual)$, which is just the scalars.

\begin{proposition}[Characterization of filtration on polynomial maps]
	\label{prop:characterization_of_filtration_on_polynomial_maps}
	The following properties are equivalent for a polynomial map $F\in \Poly(V,W)$:\index{$P$@$\Poly(V,W)^{\leq \alpha}$}
	\begin{enumerate}
		\item We have $F\in \Poly(V,W)^{\leq \alpha}$.
		\item For any linear function $\xi\in
                  (W^{\dual})^{\leq \mu}$ we have that its pullback satisfies $F^* \xi\in
                  \Sym^{\bullet}(V^\dual)^{\leq \mu+\alpha}$.
		\item For any $d\geq 1$ and any homogeneous polynomial function $f\in \Sym^{d}(W^{\dual})^{\leq \mu}$ we have that $F^* f\in \Sym^{\bullet}(V^\dual)^{\leq \mu+d\cdot\alpha}$.
	\end{enumerate}
\end{proposition}
\begin{proof}
	It is clear that (iii) implies (ii), and (ii) implies (iii) since any homogeneous polynomial function $f$ is a linear combination of products of linear functions whose weights add up to at most the weight of $f$.

	Write now a polynomial map as $F=\sum p_i\otimes w_i$ with $p_i$ polynomial functions on $V$ and $w_i$ a basis of $W$.
	We can choose the basis $\{w_i\}$ compatible with the filtration, i.e. $wt(w_i)\leq wt(w_{i+1})$.

	To check that (i) implies (ii), suppose $wt(F)\leq \alpha$, i.e. $wt(p_i)+wt(w_i)\leq \alpha, \forall i$.
	The pullback of $\xi\in \left(W^\dual\right)^{\leq \mu}$ is the polynomial function $\sum p_i \cdot \xi(w_i)$.
	Now $\xi(w_i)\neq 0$ implies $wt(w_i)\geq -\mu$ (by definition of the filtration on the dual) and so $wt(p_i)\leq \alpha+\mu$ as desired.

	To check that (ii) implies (i), pick a dual to $w_i$ basis $\xi_i\in \cW^{\dual}$, in particular with the property that $wt(\xi_i)= -wt(w_i)$.
	Applying the assumption in (ii) to each $\xi_i$ gives that $wt(p_i)\leq \alpha - wt(w_i)$ which yields the desired conclusion.
\end{proof}

\begin{definition}[Subresonant polynomial maps]
	\label{def:subresonant_polynomial_maps}
	A polynomial map $F\in \Poly(V,W)$ will be called \emph{subresonant} if it is in non-positive weight, i.e. $F\in \Poly(V,W)^{\leq 0}$, and such maps will also be denoted \index{$P$@$\Poly^{sr}(V,W)$}$\Poly^{sr}(V,W)$.
\end{definition}
For example, linear maps that preserve the filtrations are subresonant.

\begin{proposition}[Properties of subresonant polynomial maps]
	\label{prop:properties_of_subresonant_polynomial_maps}
	Let $U,V,W$ be filtered vector spaces and $V\xrightarrow{F} W$ a subresonant polynomial map.
	\begin{enumerate}
		\item For any subresonant polynomial map $G\in \Poly^{sr}(U,V)$ the composition $F\circ G$ is also subresonant.
		\item For any $x\in V$ a reference point, the differential of $F$ at $x$, viewed as a linear map $D_x F\colon V\to W$, is a linear subresonant map.
		\item Suppose $S\subseteq V$ is a linear subspace, with the induced filtration.
		Then the restricted polynomial map $F\vert_{S}\colon S\to W$ is a subresonant polynomial map.
		\item Suppose $S=V^{\leq -\lambda}\subseteq V$ is an element of the filtration of $V$, or an affine translate thereof.
		Then the image $F(S)$ is contained in $\cW^{\leq -\lambda}$, or an affine translate thereof.
		\item Suppose $S=V^{\leq -\lambda}\subseteq V$ is an element of the filtration of $V$.
		There is an induced subresonant map $[F]$ fitting into the diagram\index{$F$@$[F]$}
		\begin{equation}
			\label{eqn_cd:quotient_subresonant map}
		\begin{tikzcd}
			V
			\arrow[r, "F"]
			\arrow[d]
			& 
			W
			\arrow[d]
			\\
			\rightquot{V}{V^{\leq -\lambda}}
			\arrow[r, "{[F]}"]
			&
			\rightquot{W}{\cW^{\leq - \lambda}}
		\end{tikzcd}
		\end{equation}
	\end{enumerate}
\end{proposition}
\begin{proof}
	Part (i) follows from \autoref{prop:characterization_of_filtration_on_polynomial_maps} and the characterization in terms of pullback of functions.

	For part (ii), it suffices to check it for $x=0\in V$ since translations on the source or target preserve the set of subresonant maps.
	But the pullback of linear functions by $F$ does not increase their weight, only adds higher degree functions of lower weight.
	So the derivative is a subresonant linear map.

	Part (iii) follows because restriction of polynomial functions $\Sym^{\bullet}\left(V^{\dual}\right)\to \Sym^{\bullet}\left(S^\dual\right)$ preserves the filtrations.

	To proceed further, write explicitly $F=\sum p_i \otimes w_i$ where $\{w_i\}$ is a basis of $W$ with $wt(w_i)\leq wt(w_{i+1})$, and $wt(p_i)+wt(w_i)\leq 0$.

	For part (iv), we can assume again (by translating on the source and target) that $F(0)=0$ and $S=V^{\leq-\lambda}$.
	But the last condition implies by \autoref{eq:vanishing_sres_polynomials} that for for a polynomial function $p$ (without constant term) we have $p(S)=0$ if $wt(p)< \lambda$.
	So the terms $p_i\otimes w_i$ contribute only if $wt(p_i)\geq \lambda$, and hence $wt(w_i)\leq -\lambda$.
	This implies that $F(S)\subset \cW^{\leq -\lambda}$, as required.

	For part (v), the map $[F]$ is set-theoretically well-defined by part (iv).
	It remains to check it is subresonant.
	Denote by $[w_i]$ the image of $w_i$ in $W/\cW^{\leq -\lambda}$.
	Then the terms $p_i\otimes [w_i]$ are nonzero only if $wt(w_i)>-\lambda$, and hence $wt(p_i)<\lambda$.
	Again by \autoref{eq:vanishing_sres_polynomials} we know that $p(v+v')=p(v)$ when $wt(p)< \lambda$ and $v'\in V^{\leq -\lambda}$.
	It follows that each $p_i$ descends to a polynomial map $[p_i]$ on $V/V^{\leq -\lambda}$.
	Because taking quotients does not affect the weights, it follows that the descended map $[F]=\sum [p_i]\otimes [w_i]$ is also subresonant.
\end{proof}

\subsubsection{Veronese embedding for subresonant spaces}
	\label{sssec:veronese_embedding_for_subresonant_spaces}
Let \index{$P_V$}$P_V:=\Sym^{\bullet}\left(V^\dual\right)^{\leq \lambda_1}$ be the space of polynomial functions that have weight at most $\lambda_1$.
Then entries in any subresonant polynomial map can be only from $P_{V}$ (and the span of all such entries equals $P_{V}$).

The action of subresonant maps by pullback preserves $P_{V}$ (by \autoref{prop:characterization_of_filtration_on_polynomial_maps}) and acts on the right, i.e. $(g_1g_2)^*p=g_2^* (g_1^*p)$ (since $(g_1g_2)^* p(x)=p(g_1g_2x)$). We take the action on its dual $P_{V}^\dual$ to then be a left action, as is natural.

We have the tautological embedding given by evaluating at a point.\index{$ev$}
\begin{align*}
	V & \xrightarrow{\quad ev\quad} P_{V}^{\dual}=:L_V\\
	v & \quad \mapsto \quad ev(v)(p)=p(v)
\end{align*}
If we would have taken all polynomial functions up to some degree, this map is also known as the Veronese embedding.
Note also that we included the constant functions, therefore the image of $V$ is contained in a codimension $1$ affine subspace.

To ease notation we will use the shorthand $L_V:=P_V^{\dual}$.
As we will see below in \autoref{prop:linearization_of_subresonant_maps}, the purpose of the evaluation map $V\xrightarrow{ev}L_V$ is to linearize the action of the subresonant maps.

\subsubsection{Subresonant structures on affine spaces}
	\label{sssec:subresonant_structures_on_affine_spaces}
Suppose now that $E$ is an affine space, with associated vector space of translations $V$ (we will say that $E$ is affine over $V$).
Suppose that $V$ has a filtration $V^{\leq \bullet}$, and analogously for its dual $V^\dual$.
Then the space of polynomial functions $\Poly(E)$ on $E$ is itself filtered, with constants in weight $0$.
Note that for an affine space, we now have a \emph{filtration} by degree and not a direct sum decomposition.
The discussion so far generalizes naturally to the case of an affine space $E$ and we will denote by $P_E:=\Poly(E)^{\leq \lambda_1}$ and \index{$L_E$}$L_E:=P_{E}^{\dual}$, with evaluation embedding $E\xrightarrow{ev}L_{E}$ as before.

\begin{definition}[Subresonant structure on a manifold]
	\label{def:subresonant_structure_on_a_manifold}
	A \emph{subresonant structure} on an $n$-dimensional manifold $M$, with weights $-\lambda_1<\cdots < -\lambda_l$ is a collection of charts
	\[
		\phi_{i}\colon V_i \to U_i
	\]
	where $\bigcup U_i = M$ and $V_i$ are filtered vector spaces with the same weights, such that the transition maps on overlaps $\phi_j^{-1}\circ\phi_i$ agree with restrictions of subresonant polynomial maps in $\Poly^{sr}(V_i,V_j)$.
	A \emph{subresonant map} map between manifolds with subresonant structures is a map that is subresonant in respective coordinate charts.
\end{definition}

\subsubsection{Linearizing local systems}
	\label{sssec:linearizing_local_systems}
	Given a manifold $M$ with a subresonant structure, let \index{$P_M$}$P_M$ denote the functions on $M$ which are pullbacks of the corresponding functions in $P_{V_i}$ from various charts.
	Because subresonant maps respect the weights, the elements of $P_M$ give well-defined functions on $M$.
	On any open set in $M$ these form a finite-dimensional vector space, and if the open set is sufficiently small (simply-connected suffices) the dimension of the vector space is one and the same constant.
	In general, the elements of $P_M$ give a locally constant family of vector spaces (see \autoref{sssec:local_systems_and_holonomies} and \autoref{sssec:example_of_affine_structure}).
	Because in all our applications $M$ will be simply-connected, in those cases we will speak of a global vector space $P_M$ of functions on $M$.
	Analogously to the earlier discussion, let \index{$L_M$}$L_M:=P_{M}^{\dual}$ denote its dual.

	To distinguish a vector space of functions such as $P_U$ for $U\subset M$ a simply-connected open chart in $M$, and its dual $L_U$, we will denote the corresponding local system by \index{$L$@$\bL_M$}$\bL_M\to M$, which in local trivializations is identified as $\bL_U\isom L_U\times U$.
	Let now $M\xrightarrow{ev_{M}}\bL_M$ denote the tautological evaluation map.
	In local charts, or when $M$ is simply connected, we identify $\bL_M$ with $L_M\times M$, and conflate the evaluation map \index{$ev_M$}$ev_M$ with the induced map $M\to L_M$.

\begin{proposition}[Linearization of subresonant maps]
	\label{prop:linearization_of_subresonant_maps}
	Let $V$ be a filtered vector space with subresonant structure and \index{$G$@$\bbG^{sr}(V)$}$\bbG^{sr}(V)$ be the group of invertible subresonant maps.
	\begin{enumerate}
		\item 
		The evaluation map \index{$ev$}$ev\colon V\to L_V$ is injective and equivariant for the (left) action of subresonant maps.
		It gives a faithful group representation
		\[
			\index{$\rho_L$}\rho_L\colon \bbG^{sr}(V)\to \GL\left(L_V\right)
		\]
		The linear span of $ev(V)$ inside $L_V$ equals $L_V$.
		\item Let $M_1,M_2$ be two manifolds with subresonant structures and $ M_1\xrightarrow{f} M_2$ a subresonant map.
		Then there is an associated map of local systems $\bL_{M_1}\xrightarrow{L_f}\bL_{M_2}$, covering $f$ and linear on the fibers, such that the diagram commutes
		\begin{equation}
			\label{eqn_cd:linearization_commutes}
		\begin{tikzcd}
			M_1
			\arrow[r, "f"]
			\arrow[d, "ev_{M_1}"]
			& 
			M_2
			\arrow[d, "ev_{M_2}"]
			\\
			\bL_{M_1}
			\arrow[r, "L_f"]
			&
			\bL_{M_2}
		\end{tikzcd}
		\end{equation}
	\end{enumerate}	
\end{proposition}
\begin{proof}
	For part (i), observe that the affine functions are included among all the polynomials in $P_V$.
	This gives the injectivity of $ev_{V}$, as well as the faithfulness of the representation (since if a polynomial map fixes all affine functions, it must be the identity).
	By the definition of subresonant maps, they preserve $P_V$ and clearly act by linear transformations on it, so the map is a linear representation indeed.

	For part (ii), observe that pullback of functions induces a map in local charts $U_i$ on $M_i$ as $P_{U_2}\xrightarrow{f^*}P_{U_1}$, so the induced transpose map on dual vector spaces $L_{U_1}\xrightarrow{L_f}L_{U_2}$ is well-defined.
	Compatibility with the evaluation map is then clear.
\end{proof}

\subsubsection{Subresonant vector bundles}
	\label{sssec:subresonant_vector_bundles}
Let $V$ be a vector space with a subresonant structure and weights $-\lambda_1<\cdots<-\lambda_l<0$, and let $E_0$ be another vector space equipped with a filtration
\[
	E_0=E_0^{\leq \mu_1}\supsetneq E_0^{\leq \mu_2}\cdots E_0^{\leq \mu_{k}}\supsetneq 0
	\quad\text{and arbitrary weights }
	\mu_1>\cdots>\mu_k.
\]
We will take $E:=E_0\times V$ as a standard model of vector bundle with subresonant structure (cf. \autoref{def:subresonant_maps_structure_on_bundle} below).

\begin{definition}[Subresonant maps and structures on bundle]
	\label{def:subresonant_maps_structure_on_bundle}
	Let $E=E_0\times V\to V$ and $E'=E_0'\times V'\to V'$ be filtered vector bundles as in \autoref{sssec:subresonant_vector_bundles}.
	A \emph{subresonant map} between them is a polynomial map $F\colon E\to E'$, covering a subresonant map $f\colon V\to V'$, such that $F$ is linear on the fibers, and of total weight required to be nonpositive (see \autoref{rmk:on_weight_condition_for_subresonant_maps_of_vector_bundles}).

	Suppose $M$ is a manifold with subresonant structure and $E\to M$ is a filtered vector bundle.
	A subresonant structure on $E$ is a collection of trivializations
	\[
		c_i\colon E_{u_i}\times V_i \toisom  E\vert_{U_i}
	\]
	where $\cup U_i = M$ and the maps $V_i\toisom U_i$ give the subresonant structure on $M$, such that the gluing maps $c_{j}^{-1}\circ c_i$ are subresonant maps of vector bundles.

	More generally, given bundles $E_i\to M_i, i=1,2$ with subresonant structure, a subresonant map between them is one which is a subresonant map of bundles in each trivializing chart.

\end{definition}

\begin{remark}[On weight condition for subresonant maps of vector bundles]
	\label{rmk:on_weight_condition_for_subresonant_maps_of_vector_bundles}
	To clarify the weight condition in the above \autoref{def:subresonant_maps_structure_on_bundle}, observe that a polynomial map $F$ between trivialized vector bundles splits, using the direct product structure of the bundles, as $F=(F_E,f)$ where $f\in \Poly(V,V')^{\leq 0}$ and $F_E\in \Big[\Poly(V)\otimes \Hom(E_{0},E_{0}')\Big]^{\leq 0}$.
	See also \autoref{sssec:example_of_subresonant_bundle_map} for an example.
\end{remark}

\begin{proposition}
	\label{prop:properties_of_bundles_with_subresonant_structure}
	Let $E_V=E_0\times V\to V$ and $E_{V'}=E_{0}'\times V'\to V'$ be trivialized bundles with subresonant structure and weights $\mu_1>\cdots >\mu_k$ on the bundle and $-\lambda_1<\cdots<-\lambda_l$ on the base.
	Let $F\colon E_V\to E_{V'}$ be a subresonant bundle map.
	Set $D:=\mu_1-\mu_k$.
	\begin{enumerate}
		\item  Consider the space $P_{E_{V}}:=\Poly(V)^{\leq D}\otimes (E')^{\dual}$ of functions on $E_{V}$, which are linear on the fibers and depend polynomially on the base (with $P_{E_{V'}}$ defined analogously).
		Then $F^* P_{E_{V'}}\subseteq P_{E_{V}}$, i.e. pullback of functions preserves these spaces.

		More generally, $F^*$ preserves the weight filtration on functions on the total spaces of the bundles.
		\item The map $F$ respects the vector bundle filtrations on $E_V$ and $E_{V'}$ and induces a constant map on the associated graded bundles.
	\end{enumerate}
\end{proposition}

\begin{proof}
	Let us write the map $F$ as $F=(F_E,f)$ using the trivializations, with $f\in \Poly(V,V')^{\leq 0}$ and $F_E\in \Poly(V)\otimes \Hom(E_{0},E_{0}')^{\leq 0}$.
	Choose now a basis $\xi_i$ for $E_0^{\dual}$, $e_j'$ for $E_{0}'$, and $\xi_j'$ a dual basis for $(E_0')^{\dual}$.
	We choose the bases such that $\xi_l'(e_j')\neq 0$ implies $wt(\xi_l')+wt
	(e_j')\geq 0$ (and is equivalent to it).
	Let us write explicitly
	\[
		F_E = \sum p_{i,j} \otimes \xi_i \otimes e_{j}'
		\text{ with }
		wt(p_{i,j}) + wt(\xi_i) + wt(e_j') \leq 0
	\]
	For a linear functional $\xi_l'$ on $E_0'$, we can compute
	\begin{align}
		\label{eq:pullback_vector_bundle_linear_function}
		F^*\xi_{l}' 
		=
		\sum p_{i,j}\otimes \xi_i \cdot \xi_l'(e_j')
	\end{align}
	and we see that only terms with $wt(e_j')\geq -wt(\xi_l')$ contribute, i.e. we must have $wt(p_{i,j}) + wt(\xi_i)\leq wt(\xi_l')$.
	Therefore $wt(F_E^*\xi_l')\leq wt(\xi_l')$, and note also that $wt(p_{i,j})\leq wt(\xi_l')-wt(\xi_i)\leq D$.
	So at least for linear functions, part (i) is verified.
	If we have an element of $P_{E'}$ of the form $p' \cdot \xi'$ with $p'$ a polynomial function on $V'$ and $\xi'$ a linear functional on $E_0'$, then $F^*(p'\cdot \xi')=f^*(p')\cdot F_E^*(\xi')$ and it is clear again that the weight is respected (since $f^*$ respects it).

	For (ii), observe that the expression for $F^*\xi_{l}'$ from \autoref{eq:pullback_vector_bundle_linear_function} implies that if a term $p_{i,j}\otimes \xi_i$ appears, then $wt(\xi_i)\leq wt(\xi_l')$ since $wt(p_{i,j})\geq 0$.
	Furthermore these last two inequalities are strict unless $p_{i,j}$ is a constant.
	Therefore the map $F$ respects the grading of the bundles pointwise, and induces constant maps on the gradings.
\end{proof}

\subsubsection{Linearization of subresonant bundles}
	\label{sssec:linearization_of_subresonant_bundles}
We now explain how to produce canonical vector spaces and linearizations associated to vector bundles with subresonant structure, analogously to \autoref{prop:linearization_of_subresonant_maps}.
Let $M$ be a manifold with subresonant structure and $E\to M$ a bundle with subresonant structure.
Let \index{$P$@$\bP_E$}$\bP_E$ be the local system of functions on $E$ that in local trivializing charts are in $\Poly^{\leq D}\otimes E_0^{\dual}$.
By \autoref{prop:properties_of_bundles_with_subresonant_structure}, the space $\bP_E$ is independent of the choice of trivializing charts and in sufficiently small (simply-connected suffices) open sets on $M$, it is a vector space of fixed dimension.
Let then \index{$L$@$\bL_E$}$\bL_E:=\bP_E^\dual$ be its dual.

We have as before a tautological evaluation map
\begin{align}
	\label{eq:evaluation_embedding_sresn_vector_bundle}
	\begin{split}
	E & \xrightarrow{ev} \bL_E\\
	e & \mapsto ev(e)(p):=p(e) \quad \text{for }p\in \bP_E.
	\end{split}
\end{align}
Again for every subresonant map of vector bundles $F\colon E\to E'$ covering a map $f\colon M\to M'$ there exists a map of local systems $\rho_L(F)\colon \bL_E\to \bL_{E'}$ such that
\begin{equation}
	\label{eqn_cd:subresonant_bundle_linearization}
\begin{tikzcd}
	E
	\arrow[r, "F"]
	\arrow[d, "ev_E ", swap]
	& 
	E'
	\arrow[d, "ev_{E'}"]
	\\
	\bL_E
	\arrow[r, "\rho_L(F)"]
	&
	\bL_{E'}
\end{tikzcd}
\end{equation}
commutes.

\subsubsection{Local systems and ``holonomies''}
	\label{sssec:local_systems_and_holonomies}
The linearization constructions associated to a subresonant manifold $M$ from \autoref{prop:linearization_of_subresonant_maps}, and a subresonant vector bundle $E\to M$ from \autoref{sssec:linearization_of_subresonant_bundles} can be phrased as follows.
The local systems $\bL_M$ and $\bL_E$ on sufficiently small open sets (simply-connected suffices) provide a vector space, and on overlaps the transition maps between the vector spaces are constant.
In most situations of dynamical interest, these considerations are applied to small open sets which are simply connected, so we can speak of a single vector space.

We can alternatively say $\bL_M$ and $\bL_E$ give vector bundles over $M$ equipped with a flat ``Gauss--Manin'' connection.
The parallel transport maps from $x$ to $y$ along a path in $M$ will be denoted \index{$P$@$\cP_{GM}(x,y)$}$\cP_{GM}(x,y)$, with subscripts $M$ or $E$ when we'll need to distinguish between $L_M$ and $L_E$.
The maps have the property
\[
	\cP_{GM}(x,z)=\cP_{GM}(y,z)\circ\cP_{GM}(x,y)
\]
and when such maps are obtained from dynamical considerations, they are also called ``holonomies'' in the dynamical literature.

\subsubsection{Example of subresonant maps}
	\label{sssec:example_of_subresonant_maps}
Assume that we have $3$ Lyapunov exponents $\lambda_1>\lambda_2>\lambda_3>0$, with the filtration of $V$ given as
\[
	0\subsetneq V^{\leq-\lambda_1}\subsetneq V^{\leq-\lambda_2} \subsetneq V^{\leq -\lambda_3}=V
\]
and the filtration of the dual:
\[
	0\subsetneq \left(V^\dual\right)^{\leq\lambda_3}\subsetneq 
	\left(V^{\dual}\right)^{\leq\lambda_2} \subsetneq \left(V^{\dual}\right)^{\leq \lambda_1}=V^{\dual}.
\]
Let us simplify further and assume that $\dim V = 3$ and each subquotient has dimension $1$.

Choose a basis $e_x,e_y,e_z$ of $V$ and dual basis $x,y,z\in V^\dual$ of coordinate functions such that filtration on $V$ is $span(e_x)\subseteq span(e_x,e_y)\subseteq span(e_x,e_y,e_z)$ and the dual filtration on $V^\dual$ is $span(z)\subseteq span(z,y)\subseteq span(z,y,x)$.
Then the general subresonant map from $V$ to itself is of the form
\begin{align}
	\label{eq:subresonant_map_example}
	\begin{bmatrix}
		x\\
		y\\
		z
	\end{bmatrix}
	\xrightarrow{F}
	\begin{bmatrix}
		p\cdot x & + \sum\limits_{i\cdot \lambda_2 + j\cdot \lambda_3 \leq \lambda_1} a_{ij} y^i z^j\\
		q\cdot y & + \sum\limits_{k\cdot \lambda_3\leq \lambda_2}b_k z^k\\
		r\cdot z &
	\end{bmatrix}
\end{align}
Let us make explicit the linearization from \autoref{sssec:veronese_embedding_for_subresonant_spaces}.
To simplify notation, we omit the constant functions; including them would add an extra coordinate at the bottom.
\begin{align}
	\label{eq:subresonant_map_linearization_example}
	\rho(F)=
	\begin{bmatrix}
		p & \cdots  a_{i,j}\cdots & a_{1,0} & \cdots & a_{0,k} & \cdots &a_{0,1} \\
		 & \vdots &  &\\
		0 & q^ir^j & \cdots & c_{i,j,l,k} & \cdots \\
		0 & 0 & \vdots\\
		0 & 0 & q & \cdots & b_{k} & \cdots & b_1\\
		0 & 0 & 0 & \ddots & 0 & 0 & 0\\
		0 & 0 & 0 & 0 & r^k & 0 & 0 \\
		0 & 0 & 0 & 0 & 0 & \ddots & 0\\
		0 & 0 & 0 & 0 & 0 & 0 & r
	\end{bmatrix}
	\begin{bmatrix}
		x\\
		y\\
		z
	\end{bmatrix}
	\to
	\begin{bmatrix}
		x\\
		\vdots\\
		y^iz^j\\
		\vdots\\
		y\\
		\vdots\\
		z^k\\
		\vdots\\
		z
	\end{bmatrix}
\end{align}
where the coefficients $c_{i,j,l}$ are computed from $y^iz^j\mapsto (y+\sum b_k z^k)^iz^j = \sum c_{i,j,l,k}y^lz^{j+k(i-l)}$


\subsubsection{Example of subresonant bundle map}
	\label{sssec:example_of_subresonant_bundle_map}
For simplicity, we will take the base subresonant manifold to be $2$-dimensional with coordinate functions $(x,y)$ and respective weights $-\lambda_1<-\lambda_2<0$.
Let $E$ be a $3$-dimensional vector space with basis $e_1,e_2,e_3$ and weights $\mu_1>\mu_2>\mu_3$ and dual basis $\xi_1,\xi_2,\xi_3$.
A general bundle map (with no translation and identity derivative at the origin on the base) is of the form
\begin{align*}
	\label{eq:subresonant_bundle_map}
	\everymath{\displaystyle}
	\begin{bmatrix}
		\xi_1
		\vphantom
		{\left(\sum_{i\lambda_1+j\lambda_2\leq \mu_1-\mu_2}b_{i,j}^{(1,2)}x^iy^j\right)}
		\\
		\xi_2
		\vphantom
		{\left(\sum_{i\lambda_1+j\lambda_2\leq \mu_1-\mu_2}b_{i,j}^{(1,2)}x^iy^j\right)}
		\\
		\xi_3\\
		x
		\vphantom{\sum_{i\lambda_2\leq \lambda_1}a_iy^i}
		\\
		y
	\end{bmatrix}
	\mapsto
	\begin{bmatrix}
		\xi_1
		 & +
		\xi_2\cdot 
		\left(\sum_{i\lambda_1+j\lambda_2\leq \mu_1-\mu_2}b_{i,j}^{(1,2)}x^iy^j\right)
		& +
		\xi_3\cdot
		\left(\sum_{i\lambda_1+j\lambda_2\leq \mu_1-\mu_3}b_{i,j}^{(1,3)}x^iy^j\right)\\
		\xi_2
		& +
		\xi_3\cdot
		\left(\sum_{i\lambda_1+j\lambda_2\leq \mu_2-\mu_3}b_{i,j}^{(2,3)}x^iy^j\right)\\
		\xi_3\\
		x & +\sum_{i\lambda_2\leq \lambda_1}a_iy^i\\
		y &
	\end{bmatrix}
\end{align*}
The embedding corresponding to \autoref{eqn_cd:subresonant_bundle_linearization} that linearizes the vector bundle maps takes the point with coordinates $[\xi_1\,\, \xi_2\,\, \xi_3\,\,  x\,\, y]^t$ to
\[
	[\xi_1\,\,  \cdots 
	(x^iy^j \xi_2) \cdots \xi_2 
	\cdots
	(x^k y^l\xi_3)
	\cdots
	\xi_3
	\,\,
	x\,\, y
	]^t
\]
and makes the corresponding bundle maps upper-triangular and constant.
If desired, we could also linearize the dynamics on the base.

\begin{example}[Affine structure as subresonant structure]
	\label{eg:affine_structure_as_subresonant_structure}
	\label{sssec:example_of_affine_structure}
	We end with the simplest example of a subresonant structure, which should also clarify the ``local system'' aspect of the discussion in \autoref{sssec:local_systems_and_holonomies}.
	Let $A$ be an affine space over a vector space $V$, i.e. we have a free and transitive action of $V$ on $A$.
	First, observe that $A$ can be canonically realized as an (affine) hyperplane, such that affine maps of $A$ act canonically as linear maps.
	Indeed, let $P_A$ denote the vector space of affine functions on $A$, i.e. $\xi\in P_A$ if for any $a\in A$ the function $v\mapsto \xi(a+v)-\xi(a)$ (for $v\in V$) is a linear function on $V$.
	Then $P_A$ has dimension $\dim V + 1$ (since it includes the constants) and if we denote its dual by $L_A$, we have a canonical embedding $ev\colon A\into L_A$ by evaluation of functions.
	The image of $A$ is the affine hyperplane of functionals that take the constant function $1\in P_A$ to $1$.

	Finally, let $\Lambda\subset V$ be a lattice (i.e. discrete cocompact group).
	Then we can form the quotient $T:=A/\Lambda$.
	Locally on $T$ we have the notion of affine function, but globally there are no affine functions.
	Therefore $\bP_T$ can be viewed as a vector bundle with flat connection, or a locally constant family of vector spaces.
\end{example}

The next result will turn out useful in constructing interpolation maps:
\begin{lemma}[Polynomial maps from linear maps]
	\label{lem:polynomial_maps_from_linear_maps}
	Suppose $V$ is a vector space, equipped with a filtration with exponents, and associated subresonant structure.
	Let $ev\colon V\to L_V$	be the linearization map from \autoref{sssec:veronese_embedding_for_subresonant_spaces}.
	Suppose that $g\in \GL(L_V)$ is a linear map, preserving $ev(V)$ set-theoretically, and respecting the weight filtration on $L_V$.

	Then there exists a subresonant polynomial map $\cG\colon V\to V$ such that it's linearization agrees with $g$.
\end{lemma}

\begin{corollary}[Polynomial maps between different subresonant manifolds]
	\label{cor:polynomial_maps_between_different_subresonant_manifolds}
	Suppose that $V_1,V_2$ are vector spaces with filtrations, equipped with subresonant structures, and with linearizations $ev_{i}\colon V_i\to L_{V_i}$.
	Suppose also that there exists some subresonant isomorphism $\phi\colon V_1\to V_2$.

	Suppose now that $g\colon L_{V_1}\to L_{V_2}$ is a linear map preserving the weight filtrations and taking $ev\left(V_1\right)$ to $ev\left(v_2\right)$.
	Then there exists a subresonant map $\cG\colon V_1\to V_2$ whose linearization is $g$.
\end{corollary}
\begin{proof}
	Let $L\phi\colon L_{V_1}\to L_{V_2}$ be the linearization of the isomorphism $\phi$.
	Apply \autoref{lem:polynomial_maps_from_linear_maps} to $h:=L\phi^{-1}\circ g\colon L_{V_1}\to L_{V_1}$ to obtain a subresonant map $\cH\colon V_1\to V_1$.
	Then set $\cG:=\phi\circ \cH\colon V_1\to V_2$.
\end{proof}

\begin{proof}[Proof of \autoref{lem:polynomial_maps_from_linear_maps}]
	Let $g^t\in \GL(L_{V}^{\dual})$ be the transpose of the given linear map.
	Recall from \autoref{sssec:veronese_embedding_for_subresonant_spaces} that $L_{V}^\dual$ is a space of polynomial functions on $V$.
	Fix a set of (say, linear) coordinate functions $x_1,\ldots,x_n$ on $V$, which are in particular elements of $L_V^{\dual}$.
	Then $g^t(x_1),\ldots,g^t(x_n)$ can be taken as the coordinate functions of a polynomial map $\cG\colon V\to V$.

	Because $g$, and hence $g^t$, respects the weight filtration, it follows that the weights of $g^t(x_i)$ are at most those of $x_i$, hence that $\cG$ is a subresonant polynomial map, with linearization $L\cG\colon L_V\to L_V$.
	We now consider the linear map $g^{-1}\circ L\cG$ of $L_V$, which still preserves $ev(V)$.

	Let $x_{I}^{\alpha_{I}}$ be the list of monomials that form a basis of $L_{V}^{\dual}$, with $I,\alpha_I$ appropriate multi-indexes.
	Note that for any subresonant map $\cG$, the action of its linearization on $L_V$ in this basis is given by the pullback $\cG^{*}x_I^{\alpha_I}$ for the various monomials.
	In particular, the action is compatible whether we view $x_I^{\alpha_I}$ as a basis element of $L_V^{\dual}$ or as a function on $V$, and the compatibility is obtained by $ev\colon V\into L_V$.

	It follows that $g^{-1}\circ L\cG$ acts on $x_1,\ldots,x_n$ as the identity map.
	It then follows that it acts as the identity on all monomials $x_I^{\alpha_I}$ and hence that it is the identity map.
\end{proof}



\subsection{Resonant maps and groups of (sub)resonant transformations}
	\label{ssec:resonant_maps_and_}

\subsubsection{Setup}
	\label{sssec:setup_resonant_maps_and_subresonant_transformations}
We continue our study of subresonant maps and introduce the notion of resonant maps.
We then study the associated groups of resonant, subresonant, and strictly subresonant maps.

\subsubsection{Splittings and resonant maps}
	\label{sssec:splittings_and_resonant_maps}
Given a filtered vector space $V$ with weights $\{-\lambda_i\}$, let $-\Lambda\colon V\to V$ be a linear map that preserves the filtration and acts by $-\lambda_i$ on $V^{\leq -\lambda_{i}}/V^{\leq -\lambda_{i-1}}$.
Then this operator induces a \emph{splitting of the filtration}, i.e. a decomposition $V=\oplus V^{-\lambda_i}$ given by the eigenspaces of $-\Lambda$.
Conversely, any such splitting determines a unique operator $-\Lambda$.

Note that an initial splitting of the filtered vector space $V$ determines a splitting of all corresponding linear algebra constructions $V^\dual$, $\Sym^d V$ etc., and corresponding operators $\Lambda^t$ on $V^\dual$, $\Sym^d(-\Lambda)$, etc.
We maintain our standing assumption that on vector spaces $V,W$, etc. the weights are strictly negative.

\begin{definition}[Resonant maps]
	\label{def:resonant_maps}
	Let $V,W$ be filtered vector spaces, with splittings determined by $-\Lambda_V,-\Lambda_W$.
	A polynomial map $F\in \Poly(V,W)$ is called \emph{resonant} (relative to the fixed splittings) if it is in weight $0$, i.e. if \index{$P$@$\Poly^0(V,W)$}$F\in \Poly^0(V,W)$.
\end{definition}

Let us observe first that nontrivial translations \emph{cannot} be among the resonant maps, as they must contain terms with strictly negative weights.
The properties of subresonant maps from \autoref{prop:properties_of_subresonant_polynomial_maps} translate immediately to analogous ones for resonant maps.

\begin{proposition}[Resonant and subresonant map characterizations]
	\label{prop:resonant_and_subresonant_map_characterizations}
	Let $V,W$ be filtered vector spaces with fixed splittings given by $-\Lambda_V,-\Lambda_W$.
	\begin{enumerate}
		\item A polynomial map $F\in \Poly(V,W)$ is resonant if and only if it commutes with the (exponentiated) splitting operators, i.e.
		\[
			F\circ e^{\Lambda_V} = e^{\Lambda_W} \circ F \text{ and similarly for }-\Lambda_V,-\Lambda_W
		\]
		\item A polynomial map $F\in \Poly(V,W)$ is subresonant if and only if the following limit exists and:
		\[
			\lim_{n\to +\infty} e^{-n\Lambda_W}\circ F \circ e^{n\Lambda_V}= F^{(0)}
		\]
		where $F^{(0)}$ is the resonant map which is the projection to the $0$th weight of $F$ in $\Poly(V,W)$.
	\end{enumerate}
\end{proposition}
The characterizations provided by this proposition show that our notions of resonant and subresonant maps agree with those used in \cite[Def. 2.2]{KalininSadovskaya2017_Normal-forms-for-non-uniform-contractions}, see also \cite[\S3]{Feres_A-differential-geometric-view-of-normal-forms-of-contractions}.

\begin{proof}
	Using the splittings, we can write any polynomial map $F$ as $\sum_{i,j} p_{i,j}^{(\alpha_{i,j})}\otimes w_i$ where $p_{i,j}^{(\alpha_{i,j})}$ are polynomial functions on $V$ of weight $\alpha_{i,j}\geq 0$ and $w_i$ are a basis of $W$ with weights ${\phi(i)}$ where $\phi\colon \{1,\ldots,\dim W\}\to \{-\lambda_1,\ldots,-\lambda_l\}$ is an order-preserving map.
	Then we have
	\begin{align*}
		F\circ e^{\Lambda_V} & = \sum_{i,j}e^{\alpha_{i,j}}p_{i,j}^{(\alpha_{i,j})}\otimes w_i\\
		e^{\Lambda_W} \circ F & = \sum_{i,j}e^{-\phi(i)}p_{i,j}^{(\alpha_{i,j})}\otimes w_i
	\end{align*}
	so part (i) follows, i.e. we must have $\alpha_{i,j}+\phi(i)=0,\forall i,j$.

	For part (ii) an analogous calculation gives
	\[
		e^{-n\Lambda_W}F\circ e^{n\Lambda_V} = 
		\sum_{i,j}
		e^{n(\alpha_{i,j} + \phi(i))}
		p_{i,j}^{(\alpha_{i,j})}\otimes w_i
	\]
	so in order to have a limit at all, we must have $\alpha_{i,j}+\phi(i)\leq 0,\forall i,j$, i.e. the map is subresonant, and in the limit as $n\to +\infty$ only the terms of weight $0$ survive.	
\end{proof}

\subsubsection{Further decompositions of polynomial maps}
	\label{sssec:further_decompositions_of_polynomial_maps}
The splitting $V = \oplus V^{-\lambda_i}$ induces corresponding splittings on polynomial maps, and the ones relevant for our discussion are those of nonpositive weight:
\[
	\Poly(V,V)^{\leq 0} = \bigoplus_{-\lambda_i+\sum a_j\lambda_j\leq 0}\Poly(V,V)^{\sum a_j\lambda_j,-\lambda_i}
\]
We also have the decomposition into weight $0$ and strictly negative weight, and a further refinement:
\begin{align*}
	\Poly(V,V)^{\leq 0} & = \Poly(V,V)^{0}\oplus \Poly(V,V)^{<0} \\
	\Poly(V,V)^{0} & = \left(\bigoplus_{\lambda_i} \Poly(V,V)^{\lambda_i,-\lambda_i}\right)
	\oplus
	\underbrace{\left(\bigoplus_{\lambda_i =\sum_{j>i} a_j \lambda_j}\Poly(V,V)^{\sum a_j \lambda_j,-\lambda_i}\right)}_{\Poly(V,V)^{0}_{nil}}
\end{align*}
where the summand $\Poly(V,V)^{\lambda_i,-\lambda_i}$ can also be identified with the linear maps $\End\left(V^{-\lambda_i}\right)$.
The maps denoted by \index{$P$@$\Poly(V,V)^{0}_{nil}$}$\Poly(V,V)^{0}_{nil}$ form a nilpotent Lie algebra under composition (and the summation is over tuples of $a_j$ such that $\lambda_i = \sum a_j \lambda_j$ and $j>i$).

Recall that the linearization map $V\xrightarrow{ev}L_V$ induces an injection compatible with composition of maps:
\[
	\Poly(V,V)^{\leq 0}\xrightarrow{L_\rho}\End(L_V)^{\leq 0}
\]
where $\End(L_V)^{\leq 0}$ are the linear maps preserving the weight filtration.

\begin{definition}[Groups of resonant and strictly subresonant maps]
	\label{def:groups_of_resonant_and_strictly_subresonant_maps}
	Let \index{$G$@$\bbG^{sr}(V)$}$\bbG^{sr}(V)$ denote the set of subresonant maps which have invertible derivative at the origin.

	Let also \index{$G$@$\bbG^{ssr}(V)$}$\bbG^{ssr}(V)$ denote the subset of $\bbG^{sr}(V)$ consisting of elements that, under the map to $\End(L_V)^{\leq 0}$, act as the identity on the associated graded $\gr_\bullet L_V$.
	These will be called \emph{strictly subresonant maps}.

	Define also the set of maps \index{$G$@$\bbG^{r}$}$\bbG^{r}:=\bbG^{sr}\cap \Poly(V,V)^0$ consisting of elements in $\bbG^{sr}$ that are also resonant maps.
\end{definition}
Note that the sets $\bbG^{sr}(V)$ and $\bbG^{ssr}(V)$ only depend on the weight filtration, and not the decomposition of $V$.
In particular, from \autoref{prop:group_structure_on_subresonant_maps} below it will follow that such groups can be associated to manifolds with subresonant structure, i.e. these objects are independent of subresonant charts that we choose.


\begin{proposition}[Group structure on (sub)resonant maps]
	\label{prop:group_structure_on_subresonant_maps}
	\leavevmode
	\begin{enumerate}
		\item 
		\label{item:ssresn_are_a_group}
		The strictly subresonant maps can be expressed as
		\[
			\bbG^{ssr}(V) = \id_V + \Poly(V,V)^{<0}
		\]
		and in particular form a unipotent algebraic group.
		\item
		\label{item:resn_gp_semidir}
		The resonant maps can be expressed as a semidirect product
		\[
			\bbG^r(V) = \left(\prod \GL\left(V^{-\lambda_i}\right)\right)\ltimes\left(\id_V + \Poly(V,V)^{0}_{nil}\right)
		\]
		and in particular form an algebraic group.
		\item
		\label{item:sresn_gp_semidir}
		 All the subresonant maps can be expressed as a semidirect product
		\[
			\bbG^{sr}(V) = \bbG^r(V)\ltimes \bbG^{ssr}(V)
		\]
		and in particular form an algebraic group.
		The group $\bbG^{ssr}(V)$ is a normal subgroup of $\bbG^{sr}(V)$.
	\end{enumerate}
\end{proposition}
\begin{proof}
	For part \autoref{item:ssresn_are_a_group}, by the definition of $\bbG^{ssr}(V)$, these are subresonant maps acting as the identity on the associated graded of $L_V$.
	So they cannot have a nontrivial component in $\Poly(V,V)^0$ besides the identity map $\id_V$, otherwise they would pull back some linear function to a nontrivial polynomial of the same weight.
	But since $\Poly(V,V)^{<0}$ form a nilpotent Lie algebra under composition, existence of an inverse and the group structure follow.

	For part \autoref{item:resn_gp_semidir}, it is clear that post-composing with an element of $\prod \GL\left(V^{-\lambda_i}\right)$ still gives resonant transformations, and in particular such elements are resonant transformations.
	Thus, given an arbitrary resonant map $F$, by post-composing with a linear map, we can assume that $F= \id_V + F_{nil}$ where $F_{nil}\in \Poly(V,V)^0_{nil}$.
	Again we have a nilpotent Lie algebra $\Poly(V,V)^0_{nil}$ and invertibility follows.
	The semidirect product structure also immediately follows by checking that $\prod \GL(V^{-\lambda_i})$ normalizes $\Poly(V,V)^0_{nil}$.

	For part \autoref{item:sresn_gp_semidir}, it follows from the definitions that $\bbG^{ssr}$ is the kernel of the map $\bbG^{sr}\to \prod_{\mu_j}\GL(\gr_{\mu_j}L_V)$.
	Given any element $F$ in $\bbG^{sr}$ (so it has invertible derivative at the origin), we can first compose it with an element of $\prod \GL(V^{-\lambda_i})$ to make the derivative at the origin strictly subresonant (i.e. unipotent for the filtration $V^{\leq -\lambda_i}$).
	Then by applying another resonant transformation of the form $\id_{V}+F_{nil}$ we can ensure that $F-\id_V\in \Poly(V,V)^{<0}$, i.e. $F$ is strictly resonant.
	It follows that $\bbG^r$ surjects onto and is isomorphic to the quotient $\bbG^{sr}/\bbG^{ssr}$ and the result follows.
\end{proof}

Because strictly subresonant maps form a normal subgroup of all subresonant maps, we obtain:
\begin{corollary}[Strictly subresonant maps on subresonant manifolds]
	\label{cor:strictly_subresonant_maps_on_subresonant_manifolds}
	Suppose that $M$ is a manifold equipped with a subresonant structure in the sense of \autoref{def:subresonant_structure_on_a_manifold}, which is furthermore isomorphic as a subresonant manifold to a vector space with a filtration.
	We will refer to the set of such isomorphisms as ``global charts'' on $M$.

	Then the set of maps of $M$ induced by strictly subresonant ones on global charts are independent of the choice of global charts, and will be called the strictly subresonant maps of $M$.
\end{corollary}
This will apply, in particular, to the stable and unstable manifolds with the subresonant structure (and global charts) from \autoref{thm:subresonant_normal_form_suspended_systems}.

\begin{remark}[On resonant maps between manifolds]
	\label{rmk:on_resonant_maps_between_manifolds}
	Given a manifold with subresonant structure, it does not make sense to speak of endowing it with a resonant structure, since the group of resonant transformations of a filtered vector space with splitting does not act transitively on the vector space.
	However, given two pointed manifolds $(M_i,p_i)$ with subresonant structure, as well as charts $U_i\to M_i$ with $0\mapsto p_i$ and $U_i$ open sets in a filtered vector space with a splitting, it makes sense to ask if a map $F\colon M_1\to M_2$ with $F(p_1)=p_2$ is resonant, by considering it in the corresponding charts.
	
	Note that one has to specify not just a splitting of the tangent spaces $T_{p_i}M_i$, but in addition charts and in particular splittings of the linearization spaces $L_{M_i}$.
	In the case of interest in dynamics, the resonant normal form coordinates arise by considering the splittings of the linearization spaces of stable/unstable manifolds $L\cW^{s/u}$ that come from the Oseledets decomposition, since the spaces $L\cW^{s/u}$ form a cocycle for the dynamics.
\end{remark}





\section{Normal forms on manifolds and cocycles}
	\label{appendix:normal_forms_on_manifolds_and_cocycles}

\paragraph{Outline of section}
We state the main results in the discrete time setting, the generalization to continuous time being immediate, since one can apply the discrete result to the time-$1$ map of a flow and make use of the uniqueness properties of the normal forms coordinates.

The normal forms coordinates on stable manifolds are discussed in \autoref{ssec:subresonant_normal_forms_on_stable_manifolds}, and for cocycles on stable manifolds in \autoref{ssec:cocycle_normal_forms}.

\subsubsection*{References}
	\label{sssec:references}
Most of the results that we need can be extracted from the literature, with the most suitable for our purposes being the work of Kalinin--Sadovskaya \cite{KalininSadovskaya2017_Normal-forms-for-non-uniform-contractions}, see also Melnick \cite{Melnick_Non-stationary-smooth-geometric-structures-for-contracting-measurable}, Feres \cite{Feres_A-differential-geometric-view-of-normal-forms-of-contractions}, Katok--Rodriguez Hertz\cite[App.~A]{KatokRodriguez-Hertz2016_Arithmeticity-and-topology-of-smooth-actions-of-higher-rank-abelian}, and Li--Lu \cite{LiLu2005_Sternberg-theorems-for-random-dynamical-systems}.


\subsection{Subresonant normal forms on stable manifolds}
	\label{ssec:subresonant_normal_forms_on_stable_manifolds}

\subsubsection*{Lyapunov exponents of a dual cocycle}
	\label{sssec:lyapunov_exponents_of_a_dual_cocycle}
Suppose $E\to Q$ is a cocycle over a dynamical system on $Q$, and the dynamical system is not necessarily invertible, but the cocycle matrices are.
Then the dual cocycle is defined by
\[
	\index{$E_{x}^{\dual}$}E_{x}^{\dual} \xrightarrow{(A_x^{-1})^t} E_{Tx}^\dual
\]
where \index{$\norm$@$(\bullet)^\dual$}$(\bullet)^\dual$ denotes the dual vector space, \index{$\norm$@$(\bullet)^t$}$(\bullet)^t$ denotes the transpose map, $T$ denotes the dynamics (discrete or continuous) and $A_x\colon E_x\to E_{Tx}$ the original cocycle.
Note that
\begin{align*}
	\text{if the Lyapunov exponents of }E & \text{ are }
	\lambda_1>\lambda_2>\cdots\\
	\text{then the Lyapunov exponents of }E^\dual & \text{ are }
	-\lambda_1< -\lambda_2 < \cdots
\end{align*}

\subsubsection{Assumptions for normal forms}
	\label{sssec:assumptions_for_normal_forms_single_diffeo}
Suppose that $Q$ is a manifold, \index{$f$}$f\colon Q\to Q$ is a diffeomorphism preserving an ergodic measure $\nu$, such that furthermore the tangent cocycle $Df$ satisfies
\[
	\log \max \left(\norm{Df},\norm{Df^{-1}}\right)\in L^1(\nu).
\]
Suppose further that for any $k\geq 1$, there exists a measurable tempered function \index{$r_k$}$r_k\colon Q\to \bR_{>0}$, such that the $k$-th derivative cocycles $D^{(k)}f$ satisfy
\[
	q\mapsto \norm{f}_{C^k(q,r_k(q))} \text{ is tempered.}
\]
Here \index{$C^k(q,r_k(q))$}$C^k(q,r_k(q))$ denotes the space of $k$-times differentiable functions on the ball of radius $r_k(q)$ centered at $q\in Q$, equipped with the $C^k$-norm.

Let $\cW^{s}[q]$ denote the stable manifolds defined for $\nu$-a.e. $q\in Q$, as in \autoref{ssec:lyapunov_charts_and_stable_manifolds}.
Recall also that cocycles smooth along unstables are defined in \autoref{def:regularity_along_stables}.
In this setting, the main goal of this section is:

\begin{theorem}[Linearization of stable dynamics]
	\label{thm:linearization_of_stable_dynamics_single_diffeo}
	With the setup as above, there exists full $\nu$-measure sets \index{$Q_{00}$}$Q_{00}\subset Q_0$ such that the points in $Q_{00}$ are Oseledets-biregular, while \index{$Q_0$}$Q_0$ is the saturation of $Q_{00}$ by stable manifolds, and a cocycle $L\cW^s\to Q_0$ which is smooth along stable manifolds, with the following properties.
	\begin{enumerate}
		\item \textbf{Linearization:} For $q\in Q_0$ there
                  exist smooth embeddings\index{$L_q$}
		\[
			\cW^s[q] \xrightarrow{L^-_q}L\cW^s(q)
		\]
		of the unstable manifolds into the linear spaces $L\cW^s(q)$, and the embeddings are compatible with the (nonlinear) dynamics on $\cW[q]$ and linear cocycle dynamics on \index{$L\cW^s(q)$}$L\cW^s(q)$, as in \autoref{eqn_cd:linearization_commutes}.
		\item \textbf{Flat connection:} The bundle $L\cW^s$ is equipped with a smooth flat connection \index{$P_{GM}$}$P_{GM}$ along the stable manifolds.
		The embeddings $L^-_q$ are compatible with the flat connection, i.e.
		\[
			P_{GM}(q,q')\circ L^-_q  = L^-_{q'} \text{ as maps on }\cW^s[q]=\cW^s{[q']}
		\]
		when $q'\in \cW^s[q]$ and $q,q'\in Q_0$.
		\item \textbf{Holonomy:} The flat connection $P_{GM}$ induces a holonomy on the bundle $L\cW^s$, in the sense of \autoref{sssec:holonomies_discussion}.
		It agrees with the standard measurable connection $P^-$ (see \autoref{sssec:standard_measurable_connection}) on the associated graded bundle, and therefore differs from the standard measurable connection on $L\cW^s$ itself by a nilpotent map that strictly lowers the Lyapunov filtration.
	\end{enumerate}
\end{theorem}

\begin{proof}
	The results are a consequence of \autoref{thm:subresonant_normal_form_suspended_systems} which equips the stable manifolds with a subresonant manifold structure as in \autoref{def:subresonant_structure_on_a_manifold}.
	The properties developed in \autoref{prop:linearization_of_subresonant_maps} yield the compatibility properties.

	To deduce compatibility of holonomies, observe that the associated graded of $L\cW^s$ can be expressed in terms of tensor constructions on the associated graded of the tangent cocycle $W^s$, and tensor constructions are compatible with taking standard holonomies.
\end{proof}

\begin{remark}
	\label{rmk:equivariant_section}
	\leavevmode
	\begin{enumerate}
		\item 
		We can equivalently state the linearization part of \autoref{thm:linearization_of_stable_dynamics_single_diffeo} as follows.
		There exists an equivariant section \index{$\iota_s$}$\iota_s\colon Q_0\to L\cW^s$ which is smooth along stable manifolds.
		It's relation to the map $L_q^-$ is that
		\[
			L^-_q(q') = P_{GM}(q',q)\big(\iota_s(q') \big)
		\]
		Note that by construction, the image of $\iota_s$ will be contained in an affine hyperplane and in particular $\iota_s$ will never vanish.

	\item In \autoref{thm:subresonant_normal_form_suspended_systems} we will describe parametrizations $\cN_q$ of $\cW^s[q]$ by the tangent space $W^s(q)$.
	The embeddings $L^-_q$ are polynomial maps when precomposed with these parametrizations.
	\end{enumerate}
\end{remark}

In the main body of the text, the primary application will be to the case of \emph{unstable}\index{$L\cW^u(q)$}\index{$L_q$} manifolds and we will denote the corresponding maps
\begin{align}
	\label{eqn:def_of_L_unstable_linearization}
	\cW^u[q]\xrightarrow{L_q} L\cW^u(q).
\end{align}
Note that there is no superscript.

\subsubsection{Skew-product variant}
	\label{sssec:assumptions_normal_forms_skew_product_variant}
Let $(\Omega,\nu_{\Omega},T)$ be an ergodic probability measure-preserving system.
Suppose $Q$ is a manifold and \index{$Q$@$\hat{Q}$}$\hat{Q}:=\Omega\times Q$ is equipped with a map \index{$T$@$\hat{T}$}$\hat{T}$ for which the map $T$ on $\Omega$ is a factor, and furthermore on the fibers the maps \index{$T$@$\hat{T}_\omega$}$\hat{T}_\omega\colon Q\to Q$ are diffeomorphisms.
Let \index{$\nu$@$\hat{\nu}$}$\hat{\nu}$ be a $\hat{T}$-invariant ergodic probability measure.

Suppose, furthermore, that the derivative cocycles in the $Q$-direction, denoted \index{$D^{(k)}\hat{T}_\omega(q)$}$D^{(k)}\hat{T}_\omega(q)$, satisfy the assumptions analogous to \autoref{sssec:assumptions_for_normal_forms_single_diffeo}.
Let then \index{$W$@$\cW^s[\hat{q}]$}$\cW^s[\hat{q}]:=\cW^{s}[\omega,q]\subset Q$ be the associated stable manifolds.
We then have a result analogous to \autoref{thm:linearization_of_stable_dynamics_single_diffeo}, with an extra property in case $\Omega=S^{\bZ}$ is a Bernoulli system and $\hat{T}_{\omega}$ depends only on the $0$-th coordinate (and we follow the notation as in \autoref{ssec:random:skew:product}):

\begin{theorem}[Linearization of dynamics for skew products]
	\label{thm:linearization_of_dynamics_for_skew_products}
		There exists full $\hat{\nu}$-measure sets $\hat{Q}_{00}\subset\hat{Q}_0$, such that points in \index{$Q$@$\hat{Q}_{00}$}$\hat{Q}_{00}$ are Oseledets-biregular and \index{$Q$@$\hat{Q}_0$}$\hat{Q}_0$ is the saturation of $\hat{Q}_0$ by the stable manifolds $\cW^{s}[\hat{q}]$, and a cocycle \index{$L$@$L\cW^s$}$L\cW^s\to \hat{Q}_{0}$ which is smooth along stable manifolds, with the following properties.
	\begin{enumerate}
		\item \textbf{Linearization:} For $q\in \hat{Q}_{0}$ there exist smooth embeddings
		\[
			\cW^s[\hat{q}] \xrightarrow{L^-_{\hat{q}}}L\cW^s(\hat{q})
		\]
		of the unstable manifolds into the linear spaces $L\cW^s(\hat{q})$, and the embeddings are compatible with the (nonlinear) dynamics on $\cW[\hat{q}]$ and linear cocycle dynamics on $L\cW^s(\hat{q})$, as in \autoref{eqn_cd:linearization_commutes}.
		\item \textbf{Flat connection:} The bundle $L\cW^s$ is equipped with a smooth flat connection $P_{GM}$ along $\cW^s[\hat{q}]$.
		The embeddings \index{$L^-_{\hat{q}}$}$L^-_{\hat{q}}$ are compatible with the flat connection, i.e. if $\hat{q}=(\omega,q)$ and $\hat{q}'=(\omega,q')$ are in $\hat{Q}_{0}$ with $q'\in\cW^{s}[\hat{q}]$ then:
		\[
			P_{GM}(\hat{q},\hat{q}')\circ L^-_{\hat{q}}  = L^-_{\hat{q}'} \text{ as maps on }\cW^s[\hat{q}]=\cW^s{[\hat{q}']}.
		\]
		\item \textbf{Holonomy:} The flat connection $P_{GM}$ induces a holonomy on the bundle $L\cW^s$, in the sense of \autoref{sssec:holonomies_discussion}.
		It agrees with the standard measurable connection $P^-$ (see \autoref{sssec:standard_measurable_connection}) on the associated graded bundle, and therefore differs from the standard measurable connection on $L\cW^s$ itself by a nilpotent map that strictly lowers the Lyapunov filtration.
		\item \textbf{Changing the past:} Suppose that $\Omega=S^{\bZ}$ and let \index{$W$@$\cW^{s}_{loc}[\omega]$}$\cW^{s}_{loc}[\omega]$ denote the set $S^{\bZ_{<0}}\times \omega_+$ when $\omega=(\omega_-,\omega_+)\in S^{\bZ}$.
		Then if $\hat{q}=(\omega,q)$ and $\hat{q}'=(\omega',q)$ with $q'\in \cW^s_{loc}[\omega]$ we have that $\cW^{s}[\hat{q}]=\cW^s[\hat{q}']$ as submanifolds of $Q$, and the Gauss--Manin connection and linearization are compatible with this identification.
		\end{enumerate}
\end{theorem}
The proof of this result is analogous to \autoref{thm:linearization_of_stable_dynamics_single_diffeo} and follows from the formalism that we describe next.

\subsubsection{Statement of subresonant normal form}
	\label{sssec:statement_of_subresonant_normal_form}
The next result is an extended version of \cite[Thm. 2.5]{KalininSadovskaya2017_Normal-forms-for-non-uniform-contractions} but stated for a suspended dynamical system; in loc. cit. the base $B$ is a point.
Note that in our statement, we parametrize the manifold by the tangent space, whereas \cite{KalininSadovskaya2017_Normal-forms-for-non-uniform-contractions} use the inverse map.

Keep the notation and assumptions of \autoref{sssec:assumptions_normal_forms_skew_product_variant}, with $\cW^s[\hat{q}]$ denoting the fiberwise stable manifolds, with tangent space \index{$W^s(\hat{q})$}$W^s(\hat{q})\subset T_q Q$.
Equip the tangent space $W^s(\hat{q})$ with the subresonant structure coming from the Lyapunov filtration and exponents.
For general $B$ set $\cW^s_{loc}[b]:=b$ while for $B=S^{\bZ}$ set $\cW^s_{loc}[b]:=S^{\bZ_{<0}}\times b_+$ as in \autoref{thm:linearization_of_dynamics_for_skew_products}.

\begin{theorem}[Subresonant normal form, suspended systems]
	\label{thm:subresonant_normal_form_suspended_systems}
	There exists set $\hat{Q}_{00}\subset \hat{Q}_0$ of full $\hat{\nu}$-measure, such that points in $\hat{Q}_{00}$ are Oseledets-biregular while $\hat{Q}_0$ is its saturation by $\cW^s[\hat{q}]$-manifolds, as well as a tempered radius functions $r_k(\hat{q})$, such that for $\hat{q}\in \hat{Q}_{00}$ there exist smooth maps\index{$N$@$\cN_{\hat{q}}$}
	\[
		\cN_{\hat{q}}\colon W^s(\hat{q}) \to \cW^s{[\hat{q}]}
	\]
	which have derivative at the origin equal to the identity and for any $k\geq 1$ the derivative norms $\norm{\cN_{q}}_{C^{k}(B(0,r_{k}(\hat{q})))}$ are tempered.

	Equip $W^s(q)$ with a subresonant structure coming from the Oseledets filtration and Lyapunov exponents.
	Then the charts $\cN_{\hat{q}}$ induce a subresonant structure on $\cW^s{[\hat{q}]}$ in the sense of \autoref{def:subresonant_structure_on_a_manifold}, compatible with the dynamics in the following sense:
	\begin{enumerate}
		\item \textbf{Subresonant polynomial dynamics:} For any $\hat{q}\in \hat{Q}_{00}$ the composition
		\[
			\cN_{\hat{f}(\hat{q})}^{-1} \circ \hat{f} \circ \cN_{\hat{q}} 
			\colon
			W^s{\left(\hat{q}\right)}
			\to
			W^s{\left(\hat{f}(\hat{q})\right)}
		\]
		is a subresonant polynomial map.
		\item \textbf{Consistency of charts:} For any $\hat{q}=(\omega,q)$ and $\hat{q}'=(\omega',q')$ (with $\hat{q},\hat{q}'\in \hat{Q}_{00}$) such that $\omega'\in \cW^s_{loc}[\omega]$ and $q'\in \cW^s{[\hat{q}]}$, the composition:
		\[
			\cN_{\hat{q}'}^{-1}\circ \cN_{\hat{q}}
			\colon
			\cW^s{(\hat{q})}\to
			\cW^{s}{(\hat{q}')}
		\]
		is a subresonant polynomial map.
		\item \textbf{Uniqueness:} Any other measurable family of maps that satisfies the same assumptions as $\cN_{\hat{q}}$ differs from it by precomposition with a measurable family of subresonant automorphisms of $W^s{(\hat{q})}$.
		Conversely, precomposing with any measurable family of subresonant maps, with identity derivative at the origin, yields a new family of $\cN_{\hat{q}}$'s.
	\end{enumerate}
\end{theorem}
\begin{proof}
	The existence of the maps $\cN_{\hat{q}}$ is a direct consequence of \cite[Thm.~2.3]{KalininSadovskaya2017_Normal-forms-for-non-uniform-contractions}, which gives maps that conjugate the dynamics along stables to subresonant polynomial maps, and also yields the uniqueness.
	The only additional property that we need to verify is consistency of charts.
	This can proved analogously to \cite[Thm.~2.5]{KalininSadovskaya2017_Normal-forms-for-non-uniform-contractions}, which involves the case of a single diffeomorphism, but the proof in the suspended case is virtually identical, since all estimates are on individual stable manifolds.
\end{proof}



\subsection{Measurable connections and subresonant maps}
	\label{ssec:measurable_connections_and_subresonant_maps}

\subsubsection{Setup}
	\label{sssec:setup_measurable_connections_and_subresonant_maps}
We keep the notation and assumptions as in \autoref{ssec:subresonant_normal_forms_on_stable_manifolds}, but for simplicity treat just the case of a single diffeomorphism $f\colon Q\to Q$.
The statements below in the case of a suspended system only require the addition of $\hat{\bullet}$.
For compatibility with the main text, we now switch to \emph{unstable} manifolds and normal forms on them, so we consider the linearizations\index{$L\cW^u(q)$}\index{$L_q$}
\[
	\cW^u[q]\xrightarrow{L_q} L\cW^u(q).
\]

\begin{proposition}[Subresonant maps from holonomies]
	\label{prop:subresonant_maps_from_holonomies}
	Let \index{$P^-$}$P^-$ and \index{$P^+$}$P^+$ denote the measurable stable and unstable connections on $L\cW^u$ defined in \autoref{sssec:standard_measurable_connection}.
	Suppose that $q,q',y\in Q_{00}$, and that $y\in \cW^u[q],q'\in \cW^s[q]$.
	\begin{enumerate}
		\item The map $P^{+}(y,q)$ is induced by a strictly subresonant map \index{$P$@$\wp^+(y,q)$}$\wp^+(y,q)\in \bbG^{ssr}(\cW^u[q])$ of the unstable $\cW^u[q]=\cW^u[y]$, taking $y$ to $q$.
		\item The map $P^-(q,q')$ is induced by a resonant map \index{$P$@$\wp^-(q,q')$}$\wp^-(q,q')$ between $\cW^u[q]$ and $\cW^u[q']$, where the resonant structures are induced by the Lyapunov splittings at $q$ and $q'$ respectively.
	\end{enumerate}
\end{proposition}
\begin{proof}
	By Ledrappier's invariance principle, the measurable holonomies $P^+$ and $P^-$ respect the algebraic hulls of the cocycles.
	In particular, both induce subresonant maps $\wp^{\pm}$ between corresponding manifolds, by \autoref{cor:polynomial_maps_between_different_subresonant_manifolds}.
	In fact, for the same reason both maps are resonant, for the resonant structures induced by the Oseledets decompositions at the respective endpoints.

	To see that $\wp^+(y,q)$ is strictly subresonant, note that $P^+$ acts trivially on the associated graded of $L\cW^u$ by construction, since it agrees with the Gauss--Manin connection there.
	It follows $\wp^+(y,q)$ is strictly subresonant, by the semi-direct product structure of \autoref{prop:group_structure_on_subresonant_maps}\autoref{item:sresn_gp_semidir}.
\end{proof}

We will make use of the following compatibility between a family of invariant subgroups as in \autoref{def:compatible_family_of_subgroups}, and the measurable connection $P^+$ on $L\cW^u$ defined in \autoref{sssec:standard_measurable_connection}.
We will denote by \index{$P^+_{ad}$}$P^+_{ad}$ the induced (adjoint) action on $\End(L\cW^u)$ by conjugation, and so on the cocycle $\frakg^{sr}\subset \End(L\cW^u)$ of Lie algebras of subresonant maps.

\begin{proposition}[Compatibility of measurable connection and the family of subgroups]
	\label{prop:compatibility_of_measurable_connection_and_the_family_of_subgroups}
	There exists a set \index{$X_{00}$}$X_{00}\subset X$ of full measure, such that $U^+$ is defined on $X_{00}$ and with the following property.
	Suppose that $q,q'\in X_{00}$ and $q'\in \cW^u[q]$.

	Then the strictly subresonant map $\wp^+(q,q')\in \bbG^{ssr}(\cW^u[q])$ from \autoref{prop:subresonant_maps_from_holonomies} satisfies $\wp^+(q,q')U^{+}[q]=U^+[q']$.
	Analogously, the measurable connection satisfies
        $P^+_{ad}(q,q')\fraku^+(q)=\fraku^+(q')$.
\end{proposition}
\begin{proof}
	The existence of the compatible family of subgroups $U^+(q)\subset \bbG^{ssr}(q)$ implies that it is preserved pointwise by the algebraic hull of the cocycle on $L\cW^u$, viewed as a subgroup of $\bbG^{sr}(q)$.
	The claim then follows again from the Ledrappier invariance principle \autoref{thm:ledrappier_invariance_principle}.
\end{proof}

We record, for future use, the following elementary consistency result for the structures constructed above.

\begin{proposition}[Agreement of linearizations and compatible family of subgroups]
	\label{prop:agreement_of_linearizations_and_compatible_family_of_subgroups}
	Fix $n_u:=\dim W^u$ and subresonant structure on $(\bR^{n_u},0)$, with the same Lyapunov exponents as on $W^u$ with $\bR^{n_u}=\oplus \bR^{n_{\lambda_i}}$, and equip $\bR^{n_u}$ with its standard Euclidean norm.
	Consider its linearization $L_0\colon \bR^{n_u}\into L\bR^{n_u}$, and the associated direct sum decomposition on $L\bR^{n_u}$.
	Then there exists a set of full measure $Q_{00}$ with the following properties:
	\leavevmode
	\begin{enumerate}
		\item There exists a linear map $T(q)\colon L\cW^u(q)\toisom L\bR^{n_u}$ and resonant map $\tau(q)\colon \cW^u[q]\to \bR^{n_u}$, depending measurably on $q$, such that $T(q)$ is the linearization of $\tau(q)$, such that $T(q)$ maps $L_q(\cW^u[q])$ to $L_0(\bR^{n_u})$ compatibly with the linearization maps (as in \autoref{eqn_cd:linearization_commutes}), and such that $T(q)$ maps the Oseledets decomposition of $L\cW^u(q)$ to the corresponding decomposition on $L\bR^{n_u}$.
		\item If, in addition, $Q$ admits a family $U^+$ of subgroups compatible with the measure as in \autoref{def:compatible_family_of_subgroups} on the finite cover $X$, then $\tau$ can be lifted to $X_{00}\subset X$ and there exists $U_0^+\subset \bbG^{ssr}(\bR^{n_u})$ such that $\tau(q)(U^+(q))=U_0^+$ for $q\in X_{00}$.
	\end{enumerate}
	In particular, for any $\delta>0$ there exists a compact set $K\subset X$ of measure at least $1-\delta$ such that $\norm{T}\leq C(\delta)$ on $K$.
\end{proposition}
\begin{proof}
	The existence of $T(q)$ follows from \autoref{prop:algebraic_hull_and_reduction_of_structure}, and the existence of $\tau(q)$ is a consequence of the fact that the algebraic hull on $L\cW^u$ is contained in the group of resonant maps $\bbG^{r}(\cW^u[q])$, with resonant structure induced by the Oseledets decomposition at $q$.
	Since the algebraic hull preserves the family $U^+$, it follows that $T,\tau$ preserve it as well under the measurable trivialization adapted to the algebraic hull.
\end{proof}




\subsection{Cocycle normal forms}
	\label{ssec:cocycle_normal_forms}

\subsubsection{Setup}
\label{sssec:setup_cocycle_normal_forms}
We keep the notation and terminology from \autoref{ssec:subresonant_normal_forms_on_stable_manifolds}, and refer to \autoref{appendix:generalities_on_cocycles} for general notions related to cocycles.
The construction of holonomies is carried out in \autoref{cor:holonomies_on_graded}.
The existence of cocycle normal forms follows an analogous scheme, but taking into account higher order Taylor expansions.

\subsubsection{Cocycles smooth along stable manifolds}
	\label{sssec:cocycles_smooth_along_stable_manifolds}
Fix a parametrization of the stable manifolds, say the one given by normal forms in \autoref{thm:subresonant_normal_form_suspended_systems}.
A measurable cocycle $E\to Q$ will be called \emph{smooth along stable manifolds} if there exists a set of $Q'$ of full measure (which can be assumed invariant by the dynamics), and trivializations of $E$ on countably many measurable sets $Q_\alpha'$ with $\cup_\alpha Q_{\alpha}'=Q'$, with the following properties.
The sets $Q_{\alpha}'$ are $\cW^{s}_{loc}$-saturated, i.e. if $q\in Q_{\alpha}'$ then $\cW^s[q]_{loc}\subset Q_{\alpha}'$.
The gluing maps of the cocycle, and the cocycle maps, when restricted to any local stable manifold, are smooth maps with bounds on derivatives that depend only on the index of chart.

In particular, we can consider the induced family of trivializations:
\[
	\cE_q\colon E(q)\times \cW^s_{loc}(q)\to E\vert_{\cW^s_{loc}[q]}
\]
varying measurably for $q\in Q'$, with the following additional properties.
If $q'\in \cW^s_{loc}[q]\cap Q'$ then on the overlap we have the estimate
\[
	\norm{\cE_{q'}^{-1}\circ \cE_{q}}_{C^k(B(q,r(q))\cap B(q',r(q')))} \leq c_k(q)
\]
where $c_k(q)$ is a tempered in $q$ measurable function.
Additionally, we assume that the induced cocycle maps
\begin{align*}
	A_q &\colon E(q)\times \cW^s_{loc}(q)
	\to
	E{(\hat{f}(q))}\times \cW^s_{loc}(\hat{f}(q))&& \text{obtained from}\\
	\wtilde{\cA}_q &\colon E\vert_{\cW^s_{loc}[q]}
	\to
	E_{\hat{f}(q)}\vert_{\cW^s_{loc}[\hat{f}(q)]}
\end{align*}
satisfy $\norm{A_q}_{C^k(B(q,r(q)))}\leq a_k(q)$ for a tempered measurable function $a_k$.

\subsubsection{Smoothness of Oseledets filtration along stables}
	\label{sssec:smoothness_of_oseledets_filtration_along_stables}
Let us recall that by \cite[Rmk. 5.2b)]{Ruelle1979_Ergodic-theory-of-differentiable-dynamical-systems} as well as Thm. 5.1 in op. cit., the stable Oseledets filtrations of smooth and natural cocycles, in the sense of \autoref{def:smooth_and_natural_cocycles}, are smooth along stable manifolds.

Below, \emph{subresonant maps of vector bundles} are understood in the sense of \autoref{sssec:subresonant_vector_bundles}, and for the Lyapunov exponents as described in \autoref{sssec:preliminary_notation_for_the_proof}.

\begin{theorem}[Subresonant normal forms, cocycles]
	\label{thm:subresonant_normal_forms_cocycles}
	Let $E\to Q$ be a cocycle smooth along the stable manifolds and let $\cN_q$ be normal form coordinates as in \autoref{thm:subresonant_normal_form_suspended_systems}.
	There exists a set $Q_0$ of full measure, such that for any $q\in Q_0$ there exists a measurable family of smooth trivializations
	\[
		\cN\cE_q \colon E(q)\times \cW^s_{loc}{(q)}\to E\vert_{\cW^s_{loc}[q]}
	\]
	which are linear on the fibers (and the identity on the fiber at $0$), cover the normal form coordinates given by $\cN_q$, and have the following properties.
	\begin{enumerate}
		\item Write $\wtilde{\cA}_q$ for the cocycle maps and let $NA_q$ be defined by the diagram:
		\begin{equation}
			\label{eqn_cd:cocycle_normal_forms}
		\begin{tikzcd}
			E(q)\times \cW^{s}_{loc}{(q)}
			\arrow[r, "\cN\cE_{q}"]
			\arrow[d, "NA_q", dashed, swap ]
			& 
			E\vert_{\cW^s_{loc}{[q]}}
			\arrow[d, "\wtilde{\cA}_q"]
			\\
			E({\hat{f}(q)})\times \cW^s_{loc}{(\hat{f}(q))}
			\arrow[r, "\cN\cE_{\hat{f}(q)}"]
			&
			E\vert_{\cW^s_{loc}{[\hat{f}(q)]}}
		\end{tikzcd}
		\end{equation}
		i.e. by the formula $NA_q=\cN\cE_{\hat{f}(q)}^{-1}\circ \wtilde{\cA}_q \circ \cN\cE_q$.
		Then $NA_q$ is a subresonant map of vector bundles.
		\item For any $q'\in Q_0$ with $q'\in \cW^s_{loc}[q]$ the change of coordinates map
		\[
			\cN\cE_{q'}^{-1}\circ \cN\cE_q
			\colon
			E(q)\times \cW^s_{loc}(q)
			\longrightarrow
			E(q')\times \cW^s_{loc}(q')
		\]
		is a subresonant map of vector bundles, on the subset where it is defined.
		\item Given any other measurable family $\cN\cE_q'$ with the above two properties, the maps $\cN\cE_q^{-1}\circ \cN\cE_{q}'$ are subresonant maps of vector bundles, with tempered $C^k$ norms in $B(0,r(q))$, for any $k$.
		Additionally, if we would use the smooth trivializations of $\cE$ along stables instead of $\cN\cE_q^{'}$, then the resulting coordinate changes would also have tempered $C^k$-norms in $B(0,r(q))$.
	\end{enumerate}
	Therefore, when the stable manifolds are equipped with their subresonant structure as in \autoref{thm:subresonant_normal_form_suspended_systems}, the total space of the cocycle has a natural structure of subresonant vector bundle in the sense of \autoref{sssec:subresonant_vector_bundles}.
\end{theorem}

\begin{corollary}[Holonomy linearization of arbitrary cocycles]
	\label{cor:holonomy_linearization_of_arbitrary_cocycles}
	Let $E\to Q$ be a cocycle smooth along stable manifolds.
	Then there exists an associated cocycle \index{$LE$}$LE\to Q$, also smooth along the unstable, with the following properties.
	\begin{enumerate}
		\item There is a natural embedding
		\[
			E\xrightarrow{ev}LE
		\]
		which is linear on the fibers and is compatible with the dynamics of the cocycles.
		In particular, it is compatible with the Lyapunov filtration.
		\item The bundle $LE$ admits a flat connection \index{$P_{GM}$}$P_{GM}$ along the stable manifolds, which in particular provides smooth holonomies in the sense of \autoref{sssec:holonomies_discussion}.
	\end{enumerate}
\end{corollary}

Note that when $E$ has a single Lyapunov exponent, or is fiber bunched, the above embedding into $LE$ is an isomorphism and the flat connection $P_{GM}$ is the holonomy map usually constructed in such a setting.

\begin{proof}[Proof of \autoref{cor:holonomy_linearization_of_arbitrary_cocycles}]
	The bundle $LE$ is given by the construction in \autoref{sssec:linearization_of_subresonant_bundles} and the embedding is given by \autoref{eq:evaluation_embedding_sresn_vector_bundle}.
	The required properties are a consequence of the subresonant bundle structure on $E$ provided by \autoref{thm:subresonant_normal_forms_cocycles}.
\end{proof}

\begin{remark}[Compatibility of holonomies]
	\label{rmk:compatibility_of_holonomies}
	Because the bundle $LE$ admits a flat connection on stable manifolds, which is compatible with the dynamics, it follows that the standard measurable connection (\autoref{sssec:standard_measurable_connection}) on the associated graded $\gr_\bullet LE$ agrees with the flat connection induced on it.
	Note that the standard measurable connection
    $P^+_{LE}$ on $LE$ is usually \emph{different} from the Gauss--Manin connection $P_{GM}$, since $P^+_{LE}$ involves the full Lyapunov decomposition, not just the filtration.
	However, the difference of maps
	\[
		P_{GM}(q,q')-P^+_{LE}(q,q') \colon LE_q \to LE_{q'}
	\]
	whenever defined, takes a term of the Lyapunov filtration to a strictly lower one.
	This holds because the maps respect the filtrations and agree on the associated graded.

	Returning to the embedding $E\xrightarrow{ev} LE$, it is compatible with the Lyapunov filtration and so it is also compatible with the standard holonomies, as a consequence, for instance, of Ledrappier's invariance principle (it also follows from making explicit the construction of $LE$).
	Note, however, that $E$ is \emph{not} invariant under the flat connection $P_{GM}$ on $LE$.
\end{remark}

\subsubsection{Preliminary notation for the proof}
	\label{sssec:preliminary_notation_for_the_proof}
We will denote the Lyapunov exponents by
\begin{align*}
	-\lambda_1<\cdots < -\lambda_l < 0 && \text{ on the stable manifold}\\
	\mu_1 > \cdots > \mu_k && \text{ on the bundle }E.
\end{align*}
Fix a natural number $D$ and $\ve>0$ sufficiently small, satisfying:
\[
	D(-\lambda_l + \ve) + (\mu_1-\mu_k+2\ve)<-10\ve <0.
\]

\subsubsection{Proof of \autoref{thm:subresonant_normal_forms_cocycles}}
	\label{sssec:proof_of_cocycle_normal_forms}
The proof naturally breaks down into several steps.
Step 1 is the case when $E$ has only one Lyapunov exponent.
Then the statement of the theorem is the existence of holonomies for smooth conformal cocycles, which is well-known, see for instance \cite{AvilaSantamariaViana2013_Holonomy-invariance:-rough-regularity-and-applications-to-Lyapunov-exponents}.

Step 2, which will be the main step, is to assume that the given trivialization of $E$ is already block-diagonal, with cocycles $A_i$, $i=1,\ldots,k$ for each diagonal Lyapunov block, and that additionally all upper-triangular entries of the cocycle are in normal form, except for the deepest one, denoted $A_{1,k}$.
See \autoref{sssec:step_2_deepest_lyapunov_block} for more precise notation.
Then we will show that by a smooth change of coordinates on the cocycle, $A_{1,k}$ is a matrix with polynomial entries of degree bounded by $D$ (see \autoref{sssec:preliminary_notation_for_the_proof} for the defining property of $D$).

Step 2', detailed in \autoref{sssec:step_2_removing_the_finitely_many_terms} will take as input the conclusion of Step 2 and, by a polynomial change of coordinates, bring the cocycle to the desired normal form.

Step 3 consists of applying sequentially steps 2 and 2' to each Lyapunov block of a general cocycle, ordering the blocks $A_{ij}$ according to the lexicographic order on the indexing $(i-j,j)$.
\hfill \qed

\subsubsection{Step 2: deepest Lyapunov block}
	\label{sssec:step_2_deepest_lyapunov_block}
As stated in the proof outline, we assume that our cocycle is in normal form coordinates, except for the block piece $A_{1k}$ which is a smooth function.
By a Taylor expansion of this smooth function, decompose the cocycle matrices in the form:
\begin{align}
	\label{eq:cocycle_decomposition_smooth_error}
	A = PA + R_{1k}
\end{align}
where $PA$ is an upper-triangular matrix with polynomial entries, with all blocks except $PA_{1k}$ in subresonant normal form, and the block $PA_{1k}$ consisting of the Taylor expansion of $A_{1k}$ of degree $\leq D$.
Then $R_{1k}$ is a smooth matrix-valued function satisfying
\begin{align}
	\label{eq:Taylor_estimate_C0}
	\norm{R_{1k}(v)}\leq c_0(q) \norm{v}^{D} \quad \forall v\in \cW^s_{loc}(q)
\end{align}
where $c_0$ is a tempered in $q$ function.
We also have the derivative estimates
\[
	\norm{R_{1k}}_{C^n(B(0,r(q)))}\leq c_n(q) \norm{v}^{\max(D-n,0)}
\]
which will be used to show that the resulting change of coordinates satisfies the claimed differentiability.

\subsubsection{Further notation}
	\label{sssec:further_notation_cocycle_NFs}
Our dynamics is given by $f\colon Q\to Q$ and we will denote by $F_q\colon \cW^s(q)\to \cW^s(f(q))$ the subresonant normal form dynamics provided by \autoref{thm:subresonant_normal_form_suspended_systems}.
Let also $F^{\circ n}_q$ denote its $n$-fold iterate, i.e. $F^{\circ (n+1)}_q:=F_{f^n(q)}\circ F^{\circ n}_q$.
Similarly, we have the fiber-wise linear functions
\[
	A_q^{\circ n}\colon \cW^s(q)\times E(q)\to \cW^s(f^{n}(q))\times E(f^n(q))
\]
and for all the other matrices defined earlier.

We will show that the sequence of transformations
\begin{align}
	\label{eq:definition_cocycle_conjugacy_time_n}
	\cN^{(n)}_{q}(v):= \left(PA_{q}^{\circ n}(v)\right)^{-1}
	\circ
	A_{q}^{\circ n}(v)
\end{align}
converges, exponentially fast, to a $\GL(E_{q})$-valued map $\cN_{q}(v)$ which is smooth in $v$ (and with tempered $C^k$ norms in $B(0,r(q))$).
Assuming convergence, by construction the limit satisfies
\[
	\cN_{{f}(q)}(F_{q}(v)) \circ A_{q}(v)
	=
	PA_{q}(v) \circ \cN_{a}(v)
\]
i.e. it conjugates the cocycle $A$ to the polynomial form $PA$, which by construction is a matrix with subresonant entries everywhere except in the upper-right corner, where it is has a polynomial of degree bounded above by $D$.

To evaluate \autoref{eq:definition_cocycle_conjugacy_time_n}, observe that in the decomposition of the cocycle $A=PA + R_{1k}$ from \autoref{eq:cocycle_decomposition_smooth_error}, when we substitute this into the definition of $A_q^{\circ n}$, and expand, only terms with at most one occurrence of $R_{1k}$ survive, since the image of $R_{1k}$ is the lowest Lyapunov subbundle, its kernel is the second-highest subbundle in the Lyapunov filtration, so the kernel contains the image.
This leads to the expression:
\[
	\cN^{(n)}_{q}(v) = 1 + 
	\sum_{i=1}^n 
	PA^{\circ -(i+1)}_{f^{(i+1)}(q)}\left(F_q^{i+1}(v)\right)\circ
	R_{f^i(q)}\left(F^i_q(v)\right) \circ
	PA^{\circ i}_{q}(v)
\]
To establish the needed estimates, it suffices to show that the individual elements in the sum above decay exponentially fast, and the same for their derivatives to obtain differentiability of the change of coordinates.

A crucial observation which simplifies many calculations is that the term
\[
	\tau_i:=PA^{\circ -(i+1)}_{f^{(i+1)}(q)}\left(F_q^{i+1}(v)\right)\circ
	R_{f^i(q)}\left(F^i_q(v)\right) \circ
	PA^{\circ i}_{q}(v)
\]
only depends on $v$ in the $R$-part, since the diagonal part of the cocycle has been already brought to normal forms and in particular is independent of $v$.
We can thus rewrite the above as
\[
	\tau_i:=
	A_k^{\circ -(i+1)}\left(f^{(i+1)}(q)\right)
	\circ
	R_{f^i(q)}\left(F^i_q(v)\right)
	\circ
	A_1^{\circ i}(q)
\]
where $A_j^{\circ n}(q)$ denotes the $j$-th block on the diagonal of the cocycle, iterated $n$ times starting at $q$.

\subsubsection{Estimating the absolute norm}
	\label{sssec:estimating_the_absolute_norm}
It is convenient to use the Lyapunov norms for $E$ along the orbit of $q$, for which we do not introduce special notation but recall that they depend on a parameter $\ve>0$ which we choose sufficiently small depending on the Lyapunov spectrum.
Then for the individual terms above we have the estimate
\[
	\norm{\tau_i}\leq c_1(q)e^{i\cdot(\mu_1-\mu_k+2\ve)} \norm{R_{f^i(q)}\left(F^i_q(v)\right)}
\]
Combined with the basic contraction estimate
\[
	\norm{F^i_q(v)}\leq c_2(q)e^{-i\cdot (\lambda_l-\ve)}\norm{v}
\]
and the Taylor remainder estimate \autoref{eq:Taylor_estimate_C0} yields, for appropriate slowly-varying constants $c_i(q)$:
\[
	\norm{\tau_i}\leq
	c_1(q)c_3\left(f^i(q)\right)c_2(q)^D 
	e^{i\left(\mu_1-\mu_2 + 2\ve - D(\lambda_l-\ve)\right)} \norm{v}^D
\]
Now since $c_3(q)$ is $\ve$-slowly varying, and by the choice of $D$ in \autoref{sssec:preliminary_notation_for_the_proof}, the $C^0$-convergence of the change of variables follows.

\subsubsection{Estimating the derivatives}
	\label{sssec:estimating_the_derivatives}
In order to estimate $\norm{\tau_i}_{C^n}$ for arbitrary $n$, recall that the dynamics of $F^{i}_q(v)$ is already in polynomial normal forms.
In particular we have an estimate
\[
	\norm{D^{(n)}F^i_q(v)}\leq c_n(q)e^{-i\cdot (n\lambda_l + \ve)}
\]
as soon as $n\geq 1$, where $D^{(n)}$ represents the vector of all $n$-th derivatives, and in fact the derivative vanishes if $n$ is larger than the degree of normal form coordinates on the stable manifold.

Applying now the chain rule to the definition of $\tau_i$, the desired differentiability follows and ends the proof.

\subsubsection{Step 2': removing the finitely many terms}
	\label{sssec:step_2_removing_the_finitely_many_terms}
After the previous step, we are left with a remainder $R_{1k}$ which is a polynomial of large degree $D$.
Choose now a splitting $R_{1k}=R_{1k}^{\leq 0} + R_{1k}^{>0}$ into polynomial functions of respective Lyapunov exponent, for instance by using the Lyapunov splitting at biregular points.
Then the term $R_{1k}^{\leq 0}$ is part of the subresonant normal form, while the exact same process as in the previous step shows that using $R_{1k}^{>0}$ as the remainder term yields again exponentially fast convergence of the conjugation.

\section{Estimates for smooth cocycles and their normal forms}
	\label{sec:estimates_for_smooth_cocycles_and_their_normal_forms}

\subsection{Cocycles: smooth, forward dynamically defined}
	\label{ssec:cocycles_smooth_forward_dynamically_defined}

\subsubsection{Setup}
	\label{sssec:setup_cocycles_smooth_forward_dynamically_defined}

Let $E\to Q$ be a vector bundle, equipped with a lift of the dynamics from $Q$ (e.g. flow, group action, etc.) so it gives a cocycle.
We also assume that $Q$ is equipped with \emph{stable and unstable manifolds} for the dynamics of $g_t$, which will be denoted $\cW^s[q]$ and $\cW^u[q]$ respectively, for $q\in Q$.
Note that in the suspended dynamics case, the stable/unstable manifolds are by convention contained in the fiber direction, so we ignore the stable/unstable manifolds on the base.

\begin{definition}[Smooth and natural cocycles]
	\label{def:smooth_and_natural_cocycles}
	We call a cocycle \emph{smooth and natural} if it is obtained by standard tensor operations (tensor product, direct sum, duality, tensor decompositions) starting from $k$-jets of functions on $Q$, for arbitrary $k \in \bN$, and their subcocycles induced by order of vanishing of a function.

	In addition to these cocycles, we can add to this list cocycles that are smooth on $Q$ and satisfy bounds analogous to those stated in \autoref{ssec:assumptions_on_the_total_space}.
\end{definition}
For instance the $1$-jets of functions, modulo constants, is just the cotangent bundle of $Q$, which is therefore smooth and natural.
We equip smooth and natural cocycles with metrics using the Riemannian metric on $Q$ and the associated Riemannian Levi--Civita connection.
Note that the Lyapunov exponents of any smooth and natural cocycle are determined by those of $TQ$.	

\subsubsection{Local trivializations of smooth and natural cocycles}
	\label{sssec:local_trivializations_of_smooth_and_natural_cocycles}
We will make use of two classes of charts for a smooth and natural cocycle $E$.
First, recall from \autoref{sssec:charts_with_bounds_on_jets} that for every $q\in Q$ we have charts
\[
	\exp_q \colon B(q;r_a(q))\to Q\quad B(q;r_a(q))\subset TQ(q)
\]
subject to the properties listed in \autoref{sssec:systole_function}.
In particular, by possibly shrinking the charts and increasing the function $\phi_a$, we can assume that for every $k\in \bN$, if $\exp_q$ and $\exp_{q'}$ have overlapping image, the $C^k$-norm of the transition map is bounded by $\max(\phi_a(q),\phi_a(q'))^{e'_k}$ for some $e'_k\geq 1$.

Note that such charts induce trivializations of all jet bundles of $Q$, and hence on all smooth and natural cocycles.
We will denote these charts by\index{$e$@$\exp_{q}^E$}
\begin{align}
	\label{eqn:smooth_natural_charts_covering_exp_q}
	\exp_{q}^E\colon B(q;r_a(q))\times E(q) \toisom E\vert_{\exp_q(B(q;r_a(q)))}.
\end{align}
We also extract a locally finite cover of $Q$ from these charts and denote by \index{$\chi_{E,i}$}$\chi_{E,i}$ the associated charts, where $i\in I$ is some indexing set.

\subsubsection{Forward Pesin sets of smooth and natural cocycles}
	\label{sssec:forward_pesin_sets_of_smooth_and_natural_cocycles}
Let $E$ be a smooth and natural cocycle.
For a $k$-tuple $\Lambda=\{\lambda_1>\dots >\lambda_k\}$ and $N,\ve>0$ denote by \index{$Q^{+,\Lambda}_{N,\ve}(E)$}$Q^{+,\Lambda}_{N,\ve}(E)$ the set of $q\in Q$ satisfying for $t\geq 0$:
\begin{align}
	\label{eq:Pesin_set_smooth_natural}
	\begin{split}
	\exists \text{orthog. decomposition }E(q)=\bigoplus_{i=1}^k E_i(q) \text{ such that:}\\
	\forall t\geq N,\quad \forall i\in \{1,\cdots, k\},\quad \forall v\in E_i(q) \\
	e^{t(\lambda_i-\ve)}\tnorm{v}\leq \tnorm{g_t v}\leq e^{t(\lambda_i+\ve)}\tnorm{v}
	\text{ and }
	\lim_{t\to +\infty}\tfrac 1t \log\tnorm{g_t v}=\lambda_i\\
	\text{and }
	\lambda_1 + \dots +\lambda_k = \lim_{t\to +\infty} \tfrac 1t \log\tnorm{g_t\vert_{\det E}}
	\end{split}
\end{align}
Recall that $\tnorm{-}$ denotes the ambient norm on $E$ inherited from the Riemannian metric and natural jet bundle constructions on $Q$.
Let us note that the set of conditions in \autoref{eq:Pesin_set_smooth_natural} is the standard assumption on forward regular points (see \cite[\S1.3.3-7]{BarreiraPesin2007_Nonuniform-hyperbolicity}) and ensures, in particular, the subexponential decay of angles between the images of $E_i$.

We next define:\index{$Q^{+,\Lambda}(E)$}
\[
	Q^{+,\Lambda}(E):=\bigcap_{\substack{\ve\searrow 0}} \bigcup_{N\nearrow +\infty}Q^{+,\Lambda}_{N,\ve}(E)
\]
in other words $q\in Q^{+,\Lambda}(E)$ if for any $\ve>0$ there exists $N(q,\ve)$ such that $q\in Q^{+,\Lambda}_{N(q,\ve),\ve}(E)$.

\subsubsection{Forward Lyapunov bundles}
	\label{sssec:forward_lyapunov_bundles}
For $q\in Q^{+,\Lambda}$, the set of vectors $v\in E(q)$ that satisfy
\[
	\lim_{t\to +\infty} \tfrac 1t \log \tnorm{g_t v} \leq \lambda_i
\]
form a linear subspace, which we denote by \index{$E^{\leq \lambda_i}(q)$}$E^{\leq \lambda_i}(q)$ and will refer to as a forward (or stable) Lyapunov subspace.
For a fixed $\lambda_i$, the set of $E^{\leq \lambda_i}(q)$ over $q\in Q^{+,\Lambda}$ will be called a forward Lyapunov subbundle.
Note that $Q^{+,\Lambda}$ is $g_t$-invariant and the forward Lyapunov subbundles are $g_t$-equivariant.

\begin{proposition}[Stable saturation of Pesin sets]
	\label{prop:stable_saturation_of_pesin_sets}
	Suppose $E\to Q$ is a smooth and natural cocycle and $\nu$ is a $g_t$-invariant ergodic measure, with Lyapunov exponents $\Lambda$ on $E$.

	Then there exists a $g_t$-invariant set of full $\nu$-measure $Q_0$, such that its saturation by stables $\cW^s[Q_0]$ is contained in $Q^{+,\Lambda}$.
	Furthermore, given $k\geq 2$ and sufficiently small $\ve>0$, for $q\in Q_0$ and trivialization of $E$ in the Lyapunov charts $\cL_{k,\ve}[q]$, the composed map
	\[
		W^s_{k,\ve}(q;r_0)\xrightarrow{\exp_q(\id,h^s_q)}\cW^s_{k,\ve}[q]
		\xrightarrow{E^{\leq \lambda_i}}\Gr\left(k_i;E\right)
		\xrightarrow{\left(\exp_{q}^{E}\right)^{-1}}\Gr\left(k_i;\bR^{\rk E}\right)
	\]
	is smooth, with $C^k$-norm bounded by some $C(E,k,\ve)$ (using
        the Lyapunov norm on $W^s(q)$).
\end{proposition}
In the above statement, recall that $\exp_q(\id,h^s_q)$ is defined in \autoref{prop:stableman} as a parametrization of the stable manifold.

\begin{proof}
	The smoothness of the Lyapunov subbundles along stable manifolds is a standard fact, see e.g. \cite[Rmk.~5.2(b)]{Ruelle1979_Ergodic-theory-of-differentiable-dynamical-systems} which also gives the required bounds (recall that in Lyapunov charts with the Lyapunov norms, the dynamics is uniformly hyperbolic, see \autoref{prop:rescaling_metrics_to_size_1}).
	Additionally, for a set of full measure $Q_0$, the conditions defining $Q^{+,\Lambda}_{N,\ve}$ are satisfied for any $q\in Q_0$, but also for any $q'\in \cW^{s}_{loc}[q]\cap \cL_{k,\ve}[q]$ (after possibly increasing $N,\ve$).
	Full saturation by unstables follows from the $g_t$-invariance of $Q^{+,\Lambda}$.
\end{proof}

\renewcommand{\labelenumi}{(\alph{enumi})}

\begin{definition}[Forward dynamically defined cocycle]
		\label{def:forward_dynamically_defined_cocycle}
	A measurable cocycle $F\to Q^{+,\Lambda}$ is called \emph{forward dynamically defined} (abbreviated fdd), and similarly an equivariant map between such cocycles will be called \emph{forward dynamically defined}, if it can be obtained successively from the following starting data and operations:
	\begin{enumerate}
		\item Any smooth and natural cocycle is fdd, and the natural inclusion and quotient maps are fdd.
		\item Any forward (i.e. stable) Lyapunov subbundle of an fdd cocycle is also fdd, and the inclusion map of the Lyapunov subbundle is fdd.
		\item Any standard tensor operation on fdd cocycles is also fdd, and the corresponding maps are fdd.
		This includes tensor products, duality, direct sum, as we all images and quotients.
	\end{enumerate}
	More generally, given a principal $G$-bundle cocycle $P\to Q$, for an affine algebraic group $G$, we will say that $P$ is fdd if, for any linear representation $G\xrightarrow{\rho}\GL(V)$ the linear cocycle $V_\rho:=P\times_G V$ is fdd.
	For a subgroup $H\subset G$, we can then consider the cocycle $P/N\to Q$ whose fibers are isomorphic to $G/H$, which we will also refer to as an fdd fiber bundle cocycle.
\end{definition}

\renewcommand{\labelenumi}{(\roman{enumi})}

\begin{remark}[Basic observations on fdd cocycles]
	\label{rmk:basic_observations_on_fdd_cocycles}
	\leavevmode
	\begin{enumerate}
		\item The construction of fdd cocycles naturally endows them with metrics for which the cocycles are integrable and have Lyapunov exponents (determined by the exponents of the tangent bundle and the construction of the cocycle).
		\item Familiar linear algebra constructions are included in the definition of fdd cocycles.
		For instance the span of two fdd subcocycles $S_i\subset F$ is the image of $S_1\oplus S_2\subset F\oplus F$ under the summation map $F\oplus F\to F$.
		Intersections can be defined using kernels and quotients.
		\item The next remark does not play any role in our proofs: the definition of smooth and natural, as well as fdd cocycles, can be looked at from the point of view of the tannakian formalism.
		The smooth and natural cocycles form a rigid abelian tensor category and the associated pro-algebraic group is the group of formal germs of diffeomorphisms of $(\mathbf{0},\bR^{\dim Q})$.
		The fdd cocycles can be organized in an analogous way, and the corresponding Tannaka group is a subgroup of that of formal germs of diffeomorphisms, determined by the Lyapunov spectrum $\Lambda$.
	\end{enumerate}
\end{remark}

\subsubsection{Aside on Lyapunov filtrations}
	\label{sssec:aside_on_lyapunov_filtrations}
Recall that a linear map of filtered vector spaces $\left(V_1,V_1^{\leq \bullet}\right)\xrightarrow{l} \left(V_2,V_2^{\leq \bullet}\right)$ is \emph{strict} if $l\left(V_1^{\leq \lambda}\right) = l(V_1)\cap V_2^{\leq \lambda}, \forall \lambda$.
For strict linear maps, the kernels and cokernels have well-defined filtrations and the corresponding induced maps are also strict.

Recall also that a linear map is strict if there exist decompositions $V_i = \oplus V_{i}^{\lambda_j}$ compatible with the filtration, such that $l\left(V_1^{\lambda_j}\right)\subseteq V_{2}^{\lambda_j}$.
As a consequence, in the presence of an invariant measure, all dynamics-equivariant maps between cocycles, and associated Lyapunov filtrations, are strict.
In particular Lyapunov filtrations on sub and quotient cocycles are unambiguously defined and are mapped accordingly.

We record for later use:
\begin{proposition}[Fdd cocycles and associated graded]
	\label{prop:fdd_cocycles_and_associated_graded}
	\leavevmode
	Assume that $S,F$ are fdd cocycles and $S\xrightarrow{f}F$ is an fdd map of cocycles, and $\nu$ is an ergodic $g_t$-invariant measure.
	Then there exists a set of full measure $Q_0$ and on its stable saturation $\cW^s[Q_0]$ the following hold.

	Then there exists an induced map on the associated graded cocycles for the filtration by Lyapunov subspaces:
	\[
		gr_{\leq \bullet}S\xrightarrow{gr_{\leq\bullet}f} \gr_{\leq \bullet}F
	\]
	Additionally, this association is exact, in the sense that if $f$ is injective (resp. surjective) then so is $gr_{\leq \bullet}f$.
	Furthermore, the associated graded cocycles and map are also fdd.
\end{proposition}

\begin{proof}
	All equivariant morphisms of cocycles are fiberwise strict for the forward (i.e. stable) filtration, since the underlying filtrations come from a direct sum decomposition of each vector space (into Lyapunov subspaces), and the equivariant morphism of cocycles respects this decomposition.

	All equivariant maps of cocycles are therefore strict and hence induce maps on associated graded pieces.
	Furthermore, for an fdd cocycle $S$ the associated graded cocycle $\gr_{\bullet \leq}S$ is naturally fdd since it comes from the Lyapunov filtration, so each summand is fdd.
	Furthermore, the map induced between summands with the same Lyapunov exponent is itself fdd, again by strictness of the filtration.
\end{proof}

We record an elementary property of the construction of fdd cocycles:

\begin{proposition}[Subquotient representation of fdd cocycles]
	\label{prop:subquotient_representation_of_fdd_cocycles}
	Let $F$ be an fdd cocycle.
	Then there exists a smooth and natural cocycle $E\to Q$, called the smooth envelope of $F$, and two fdd subcocycles $S_1\subseteq S_2\subseteq E$ with an fdd isomorphism $F\toisom S_2/S_1$.

	In particular, $F$ is well-defined on the set $Q^{+,\Lambda}\left(E\right)$.
\end{proposition}

\begin{proof}
	The claim follows by induction on the number of steps from \autoref{def:forward_dynamically_defined_cocycle} needed to construct $F$.
	First, we note that for Lyapunov subbundles of smooth and natural cocycles, the claim is immediate.

	Suppose we already have an fdd cocycle $F=S_2/S_1$ with smooth envelope $E$.	
	Then the Lyapunov subbundle $F^{\leq \lambda_k}$ is the quotient $\left(S_2^{\leq \lambda_k}+S_1\right)/S_1$ while $S_2^{\leq \lambda_k} = S_2 \cap \left(E^{\leq \lambda_k}\right)$.
	So $F^{\leq \lambda_k}$ satisfies the assumptions, with the same smooth envelope and fdd subcocycles as described.
	All other cases of linear-algebraic operations on fdd cocycles are checked similarly, using that all fdd maps between fdd cocycles are strict by \autoref{prop:fdd_cocycles_and_associated_graded}.
\end{proof}

\subsubsection{Norms on fdd cocycles}
	\label{sssec:norms_on_fdd_cocycles}
Given an fdd cocycle $F$, let $S_1\subset S_2\subseteq E$ be a presentation of $F$ as in \autoref{prop:subquotient_representation_of_fdd_cocycles}.
Then $S_1,S_2$ inherit a norm from $E$, whereas $F$ is equipped with the norm obtained by taking the orthogonal complement of $S_1$ inside $S_2$.
This norm will be denoted by \index{$\norm{\bullet}_F$}$\norm{\bullet}_F$ when it needs to be distinguished from other norms, and referred to as the ``fdd norm'' of $F$.

Note that this norm, as well as some of the constructions in \autoref{ssec:estimates_for_smooth_and_fdd_cocycles}, depend on the particular presentation of $F$.
For each fdd cocycle, we fix once and for all one such presentation.
We will refer to $E$ as the \emph{smooth envelope} of $F$.

Let us note, finally, that we can equip $F$ with an $\ve$-Lyapunov norm (see \autoref{sssec:lyapunov_norms}) which is defined on a set of full measure $Q_0$.
On the set $Q_0$, the fdd norm and the $\ve$-Lyapunov norm are comparable up to a tempered function.



\subsection{Estimates for smooth and fdd cocycles}
	\label{ssec:estimates_for_smooth_and_fdd_cocycles}


We now proceed to establish useful estimates for the class of cocycles defined in \autoref{ssec:cocycles_smooth_forward_dynamically_defined}.

\subsubsection{Notation}
	\label{sssec:notation_for_fdd_cocycles}
We set\index{$Q_{k,\ve}[q,r]$}\index{$Q_{a}[q,r]$}
\begin{align*}
		Q_{k,\ve}[q,r] & :=\exp_{q}\left(B_{k,\ve}(q,r)\right)\subset Q \\
		Q_{a}[q,r] & :=\exp_{q}\left(B(q,r)\right)\subset Q
\end{align*}
where $B_{k,\ve}(q,r)\subset TQ(q)$ is the ball of radius $r$ in the $(k,\ve)$-Lyapunov norm, and $B(q,r)$	is the ball in the usual norm.

For a smooth and natural cocycle $E$ and $q\in Q$, we set\index{$Q$@$\wtilde{Q}^{+,\Lambda}_{N,k,\ve}\left( q, r\right)$}
\begin{align*}
	\index{$Q^{+,\Lambda}_{N,k,\ve}[q,r]$}Q^{+,\Lambda}_{N,k,\ve}[q,r] & :=
	Q_{k,\ve}[q,r] \cap Q^{+,\Lambda}_{N,\ve}
	\text{ and }
	\wtilde{Q}^{+,\Lambda}_{N,k,\ve}\left( q, r\right) :=
	\exp_{q}^{-1}\left(Q^{+,\Lambda}_{N,k,\ve}[q,r]\right)
\end{align*}
Moving forward, we will write $F$ for an fdd cocycle and $E$ for its smooth envelope.
The sets $Q^{+,\Lambda}$ and $Q^{+,\Lambda}_{N,\ve}$ will always be associated to $E$.


\begin{definition}[Trivializations]
	\label{def:trivializations}
A \emph{trivialization} of a cocycle $F$ (not necessarily fdd) is the data of a measurable set $Q_0\subset Q$, as well as a measurable family of maps, for all $q\in Q_0$:\index{$\chi_{q}$}
\[
	\chi_{q}\colon F(q)\times S(q) \to F\vert_{S[q]} \quad S(q)\subset TQ(q)
\]
for a (measurable) family of subsets $S(q)\subset TQ(q)$ and $S[q]:=\exp_q(S(q))$, such that $\chi_q$ covers $\exp_q$ on the second factor, and is a linear isomorphism on fibers.

Given a trivialization $\chi$, we will write
\begin{align}
	\label{eqn:G_t_cocycle_trivialization_definition}
	\begin{split}
		\index{$G_{t,q}^{\chi}(s)$}G_{t,q}^{\chi}(s) & \colon F(q) \to F(g_t q)\qquad \quad s\in S(q)\\
		\index{$G_{t,q}^{\chi}(s)v$}G_{t,q}^{\chi}(s)v &:=\chi_{g_tq}^{-1} \circ  g_t\circ \chi_q(v, s) \quad v\in F(q)
	\end{split}
\end{align}
provided that $g_t(\exp_q(s))\in S[g_t q]$.
\end{definition}
In all our constructions, the set $Q_0$ will have full measure.

Recall also that our local stable manifolds are denoted by $\cW^s_{loc}[q]\subset Q$, and we will write \index{$W$@$\wtilde{\cW}^s_{loc}(q)$}$\wtilde{\cW}^s_{loc}(q)\subset TQ(q)$ for their preimage under $\exp_q$.
We also denote by \index{$W$@$\cW^{s}_{k,\ve}[q]$}$\cW^{s}_{k,\ve}[q]$ and \index{$W$@$\wtilde{\cW}^s_{k,\ve}(q)$}$\wtilde{\cW}^s_{k,\ve}(q)$ the corresponding versions in Lyapunov charts.
We recall from \autoref{prop:stableman} that $\wtilde{\cW}^{s}_{loc}(q)$ is identified with a subset of $W^s(q)$ using a graph $v\mapsto (v,h_q^s(v))$.

\begin{definition}[Equivalence of trivializations]
	\label{def:holder_equivalence_of_trivializations}
Suppose that $F$ is an fdd or smooth and natural cocycle, equipped with two trivializations $\chi,\xi$.

We will say that $\chi,\xi$ are \emph{\Holder-equivalent} if there exists $\beta>0$, and for every $N,k,\ve$ a tempered $C_{N,k,\ve}$ such that for a set of full measure $Q_0$, we have for $q\in Q_0$ that $\chi_q,\xi_q$ are defined on $\wtilde{Q}^{+,\Lambda}_{N,k,\ve}(q;r_0)$ and the maps
\[
	\chi_{q}^{-1}\circ \xi_q \colon \wtilde{Q}^{+,\Lambda}_{N,k,\ve}(q;r_0) \to \GL\left(E(q)\right)
\]
have $\beta$-\Holder norms bounded by $C_{N,k,\ve}(q)$, in the $(k,\ve)$-Lyapunov norms.

We will similarly say that the two trivializations are \emph{smoothly equivalent along stables} if the same maps, restricted to $\wtilde{\cW}^s_{k,\ve}(q)$, have $C^k$-norms bounded by $C'_{N,k,\ve}(q)$.
\end{definition}

We recall that the pointwise product of two $\beta$-\Holder functions is $\beta$-\Holder.

\begin{definition}[Regularity along stables]
	\label{def:regularity_along_stables}
	Fix a $g_t$-invariant ergodic measure $\nu$.
	A measurable cocycle $F\to Q$ is called \emph{smooth along stables} if the following hold.
	There exists a set of full measure $Q_0$, such that $F$ is defined on its $\cW^s$-saturation $\cW^s[Q_0]$.

	Furthermore, for every $k\geq k_0(F)$, $\ve\in (\ve_0(F),0)$ there exists  $C=C_{F,k,\ve}$ together with a trivialization $\chi:=\chi_{k,\ve}$, defined for all $q\in Q_0$, with the following property.
	The sets $S(q):=S_{k,\ve}(q)$ contain $\wtilde{\cW}^s_{k,\ve}(q)$, and the maps $G_{t,q}^{\chi}$ restricted to $\wtilde{\cW}^s_{k,\ve}(q)$ have $C^k$-norm bounded by $C$ (for $t\in [-1,1]$).
	To measure the $C^k$-norm, we parametrize $\wtilde{\cW}^s_{k,\ve}(q)$ by $W^s_{k,\ve}(q)$ as in \autoref{prop:stableman}.
\end{definition}


\begin{proposition}[\Holder properties of fdd cocycles]
	\label{prop:holder_properties_of_fdd_cocycles}
	Suppose $F$ is an fdd cocycle, and let $S_1,S_2, E$ be the fdd, resp. smooth and natural, cocycles provided by \autoref{prop:subquotient_representation_of_fdd_cocycles}.
	Then, there exists a set $Q_0$ of full measure, and $k_0(E),\ve_0(E)>0$ with the following properties.
	\begin{enumerate}
		\item
		\label{item:Holder_fdd_cocycles}
		There exist $\beta_1>0$, such that for every $k\geq k_0$ and $\ve\in (0, \ve_0)$, there exist tempered functions $C=C_{S_1,E,N,k,\ve}, r=r_{S_1,E,N,k,\ve}$ such that for every $q\in Q_0$, the map
		\[
			\exp_{q}^{E,-1}
			\left(\exp_q(-),S_1(\exp_q(-))\right)
			\colon \wtilde{Q}^{+,\Lambda}_{N,k,\ve}(q;r(q))\to \Gr\left(\dim S_1;E(q)\right)
		\]
		is \Holder-continuous, with exponent $\beta_1$ and \Holder constant $C(q)$ for the $(k,\ve)$-Lyapunov norm (see \autoref{eqn:smooth_natural_charts_covering_exp_q} for the trivialization $\exp_q^E$).
		
		\item For every $k\geq k_0,\ve\in (0, \ve_0)$ there exists $C=C_{S_1,E,k,\ve}$ such that the map
		\[
			\exp_{q}^{E,-1}
			\left(\exp_q(-),S_1(\exp_q(-))\right)
			\colon \wtilde{W}^s_{k,\ve}(q) 
			\to \Gr\left(\dim S_1;E(q)\right)
		\]
		is bounded in $C^k$-norm by $C(q)$ (using the $(k,\ve)$-Lyapunov norm).
	\end{enumerate}
	The analogous claims hold for $S_2$.
\end{proposition}
Note that the same bounds (with different tempered functions) hold if using the smooth and natural norm on $E(q)$ rather than the $(k,\ve)$-Lyapunov norm, since the two norms are comparable up to a tempered function.
\begin{proof}
	Since we work with the Lyapunov norm, for which the hyperbolicity estimates are uniform, the claim in part (i) when $S_1$ is a forward Lyapunov subbundle follows in a standard manner as in \cite[Thm.~5.3.1, Lemma~5.3.3]{BarreiraPesin2007_Nonuniform-hyperbolicity}.
	When $S_1$ is a smooth and natural subcocycle, the trivialization $\exp^E$ yields a smooth map with control on derivatives on Lyapunov charts.
	For a general $S_1$, one proceeds by induction on the number of steps involved in the construction.

	Part (ii), i.e. smoothness, follows as in \autoref{prop:stable_saturation_of_pesin_sets} from \cite[Rmk.~5.2(b)]{Ruelle1979_Ergodic-theory-of-differentiable-dynamical-systems} for the case when $S_1$ is a Lyapunov subbundle, and in general again by induction.
\end{proof}



\subsubsection{Construction: trivialization $\chi_F$ of an fdd cocycle}
	\label{sssec:construction_trivialization_chi_f_of_an_fdd_cocycle}
We construct first a trivialization $\chi_F$ of any fdd cocycle $F$, using the trivialization of its smooth envelope $E$, and then later in \autoref{eqn:trivializations_fdd} we will define another trivialization of $F$, adapted to the Oseledets filtration.


Suppose that $S_1,S_2,E,F$ are as in \autoref{prop:holder_properties_of_fdd_cocycles}.
Take the orthogonal decomposition $E(q)=S_1(q)\oplus S_{1}^c(q)\oplus S_2^c(q)$ where $S_1(q)\oplus S_1^c(q)=S_2(q)$.
Then in the corresponding flag manifold parametrizing configurations $S_1\subset S_2\subset E(q)$, an open neighborhood of the reference one is given by
\[
	(\phi_1,\phi_2)\in \Hom\left(S_1(q),S_1^c(q)\right) \oplus \Hom\left(S_2(q),S_2^c(q)\right)
\]
where we define $S_2(\phi_1,\phi_2):=\Gamma_{\phi_2}(S_2(q))$ as the graph of $\phi_2$, and then $S_1(\phi_1,\phi_2)$ as $\Gamma_{\phi_1}(\Gamma_{\phi_2}(S_1(q)))$, where as usual \index{$\Gamma_\phi(x)$}$\Gamma_\phi(x):=x\oplus \phi(x)$.

Suppose that $x\in Q^{+,\Lambda}\cap Q[q;r_0]$ with $\exp_q(\wtilde{x})=x$ \emph{and} in the trivialization $\exp^E$ of $E$, the subspaces $\wtilde{S}_i(\wtilde{x})\subset E(q)$ are in the chart constructed above, so that there exist $\phi_i(\wtilde{x})$ with the corresponding properties.
We then use the image $\Gamma_{\phi_2\left(\wtilde{x}\right)}(S_1^{c}(q))$ to identify of $F(x)$ with $F(q)$.

The analogous construction, together with orthogonal projection, yields trivializations
\begin{align}
	\label{eqn:trivializations_E_adapdated_F}
	\index{$I^{E,F}_{q_0,q}$}I^{E,F}_{q_0,q}\colon E(q_0)\to E(q)
	\text{ such that }
	I^{E,F}_{q_0,q}(S_i(q_0))=S_i(q) \quad i=1,2.
\end{align}
which are \Holder and smoothly along stables equivalent to $\exp_q^E$.

\begin{proposition}[Lyapunov-adapted trivializations]
	\label{prop:lyapunov_adapted_trivializations}
	Suppose that $F$ is an fdd cocycle and $\nu$ is a $g_t$-invariant ergodic measure.
	Then there exists a set of full measure $Q_0$ and a trivialization \index{$\chi_q$}$\chi_q$, defined for $q\in Q_0$, and an exponent $\beta>0$, with the following properties.

	There exists $k_0,\ve_0$ such that if $k\geq k_0, \ve \in (0,\ve_0)$, then there exist tempered functions $C=C_{F,N,k,\ve},r=r_{F,N,k,\ve}$ with the following properties.
	\begin{description}
		\item[Forward-regular points are in] 
		The sets $S(q)\subset TQ(q)$ on which the trivialization is given contain $\wtilde{Q}^{+,\Lambda}_{N,k,\ve}(q;r(q))$.

		\item[Lyapunov flag compatibility] We have that
		\[
			\chi\left(F^{\leq \lambda_i}(q)\times s\right)
			= F^{\leq \lambda_i}\left(\exp_q(s)\right)
			\quad
			\forall s\in \wtilde{Q}^{\Lambda,+}_{N,k,\ve}(q,r(q)).
		\]

		\item[\Holder bound] For $s\in [-1,1]$ we have that
		\[
			\norm{G_{s,q}^{\chi}(v_1)-G_{s,q}^{\chi}(v_2)} \leq C(q) \norm{v_1-v_2}^{\beta}
			\quad \forall v_i\in \wtilde{Q}^{\Lambda,+}_{N,k,\ve}(q,r(q)).
		\]
		\item[Norm bound] The map $\chi_q(v)\colon F(q)\to F\left(\exp_q(v)\right)$ has norm bounded by $C(q)$.

		\item[Smoothness along stables] For $s\in[-1,1]$ the induced maps $G_{s,q}^{\chi}$ have $C^k$-norm on $\wtilde{\cW}^s_{k,\ve}(q)$ bounded by $C(q)$.

		\item[Compatibility] The trivializations $\chi$ and $\chi_F$ are \Holder and smoothly along stables equivalent.
	\end{description}
\end{proposition}
\begin{proof}
	The construction is analogous to that of the trivialization $\chi_F$ and $I^{E,F}$ in \autoref{sssec:construction_trivialization_chi_f_of_an_fdd_cocycle}, except that we need to ensure compatibility of forward Lyapunov flags.
	Indeed, we have a filtration of the smooth envelope $E$ of $F$ as
	\[
		S_1\subset SF^{\leq \lambda_n}\subset \dots \subset SF^{\leq \lambda_1}=S_2\subseteq E
	\]
	where $SF^{\leq \lambda_i}$ induce the forward Lyapunov filtration on $F:=S_2/S_1$.
	Now \autoref{prop:holder_properties_of_fdd_cocycles} applies to each of $SF^{\leq \lambda_i}$, so the required transversality holds for points in $Q^{+,\Lambda}_{N,k,\ve}[q;r(q)]$ for an appropriate tempered $r(q)$.
	Furthermore, \autoref{prop:holder_properties_of_fdd_cocycles} yields the required estimates from the construction of $\chi$.
\end{proof}

Going forward, we denote by\index{$I_{x,y}^F$}
\begin{align}
	\label{eqn:trivializations_fdd}
	I_{x,y}^F := \chi\left(-,\exp_x^{-1}(y)\right) \colon F(x)\to F(y)
\end{align}
the linear map obtained from the trivialization provided by \autoref{prop:lyapunov_adapted_trivializations}.
When the cocycle is clear from the context, we will omit the superscript $F$.


We record the following elementary property:
\begin{proposition}[Closeness to the identity]
	\label{prop:closeness_to_the_identity}
	Suppose that $R(x)\colon F(x)\to \bR^n$ is a measurable trivialization of the fdd cocycle $F$ over a set of full measure $Q_0$.
	Then for any $\delta>0$ there exists a compact set $K$ of measure at least $1-\delta$ such that if $x,y\in K$ and $I_{x,y}^F$ is defined then the measurable function
	\[
		I_{x,y}^{R}:=R(y)\circ I_{x,y}^F\circ R(x)^{-1}\in \GL_n(\bR)
	\]
	satisfies $I_{x,y}^{R}\to \id$ as $y\to x$ with all points in $K$.
\end{proposition}
\begin{proof}
	One can extend the measurable trivialization $R$ of $F$ to a compatible trivialization $\wtilde{R}$ of the bundle $E$ used to construct $I_{x,y}$ in \autoref{eqn:trivializations_fdd}.
	Then we can pick the compact set $K$ such that the bundles $S_1,S_2$, as well as the Lyapunov flags of $E$, are given by continuous functions to the Grassmannian of $\bR^n$, and such that the corresponding inner products are continuous as well.
	The claimed property then follows from the construction of $I_{x,y}$.
\end{proof}

%

\subsubsection{Slowly divergent points}
	\label{sssec:slowly_divergent_points}
For the next propositions, until the end of proof of \autoref{distortionest}, we switch to backward dynamically defined (bdd) rather than fdd cocycles, for correspondence with the main text.
Suppose that $Q_0$ is a full measure set such that Lyapunov charts are defined for all points in $Q_0$.
Given $y\in Q_0$ and $z\in \cL_{k,\ve}[y]$, we will say that $y,z$ are $(\beta, \ve,T)$-slowly diverging (with $\beta,\ve>0$) if, introducing the notation $y_t:=g_t y, z_t:=g_t z$ we have the estimates:
\begin{itemize}
	\item $D:=d^{\cL}_{y}(y,z)\leq e^{-\beta T}$.
	\item $d^{\cL}_{y_t}(y_t, z_t )\leq D e^{\ve t}$ for all $t\in [0,T]$
	\item $\ve\leq \beta/10$.
\end{itemize}
In our applications $\beta$ will stay fixed, $\ve>0$ can be chosen arbitrarily small, and $T\to +\infty$.

For the next statement, recall that we reversed time direction and $z$ will now be assumed \emph{backward}-regular.

\begin{proposition}[Slow divergence estimates]
	\label{prop:slow_divergence_estimates}
	Suppose that $y,z$ are $(\beta,\ve,T)$-slowly diverging, and in addition $z\in Q^{-,\Lambda}_{N,k,\ve}$; set $\wtilde{z}_t:=\exp_{y_t}^{-1}(z_t)\in TQ(y_t)$.
	
	Suppose that $E$ is a smooth and natural cocycle and $S\subset E$ is a bdd subcocycle, while $F$ is a bdd cocycle.
	Then, provided $\ve$ is sufficiently small, for each of the assertions below there exists $\beta_1>0$ satisfying:
	\begin{enumerate}
		\item 
		\label{item:slow_divergence_E}
		There exists a tempered $C=C_{E,S,N,k,\beta,\ve}$ such that in the smooth trivialization $\exp^E_{y_t}$ of $E$ at $y$ we have
		\[
			d\left(S(y_t), \wtilde{S}(\wtilde{z}_t)\right)\leq C(y) e^{-\beta_1 T}
			\quad \forall t\in[0,T]
		\]
		for the spherical distance on the Grassmannian of $E(y_t)$, equipped with the $(k,\ve)$-Lyapunov norm.
		\item
		\label{item:slow_divergence_F}
		There exist functions tempered functions $e^{\ell_0},C=C_{F,N,k,\beta,\ve}$ such that if $T\geq \ell_0(y)$, then the trivialization $\chi_F$ from \autoref{sssec:construction_trivialization_chi_f_of_an_fdd_cocycle} is well-defined between $F(y_t)$ and $F(g_t z)$, for all $t\in[0,T]$.
		Furthermore, in the bdd norm on $F$, the corresponding cocycle maps satisfy the bounds
		\[
			\norm{G_{s,y_t}^{\chi}(0)-G_{s,y_t}^{\chi}(\wtilde{z}_t)} \leq 
			C(y) e^{-\beta_1 T}
			\quad \forall t\in[0,T] \quad \forall s\in[-1,1]
		\]
		with notation as in \autoref{eqn:G_t_cocycle_trivialization_definition}.
		Additionally, the induced map $\chi_{y_t}(\wtilde{z}_t)\colon F(y_t)\to F(z_t)$ satisfies
		\[
			\norm{\chi_{y(t)}(\wtilde{z}_t)}\leq 1 + C(y)e^{-\beta_1 T}
			\quad
			\norm{\chi_{y(t)}^{-1}(\wtilde{z}_t)}\leq 1 + C(y)e^{-\beta_1 T}
		\]
	\end{enumerate}
\end{proposition}
\begin{remark}[On estimates from slow divergence]
	\label{rmk:on_estimates_from_slow_divergence}
	Said differently, \autoref{prop:slow_divergence_estimates} asserts estimates for $g_t z$ as if it were in a uniform backwards Pesin set.
	Note that with the assumptions, we have that $\norm{\wtilde{z}_t}\leq e^{-(\beta-\ve) T}$.

	Moreover, we can bound from above $C(y)$ by an $\ve$-slowly varying $\wtilde{C}(y)$ and hence $\wtilde{C}(y_t)\leq \wtilde{C}(y)e^{\ve t}$, so it is immaterial whether we use $C(y)$ or $C(y_t)$ in the statements above.
	Furthermore, the estimates hold for both the $(k,\ve)$-Lyapunov norms and the norm from \autoref{sssec:norms_on_fdd_cocycles}, as we can pass between them at $y$ using a tempered function.
\end{remark}

\begin{proof}[Proof of \autoref{prop:slow_divergence_estimates}]
	We observe deduce the assertion in \autoref{item:slow_divergence_F} from that in \autoref{item:slow_divergence_E}.
	Indeed, fix $S_1\subset S_2\subset E$ and a bdd isomorphism $F\toisom S_2/S_1$.
	Then, if $T$ is sufficiently large (depending on the tempered functions arising in \autoref{item:slow_divergence_E}), it follows that $\wtilde{S}_i(\wtilde{z}_t)$ are sufficiently close to $S_i(y_t)$ that the trivialization $\chi_F$ from \autoref{sssec:construction_trivialization_chi_f_of_an_fdd_cocycle} is well-defined, and furthermore all the spaces involved in the construction satisfy a \Holder bound between $y_t$ and $z_t$.
	Furthermore, since the cocycle on $E$ is smooth, it follows that the induced cocycle on $F$ satisfies the stated \Holder bound.

	For the proof of \autoref{item:slow_divergence_E}, it suffices to prove it for the case when $S$ is given by an element of the Lyapunov filtration, or a smooth and natural subcocycle of $E$.

	If $S$ is a smooth and natural subcocycle of $E$, then in the trivialization $\exp_{y_t}^{E}$ it gives a smoothly varying subspace of $E(y_t)$, and therefore the bound holds with $\beta_1=1$.
	It remains to treat the case when $S\subset E$ is part of the forward Lyapunov filtration.
	Let $\rho_1$ be the smallest Lyapunov exponent in $S$, and $\rho_2<\rho_1$ the largest Lyapunov exponent outside of $S$.

	Let $\cG\to Q$ be the Grassmannian bundle of subspaces of $E$ of dimension $\dim S$, so that $S\colon Q^{-,\Lambda}\to \cG$ is an equivariant section.
	We first check that over $Q_0$, the cocycle is a strict contraction on the fibers, in a neighborhood of $S(x)$, when viewed with the spherical metric on $\cG(x)$ induced by the Lyapunov norm.
	Indeed, an open neighborhood of $S(x)$ in $\cG(x)$ is parametrized by $\Hom(S(x),S(x)^{\perp})$, and the induced cocycle map on this space (identifying $S^{\perp}\isom E/S$) has largest Lyapunov exponent $-\rho_1+\rho_2=:-\delta' <0$.
	Therefore, in the Lyapunov norm for time $1$ the cocycle is a contraction by a factor of $e^{-(\delta'-2\ve)}$.
	If $\ve$ is sufficiently small then $\delta'':=\delta'-2\ve>0$.
	Now the spherical metric on $\cG(x)$, and the ``flat'' metric in the cell $\Hom(S(x),S(x)^{\perp})$, are Riemannian and agree on the tangent space, so there exists $r_0>0$ such that in the ball of spherical radius $r_0$ the metrics are $e^{\delta''/2}$-Lipschitz-equivalent.
	This implies that in the ball of radius $r_0$ in the spherical metric around $S(x)$, the cocycle is uniformly contracting for the spherical metric.

	Next, we use that in the trivialization $\exp_{y}^E$, since $z\in Q^{-,\Lambda}_{N,k,\ve}$, we have from \autoref{prop:holder_properties_of_fdd_cocycles} that
	\[
		d\left(S(y_0),\wtilde{S}\left(\wtilde{z}_0\right)\right)
		\leq C_2(y_0)\norm{\wtilde{z}_0}^{\beta_2}
	\]
	Provided $T$ is sufficiently large (depending on $C_2(y_0)$) so that $\norm{\wtilde{z}_0}$ is sufficiently small, the subspace $\wtilde{S}(\wtilde{z}_0)$ is in the $r_0$-neighborhood of $S(y_0)$.
	We will now establish analogous estimates for discrete times $t\in [0,T]$.

	Set $M_t:=d\left(S(y_t),\wtilde{S}(\wtilde{z}_t)\right)$
	To simplify notation, we write $A_1\left(y_t\right):=G_{y_t,1}(0)$ and $A_{1}\left(z_t\right):=G_{y_t,1}(\wtilde{z}_t)$ for the relevant cocycle maps $E\left(y_t\right)\to E\left(y_{t+1}\right)$.
	From the smoothness of the cocycle on $E$, and the uniformity of the dynamics in the Lyapunov norm, we have the estimate
	\[
		\norm{A_1(y_t)-A_1(z_t)}
		\leq C_3\norm{\wtilde{z}_t}
		\quad \forall s\in [-1,1]
	\]
	for a constant $C_3$.
	We can now compute, using the contraction property of the cocycle along $y_t$:
	\begin{align*}
		M_{t+1} & = 
		d\left(A_1(y_t) S(y_t), A_1(z_t)\wtilde{S}(\wtilde{z}_t)\right)\\
		& \leq 
		d\left(A_1(y_t) S(y_t), A_1(y_t)\wtilde{S}(\wtilde{z}_t)\right) \\
		& + 
		d\left(A_1(y_t)\wtilde{S}(\wtilde{z}_t), A_1(z_t)\wtilde{S}(\wtilde{z}_t)\right)\\
		& \leq e^{-\delta}M_t + 
		\norm{A_1(y_t)-A_1(z_t)}
		\leq e^{-\delta} M_t + C_3\norm{\wtilde{z}_t}
	\end{align*}
	Summing, and assuming $t$ is an integer, we have that:
	\begin{align*}
		M_{t+1}
		& \leq 
		\left(\sum_{k=0}^t
		e^{-k\delta} C_3 \norm{\wtilde{z}_{t-k}}\right)
		+e^{-(t+1)\delta}M_0\\
		& \leq 
		C_3 e^{-\beta T} e^{\ve t}\left(\sum_{k=0}^{t}e^{-k(\delta+\ve)}\right)
		+ C_2(y_0)e^{-\beta_3 T}\\
		& \leq C_4(y_0) e^{-\beta_1 T}
	\end{align*}
	as required.
\end{proof}

\begin{proposition}[Distortion Estimate]
	\label{distortionest}
	Suppose $F$ is a bdd cocycle.
	For any $\ve_1,\beta_1>0$ there exists $\beta_0,\ve_0>0$, and tempered functions $e^{\ell_0}, C = C_{F,N,k,\ve}$ with the following properties.
	First, $\beta_0$ is bounded away from $0$ if $\ve_1\to 0$ and $\beta_1$ is fixed.
	Next, assume $k\geq k_0(F), \ve\in (0,\ve_0)$, and
	suppose that $y\in Q_0$, that $z\in Q^{-,\Lambda}_{N,k,\ve}$, and that $y,z$ are $(\beta_1,\ve,T)$-slowly diverging, and furthermore $T\geq \ell_0(y)$.

	Suppose that a linear map $\phi\colon F(y)\to F(z)$ satisfies
	\[
		\norm{\gr_{\bullet}^{\geq}\phi - \gr_{\bullet}^{\geq}I_{y,z}}
		\leq B.
	\]
	Then we have
	\[
		\norm{g_{T} \circ \gr_{\bullet}^{\geq}\phi \circ g_{-T} - \gr_{\bullet}^{\geq}I_{g_T y,g_T z}}
		\leq C(y)\left(e^{-\beta_0 T}+Be^{\ve_1 T}\right)
	\]
	where the maps inside the norm are between $\gr_{\bullet}^{\geq}F(g_T y)$ and $\gr_{\bullet}^{\geq}F(g_T z)$.
\end{proposition}
The above statement is for the bdd norm introduced in \autoref{sssec:norms_on_fdd_cocycles}.
\begin{proof}
	We note first that it suffices to prove the stated bound using Lyapunov norms, as follows.
	We write for convenience $y_t:=g_t y, z_t:=g_t z$ and $\phi_t:=g_t \circ \phi \circ g_{-t}$.
	Write $\wtilde{\phi}:=I_{y,z}^{-1}\circ \phi\colon F(y)\to F(y)$ and more generally (for $x_s\in \{y_s,z_s\}$):
	\begin{align*}
		\wtilde{\phi}_{t} := I_{y_t, z_t}^{-1}\circ \phi_t & \colon F(y_t)\to F(y_t)\\
		A_t(x_s) := I_{y_{s+t},x_{s+t}} \circ g_t \circ  I_{y_s,x_s}^{-1}\circ & \colon F(y_{s})\to F(y_{s+t})
	\end{align*}
	so that $\wtilde{\phi}_t = A_t(z) \circ \wtilde{\phi} \circ A_{-t}(y_t)$.
	Note that $I_{y_t,y_t}=\id$.
	Now the assumption implies that for a tempered function $C_{k,\ve}$ and $(k,\ve)$-Lyapunov norm we have:
	\[
		\norm{\gr_{\bullet}^{\geq}\wtilde{\phi}-\id}_{k,\ve}\leq B\cdot C_{k,\ve}(y)
	\]
	and if we prove the bound
	\[
		\norm{\gr_{\bullet}^{\geq}\wtilde{\phi}_t-\id}_{k,\ve}
		\leq 
		Be^{\ve_1 T} + e^{-\beta_0 T}
	\]
	then the conclusion follows using the norm bound from \autoref{prop:slow_divergence_estimates}\autoref{item:slow_divergence_F}.
	
	Furthermore, it suffices to assume that $F$ has a single Lyapunov exponent, since each piece of the associated graded is itself bdd, and it suffices to establish the bound on each piece individually.
	In particular, we have
	\[
		e^{(t+s)(\lambda-\ve)} \leq \norm{A_t(y_s)}_{k,\ve} \leq e^{(t+s)(\lambda+\ve)}
	\]
	for some $\lambda\in \bR$.

	We drop the prefix $\gr_{\bullet}^{\geq}$ for the rest of the proof, as well as the subscript $(k,\ve)$ for norms: all norms below are $(k,\ve)$-Lyapunov norms.
	It follows from the estimate in \autoref{prop:slow_divergence_estimates}\autoref{item:slow_divergence_F} that
	\[
		\norm{A_1(z_n)-A_1(y_n)}\leq C'(y)e^{-\beta' T}
	\]
	and therefore that
	\begin{align}
		\label{eqn:dist_est_norm_bound_slow_divergence}
		\begin{split}
		\norm{A_1(z_n)} & \leq \norm{A_1(y_n)}\left(1+C'(y)e^{-\beta' T}\right)\\
		& \leq
		e^{-(\lambda-\ve)}\left(1+C'(y)e^{-\beta' T}\right)
		\end{split}
		\intertext{as well as}
		\label{eqn:AA_inverse_bound}
		\begin{split}
		\norm{A_1(z_n)A_1(y_n)^{-1} - \id}
		\leq C''(y)e^{-\beta' T}.
		\end{split}
	\end{align}
	We evaluate the cocycle at discrete times for convenience and our interest is in estimating $b_n:=\norm{\wtilde{\phi}_n-\id}_{k,\ve}$, and we compute directly:
	\begin{align*}
		\norm{\wtilde{\phi}_{n+1} - \id} & = 
		\norm{A_{1}(z_n) \wtilde{\phi}_{n} A_{-1}(y_{n+1}) -\id}
		\\
		& = \norm{A_{1}(z_n) \left(\wtilde{\phi}_{n}-\id\right) A_{-1}(y_{n+1}) +
		\left[A_{1}(z_n) A_{-1}(y_{n+1}) - \id\right]}\\
		& \leq \norm{A_{1}(z_n)}\norm{\wtilde{\phi}_{n}-\id}\norm{A_{-1}(y_{n+1})} +
		C''(y) e^{-\beta' T}\\
		& \leq b_n \norm{A_{1}(z_n)}\norm{A_{-1}(y_{n+1})} +
		C''(y) e^{-\beta' T}\\
		& \leq b_n e^{3\ve} +
		C''(y) e^{-\beta' T}
	\end{align*}
	using \autoref{eqn:dist_est_norm_bound_slow_divergence} and \autoref{eqn:AA_inverse_bound}, provided $T$ is large enough.
	So we find that $b_n$ satisfies the recurrent bounds:
	\[
		b_{n+1}\leq b_n e^{3\ve} + C''(y) e^{-\beta' T}
	\]
	from which we conclude that
	\[
		b_{T}\leq b_0 e^{3\ve T} + C''(y)e^{-(\beta'-4\ve)T}
		\leq B e^{3\ve T} + e^{-\beta'' T}
	\]
	if $T$ is large enough.
\end{proof}

\subsubsection{Holonomy construction}
	\label{sssec:holonomy_construction}
We return to a fixed fdd cocycle $F$, with smooth envelope $E$ and associated subcocycles $S_1,S_2$.
The associated graded in the next construction is taken for the forward Lyapunov subbundles:\index{$g$@$\gr_{\bullet}F$}\index{$g$@$\gr_{\lambda_i}F$}
\[
	\gr_{\bullet}F:=\bigoplus_{\lambda_i} \gr_{\lambda_i}F \quad
	\text{ and }\quad
	\gr_{\lambda_i}F :=F^{\leq \lambda_i} / F^{<\lambda_i}.
\]
For the next statement, we write $I_{x,y}$ rather than $I_{x,y}^F$ as in \autoref{eqn:trivializations_fdd}.

\begin{proposition}[Holonomy estimates and trivialization]
	\label{prop:holonomy_estimates_and_trivialization}
	There exists $\kappa>0$, $k_0>0$, such that for all $k\geq k_0$ there exists $\ve_k>0$ and tempered functions $C_0=C_{F,\ve},C_k=C_{F,k,\ve},r_k=r_{F,k,\ve}$, such that if $\ve\in(0,\ve_k)$ then the following properties hold.

	There exists a full measure set $Q_0$, such that if $x\in Q_0$ and $y\in \cW^{s}_{k,\ve}[x;r_k(x)]$ we have:
	\begin{enumerate}
		\item For $t\geq 0$ the map\index{$I$@$\wtilde{I}_{x,y}^{t}$}
		\begin{align*}
			\wtilde{I}_{x,y}^{t} & :=
			g_{-t} \circ  I_{g_t x,g_t y}
			\circ g_t 
			\colon
			F(x) \to F(y)
		\end{align*}
		is well-defined (i.e. $I_{g_t x,g_t y}$ is) and for $t\geq 1$ we have the estimate
		\[
			\norm{\gr_{\bullet}\wtilde{I}_{x,y}^t - \gr_{\bullet}\wtilde{I}_{x,y}^{t-1}}
			\leq C_0(x) e^{-{\kappa}t}
		\]

		\item The map
		\[
			y\mapsto \left(\gr_{\bullet}{I}_{x,y}\right)^{-1}
			\circ
			\gr_{\bullet}\wtilde{I}_{x,y}^{t}
		\]
		has $C^k$-norm bounded by $C_k(x)$ on $\cW^{s}_{k,\ve}[x;r_k(x)]$.
		\item The limit map
		\[
			H_{x,y}:= \lim_{t\to +\infty} \gr_{\bullet} I_{x,y}^{t}
		\]
		exists and has the properties listed in \autoref{cor:holonomies_on_graded}.
	\end{enumerate}
\end{proposition}
\begin{proof}
	It suffices to establish the estimates on each individual piece of the associated graded, i.e. we can assume that $F$ has a single Lyapunov exponent $\lambda$.
	We work exclusively with $\ve$-Lyapunov norms on $F(x_\bullet)$, so we have uniform hyperbolicity estimates, and recall again that the fdd and Lyapunov norms are equivalent up to a tempered factor.

	Let $\kappa>0$ be such that for $x$ in a set of full measure $Q_0$:
	\[
		d^{\cL}_{g_t x}(g_t x, g_t y) \leq e^{-\kappa t}d^{\cL}_{x}(x,y)
		\quad \forall y\in \cW^s_{k,\ve}[x], t\geq 0.
	\]
	As in the proof of \autoref{distortionest}, we work with discrete time.
	and introduce the notation $x_n:=g_n x$ and $\exp_{x_n}(\wtilde{y}_n)=y_n$; note that $\wtilde{x}_n=0\in TQ(x_n)$.
	We use the trivialization $I_{x,y}$, as well as the parametrization of $\cW^s_{k,\ve}$ by $W^{s}_{k,\ve}$, to write the cocycle maps:
	\[
		G_{s}(\wtilde{z}_t)\colon F(x_t)\to F(x_{t+s}) \text{ where }
		z\in \{x,y\} \text{ and } \wtilde{z}_{t}\in TQ(x_t).
	\]
	Since $\norm{\wtilde{y}_n}\leq e^{-\kappa n}d^{\cL}_{x_0}(x_0,y_0)$ we find that
	\begin{align*}
		\norm{G_{-1}(\wtilde{y}_{n+1}) \circ G_{1}(\wtilde{x}_n) - \id}
		& \leq 
		\norm{\big[G_{-1}(\wtilde{y}_{n+1}) - G_{-1}(\wtilde{x}_{n+1}) \big]G_{1}(\wtilde{x}_n) }\\
		& \leq 
		\norm{G_{-1}(\wtilde{y}_{n+1}) - G_{-1}(\wtilde{x}_{n+1}) }\norm{G_{1}(\wtilde{x}_n) }\\
		& \leq  C'({x}_{n+1}) e^{-\kappa(n+1)}d^{\cL}_{x_0}(x_0,y_0) e^{\lambda+\ve}\\
		& \leq  \wtilde{C}'(x) e^{-(\kappa-\ve)n} d^{\cL}_{x_0}(x_0,y_0)
	\end{align*}
	where we used that the cocycle is Lipschitz on the stable manifold, with Lipschitz constant a tempered $C'$, and then bounded $C'$ by an $\ve$-slowly varying $\wtilde{C}'$.

	We now also verify that $\norm{G_{-n}(\wtilde{y}_n)}\leq \tfrac 32 e^{-(\lambda - \ve)n}$, provided that $\delta:=d^{\cL}_{x_0}(x_0,y_0)$ is sufficiently small depending on the Lipschitz constant of $G$ (which is a tempered function).
	Namely, we set
	\[
		D_{n}:=\norm{G_{-n}(\wtilde{y}_n) - G_{-n}(\wtilde{x}_n)}
		\text{ so }
		\norm{G_{-n}(\wtilde{y}_n)} 
		\leq D_n + \norm{G_{-n}(\wtilde{x}_n) }
		\leq D_n + e^{-(\lambda-\ve)n}
	\]
	so it suffices to prove $D_n\leq \tfrac 12 e^{-(\lambda-\ve)n}$.
	We compute:
	\begin{multline*}
		D_{n+1}=\norm{G_{-{n+1}}(\wtilde{y}_{n+1}) - G_{-{n+1}}(\wtilde{x}_{n+1})}\\
		 = 
		\norm{
		\big[G_{-1}(\wtilde{y}_{n+1})-G_{-1}(\wtilde{x}_{n+1})\big]G_{-n}(\wtilde{y}_n)
		+
		G_{-1}(\wtilde{x}_{n+1})\big[G_{-n}(\wtilde{y}_{n})-G_{-n}(\wtilde{x}_n)\big]
		}\\
		 \leq 
		C'({x}_{n+1})e^{-\kappa(n+1)}\delta \norm{G_{-n}(\wtilde{y}_n)}
		+ e^{-(\lambda-\ve)}D_n\\
		\leq 
		C'({x}_{n+1})e^{-\kappa(n+1)}\delta \left(D_n + e^{-(\lambda-\ve)n}\right)
		+ e^{-(\lambda-\ve)}D_n
	\end{multline*}
	If we set $d_{n}:=e^{(\lambda-\ve)n}D_n$ and bound $C'({x}_{n+1})$ by an $\ve$-slowly varying $\wtilde{C}'(x_0)e^{-\ve n}$ (and absorb a factor of $e^{\kappa+\lambda}$ into it) we can rewrite the above bound as
	\begin{align*}
		d_{n+1} & \leq \wtilde{C}'(x_0) e^{-(\kappa-\ve)n}\delta(d_n + 1) + d_n\\
		&\leq 
		\left(1 + \wtilde{C'}(x_0)\delta e^{-(\kappa-\ve)n}\right)d_n + 
		\wtilde{C'}(x_0)\delta e^{-(\kappa-\ve)n}\\
		& \leq
		\left(1 + \delta' e^{-\kappa' n}\right)d_n + \delta' e^{-\kappa'n}
	\end{align*}
	where $\delta':=\delta \wtilde{C}'(x_0)$ and $\kappa':=\kappa-\ve$.
	Now $\prod_{n\geq 1} \left(1+\delta' e^{-\kappa' n}\right)$ converges, and goes to $1$ as $\delta'\to 0$, since $\delta'\sum e^{-\kappa' n}$ is finite, so we conclude that
	\[
		d_n \leq \delta' (c_{1,\kappa} d_1 + c_{2,\kappa}) \quad \forall n \geq 2
	\]
	and noting that $d_1$ can itself be estimated linearly in $\delta$, we conclude that $d_n\leq \delta c_{3,\kappa}$ for all $n\geq 1$.
	In particular, if $\delta$ is sufficiently small we have $d_n\leq \tfrac 12$.

	We are now ready to estimate the difference between $\wtilde{I}^{n+1}_{x,y}$ and $I^{n}_{x,y}$.
	In our trivialization the map $\wtilde{I}^{n}_{x,y}$ becomes:
	\[
		\iota_n = G_{-n}(\wtilde{y}_n)\circ G_{n}(\wtilde{x}_0)\colon F(x)\to F(x)
	\]
	so we can estimate the difference as follows:
	\begin{align*}
		\norm{\iota_{n+1}-\iota_n} & =
		\norm{
		G_{-n}(\wtilde{y}_n) 
		\left[G_{-1}(\wtilde{y}_{n+1}) \circ G_{1}(\wtilde{x}_n) - \id\right]
		 G_{n}(\wtilde{x}_0)}
		\\
		& \leq
		\norm{G_{-n}(\wtilde{y}_n) } 
		\norm{G_{-1}(\wtilde{y}_{n+1}) \circ G_{1}(\wtilde{x}_n) - \id}
		\norm{G_{n}(\wtilde{x}_0)}
		\\
		& \leq 
		\tfrac 32 e^{-(\lambda-\ve)n} \wtilde{C}'(x) e^{-(\kappa-\ve) n} d^{\cL}_{x_0}(x_0,y_0) e^{(\lambda+\ve)n}\\
		& \leq
		\tfrac 32 \wtilde{C}'(x)d^{\cL}_{x_0}(x_0,y_0)e^{-(\kappa -5\ve)n}
	\end{align*}
	Using analogous estimates and the chain rule, combined with the fact that the dynamics on stable leaves is uniformly contracting, uniform $C^k$-bounds on $\iota_n$ follow analogously, see e.g. \cite[\S3.2]{AvilaSantamariaViana2013_Holonomy-invariance:-rough-regularity-and-applications-to-Lyapunov-exponents}.
\end{proof}

\begin{corollary}[Holonomies on graded]
	\label{cor:holonomies_on_graded}
	There exists a set of full measure $Q_0$, such that for every $x\in Q_0$, for every $y\in \cW^s[x]$ there is a linear map\index{$H(x,y)$}
	\[
		H(x,y)\colon \gr_{\bullet} F(x)\to \gr_{\bullet} F(y)
	\]
	with the following properties:
	\begin{description}
	\item[Equivariance] We have $H({g_tx},g_t y) \circ  g_t=   g_t\circ H(x,y)$ for all $t\in \R$.

	\item[Normalization] We have $H(x,x)=\id$.

	\item[Boundedness] For every $k$ there is a tempered function $C_k$ such that restricted to $\cW^s_{k,\ve}[x]$, the map
	$$y\mapsto (\gr_\bullet   I_{x,y})^{-1}\circ H(x,y)$$
	is $C^k$ with $C^k$-norm bounded above by $C_k(x)$.

	\item[Uniqueness] Suppose that there exist tempered functions $\wtilde{r},\wtilde{C}$ and another such family of maps $\wtilde{H}_{x}(y)$, defined for $y\in \cW^s[x;\wtilde{r}(x)]$, such that
	\[
		\norm{H(x,y)^{-1}\circ \wtilde{H}(x,y)} \leq \wtilde{C}(x).
	\]
	Then on a further subset of full measure the maps coincide.

	\item[Symmetry] If $x,y\in Q_0$ and $y\in \cW^s[x]$ then $H(x,y) = H(y,x)^{-1}$.
	\end{description}
\end{corollary}
\begin{proof}
	We take
	\[
		H(x,y):=\lim_{s\to+\infty}\wtilde{I}_{x,y}^{s}
	\]
	for $y\in \cW^{s}_{k,\ve}[x]$ and extend it by dynamical equivariance to all $y\in \cW^s[x]$.
	The equivariance, boundedness, uniqueness, and normalization properties follow from the construction and the estimates in \autoref{prop:holonomy_estimates_and_trivialization}.

	Symmetry follows from uniqueness.
	Indeed, we can take a stable Markov partition $\frakS$ and consider the associated conditional measures.
	If they are atomic, there is nothing to prove.
	Otherwise, by a measurable selection theorem we can choose for each atom $\frakS[x]$ a ``center'' $x_*$ and use it to define another family of fiber identifications $\wtilde{H}(x,y):=H({x_*},y)H({x_*},x)^{-1}$ which, by uniqueness, must agree with $H(x,y)$.
\end{proof}

We record the following soft continuity property for holonomy maps, with notation matching its application in \autoref{prop:convergence_at_tilde_points}.
\begin{proposition}[Continuity estimate for Holonomy]
	\label{prop:continuity_estimate_for_holonomy}
	Let $F$ be an fdd cocycle, with holonomy $H$ and trivialization $I_{x,y}$.
	For any $\delta>0$ there exists a compact set $K$ of measure at least $1-\delta$ with the following property.
	Suppose given points $z,\tilde{q},\tilde{q}',q'$ such that $\tilde{q},\tilde{q}',q'\in K$, and for any of these three points, the remaining ones (including $z$) are in its $(k,\ve)$-Lyapunov chart.
	Suppose further that $z\in Q^{+,\Lambda}_{N,\ve}$, and $z\in \cW^{s}_{loc}[q']$, $\tilde{q}\in \cW^s_{loc}[\tilde{q}']$.
	Then as $d^{Q}(z,\tilde{q})+d^{Q}(q',\tilde{q}')\to 0$ we have
	\begin{align}
		\label{eqn:holonomy_graded_estimate}
		\norm{ H\left(\tilde{q},\tilde{q}'\right) -   \gr_\bullet I_{q',\tilde{q}'} \circ H(z,q')\circ \gr_{\bullet} I_{\tilde{q},z}}
		\to 0.
	\end{align}
\end{proposition}
\begin{proof}
	We first note that it suffices to prove the property on each piece of the associated graded, so we can assume that $F=\gr_{\bullet}F$ and drop $\gr_{\bullet}$ going forward.
	Let $E$ be a smooth envelope of $F$ with $F\isom S_2/S_1$ and $S_1\subset S_2 \subset E$.
	We identify $F(q)$ with $S_1^{\perp}(q)\cap S_2(q) \subset E(q)$ in what follows.
	
	Next, as a first step, using the locally finite cover of $Q$ from \autoref{sssec:local_trivializations_of_smooth_and_natural_cocycles}, we can assume that we have finitely many charts covering $K$, and that our points are all in the same chart $U_i=:U\subset Q$, with $\chi_i\colon U\to \hat{U}\subset \bR^{\dim Q}$ denoting the chart map.
	We will write $\chi_E\colon E\vert_{U}\to \hat{U} \times \bR^{m_e}$ for the trivialization of $E$, and $\hat{U}_0:=\chi_i(Q_0\cap U)$, as well as $\hat{U}_0^s:=\chi_i\left(\cW^{s}_{loc}[Q_0]\cap U\right)$.
	For a point $x\in U$ we will write $\hat{x}:=\chi_{i}(x)\in \hat{U}$.

	Fix a measurable family of Lyapunov charts $\gamma\colon \hat{U}_0 \to C^0(B^{cu}\times B^s; \hat{U})$ where $B^{cu},B^s$ are fixed euclidean balls, and the stable manifold through $\hat{x}$ is $\gamma_{x}\left(\mathbf{0}\times B^s\right)$.

	
	Then we can write the above family of subspaces as:
	\begin{align*}
		\hat{F}^s\colon \hat{U}_0 & 
		\to 
		C^1\left(B^s;\Gr(m_f;\bR^{m_e})\right)
		\\
		\hat{F}^s(x) & =
		\Big[
		b\mapsto 
		\chi_{E} \circ {F}\circ \chi_i^{-1}(\gamma_{x}(\mathbf{0},b))
		\Big].
	\end{align*}
	We can assume that $\hat{F}$ is continuous on $K$, so the subspaces are all in a single chart on the Grassmannian, so $\hat{F}^s(x)\in \Hom\left(\bR^{m_f};\bR^{m_f,\perp}\right)$ (for an appropriate choice of $\bR^{m_f}\subset \bR^{m_e}$).
	In particular, we also have a fixed (measurable) identification $F(x)\to \bR^{m_f}$.

	With this data, the charts $I_{x,\bullet}$ correspond to maps
	\[
		\hat{I}_{x,\bullet}\colon \chi_i\left(U\cap Q^{+,\Lambda}_{N,\ve}\right)
		\to \GL\left(\bR^{m_f}\right)
	\]
	that are \Holder-continuous by \autoref{prop:holder_properties_of_fdd_cocycles}, with uniform constant (by the appropriate choice of compact set $K$) and uniform \Holder exponent.
	In particular, $\hat{I}_{x,y}\to \id$ as $y\to x$, provided $y$ stays in a uniform Pesin set.
	In particular, expressing \autoref{eqn:holonomy_graded_estimate} in these charts, the maps $\hat{I}_{\hat{\tilde{q}},\hat{z}}\to \id$ and $\hat{I}_{\hat{q}',\hat{\tilde{q}'}}\to \id$.

	We can also express the holonomy maps using the measurable function:
	\begin{align*}
		\hat{H} \colon \hat{U}_0 & \to C^1(B^s ; \GL\left(\bR^{m_f}\right))\\
		\hat{H}(\hat{x}) & := 
		\Big[
		b\mapsto 
		\chi_E\circ H\left(x,\chi_i^{-1}\gamma_{\hat{x}}(\mathbf{0},b)\right)\circ \chi_{E}^{-1}
		\Big].
	\end{align*}
	We can choose the compact set $K$ such that $\hat{H}$ is continuous on it, so $\hat{H}(\hat{x})\to \hat{H}(\hat{y})$ as $\hat{x}\to \hat{y}$.
	Furthermore, if we denote by $b_z,b_{\tilde{q}}\in B^s$ the points such that $\gamma_{\hat{q}'}(\mathbf{0},b_z)=z$ and $\gamma_{\hat{\tilde{q}}'}(\mathbf{0},b_{\tilde{q}})=\tilde{q}$, then $b_z\to b_{\tilde{q}}$ as $z\to \tilde{q}$.
	This implies that
	\[
		\hat{H}\left(\hat{{q}}'\right)(b_{{z}})
		\to 
		\hat{H}\left(\hat{\tilde{q}}'\right)(b_{\tilde{q}})
	\]
	which yields the claim.
\end{proof}

Recall that given a cocycle $F$ with Lyapunov decomposition $F=\oplus F^{\lambda_i}$, at a biregular point $q$ we have an equivariant isomorphism\index{$\Lambda_F^{\leq}$}
\begin{align}
	\label{eqn:isomorphism_graded_and_cocycle_at_biregular}
	\Lambda_F^{\leq}\colon \gr_{\bullet}^{\leq}F\toisom F
\end{align}
given by the Oseledets decomposition.
We can alternatively view $\Lambda_F^{\leq}$ as a measurable equivariant section of $\Hom\left(\gr_{\bullet}^{\leq} F, F\right)$.
We introduce the notation\index{$Q^{+,-,\Lambda}_{N,\ve}$}
\begin{align}
	\label{eqn:biregular_points_notation}
	Q^{+,-,\Lambda}_{N,\ve}:=Q^{+,\Lambda}_{N,\ve}\cap Q^{-,\Lambda}_{N,\ve}
\end{align}
for the following statement:

\begin{proposition}[\Holder equivalence of cocycle and graded]
	\label{prop:holder_equivalence_of_cocycle_and_graded}
	Suppose that $F$ is an fdd cocycle and $I_{x,y}^{F}$ the trivialization from \autoref{eqn:trivializations_fdd}
	and $I_{x,y}^{\gr_\bullet^{\leq} F}$ the induced trivialization on the associated graded.

	Then there exists a set of full measure $Q_0$ and $\beta>0$, as well as $k_0,\ve>0$, such that if $k\geq k_0,\ve\in(0,\ve_0)$, we have the following properties.
	There exists tempered functions $r, C$ depending on ${F,N,k,\ve}$ such that, for every $q\in Q_0$, the function
	\[
		I_{q,\exp_q(v)}^{F,-1}\circ \Lambda_{F}^{\leq}(\exp_q(v))\circ I_{q,\exp_q(v)}^{\gr_{\bullet}^{\leq}F}
		\colon \gr_{\bullet}^{\leq }F(q) \to F(q)
	\]
	is defined on $Q_0\cap Q^{+,-,\Lambda}_{N,k,\ve}[q;r(q)]$ and is \Holder-continuous with exponent $\beta$ and constant $C(q)$.
\end{proposition}
\begin{proof}
	Let $E$ be the smooth envelope of $F$ with $F\isom S_2/S_1$ and $S_1\subset S_2\subseteq E$.
	On the set in question, the forward and backward Lyapunov filtrations of $E$ as well as $S_1$ and $S_2$ are \Holder-continuous.
	The Lyapunov decomposition of $F$ is constructed in terms of these spaces, which , which implies the \Holder continuity of $\Lambda^{\leq}_F$.
\end{proof}

We record the following consequence:
\begin{corollary}[\Holder properties of measurable holonomy]
	\label{cor:holder_properties_of_measurable_holonomy}
	Suppose $E$ is a smooth and natural cocycle.
	Then there exists a set of full measure $Q_0$ and a constant $\beta>0$, such that for every $k\geq k_0(E)$ and $\ve\in(0,\ve_0(E)), N\geq 1$, there exist tempered functions $r,C$ depending on ${E,N,k,\ve}$, with the following property.

	For every $q\in Q_0$, $Q_0\cap Q^{+,-,\Lambda}_{N,k,\ve}[q;r(q)]$ and $x_i = \exp_q(\wtilde{x}_i)$, such that $x_2\in \cW^{s}[x_1]$ the function
	\[
		P^{-}_{q}(x_1,x_2):=\exp_{q}^{E,-1}(\wtilde{x}_2)\circ P^{-}(x_1,x_2)\circ \exp_{q}^{E}(\wtilde{x}_1,-)\colon E(q)\to E(q)
	\]
	satisfies the bound
	\[
		\norm{P^-_q(x_1,x_2)-\id}\leq C(q)d^{\cL}_q(x_1,x_2)^{\beta}.
	\]
	The analogous bound for $P^+$ holds.
\end{corollary}

\subsection{Fdd properties of stable manifolds}
	\label{ssec:fdd_properties_of_stable_manifolds}


\subsubsection{Parametrized stable manifolds}
	\label{sssec:parametrized_stable_manifolds}
	To illustrate the intended application of the notion of fdd cocycle, we now explain why the cocycles arising from normal form coordinates on stable manifolds are fdd.

	For any manifold $M$ and $m\in M$ let \index{$J$@$\cJ^{(k)}_mM$}$\cJ^{(k)}_mM$ denote the (finite-dimensional) vector space of $k$-jets of smooth functions on $M$ which vanish at $m$.
	These assemble into a smooth and natural cocycle denoted \index{$J$@$\cJ^{(k)}M$}$\cJ^{(k)}M$.
	For example $\cJ^{(1)}M$ is the cotangent bundle of $M$, and $\cJ^{(k+1)}$ is an extension of $\cJ^{(k)}$ by \index{$S$@$\Sym^{(k+1)}\cJ^{(1)}M$}$\Sym^{(k+1)}\cJ^{(1)}M$, i.e. if the $(k+1)$-jet of a function vanishing to order $k$ at $m$ determines a polynomial of degree $(k+1)$ on the tangent space at $m$.
	The jet bundles $\cJ^{(k)}M$ have subbundles $\cK^{(i,k)}M$ of functions that vanish to order $i$.
	We will use this notation both for the ambient manifold $Q$ and the stable manifolds $\cW^s$.

	By the Oseledets theorem, we have the stable subcocycle $W^s\subset TQ$ of strictly negative Lyapunov exponents.
	By definition it is fdd.
	We now explain how to construct the jet of the stable manifold $\cW^s$, and associated cocycles such as $L\cW^s$, using the allowed operations for fdd cocycles.

	Specifically, define\index{$Maps^{(k)}_q(W^s,Q)$}
	\begin{align*}
		Maps^{(k)}_q(W^s,Q):=&\left\lbrace k\text{-jets of maps from }W^s\text{ to }Q\right.\\
		 & \text{ taking }0\in W^s(q)\text{ to }q\in Q\\
		 &\left. \text{and with $1$-jet equal to the identity on }W^s(q)\right\rbrace
	\end{align*}
	In particular, we only consider maps that are already tangent to the desired tangent space, with derivative the identity map there.
	Note that $W^s$ is regarded as a linear space, so for example the jets of functions on it have natural splittings according to degree.
	The spaces $Maps^{(k)}_q$ assemble into a fiber bundle over $Q$, however it is not linear.
	To linearize it, consider the dynamically equivariant embedding
	\[
		Maps^{(k)}_q(W^s,Q)\into \Hom_{filtr}\left(
		\cJ^{(k)}Q, \cJ^{(k)}W^s
		\right)
	\]
	where the jet of a map induces a homomorphism on jet bundles by pullback of jets of functions from $Q$ to $W^s$.
	We only consider maps which respect the natural filtrations on jets by order of vanishing, i.e. those taking $\cK^{(i,k)}Q$ to $\cK^{(i,k)}W^s$.
	Note that $Maps^{(k)}(W^s,Q)$ is not a linear cocycle, but it is a subset of the corresponding linear space cut out by a collection of polynomial equations.
	We can therefore encode this fiber bundle by the ideal of functions cutting it out in the appropriate cocycle of polynomial functions on the linear space, in particular we can speak about the conclusions of \autoref{prop:subquotient_representation_of_fdd_cocycles} and \autoref{prop:holder_properties_of_fdd_cocycles}.

	Define now:
	\begin{align}
		\label{eq:stable_maps}
		Maps^{(k)}(W^s,Q)^{s}:= Maps^{(k)}(W^s,Q) \cap \Hom_{filtr}\left(
		\cJ^{(k)}Q, \cJ^{(k)}W^s
		\right)^{\leq 0}
	\end{align}
	These are maps which land in the nonpositive Lyapunov subspace of the linearizing cocycle.
	We clearly have forgetful maps by lowering the order of jets, in particular we have:
	\[
		Maps^{(k+1)}(W^s,Q)^{s} \xrightarrow{\pi_{k+1}} Maps^{(k)}(W^s,Q)^{s}.
	\]
	Note also that we are taking exponents $\leq 0$ on the right-hand side, since the identity map (when $k=1$) will be there.

\begin{proposition}[Stable parametrizations of the stable]
	\label{prop:stable_parametrizations_of_the_stable}
	Suppose that $k$ is sufficiently large depending on the Lyapunov spectrum.
	
	Then the forgetful maps $\pi_{k}$ are bijections, in particular $Maps^{(k)}(W^s,Q)^{s}$ becomes independent of $k$.

	Additionally, the action by strictly subresonant maps on $W^s(q)$ of $\bbG^{ssr}(W^s(q))$ induces an action on $Maps^{(k)}(W^s,Q)^s$, which is free and transitive.
\end{proposition}
\begin{proof}
	For $k=1$, $Maps^{(1)}(W^s,Q)$ consists of only one element, the identity map, which is in the $0$-th Lyapunov subspace and in particular $Maps^{(1)}(W^s,Q) = Maps^{(1)}(W^s,Q)^s$.
	Now suppose that $k\geq 2$.

	Because we consider \emph{filtered} maps in the homomorphisms of jets, we have a forgetful map
	\[
		\Hom_{filtr}\left(
		\cJ^{(k+1)}Q, \cJ^{(k+1)}W^s
		\right)
		\to 
		\Hom_{filtr}\left(
		\cJ^{(k)}Q, \cJ^{(k)}W^s
		\right)
	\]
	whose kernel is the linear subcocycle $\Hom(\cJ^{(k+1)}Q,\cK^{(k,k+1)}W^s)$.
	This projection is compatible with the projection $\pi_{k+1}$ on the space of maps.
	Given an element $f\in Maps^{(k)}(W^s,Q)$, any two lifts of it to $Maps^{(k+1)}(W^s,Q)$ differ by some element of $\Hom(\cJ^{1}Q,\cK^{(k,k+1)}W^s)$.
	Indeed, this space can be identified with the subspace of elements in $\Hom(\cJ^{(k+1)}Q,\cK^{(k,k+1)}W^s)$ that vanish on $\cJ^{(2,k+1)}Q$, and the identification with the fiber of $f\in Maps^{(k)}(W^s,Q)$ follows because for two lifts $\tilde{f}_1,\tilde{f}_2$, the induced pullbacks on jets are in fact algebra homomorphisms for multiplication of functions.
	In other words, if $\tilde{f}_1-\tilde{f}_2$ only involves polynomial functions of degree $k+1$ in the $W^s$-variables, then the pullback from $Q$ of any function jet vanishing at $q$ to degree $2$ or higher is the same for both $\tilde{f}_1$ and $\tilde{f}_2$.
	So the difference $\tilde{f}_1-\tilde{f}_2$ is determined by where it maps the jets of order $1$ on $Q$.

	Recall that $\cK^{(k,k+1)}W^s=\Sym^{k+1}\cJ^1W^s$, i.e. $(k+1)$-jets of functions that vanish to order $k$ can be naturally identified with polynomials of degree $k+1$.
	To identify the relevant Lyapunov exponents, denote the exponents on $TQ$ by
	\[
		\mu_1 >\cdots > \mu_n> 0 > -\lambda_m>\cdots >-\lambda_1
	\]
	where $0$ might, or might not, be an exponent.
	Recall also that $\cJ^{(1)}Q\isom TQ^{\dual}$, that $\Hom(E_1,E_2)=E_1^{\dual}\otimes E_2$, and that when dualizing exponents flip sign.
	So the Lyapunov exponents of $\Hom(\cJ^{(1)}Q,\cK^{(k,k+1)}W^s)=\Hom(\cJ^{(1)}Q,\Sym^{k+1}\cJ^1W^s)$ are of the form
	\begin{align}
		\label{eqn:resonance_equation_stable}
		\xi_{i_0} + \lambda_{i_1}+\cdots + \lambda_{i_{k+1}}
	\end{align}	
	where $\xi_{i_0}$ is some exponent of $TQ$, while $\lambda_{i_l}$ are the negative of the exponents on $W^s$ (with multiplicities allowed), so they contribute \emph{positively} (with indexing of exponents as above, we have $i_{k}\leq m, \forall k$).
	The only way for the above expression to yield a non-positive exponent is thus precisely a resonance in the sense of \autoref{appendix:subresonant_linear_algebra}, i.e. if we have $\xi_{i_0}=-\lambda_{i_0}$ for some negative Lyapunov exponent, and $\lambda_{i_0}\geq \lambda_{i_1}+\cdots + \lambda_{i_{k+1}}$.

	Recall that $Maps^{(k)}(W^s,Q)^s$ was defined as the intersection of $Maps^{(k)}(W^s,Q)$ with the strictly negative exponents.
	From the above calculation of Lyapunov exponents it follows that any two lifts of $f\in Maps^{(k)}(W^s,Q)^{s}$ to $Maps^{(k+1)}(\cW^s,Q)^{s}$ differ by a pre-composition by the $k+1$-jet of a strictly subresonant map.
	Furthermore, the lifts are unique as soon as there are no resonances, i.e. $(k+1)\cdot \lambda_m>\lambda_1$.
\end{proof}

Denote by
\[
	Maps(W^s,Q)^s := Maps^{(k)}(W^s,Q)^s \quad \text{ for }k\gg 1
\]
which by \autoref{prop:stable_parametrizations_of_the_stable} is finite-dimensional and well-defined as soon as $k$ is sufficiently large.
The content of the same proposition implies:

\begin{corollary}[Stable parametrizations are fdd]
	\label{cor:stable_parametrizations_are_fdd}
	The fiber bundle cocycle $Maps(W^s,Q)^s$ is forward dynamically defined (fdd).
\end{corollary}

We can now establish:
\begin{proposition}[Linearization cocycle is fdd]
	\label{prop:linearization_cocycle_is_fdd}
	The cocycle $L\cW^s$ constructed in \autoref{thm:linearization_of_stable_dynamics_single_diffeo} is fdd.
\end{proposition}
\begin{proof}
	Let $k\in \bN$ be sufficiently large depending on the Lyapunov spectrum.
	Recall that $L\cW^s(q)$ consists of the \emph{dual} of an appropriate space of functions on the stable manifold $\cW^s[q]$ through $q$.
	Therefore it suffices to exhibit the cocycle $F\cW^s$ of functions which is dual to $L\cW^s$.
	To do so, we will show it is naturally isomorphic to a quotient $\tilde{F}\cW^s/\tilde{K}\cW^s$ of two fdd subcocycles in $\cJ^{(k)}Q$, the $k$-jets of functions on $Q$.

	Specifically, we have $Maps(W^s,Q)^s\subset \Hom(\cJ^{(k)}Q,\cJ^{(k)}W^s)$ which is fdd, and which is a single orbit for the free action of strictly subresonant maps $\bbG^{ssr}(W^s)$ on $W^s$.
	Inside $\cJ^{(k)}W^s$ we have the space of functions $FW^s$ which is used to define the linearization, see \autoref{sssec:veronese_embedding_for_subresonant_spaces}, which is invariant under $\bbG^{ssr}(\cW^s)$ and is fdd.
	Therefore, we can consider $\tilde{F}W^s\subset \cJ^{(k)}Q$ to be the preimage of $FW^s$ under any element $\xi\in Maps(W^s,Q)^s$.
	This subspace is independent of the choice of $\xi$.
	Similarly, let $\tilde{K}W^s\subset \cJ^{(k)}Q$ be the kernel of any such $\xi$, again independent of $\xi$.
	So $\tilde{F}W^s$ represents $k$-jets of functions on $Q$ which restrict to the jets on stable manifolds, while $\tilde{K}W^s$ represents $k$-jets of functions that restrict to zero on stable manifolds.
	It is now clear that $\cF \cW^s = \tilde{F}W^s/\tilde{K}W^s$ and all the involved bundles are fdd by construction.
\end{proof}

\begin{remark}[Uniqueness of the jet of the image]
	\label{rmk:uniqueness_of_the_jet_of_the_image}
	In both \autoref{prop:stable_parametrizations_of_the_stable} above and \autoref{prop:jets_of_center_stable_parametrizations} below, the space of parametrizations of the center-stable, or stable manifold, is acted upon transitively by a group which acts by precomposition on the source.
	Therefore, the jet of the image manifold is well-defined and independent of the choice of parametrization.
\end{remark}

\subsubsection{A variant for center-stable manifolds}
	\label{sssec:a_variant_for_center_stable_manifolds}
We now establish the analogous results for center-stable manifolds.
Let therefore $W^{cs}\subset TQ$ denote the linear subcocycle with Lyapunov exponents $\leq 0$.
We have the subcocycle $W^{s}\subset W^{cs}$, although note that the complementary central part is not fdd, so we will not mention it.
Let \index{$Maps^{(k)}(W^{cs},Q)$}$Maps^{(k)}(W^{cs},Q)$ be the analogous set of $k$-th order parametrizations, mapping $0\mapsto q$ and tangent to the identity to order $1$, which again embeds in $\Hom_{filtr}\left(\cJ^{(k)}Q,\cJ^{k}W^{cs}\right)$, and denote by \index{$Maps^{(k)}(W^{cs},Q)^{cs}$}$Maps^{(k)}(W^{cs},Q)^{cs}$ the subset obtained by intersecting with the Oseledets subspace of exponents $\leq 0$.

\begin{proposition}[Jets of center-stable parametrizations]
	\label{prop:jets_of_center_stable_parametrizations}
	For $k\geq 2$, the fibers of the forgetful map
	\[
		Maps^{(k+1)}(W^{cs},Q)^{cs}\to Maps^{(k)}(W^{cs},Q)^{cs}
	\]
	are acted upon transitively by the group of $(k+1)$-jets of polynomial automorphisms of $W^{cs}$, that are tangent to order $k$ to the identity map, and that map $W^s$-translates to $W^s$ translates and do so by subresonant maps.
\end{proposition}
Note that since the maps of $W^{cs}$ in question preserve the origin and the foliation by $W^s$-translates, they must take $W^s$ to itself.
Furthermore, once $k$ is sufficiently large, strictly subresonant maps of $W^{s}$ which are tangent to the identity to order $k$, must equal the identity.

Note finally that the above group of transformations of $W^{cs}$ can be seen as strictly subresonant maps of $W^{cs}$, where we take the center direction to have negative, but arbitrarily small, exponent.

\begin{proof}
	Up to \autoref{eqn:resonance_equation_stable}, the proof is exactly the same as that of \autoref{prop:stable_parametrizations_of_the_stable} (with $s$ replaced by $cs$) but now we have a larger kernel.
	Consider again the quantity
	\[
		\xi_{i_0} + \lambda_{i_1}+\cdots + \lambda_{i_{k+1}}
	\]
	where $\lambda_i\geq 0$ but not necessarily strictly so.
	In order for the above expression to be $\leq 0$, we must have $\xi_{i_0}\leq 0$.
	If $\xi_{i_0}=0$, then necessarily we must have $\lambda_{i_j}=0,\forall j=1,\ldots,k+1$, and this corresponds to the reparametrizations in the center direction, which preserve a decomposition $W^{cs}=W^{s}\oplus W^{c}$ and act nontrivially only in the $W^c$-direction.
	If $\xi_{i_0}<0$, then we might have some $\lambda_{i_k}>0$, but only up to finitely many allowed subresonances, and all the other exponents must be $0$.
	This corresponds to maps which can ``shear'' $W^c$ along $W^s$.
\end{proof}



\subsection{Fdd cocycles and normal form coordinates}
	\label{ssec:fdd_cocycles_and_normal_form_coordinates}
In this section, we connect the notion of fdd cocycle with subresonant structures on stable manifolds.
In particular, we obtain an alternative approach to the construction of normal forms, see \autoref{sssec:subresonant_normal_forms_in_dynamics}.
This approach allows us to obtain quantitative estimates on the variation of different parametrizations, see \autoref{cor:integrability_with_parameters}.

Some related constructions appear in works of Feres and Melnick \cite{Feres_A-differential-geometric-view-of-normal-forms-of-contractions,Melnick_Non-stationary-smooth-geometric-structures-for-contracting-measurable}, but it seems to us that our point of view and results are different.
A basic reference for standard facts we will use on principal bundles and various flavors of connections is \cite{Sharpe1997_Differential-geometry}.

\subsubsection{Klein Geometry as $(G,X)$-structure}
	\label{sssec:klein_geometry_as_gx_structure}
Given a Lie group $G$ and a subgroup $H\subset G$, we recall that a \emph{Klein geometry} modeled on $G/H$ (see \cite[Ch.~4]{Sharpe1997_Differential-geometry}) on a manifold $M$ consists of a collection of charts to $X:=G/H$, such that the transition maps are given by elements of $G$ acting by left translation.
This is also known as a \emph{$(G,X)$-structure}.

\subsubsection{Klein geometry from differential geometry}
	\label{sssec:klein_geometry_from_differential_geometry}
The equivalent differential-geometric description is as follows (see \cite[Ch.~4-5]{Sharpe1997_Differential-geometry}), where $\frakg=\Lie G, \frakh=\Lie H$.
Assume given:
\begin{description}
	\item[Principal bundle] A manifold $P$ with surjection $P\onto M$ and a proper, free and transitive right action of $H$ on the fibers;
	The action will be denoted $p\mapsto ph$ or $p\mapsto R_h(p)$.
	\item[Cartan Connection] The data of:
	\begin{itemize}
		\item a $\frakg$-valued $1$-form $\omega\colon TP\to \frakg$ such that $\omega(p)\colon TP(p)\toisom \frakg$ is an isomorphism for all $p\in P$
	\item which is $H$-equivariant:
	\[
		\omega(ph) \circ DR_{h,*} = \Ad_{h}^{-1} \circ \omega(p)\text{ as maps } TP(p)\to \frakg, \forall p\in P, \forall h\in H
	\]
	\item and which is also $H$-compatible, i.e. for $X\in \frakh$ inducing the vector field $X_P$ on $P$, require
	\[
		\omega(p)\left(X_P(p)\right) = X \quad \forall p\in P.
	\]
	\end{itemize}
	\item[Flatness] Require $\omega$ to satisfy the Maurer--Cartan equation
	\[
		d\omega + \tfrac12 [\omega,\omega]=0.
	\]
\end{description}
Without the flatness assumption the above data is called a \emph{Cartan geometry} modeled on $G/H$.

\subsubsection{Maurer--Cartan form and Curvature}
	\label{sssec:maurer_cartan_form_and_curvature}
We recall that given a pointed manifold $(p,P)$, to give a (differentiable) map to the Lie group $G$ with $p\mapsto \id$ is equivalent to giving a differential $1$-form $\omega\colon TP\to \frakg$ satisfying the Maurer--Cartan equation $d\omega + \tfrac12 [\omega,\omega]=0$ (see \cite[Ch.~3, \S6]{Sharpe1997_Differential-geometry});
note that $P$ and $G$ need not have the same dimension for this statement.

In the setting of a Cartan connection $\omega$ as in \autoref{sssec:klein_geometry_from_differential_geometry}, its curvature $2$-form is defined as:
\[
	\Omega:=d\omega + \tfrac12 [\omega,\omega]\colon
	\Lambda^{2}(TP) \to \frakg
\]
The curvature vanishes on vertical directions \cite[Ch.~4, Cor.~3.10]{Sharpe1997_Differential-geometry}, and therefore the ``curvature function'' $K$ (see \cite[Ch.~4, Lemma~3.23]{Sharpe1997_Differential-geometry} and subsequent discussion) descends to a section of the vector bundle over $M$ associated to the $H$-representation $V:=\Hom(\Lambda^2(\frakg/\frakh),\frakg)$, namely to a section of $P\times_{H}V$.

Two elementary properties will be needed later.
First, the $2$-form $\Omega$ vanishes if and only if $K$ vanishes.
Second, an isomorphism of two Cartan geometries (i.e. fibered isomorphism of principal bundles preserving the Cartan connection) identifies the curvature functions.
This will imply that the Cartan geometries on stable manifolds modeled on subresonant structures have, in fact, vanishing curvature and are therefore Klein geometries.
The reason is that the dynamics is contracting on the leaves, and recurrent, so the norm of the curvature function has to vanish.

\subsubsection{Reconstructing the Klein geometry from the differential-geometric data}
	\label{sssec:reconstructing_the_klein_geometry_from_the_differential_geometric_data}
To reconstruct the more familiar description of the $(G,X)$-structure via charts, we note that given any $p\in P$, by flatness of $\omega$ there exists a unique (locally on $M$) map $\chi_{p}\colon P\to G$ such that $\chi_p(p)=\id$ and $\omega=\chi_{p}^{*}\omega_{MC}$ where $\omega_{MC}$ is the Maurer--Cartan form on $G$.
The assumptions imply that $\chi_p$ is $H$-equivariant, so we can descend $\chi_p$ to maps $\chi_{m}$ from small open sets on $M$ to $G/H$.
For another $p'\in P$ there exists a unique $g=g(p,p')$ such that $\chi_{p'}=g\cdot \chi_{p}$, so the charts on $M$ have the required properties.

We note that in the description of \autoref{sssec:klein_geometry_from_differential_geometry}, only the Lie algebra $\frakg$ appears and not the group $G$.
One can work with the unique simply-connected Lie group associated to $\frakg$, but this is not strictly necessary as all our arguments will be local on $M$.

\begin{remark}[Klein geometries and connections]
	\label{rmk:klein_geometries_and_connections}
	When there exists an $H$-invariant decomposition $\frakg=\frakh\oplus \fraks$ we obtain a more familiar Ehresmann connection on $P$ by declaring the horizontal planes on $TP$ to be $\omega^{-1}(\fraks)$.
	Note that $\dim \fraks=\dim M$ but that $\fraks$ need not be a Lie subalgebra of $\frakg$ (it must be, however, an $\frakh$-module).
	Recall that an Ehresmann connection is an $H$-invariant family of planes on $TP$, of the same dimension as $M$, that are transverse to the vertical planes, i.e. to the kernel of the projection to $M$.
	Equivalently, it is an $H$-equivariant and $H$-compatible map $\omega\colon TP\to \frakh$ (take the planes to be the kernel of $\omega$).

	We note also that we can pick any transverse subspace for $\frakh$ inside $\frakg$, but this will not lead to an $H$-equivariant family of planes unless the transversal is $H$-invariant.
\end{remark}

\subsubsection{Subresonant structures as Klein geometries}
	\label{sssec:subresonant_structures_as_klein_geometries}
A subresonant structure as in \autoref{def:subresonant_structure_on_a_manifold} is the same as a $(G,X)$-structure modeled on $X=G/H$, where $G$ is the group of subresonant maps of some $\bR^n$, and $H$ is the subgroup fixing the origin; see \autoref{rmk:on_subresonant_vs_strictly_subresonant_structures} on the question of subresonant vs strictly subresonant structure.
The above discussion shows how to give a differential-geometric description of the same data.
Our goal is to explain how the data in \autoref{sssec:klein_geometry_from_differential_geometry} naturally arises from dynamical considerations.

Fix a manifold $\cW$ and a vector space $W$ each of dimension $n$, and equip $W$ with a subresonant structure (i.e. filtration and Lyapunov exponents).
We will use $W$ as our model space.
Set $\bbG:=\bbG^{sr}(W)$ to be the group of subresonant maps of $W$ (including translations) and let $\bbG_0\subset\bbG$ be the subgroup of points fixing the origin, so we have that $W\isom \bbG/\bbG_0$.
Set $\frakg$ and $\frakg_0$ to be the corresponding Lie algebras.

Using the vector space structure on $W$ we could define $\frakt\subset \frakg$ to be the Lie subalgebra generated by translations.
It is transverse to $\frakg_0$, but we emphasize that the decomposition $\frakg = \frakg_0\oplus \frakt$ is \emph{not} $\frakg_0$-invariant in general and for this reason we will avoid it.

\subsubsection{Jets of changes of coordinates}
	\label{sssec:jets_of_changes_of_coordinates}
Denote by $\GL^{(r)}(W)$ the group of jets of all polynomial maps of $W$ that fix $\mathbf{0}$, modulo monomials vanishing to order $r+1$ or higher.
For example $\GL^{(1)}(W)=\GL(W)$ 
We will use the abbreviation $\GL_{n}^{(r)}:=\GL^{(r)}(W)$ for the group and $\frakgl_n^{(r)}$ for its Lie algebra.
For each $r\geq 2$ we have an exact sequence:
\[
	0 \to \bbK^{(r)}_n \into \GL^{(r)}_n \onto \GL^{(r-1)}_n\to 1
	\text{ with }
	\bbK^{(r)}_n\isom
	\Hom\left(\Sym^{r}W, W\right).
\]

\subsubsection{Jets of vector fields}
	\label{sssec:jets_of_vector_fields}
More generally, denote by $\frakaff^{(r)}_n$ the vector space of jets of vector fields with polynomial coefficients modulo polynomials of degree $\geq r+1$.
Note that when a general element of $\frakaff^{(r)}_n$ acts on a jet of functions of order $s$, the result is a jet of function of order $s-1$.
Concretely, differentiating by a constant coefficient vector field reduces the order to which the jet of a function is well-defined by $1$.
For this reason $\frakaff^{(r)}_n$ does \emph{not} form a Lie algebra under commutator.

We do have the inclusion $\frakgl_{n}^{(r)}\subset \frakaff^{(r)}_n$ and this is a Lie algebra under commutator of vector fields.
Furthermore, $\frakaff^{(r)}_{n}$ is a module over $\frakgl_{n}^{(r)}$, with quotient $\frakaff^{(r)}_{n}/\frakgl^{(r)}_n\isom W$.

\subsubsection{Jet bundles and formal frames}
	\label{sssec:jet_bundles_and_formal_frames_and_diffeomorphisms}
For $r\in \bN$ let $\cF^{(r)}\cW$ denote the bundle of $r$-frames on $\cW$, i.e. the fiber at $x\in \cW$ consists of $r$-jets of smooth diffeomorphisms $\chi\colon(\mathbf{0},W)\to (x,\cW)$.
Then $\cF^{(r)}\cW$ is a right principal $\GL^{(r)}_n$-bundle over $\cW$.

When there is no risk of confusion, we will write $\pi^{r+1}_{r}$ for the projection maps $\cF^{(r+1)}\cW\onto \cF^{(r)}\cW$ and $\GL_{n}^{(r+1)}\onto \GL_{n}^{(r)}$.

\subsubsection{Holonomic sections and planes}
	\label{sssec:holonomic_sections_and_planes}
We would like to now address a subtle point regarding the forgetful map $\cF^{(r+1)}\cW\to \cF^{(r)}\cW$.
Let $\cJ^1\left(\cF^{(r)}\cW;\cW\right)$ denote the space of $1$-jets of sections $\cW\to \cF^{(r)}\cW$, which is the same as tangent $n$-planes of $\cF^{(r)}\cW$ that are transverse to the vertical ones (i.e. to kernel of the projection).
It is stated in the third paragraph of \cite[\S5.1]{Feres_A-differential-geometric-view-of-normal-forms-of-contractions}, as well as the second paragraph of \cite[\S1.1.3]{Benoist1997_Orbites-des-structures-rigides-dapres-M.-Gromov} that there exists a ``canonical'' map
\begin{align}
	\label{eqn:jet_splitting_by_translations}
	\cF^{(r+1)}\cW\into \cJ^1\left(\cF^{(r)}\cW/\cW\right)
\end{align}
that assigns to the $(r+1)$-jet of chart $\chi\colon (\mathbf{0},W)\to (x,\cW)$ the tangent $n$-plane of the family of charts $y\mapsto \chi_{y}$ for $y\in \cW$ where $\chi_{y}=\chi\circ t_{-\chi^{-1}(y)}$ and $t_w$ is translation by $w$ in $W$.

We recall that the group $\GL^{(r)}_n$ acts on $\cF^{(r)}\cW$ and therefore on $\cJ^1\left(\cF^{(r)}\cW;\cW\right)$ since the fibers of the projection to $\cW$ are preserved.
The map in \autoref{eqn:jet_splitting_by_translations} is injective, therefore it \emph{cannot} be equivariant for the homomorphism $\GL^{(r+1)}_n\to \GL^{(r)}_n$ and corresponding group actions.
This boils down to the fact that translations are not in any way ``canonical'' on a smooth manifold; this persists even after we reduce the structure to a subresonant one (with some exceptions when the subresonant structure is equivalent to an affine structure).

\subsubsection{Equivariant maps}
	\label{sssec:equivariant_maps}
The above lack of equivariance can be corrected by considering instead the linear map
\begin{align}
	\label{eqn:equivariant_jet_budle_map}
	\begin{split}
	\eta\colon \cF^{(r+1)}\cW \times \frakaff_n^{(r)} & \to T\cF^{(r)}\cW\\
	(\chi,X) & \mapsto \tfrac{d}{dt}\vert_{t=0}\left(\chi\circ \phi_{t\wtilde{X}}\right)
	\end{split}
\end{align}
where $\wtilde{X}$ is a vector field representing the jet $X$ and $\phi_{t\wtilde{X}}$ is the flow generated by $\wtilde{X}$.
We emphasize that the map $\eta$ involves $\frakaff_n^{(r)}$ and not $\frakaff^{(r+1)}_{n}$ since differentiating along vector fields with vanishing $r$-jet will result in a vanishing tangent vector.

The map $\eta$ is now equivariant for $\pi^{r+1}_{r}\colon \GL_{n}^{(r+1)}\onto \GL_{n}^{(r)}$, namely
\[
	\eta\left(\chi\circ g, \Ad_{g}^{-1}X\right) = \eta\left(\chi,X\right)\circ \pi^{r+1}_r(g) \quad \forall g\in \GL_{n}^{(r+1)} 	
\] 
and we could quotient by the kernel $\bbK^{(r+1)}_n$ to obtain a $\GL_n^{(r)}$-equivariant embedding of vector bundles over $\cF^{(r)}\cW$.

Most importantly, we note that for $X\in \frakgl_{n}^{(r)}\subset \frakaff^{(r)}_{n}$ we do have
\[
	\eta\left(\chi, X\right) = X\in T\cF^{(r)}\cW(\pi^{r+1}_r(\chi))
\]
where on the right-hand side $X$ denotes the vertical vector field induced by the action of $\GL_n^{(r)}$ on the fibers.

\subsubsection{Sections and connections}
	\label{sssec:sections_and_connections}
Given a section $\sigma\colon \cF^{(r)}\cW\to \cF^{(r+1)}\cW$ we obtain an isomorphism
\[
	\eta_\sigma:=\eta\circ \left(\sigma\times \id\right)
	\colon
	\cF^{(r)}\cW\times \frakaff^{(r)}_{n}\toisom T\cF^{(r)}\cW
\]
with the property that $\eta_\sigma(\chi,X)=X\in T\cF^{(r)}\cW(\chi)$ for $X\in \frakgl_n^{(r)}$.
We will use this isomorphism to construct the Cartan connection in \autoref{prop:data_for_subresonant_structure}.

\subsubsection{Subresonant structures}
	\label{sssec:subresonant_structures_dg_notation}
We can now return to subresonant structures.
Let $\bbG^{(r)}_0\subset \GL_{n}^{(r)}$ be the image of the subgroup $\bbG_0$ of subresonant maps that fix $\mathbf{0}$ given by taking $r$-jets.
For some $r_0$ depending on the Lyapunov spectrum, if $r\geq r_0$ the forgetful map $\bbG_0\to \bbG^{(r)}_0$ is an isomorphism and the projection maps $\bbG^{(r+1)}\to \bbG^{(r)}$ are isomorphisms.
We always assume from now on that $r\geq r_0$.

At the level of Lie algebras, we also have $\Lie \bbG=\frakg \into \frakaff^{(r)}_n$ with the property that $\frakg_0=\frakg\cap\frakgl^{(r)}_n$.
With this remark we can now proceed to a basic construction:

\begin{proposition}[Data for subresonant structure]
	\label{prop:data_for_subresonant_structure}
	\leavevmode
	\begin{enumerate}
		\item Suppose that $\cW$ is equipped with a reduction of the structure group of $\cF^{(r+1)}\cW$ from $\GL_n^{(r+1)}$ to $\bbG_0^{(r+1)}$, i.e. we have a section $\delta\colon \cW\to \cF^{(r+1)}\cW/\bbG_0^{(r+1)}$.
		Then the map $\eta$ from \autoref{eqn:equivariant_jet_budle_map} induces on $\cW$ a Cartan geometry in the sense of \autoref{sssec:klein_geometry_from_differential_geometry} modeled on $\bbG/\bbG_0$.
		\item Suppose that the associated Cartan connection has vanishing curvature.
		Then this endows $\cW$ with a Klein geometry modeled on $\bbG/\bbG_0$, i.e. a subresonant structure.
	\end{enumerate}
\end{proposition}
\begin{proof}
	Vanishing curvature implies that the Cartan geometry is in fact a Klein geometry by \autoref{sssec:maurer_cartan_form_and_curvature}, so it remains to explain the construction of the connection.

	Denote by $P^{(r+1)}\subset \cF^{(r+1)}\cW$ the principal $\bbG^{(r+1)}_0\isom \bbG_0$-bundle associated to the reduction.
	Since $\bbG^{(r+1)}\isom \bbG^{(r)}$ and the isomorphism is established by the natural projection, we set $P^{(r)}:=\pi^{r+1}_r\left(P^{(r+1)}\right)$ to be the image of the reduction, which is still a principal $\bbG_0$-bundle since $r\geq r_0$.
	We then have an inverse isomorphism $\sigma\colon P^{(r)}\toisom P^{(r+1)}$.

	Recall that $\frakg\subset \frakaff^{(r)}_n$, so we can now plug in this isomorphism into the construction from \autoref{sssec:sections_and_connections} to obtain the Cartan connection:
	\begin{align*}
		\omega^{-1}\colon P^{(r)}\times \frakg & \toisom TP^{(r)} \\
		(p, Y) & \mapsto \eta\left(\sigma(p), Y\right)
		\quad \forall p\in P^{(r)}, Y\in \frakg
	\end{align*}
	Note that the Cartan connection $\omega$ is inverse of the above isomorphism.
	The required $\bbG_0$-equivariance follows by construction, as does the $\bbG_0$-compatibility, i.e. if $Y\in \frakg_0$ then $\omega^{-1}(p,Y)=Y_{P}(p)$ for the corresponding vector field $Y_{P}$ on $P^{(r)}$.	
\end{proof}

For applications in \autoref{ssec:estimates_for_linearized_unstables}, we will need a version of the above construction in families:
\begin{corollary}[Integrability with parameters]
	\label{cor:integrability_with_parameters}
	Suppose $x_0\in \cW$ is fixed, and for some metric space $Y$ we have a family of sections
	\[
		\delta_{y}\colon \cW \to \cF^{(r+1)}\cW/\bbG_0 \quad \forall y\in Y
	\]
	with the property that $\delta_y(x_0)=\delta_{y'}(x_0)$ for all $y,y'\in Y$ and the Cartan connections $\omega_y$ provided by \autoref{prop:data_for_subresonant_structure} are all flat.

	Fix $p_0\in \delta_y(x_0)$ and denote by $P_y\subset \cF^{(r)}\cW$ the principal $\bbG_0$-bundles determined by $\delta_y$ (after projecting from $\cF^{(r+1)}$ to $\cF^{(r)}$).
	Let
	\[
	 	\chi_y\colon \cW \to \bbG/\bbG_0
 	\]
 	be the family of charts as in \autoref{sssec:reconstructing_the_klein_geometry_from_the_differential_geometric_data} coming from the maps $\tilde{\chi}_{y}\colon P_y\to \bbG$, unique with the property that $\tilde{\chi}_y(p_0)=\id$.

 	Then the family $\chi_y$ has the same regularity as $\delta_y$, locally on $\cW$.
 	Namely, if $\delta_y$ is a \Holder in $y$ family of smooth sections, then $\chi_y$ is a \Holder in $y$ family of smooth maps, with same exponent and constant within a uniform multiplicative factor that only depends on a fixed smooth chart at $x_0$.
 	More generally, the analogous result holds if $Y$ is a smooth manifold and $\delta_y$ has a $C^k$-bound.
\end{corollary}

\begin{proof}
	The result is a consequence of the standard existence and uniqueness theorem for ODEs with parameter and the description of charts in \autoref{sssec:reconstructing_the_klein_geometry_from_the_differential_geometric_data}, see also \autoref{sssec:maurer_cartan_form_and_curvature}.
	We recall (see \cite[Ch.~4, Thm.~6.1]{Sharpe1997_Differential-geometry} for details) that the maps $\tilde{\chi}_{p}\colon P\to G$ are constructed by taking on $P\times G$ the $\frakg$-valued $1$-form $\phi:=\pi_{\cW}^{*}\omega_P-\pi^{*}_G\omega_G$ (where $\omega_G$ is the Maurer--Cartan form on $G$ and $\omega_P$ is the Cartan connection on $P$).
	Then the kernel of $\phi$ determines a plane field, and the flatness (i.e. Maurer--Cartan equations) imply by the Frobenius theorem that it is integrable.
	This construction works with parameters, as long as the initial conditions are fixed.
\end{proof}

\subsubsection{Subresonant normal forms in dynamics}
	\label{sssec:subresonant_normal_forms_in_dynamics}
We note that the above formalism provides an alternative proof of the existence of normal form coordinates as in \cite[Thm.~2.3-2.5]{KalininSadovskaya2017_Normal-forms-for-non-uniform-contractions}, as well as their compatibility and preservation by centralizers.

Indeed, the constructions in \autoref{ssec:fdd_properties_of_stable_manifolds} show that stable manifolds admit reductions of their higher frame bundles to the group of subresonant maps, and do so via fdd cocycles that are in particular smooth along stables \emph{and} independent of the choice of point in a full measure set on a particular unstable.
This provides at least a Cartan geometry modeled on $\bbG/\bbG_0$, which is preserved by the dynamics as well as any map that commutes with it, by construction.
However the curvature vanishes (see \autoref{sssec:maurer_cartan_form_and_curvature}) since it is clearly invariant under the dynamics and the dynamics is contracting on the stable manifolds.

\begin{remark}[On subresonant vs. strictly subresonant structures]
	\label{rmk:on_subresonant_vs_strictly_subresonant_structures}
	We could have alternatively worked with the groups $\bbG':=\bbG^{ssr}(W)\subset \bbG$ of strictly subresonant maps and the stabilizer of the origin $\bbG'_0$.
	Indeed $\bbG'/\bbG'_0\isom W \isom \bbG/\bbG_0$ and $\bbG'\normal \bbG, \bbG'_0\normal\bbG_0$ are normal subgroups.

	We also recall that an automorphism of a Klein geometry $(G,X:=G/H)$ is a group automorphism $\phi\colon G\to G$ such that $\phi(H)=H$ and we have the induced map on $X$.
	We note that conjugation by elements of $H$ gives such automorphisms, but in general there can be others.
	Let us write $W':=\bbG'/\bbG_0$ and $W:=\bbG/\bbG_0$, and there is an identification $W'\isom W$.
	Then \autoref{cor:polynomial_maps_between_different_subresonant_manifolds} implies that,
	at least when restricted to automorphisms that are represented by smooth functions, the group of automorphisms of the Klein geometries $(\bbG,W)$ and $(\bbG',W')$ agree and are equal to $\bbG_0=\bbG^{sr}_0(W)$.

	We do note that reduction of the frame bundle to $\bbG'_0$ gives a reduction to $\bbG_0$, but not conversely.
\end{remark}



\subsection{Estimates for linearized unstables}
	\label{ssec:estimates_for_linearized_unstables}

\subsubsection{Setup}
	\label{sssec:setup_estimates_for_linearized_unstables}
In this subsection, to be consistent with the intended application in the main text, we switch to unstable manifolds, backwards-regular Pesin sets, and backwards dynamically defined (bdd) cocycles.

Let \index{$L_x$}$L_x$ denote the linearization map $\cW^u[x] \to L\cW^u(x)$. 

\begin{lemma}[BiLipschitz bound on linearization]
	\label{lemma:L:is:bilipshitz}
	There exists a set of full measure $Q_0$, as well as constants $k_0,\ve_0$ such that if $k\geq k_0,\ve\in(0,\ve_0)$, there exists a tempered function $C=C_{k,\ve}$ with the following property.
	Suppose that $y\in Q_0$ and $x,x'\in \cW^{u}_{k,\ve}[y]$.
	Then
	\[
		C(q)^{-1}d^{Q}(x,x') \leq
		\norm{L_y(x)-L_y(x')}
		\leq C(q) d^{Q}(x,x').
	\]
\end{lemma}

\begin{proof}
The parametrization $W^{u}_{k,\ve}(y)\to \cW^{u}_{k,\ve}[y]$ from \autoref{prop:stableman} is biLipschitz by construction, uniformly so in the Lyapunov metric and with a tempered constant for $d^Q$.
Furthermore, for normal forms parametrizations $Maps(W^s,Q)^{s}$, we can measurably choose a section, as in \autoref{thm:linearization_of_stable_dynamics_single_diffeo}, for example by removing the super-resonant terms from the Taylor expansion of the parametrization provided by \autoref{prop:stableman}.

Now the estimate, with $x,x'\in W^u_{k,\ve}(y)$ and $d^Q(x,x')$ replaced by $\norm{x-x'}$ is clear.
Since the choice of parametrizations is biLipschitz, up to tempered multiplicative constants, the claim follows.
\end{proof}

For the following statements, recall that $I^{L\cW^u}_{y,z}\colon L\cW^u(y) \to L\cW^u(z)$ denoted the map constructed in
\autoref{eqn:trivializations_fdd} for the bdd cocycle $L\cW^u$.
We will make use of its inverse, which is usually \emph{distinct} from $I^{L\cW^u}_{z,y}$ (which might not even be defined).
We note that $\|I^{LW}_{y,z}\|$ is uniformly bounded on Pesin sets by \autoref{prop:lyapunov_adapted_trivializations}.

\begin{lemma}[Common parametrization of linearizations]
	\label{lem:common_parametrization_of_linearizations}
	There exists a set of full measure $Q_0$, as well as constants $\beta > 0$, and $k_0,\ve_0>0$ with the following properties.

	For any $k\geq k_0,\ve\in(0,\ve_0)$, there exists tempered functions $C=C_{N,k,\ve}, r = r_{N,k,\ve}$ such that for every $y\in Q_0$, there exists a function
	\[
		h_y\colon Q^{-,\Lambda}_{N,k,\ve}[y;r(y)] \to \GL\left(L\cW^u(y)\right)
	\]
	satisfying, for $z\in Q^{-,\Lambda}_{N,k,\ve}[y;r(y)]$:
	\begin{align*}
			h_{y}(z)\cdot
	        \left(I^{L\cW^u}_{y,z}\right)^{-1}
	        \left(L_{z}(\cW^u[z])\right) & =
	        L_{y}\left(\cW^u[y]\right) && \text{on linearized unstables}\\
	        h_{y}(z)\cdot
	        \left(I^{L\cW^u}_{y,z}\right)^{-1}
	        \left(L\cW^u(z)^{\geq \bullet}\right)
	        &=L\cW^u(y)^{\geq \bullet} && \text{on flags}
	\end{align*}
	as well as $h_y(y)=\id$ and
	\begin{displaymath}
	\norm{h_y(z_1)-h_y(z_2)}\leq C(y) d^Q(z_1,z_2)^\beta
	\quad \forall z_1,z_2\in Q^{-,\Lambda}_{N,k,\ve}[y;r(y)].
	\end{displaymath}
\end{lemma}

Note that because $z\in Q^{-,\Lambda}_{N,k,\ve}$, it means, in particular, that there exist subresonant isomorphisms between $\cW^{u}[z]$ and $\cW^{u}[y]$.
Indeed, the fact that the cocycle $L\cW^{u}$ is defined at $z$ implies, in particular, that so is the unstable manifold, and that the Lyapunov spectrum agrees with that of the background measure.

Furthermore, by \autoref{cor:polynomial_maps_between_different_subresonant_manifolds}, the map $h_{y}(z)\circ \left(I^{L\cW^u}_{y,z}\right)^{-1}$ is the linearization of a subresonant isomorphism $\cW^{u}[z]\to \cW^{u}[y]$.

\begin{proof}
	The existence of an $h_y(z)$, without estimates on the $z$-dependence, follows from the trivializations constructed in \autoref{prop:agreement_of_linearizations_and_compatible_family_of_subgroups}.
	
	To construct $h_y(z)$ satisfying a \Holder estimate, we note first that \index{$M$@$\cM(x)$}$\cM(x):=L_x(\cW^u[x])\subset L\cW^u[x]$ is a Zariski-closed algebraic subset, with a transitive action of the unipotent algebraic group $\bbG^{ssr}(x)$; the Zariski-closed property holds since all orbits of unipotent groups are Zariski-closed.
	The construction in \autoref{prop:stable_parametrizations_of_the_stable} shows that $\cM\to \cW^u[Q_0]$ is bdd in the sense of \autoref{def:forward_dynamically_defined_cocycle}.
	Therefore, if we let $\cJ_{\leq d}(x)$ denote the polynomials of degree at most $d$ that vanish on $\cM(x)$, we have that the induced $I_{y,z}^{-1}\left(\cJ_{\leq d}(z)\right)$ and $\cJ_{\leq d}(y)$ are \Holder-close in the corresponding Grassmannian.
	The same applies to the Lyapunov flags.
	We also know that, in the product of the corresponding Grassmannians, the subspaces are in the same $\GL(L\cW^u(y))$-orbit, and for large enough $d$ the subspace $\cJ_{\leq d}(x)$ determines $\cM(x)$ (and the ideal of all polynomials vanishing on it).

	We next recall that for an algebraic group $\bbG$ acting on a variety $\bbV$ and point $p\in \bbV(\bR)$, the orbit $\bbO:=\bbG\cdot p$ is locally closed in $\bbV$ in the Zariski topology.
	Let $\bbS\subset \bbG$ be the stabilizer of $p$, and let $\fraks\subset \frakg$ denote the corresponding Lie algebras, and let $\fraks^{\perp}$ denote the orthogonal complement of $\fraks$ in $\frakg$; note that we apply this with $\bbG=\GL_N$ so we can use the natural Killing form.
	Since $\bbG(\bR)\to \bbO(\bR)$ has surjective differential at $\id\to p$, by the implicit function theorem, for an open ball $\fraks^{\perp}_\ve\subset\fraks$ around $0$, the map $X\mapsto \exp(X)\cdot p$ is biLipschitz between $\fraks^{\perp}_{\ve}$ and an open neighborhood of $p$ in $\bbO(\bR)$.
	The existence of $h$ with the required bound now follows.
\end{proof}

\begin{proposition}[Linearization varies \Holder-continuously]
	\label{prop:linearization_varies_holder_continuously}
	There exists a set of full measure $Q_0$, such that for every $k_1$ there exist constants $\beta > 0$, $k_0,\ve_0>0$ with the following properties.

	For any $k\geq k_0,\ve\in(0,\ve_0)$, there exists tempered functions $C=C_{N,k,\ve}$, $ r = r_{N,k,\ve}$, $r_0=r_{0,N,k,\ve}>0$ with the following property.
	For $y\in Q_0$ let $B_r(y;\cW^u)$ denote by intersection of $L_y(\cW^u[y])$ with the ball of radius $r$ around $L_y(y)$.
	Then the map:
	\begin{align*}
		\xi_y\colon Q^{-,\Lambda}_{N,k,\ve}[y;r(y)] & \to 
		C^{k_1}(B_{r_0(y)}(y;\cW^u);TQ(y))\\ 
		z & \mapsto \exp_y^{-1}\circ L_z^{-1} \circ I_{y,z}^{L\cW^u}\circ h_y^{-1}(z)
	\end{align*}
	satisfies $\norm{\xi_y(z_1)-\xi_y(z_2)}_{C^{k_1}}\leq C(y)d^Q(z_1,z_2)^{\beta}$.
\end{proposition}
Informally, the above result says that combining the canonical linearization $L_z$ together with the parametrization $h_y(z)$ from \autoref{lem:common_parametrization_of_linearizations} yields a \Holder-continuous family of maps for the $C^{k_1}$-norm.

\begin{proof}
	We consider first the parametrizations of unstable manifolds that arise from \autoref{prop:stableman}:
	\begin{align*}
		\zeta_y\colon Q^{-,\Lambda}_{N,k,\ve}[y;r(y)] & \to 
		C^{k_1}\left(W^u(y;r_0'(y));TQ(y)\right)
	\end{align*}
	where $W^u(y;r)$ is the ball of radius $r$ in $W^u(y)$ (in the $(k,\ve)$-Lyapunov norm).
	The defining property of $\zeta_y(z)$ is that $\cW^{u}_{k,\ve}[z]\cap \cL_{k,\ve}[y]$ is equal to $\exp_y$ applied to the graph of $\zeta_y(z)$.

	Then $\zeta_y$ is a \Holder function, with exponent depending on $k,k_1,\ve$, see for example \cite[\S8.1.5]{BarreiraPesin2007_Nonuniform-hyperbolicity} or \cite[Prop.~3.12]{Ben-Ovadia2018_Symbolic-dynamics-for-non-uniformly-hyperbolic-diffeomorphisms-of-compact}.
	Note that both of these references make the statement only for $k_1=1$, but the general case follows by either of the following two approaches.
	One can go through the proof of the stable manifold theorem and verify that the estimate holds (essentially the approach in \cite{Ben-Ovadia2018_Symbolic-dynamics-for-non-uniformly-hyperbolic-diffeomorphisms-of-compact}).
	Alternatively, one first notes that unstable holonomies are \Holder-continuous relative to the family of transverse coordinate planes, uniformly over Pesin sets.
	Then $\exp_{y}^{-1}(\cW^{u}_{k,\ve}[z])$ is a graph of a smooth function over $W^{u}(y)$, and the \Holder estimates for the jets of $\cW^u[z]$ from \autoref{prop:holder_properties_of_fdd_cocycles} imply the claimed \Holder estimate on $\zeta_y$.

	Set now
	\begin{align*}
		\tilde{\xi}_y\colon Q^{-,\Lambda}_{N,k,\ve}[y;r(y)] & \to 
		C^{k_1}(B_{r_0(y)}(y;\cW^u);TQ(y))\\
		z & \mapsto \zeta_y(z) \circ \phi_y  
	\end{align*}
 	where $\phi_y$ is the unique subresonant isomorphism $\phi_y\colon W^u(y)\toisom L_y(\cW^{u}[y])$ such that
 	\[
 		\exp_y\circ \tilde{\xi}_y(y)
 		=
 		L_{y}^{-1}
 		\text{ as maps }
 		L_y(\cW^u[y]) \toisom \cW^{u}[y].
 	\]
 	Then by construction $\tilde{\xi}_y$ is a \Holder map.

 	We recall also that $h_{y}(z)\circ \left(I^{L\cW^u}_{y,z}\right)^{-1}$ is induced by the linearization of a subresonant isomorphism $\cW^{u}[z]\to \cW^{u}[y]$.
 	Therefore, $\tilde{\xi}_y(z)$ differs from $\tilde{\xi}_y(z)$ as defined in the statement by precomposition with a diffeomorphism $\rho_y(z)$ of $L_y(\cW^u[y])$, such that the induced map gives the subresonant isomorphism of unstables.

 	We now note that $\tilde{\xi}_y$ induces a \Holder family of subresonant structures on $L_y(\cW^u[y])$, with parameter space $Q^{-,\Lambda}_{N,k,\ve}[y;r(y)]$, so we can apply \autoref{cor:integrability_with_parameters} to obtain \Holder variation of the reparametrizing maps.
\end{proof}

\begin{proposition}[Measuring distances in linearization]
	\label{prop:measuring_distances_in_linearization}
	With the same assumptions and notation as in \autoref{lem:common_parametrization_of_linearizations},
	there exists $\beta_2>0$ and tempered $C=C_{N,k,\ve}$ such that, for $x_z\in \cW^u_{k,\ve}[z]\cap \cL_{k,\ve}[y]$ and $x_{y}\in \cW^u_{k,\ve}[y]$, we have the estimate
	\[
		d^{\cL}_y(x_z,x_y) \leq C(y)\left(\norm{I_{y,z}^{-1}\circ L_z(x_z) - L_y(x_y)} 
				+ d^{\cL}_y(z,y)\right)^{\beta_2}
	\]
	where $d^{\cL}_{y}$ is the Lyapunov distance based at $y$ (see \autoref{sssec:lyapunov_charts_and_distance}).
\end{proposition}

\begin{proof}
	We first reduce to the case when $x_y$ has the property that $L_y(x_y) = h_{y}(z)\circ I_{y,z}^{-1}\circ L_z(x_z)$.
	Indeed, set $x_y'\in \cW^u_{k,\ve}[y]$ to be the unique point such that $L_y(x_y')=h_y(z)\circ I_{y,z}^{-1}\circ L_z(x_z)$.
	Then since linearization is locally Lipschitz we have:
	\begin{align*}
		d^{\cL}_y(x_y',x_y)& \leq  C(y)\norm{L_y(x_y') - L_y(x_y)}	\\
		& = C(y) \norm{h_y(z)\circ I_{y,z}^{-1}\circ L_z(x_z) - L_y(x_y)}\\
		& \leq C(y) \big[\norm{h_y(z)-\id} \norm{I_{y,z}^{-1}\circ L_z(x_z)} + 
		\norm{I_{y,z}^{-1}\circ L_z(x_z) - L_y(x_y)} \big]\\
		& \leq C_1(y)\left(d^{\cL}_y(y,z)^{\beta} + \norm{I_{y,z}^{-1}\circ L_z(x_z) - L_y(x_y)}\right)
	\end{align*}
	So it suffices to treat the case when $x_y=x_y'$.

	In this case, set $\tilde{x}_y:=L_y(x_y)\in L_y(\cW^u[x_y])$.
	Then, for $\xi_y$ as constructed in \autoref{prop:linearization_varies_holder_continuously}, we have that
	\begin{align*}
		x_z & = \exp_y(\xi_y(z)(\tilde{x}_y))\\
		x_y & = \exp_y(\xi_y(y)(\tilde{x}_y))
	\end{align*}
	and so
	\[
		d^{\cL}_y(x_z,x_y)\leq C(y) d_{C^0}(\xi_y(z),\xi_y(y))
		\leq
		C(y) d^{\cL}_y(y,z)^{\beta}
	\]
	where the last inequality is a consequence of \autoref{prop:linearization_varies_holder_continuously}.	
\end{proof}

\bibliographystyle{main_text_bibstyle}
\bibliography{main_text}

\printindex

\end{document}